\documentclass[11pt, eqno]{article}

\usepackage{amsmath,amssymb,latexsym}
\newcommand{\1}{{{\mathchoice {\rm 1\mskip-4mu l} {\rm 1\mskip-4mu l}
{\rm 1\mskip-4.5mu l} {\rm 1\mskip-5mu l}}}}

\newcommand{\N}{{\mathbb N}}
\newcommand{\C}{{\mathbb C}}
\newcommand{\R}{{\mathbb R}}
\newcommand{\Q}{{\mathbb Q}}
\newcommand{\Z}{{\mathbb Z}}
\newcommand{\CP}{{\mathbb CP}}

\newtheorem{theorem}{Theorem}[section]

\newtheorem{corollary}[theorem]{Corollary}

\newtheorem{definition}[theorem]{Definition}
\newtheorem{example}[theorem]{Example}
\newtheorem{remark}[theorem]{Remark}

\newtheorem{lemma}[theorem]{Lemma}
\newtheorem{claim}[theorem]{Claim}

\newtheorem{proposition}[theorem]{Proposition}
\newtheorem{assumption}[theorem]{Assumption}

\begin{document}
\title{Constructing virtual Euler cycles and classes}
\author{Guangcun Lu\thanks{The first author is partially supported by the NNSF 10371007 of China and
the Program for New Century Excellent Talents of the Education Ministry of China.}\\
Department of Mathematics,\\
Beijing Normal University, China\\
(gclu@bnu.edu.cn)\vspace{4mm}\\
Gang Tian\\
Department of Mathematics,\\
 Princeton University and Beijing University\\
(tian@math.princeton.edu)}

%\subjclass[2000]{58F05, 53C15}

\date{First version May 10, 2006\\
Revisied July 5, 2007} \maketitle \vspace{-0.1in}

\abstract{The constructions of the virtual Euler (or moduli)
cycles and their properties are explained and developed
systematically in the general abstract settings. }
 \vspace{-0.1in}
\medskip
\tableofcontents

\section{Introduction}

In many mathematical studies, one encounters the following {\bf
Moduli Problem}: $E\to X$ is a vector bundle and $S:X\to E$ is a
section, the zero locus $Z(S):=S^{-1}(0)$ contains a lot of
information of about the triple $(X, E, S)$.

{\bf A classical differential topology fact is}: If $X$ is a
closed smooth manifold of finite
 dimension and $E\to X$ is a smooth vector bundle, then a generic smooth
 section $S:X\to E$ is transversal to the zero section and
 $Z(S)$ is a closed submanifold of $X$ which is a
 representative cycle of the Poincar\'e dual of the Euler class of
 $E$.

For the general moduli problem above it is expected to get some
analogies. A nice historical discussion of work related to this
question can be found in the introduction of \cite{CMSa}. The goal
of this paper is to explain and study constructions and properties
of (virtual) Euler chains and classes of the following four kinds
models of moduli problems.

\begin{description}
\item[I.]   $(X, E, S)$ is a {\bf Banach Fredholm bundle} of index
$r$ and with compact zero locus. Namely, $X$ is a separable Banach
manifold, $E\to X$ is a Banach vector bundle and $S$ is a Fredholm
section of index $r$ and with compact zero locus. When the
determinant bundle ${\rm det}(S)\to Z(S)$ is oriented, i.e., it is
trivializable and is given a continuous section nowhere zero, $(X,
E, S)$ is said to be {\bf oriented}.

\item[II.]   $(X, E, S)$ is a {\bf Banach Fredholm orbibundle} of
index $r$ and with compact zero locus. Namely, $X$ is a separable
Banach orbifold, $E\to X$ is a Banach orbibundle and $S$ is a
Fredholm section of index $r$ and with compact zero locus.   $(X, E,
S)$ is called {\bf oriented} if the determinant bundle ${\rm
det}(S)\to Z(S)$  is oriented, (cf. \S2.1 for precise definition).

\item[III.]  Roughly speaking, $(X, E, S)$ consists of
a separable PS (partially smooth) Banach manifold
(Def.\ref{def:3.1}), a PS Banach bundle  $E\to X$
(Def.\ref{def:3.4}) and a PS section $S:X\to E$ which has a compact
zero locus $Z(S)$ and whose restriction to each stratum is Fredholm
and  restrictions to the lower strata have less indexes. (See
Def.\ref{def:3.5} for precise definition). The bundle $E$ is also
required to have a class of PS sections ${\cal A}$ which is rich
near $Z(S)$ (Def.\ref{def:3.8}) and $S$ is $k$-good relative to the
class ${\cal A}$ (Def.\ref{def:3.13}).  $(X, E, S)$ is said to be
{\bf oriented} if the restriction of the determinant bundle ${\rm
det}(S)\to Z(S)$ to the top strata $Z(S)\cap X_0$ is oriented.

\item[IV.]   $(X, E, S)$ is, roughly speaking, consisting of
 a PS Banach orbibundle  $E\to X$ over a separable PS Banach
 orbifold $X$ and a PS stratawise Fredholm section $S:X\to E$
 with similar properties as in III.
  \end{description}

Clearly, the Model IV includes the first three models as a special
case. For each model we shall give a detailed discussions since
their precise forms can be more suitable for different problems and
we  can get very refined results. For Model II and IV we shall adopt
the method developed by Liu-Tian in \cite{LiuT1}-\cite{LiuT2}. The
study of some properties is motivated by \cite{Lu2}-\cite{Lu3}.
Related works were done by Hofer et al in the language of polyfolds
\cite{Ho}, by Cieliebak-Mundet i Riera-Salamon for an oriented
$G$-equivariant Fredholm section of a Hilbert space bundle over a
Hilbert manifold \cite{CMSa}\footnote{Here $G$ is a compact oriented
Lie group such that the isotropy group of $G$ at each point of the
zero locus is finite},and by McDuff in the language of groupoids
\cite{Mc3}. In \cite{Mc3}, McDuff proposed an intrinsic definition
of a weighted branched manifold and used it to construct the virtual
moduli cycle in a finite dimensional setting. She also adapted the
method of Liu-Tian and expected to carry this construction to
certain infinite dimensional setting which covers at least
Gromov-Witten invariants. It will be interesting to compare the
construction here to hers.

Section 1 will explain the constructions of the Euler classes of
Banach Fredholm bundles though many of them are known. We here
present them so that to have a good understanding and guide for
other complex cases. The readers only need to skim over it.
 For a Banach Fredholm bundle $(X, E, S)$ of index $r$ as above,
  one can give it a small perturbation section $\sigma$ so
that the section $S+\sigma$ is still a Fredholm one with index $r$
and is transversal to the zero section. Then $Z(S+\sigma)$ is a
compact manifold of dimension $r$ and is called a {\bf Euler chain}
of $(X,E,S)$. If $Z(S+\sigma)$ has no part of codimension one then
it is a cycle in $X$ of dimension $r$, determines a homology class
in $H_r(X,\Z_2)$ (resp. $H_r(X, \Z)$ if $S$ is orientable), denoted
by $e(E,S)$ and called a {\bf Euler cycle} of $(X,E,S)$. In this
case its homology is independent of such a generic small
perturbation section $\sigma$, called a {\bf Euler class} of the
triple $(X,E,S)$. This classical result and properties of Euler
classes have a very important position in the global analysis,
geometry and topology.

In the second section we shall consider the case of the Model II.
That is, for Banach Fredholm orbibundles we develop a
corresponding construction theory with Section 1. However, it is
not hard to see that the usual arguments for Banach Fredholm
bundles cannot be directly applied to the Banach Fredholm
orbibundle $(X,E,S)$ though one can prove that the zero set $Z(S)$
is still an orbifold when $S$ is transversal to the zero section.
In most cases, $S$ is not transversal to the zero section,
moreover, one cannot add an arbitrarily small section $\sigma$ to
$S$ to get a section $S+\sigma$ which is transversal to the zero
section. For example, let $\widetilde X$ be a closed smooth
manifold of finite
 dimension, and  $\Gamma$ a finite group which has not only an automorphism
 representation on $\widetilde X$ but also a linear one on $\R^k$. Assume that
 $\gamma$ acts on $\widetilde X$ effectively. Denote by
 $X=\widetilde X/\Gamma$. Then the natural projection
 $p:E=(\widetilde X\times\R^k)/\Gamma\to \widetilde X/\Gamma$
 is an orbibundle of rank $k$ over orbifold $X=\widetilde X/\Gamma$.
A smooth section $S$ of $p$ may be identified with  an equivariant
smooth section $\tilde S:\widetilde X\to \widetilde X\times\R^k$,
i.e. $\tilde S(g\cdot x)=g\cdot\tilde S(x)$ for any $x\in
\widetilde X$ and $g\in\Gamma$. If such a section $\tilde S:\tilde
S:\widetilde X\to \widetilde X\times\R^k$ is not transversal to
the zero section it follows from Theorem~\ref{th:1.5}(B) that
there exists a sufficiently small section $\tilde
\sigma:\widetilde X\to\widetilde X\times\R^k$ such that $\tilde
S+\tilde \sigma$ is transversal to the zero section. But we can
not guarantee that the section $\tilde \sigma$ and thus $\tilde
S+\tilde\sigma$ is equivariant. So $\tilde S+\tilde\sigma$ can not
descend to a single value section of the orbibundle $E\to X$ in
general. So to deal with the Banach Fredholm bundles one needs to
develop new methods. See  \cite{FuO}, \cite{LiT},
\cite{LiuT1}-\cite{LiuT3},
 \cite{R2} and \cite{Sie} for some previous works. These works concern
 the important case of Gromov-Witten invariants. We are more interested in
 constructing virtual classes in the most general cases. More precisely,
 {\bf for an oriented Banach Fredholm orbibundle $(X,E,S)$ with
compact zero locus, if $X$ satisfies Assumption~\ref{ass:2.45},
i.e., each $x\in X^{sing}\cap Z(S)$ has the singularity of
codimension at least two, we shall construct its (virtual rational)
Euler class.} These will be completed in
Sections~\ref{sec:2.3}-\ref{sec:2.8} after we give an overall
strategy in a special case in Section~\ref{sec:2.2}.

\S 3 deals with  the case of the Model III.  We abstract two
essential notions in Definition~\ref{def:3.10} and
Definition~\ref{def:3.12} and use them to develop the
corresponding theory with \S 1.

 Finally, in \S4 we
combine the methods in \S 2 and \S 3 to construct the virtual
Euler chains and classes for  the PS Banach Fredholm orbibundles
in the Model IV, and study their properties. Our settings are
motivated by constructing the Floer homology and Gromov-Witten
invariants on general\\ symplectic manifolds, but more general. We
firstly discuss the general constructions, and then divide into
two concrete frameworks. One of them, Framework I, is designed for
the closed string Gromov-Witten invariants, under the
Assumption~\ref{ass:4.7} we construct the virtual Euler chains and
classes and also get all corresponding results with those in \S2.
Framework II is designed for the Floer homology and the open
string Gromov-Witten invariants, Theorem~\ref{th:4.20}  is enough
for the Floer homology.  In actual applications to them the
goodness conditions relative to ${\cal A}$ can be checked by the
gluing arguments. Our attention is how to complete global
constructions under the least local assumptions.

In this paper we shall {\bf always assume}: our Banach manifolds
(resp. orbifolds) are separable if we need to use Sard-Smale
theorem, and admit smooth cut-off functions if we need them (for
example, Hilbert manifolds (resp. orbifolds) and those Banach
manifolds (resp. orbifolds) modelled and Sobolev spaces $L^p_k$
with $p$ an even integer).

{\bf A guide for the reader}. The second section is the core of this
paper.  In that section, we explain ideas and methods of overcoming
the difficulties of orbifolds in our construction in details as much
as possible. In particular, the reader should be able to get a clear
overall strategy for the construction. Having a good understanding
of Section 2, one can easily understand Section 4 in which the
related constructions are presented in a more general category
including one for GW-invariants and Floer homologies.

{\bf Acknowledgements}.\quad The first author would like to thank
the Abdus Salam International Centre for Theoretical Physics and the
Institut des Hautes \'Etudes Scientifiques for their hospitality and
financial support where much of the work was done. We also warmly
thank the referee for chasing an uncountable number of mistakes, and
suggesting many improvements.

\section{The Euler cycle of Banach Fredholm bundles}\label{sec:1}

We first construct the Euler class of a Banach Fredholm bundle with
compact zero locus, and then give  two localization formulas and
some properties. They are useful for us constructing and
understanding the virtual Euler chain (or class). The basic
references of our arguments are [AS], [Bru], \S4.3 in [DoKr], [LiT]
and [R].  We here adopt a way that is closely related to our
generalization in the following sections though the method appears
to be more complex.

\subsection{The construction of the Euler cycles}\label{sec:1.1}

 Let $p:E\to X$ be a Banach vector
bundle over  a Banach manifold  $X$ with corners, and $S:X\to E$ be
a $C^1$-smooth section  with {\bf zero locus} $Z(S):=\{x\in X\,|\,
S(x)=0\}$. If $Z(S)\ne\emptyset$, for each $x\in Z(S)$ let
$DS(x):T_xX\to E_x$ is the composition of the projection $\Pi(x):
T_{(x,0)}E\to E_x$ and the differential $dS(x):T_xX\to T_{(x,0)}E$.
We call $DS(x)$ the {\bf vertical differential} of $S$ at $x$. Using
a connection $\nabla$ on $E$ ( which in general may not exist on a
Banach bundle) we can extend $DS$ from $Z(S)$ to $M$ by $\nabla S$.
In the literature one often meets the following two definitions of
the Fredholm section of index $r$.

$\bullet$ For each $x\in Z(S)$, $DS(x):T_xX\to E_x$  is a Fredholm
operator of index $r$.

$\bullet$ In each local trivialization  $S$ may be represented by
a Fredholm map of index $r$ from the base to the fibre. (cf. page
137 in [DoKr].)

Clearly, the second definition seems to be stronger than the first
one, and is always well-defined whether  $Z(S)\ne\emptyset$ or not.
Later we say the section satisfying the second definition to be {\bf
strongly Fredholm}. It is not hard to see that a strongly Fredholm
section is stable under a small perturbation near a compact subset.
The following result shows that these two definitions have no
essential differences in the most cases.

\begin{lemma}\label{lem:1.1} Let $p:E\to X$ be any Banach bundle over a
Banach manifold $X$ with corners, and $S:X\to E$ be a Fredholm
section. If the zero locus $Z(S)$ is nonempty  then there exists an
open neighborhood $U(Z(S))$ of $Z(S)$ such that the restrictions
$S|_{U(Z(S))}$ is strongly Fredholm.
\end{lemma}

\noindent{\bf Proof.}\quad For $x\in Z(S)$ let $O_x$ be an open
neighborhood of $x$ in $X$ and
$$
\psi: O_x\times H\to E|_{O_x},\;(y,\xi)\mapsto\psi_y(\xi)
$$
 be a local trivialization. Here $\psi_y:H\to E_y$ is the topological linear
 isomorphism from a Banach space $H$ to $E_y$. Then $\psi^{-1}\circ
(S|_{O_x}):O_x\to O_x\times H$ may be written as the form
$\psi^{-1}\circ (S|_{O_x})(y)=(y, S_\psi(z))\,\forall y\in O_x$,
where
$$S_\psi:O_x\to H,\, y\mapsto \psi_y^{-1}(S(y))$$
 is called the {\bf local
representation} of $S$ under the trivialization $\psi$, or a local
representation of $S$ near $x$. For $x\in Z(S)$ one easily proves
that $DS(x)$ is Fredholm if and only if the differential
$dS_\psi(x): T_xX\to H$ is Fredholm for some (and thus any) local
representative $S_\psi$ of $S$ near $x$. Since the Fredholmness of
bounded linear operators between Banach spaces is stable under small
perturbations we may assume that the differential $dS_\psi(y):T_y
X\to H$ is also Fredholm and has the same index as $dS_\psi(x)$ for
each $y\in O_x$ by shrinking $O_x$ if necessary. Setting
$U(Z(s))=\cup_{x\in Z(S)}O_{x}$
 it is easy to see that the
restriction of $S$ to $U(Z(S))$ is strongly Fredholm.
Lemma~\ref{lem:1.1} is proved. \hfill$\Box$\vspace{2mm}

 Assume that $p:E\to X$ is a Hilbert vector
bundle over a Hilbert manifold. Since the tangent space of a fibre
$E_x$ at any point $\xi\in E_x$ may be identified with $E_x$,
using the orthogonal projection from $T_\xi E$ to $E_x$, denoted
by ${\bf P}_\xi$, we have a ``connection" $D$ defined by
\begin{equation}\label{e:1.1}
 DS(x):={\bf P}_{S(x)}\circ dS(x): T_x
X\to E_x.
\end{equation}
 It is still called the {\bf vertical
differential} of $S$ at $x$ without occurrence of confusions. Here
we put the word connection between quotation marks because  the
projection ${\bf P}_\xi$ can {\bf only depend on $\xi$ continuously}
if both the rank of $E$ and the dimension of $X$ are infinite. So we
can only know that {\bf $DS(x)$ is continuous in $x$} even if $S$ is
smooth. However, it is easily checked that for any local
representative $S_\psi$, $DS(x)$ and $dS_\psi(x)$ have the same
linear functional analysis properties. For example, $DS(x)$ is
Fredholm (resp. onto) if and only if $dS_\psi(x)$ is Fredholm (resp.
onto) and have the same index. These are sufficient for our many
arguments.

 A Fredholm section $S:X\to E$
is called {\bf orientable} if  the real (determinant) line bundle
${\rm det(Ind(DS))}$ is trivial over $Z(S)$, i.e., it admits a
continuous nowhere zero section on $Z(S)$. Clearly, such a given
section determines an orientation of $S$. In this case we say the
Banach Fredholm bundle $(X, E, S)$ to be {\bf oriented}.

The following two lemmas are elementary functional analysis
exercises.

\begin{lemma}\label{lem:1.2}
Let both $A$ and $B$ be two continuous linear operators from
Banach spaces $Y$ to $Z$. Suppose that $A$ is onto. Then there
exists a $\epsilon>0$ such that $B$ is also onto as
$\|A-B\|<\epsilon$.
\end{lemma}

\begin{lemma}\label{lem:1.3}
Let $X, Y$ and $Z$ be Banach spaces, and  both $A:X\to Z$ and
$B:Y\to Z$ be continuous linear operators. Define
$$A\oplus B: X\times Y\to Z,\;(x,y)\mapsto Ax+ By.$$
{\rm (i)} If $A$ is Fredholm and $A\oplus B$ is onto, then the
restriction of the natural projection $X\times Y\to Y$ to ${\rm
Ker}(A\oplus B)$,   $\Pi:{\rm Ker}(A\oplus B)\to Y$
is a Fredholm operator with ${\rm Index}(\Pi)={\rm Index}(A)$;\\
{\rm (ii)} If $\dim Y<+\infty$ then $A\oplus B$ is Fredholm if and
only if $A$ is Fredholm, and in this case ${\rm Index}(A\oplus
B)={\rm Index}(A)+ \dim Y$.
\end{lemma}

\noindent{\bf Proof.} We only prove (ii). Since $\dim Y<+\infty$ we
have a direct sum decomposition $Y=Y_0\oplus Y_1$, where
 $Y_0={\rm Ker}(B)$. Then $B|_{Y_1}:Y_1\to Z$ is an injection. Set
$Y_{10}=(B|_{Y_1})^{-1}({\rm Im}(A))$ and decompose $Y_1$ into
$Y_{10}\oplus Y_{11}$. Then $B|_{Y_1}(Y_{11})\cap{\rm Im}(A)=\{0\}$
and ${\rm Im}(A\oplus B)={\rm Im}(A)\oplus B|_{Y_1}(Y_{11})$. This
and ${\rm Ker}(A)\subset {\rm Ker}(A\oplus B)$ imply that $A$ is
Fredholm if $A\oplus B$ is so.  Conversely, if $A$ is Fredholm we
can decompose $X=X_0\oplus X_1$, where $X_0={\rm Ker}(A)$. So
\begin{eqnarray*}
{\rm Ker}(A\oplus B)\!\!\!\!\!\!\!\!\!\!&&=({\rm
Ker}(A)\times\{0\})\oplus(\{0\}\times{\rm Ker}(B))\\
&&\oplus\{(x,y)\in X_1\times Y_{10}\;|\;A|_{X_1}x+ B|_{Y_1}y=0\}.
\end{eqnarray*}
Since  $A|_{X_1}:X_1\to Y$ is injective and $B|_{Y_1}(Y_{10})\subset
{\rm Im}(A)={\rm Im}(A|_{X_1})$ we get
\begin{eqnarray*}
\dim{\rm Ker}(A\oplus B)\!\!\!\!\!\!\!\!\!\!&&=\dim {\rm Ker}(A)+
\dim{\rm Ker}(B)+ \dim
Y_{10}\\
&&=\dim {\rm Ker}(A)+ \dim{\rm Ker}(B)+ \dim Y-\dim Y_{11}-\dim Y_0\\
&&=\dim {\rm Ker}(A)+  \dim Y-\dim Y_{11}.
\end{eqnarray*}
Hence
  $\dim{\rm Ker}(A\oplus B)<+\infty$ and thus $A\oplus B$ is Fredholm
(using ${\rm Im}(A\oplus B)={\rm Im}(A)\oplus B|_{Y_1}(Y_{11})$
again). To get the final index equality we decompose $Z={\rm
Im}(A\oplus B)\oplus Z_2={\rm Im}(A)\oplus B|_{Y_1}(Y_{11})\oplus
Z_2$. Then ${\rm Coker}(A\oplus B)\cong Z_2$ and ${\rm
Coker}(A)\cong B|_{Y_1}(Y_{11})\oplus Z_2$. It follows that
\begin{eqnarray*}
{\rm Ind}(A\oplus B)\!\!\!\!\!\!\!\!\!\!&&=\dim{\rm Ker}(A\oplus
B)-\dim{\rm Coker}(A\oplus B)\\
&&=\dim {\rm Ker}(A)+  \dim Y-\dim Y_{11}-\dim Z_2\\
 &&=\dim {\rm Ker}(A) + \dim Y-(\dim B|_{Y_1}(Y_{11})+ \dim Z_2)\\
&&=\dim {\rm Ker}(A) + \dim Y-\dim{\rm Coker}(A)\\
 &&={\rm Index}(A)+ \dim Y.
\end{eqnarray*}
The desired result is proved.\hfill$\Box$\vspace{2mm}

Let $(X, E, S)$ be a Banach Fredholm bundle of index $r$ and with
compact zero locus $Z(S)$. By Lemma~\ref{lem:1.1} we may assume that
{\bf the section $S:X\to E$ is strongly Fredholm}. For each $x\in
Z(S)$ we can take a small open neighborhood $O_x$ of $x$ in $X$ and
a {\bf local trivialization} $\psi^{(x)}: O_x\times E_x\to
E|_{O_x}$. Let $S_{\psi^{(x)}}: O_x\to E_x$ be the corresponding
trivialization representative of $S|_{O_x}$, i.e., it is defined by
$\psi^{(x)}(y, S_{\psi^{(x)}}(y))=S(y)\;\forall y\in O_x$. Then it
is a Fredholm map. So there exist nonzero $v_{x1},\cdots,v_{xk}\in
E_x$, such that
$$
dS_{\psi^{(x)}}(x)(T_x X)+ {\rm span}(\{v_{1x},\cdots,
v_{xk}\})=E_x.
$$
 By Lemma~\ref{lem:1.2}  we may shrink $O_x$ and
assume that
\begin{equation}\label{e:1.2}
dS_{\psi^{(x)}}(y)(T_y X)+ {\rm span}(\{v_{1x},\cdots,
v_{xk}\})=E_x,\;\forall y\in O_x.
\end{equation}
 Take a smooth cut-off function $\beta_x: O_x\to [0, 1]$ such that it is
equal to $1$ near $x$. Denote by $O_x^0:=\{y\in O_x\,|\,
\beta_x(y)>0\}$.  ({\bf This is where we need to require $X$ having
smooth cut-off functions}.) Then (2) implies that
\begin{equation}\label{e:1.3}
dS_{\psi^{(x)}}(y)(T_y X)+ {\rm span}(\{\beta_x(y)\cdot
v_{1x},\cdots, \beta_x(y)\cdot v_{xk}\})=E_x
\end{equation}
for any $y\in O^0_x$. For $i=1,\cdots, k$, defining $\sigma_{xi}:
X\to E$ by
$$
\sigma_{xi}(y)=\left\{\begin{array}{ll}
\beta_x(y)\cdot\psi^{(x)}(y, v_{xi})\;&{\rm if}\;y\in O_x,\\
0&{\rm if}\; y\notin O_x,
\end{array}\right.
$$
 they have trivialization representatives
$\beta_x\cdot v_{xi}$ under $\psi^{(x)}$, and
$$
 {\rm supp}(\sigma_{xi})\subset
Cl(O_x^0)\subset O_x,\;i=1,\cdots, k.
$$
It follows from (\ref{e:1.2}) that for $y\in O^0_x\cap Z(S)$,
\begin{equation}\label{e:1.4}
DS(y)(T_y X)+ {\rm span}(\{\sigma_{x1}(y),\cdots,
\sigma_{xk}(y)\})=E_x.
\end{equation}
Since $Z(S)$ is compact it can be covered by finitely many such
open subsets $O_{x_i}$, $i=1,\cdots,n$. Setting ${\mathcal
O}(Z(S))=\cup^n_{i=1}O_{x_i}$ we easily get:

\begin{lemma}\label{lem:1.4}
 There exist an open neighborhood
${\mathcal O}(Z(S))$ of $Z(S)$ in $X$  and finitely many smooth
sections $\sigma_1, \sigma_2,\cdots, \sigma_m$ of the bundle $E\to
X$ such that:
\begin{description}
 \item[(i)] $DS(y)(T_y
X)+ {\rm span}(\{\sigma_1(y), \cdots,\sigma_m(y) \})=E_y$ for any
$y\in Z(S)$.

\item[(ii)] ${\rm supp}(\sigma_i)\subset{\mathcal O}(Z(S))$ for
$i=1,2,\cdots, m$.
 \end{description}
 \end{lemma}

Let ${\bf t}=(t_1,\cdots, t_m)\in\R^m$.  Consider smooth sections
\begin{eqnarray}
&&\Phi:X\times\R^m\to\Pi_1^\ast E,\;(y,{\bf t})\mapsto
S(y)+\sum^m_{i=1}t_i\sigma_i(y),\label{e:1.5}\\
&&\Phi_{\bf t}:X\to E,\;y\mapsto
S(y)+\sum^m_{i=1}t_i\sigma_i(y),\label{e:1.6}
 \end{eqnarray}
  where $\Pi_1$ is the projection to the first factor
 of $X\times\R^m$. Lemma~\ref{lem:1.4}(ii)
 implies that
$$
Z(S)\times\{0\}\subset\Phi^{-1}(0)\subset{\mathcal O}(Z(S))\times
\R^m.
 $$

\begin{theorem}\label{th:1.5} There  exist an open
neighborhood ${\cal W}\subset {\cal O}(Z(S))$ of $Z(S)$ and  a
small $\varepsilon>0$ such that:
\begin{description}
\item[(A)] The zero locus of $\Phi$ in $Cl({\cal W}\times
B_\varepsilon(\R^m))$ is compact. Consequently,  for any given
small open neighborhood ${\mathcal U}$ of $Z(S)$ there exists a
$\epsilon\in (0, \varepsilon]$ such that $Cl({\cal
W})\cap\Phi^{-1}_{{\bf t}}(0)\subset{\mathcal U}$ for any ${\bf
t}\in B_\epsilon(\R^m)$. In particular,  each set ${\cal
W}\cap\Phi^{-1}_{\bf t}(0)$ is {\rm compact} for ${\bf t}\in
B_\varepsilon(\R^m)$ sufficiently small.

\item[(B)] The restriction of $\Phi$ to ${\cal W}\times
B_\varepsilon(\R^m)$ is (strongly) Fredholm and also transversal to
the zero section. So
$$ U_\varepsilon:=\{(y,{\bf t})\in {\cal W}\times
B_\varepsilon(\R^m)\,|\, \Phi(y,{\bf t})=0\}
$$
is a smooth manifold of dimension $m+{\rm Ind}(S)$, and  for
 ${\bf t}\in B_\varepsilon(\R^m)$ the section
 $\Phi_{\bf t}|_{\cal W}: X\to E$
 is transversal to the zero section if and only if ${\bf t}$ is
a regular value of the (proper) projection
 $$P_\varepsilon: U_\varepsilon\to B_\varepsilon(\R^m),\;(y,{\bf t})\mapsto{\bf t},$$
 and $\Phi_{\bf t}^{-1}(0)\cap{\cal W}=P_\varepsilon^{-1}({\bf t})$.
(Specially, ${\bf t}=0$ is a regular value of $P_\varepsilon$ if
$S$ is transversal to the zero section).
 Then the Sard theorem yields  a residual subset
  $B_\varepsilon(\R^m)_{res}\subset B_\varepsilon(\R^m)$  such that:
 \begin{description}
 \item[(B.1)] For each
 ${\bf t}\in B_\varepsilon(\R^m)_{res}$ the section $\Phi_{\bf t}|_{\cal
 W}$ is a Fredholm section of index ${\rm Ind}(S)$ and
 the set
  $(\Phi_{\bf t}|_{\cal W})^{-1}(0)\thickapprox(\Phi_{\bf t}|_{\cal W})^{-1}(0)\times\{{\bf t}\}
  =P_\varepsilon^{-1}({\bf t})$
  is a compact smooth manifold of dimension
  ${\rm Ind}(S)$ and  all $k$-boundaries
$$\partial^k(\Phi_{\bf t}|_{\cal W})^{-1}(0)=(\partial^k X)\cap (\Phi_{\bf t}|_{\cal W})^{-1}(0)$$
for $k=1, 2,\cdots$. Specially, if $Z(S)\subset {\rm Int}(X)$ one
can shrink $\varepsilon>0$ so that $(\Phi_{\bf t}|_{\cal
W})^{-1}(0)$ is a closed manifold for each ${\bf t}\in
B_\varepsilon(\R^m)_{res}$.

\item[(B.2)] If the Banach Fredholm bundle $(X,E, S)$ is {\bf
oriented}, i.e., the determinant  bundle ${\rm det}(DS)\to Z(S)$
is given  a nowhere vanishing continuous section over $Z(S)$, then
it determines an orientation on $U_\varepsilon$. In particular, it
induces a natural orientation on every $(\Phi_{\bf t}|_{\cal
W})^{-1}(0)$ for ${\bf t}\in B_\varepsilon(\R^m)_{res}$.

 \item[(B.3)]   For any $l\in\N$ and two different  ${\bf t}^{(1)},
  {\bf t}^{(2)}\in B_\varepsilon(\R^m)_{res}$ the smooth
  manifolds $(\Phi_{{\bf t}^{(1)}}|_{\cal W})^{-1}(0)$ and $(\Phi_{{\bf
  t}^{(2)}}|_{\cal W})^{-1}(0)$ are cobordant in the sense that for a generic
  $C^l$-path $\gamma:[0,1]\to
B_\varepsilon(\R^m)$ with $\gamma(0)={\bf t}^{(1)}$ and
$\gamma(1)={\bf t}^{(2)}$ the set
$$\Phi^{-1}(\gamma):=\cup_{t\in
[0,1]}\{t\}\times (\Phi_{\gamma(t)}|_{\cal W})^{-1}(0)
$$
is a compact smooth manifold with boundary
$$\{0\}\times
(\Phi_{{\bf t}^{(1)}}|_{\cal
W})^{-1}(0)\cup(-\{1\}\times(\Phi_{{\bf t}^{(2)}}|_{\cal
W})^{-1}(0)).$$
 In particular,  if $Z(S)\subset {\rm Int}(X)$ and
 $\varepsilon>0$ is suitably shrunk so that $(\Phi_{\bf t}|_{\cal W})^{-1}(0)\subset {\rm Int}(X)$
for any ${\bf t}\in B_\varepsilon(\R^m)$ then $\Phi^{-1}(\gamma)$
has no corners.

\item[(B.4)] The cobordism class of the manifold $(\Phi_{\bf
t}|_{\cal W})^{-1}(0)$ above is independent of all related
choices.
  \end{description}
  \end{description}
  \end{theorem}

When a smooth map $f:M\to N$ is transversal to a submanifold
$S\subset N$, the inverse image $f^{-1}(S)$ is either a manifold or
an empty set. In this paper we always follow the usual {\bf
convention}: not mentioning the second case. \vspace{2mm}

\noindent{\bf Proof of Theorem~\ref{th:1.5}.}\quad The direct
computation shows that the vertical differential
\begin{equation}\label{e:1.7}
D\Phi(y,0):T_{(y,0)}(X\times\R^m)\to (\Pi_1^\ast E)_{(y,0)}=E_y
\end{equation}
 of $\Phi$ at any
point $(y,0)\in Z(S)\times\R^m$ is given by
$$
D\Phi(y,{\rm 0})(\xi,{\rm u})=DS(y)(\xi)+\sum^m_{k=1} u_k\cdot
\sigma_k(y)
$$
 for $\xi\in T_y X$ and ${\bf u}=(u_1, u_2, \cdots, u_m)\in \R^m$.
By Lemma~\ref{lem:1.4}(i), the linear continuous operator
$D\Phi(y,0)$ is surjective for each $(y,0)\in
Z(S)\times\{0\}\subset\Phi^{-1}(0)$. Lemma~\ref{lem:1.3}(ii) also
implies that $D\Phi(y,0)$ is a Fredholm operator of index ${\rm
Ind}(S)+ m$.
 From the compactness of the subset
$Z(S)\times\{0\}\subset X\times\R^m$ and Lemma~\ref{lem:1.1} it
follows that in some open neighborhood of $Z(S)\times\{0\}$ the
section $\Phi$ is (strongly) Fredholm and also transversal to the
zero section. Using the properness of the Fredholm map again we may
derive that there exist an open neighborhood ${\cal W}\subset {\cal
O}(Z(S))$ of $Z(S)$ and a small $\varepsilon>0$ such that:

$\bullet$ The restriction of $\Phi$ to ${\cal W}\times
B_\varepsilon(\R^m)$ is a (strongly) Fredholm section of index ${\rm
Ind}(S)+ m$ and is also transversal to the zero section.

$\bullet$ The zero locus of $\Phi$ in $Cl({\cal W}\times
B_\varepsilon(\R^m))$ is compact.

\noindent{\bf Step 1}. Let us prove that the first claim in
Theorem~\ref{th:1.5}(A) implies others.  Suppose that there exist
a small neighborhood ${\cal U}$ of $Z(S)$ in ${\cal W}$ and a
sequence of zero points $\{(y_i,{\bf t}^{(i)})\}$ of $\Phi$ in
${\cal W}\times B_\varepsilon(\R^m)$ such that ${\bf t}^{(i)}\to
0$ as $i\to \infty$ and $y_i\notin{\mathcal U}$ for
$i=1,2,\cdots$. Since $\Phi^{-1}(0)\cap Cl({\cal W}\times
B_\varepsilon(\R^m))$ is compact there exists a subsequence of
$\{y_i\}$, still denoted by $\{y_i\}$, converging to some $y_0\in
Cl({\mathcal W})$. So $\Phi(y_0, 0)=0$ or $S(y_0)=0$. That is,
$y_0\in Z(S)$. This implies that $y_i\in{\mathcal U}$ for very
large $i$, which leads to a contradiction. Next we prove the
second claim. Take a neighborhood ${\cal V}$ of $Z(S)$ in ${\cal
W}$ such that $Cl({\cal V})\subset{\cal W}$. By the first claim,
after shrinking $\varepsilon>0$ we get that $Cl({\cal
W})\cap\Phi^{-1}_{{\bf t}}(0)\subset{\mathcal V}$ for any ${\bf
t}\in B_\varepsilon(\R^m)$. For any sequence $\{y_k\}\subset{\cal
W}\cap\Phi^{-1}_{\bf t}(0)$, as above we can derive from the first
claim that $\{y_k\}$ has a subsequence $\{y_{k_i}\}$ converging to
$y_0\in Cl({\cal V})\subset{\cal W}$. Specially, $\Phi(y_0, {\bf
t})=0$. Therefore $y_0\in{\cal W}\cap\Phi_{\bf t}^{-1}(0)$. The
desired conclusion is proved.

We only need to prove  (B.2), (B.3) and (B.4) since other claims in
Theorem~\ref{th:1.5}(B) are clear. (In fact, the claim above Step 1
implies that $U_\varepsilon$ is a smooth manifold of dimension $m+
{\rm Ind}(S)$. So the projection $P_\varepsilon$ has Fredholm index
${\rm Ind}(S)$.
 Then the first claim in (B.1) follows from the fact that the kernel ${\rm Ker}(D\Phi_{\bf t}(x))={\rm
Ker}(DP_\varepsilon(x))$ has dimension ${\rm Ind}(S)$ and $\dim{\rm
Coker}(D\Phi_{\bf t}(x))=\dim{\rm Coker}(DP_\varepsilon(x))=0$ for
any $x\in\Phi^{-1}_{\bf t}(0)\cap{\cal W}$.)
 The following proofs
are in three steps.

\noindent{\bf Step 2}.  In order to prove (B.2) let us consider
the homotopy sections
$$
\Psi^t:{\cal W}\times B_\varepsilon(\R^m)\to\Pi_1^\ast E,\;(y,{\bf
t})\mapsto S(y) + t\sum^m_{i=1}t_i\sigma_i(y)
 $$
for $0\le t\le 1$. By shrinking $\varepsilon>0$ we can assume that
they are all Fredholm and thus have the same index as $\Phi$. By
the construction of $U_\varepsilon$ above,
$$
T_{(y,{\bf t})}U_\varepsilon={\rm Ker}(D\Phi)(y,{\bf t})\quad{\rm
and}\quad {\rm Coker}(D\Phi)(y,{\bf t}))=\{0\}.
$$
So an orientation of $U_\varepsilon$ is equivalent to giving a
continuous nowhere zero section of ${\rm det}(D\Phi)$.  Note that
$Z(S)\times\{0\}\subset (\Psi^t)^{-1}(0)$ for any $t\in [0,1]$. We
have
\begin{eqnarray*}
{\rm det}(D\Phi)|_{Z(S)\times\{0\}}&=&{\rm
det}(D\Psi^1)|_{Z(S)\times\{0\}}\\
 &=&{\rm det}(D\Psi^0)|_{Z(S)\times\{0\}}\\
 &=&({\rm det}(DS)\otimes\det(\R^m))|_{Z(S)\times\{0\}}.
\end{eqnarray*}
Let $o(S)$ be the given continuous nowhere zero section of ${\rm
det}(DS)|_{Z(S)}$. Then $o(S)\otimes{\bf 1}$ is such an section of
$({\rm det}(DS)\otimes\det(\R^m))|_{Z(S)\times\{0\}}$ and thus
gives that of ${\rm det}(D\Phi)|_{Z(S)\times\{0\}}$.

Since $E$ is local trivial, so is ${\rm
det}(D\Phi)|_{U_\varepsilon}$.
 It follows that  the nowhere zero section
$o(S)\otimes{\bf 1}$  may be extended into a continuous nowhere
zero section of  ${\rm det}(D\Phi)|_{U_\varepsilon}$. The latter
naturally restricts to such a section of ${\rm det}(D\Phi_{\bf
t})$ on  $(\Phi_{\bf t}|_{\cal W})^{-1}(0)$ for each ${\bf t}\in
B_\varepsilon(\R^m)_{res}$.

\noindent{\bf Step 3}.  Now let us to prove (B.3).
 Let ${\cal P}^l({\bf t}^{(1)}, {\bf t}^{(2)})$
denote the space of all $C^l$-smooth paths $\gamma:[0,1]\to
B_\varepsilon(\R^m)$ with $\gamma(0)={\bf t}^{(1)}$ and
$\gamma(1)={\bf t}^{(2)}$. Consider the section
\begin{eqnarray*}
{\cal F}: {\cal W}\times [0, 1]\times {\cal P}^l({\bf t}^{(1)},
{\bf
t}^{(2)})\to \Pi_1^\ast E,\\
(x, \tau,\gamma)\mapsto
S(x)+\sum^m_{i=1}\gamma_i(\tau)\sigma_i(x),
\end{eqnarray*}
where $\Pi_1: {\cal W}\times [0, 1]\times {\cal P}^l({\bf
t}^{(1)}, {\bf t}^{(2)})\to {\cal W}$ is the projection and
$\gamma=(\gamma_1,\cdots, \gamma_m)\in {\cal P}^l({\bf t}^{(1)},
{\bf t}^{(2)})$. The (vertical) differential of ${\cal F}$ at any
zero $(x,\tau,\gamma)$ of it is given by
\begin{eqnarray*}
&&D{\cal F}(x,\tau,\gamma): T_xX\times\R\times T_\gamma{\cal
P}^l({\bf t}^{(1)}, {\bf t}^{(2)})\to E_x,\\
&&\hspace{20mm} (\xi, e, \alpha)\to D\Phi(x, \gamma(\tau))(\xi,
\alpha(\tau)+ e\dot\gamma(\tau)).
\end{eqnarray*}
For $\tau=0$ and $\tau=1$ it is surjective as a function of $\xi$
alone since the sections ${\cal F}(\cdot, 0,\gamma)=\Phi_{\bf
t}^{(1)}|_{\cal W}$ and ${\cal F}(\cdot, \tau,\gamma)=\Phi_ {\bf
t}^{(2)}|_{\cal W}$ are transversal to the zero section. For
$\tau\in [0, 1]\setminus\{0,1\}$ we can choose $\alpha(\tau)$
arbitrarily and thus derive that this operator is also surjective as
a function of $\xi$ and $\alpha$. Hence ${\cal F}$ is transversal to
the zero section. Let $P$ be the projection from ${\cal F}^{-1}(0)$
to the third factor. For $(x,\tau,\gamma)\in {\cal F}^{-1}(0)$ we
have
\begin{eqnarray*}
&&T_{(x,\tau,\gamma)}{\cal F}^{-1}(0)={\rm Ker}(D{\cal
F}(x,\tau,\gamma))\\
&&=\{(\xi, e, \alpha)\in T_xX\times\R\times T_\gamma{\cal P}^l({\bf
t}^{(1)}, {\bf t}^{(2)})\,|\\
&&\hspace{40mm} D\Phi(x, \gamma(\tau))(\xi, \alpha(\tau)+
e\dot\gamma(\tau))=0\}.
\end{eqnarray*}
So $(\xi, e, \alpha)\in{\rm Ker}(DP(x,\tau,\gamma))$ if and only if
$\alpha=0$ and
$$D\Phi(x, \gamma(\tau))(\xi, e\dot\gamma(\tau))=0.$$
Note that $0={\cal F}(x,\tau,\gamma)=\Phi(x,\gamma(\tau))$. This
means that $(\xi, e\dot\gamma(\tau))\in
T_{(x,\gamma(\tau))}U_\varepsilon$. It follows that $\dim {\rm
Ker}(DP(x,\tau,\gamma))={\rm Ind}(S)+ 1$. Hence $P$ is a Fredholm
operator of index ${\rm Ind}(S)+ 1$. One also easily prove that
$\gamma\in {\cal P}^l({\bf t}^{(1)}, {\bf t}^{(2)})$ is a regular
value of $P$ if and only if the section
$$
{\cal F}_\gamma: {\cal W}\times [0, 1]\to \Pi_1^\ast E, (x,
\tau)\mapsto S(x)+\sum^m_{i=1}\gamma_i(\tau)\sigma_i(x),
$$
is  transversal to the zero section. (This also implies ${\cal
F}_\gamma$ to be a Fredholm section of index ${\rm Ind}(S)+1$.)
Hence
 the Sard-Smale theorem  yields a residual subset ${\cal P}^l_{reg}({\bf t}^{(1)},
{\bf t}^{(2)})\subset {\cal P}^l({\bf t}^{(1)}, {\bf t}^{(2)})$
  such that each $\gamma\in {\cal P}^l_{reg}({\bf
t}^{(1)}, {\bf t}^{(2)})$ gives a compact  cobordsim
$P^{-1}(\gamma)={\cal F}_\gamma^{-1}(0)$ between $(\Phi_{{\bf
t}^{(1)}}|_{\cal W})^{-1}(0)$ and $(\Phi_{{\bf t}^{(2)}}|_{\cal
W})^{-1}(0)$.  The second claim easily follows from the second one
in (B.1). (B.3) is proved.

 \noindent{\bf Step 4}. Finally we prove (B.4). As above assume that
$\varepsilon'>0$, an open neighborhood ${\cal W}'$ of $Z(S)$ and
 another group of sections of $E$, $\sigma'_1,\cdots,\sigma'_{m'}$  such that
 the   section
\begin{equation}\label{e:1.8}
\Psi: {\cal W}'\times B_{\varepsilon'}(\R^{m'})\to \Pi^\ast_1E,\;
(y, {\bf t}')\mapsto S(y)+ \sum^{m'}_{i=1}t'_i\sigma'_i(y),
\end{equation}
 is Fredholm and transversal to the zero section and that  the  set $\Psi^{-1}_{{\bf
t}'}(0)$ is compact for each ${\bf t}'\in
B_{\varepsilon'}(\R^{m'})$. Here the section $\Psi_{{\bf t}'}:
{\cal W}'\to E$ is given by $\Psi_{{\bf t}'}(y)=\Psi(y,{\bf t}')$.
Let $B_{\varepsilon'}(\R^{m'})_{res}\subset
B_{\varepsilon'}(\R^{m'})$ be the corresponding residual subset
such that for each ${\bf t}'\in B_{\varepsilon'}(\R^{m'})_{res}$
the section $\Psi_{{\bf t}'}$ is transversal to the zero section
and that any two ${\bf t}', {\bf s}'\in
B_{\varepsilon'}(\R^{m'})_{res}$ yield cobordant manifolds
$(\Psi_{{\bf t}'})^{-1}(0)$  and $(\Psi_{{\bf s}'})^{-1}(0)$.

For $0<\eta<\min\{\varepsilon,\varepsilon'\}$ let us consider the
section
\begin{eqnarray}
&&\Theta:{\cal W}\cap{\cal W}'\times B_\eta(\R^m)\times
B_\eta(\R^{m'})\times [0,1]\to\Pi_1^\ast
E,\label{e:1.9}\\
&&\qquad (y, {\bf t}, {\bf t}', \tau)\mapsto S(y)+
(1-\tau)\sum^m_{i=1} t_i \sigma_i(y)+\tau\sum^{m'}_{j=1} t'_j
\sigma^\ast_j(y).\nonumber
\end{eqnarray}
Then $Z(S)\times\{0\}\times\{0\}\times [0,1]\subset\Theta^{-1}(0)$
and the vertical differential
$$D\Theta(y,0,0,\tau): T_yX\times\R^m\times\R^{m'}\times\R\to E_y$$
is Fredholm and surjective for any $(y,0,0,\tau)\in
Z(S)\times\{0\}\times\{0\}\times [0,1]$.   The standard arguments
lead to:

\begin{claim}\label{claim:1.6}
 There exist an open neighborhood ${\cal W}^\ast\subset{\cal W}\cap{\cal W}'$
 of $Z(S)$ and a small $\eta>0$  such that:
\begin{description}
\item[(ii)] The set $\Theta^{-1}(0)\cap Cl({\cal W}^\ast\times
B_{\eta}(\R^m)\times B_{\eta}(\R^{m'})\times [0,1])$ is compact.
So for any given small open neighborhood ${\mathcal V}$ of $Z(S)$
there exists a small $\epsilon\in (0,\eta]$ such that
 for any point $(y, {\bf t}, {\bf
t}',\tau)$ in $\Theta^{-1}(0)\cap Cl({\cal W}^\ast\times
B_{\epsilon}(\R^m)\times B_{\epsilon}(\R^{m'})\times [0,1])$ one
has $y\in{\mathcal V}$.

\item[(ii)]  The restriction of $\Theta$ to ${\cal W}^\ast\times
B_{\eta}(\R^m)\times B_{\eta}(\R^{m'})\times [0,1]$
 is transversal to the zero section.
So there exists a residual subset $(B_{\eta}(\R^m)\times
B_{\eta}(\R^{m'}))_{res}\subset B_{\eta}(\R^m)\times
B_{\eta}(\R^{m'})$ such that for any $({\bf t},{\bf t}')\in
(B_{\eta}(\R^m)\times B_{\eta}(\R^{m'}))_{res}$ the section
$$
\Theta_{({\bf t}, {\bf t}')}: {\cal W}^\ast\times[0,1]\to
\Pi^\ast_1(E|_{{\cal W}^\ast}),\;(y, \tau)\mapsto \Theta(y, {\bf
t}, {\bf t}',\tau)
$$
is Fredholm and transversal to the zero section, and thus
$\Theta_{({\bf t},{\bf t}')}^{-1}(0)$ is a compact manifold of
dimension $r+1$ and with boundary (and corners).
\end{description}
\end{claim}

To finish the final proof we also need the following lemma.

\begin{lemma}\label{lem:1.7}\cite{LeO}\quad Let $X$ and $Y$ be
metric spaces which satisfy the second axiom of countability.
Suppose that $S$ is a countable intersection of open dense subsets
in the product space $X\times Y$. Consider the space $X_S$
consisting of those $x\in X$ such that $S\cap\{x\}\times Y$ is a
countable intersection of open dense subsets in $\{x\}\times Y$.
Then $X_S$ is a countable intersection of open dense subsets in $X$.
\end{lemma}

Take an open neighborhood ${\cal V}$ of $Z(S)$ in $X$ so that
$Cl({\cal V})\subset {\cal W}^\ast$.  By Claim~\ref{claim:1.6} and
Lemma~\ref{lem:1.7} we may take
$$
({\bf t},{\bf t}')\in \bigl((B_{\eta}(\R^m)\times
B_{\eta}(\R^{m'}))_{res}\bigr)\cap
\bigl(B_{\varepsilon}(\R^m)_{res}\times
B_{\varepsilon'}(\R^{m'})_{res}\bigr)
$$
so small that $y\in{\cal V}$ for any $(y, \tau)\in\Theta_{({\bf
t},{\bf t}')}^{-1}(0)$. Note  that $\Theta(\cdot, {\bf t}, {\bf
t}', 0)=\Phi_{\bf t}|_{\cal W}$ and $\Theta(\cdot, {\bf t}, {\bf
t}', 1)=\Psi(\cdot,{\bf t}')$. The manifold $\Theta_{({\bf t},{\bf
t}')}^{-1}(0)$  forms a cobordsim between $(\Phi_{\bf t}|_{\cal
W})^{-1}(0)\thickapprox(\Phi_{\bf t}|_{\cal
W})^{-1}(0)\times\{0\}$ and $\Psi_{{\bf
t}'}^{-1}(0)\thickapprox\Psi^{-1}_{{\bf t}'}(0)\times\{1\}$,
 i.e. $\partial\Theta_{({\bf t},{\bf t}')}^{-1}(0)=(\Phi_{\bf t}|_{\cal
W})^{-1}(0)\times\{0\}\cup(-\Psi^{-1}_{{\bf t}'}(0)\times\{1\})$. If
$(X,E,S)$ is oriented then it is also oriented cobordsim. From this
and (B.2) it follows that for any ${\bf t}\in
B^{reg}_{\varepsilon}(\R^m)$ and ${\bf t}'\in
B^{reg}_{\varepsilon'}(\R^{m'})$ the corresponding manifolds
$(\Phi_{\bf t}|_{\cal W})^{-1}(0)$ and $\Psi^{-1}_{{\bf t}'}(0)$ are
cobordism. In summary we have shown that the cobordant class of the
manifold $(\Phi_{\bf t}|_{\cal W})^{-1}(0)$  is independent of all
related choices. Theorem~\ref{th:1.5} is proved.
\hfill$\Box$\vspace{2mm}

Each $(\Phi_{\bf t}|_{\cal W})^{-1}(0)$ is called the {\bf Euler
chain} of $(X, E, S)$. If it has no boundary $(\Phi_{\bf t}|_{\cal
W})^{-1}(0)$ is a cycle, called the {\bf Euler cycle} of $(X, E,
S)$. In this case the homology class $[(\Phi_{\bf t}|_{\cal
W})^{-1}(0)]$ in $H_r(X,\Z_2)$ (resp. $H_r(X, \Z)$ if $(X,E,S)$ is
oriented) is only dependent on $(X,E,S)$, called the {\bf Euler
class} of it and denoted by $e(E,S)$.

Some authors, ex. \cite{R1}, \cite{R2} and \cite{CMSa}, like to view
$e(E,S)$ as a homomorphism $\mu_{(E,S)}$ from $H^\ast (X,\R)$ to
$\R$ by
$$
\mu_{(E,S)}(\alpha)=\langle\alpha, e(E,S)\rangle\;{\rm
for}\;\alpha\in H^r(X,\R),
$$
and $\mu_{(E,S)}(\alpha)=0$ for other $\alpha$. Let
$P_\varepsilon: U_\varepsilon\to B_\varepsilon(\R^m)$ be as in
Theorem~\ref{th:1.5}(B). Take a Thom form $\tau$ on $\R^m$ with
support in $B_\varepsilon(\R^m)$. Then
$$\mu_{(E,S)}(\alpha)=\int_{U_\varepsilon}(R_\varepsilon^\ast\alpha)\wedge(P^\ast_\varepsilon\tau),$$
where $R_\varepsilon$ is the composition of the projection
$U_\varepsilon\to{\cal W}$ and the inclusion ${\cal
W}\hookrightarrow X$, and we have used the same notation
$R_\varepsilon^\ast\alpha$ to denote its closed form
representative. This sometimes is convenient.

\begin{remark}\label{rem:1.8}
{\rm By the arguments in lemma~\ref{lem:1.4} we may require that at
least one of $\sigma_1(z)$, $\sigma_2(z),\cdots, \sigma_m(z)$ be
nonzero at each point $z\in Z(S)$ even if $DS(y)(T_yX)=E_y$ for some
$y\in Z(S)$. In this case the linear map
$$L_x:\R^m\to {\rm Span}(\{\sigma_1(x), \cdots,
\sigma_m(x)\}),\;{\bf t}\mapsto\sum^m_{i=1}t_i\sigma_i(x).
$$
is nonzero for any $x\in Z(S)$. So the subspace ${\rm
Ker}(L_x)\subset\R^m$ has at least codimension one. Denote by $S^m$
the unit sphere in $\R^m$. Take a small open neighborhood ${\cal
U}(S^m\cap{\rm Ker}(L_x))$ of $S^m\cap{\rm Ker}(L_x))$ in $S^m$.
Then $L_x({\bf t})\ne 0$ for any ${\bf t}\in S^m\setminus {\cal
U}(S^m\cap{\rm Ker}(L_x))$. Note that $S^m\setminus {\cal
U}(S^m\cap{\rm Ker}(L_x))$ is a compact subset in $S^m$. By a
contradiction arguments one easily shows that there exists a small
open neighborhood $O_x$ of $x$ in $X$ such that $L_z({\bf
t})=\sum^m_{i=1}t_i\sigma_i(z)\ne 0$ for any $z\in O_x$. Let
$\angle{\cal U}(S^m\cap{\rm Ker}(L_x))$ be the open cone spanned by
${\cal U}(S^m\cap{\rm Ker}(L_x))$. Then
$B_\varepsilon(\R^m)_{res}\setminus\angle{\cal U}(S^m\cap{\rm
Ker}(L_x))$ is a residual subset in
$B_\varepsilon(\R^m)\setminus\angle{\cal U}(S^m\cap{\rm Ker}(L_x))$
and for any ${\bf t}\in
B_\varepsilon(\R^m)_{res}\setminus\angle{\cal U}(S^m\cap{\rm
Ker}(L_x))$ the manifold $(\Phi_{\bf t}|_{\cal W})^{-1}(0)$
satisfies
$$(\Phi_{\bf t}|_{\cal W})^{-1}(0)\cap(O_x\cap Z(S))=\emptyset.$$
So for a given finite subset $F\subset Z(S)$ one can choose
$(\Phi_{\bf t}|_{\cal W})^{-1}(0)$ so that $(\Phi_{\bf t}|_{\cal
W})^{-1}(0)\cap F=\emptyset$.

 \noindent{\bf Question}.\quad Is there ${\bf t}\in
B_\varepsilon(\R^m)_{res}$ such that $(\Phi_{\bf t}|_{\cal
W})^{-1}(0)\cap Z(S)=\emptyset$ or $(\Phi_{\bf t}|_{\cal
W})^{-1}(0)\cap Z(S)$ is nowhere dense in $(\Phi_{\bf t}|_{\cal
W})^{-1}(0)$? (Some attempts: Using the compactness of $Z(S)$ we can
find finitely many subspaces of codimension at least one in $\R^m$,
saying $H_1,\cdots, H_r$, and a very thin open neighborhood ${\cal
U}(\cup^r_{i=1}S^m\cap H_i)$ of $\cup^r_{i=1}S^m\cap H_i$ in $S^m$
such that in some open neighborhood of $Z(S)$,
$\sum^m_{i=1}t_i\sigma_i(y)\ne 0$ for each ${\bf
t}\in\R^m\setminus\{0\}$ with ${\bf t}/|{\bf t}|\in {\cal
U}(\cup^r_{i=1}S^m\cap H_i)$. Since we have assumed $X$ to be a
separable Banach manifold it possess a countable base. It follows
that the compact subspace $Z(S)$ has a countable base and thus is
also separable.) }
\end{remark}

\begin{remark}\label{rem:1.9}{\rm Carefully checking the above
construction one easily sees that for any compact subset
$\Lambda\subset Z(S)$ (even if $Z(S)$ is not compact) we can find
  open neighborhoods ${\cal
W}(\Lambda)\subset{\cal O}(\Lambda)$ of $\Lambda$, smooth sections
$\sigma_i$ with supports in ${\cal O}(\Lambda)$, $i=1,\cdots,l$,
and $\epsilon>0$ such that:
\begin{description}
\item[(a)] The section
$$\Phi^\Lambda: {\cal W}(\Lambda)\times B_\epsilon(\R^l)\to\Pi^\ast_1E,\,
(y, {\bf t})\mapsto S(y)+\sum^l_{i=1}t_i\sigma_i(y)
$$
is Fredholm and transversal to the zero section.

\item[(b)] The zero locus $(\Phi^\Lambda)^{-1}(0)$ has compact
closure in $Cl({\cal W}(\Lambda)\times B_\varepsilon(\R^l))$ and
thus $(\Phi^\Lambda)^{-1}(0)\cap(\Lambda\times\{0\})$ is compact.
\end{description}
Then there exists a residual subset $B_\epsilon(\R^l)_{res}\subset
B_\epsilon(\R^l)$ such that for each ${\bf t}\in
B_\epsilon(\R^l)_{res}$ the section
$$
\Phi^\Lambda_{{\bf t}}: {\cal W}(\Lambda)\to
E,\;y\mapsto\Phi^\Lambda(y, {\bf t})
$$
is transversal to the zero section and each intersection
$(\Phi^\Lambda_{{\bf t}})^{-1}(0)\cap\Lambda$ is compact. As before
any two different  ${\bf t}, {\bf t}'\in B_\epsilon(\R^l)_{res}$
give smooth cobordant
  manifolds $(\Phi^\Lambda_{{\bf t}})^{-1}(0)$ and $(\Phi^{\Lambda}_{{\bf
  t}'})^{-1}(0)$ in the sense that there exist generic paths $\gamma:[0,1]\to
B_\epsilon(\R^l)$ with $\gamma(0)={\bf t}$ and $\gamma(1)={\bf
t}'$ such that
$$(\Phi^\Lambda)^{-1}(\gamma):=\cup_{t\in
[0,1]}\{t\}\times(\Phi^\Lambda_{\gamma(t)})^{-1}(0)
$$
is a smooth manifold with boundary
$\{0\}\times(\Phi^{\Lambda}_{\bf
t})^{-1}(0)\cup(-\{1\}\times(\Phi^\Lambda_{{\bf t}'})^{-1}(0))$
which is relative compact in $[0,1]\times Cl({\cal W}(\Lambda))$.
Note also that each intersection $\Lambda\cap(\Phi^\Lambda_{\bf
t})^{-1}(0)$ is compact since $\Lambda$ is compact.  Later each
$(\Phi_{\bf t}^\Lambda)^{-1}(0)$ is also called a {\bf local Euler
chain} of $(X, E, S)$ near $\Lambda$. Moreover, for any given open
neighborhood $U(\Lambda)$ of $\Lambda$ in $X$ we can shrink
$\epsilon>0$ so that each $(\Phi^\Lambda_{\bf t})^{-1}(0)$ is
contained in $U(\Lambda)$ for any ${\bf t}\in B_\epsilon(\R^l)$.}
\end{remark}

\begin{remark}\label{rem:1.10}
{\rm If $Z(S)$ is not  compact we can extend Theorem~\ref{th:1.5} to
the case that $Z(S)$ is $\sigma$-compact, i.e.,  it is the union of
countable compact sets. In this case we can similarly define  the
{\bf Euler class} of the triple $(X, E, S)$ sitting in the homology
of the second kind $H_r^{II}(X,\Z_2)$ (resp. $H_r^{II}(X,\Z)$ if the
real (determinant) line bundle ${\rm det(Ind(DS))}$ is trivial over
$Z(S)$).}
\end{remark}

\subsection{Localization formula}\label{sec:1.2}

In applications we often meet the following {\bf model}:
 Let $(X,E,S)$ be a Banach Fredholm
bundle of index $r$ and with compact zero locus $Z(S)$. Assume
that it has the Euler class $e(E,S)$ as above. Let $P$ be a
manifold of finite dimension, and $f:X\to P$ be a smooth map. For
$\alpha\in H^r(P, \R)$ let $\alpha^\ast$ be a differential form
representative of it.  If $(\Phi_{\bf t}|_{\cal W})^{-1}(0)$ is a
Euler cycle of $(X,E,S)$ that is a closed manifold as above then
\begin{equation}\label{e:1.10}
\langle f^\ast\alpha, e(E,S)\rangle=\langle\alpha, f_\ast
e(E,S)\rangle =\int_{(\Phi_{\bf t}|_{\cal
W})^{-1}(0)}f^\ast\alpha^\ast.
\end{equation}

\begin{proposition}\label{prop:1.11}
({\bf First localization  formula}). If $\alpha$ has a
representative form $\alpha^\ast$ such that $f^\ast\alpha^\ast$
has support ${\rm supp}(f^\ast\alpha^\ast)$ contained in a compact
subset $\Lambda\subset Z(S)$ then for the family
$\{(\Phi^\Lambda_{\bf t})^{-1}(0)\,|\,{\bf t}\in
B_\epsilon(\R^l)_{res}\}$ in Remark~\ref{rem:1.9} there exists a
residual subset $B_\epsilon(\R^l)^\star_{res}\subset
B_\epsilon(\R^l)_{res}$ such that
$$
\langle e(E,S), f^\ast\alpha\rangle=
\int_{(\Phi^\Lambda_{\bf t})^{-1}(0)}f^\ast\alpha^\ast
$$
for ${\bf t}\in B_\epsilon(\R^l)^\star_{res}$.
\end{proposition}

\noindent{\bf Proof.}\quad We always assume that $Z(S)$ is compact
in the following proof. Assume that the family $\{(\Phi^\Lambda_{\bf
t})^{-1}(0)\,|\,{\bf t}\in B_\epsilon(\R^l)_{res}\}$ is as in
Remark~\ref{rem:1.9}. Fix a small open neighborhood ${\cal
W}(\Lambda)_0$ of $\Lambda$ in ${\cal W}(\Lambda)$ so that
$$
{\cal W}(\Lambda)_0\subset\subset{\cal W}(\Lambda).
$$
{\bf Hereafter} we write $A\subset\subset B$ to mean that $A$ has
closure contained in $B$ without special statements.
 Note that
$Z(S)\setminus{\cal W}(\Lambda)$ is compact. We furthermore take
finitely many points $x_i\in Z(S)\setminus{\cal W}(\Lambda)$ and
their open neighborhoods $O_i$ in $X$, $i=k+1,\cdots, n$, and smooth
sections $\sigma_j$
 of $E$ with supports in ${\cal
O}(\Lambda)_2:=\cup^n_{i=k+1}O_i$, $j=l+1,\cdots, m$, satisfying
$$
Z(S)\setminus{\cal W}(\Lambda)\subset{\cal O}(\Lambda)_2\quad{\rm
and}\quad {\cal O}(\Lambda)_2\cap{\cal W}(\Lambda)_0=\emptyset.
$$
As before one can find an open neighborhood ${\cal W}\subset{\cal
W}(\Lambda)\cup{\cal O}(\Lambda)_2$ and $\varepsilon\in (0,
\epsilon]$ such that:
\begin{description}
\item[(I)]  The section
$$\Phi: {\cal W}\times B_\varepsilon(\R^m)\to\Pi^\ast_1E,\,
(y, {\bf s})\mapsto S(y)+\sum^m_{i=1}s_i\sigma_i(y)
$$
is Fredholm and transversal to the zero section.

\item[(II)] The zero set $\Phi^{-1}(0)$ has compact closure in
$Cl({\cal W}\times B_\varepsilon(\R^m))$.
\end{description}
Then one has a residual subset $B_\varepsilon(\R^m)_{res}\subset
B_\varepsilon(\R^m)$ such that the section
$$\Phi_{\bf s}: {\cal W}\to\Pi^\ast_1E,\,
y\mapsto S(y)+\sum^m_{i=1}s_i\sigma_i(y)
$$
is transversal to the zero section for each ${\bf s}\in
B_\varepsilon(\R^m)_{res}$.

Let $\Pi^m_l$ denote the projection from $\R^m$ onto the former
$l$ coordinates. By Lemma~\ref{lem:1.7},
$\Pi^m_l(B_\varepsilon(\R^m)_{res})$ is still a residual subset in
$B_\varepsilon(\R^l)$. So is the intersection
$$
B_\varepsilon(\R^l)^\star_{res}:=B_\epsilon(\R^l)_{res}\cap
\Pi^m_l(B_\varepsilon(\R^m)_{res}).
$$
Note that for any ${\bf s}'=(s_1,\cdots,s_l)\in
B_\epsilon(\R^l)^\star_{res}$ we can always choose
$s_{l+1},\cdots, s_m$ in $\R$ such that  ${\bf s}=(s_1,\cdots,
s_m)\in B_\varepsilon(\R^m)_{res}$. For such ${\bf s}$ it is
easily checked that $\Phi_{\bf s}^{-1}(0)\cap{\cal W}(\Lambda)_0$
agrees with $(\Phi^\Lambda_{{\bf s}'})^{-1}(0)\cap{\cal
W}(\Lambda)_0$ because the sections $\sigma_i$, $i=l+1,\cdots, m$
have supports outside ${\cal W}(\Lambda)_0$. It follows from ${\rm
supp}(f^\ast\alpha^\ast)\subset \Lambda$ that
\begin{eqnarray*}
\int_{(\Phi^\Lambda_{{\bf s}'})^{-1}(0)}f^\ast\alpha^\ast
&=&\int_{(\Phi^\Lambda_{{\bf s}'})^{-1}(0)\cap{\cal
W}(\Lambda)_0}f^\ast\alpha^\ast\\&=&\int_{(\Phi_{\bf
s})^{-1}(0)\cap{\cal
W}(\Lambda)_0}f^\ast\alpha^\ast\\&=&\int_{(\Phi_{\bf
s})^{-1}(0)}f^\ast\alpha^\ast\\&=&\langle e(E,S),
f^\ast\alpha\rangle.
\end{eqnarray*}
The desired result is proved. \hfill$\Box$\vspace{2mm}

If $\Lambda$ is a connected component of $Z(S)$, it is also
compact and we can require that the above small neighborhood
${\cal W}(\Lambda)$ of $\Lambda$ is disjoint with
$Z(S)\setminus\Lambda$. In this situation
 $\{(\Phi^\Lambda_{\bf t})^{-1}(0)\,|\,{\bf t}\in
B_\epsilon(\R^l)_{res}\}$ is a family of cobordant cycles and thus
determines a homology class in $H_r(X)$, denoted by
$e(E,S)_\Lambda$. Let $\Lambda_i$, $i=1,\cdots, p$, be all
connected components of $\Lambda$. Then
$$
e(E,S)=\sum^p_{i=1}e(E,S)_{\Lambda_i}.
$$

 The second  localization formula is as follows:

\begin{proposition}\label{prop:1.12}
({\bf Second localization  formula}). Let $(X,E,S)$ be an oriented
Banach Fredholm bundle of index $r$ and with compact zero locus
$Z(S)$. Suppose that a closed subset $Y\subset X$ is a Banach
submanifold of finite positive codimension.   Then $(Y, E|_Y,
S|_Y)$ is still a Banach Fredholm bundle with a natural induced
orientation and there exist the Euler chains $M$ of $(X, E, S)$
and $N$ of $(Y, E|_Y, S|_Y)$ such that $M\cap Y=N$. Moreover, $M$
and $N$ can be chosen as closed manifolds if both $e(E,S)$ and
$e(E|_Y, S|_Y)$ exist. In particular, if $P$ is an oriented smooth
manifold of finite dimension, $Q\subset P$ is an oriented closed
submanifold and a smooth map $f:X\to P$ is transversal to $Q$,
then for $Y:=f^{-1}(Q)\subset X$ ( a Banach submanifold satisfying
the above requirements!), when both $e(E,S)$
 and $e(E|_Y, S|_Y)$ exist, it holds that
 $$I_{Q\ast}(a)\cdot_P f_\ast(e(E,S))
 =a\cdot_Q
f_\ast(e(E|_Y, S|_Y))
$$
for any class $a\in H_\ast(Q, \R)$ of codimension $r$. (This
should imply that $f_\ast(e(E|_Y, S|_Y)=f_\ast(e(E,S))\cap
PD(Q)$.)  Here $I_Q:Q\to P$  is the inclusion map, $\cdot_P$
(resp. $\cdot_Q$) is the  intersection product in $P$ (resp. $Q$).
\end{proposition}

\noindent{\bf Proof.}\quad By Proposition~\ref{prop:C.1}, the
restriction of a bounded linear Fredholm operator to a closed
subspace of finite codimension is also Fredholm. So from the fact
that $Y$ has finite codimension in $X$ it easily follows that $(Y,
E|_Y, S|_Y)$ is also a Banach Fredholm bundle. As before we can take
an open neighborhood ${\cal W}$ of $Z(S)$ in $X$, $\varepsilon>0$
and smooth sections $\sigma_i:X\to E$, $i=1,\cdots,m$, such that for
each ${\bf t}\in B_\varepsilon(\R^m)_{res}$ both sections
\begin{eqnarray*}
&&\Phi_{\bf t}|_{\cal W}: {\cal W}\to E,\; y\mapsto
S(y)+\sum^m_{i=1}t_i\sigma_i(y),\\
&&\Phi^Y_{\bf t}: {\cal W}\cap Y\to E|_Y,\; y\mapsto
S|_Y(y)+\sum^m_{i=1}t_i\sigma_i|_Y(y),
\end{eqnarray*}
are transversal to the zero section. Note that ${\cal W}\cap Y$ is
an open neighborhood of $Z(S|_Y)=Z(S)\cap Y$ in $Y$. Then
$(\Phi_{\bf t}|_{\cal W})^{-1}(0)\cap Y=(\Phi^Y_{\bf t})^{-1}(0)$
and thus $M=(\Phi_{\bf t}|_{\cal W})^{-1}(0)$ and $N=(\Phi^Y_{\bf
t})^{-1}(0)$ are the desired Euler chains.

In particular, suppose that $Y=f^{-1}(Q)$ is as above. Let
$a^\ast_Q$ be a differential form representative of the Poincar\'e
dual $a\in H_\ast(Q, \R)$ and $R_Q$ be a retraction from a tubular
neighborhood of $Q$ in $P$ onto $Q$. If $Q^\ast$ is a representative
form of the Poincar\'e dual of $Q$ in $P$ whose support is contained
in the above tubular neighborhood, then
$a^\ast:=R^\ast_Q(a^\ast_Q)\wedge Q^\ast$ is a representative form
of the Poincar\'e dual of $a\in H_\ast(Q, \R)\subset H_\ast(P, \R)$
in $P$. So
\begin{eqnarray*}
 a\cdot_Q f_\ast(e(E|_Y,S|_Y))\!\!\!\!\!\!\!\!&&=\int_{f(N)}a^\ast_Q\\
 &&=\int_{f(M)\cap Q}a^\ast_Q\\
 &&=\int_{R_Q(f(M)\cap Q)}a^\ast_Q\\
 &&=\int_{f(M)\cap Q}R^\ast_Q(a^\ast_Q)\\
  &&=\int_{f(M)}R^\ast_Q(a^\ast_Q)\wedge Q^\ast\\
  &&=\int_{f(M)}a^\ast=a\cdot_P f_\ast(e(E, S)).
\end{eqnarray*}
 Here the fifth equality is obtained by jiggling $f|_M$ so that it is transversal to $Q$.
  The desired result is proved. \hfill$\Box$\vspace{2mm}

\subsection{Properties }\label{sec:1.3}

 Our first property is a stability result, which comes from  Proposition 14 in [Bru].
 But we here requires   slightly strong assumptions because
 our proof actually give the relations of the corresponding
  Euler cycles.

\begin{proposition}\label{prop:1.13}
({\bf Stablity}) Let $(X, E, S)$ be a Banach Fredholm bundle of
index $r$ and with the Euler class $e(E,S)$.  Assume that there
exists a smooth direct sum decomposition of Banach bundles
$E=E'\oplus E''$ over a neighborhood ${\cal W}$ of $Z(S)$. Let
${\rm P}_{E'}$ and ${\rm P}_{E''}$ be the  fibrewise projections
onto $E'$ and $E''$ respectively. For the natural sections
$$S'={\rm P}_{E'}\circ S:{\cal W}\to E'|_{\cal W}\quad{\rm and}\quad
S''={\rm P}_{E''}\circ S: {\cal W}\to E''|_{\cal W}
$$
($Z(S)=Z(S')\cap Z(S'')$ and $Z(S)=Z(S'|_{Z(S'')})$ are clear),
 suppose that $DS''|_{Z(S)}$ is surjective, i.e. for any $x\in Z(S)$ the
vertical differential $DS''(x): T_xX\to E''_x$ is onto. Then
\begin{description}
\item[(a)] for some small open neighborhood ${\cal U}$ of $Z(S)$
in $X$ the intersection $Z(S'')^\star:=Z(S'')\cap{\cal U}$ is a
smooth manifold;

\item[(b)] $D(S'|_{Z(S'')}):(TZ(S''))|_{Z(S)} \to E'|_{Z(S)}$ is a
Fredholm bundle map of index ${\rm Ind}(D(S'|_{Z(S'')})={\rm
Ind}(DS)$ and there exists a closed manifold representative $N$ of
$e(E'|_{Z(S'')^\star}, S'|_{Z(S'')^\star})$ and that $M$ of
$e(E,S)$ such that $M=N$. In particular it implies
$$
e(E,S) = (i_{Z(S'')^\star})_\ast e(E'|_{Z(S'')^\star},
S'|_{Z(S'')^\star}),
$$
where $(i_{Z(S'')^\star})_\ast$ is the homomorphism between the
homology groups induced by the inclusion $i_{Z(S'')^\star}:
Z(S'')^\star\hookrightarrow X$.
\end{description}
\end{proposition}

\noindent{\bf Proof.}\quad The conclusion (a) is obvious.
  For $x\in Z(S)=Z(S')\cap
Z(S'')$, it is easily checked that the linear Fredholm operator
$DS(x): T_xX\to E_x$ may be decomposed into $DS(x)=DS'(x)\oplus
DS''(x): T_xX\to E_x'\oplus E''$. So the induced operator
$\widehat{DS(x)}:T_xX\to E_x/E_x'$ ( as the composition of $DS(x)$
and the quotient map $E_x\to E_x/E_x'$ ) may be identified with
$DS''(x)$ since we identify $E_x''\equiv E_x/E_x'$. Then we can
identify ${\rm Ker}(\widehat{DS(x)})$ with ${\rm Ker}(DS''(x))=
(DS(x))^{-1}(E_x')$ and the map
$$
\widetilde{DS(x)}: {\rm Ker}(\widehat{DS(x)})\to E_x',\; v\mapsto
DS(x)(v),
$$
with one
$$D(S'|_{Z(S'')})(x): T_xZ(S'')\to E_x'$$
because the surjectivity of $DS''(x)$ implies that $T_xZ(S'')={\rm
Ker}(DS''(x))$. By Proposition~\ref{prop:C.3} in
Appendix~\ref{app:C} we get that $D(S'|_{Z(S'')})(x)$ is a Fredholm
operator with index ${\rm Ind}(D(S'|_{Z(S'')})(x))={\rm
Ind}(DS(x))$.

Next we prove the second claim in (b). Since $Z(S)$ is compact  we
may take finitely many points $x_1,\cdots, x_m$ in $Z(S)$ and their
open neighborhoods $O_j$ in $X$ which can be required in ${\cal U}$,
$j=1,\cdots,m$ such that $Z(S)\subset {\cal
O}(Z(S)):=\cup^m_{j=1}O_j$ and that there exist trivializations
$\psi_j: O_j\times E_{x_j}\to E|_{O_j}$ so that for any $y\in O_j$
the maps $\psi_{jy}=\psi_j(y,\cdot): E_y\to E_{x_j}$ preserve  the
splitting $E_y=E_y'\oplus E_y''\to E_{x_j}=E_{x_j}'\oplus
E_{x_j}''$, i.e., $\psi_{jy}(E_y')=E_{x_j}'$ and
$\psi_{jy}(E_y'')=E_{x_j}''$. Let $S'_j:O_j\to E'_{x_j}$ (resp.
$S''_j: O_j\to E''_{x_j}$) be trivialization representatives of
$S'|_{O_j}$ (resp. $S''|_{O_j}$), i.e., $\psi_j(y, S'_j(y))=S'(y)$
(resp. $\psi_j(y, S''_j(y))=S''(y)$, $j=1,\cdots,m$. Then
$S_j:=S_j'+ S_j''$ are trivialization representatives of $S|_{O_j}$,
i.e., $\psi_j(y, S_j(y))=S(y)$, $j=1,\cdots,m$. Note that
$DS''|_{Z(S)}$ is surjective. As before we may choose finitely many
sections $\sigma_1',\cdots,\sigma_k'$ of $E'$ such that their
trivialization representatives under $\psi_j$,
$$
\sigma'_{ij}, i=1,\cdots,k,\; j=1,\cdots, m,
$$
 satisfy:
\begin{eqnarray}
&&\left.\begin{array}{ll} dS'_j(y)(T_yZ(S''))+ {\rm
span}(\{\sigma'_{1j}(y),\cdots,
\sigma'_{kj}(y)\}=E'_{x_j}\\
\qquad\forall y\in Z(S'')\cap O_j,\;j=1,\cdots,m,
\end{array}\right\}\label{e:1.11}\\
 &&\left.\begin{array}{ll}
 dS_j(y)(T_yX)+ {\rm span}(\{\sigma'_{1j}(y),\cdots,
\sigma'_{kj}(y)\}=E_{x_j}\\
 \qquad\forall y\in O_j,\;
j=1,\cdots,m.\end{array}\right\}\label{e:1.12}
\end{eqnarray}
 For ${\bf t}\in\R^k$ consider the
sections
\begin{eqnarray*}
&&\Phi'_{{\bf t}}: {\cal O}(Z(S))\to E'|_{{\cal O}(Z(S))},
 \quad x\mapsto S'(x)+\sum^k_{i=1}t_i\sigma'_i(x),\\
&&\Phi_{\bf t}: {\cal O}(Z(S))\to E|_{{\cal O}(Z(S))}, \quad
x\mapsto S(x)+\sum^k_{i=1}t_i\sigma'_i(x).
\end{eqnarray*}
As before it follows from (\ref{e:1.11}) and (\ref{e:1.12}) that
exist a small open neighborhood ${\cal W}\subset{\cal O}(Z(S))$ of
$Z(S)$ in $X$ and  $\varepsilon>0$ such that for the open
neighborhood ${\cal N}:={\cal O}(Z(S))\cap{\cal W}$ of $Z(S)$ in
$Z(S'')^\star$ and each ${\bf t}\in B_\varepsilon(\R^k)_{res}$ the
sections
$$
\Phi_{\bf t}|_{\cal W}: {\cal W}\to E|_{\cal W}\quad{\rm and}\quad
\Phi'_{\bf t}|_{\cal N}: {\cal N}\to E'|_{\cal N}
$$
are all transversal to the zero section. Let $R_{\bf t}=\Phi'_{\bf
t}|_{\cal N}$.  Note that
\begin{equation}\label{e:1.13}
{\rm P}_{E'}\circ\Phi_{\bf t}=\Phi_{\bf t}',\;\;{\rm
P}_{E''}\circ\Phi_{\bf t}=S''\;\;{\rm and}\;\; {\rm
P}_{E'}\circ\Phi_{\bf t}+{\rm P}_{E''}\circ\Phi_{\bf t}=\Phi_{\bf
t}.
\end{equation}
We claim that for ${\bf t}\in B_\varepsilon(\R^k)_{res}$ small
enough,
\begin{equation}\label{e:1.14}
(\Phi_{\bf t}|_{\cal W})^{-1}(0)=R_{\bf t}^{-1}(0).
\end{equation}
Then the desired result follows from it.

 To prove (\ref{e:1.14})  note that
 ${\rm P}_{E'}\circ\Phi_{\bf
t}(y)=S'(y)+ \sum^m_{i=1}t_i\sigma_i'(y)=R_{\bf t}(y)$ for $y\in
Z(S'')\cap{\cal W}$. So ${\rm P}_{E'}\circ\Phi_{\bf t}(y)=0$ for any
$y\in R_{\bf t}^{-1}(0)$. It follows from (\ref{e:1.13}) that for
any $y\in R_{\bf t}^{-1}(0)\subset Z(S'')\cap{\cal W}$, $\Phi_{\bf
t}(y)=0$ and thus
$$
R_{\bf t}^{-1}(0)\subset(\Phi_{\bf t}|_{\cal W})^{-1}(0).
$$
 On the other hand,  (\ref{e:1.13}) implies that
\begin{eqnarray*}
(\Phi_{\bf t}|_{\cal W})^{-1}(0)\!\!\!\!\!\!&&\subset  ({\rm
P}_{E'}\circ\Phi_{\bf t}|_{\cal W})^{-1}(0)\cap({\rm P}_{E''}\circ
S)^{-1}(0)\\
&&=({\rm P}_{E'}\circ\Phi_{\bf t}|_{\cal W})^{-1}(0)\cap Z(S'').
\end{eqnarray*}
By the definition of $R_{\bf t}$ we get $(\Phi_{\bf t}|_{\cal
W})^{-1}(0)\subset R_{\bf t}^{-1}(0)$ and thus (\ref{e:1.14}). The
desired result is proved.
 \hfill$\Box$\vspace{2mm}

In the proof above it is important for us to assume that there
exists a smooth direct sum decomposition of Banach bundles
$E=E'\oplus E''$ over a neighborhood ${\cal W}$ of $Z(S)$. As a
consequence we get the following special case of Proposition 2.8
in [R1].

\begin{corollary}\label{cor:1.14}
 Let $(X, E, S)$ be as in Proposition~\ref{prop:1.13}. Assume:
\begin{description}
\item[(i)] $\dim{\rm Coker}DS(y)=k$ on $Z(S)$, $Z(S)$ is a closed
smooth manifold of dimension $r+ k$ and thus ${\rm Coker}DS$ forms
an obstruction bundle ${\mathcal E}$ over $Z(S)$.

\item[(ii)] There exist an open neighborhood ${\cal W}$ of $Z(S)$
and a decomposition of the direct sum of Banach bundles
$$E|_{{\cal W}}=F\oplus F^c$$
such that $F^c|_{Z(S)}$ is isomorphic to ${\cal E}$ and that the
Banach subbundle $F$ of corank $k$  restricts to ${\rm Im}(DS)$ on
$Z(S)$, i.e., $F_y={\rm Im}(DS(y))$ for any $y\in Z(S)$.
\end{description}
Then there exist  closed manifold representatives $M$ of $e(E,S)$
and $N$ of $e({\cal E})$ such that $M=N\cap Z(S)$, and in
particular $ e(E,S)=e({\mathcal E})\cap [Z(S)]$.
 \end{corollary}

\noindent{\bf Proof.}\quad Note first that $e(F^c|_{Z(S)})=e({\cal
E})$ because $F^c|_{Z(S)}$ is isomorphic to ${\cal E}$. Let ${\rm
P}_F: E|_{{\cal W}}\to F$ and ${\rm P}_{F^c}: E|_{{\cal W}}\to
F^c$ be clear bundle projections.
  By the  assumption above, the section ${\rm P}_F\circ S: {\cal W}\to F$
 is transversal to the zero section at each point $y\in Z(S)$. The desired
 conclusions follow from Proposition~\ref{prop:1.13}.
 \hfill$\Box$\vspace{2mm}

\begin{proposition}\label{prop:1.15}
({\bf Product}) For $i=1,2$ let $(X_i, E_i, S_i)$ be
 Banach Fredholm bundles of index $r_i$ and
with compact zero locus $Z(S_i)$. Then the natural product
$(X_1\times X_2, E_1\times E_2, S_1\times S_2)$ is such a bundle
of index $r=r_1+ r_2$ and with zero locus $Z(S_1\times
S_2)=Z(S_1)\times Z(S_2)$ and  for $i=1,2$ there exist  compact
submanifolds $M_i$ of $X_i$, which are the Euler chains of $(X_i,
E_i, S_i)$,  such that the product manifold $M=M_1\times M_2$ is a
Euler chain of $(X_1\times X_2, E_1\times E_2, S_1\times S_2)$. If
 both $e(E_1, S_1)$ and $e(E_2, S_2)$ exist then
$e(E_1\times E_2, S_1\times S_2)$ also exist and the compact
manifolds $M_i$ can be chosen to be without boundary.
\end{proposition}

\noindent{\bf Proof.}\quad Let $S:=S_1\times S_2$. Then
$S((x_1,x_2))=(S_1(x_1), S_2(x_2))$ for $x=(x_1, x_2)\in X_1\times
X_2$.  One easily checks that $S$ is also a Fredholm section and
${\rm Ind}(S)={\rm Ind}(S_1)+ {\rm Ind}(S_2)$.

 Let open neighborhoods ${\cal W}\subset{\cal O}(Z(S))$ of $Z(S)$ in $X_1\times X_2$, $\varepsilon>0$
 and the sections $\sigma_j=\sigma_j^{(1)}\times\sigma_j^{(2)}: X_1\times X_2\to
E_1\times E_2$ with supports in ${\cal O}(Z(S))$, $j=1,\cdots, m$,
 be such that Theorem~\ref{th:1.5} holds for them. Let
 $$\Phi^{(i)}: X_i\times B_\varepsilon(\R^m)\to\Pi^\ast_1E_i,\; x_i\mapsto
 S_i(x_i)+ \sum^m_{j=1}t_j\sigma^{(i)}_j(x_i).$$
Then $\Phi=\Phi_1^{(1)}\times\Phi_2^{(2)}$ and
$Z(\Phi)=Z(\Phi^{(1)})\times Z(\Phi^{(2)})$. For $i=1,2$ let us take
open neighborhoods ${\cal W}_i$ of $Z(S_i)$ in $X_i$ so that ${\cal
W}_1\times{\cal W}_2\subset{\cal W}$. Note that for $x=(x_1,x_2)\in
Z(\Phi)$ the vertical differential $D\Phi(x): T_x(X_1\times X_2)\to
E_{1x_1}\times E_{2x_2}$ is onto if and only if both
$D\Phi^{(1)}(x_1): T_{x_1}X_1\to E_{1x_1}$ and $D\Phi^{(2)}(x_2):
T_{x_2}X_2\to E_{2x_2}$ are onto. By Theorem~\ref{th:1.5}(A), by
shrinking $\varepsilon>0$ we can assume that $(\Phi_{\bf t}|_{\cal
W})^{-1}(0)\subset{\cal W}_1\times{\cal W}_2$ for any ${\bf t}\in
B_\varepsilon(\R^m)$. So for each ${\bf t}\in
B_\varepsilon(\R^m)_{res}$ the sections $\Phi_{\bf t}|_{\cal W}$ and
$\Phi_{\bf t}^{(i)}|_{{\cal W}_i}$, $i=1,2$, are all transversal to
the zero section. In particular we get
$$(\Phi_{\bf t}|_{\cal W})^{-1}(0)=(\Phi_{\bf t}^{(1)}|_{{\cal W}_1})^{-1}(0)\times
(\Phi_{\bf t}^{(2)}|_{{\cal W}_2})^{-1}(0).
$$
This completes the proof of Proposition~\ref{prop:1.15}.
  \hfill$\Box$\vspace{2mm}

Two Banach Fredholm bundles $(X, E^{(i)}, S^{(i)})$, $i=0,1$,
 of index $r$ and with compact zero locus  are called {\bf homotopy}
 if there exists a Banach Fredholm bundle $(X\times [0,1], E, S)$
 of index $r+1$ and with compact zero locus such that
 $E^{(i)}=E|_{\{i\}\times X}$ and $S^{(i)}=S|_{\{i\}\times X}$ for
 $i=0, 1$. The homotopy is called an {\bf oriented homotopy} if
 $(X\times [0,1], E, S)$ and $(X, E^{(i)}, S^{(i)})$, $i=0,1$, are also oriented
 and the orientation of $(X\times [0,1], E, S)$ induces those of
$(X, E^{(i)}, S^{(i)})$, $i=0,1$.

\begin{proposition}\label{prop:1.16}
({\bf Homotopy}) If two (oriented) Banach Fredholm bundles $(X,
E^{(i)}, S^{(i)})$, $i=0,1$,
 of index $r$ and with compact zero locus  are (oriented)
 homotopic,
 then they have the same Euler class.
 \end{proposition}

\noindent{\bf Proof.}\quad By the arguments before
Theorem~\ref{th:1.5} we may choose  points
\begin{eqnarray*}
&&(x_{0i}, 0)\in Z(S),\; i=1,\cdots, m_0,\\
&&(x_{1i}, 1)\in Z(S),\; i=1,\cdots, m_1,\\
&&(y_j, t_j)\in Z(S)\cap (X\times (0, 1)), j=1,\cdots, m_2
\end{eqnarray*}
and their open neighborhoods in $X\times [0, 1]$, $O_{0i}$,
$i=1,\cdots, m_0$, $O_{1i}$, $i=1,\cdots, m_1$ and $O_j$,
$j=1,\cdots, m_2$, and smooth sections $s_{0i}$, $i=1,\cdots,
n_0$, $s_{1i}$, $i=1,\cdots, n_1$, $s_j$, $j=1,\cdots, n_2$ such
that:
\begin{description}
\item[(i)] $Z(S)\subset
(\cup^{m_0}_{i=1}O_{0i})\cup(\cup^{m_1}_{k=1}O_{1k})\cup(\cup^{m_2}_{j=1}O_j)$,

\item[(ii)] $Z(S^{(0)})\subset
\cup^{m_0}_{i=1}O_{0i}\cap(X\times\{0\})$,

\item[(ii)] $Z(S^{(1)})\subset
\cup^{m_1}_{i=1}O_{1i}\cap(X\times\{1\})$,

\item[(iii)] each $O_j$ is contained in $X\times (0, 1)$,
$j=1,\cdots, m_2$,

\item[(iv)] the support of each $s_{0i}$ (resp. $s_{1i}$, $s_j$)
is contained in some $O_{0k}$ (resp. $O_{1l}$, $O_s$),

\item[(v)] ${\rm span}\{s_{01}(p),\cdots, s_{0n_0}(p),
s_{11}(p),\cdots, s_{1n_1}(p), s_1(p),\cdots, s_{n_2}(p)\}\\
+ {\rm Im}(DS(p))=E_p$ for any $p\in Z(S)$,

\item[(vi)] ${\rm span}\{s_{01}(x),\cdots, s_{0n_0}(x)\}+{\rm
Im}(DS^{(0)}(x))=E^{(0)}_x$ for any $x\in Z(S^{(0)})$,

\item[(vii)] ${\rm span}\{s_{11}(x),\cdots, s_{1n_1}(x)\}+{\rm
Im}(DS^{(1)}(x))=E^{(1)}_x$ for any $x\in Z(S^{(1)})$.
\end{description}
It follows that there exists $\varepsilon>0$, the residual subsets
$B_\varepsilon(\R^{n_0})_{res}\subset B_\varepsilon(\R^{n_0})$,
$B_\varepsilon(\R^{n_1})_{res}\subset B_\varepsilon(\R^{n_1})$,
$B_\varepsilon(\R^{n})_{res}\subset B_\varepsilon(\R^{n})$ with
$n=n_0+ n_1 + n_2$, and a small open neighborhood ${\cal W}$ of
$Z(S)$ in $X\times [0, 1]$ such that for any
\begin{eqnarray*}
&&{\bf t}^0=(t_{01},\cdots, t_{0n_0})\in
B_\varepsilon(\R^{n_0})_{res},\\
&&{\bf t}^1=(t_{11},\cdots, t_{1n_1})\in
B_\varepsilon(\R^{n_1})_{res},\\
&&{\bf r}=(r_{01},\cdots, r_{0n_0}, r_{11},\cdots, r_{1n_1},
r_1,\cdots, r_{n_2})\in B_\varepsilon(\R^{n})_{res},
\end{eqnarray*}
the sections
\begin{eqnarray*}
\Phi_{\bf r}: {\cal W}\to E|_{\cal W},\;p\mapsto
S(p)+\sum^{n_0}_{i=1}r_{0i}s_{0i}(p)+
\sum^{n_1}_{k=1}r_{1k}s_{1k}(p)+\sum^{n_1}_{j=1}r_{j}s_{j}(p),\\
\Phi^{(0)}_{{\bf t}^0}: {\cal W}^0={\cal W}\cap(X\times\{0\})\to
E^{(0)}|_{{\cal W}^0},\;x\mapsto
S^{(0)}(x)+\sum^{n_0}_{i=1}t_{0i}s_{0i}(x), \\
\Phi^{(1)}_{{\bf t}^1}:{\cal W}^1={\cal W}\cap(X\times\{1\})\to
E^{(1)}|_{{\cal W}^1},\;x\mapsto
S^{(1)}(x)+\sum^{n_1}_{k=1}t_{1k}s_{1k}(x)
\end{eqnarray*}
are all transversal to the zero section. By Lemma~\ref{lem:1.7} we
can take ${\bf r}^0=(r_{01},\cdots, r_{0n_0})\in
B_\varepsilon(\R^{n_0})_{res}$, ${\bf r}^1=(r_{11},\cdots,
r_{1n_1})\in B_\varepsilon(\R^{n_1})_{res}$ such that ${\bf
r}=({\bf r}^0, {\bf r}^1, r_1,\cdots, r_{n_2})\in
B_\varepsilon(\R^{n})_{res}$ some $r_1,\cdots, r_{n_2}\in\R$. Then
for this ${\bf r}$ it is easily checked that
$$
\partial(\Phi_{\bf r})^{-1}(0)=(\Phi^{(0)}_{{\bf r}^0})^{-1}(0)\cup
(-(\Phi^{(1)}_{{\bf r}^1})^{-1}(0)).
$$
This implies the desired result. \hfill$\Box$\vspace{2mm}

Finally we study the functoriality of the Euler classes of Banach
Fredholm bundles. We shall give two kinds of results in two
propositions. In some senses  they might be viewed as
generalizations of Proposition~\ref{prop:1.12}.

Let $(X, E, S)$ and $(X', E', S')$ be two oriented Banach Fredholm
bundles with compact zero loci.  A {\bf morphism} from $(X, E, S)$
to $(X', E', S')$ is a pair $(f, F)$ with following properties:
\begin{description}
\item[(i)] $f$ is a smooth embedding and $F$ is a smooth injective
bundle homomorphism covering $f$, i.e. $\forall x\in X$ the
restriction $F_x: E_x\to E'_{f(x)}$ is a continuous linear
injective map;

\item[(ii)] $S\circ f=F\circ S'$ and $Z(S')=f(Z(S))$;

\item[(iii)] For any $x\in Z(S)$ the differential $df(x):T_xX\to
T_{f(x)}X'$ and the above restriction $F_x$ induce isomorphisms
\begin{eqnarray*}
&&df(x): {\rm Ker}(DS(x))\to {\rm Ker}(DS'(f(x)))\quad{\rm and}\\
&&[F_x]: {\rm Coker}(DS(x))\to {\rm Coker}(DS'(f(x))),
\end{eqnarray*}
and the resulting isomorphism from ${\rm det}(DS)$ to ${\rm
det}(DS')$ is orientation preserving.
\end{description}

By the definition  $(X, E, S)$ and $(X', E', S')$ have the same
index.

\begin{proposition}\label{prop:1.17}
Let $(X, E, S)$ and $(X', E', S')$ be two oriented Banach Fredholm
bundles with compact zero loci, and $(f, F)$ be a morphism from
$(X, E, S)$ to $(X', E', S')$. Then there exist Euler chains $N$
of $(X, E, S)$, and $M$ of $(X', E', S')$, which are also compact
manifolds of the same dimension, such that $f(N)=M\cap f(X)$.
Moreover, if $e(E,S)$ and $e(E', S')$ exist then $M$ and $N$ can
also be chosen closed manifolds. In particular,  $\langle
f^\ast\alpha, e(E,S)\rangle=\langle\alpha, e(E', S')\rangle$ for
any $\alpha\in H^\ast(X', X'\setminus f(X);\R)$. (Since $M$ and
$N$ have the same dimension this means that any connected
component of $M$ either is disjoint with $f(X)$ or is contained in
$f(X)$. All those components contained in $f(X)$ form $f(N)$.)
\end{proposition}

\noindent{\bf Proof.}\quad By the definition of the morphism we
may assume: $X$ is a Banach submanifold of $X'$, $E$ is a
subbundle of $E'|_X$ and $S=S'|_X$. Of course, both $f$ and $F$
are inclusion maps. By (ii) and (iii) it also holds that
$Z(S')=Z(S)$, and that for any $x\in Z(S)$,
$DS'(x)|_{T_xX}=DS(x)$, ${\rm Ker}(DS(x))={\rm Ker}(DS'(x))$ and
the inclusion $E_x\hookrightarrow E'_x$ induces an isomorphism
\begin{eqnarray*}
 &&E_x/DS(x)(T_xX)\to E_x'/DS'(x)(T_xX'),\\
&&  v+ DS(x)(T_xX)\mapsto v + DS'(x)(T_xX').
\end{eqnarray*}
The final claim implies: if $v_i\in E_x, i=1,\cdots, k$ such that
$DS(x)(T_xX)+ {\rm span}\{v_1,\cdots, v_k\}=E_x$ then
$DS'(x)(T_xX')+ {\rm span}\{v_1,\cdots, v_k\}=E'_x$. So we can
choose smooth sections of $E'$, $\sigma'_1,\cdots, \sigma'_m$ with
supports near $Z(S')$ such that $\sigma_1=\sigma'_1|_X, \cdots,
\sigma_m=\sigma'_m|_X$ are smooth sections of $E$ and that
\begin{description}
\item[(I)]  the sections
\begin{eqnarray*}
&&\Phi: ({\cal W}\cap X)\times
B_\varepsilon(\R^m)\to\Pi^\ast_1E,\, (y, {\bf
t})\mapsto S(y)+\sum^m_{i=1}t_i\sigma_i(y)\\
&&\Phi': {\cal W}\times B_\varepsilon(\R^m)\to\Pi^\ast_1E',\, (y,
{\bf t})\mapsto S'(y)+\sum^m_{i=1}t_i\sigma'_i(y)
\end{eqnarray*}
are Fredholm and transversal to the zero section,

\item[(II)] the zero sets $\Phi^{-1}(0)$ and $(\Phi')^{-1}(0)$
have compact closures in $Cl(({\cal W}\cap X)\times
B_\varepsilon(\R^m))$ and $Cl({\cal W}\times B_\varepsilon(\R^m))$
respectively.
\end{description}
Here $\varepsilon>0$ is very small and ${\cal W}$ is an open
neighborhood of $Z(S')$ in $X'$. As usual there exists a residual
subset $B_\varepsilon(\R^m)_{res}\subset B_\varepsilon(\R^m)$ such
that for each ${\bf t}\in B_\varepsilon(\R^m)_{res}$ one gets
compact ${\rm Ind}(S)={\rm Ind}(S')$ dimensional manifolds
$(\Phi_{\bf t})^{-1}(0)$ and $(\Phi'_{\bf t})^{-1}(0)$. Clearly,
$(\Phi_{\bf t})^{-1}(0)\subset (\Phi'_{\bf t})^{-1}(0)$.
 Note that one of two connected closed manifolds with the same dimension
cannot contain another. When $e(E,S)$ and $e(E', S')$ exist,
$(\Phi_{\bf t})^{-1}(0)$ must consist of those connected
components of $(\Phi'_{\bf t})^{-1}(0)$ which are contained in
$X$. The desired results are proved. \hfill$\Box$\vspace{2mm}

\begin{proposition}\label{prop:1.18}
({\bf Pull-back}) Let $(X, E, S)$ be a Banach Fredholm bundle of
index $r$ and with compact zero locus $Z(S)$. If $f:Y\to X$  is a
proper Fredholm map of index $d$ from another Banach manifold $Y$
to $X$ then the natural pullback $(Y, f^\ast E, f^\ast S)$ is also
a Banach Fredholm bundle with compact zero locus $Z(f^\ast S)$.
Moreover, its index equals to $r+d$ and there exists a Euler chain
$M$ (resp. $N$) of $(X, E, S)$ (resp. $(Y, f^\ast E, f^\ast S)$)
such that $f^{-1}(M)=N$. If both Euler classes $e(E,S)$ and
$e(f^\ast E, f^\ast S)$ exist then $M$ and $N$ can be chosen as
closed manifolds.
\end{proposition}

If $d=0$ and $X$ is connected the final claim implies
$$
f_\ast(e(f^\ast E, f^\ast S))={\rm deg}(f) e(E,S)
$$
and thus $\langle f^\ast\alpha, e(f^\ast E, f^\ast S)\rangle={\rm
deg}(f)\langle\alpha, e(E,S)\rangle$ for any $\alpha\in
H^\ast(X,\R)$. Here ${\rm deg}(f)$ is understand as follows: If we
do not consider the orientation of $(X, E, S)$ and thus $e(E,S)\in
H_r(X;\Z_2)$ and $e(f^\ast E, f^\ast S)\in H_r(Y;\Z_2)$ then ${\rm
deg}(f)$ denotes the topological degree of $f$ mod-2; If $(X, E, S)$
is oriented and one considers the oriented Euler classes $e(E,S)\in
H_r(X;\Z)$ and $e(f^\ast E, f^\ast S)\in H_r(Y;\Z)$ then  the notion
of some kind ``orientation'' of $f$ shall be needed to define the
$\Z$-value topological degree ${\rm deg}(f)$ of $f$. They shall be
studied in the future. \vspace{2mm}

\noindent{\bf Proof of Proposition~\ref{prop:1.18}.}\quad Note that
$Z(f^\ast S)=f^{-1}(Z(S))$ is compact because $Z(S)$ is compact and
$f$ is proper. Moreover, for $x_0\in O\subset X$ and a
trivialization $\psi:O\times E_{x_0}\to E|_O$ one has, for a given
$y_0\in f^{-1}(x_0)$, a natural induced trivialization
$$ f^\ast\psi: f^{-1}(O)\times (f^\ast E)_{y_0}\to (f^\ast
E)|_{f^{-1}(O)},\; (y, v)\mapsto \psi(f(y), v).
$$
If $S_\psi:O\to E_{x_0}$ the representation of $S|_O$ under the
trivialization $\psi$ then $f^\ast S$ has, under $f^\ast\psi$, the
corresponding trivialization representation
\begin{equation}\label{e:1.15}
(f^\ast S)_{f^\ast\psi}: f^{-1}(O)\to (f^\ast E)_{y_0},\; y\mapsto
S_\psi(f(y)).
\end{equation}
It follows that as the composition of Fredholm maps  $(f^\ast
S)_{f^\ast\psi}$ is Fredholm and ${\rm Ind}((f^\ast
S)_{f^\ast\psi})={\rm Ind}(f)+ {\rm Ind}(S_\psi)=d+r$. That is,
$(Y, f^\ast E, f^\ast S)$ is a Banach Fredholm bundle of index
$r+d$. By (\ref{e:1.15}),
$$
d(f^\ast S)_{f^\ast\psi}(y)=dS_\psi(f(y))\circ df(y)\quad\forall
y\in f^{-1}(O).
$$
 So if there exist finite
elements $v_i\in E_{f(y_0)}=E_{x_0}$, $i=1,\cdots, m$, such that
$${\rm Im}(d(f^\ast S)_{f^\ast\psi}(y))+ {\rm span}\{v_1,\cdots,
v_m\}=E_{x_0},
$$
then it easily follows that
$${\rm Im}(dS_\psi(f(y))+ {\rm span}\{v_1,\cdots,
v_m\}=E_{x_0}
$$
because ${\rm Im}(d(f^\ast S)_{f^\ast\psi}(y))\subset {\rm
Im}(dS_\psi(f(y))$.

Let points $x_i\in Z(S)$ and their $E$-trivialized open
neighborhoods $O_i$ in $X$, $i=1,\cdots, n$, and smooth sections
 of $E$ with support in ${\cal
O}(Z(S))):=\cup^n_{i=1}O_i$, $\sigma_i$, $i=1,\cdots, m$, be such
that
\begin{equation}\label{e:1.16}
{\rm Im}(dS_{\psi_i}(x))+ {\rm span}\{\sigma_{1i}(x),\cdots,
\sigma_{mi}(x)\}=E_{x_i}\;\forall x\in O_i,
\end{equation}
where $S_{\psi_{x_i}}, \sigma_{ji}: O_i\to E|_{O_i}$ are
representations of $S$, $\sigma_j$ under the trivializations
$\psi_i:O_i\times E_{x_i}\to E|_{O_i}$, $i=1,\cdots, n$ and
$j=1,\cdots, m$.  Taking $y_i\in f^{-1}(x_i)$, $i=1,\cdots, n$
then we have
\begin{equation}\label{e:1.17}
\left.\begin{array}{ll}
 {\rm Im}(d(f^\ast S)_{f^\ast\psi_i}(y))+
{\rm span}\{(f^\ast\sigma_1)_i(y),\cdots,
(f^\ast\sigma_m)_i)(y)\}\\
=(f^\ast E)_{y_i}\;\forall y\in f^{-1}(O_i)
\end{array}
\right\}
 \end{equation}
  Here $(f^\ast S)_{f^\ast\psi_{x_i}},
(f^\ast\sigma_j)_i: f^{-1}(O_i)\to (f^\ast E)|_{f^{-1}(O_i)}$ are
representations of $f^\ast S$, $f^\ast\sigma_j$ under the above
pullback trivializations $f^\ast\psi_i: f^{-1}(O_i)\times (f^\ast
E)_{y_i}\to (f^\ast E)|_{f^{-1}(O_i)}$, $i=1,\cdots, n$ and
$j=1,\cdots, m$. Consider the sections
\begin{eqnarray*}
&&\Phi: X\times B_\varepsilon(\R^m)\to E,\; (x, {\bf t})\mapsto
S(y)+\sum^m_{i=1}t_i\sigma_i(y),\\
&&f^\ast\Phi: Y\times B_\varepsilon(\R^m)\to f^\ast E,\; (y,{\bf
t})\mapsto f^\ast S(y)+\sum^m_{i=1}t_i (f^\ast\sigma_i)(y).
\end{eqnarray*}
 Then it follows from
(\ref{e:1.16}) and (\ref{e:1.17}) that for any $x\in Z(S)$ and
$y\in Z(f^\ast S)$ the maps
\begin{eqnarray*}
&&D\Phi(x, 0): T_xX\times \R^m\to E_x\quad{\rm and}\\
&&D(f^\ast\Phi)(y,0): T_yY\times \R^m\to (f^\ast E)_y
 \end{eqnarray*}
are onto. Thus there exist an open neighborhood ${\cal W}$ of
$Z(S)$ in $X$, that ${\cal W}^\ast\subset f^{-1}({\cal W})$ of
$Z(f^\ast S)$ in $Y$,  $\varepsilon>0$, and a residual subset
$B_\varepsilon(\R^m)_{res}\subset B_\varepsilon(\R^m)$ such that
for each ${\bf t}\in B_\varepsilon(\R^m)_{res}$
 the restriction of $\Phi_{\bf t}$ to ${\cal W}$ and that of
 $(f^\ast\Phi)_{\bf t}=f^\ast(\Phi_{\bf t})$ to ${\cal
W}^\ast$ are transversal to the zero section. So the sets
$(\Phi_{\bf t}|_{\cal W})^{-1}(0)$ and $((f^\ast\Phi)_{\bf
t}|_{{\cal W}^\ast})^{-1}(0)$
 are respectively closed
manifolds of dimensions $r$ and $r+d$ for small ${\bf t}\in
B_\varepsilon(\R^m)_{res}$. Clearly,
$$((f^\ast\Phi)_{\bf t}|_{{\cal
W}^\ast})^{-1}(0)\subset f^{-1}((\Phi_{\bf t}|_{\cal W})^{-1}(0)).
$$
Moreover, since $f$ is proper it is easily proved that for any open
neighborhood ${\cal U}$ of $Z(f^\ast S)$ in $Y$ there exists an open
neighborhood ${\cal V}$ of $Z(S)$ in $X$ such that $f^{-1}({\cal
V})\subset{\cal U}$. It follows from this and
Theorem~\ref{th:1.5}(A) that $f^{-1}((\Phi_{\bf t}|_{\cal
W})^{-1}(0))\subset{\cal W}^\ast$ for ${\bf t}\in
B_\varepsilon(\R^m)_{res}$ small enough. For such a ${\bf t}$, one
easily checks that $ f^{-1}((\Phi_{\bf t}|_{\cal W})^{-1}(0))\subset
((f^\ast\Phi)_{\bf t}|_{{\cal W}^\ast})^{-1}(0)$. So  we get
\begin{equation}\label{e:1.18}
f^{-1}(\Phi_{\bf t}^{-1}(0))=(f^\ast\Phi)_{\bf t}^{-1}(0)
\end{equation}
for ${\bf t}\in B_\varepsilon(\R^m)_{res}$ small enough. Clearly,
if both $e(E, S)$ and $e(f^\ast E, f^\ast S)$ exist then
$(\Phi_{\bf t}|_{\cal W})^{-1}(0)$ and $((f^\ast\Phi)_{\bf
t}|_{{\cal W}^\ast})^{-1}(0)$ are closed manifolds for ${\bf t}\in
B_\varepsilon(\R^m)_{res}$ sufficiently small. Note that $e(f^\ast
E, f^\ast S)=[((f^\ast\Phi)_{\bf t}|_{{\cal W}^\ast}^{-1}(0)]$ and
$e(E, S)=[f^{-1}((\Phi_{\bf t}|_{\cal W})^{-1}(0))]$. The desired
results are proved.  \hfill$\Box$\vspace{2mm}

\pagebreak

\section{The virtual Euler cycles of Banach Fredholm orbibundles}\label{sec:2}
\setcounter{equation}{0}

As stated in Introduction the purpose of this section is to give a
way in which we can construct an analogue of the Euler cycle in \S1,
called {\bf virtual Euler cycle}, for the Banach Fredholm
orbibundles.

\subsection{Banach orbifolds and
orbibundles}\label{sec:2.1}

 Orbifolds were first defined and studied by Satake in
[Sat1] and [Sat2] (he used the term ``V-manifold", the term
``orbifold" is due to Thurston \cite{Thu}). Before introducing our
definition of a Banach orbifold we recall a basic result due to
Newman in the theory of compact transformation groups. If a finite
group $\Gamma$ acts effectively on a connected manifold $N$ of
finite dimension,  then Newman theorem \cite[Theorem 1]{Dr} claims
that the set of points with trivial isotropy group is open and
everywhere dense. (Here the effective action means that the
intersection of all isotropy groups, $\cap_{x\in N}\Gamma_x$, only
consists of the unit element $\1$.)
 It follows that for any $x\in N$ the restriction
action of the isotropy group $\Gamma_x$ to a $\Gamma_x$-invariant
connected open neighborhood of $x$ is also effective. However,
when $N$  is of infinite dimension, Newman theorem fails, for such
a fact we can neither know how to prove it nor find a place where
this fact is proved. Thus a notion of an induced chart cannot be
defined in the effective category. This essential difference
between finite dimension and infinite dimension shows that the
effectness condition is weaker in the infinitely dimensional case
and  will not be able to be used to complete  some important
arguments as in the case of finite dimension if we completely
imitate the definition of orbifolds of finite dimension to define
orbifolds of infinite dimension. Our strategies are to drop the
effectness condition in a definition of Banach orbifold charts of
infinite dimension and then add stronger conditions when
necessary.

\begin{definition}\label{def:2.1}
{\rm Let $X$ be a paracompact Hausdorff topological space. A {\bf
Banach} (or {\bf B-}) {\bf orbifold chart} on $X$ is a triple
$(\widetilde U, \Gamma_U,\pi_U)$, where $U$ (resp. $\widetilde U$)
is a connected open subset of $X$ (resp. a Banach manifold),
$\Gamma_U$ is a finite group which acts on $\widetilde U$ by
$C^\infty$-automorphisms of $\widetilde U$, and $\pi_U$ is a
continuous surjective map from $\widetilde U$ to $U$ such that for
any $x\in\widetilde U$ and $g\in \Gamma_U$, $\pi_U(gx)=\pi_U(x)$,
and that the induced map ${\widetilde U}/\Gamma_U\to U$ is a
homeomorphism. (One often call $\widetilde U$ as a {\bf local
cover}, $\Gamma_U$ as a {\bf local group}, $U$ the {\bf support of
the chart}, and $\pi_U$ as a {\bf local covering map}).}
\end{definition}

Note that $\pi_U$ is always  {\bf proper} by Proposition 1.7 on the
page 102 of \cite{B}. In fact, a map from a Hausdorff space into
another Hausdorff space is proper if it is closed and the pre-image
of every point is also compact.
 For any $x\in U$ and $\tilde x\in (\pi_U)^{-1}(x)$ let
$\Gamma_U(\tilde x)$ be the {\bf isotropy group} of $\Gamma_U$ at
$\tilde x$, i.e., $\Gamma_U(\tilde x)=\{g\in\Gamma_U\,|\,
g\cdot\tilde x=\tilde x\}$. Clearly, if $\tilde y\in
(\pi_U)^{-1}(x)$ is another element, then there exists
$g\in\Gamma_U$ such that $g\tilde y=\tilde x$. It follows that
$\Gamma_U(\tilde x)\to\Gamma_U(\tilde y),\;h\mapsto g^{-1}hg$, is a
group isomorphism. If
 $$
 \Gamma_{U}(\tilde{x})={\rm
Ker}(\Gamma_U,\widetilde U):=\{g\in\Gamma_U\,|\, gx=x\;\forall
x\in\widetilde U\},
$$
  $\tilde{x}$ is called a {\bf regular point},
otherwise it is called a {\bf singular point}. Denote
$\widetilde{U}^{\circ}$ (resp. $\widetilde{U}^{sing}$) by the set of
all regular (resp. singular) points of $\widetilde{U}$.  We also
denote $\Gamma_{U}\tilde{x}$ by the orbit of $\tilde{x}$, i.e.,
$\Gamma_{U}\tilde{x}=\{g(\tilde{x})\,|\,g\in \Gamma_{U} \}$.

The following lemmas come from Lemmas 1.1, 1.2 in \cite{LuW}. For
convenience we also write its proof.

\begin{lemma}\label{lem:2.2}
 Let $(\widetilde{U}, \Gamma_{U},\pi_{U})$ be a Banach orbifold chart as above.
Then for any $x\in U$ and $\tilde{x}\in(\pi_{U})^{-1}(x)$, there
exists a connected open neighborhood
$\widetilde{O}(\tilde{x})\subset\widetilde{U}$ of $\tilde{x}$ such
that \begin{description}
\item[(i)] $\widetilde{O}(\tilde{x})$ is
$\Gamma_{U}(\tilde{x})$-invariant, and $\widetilde{O}(\tilde{x})\to
\widetilde{O}(\tilde{x}),\;\tilde y\mapsto g\cdot\tilde y$ is a
homeomorphism for any $g\in \Gamma_U(\tilde x)$;
\item[(ii)]
$h(\widetilde{O}(\tilde{x}))\cap\widetilde{O}(\tilde{x})=\emptyset$
 for any $h\in \Gamma_{U}\setminus \Gamma_{U}(\tilde{x})$;
 \item[(iii)] $(\widetilde{O}(\tilde{x}), \Gamma_{U}(\tilde{x}),\pi_{U}^{x}=
\pi_{U}|_{\widetilde{O}(\tilde{x})})$ and
$(g\widetilde{O}(\tilde{x}), \Gamma_{U}(g\tilde{x}),\pi_{U}^{gx}=
\pi_{U}|_{g\widetilde{O}(\tilde{x})})$ for any $g\in\Gamma_U$ are
Banach orbifold charts on $X$; Moreover,
$h(g\widetilde{O}(\tilde{x}))\cap
g\widetilde{O}(\tilde{x})=\emptyset$
 for any $h\in \Gamma_{U}\setminus \Gamma_{U}(g\tilde{x})$.
 \end{description}
\end{lemma}

The Banach orbifold chart $(\widetilde{O}(\tilde{x}),
\Gamma_{U}(\tilde{x}),\pi_{U}^{x}=
\pi_{U}|_{\widetilde{O}(\tilde{x})})$ in Lemma~\ref{lem:2.2} is
called an {\bf induced chart} of
 $(\widetilde{U}, \Gamma_{U}, \pi_{U})$ at $\tilde{x}$ (or $x$).
It follows from Lemma~\ref{lem:2.2} that the connected component of
$\pi_U^{-1}(O(x))$ containing $\tilde x$ is
$\widetilde{O}(\tilde{x})$ because $\pi_U^{-1}(O(x))$ is equal to
the union of $\widetilde{O}(\tilde{x})$ and
$\cup_{g\in\Gamma_U\setminus\Gamma_U(\tilde
x)}g\cdot\widetilde{O}(\tilde{x})$. \vspace{2mm}

\noindent{\bf Proof of Lemma~\ref{lem:2.2}}\quad
 If $\Gamma_{U}(\tilde{x})=\Gamma_{U}$, we can take $\widetilde{O}(\tilde{x})
=\widetilde{U}$. So we may assume $\Gamma_{U}(\tilde{x})\neq
\Gamma_{U}$. For any $h\in \Gamma_{U}\setminus
\Gamma_{U}(\tilde{x}),\ h(\tilde{x})\neq\tilde{x}$,  we can always
find an open neighborhood $V_{h}$ of $\tilde{x}$  such that
$h(V_{h})\cap V_{h}=\emptyset$. Then $V_{1}:=\cap_{h\in
\Gamma_{U}\setminus \Gamma_{U}(\tilde{x})}V_{h}$ is also an open
neighborhood of $\tilde{x}$ and satisfies $h(V_{1})\cap
V_{1}=\emptyset$
 for any $h\in \Gamma_{U}\setminus \Gamma_{U}(\tilde{x})$. Since $g(\tilde{x})=\tilde{x}$
 for any $g\in \Gamma_{U}(\tilde{x})$, there must exist an open neighborhood $V_{g}$
of $\tilde{x}$  such that $g(V_{g})\subset V_{1}$. Set
$V_{2}:=\cap_{g\in \Gamma_{U}(\tilde{x})}V_{g}$, we have
$g(V_2)\subset V_1$ for any $g\in \Gamma_{U}(\tilde{x})$. Take a
connected open neighborhood $V$ of $\tilde{x}$ such that $V\subset
V_{1}\cap V_{2}$. Then $\widetilde{O}(\tilde{x}):=\cup_{g\in\
\Gamma_{U}(\tilde{x})}g(V)$ is also a connected open neighborhood of
$\tilde{x}$ and satisfies (i) obviously. To see that it also
satisfies (ii), note that for any $h\in \Gamma_{U}\setminus
\Gamma_{U}(\tilde{x})$,
\begin{eqnarray*}
h(\widetilde{O}(\tilde{x}))\cap\widetilde{O}(\tilde{x})\!\!\!\!&&=
(\bigcup\limits_{g\in \Gamma_{U}(\tilde{x})}hg(V))\bigcap
(\bigcup\limits_{g'\in \Gamma_{U}(\tilde{x})}g'(V))\\
&&=\bigcup\limits_{g\in \Gamma_{U}(\tilde{x})} \bigcup\limits_{g'\in
\Gamma_{U}(\tilde{x})}(hg(V)\cap g'(V)).
\end{eqnarray*}
 Since $hg\in \Gamma_{U}\setminus
\Gamma_{U}(\tilde{x})$, $V\subset V_1$ implies $hg(V)\cap
V_{1}=\emptyset$, and $V\subset V_{2}$ implies $g'(V)\subset
V_{1}$. Hence $hg(V)\cap g'(V)=\emptyset$, and thus
$h(\widetilde{O}(\tilde{x}))\cap\widetilde{O}(\tilde{x})=\emptyset$.

For (iii),   it is clear that
$O(x)=\pi_{U}({\widetilde{O}(\tilde{x})})$ is an open subset in $U$
(and so in $X$) since $\pi_{U}$ is an open map. Next, for any $g\in
\Gamma_{U}(\tilde{x})$ and $\tilde{y}\in\widetilde{O}(\tilde{x})$,
we have $\tilde{x}_{\tilde{y}}:=g^{-1}(\tilde{y})\in
\widetilde{O}(\tilde{x})$ and $g(\tilde{x}_{\tilde{y}})=\tilde{y}$
because $g^{-1}\in \Gamma_{U}(\tilde{x})$.  It implies
$g|_{\widetilde{O}(\tilde{x})}: {\widetilde{O}(\tilde{x})}
\longrightarrow{\widetilde{O}(\tilde{x})}$ is a surjective map and
thus a homeomorphism. Finally, we show that
$\pi_{U}^{x}:{\widetilde{O}(\tilde{x})}\to O(x)$ induces the
following homeomorphism
$$
\overline{\pi_{U}^{x}}:{\widetilde{O}(\tilde{x})}/\Gamma_{U}(\tilde{x})
\to O(x), \;[\tilde{y}]_{\Gamma_{U}(\tilde{x})}\mapsto
\pi_{U}^{x}(\tilde{y})=\pi_{U}(\tilde{y}).
$$
where $\tilde y\in\widetilde{O}(\tilde x)$,
$[\tilde{y}]_{\Gamma_{U}(\tilde{x})}:=\Gamma_{U}(\tilde x)\tilde
y$. In fact, it suffices to prove that $\overline{\pi_{U}^{x}}$ is
injective and open. Assume
$\overline{\pi_{U}^{x}}\,(\,[\tilde{y}_{1}]_{\Gamma_{U}(\tilde{x})})
=\overline{\pi_{U}^{x}}\,(\,[\tilde{y}_{2}]_{\Gamma_{U}(\tilde{x})})$
for $\tilde{y}_{1},\ \tilde{y}_{2}\in\widetilde{O}(\tilde{x})$.
Then $\pi_{U}^{x}(\tilde{y}_{1})=\pi_{U}(\tilde{y}_{1})
=\pi_{U}(\tilde{y}_{2}) =\pi_{U}^{x}(\tilde{y}_{2})$, and thus
there exists a $g_{0}\in \Gamma_{U}$ such that $\tilde{y}_{1}=
g_{0}(\tilde{y}_{2})$. Note that
$h(\widetilde{O}(\tilde{x}))\cap\widetilde{O}(\tilde{x})=\emptyset$
for any $h\in \Gamma_{U}\setminus \Gamma_{U}(\tilde{x})$. We get
$g_{0}\in \Gamma_{U}(\tilde{x})$ and
$[\tilde{y}_{1}]_{\Gamma_{U}(\tilde{x})}
=[\tilde{y}_{2}]_{\Gamma_{U}(\tilde{x})}$.

Let $A$ be the open subset of
${\widetilde{O}(\tilde{x})}/\Gamma_{U}(\tilde{x})$ and
$q_{U}^{x}:\widetilde{O}(\tilde{x})\longrightarrow\widetilde{O}(\tilde{x})
/\Gamma_{U}(\tilde{x})$ be the quotient map. Then
$(q_{U}^{x})^{-1}(A)$ is open in $\widetilde{O}(\tilde{x})$ (and so
in $\widetilde{U}$). Since $\pi_{U}:\widetilde{U}\rightarrow U$ is
an open map, $\pi_{U}((q_{U}^{x})^{-1}(A))$ is an open subset of
$U$. Note that $\pi_{U}((q_{U}^{x})^{-1}(A))=
\pi_{U}^{x}((q_{U}^{x})^{-1}(A))\subset O(x)$. We get that
$\pi_{U}^{x}((q_{U}^{x})^{-1}(A))$ is an open subset of $O(x)$.
Moreover,
$\overline{\pi_{U}^{x}}(A)=\pi_{U}^{x}((q_{U}^{x})^{-1}(A))$.
 It follows that $\overline{\pi_{U}^{x}}$ is open.

For any $g\in\Gamma_U$,  $(g\widetilde{O}(\tilde{x}),
\Gamma_{U}(g\tilde{x}),\pi_{U}^{gx}=
\pi_{U}|_{g\widetilde{O}(\tilde{x})})$ is also a Banach orbifold
chart on $X$ since $\Gamma_U(g\tilde x)=\{ghg^{-1}\,|\,
h\in\Gamma_U(\tilde x)\}$. \vspace{2mm} $\Box$

 For two Banach orbifold charts $(\widetilde U, \Gamma_U,\pi_U)$ and $(\widetilde V,
\Gamma_V,\pi_V)$ on $X$, an {\bf injection} from $(\widetilde U,
\Gamma_U,\pi_U)$ to $(\widetilde V, \Gamma_V,\pi_V)$ is a pair
$\theta_{UV}=(\tilde\theta_{UV}, \gamma_{UV})$, where
$\gamma_{UV}:\Gamma_U\to \Gamma_V$ is an injective group
homomorphism and $\tilde \theta_{UV}:\widetilde U\to \widetilde V$
is a $\gamma_{UV}$-equivariant open embedding such that
$\pi_U=\pi_V\circ\tilde \theta_{UV}$ and the following {\bf
maximality condition} holds:
$$
{\rm
Im}(\gamma_{UV})=\{g\in\Gamma_V\,|\,\tilde\theta_{UV}(\widetilde
U)\cap g\cdot\tilde\theta_{UV}(\widetilde U)\ne\emptyset\}.
$$
The last condition implies that $\gamma_{UV}$ induces an isomorphism
from ${\rm Ker}(\Gamma_U, \widetilde U)$ to ${\rm
Ker}(\Gamma_V,\widetilde V)$, where ${\rm Ker}(\Gamma_U,\widetilde
U)=\{g\in\Gamma_U\,|\, gx=x\;\forall x\in\widetilde U\}$ and ${\rm
Ker}(\Gamma_V,\widetilde V)=\{g\in\Gamma_V\,|\, gx=x\;\forall
x\in\widetilde V\}$.   In particular, if $U=V$, $\gamma_{UV}$ is a
group isomorphism, and $\widetilde\theta_{UV}$ is a diffeomorphism
from $\widetilde U$ onto $\widetilde V$, then these two B-orbifold
charts are called {\bf equivalent}, and $\theta_{UV}$ is called an
{\bf equivalence} between them. Clearly, $(\widetilde{O}(\tilde{x}),
\Gamma_{U}(\tilde{x}),\pi_{U}^{x}=
\pi_{U}|_{\widetilde{O}(\tilde{x})})$ and
$(g\widetilde{O}(\tilde{x}), \Gamma_{U}(g\tilde{x}),\pi_{U}^{gx}=
\pi_{U}|_{g\widetilde{O}(\tilde{x})})$, $g\in\Gamma_U$ are
equivalent Banach orbifold charts on $X$.

\begin{definition}\label{def:2.3}
 {\rm Let $X$ be a paracompact Hausdorff topological space.
  A {\bf  Banach orbifold atlas} on $X$ is
a collection ${\mathcal A}$ of B-orbifold charts on $X$ satisfying
the following properties:
\begin{description}
\item[(i)] The supports of all charts in ${\mathcal A}$ form a
basis for open sets in $X$.

\item[(ii)] For any two charts $(\widetilde U, \Gamma_U,\pi_U)$
and $(\widetilde V, \Gamma_V,\pi_V)$ in ${\mathcal A}$ with
$U\subset V$ there exists an   injection $\theta_{UV}=(\tilde
\theta_{UV},\gamma_{UV})$ from $(\widetilde U, \Gamma_U,\pi_U)$ to
$(\widetilde V, \Gamma_V,\pi_V)$.
\end{description}}
\end{definition}

 Let ${\mathcal A}$ and
${\mathcal A}^\prime$ be two Banach orbifold atlases on $X$. If for
each chart $(\widetilde U, \Gamma_U,\pi_U)$ in ${\mathcal A}$ there
exists a chart $(\widetilde V, \Gamma_V,\pi_V)$ in ${\mathcal
A}^\prime$ and an injection from $(\widetilde U, \Gamma_U,\pi_U)$ to
$(\widetilde V, \Gamma_V,\pi_V)$ we say ${\mathcal A}$ to be a {\bf
refinement} of ${\mathcal A}^\prime$. Such two atlases on $X$ are
said to be {\bf equivalent} if they have a common refinement. It may
be proved that this is indeed an equivalence relation among the
Banach orbifold atlases on $X$. Denote by $[{\mathcal A}]$ the
equivalent class of the atlas ${\mathcal A}$ on $X$. It is called a
{\bf Banach orbifold structure} on $X$ and the pair $(X, [{\mathcal
A}])$ is called a {\bf  Banach orbifold}, usually only denoted by
$X$ without occurring of confusions. Any atlas in
 $[{\mathcal A}]$ is called an {\bf atlas} of $(X, [{\mathcal A}])$.
We say that $(\widetilde U, \Gamma_U,\pi_U)$ is a {\bf  Banach
orbifold chart of} $(X, [{\mathcal A}])$ if it belongs to some atlas
${\mathcal B}$ in $[{\mathcal A}]$.  Later we always assume that the
atlas ${\cal A}$ is the maximal one in the equivalence class
containing it.  A Banach orbifold $(X, [{\mathcal A}])$ is called an
{\bf effective Banach orbifold} if $[{\mathcal A}]$ contains an
atlas
 ${\mathcal A}$ such that each chart $(\widetilde U,
 \Gamma_U,\pi_U)$ is {\bf effective} in the sense that $\Gamma_U$ acts on
 $\widetilde U$ effectively. In this case the atlas ${\mathcal A}$ is called
 an {\bf effective Banach orbifold atlas}. It is easily shown that every Banach
 orbifold $X$ induces an effective Banach orbifold.

 Let ${\mathcal A}$ be a Banach orbifold atlas on $X$,   $(\widetilde
U, \Gamma_U,\pi_U)$ and $(\widetilde V, \Gamma_V,\pi_V)$ be two
charts in ${\mathcal A}$ with $U\subset V$ and $\theta_{UV}=(\tilde
\theta_{UV},\gamma_{UV})$ be an injection  from $(\widetilde U,
\Gamma_U,\pi_U)$ to $(\widetilde V, \Gamma_V,\pi_V)$. For $x\in U$
and $\tilde x\in\pi_U^{-1}(x)$ the maximality condition implies that
$\gamma_{UV}$ induces an isomorphism from $\Gamma_U(\tilde x)$ to
$\Gamma_V(\tilde \theta_{UV}(\tilde x))$. Thus using
Definition~\ref{def:2.3} every $x\in X$ determines a group
$\Gamma_x$, unique up to isomorphism, which is isomorphic to the
isotropy group of any lift of $x$ in any chart of ${\mathcal A}$.
(Sometimes $\Gamma_x$ is simply called the {\bf isotropy group} of
$x$.) The orders of these groups define a locally constant function
on $X$, which is called the {\bf multiplicity function}
$\mathfrak{mul}_X:X\to\N$. Given any finite group $G$, the connected
components of the set of all $x$ such that $\Gamma_x\cong G$ are
smooth Banach manifolds. This leads to a decomposition $X=\cup_iX_i$
having the properties: $X_i$ is a connected Banach manifold and
consists of points of a fixed isotropy group, $\Gamma_i$; moreover
for $i\ne j$, $X_i$ meets the closure of $X_j$ only if $\Gamma_i$ is
an abstract subgroup of $\Gamma_i$. Each connected component of $X$
has a unique open connected stratum (called principal stratum)
$X_\ast$ on which the corresponding isotropy group $\Gamma_\ast$ is
minimal, i.e. $\sharp\Gamma_\ast=\min\mathfrak{mul}_X$. A point
$x\in X$ is said to be {\bf regular} (resp. {\bf singular}) if its
lift in some chart of ${\mathcal A}$ is regular (resp. singular).
Denote by $X^\circ$ (resp. $X^{sing}$) the set of regular (resp.
singular) points in $X$. Clearly,  $X^\circ$ is equal to the union
of all principal strata of $X$.

If each local cover $\widetilde U$ in the definition of Banach
orbifolds above is replaced by a connected open subset of some
Banach manifold with boundary (resp. corner) and the action of
$\Gamma_U$ is required to preserve the boundary (resp. corner) we
obtain the {\bf Banach orbifold with boundary (resp. corner)}. If
each local cover $\widetilde U$ in the above definition is a
connected open subset of some separable Banach manifold we say $X$
to be a {\bf separable}  Banach orbifold. It is not hard to prove
that if $X$ is a Banach orbifold with boundary then its boundary
$\partial X$ also inherits a  Banach orbifold structure from $X$ and
becomes a
 Banach orbifold.

A topological subspace $Z$ of a  Banach orbifold $X$ is called a
{\bf  Banach suborbifold} of $X$ if it is a  Banach orbifold with
respect to the induced  Banach orbifold structure obtained as
follows: For each  Banach orbifold chart $(\widetilde U,
\Gamma_U,\pi_U)$ of $X$ with $U\cap Z\ne\emptyset$ there exists a
Banach submanifold $\widetilde Z_U\subset\widetilde U$ that is not
only stable under $\Gamma_U$ but also compatible with the inclusion
maps, such that the restriction $(\widetilde Z_U,
\Gamma_U|_{\widetilde Z_U},\pi_U|_{\widetilde Z_U})$ is a  Banach
orbifold chart for $Z$.

 A  Banach orbifold $X$ of finite dimension is
called {\bf locally oriented} if it has an atlas  ${\cal A}$ such
that for each chart $(\widetilde U, \Gamma,\pi_U)$ in ${\cal A}$,
$\widetilde U$ is {\bf oriented} and each element of $\Gamma_U$ is
orientation preserving. In this case we say the chart $(\widetilde
U, \Gamma,\pi_U)$ to be {\bf oriented}.  $X$ is called {\bf
oriented} if  it has an atlas  ${\cal A}$ such that  each chart
$(\widetilde U, \Gamma,\pi_U)$ in ${\cal A}$ is oriented, and that
for any two charts  $(\widetilde U, \Gamma_U,\pi_U)$ and
$(\widetilde V, \Gamma_V,\pi_V)$ in ${\mathcal A}$ with $U\subset V$
the injection $\theta_{UV}=(\tilde \theta_{UV},\gamma_{UV})$ from
$(\widetilde U, \Gamma_U,\pi_U)$ to $(\widetilde V, \Gamma_V,\pi_V)$
is orientation preserving.

\begin{lemma}\label{lem:2.4}
Let $(X, {\mathcal A})$ be a Banach orbifold,  $(\widetilde U,
\Gamma_U,\pi_U)$ be a chart in ${\mathcal A}$, and
$(\widetilde{O}(\tilde{x}), \Gamma_{U}(\tilde{x}),\pi_{U}^{x}=
\pi_{U}|_{\widetilde{O}(\tilde{x})})$ be the induced chart at
$\tilde{x}\in\pi_U^{-1}(x)$ as in Lemma~\ref{lem:2.2}. Then\\
{\rm (i)} Any chart $(\widetilde W, \Gamma_W,\pi_W)$ in ${\mathcal
A}$ with $x\in W\subset O(x)$ has the properties: the injection
$\theta_{WU}=(\tilde \theta_{WU},\gamma_{WU})$ from $(\widetilde W,
\Gamma_W,\pi_W)$ to $(\widetilde U, \Gamma_U,\pi_U)$ induces an
isomorphism from $\Gamma_W$ to $\Gamma_U(g\tilde x)$ for some
$g\in\Gamma_U$. So the chart $(\widetilde W, \Gamma_W,\pi_W)$ is
equivalent to the chart $(\tilde\theta_{WU}(\widetilde W),
\Gamma_{U}(g\tilde{x}), \pi_{U}|_{\tilde\theta_{WU}(\widetilde W)})$
and thus to the chart $(g^{-1}\tilde\theta_{WU}(\widetilde W),
\Gamma_{U}(\tilde{x}), \pi_{U}|_{g^{-1}\tilde\theta_{WU}(\widetilde
W)})$. ({\bf Note that the last two charts are also the induced
charts of $(\widetilde U, \Gamma_U, \pi_U)$ at $g\tilde x$ and
$\tilde x$ respectively.}) Hence for any given $x\in U$, $\tilde
x\in\pi_U^{-1}(U)$ and a neighborhood $V$ of $x$ there exists an
induced chart of $(\widetilde U, \Gamma_U, \pi_U)$ at $\tilde{x}$,
$(\widetilde{O}(\tilde{x}), \Gamma_{U}(\tilde{x}),\pi_{U}^{x}=
\pi_{U}|_{\widetilde{O}(\tilde{x})})$ which is equivalent to a chart
in ${\cal A}$ and  such that $O(x)=\pi_U(\widetilde O(\tilde
x))\subset V$. This chart $(\widetilde{O}(\tilde{x}),
\Gamma_{U}(\tilde{x}),\pi_{U}^{x})$ is also effective if ${\mathcal
A}$ is an effective atlas.\\
{\rm (ii)} For finite charts in ${\cal A}$, $(\widetilde U_i,
\Gamma_i,\pi_i)$, $i=1,\cdots, n$, if $x\in\cap^n_{i=1}U_i$ and
$\tilde x_i\in\pi_i^{-1}(x)$, $i=1,\cdots,n$, then there exist
induced charts of $(\widetilde U_i, \Gamma_i,\pi_i)$ at $\tilde
x_i$, $(\widetilde O(\tilde x_i), \Gamma_i(\tilde x_i),
\pi_i|_{\widetilde O(\tilde x_i)})$, $i=1,\cdots,n$ having the
same support contained in a given neighborhood of $x$, and
equivalences from $(\widetilde O(\tilde x_1), \Gamma_1(\tilde
x_1), \pi_1|_{\widetilde O(\tilde x_1)})$ to $(\widetilde O(\tilde
x_i), \Gamma_i(\tilde x_i), \pi_i|_{\widetilde O(\tilde x_i)})$
that map $\tilde x_1$ to $\tilde x_i$, $i=2,\cdots, n$. Of course,
these induced chart are also effective if the atlas ${\mathcal A}$
is effective.
\end{lemma}

\noindent{\bf Proof}.\quad (i) Since $W\subset O(x)\subset U$, it
follows from the definition of the injection that
$\tilde\theta_{WU}(\widetilde W)\subset \pi_U^{-1}(O(x))$. By
Lemma~\ref{lem:2.2} and the arguments below Lemma~\ref{lem:2.2} we
may assume that $\tilde\theta_{WU}(\widetilde W)\subset\widetilde
O(\tilde x_u)$ for some $\tilde x_u\in\pi_U^{-1}(x)$. Let $\tilde
x_w\in\pi_W^{-1}(x)$ be such that $\tilde\theta_{WU}(\tilde
x_w)=\tilde x_u$. By the maximality condition each
$g\in\Gamma_U(\tilde x_u)$ sits in $\tilde\gamma_{WU}(\Gamma_W)$.
Moreover, for any $h\in\Gamma_W$ we have $\tilde\gamma_{WU}(h)\tilde
x_u=\tilde\gamma_{WU}(h)\tilde\theta_{WU}(\tilde
x_w)=\tilde\theta_{WU}(h\tilde x_w)\subset\widetilde O(\tilde x_u)$
and
$$
\tilde\gamma_{WU}(h)\tilde
x_u=\tilde\gamma_{WU}(h)\tilde\theta_{WU}(\tilde
x_w)=\tilde\theta_{WU}(h\tilde x_w)\in\pi_U^{-1}(x)
$$
because $\pi_W=\pi_U\circ\tilde\theta_{WU}$.
 It follows from Lemma~\ref{lem:2.2} that
$\tilde\gamma_{WU}(h)\tilde x_u=\tilde x_u$. This shows that
$\tilde\gamma_{WU}(\Gamma_W)=\Gamma_U(\tilde x_u)$ and thus
$\tilde\gamma_{WU}$ induces an isomorphism from $\Gamma_W$ to
$\Gamma_U(\tilde x_u)$. Other claims are easily seen from the
arguments above.

\noindent{(ii)} For the remainder we only need to prove the final
claim. For a given neighborhood $V$ of $x$, since all supports of
charts in ${\cal A}$ form a topology basis for $X$ we can find a
chart $(\widetilde W, \Gamma_W, \pi_W)$ in ${\cal A}$ such that
$W\subset V$. As above we can also construct induced charts of
$(\widetilde U_i, \Gamma_i,\pi_i)$ at $\tilde x_i$, $(\widetilde
O(\tilde x_i), \Gamma_i(\tilde x_i), \pi_i|_{\widetilde O(\tilde
x_i)})$, $i=1,\cdots,n$ such that they are all equivalent to
$(\widetilde W, \Gamma_W, \pi_W)$ and that the equivalence from
$(\widetilde O(\tilde x_i), \Gamma_i(\tilde x_i),
\pi_i|_{\widetilde O(\tilde x_i)})$ to $(\widetilde W, \Gamma_W,
\pi_W)$ can be chosen to map all $\tilde x_1,\cdots, \tilde x_n$
to the same $\tilde x_w\in\pi_W^{-1}(x)$. The desired claim
follows easily. \vspace{2mm} $\Box$

Now let us study the structures of $X^{sing}$ and $X^\circ$. For
points $x\in U\cap X^{sing}$ and  $\tilde{x}\in\pi_U^{-1}(x)$,   it
follows from Theorem~\ref{th:B.1} that there exists a
$\Gamma_U(\tilde x)$-equivariant diffeomorphism $F$ from a
neighborhood  $N(0_{\tilde x})$
 of the origin in $T_{\tilde x}\widetilde U$ onto a neighborhood
 of $\tilde x$ in $\widetilde U$, ${\cal N}(\tilde x)=F(N(0_{\tilde
 x}))$. That is, $F(dg(\tilde x)\tilde
 v)=g\cdot F(\tilde v)$ for any $g\in\Gamma_{U}(\tilde
 x)$ and $\tilde v\in N(0_{\tilde x})$. We may shrink
$N(0_{\tilde x})$ so that $h({\cal N}(\tilde x))\cap{\cal N}(\tilde
x)=\emptyset$ for any $h\in\Gamma_U\setminus\Gamma_U(\tilde x)$.
This implies that $\Gamma_U(\tilde y)=\{g\in\Gamma_U(\tilde x)\,|\,
g\tilde y=\tilde y\}$ for any $\tilde y\in{\cal N}(\tilde x)$.
 Let ${\nu}(dg(\tilde x))={\rm Ker}({\1}-dg(\tilde x))$. It is easily proved that
 \begin{equation}\label{e:2.1}
 \pi_U^{-1}(X^{sing})\cap {\cal N}(\tilde
 x)=\bigcup_{g\in\Gamma_{U}(\tilde x)
 \setminus{\rm Ker}(\Gamma_U, \widetilde U)}F\bigl(N(0_{\tilde x})\cap{\nu}(dg(\tilde x))\bigr).
 \end{equation}
So for each $g\in\Gamma_{U}(\tilde x)\setminus{\rm
Ker}(\Gamma_U,\widetilde U)$ the submanifold $F\bigl(N(0_{\tilde
x})\cap{\nu}(dg(\tilde x))\bigr)\subset\widetilde U$, which is
relatively closed in $\widetilde U$,  is exactly the set of fixed
points of $g$ in ${\cal N}(\tilde x)$, denoted by ${\cal N}(\tilde
x)^g$, and ${\cal N}(\tilde{x})^{sing}= \cup_{g\in \Gamma_U(\tilde
x)\setminus{\rm Ker}(\Gamma_U,\widetilde U)}{\cal
N}(\tilde{x})^{g}$. We define {\bf singular codimension} of
$\widetilde U$ near $\tilde x\in\widetilde U^{sing}$ and that of $X$
near $x\in X^{sing}$, denoted by
$$
{\rm codim}_{\tilde x}\widetilde U^{sing}\quad\hbox{and}\quad{\rm
codim}_xX^{sing},
$$
to be $\min\{{\rm codim} {\cal N}(\tilde{x})^{g}\,|\, g\in
\Gamma_U(\tilde x)\setminus{\rm Ker}(\Gamma_U,\widetilde U)\}$.
Clearly, the functions $$\tilde x\mapsto {\rm codim}_{\tilde
x}\widetilde U^{sing}\quad\hbox{and}\quad x\mapsto {\rm
codim}_xX^{sing}
$$
are upper semi-continuous. As it happens, it is possible that
$dg(\tilde x)=id$ for some $g\in\Gamma_{U}(\tilde x)
 \setminus{\rm Ker}(\Gamma_U,\widetilde U)$. So in this case $\widetilde U^{sing}$ contains an
 open neighborhood of $\tilde x$ in $\widetilde U$. This shows that
 $\widetilde U^\circ$ cannot be dense in $\widetilde U$.
When ${\rm codim}_xX^{sing}\ge 1$ at any $x\in X^{sing}$, it is
easily seen that the regular point set $X^{\circ}$, is a dense open
subset in $X$ (in fact is a complementary set of a sparse subset of
$X$). For integer $k\ge 1$, we say
 the Banach orbifold chart $(\widetilde{U},\Gamma_{U},\pi_{U})$ to be
{\bf $k$-regular} if
 the codimension of $\widetilde U^{sing}$ near any $\tilde{x}\in \widetilde U^{sing}$
 is not less than $k$. It means that the fixed point set $\widetilde
U^g$ of any $g\in \Gamma_U\setminus{\rm Ker}(\Gamma_U,\widetilde U)$
has codimension not less than $k$ near any $\tilde{x}\in \widetilde
U^g$. So for any $k$-regular Banach orbifold chart $(\widetilde{U},
\Gamma_{U},\pi_{U})$, $\widetilde{U}^{\circ}$ is a dense open subset
in $\widetilde{U}$. A Banach orbifold $X$ is called a {\bf
$k$-regular  Banach orbifold} if every orbifold chart of it is
$k$-regular. In this case $X^{sing}$ is a sparse subset of $X$ and
$X^\circ=X\setminus X^{sing}$.  Summarizing the arguments above we
get:

\begin{claim}\label{cl:2.5}
 For a $1$-regular Banach orbifold $X$, $X^\circ$ is a dense open subset in $X$ and $X^{sing}$ is the
union of finitely many suborbifolds of positive codimension {\rm
locally}.
\end{claim}

  If $X$ is a $1$-regular oriented orbifold of finite dimension then it is also $2$-regular
  because $1$ can only be an
eigenvalue of $dg(\tilde x)$ of even algebraic  multiplicity.
Moreover,
 any $1$-regular complex orbifold is $2$-regular. From the
 definitions above it is easily seen that an effective Banach orbifold
 is not necessarily $1$-regular, and conversely any $k$-regular
 Banach orbifold is not necessarily an effective Banach orbifold.
In Remark 1.4 of \cite{Sie} it was also shown how to construct a
Banach orbifold structure on a Hausdorff topological space $X$ with
a family of $2$-regular Banach orbifold charts on $X$ whose supports
cover $X$.

Another important notion is a smooth map between orbifolds.

\begin{definition}\label{def:2.6}
{\rm For $0<k\le\infty$, a continuous map $f$ from a Banach orbifold
$X$ to another Banach orbifold $X^\prime$ is said to be $C^k$-{\bf
smooth} if there exist Banach orbifold atlases ${\cal A}_X$ on $X$
and ${\cal A}_{X'}$ on $X'$ such that:
\begin{description}
\item[(i)] There is a correspondence
$$
{\cal A}_X\to {\cal A}_{X'},\;(\widetilde U, \Gamma_U,\pi_U)\to
(\widetilde U^\prime, \Gamma_{U^\prime},\pi_{U^\prime})
$$
satisfying: $f(U)\subset U^\prime$ and there exist a homomorphism
$\varphi_U: \Gamma_U\to \Gamma_{U^\prime}$ and a
$\varphi_U$-equivariant $C^k$-map $\tilde f_U: \widetilde
U\to\widetilde U^\prime$ such that $\pi_{U^\prime}\circ\tilde
f_U=f\circ\pi_U$. ( Later we often say $(\tilde
f_U,\varphi_U):(\widetilde U, \Gamma_U,\pi_U)\to (\widetilde
U^\prime, \Gamma_{U^\prime},\pi_{U^\prime})$ to be a {\bf local
representation} of $f$ near $x$. Note that for any $\tilde
x\in\widetilde U$ the homomorphism $\varphi_U$ maps the isotropy
groups $\Gamma_U(\tilde x)$ to $\Gamma_{U'}(\tilde f_U(\tilde x))$.)
\item[(ii)] If $(\tilde f_V,\varphi_V):(\widetilde V, \Gamma_V,\pi_V)\to
(\widetilde V^\prime, \Gamma_{V^\prime},\pi_{V^\prime})$ be another
local representation of $f$ as in (i), and $U\cap V\ne\emptyset$,
then for any $x\in U\cap V$, $\tilde x_u\in(\pi_U)^{-1}(x)$ and
$\tilde x_v\in (\pi_V)^{-1}(x)$ there exist equivalent induced
charts at $x$ and $f(x)$ respectively,
\begin{eqnarray*}
(\widetilde O(\tilde x_u), \Gamma_U(\tilde x_u), \pi^x_U)\quad{\rm
and}\quad (\widetilde O(\tilde x_v), \Gamma_V(\tilde x_v),
\pi^x_V),\hspace{15mm}\\
(\widetilde O(\tilde f_U(\tilde x_u)), \Gamma_{U'}(\tilde f_U(\tilde
x_u)), \pi^{f(x)}_{U'})\;{\rm and}\; (\widetilde O(\tilde f_V(\tilde
x_v)), \Gamma_{V'}(\tilde f_V(\tilde x_v)), \pi^{f(x)}_{V'}),
\end{eqnarray*}
and the corresponding equivalences $(\lambda_{UV}, {\cal A}_{UV})$
(mapping $\tilde x_u$ to $\tilde x_v$) and $(\lambda_{U'V'}, {\cal
A}_{U'V'})$ (mapping $\tilde f_U(\tilde x_u)$ to $\tilde f_V(\tilde
x_v)$) such that
$$
 \left.
\begin{array}{ll}
 \tilde f_U\bigl(\widetilde O(\tilde
x_u)\bigr)\subset \widetilde O(\tilde f_U(\tilde x_u)),\\
 \tilde
f_V\bigl(\widetilde O(\tilde x_v)\bigr)\subset \widetilde
O(\tilde f_V(\tilde x_v))\quad{\rm and}\\
\tilde f_V\circ\lambda_{UV}=\lambda_{U'V'}\circ\tilde f_U.
\end{array}\right.
$$
\end{description}}
\end{definition}

A smooth map $f$ from orbifolds $X$ to $X'$ is called an {\bf
orbifold embedding} if it is a homeomorphism to the image
$f(X)\subset X'$ and each local representation $\tilde f_U:
\widetilde U\to\widetilde U^\prime$ as above is an embedding.

Now we begin to introduce a notion of Banach orbibundles. Without
loss of generality, we always assume that a finite group $G$ acts on
a smooth Banach vector bundle $E$ by bundle automorphisms of $E$.

\begin{definition}\label{def:2.7}{\rm Let $X$ be a
Banach orbifold. A {\bf  Banach orbibundle over $X$} consists of a
 Banach orbifold $E$ and a  smooth surjective
mapping $p: E\to X$ with the additional requirements:
\begin{description}
\item[(i)] For each  Banach orbifold chart $(\widetilde U,
\Gamma_U,\pi_U)$ of $X$ there exists a  Banach orbifold chart
$(\widetilde E_U, \widetilde\Gamma_U,\Pi_U)$ on $E$ with support
$E_U=p^{-1}(U)$, such that there exist a  smooth Banach bundle
structure $\widetilde p_U:\widetilde E_U\to\widetilde U$ and
 an injective homeomorphism $\lambda_U$ from
$\Gamma_U$ into the group of bundle automorphisms of $\widetilde
E_U$ such that:
\begin{enumerate}
\item[$(a)$] For each $g\in\Gamma_U$, $\lambda_U(g)$ is a lift of
$g$ (in particular $\lambda_U(g)=\1$ if $g=\1$).

\item[$(b)$] $\widetilde\Gamma_U=\{\lambda_U(g)\,|\,
g\in\Gamma_U\}$.

\item[$(c)$] $\pi_U\circ\widetilde p_U=p\circ\Pi_U$ in $E_U$.
\end{enumerate}

\item[(ii)] If $(\widetilde U, \Gamma_U,\pi_U)$ and $(\widetilde
V, \Gamma_V,\pi_V)$ are two charts in ${\mathcal U}$ with $U\subset
V$, $\theta_{UV}=(\widetilde\theta_{UV},\gamma_{UV})$ is an
injection from $(\widetilde U, \Gamma_U,\pi_U)$ to $(\widetilde V,
\Gamma_V,\pi_V)$, then there is an injection
$\Theta_{UV}=(\widetilde\Theta_{UV}, \Gamma_{UV})$ from $(\widetilde
E_U, \widetilde\Gamma_U,\Pi_U)$ to $(\widetilde E_V,
\widetilde\Gamma_V,\Pi_V)$, such that
$\widetilde\theta_{UV}\circ\widetilde p_U=\widetilde
p_V\circ\widetilde\Theta_{UV}$ and that
$\Theta_{UV}=(\widetilde\Theta_{UV},\Gamma_{UV})$ is also a bundle
open embedding.
\end{description}
Such a  Banach orbifold chart $(\widetilde E_U,
\widetilde\Gamma_U,\Pi_U)$ on $E$ is called the {\bf Banach
orbibundle chart} corresponding with $(\widetilde U,
\Gamma_U,\pi_U)$.}
\end{definition}

\begin{remark}\label{rem:2.8}
{\rm  Later we often write $(\widetilde E_U,
\widetilde\Gamma_U,\Pi_U)$ as  $(\widetilde E_U, \Gamma_U,\Pi_U)$
and understand $g\in\Gamma_U$ to act on $\widetilde E_U$ via
$G:=\lambda_U(g)$, i.e. $g\cdot\tilde\xi=G(\tilde\xi)$ for
$\tilde\xi\in\widetilde E_U$. If $X$ is $1$-regular, it follows from
Definition~\ref{def:2.7}(i) that $E$ is $1$-regular as well as.}
\end{remark}

\begin{definition}\label{def:2.9}
{\rm Let $p:E\to X$ and $p':E'\to X'$ be two Banach orbibundles,
$f:X\to X'$ be a $C^k$-smooth map for $0<k\le\infty$ and ${\cal
A}_X$ and ${\cal A}_{X'}$ be respectively Banach orbifold atlases on
$X$ and $X'$ associated with $f$ as in Definition~\ref{def:2.6}. A
continous map $\bar f: E\to E'$ is called an $C^k$-{\bf orbibundle
map} covering $f$ if the following are satisfied:
\begin{description}
\item[(i)] For each local representation of $f$, $(\tilde
f_U,\varphi_U):(\widetilde U, \Gamma_U,\pi_U)\to (\widetilde
U^\prime, \Gamma_{U^\prime},\pi_{U^\prime})$ with $(\widetilde U,
\Gamma_U,\pi_U)\in{\cal A}_X$ and $(\widetilde U^\prime,
\Gamma_{U^\prime},\pi_{U^\prime})\in{\cal A}_{X'}$, there exist
corresponding B-orbibundle charts $(\widetilde E_U, \Gamma_U,
\Pi_U)$ on $E$, $(\widetilde E'_{U'}, \Gamma_{U'}, \Pi_{U'})$ on
$E'$ and a $\varphi_U$-equivariant bundle map $\tilde{\bar f}_U:
\widetilde E_U\to\widetilde E'_{U'}$ covering $\tilde f_U$, i.e.,
$\tilde{\bar f}_U$ being a bundle map satisfying $\tilde
p_{U'}\circ\tilde{\bar f}_U=\tilde f_U\circ\tilde p_U$ and
$\tilde{\bar f}_U(g\cdot v)=\lambda_U(g)\cdot\tilde{\bar f}_U(v)$
for any $g\in\Gamma_U$ and $v\in\widetilde E_U$,  such that
$\Pi_{U'}\circ\tilde{\bar f}_U=f|_U\circ\Pi_U$. The triple
$(\tilde{\bar f}_U, \tilde f_U,\varphi_U):(\widetilde E_U,
\widetilde U, \Gamma_U)\to (\widetilde E'_{U'},\widetilde U^\prime,
\Gamma_{U^\prime})$ is called a {\bf local representation} of $(\bar
f,f)$ or $\bar f$. For the sake of clearness we later also call the
pair $(\bar f, f)$ a bundle map.
\item[(ii)] Let $(\tilde{\bar f}_V, \tilde f_V,\varphi_V):(\widetilde E_V,
\widetilde V, \Gamma_V)\to (\widetilde E'_{V'},\widetilde V^\prime,
\Gamma_{V^\prime})$ be another local representation of $(\bar f,f)$
or $\bar f$ such that $U\cap V\ne\emptyset$. For any $x\in U\cap V$,
$\tilde x_u\in(\pi_U)^{-1}(x)$ and $\tilde x_v\in (\pi_V)^{-1}(x)$
let
\begin{eqnarray*}
&&(\lambda_{UV}, {\cal A}_{UV}):(\widetilde O(\tilde x_u),
\Gamma_U(\tilde x_u), \pi^x_U)\to (\widetilde O(\tilde x_v),
\Gamma_V(\tilde x_v),
\pi^x_V)\;{\rm and}\\
&&(\lambda_{U'V'}, {\cal A}_{U'V'}): (\widetilde O(\tilde
f_U(\tilde x_u)), \Gamma_{U'}(\tilde f_U(\tilde x_u)),
\pi^{f(x)}_{U'})\to\\
&&\hspace{50mm} (\widetilde O(\tilde f_V(\tilde x_v)),
\Gamma_{V'}(\tilde f_V(\tilde x_v)), \pi^{f(x)}_{V'})\bigr )
\end{eqnarray*}
be as in Definition~\ref{def:2.6}. By shrinking $O(x)$ if
 necessary one can find ${\cal A}_{UV}$ (resp. ${\cal A}_{U'V'}$)-
equivariant bundle isomorphisms
\begin{eqnarray*}
&& \Lambda_{UV}:\widetilde E_U|_{\widetilde O(\tilde x_u)}\to
\widetilde E_V|_{\widetilde
O(\tilde x_v)}\;{\rm covering}\;\lambda_{UV}\\
&&{\rm (resp.}\;\Lambda_{U'V'}:\widetilde E_{U'}|_{\widetilde
O(\tilde f_U(x_u))}\to
 \widetilde E_{V'}|_{\widetilde O(\tilde
f_V(x_v))}\;{\rm covering}\;\lambda_{U'V'})
\end{eqnarray*}
 such that
\begin{equation}\label{e:2.2}
\tilde{\bar f}_V\circ\Lambda_{UV}=\Lambda_{U'V'}\circ(\tilde{\bar
f}_U|_{\widetilde O(\tilde x_u)}).
 \end{equation}
 \end{description}}
\end{definition}

We say an orbibundle map ${\bar f}: E\to E^\prime$ to be an {\bf
orbibundle embedding} if it is also an orbifold embedding. (In
particular, this implies that each local representation bundle map
$\tilde{\bar f}_U: \widetilde E_U\to\widetilde E'_{U'}$ as above
is injective.) An orbibundle map ${\bar f}: E\to E^\prime$ is
called an {\bf orbibundle open embedding} (from $E\to X$ to
$E^\prime\to X^\prime$) if $\bar f:E\to E^\prime$ (and thus
$f:X\to X^\prime$) is an open embedding. ({\it Note}: It implies
that for each local representation $(\tilde{\bar f}_U, \tilde
f_U,\varphi_U):(\widetilde E_U, \widetilde U, \Gamma_U)\to
(\widetilde E'_{U'}\widetilde U^\prime, \Gamma_{U^\prime})$ of
$\bar f$ as above both $\tilde f_U$ and $\tilde{\bar f}_U$ are
smooth open embedding and the restriction of $\tilde{\bar f}_U$ to
each fibre $(\widetilde E_U)_{\tilde z}$ of $\widetilde E_U$ gives
a Banach space isomorphism to the fibre $(\widetilde
E'_{U'})_{\tilde f_U(\tilde z)}$ of $\widetilde E'_{U'}$.) In
particular, if $\bar f$ is also a diffeomorphism we call it an
(orbibundle) {\bf isomorphism} from $E$ to $E'$.

\begin{definition}\label{def:2.10}
{\rm For a Banach orbibundle $p: E\to X$ as above, a continuous
map $S:X\to E$ is called a $C^k$-{\bf smooth section} if (i)
$p\circ S=id_X$ and (ii) for each $x\in X$ there exist a Banach
orbifold chart $(\widetilde U, \Gamma_U,\pi_U)$ of $X$ near $x$, a
corresponding bundle chart $(\widetilde E_U, \Gamma_U,\Pi_U)$ of
$E$ over $\widetilde U$, and a $C^k$-smooth $\Gamma_U$-equivariant
Banach bundle section $\widetilde S_U:\widetilde U\to\widetilde
E_U$ (i.e. $\widetilde S_U(g\cdot \tilde y)=g\cdot\widetilde
S_U(\tilde y)$ for $\tilde y\in\widetilde U$ and $g\in\Gamma_U$)
such that:

$\bullet$ It is a lift of $S|_U$ on $\widetilde U$, i.e.
$\Pi_U\circ\widetilde S_U=S\circ\pi_U$.

$\bullet$ All $\widetilde S_U$ are compatible in the following
sense: For any two Banach orbibundle charts $(\widetilde E_U,
\Gamma_U,\Pi_U)$ and $(\widetilde E_V, \Gamma_V,\Pi_V)$, and
corresponding Banach orbifold charts $(\widetilde U,
\Gamma_U,\pi_U)$ and $(\widetilde V, \Gamma_V,\pi_V)$, it holds
that for any $x\in U\cap V$ there exist an equivalence
$$
(\lambda_{UV}, {\cal A}_{UV}):(\widetilde O(\tilde x_u),
\Gamma_U(\tilde x_u), \pi^x_U)\to (\widetilde O(\tilde x_v),
\Gamma_V(\tilde x_v), \pi^x_V)
$$
as in Definition~\ref{def:2.9}, and a bundle isomorphism
 $\Lambda_{UV}:\widetilde E_U|_{\widetilde O(\tilde x_u)}\to
\widetilde E_V|_{\widetilde O(\tilde x_v)}$ covering
$\lambda_{UV}$ such that
\begin{equation}\label{e:2.3}
\widetilde S_V\circ\lambda_{UV}=\Lambda_{UV}\circ\widetilde
S_U\quad{\rm on}\quad \widetilde O(\tilde x_u).
\end{equation}
We call $\widetilde S_U$ a {\bf local lift} of $S$ in charts
 $(\widetilde E_U, \Gamma_U,\Pi_U)$ and $(\widetilde U,
\Gamma_U,\pi_U)$. Note that all these $\widetilde S_U$ are
uniquely determined by $S$.}
\end{definition}

It should also be mentioned that an orbifold map $S:X\to E$
satisfying  $p\circ S=id_X$ is not necessarily to give rise to a
section of $E$.

\begin{definition}\label{def:2.11}
{\rm We say $S$ to be {\bf Fredholm} if each local lift
$\widetilde S_U$  of $S$ in charts
 $(\widetilde E_U, \Gamma_U,\Pi_U)$ and $(\widetilde U,
\Gamma_U,\pi_U)$ above is so. A triple $(X, E,S)$ consisting of a
Banach orbibundle $p:E\to X$ and a Fredholm section $S:X\to E$ is
called a {\bf Banach Fredholm orbibundle}. The index of $S$ is
called the {\bf index} of $(X, E, S)$.}
\end{definition}

For a Banach Fredholm orbibundle $(X, E, S)$ and the local lift
$\widetilde S_U$ of $S$ in the above charts $(\widetilde E_U,
\widetilde\Gamma_U,\Pi_U)$ and $(\widetilde U, \Gamma_U,\pi_U)$ one
has a determinant (real line) bundle ${\rm det}(D\widetilde S_U)$
over $Z(\widetilde S_U)$. If this determinant line bundle ${\rm
det}(D\widetilde S_U)$ is trivial and is given a continuous nowhere
zero section, i.e. $\widetilde S_U$ is oriented we say $S$ to be
{\bf oriented} on $U$ or simply {\bf locally oriented}. We call $S$
{\bf oriented} if the representative $\widetilde S_U$ in each chart
$(\widetilde E_U, \Gamma_U,\Pi_U)$ of $E$ is oriented and the bundle
open embedding ${\rm det}(\Theta_{UV})$ from $({\rm det}(D\widetilde
S_U), \widetilde\Gamma_U,\Pi_U)$ to $({\rm det}(D\widetilde S_V),
\widetilde\Gamma_V,\Pi_V)$ that is naturally induced by the bundle
open embedding $\Theta_{UV}$ from $(\widetilde E_U,
\widetilde\Gamma_U,\Pi_U)$ to $(\widetilde E_V,
\widetilde\Gamma_V,\Pi_V)$  is orientation preserving in the
following sense: If $\tilde \tau_U$ (resp. $\tilde\tau_V$ ) is the
given continuous nowhere zero section of ${\rm det}(D\widetilde
S_U)$ (resp. ${\rm det}(D\widetilde S_V)$) then for each
$x\in\widetilde V$, ${\rm det}(\Theta_{UV})(\tilde\tau_V(x))$ is a
positive multiple  of $\tilde\tau_U\bigl(\theta_{UV}(x)\bigr)$. In
the case $S$ is oriented we say $(X, E, S)$ to be {\bf oriented}.

Two oriented Banach Fredholm orbibundles $(X, E^{(i)}, S^{(i)})$
 of index $r$ and with compact zero locus,
 $i=0,1$,  are called {\bf oriented homotopic}
 if there exists an oriented Banach Fredholm orbibundle $(X\times [0,1], E, S)$
 of index $r+1$ and with compact zero locus such that
 $E^{(i)}=E|_{\{i\}\times X}$ and $S^{(i)}=S|_{\{i\}\times X}$ for
 $i=0, 1$, and that the orientation of $(X\times [0,1], E, S)$ induces those of
$(X, E^{(i)}, S^{(i)})$, $i=0,1$.

Let $S:X\to E$ and $S': X'\to E'$ be two smooth sections of Banach
orbibundles. $S$ is called a {\bf pullback} of $S'$ via an
orbibundle map $\bar f: E\to E'$ covering $f:X\to X'$ if $S'\circ
f=\bar f\circ S$. It means: (i) For any local representation of
$\bar f$ as above, $(\tilde{\bar f}_U, \tilde
f_U,\varphi_U):(\widetilde E_U, \widetilde U, \Gamma_U)\to
(\widetilde E'_{U'},\widetilde U^\prime, \Gamma_{U^\prime})$, and
the local lifting of $S$, $\widetilde S_U:\widetilde
U\to\widetilde E_{U}$, and one of $S'$, $\widetilde
S'_{U'}:\widetilde U'\to\widetilde E'_{U'}$ it holds that
\begin{equation}\label{e:2.4}
\tilde{\bar f}_U\circ\widetilde S_U= \widetilde S'_{U'}\circ\tilde
f_U.
\end{equation}
(ii) If $(\tilde{\bar f}_V, \tilde f_V,\varphi_V):(\widetilde E_V,
\widetilde V, \Gamma_V)\to (\widetilde E'_{V'},\widetilde V^\prime,
\Gamma_{V^\prime})$ is another local representation of $\bar f$ with
$U\cap V\ne\emptyset$ then for any $x\in U\cap V$, in the notations
of Definition~\ref{def:2.9} it also holds that
\begin{equation}\label{e:2.5}
\Lambda_{U'V'}\circ\tilde{\bar f}_U\circ\widetilde S_U(\tilde y)=
\widetilde S'_{V'}\circ\tilde f_V\circ\lambda_{UV}(\tilde
y)\quad\forall\tilde y\in\widetilde O(\tilde x_u).
\end{equation}
A natural question is: Under what the conditions  can the section
$S'$  be determined by $S$ and $\bar f$?

\begin{claim}\label{cl:2.12}
The answer is affirmative if each local lifting of $\bar f$
restricts to a Banach isomorphism on each fiber.
\end{claim}

\noindent{\bf Proof.}  Let $(\tilde{\bar f}_U, \tilde
f_U,\varphi_U):(\widetilde E_U, \widetilde U, \Gamma_U)\to
(\widetilde E'_{U'},\widetilde U^\prime, \Gamma_{U^\prime})$ be a
local representation of $\bar f$ as above. By the assumptions,
for any $\tilde y\in\widetilde U$,
$$
(\tilde{\bar f}_U)_{\tilde y}:(\widetilde E_U)_{\tilde
y}\to(\widetilde E_{U'})_{\tilde f_U(\tilde y)}
$$
is a Banach space isomorphism. So we can define a smooth section
\begin{equation}\label{e:2.6}
\widetilde S_U: \widetilde U\to \widetilde E_U,\;\tilde y\mapsto
(\tilde{\bar f}_U)^{-1}_{\tilde y}\bigl(\widetilde
S'_{U'}\circ\tilde f_U(\tilde y)\bigr).
\end{equation}
If $(\tilde{\bar f}_V, \tilde f_V,\varphi_V):(\widetilde E_V,
\widetilde V, \Gamma_V)\to (\widetilde E'_{V'},\widetilde
V^\prime, \Gamma_{V^\prime})$ is another local representation of
$\bar f$, we can also define another smooth section
\begin{equation}\label{e:2.7}
\widetilde S_V: \widetilde V\to \widetilde E_V,\;\tilde z\mapsto
(\tilde{\bar f}_V)^{-1}_{\tilde z}\bigl(\widetilde
S'_{V'}\circ\tilde f_V(\tilde z)\bigr).
\end{equation}
If $U\cap V\ne\emptyset$ we need to prove that (\ref{e:2.3}) holds.
Indeed, following the notations in Definition~\ref{def:2.9}, for any
$\tilde y\in O(\tilde x_u)$ let $\tilde z=\lambda_{UV}(\tilde y)$.
Then by (\ref{e:2.6}) and (\ref{e:2.7}),
\begin{eqnarray*}
\widetilde S_V\circ\lambda_{UV}(\tilde y)\!\!\!\!\!&&=\widetilde
S_V(\tilde
z)\\
&&=(\tilde{\bar f}_V)^{-1}_{\tilde z}\bigl(\widetilde
S'_{V'}\circ\tilde f_V(\tilde z)\bigr)\\
&&=(\tilde{\bar f}_V)^{-1}_{\tilde z}\bigl(\widetilde
S'_{V'}\circ\tilde f_V\circ\lambda_{UV}(\tilde y)\bigr)\\
&&=(\tilde{\bar f}_V)^{-1}_{\tilde
z}\bigl(\Lambda_{U'V'}\circ\tilde{\bar f}_U\circ\widetilde
S_U(\tilde y)\bigr)\\
&&=\Lambda_{UV}\circ\widetilde S_U(\tilde y),
\end{eqnarray*}
where the fourth and fifth equalities come from (\ref{e:2.5}) and
(\ref{e:2.2}) respectively. \hfill$\Box$\vspace{2mm}

\begin{definition}\label{def:2.13}
{\rm Let $p:E\to X$ be a Banach orbibundle and $f:Y\to X$ be a
smooth map between Banach orbifolds. A {\bf pull-back Banach
orbibundle} of $E$ over $Y$ via $f$ is a Banach orbibundle $f^\ast
E\to Y$ together with a smooth orbibundle map $\bar f: f^\ast E\to
E$ covering $f:Y\to X$ and such that each local lifting of $\bar f$
restricts to  a Banach space isomorphism on each fiber.
Claim~\ref{cl:2.12} shows that for each smooth section $S:X\to E$
one can get a smooth section of $f^\ast E\to Y$ by pullback via
$\bar f$, denoted by $f^\ast S$.}
\end{definition}

 Let $(X, E, S)$ and $(X', E', S')$ be
Banach Fredholm orbibundles, and $f:X\to X'$ be a Fredholm map. An
orbibundle map $\bar f: E\to E'$ covering $f$ is called a {\bf
Fredholm orbibundle map} from  $(X, E, S)$ to $(X', E', S')$ if
$S'\circ f=\bar f\circ S$, i.e., $S$ is the pullback of $S'$ via
$\bar f$.

 Let $E\to X$ be a Banach
orbibundle and $(X', E', S')$ be a Banach Fredholm orbibundle. Let
$\bar f:E\to E'$ be an orbibundle map covering a Fredholm map
$f:X\to X'$ and $S:X\to E$ be the pullback of $S'$ via $\bar f$.
From (\ref{e:2.4}) it is not hard to see that $S$ is not necessarily
Fredholm. Clearly, we have:

\begin{claim}\label{cl:2.14}
If $\bar f$ satisfies the assumption in Claim~\ref{cl:2.12} then $S$
is also Fredholm and ${\rm Index}(S)={\rm Index}(f)+ {\rm
Index}(S')$.
\end{claim}

Consequently,  for a Fredholm Banach orbibundle $(X, E, S)$ of index
$r$ and a Fredholm map of index $d$ from a Banach manifold $Y$ to
$X$, if $\bar f: f^\ast E\to Y$ is a pullback bundle of $E$ on $Y$
via $f$ then one has a pullback Fredholm section $f^\ast S:Y\to
f^\ast E$ of index $r+d$. The Fredholm bundle $(Y, f^\ast E, f^\ast
S)$ is called the {\bf pullback Fredholm Banach orbibundle} of $(X,
E, S)$ via $f$.

  When the word ``Banach'' above is
replaced by the word ``Hilbert'' we get the definitions of the
{\bf  Hilbert orbifold}, {\bf Hilbert orbibundle} and {\bf Hilbert
Fredholm orbibundle}.

Similarly we have also the notions of complex Banach orbifolds,
holomorphic Banach orbibundles and holomorphic sections.

\subsection{An overall strategy and an example}\label{sec:2.2}

The method developed by Liu-Tian in [LiuT1]-[LiuT3] is  beautiful
and complicated. In this subsection we  explain the overall strategy
for their constructions and then give an example to show that this
method can be used to calculate the orbifold Euler characteristic of
an orbifold. Without special statements we always {\bf assume that
all Banach orbifolds are  $1$-regular} in this section.

 Consider a Banach orbifold
$(X, {\mathcal A})$ and assume that there exist two charts in
${\mathcal A}$, $(\widetilde W_i, \Gamma_i,\pi_i)$, $i=1,2$, such
that $W_1\cup W_2=X$. By them we can understand ``orbifoldness'' of
$X$ on $W_1$ and $W_2$ respectively. This means that on the
intersection $W_{12}:=W_1\cap W_2$ one has two kinds of
understanding ways. By the definition of orbifolds (cf.
Lemma~\ref{lem:2.4}) we can know that one can get  the same
``orbifold information'' near any point $x\in W_{12}$ in these two
kinds of ways (local change of coordinates). However, we have no a
``global'' understanding on $W_{12}$.

Firstly, assume $X$ is a manifold and $(\widetilde
W_i,\Gamma_i,\pi_i)$ become the usual coordinate charts $(\widetilde
W_i,\{\1\},\pi_i)\equiv(\widetilde W_i, \pi_i)$, $i=1,2$. In this
case there exists a diffeomorphism $\varphi:\pi_1^{-1}(W_1\cap
W_2)\to\pi_2^{-1}(W_1\cap W_2)$ such that $\pi_2\circ\varphi=\pi_1$.
That is, $\varphi$ is the change of coordinates (or transition
function). It means that $\pi_1$ and $\pi_2$ give the same smooth
topology information on $W_{12}$. Consider the graph of $\varphi$ in
$\widetilde W_1\times\widetilde W_2$, ${\rm
Graph}(\varphi)=\{(\tilde x, \varphi(\tilde x))\,|\,\tilde
x\in\pi_1^{-1}(W_1\cap W_2)\subset\widetilde W_1\}$, which is a
submanifold in $\widetilde W_1\times\widetilde W_2$. It is exactly
the fiber product $\widetilde W_{12}:=\{(\tilde x_1, \tilde
x_2)\in\widetilde W_1\times\widetilde W_2\,|\,\pi_1(\tilde
x_1)=\pi_2(\tilde x_2)\}$. Note that the projection
$\pi_{12}:\widetilde W_{12}\to W_{12},\;(\tilde x_1, \tilde
x_2)=\pi_1(\tilde x_1)=\pi_2(\tilde x_2)$ is well-defined. Denote by
$\pi^{12}_i:\widetilde W_{12}\to\widetilde W_i, (\tilde x_1, \tilde
x_2)\mapsto\tilde x_i$, $i=1,2$. Then the existence of the
transition function (diffeomorphism) $\varphi$ above is equivalent
to the following claim:
\begin{center}
{\bf
 There exist a smooth manifold $\widetilde W_{12}$,  a continuous
 surjective
 $\pi_{12}: \widetilde W_{12}\to W_{12}$,  and smooth open embeddings
 $\pi^{12}_i: \widetilde W_{12}\to\widetilde W_i$  such that
$\pi_i\circ\pi^{12}_i=\pi_{12}, i=1,2.$}\end{center}

Now we consider the case that $X$ is an orbifold. In the present
case there is no a diffeomorphism such as $\varphi$ above. We cannot
expect that the claim above is true in general. However, the fiber
product $\widetilde W_{12}$ and the projections $\pi_{12}:\widetilde
W_{12}\to W_{12}$, $\pi^{12}_i:\widetilde W_{12}\to\widetilde W_i$,
$i=1,2$, are still well-defined, and also satisfy
$\pi_i\circ\pi^{12}_i=\pi_{12}$, $i=1,2$. Moreover, $\pi_{12}$ is
invariant under the obvious action
$\Gamma_{12}=\Gamma_1\times\Gamma_2$ on $\widetilde W_{12}$ and also
induces a homeomorphism from $\widetilde W_{12}/\Gamma_{12}$ to
$W_{12}$; and for $i=1,2$ the projections $\pi^{12}_i:\widetilde
W_{12}\to\widetilde W_i$, are equivariant with respect to group
homomorphisms
$$
\lambda^{12}_i:\Gamma_{12}\to\Gamma_i,\;(g_1, g_2)\mapsto g_i.
$$
But it is regrettable that $\widetilde W_{12}$ {\bf is not a
manifold in general}! What is $\widetilde W_{12}$? For any $x\in
W_{12}$ and $\tilde x_i\in\pi_i^{-1}(x)$, $i=1,2$, recall that we
have induced charts of $(\widetilde W_i,\Gamma_i,\pi_i)$ at $\tilde
x_i$, $(\widetilde O(\tilde x_i), \Gamma_i(\tilde x_i),
\pi_i|_{\widetilde O(\tilde x_i)})$, and an equivalence
$$
(\lambda_{12}, {\mathcal A}_{12}): (\widetilde O(\tilde x_1),
\Gamma_1(\tilde x_1), \pi_1|_{\widetilde O(\tilde x_1)})\to
(\widetilde O(\tilde x_2), \Gamma_2(\tilde x_2), \pi_2|_{\widetilde
O(\tilde x_2)})
$$
consisting of a group isomorphism ${\mathcal A}_{12}:\Gamma_1(\tilde
x_1)\to \Gamma_2(\tilde x_2)$ and a ${\mathcal A}_{12}$-equivariant
diffeomorphism $\lambda_{12}:\widetilde O(\tilde x_1)\to \widetilde
O(\tilde x_2)$. It follows that the graph of $\lambda_{12}$, ${\rm
Graph}(\lambda_{12})$ is a submanifold of the fiber product
$\widetilde W_{12}$ which is diffeomorphic to $\widetilde O(\tilde
x_1)$ and $\widetilde O(\tilde x_2)$. So the graph map
 ${\rm Gr}(\lambda_{12}):\widetilde O(\tilde x_1)\to \widetilde
 W_1\times\widetilde W_2,\, \tilde y\mapsto (\tilde y, \lambda_{12}(\tilde
 y))$ is a smooth open embedding with image contained $\widetilde W_{12}$.
 Moreover, for any $g\in\Gamma_2(\tilde x_2)$  the composition $({\1}\times
g)\circ{\rm Gr}(\lambda_{12})$ is also such an open embedding. Hence
we get that
$$
\pi_{12}^{-1}(O(x))=\cup_{g\in\Gamma_1(\tilde x_1)}({\1}\times
g)\circ{\rm Gr}(\lambda_{12})(\widetilde O(\tilde x_1)),
$$
that is,  an union of $\sharp\Gamma_1(\tilde x_1)$ copies of
$\widetilde O(\tilde x_1)$. This shows that near $(\tilde x_1,\tilde
x_2)$ the fiber product $\widetilde W_{12}$ is a Banach variety (the
union of finitely many Banach manifolds), rather than a Banach
manifold in general! Such a good local structure suggests that there
is a possible desingularization of $\widetilde W_{12}$. Call a
nonempty connected relative open subset in $({\1}\times g)\circ{\rm
Gr}(\lambda_{12})(\widetilde O(\tilde x_1))$ as a {\bf local
component} of $\widetilde W_{12}$ near any point of this open
subset. Two local components containing a point $\tilde
y_{12}\in\widetilde W_{12}$ are said to belong to different kinds if
their intersection is not a local component of $\widetilde W_{12}$.
In the disjoint union of all local components of $\widetilde
W_{12}$, a true Banach manifold, two points are called {\bf
equivalent} if both can be contained a local component of
$\widetilde W_{12}$ and are also same as points of $\widetilde
W_{12}$. It may be proved that this indeed gives rise to a regular
relation in the disjoint union of all local components of
$\widetilde W_{12}$. Let $\widehat W_{12}$ denote the manifold of
all equivalence classes. (Note that it is not necessarily
connected.) The compositions of the natural quotient map from
$\widehat W_{12}$ to $\widetilde W_{12}$ and $\pi_{12}$,
$\pi^{12}_i$  give the maps $\hat\pi_{12}:\widehat W_{12}\to W_{12}$
and $\hat\pi^{12}_i:\widehat W_{12}\to\widehat W_i:=\widetilde W_i$
satisfying $\hat\pi_i\circ\hat\pi^{12}_i=\hat\pi_{12}$, where
$\hat\pi_i=\pi_i$, $i=1,2$. Moreover, the action of $\Gamma_{12}$ on
$\widetilde W_{12}$ naturally lifts to an action of $\Gamma_{12}$ on
$\widetilde W_{12}$ (which is also effective if the actions of
$\Gamma_1$ and $\Gamma_2$ are effective) such that (i) $\pi^{12}_i$
are also equivariant with respect to the group homeomorphisms
$\lambda^{12}_i$ above, $i=1,2$, and (ii) $\hat\pi_{12}$ is
invariant under this action and induces a homeomorphism from
$\widehat W_{12}/\Gamma_{12}$ to $W_{12}$. In addition the set of
regular points in $\widehat W_{12}$, $\widehat W_{12}^\circ$, is
exactly $(\hat\pi_{12})^{-1}(W_{12}^\circ)$, and the restrictions of
$\hat\pi^{12}_i$ to $\widehat W_{12}^\circ$ are
$$
\frac{\prod^2_{i=1}(\sharp\Gamma_i-\sharp{\rm
Ker}(\Gamma_i,\widetilde W_i)+ 1)}{\sharp\Gamma_i-\sharp{\rm
Ker}(\Gamma_i,\widetilde W_i)+ 1}-\hbox{fold (regular) coverings}
$$
to ${\rm Im}(\hat\pi^{12}_i)\cap\widehat W_i^\circ\subset\widehat
W_i$, $i=1,2$. In summary, we have proved:
\begin{center}
 {\bf  There exist an orbifold
chart $(\widehat W_{12},\Gamma_{12}, \hat\pi_{12})$ for the open
subset $W_{12}$ (which is not necessarily in the given orbifold
structure on\\ $X$), surjective group homomorphisms
$\lambda^{12}_i:\Gamma_{12}\to\Gamma_i$ and
$\lambda^{12}_i$-equivariant smooth maps $\hat\pi^{12}_i:\widehat
W_{12}\to\widehat W_i$ such that
$\hat\pi_i\circ\hat\pi^{12}_i=\hat\pi_{12}$, $i=1,2$.}
\end{center}
It will be a substitute of the above claim in the case of manifolds.
Such a system
$$
\bigr\{(\widehat W_i, \Gamma_i,\hat\pi_i), (\widehat W_{12},
\Gamma_{12}, \hat\pi_{12}), \hat\pi^{12}_i, \lambda^{12}_i\;\bigm|\;
i=1,2\bigl\}
$$
is called a {\bf desingularization} of $X$, and $\{(\widehat W_{12},
\Gamma_{12}, \hat\pi_{12}), \hat\pi^{12}_i, \lambda^{12}_i\;|\;
i=1,2\}$ will be used to replace the transition function between
$(\widehat W_1, \Gamma_1,\hat\pi_1)$ and $(\widehat W_2,
\Gamma_2,\hat\pi_2)$ which does not exist in general. (In order to
help the reader understanding the construction of $\widehat W_{12}$
we consider a simple example. Let $\widetilde W_{12}$ be the union
of $x$-axis $l_x$ and $y$-axis $l_y$ in $\R^2$, i.e., $\widetilde
W_{12}=l_x\cup l_y$. Then each local component of $\widetilde
W_{12}$ is a relatively open subset in $x$-axis or $y$-axis. Let
$0_x$ and $0_y$ denote the origin in $x$-axis and $y$-axis
respectively. As points in $\widetilde W_{12}$ both are same, i.e.
$0_x=0_y\in\widetilde W_{12}$. A local component $C_x$ containing
$0_x$ and one $C_y$ containing $0_y$ are of different kinds if and
only if one of both is a nonempty relatively open subset in $l_x$
and another is such an open subset in $l_y$. Hence our definition of
the equivalence above shows that $\widehat W_{12}$ is the disjoint
union of $l_x$ and $l_y$.)

Similarly, if $p:E\to X$ is an orbibundle, and $\tilde
p_i:(\widetilde E_i, \Gamma_i,\Pi_i)\to (\widetilde
W_i,\Gamma_i,\pi_i)$, $i=1, 2$, are corresponding orbibundle charts,
then we can get a system of Banach bundles
$$
\bigr\{(\widehat E_i, \Gamma_i,\hat\Pi_i), (\widehat E_{12},
\Gamma_{12}, \hat\Pi_{12}), \hat\Pi^{12}_i, \lambda^{12}_i\;\bigm|\;
i=1,2\bigl\},
$$
which is not only a desingularization of $E$ as a Banach orbifold
but also there exists a Banach bundle structure $\hat
p_{12}:\widehat E_{12}\to\widehat W_{12}$ such that
$$\hat\pi_{12}\circ\hat
p_{12}=p\circ\hat\Pi_{12},\quad \hat\pi_i\circ\hat
p_i=p\circ\hat\Pi_i\quad{\rm and}\quad\hat
p_i\circ\hat\Pi^{12}_i=\hat\pi^{12}_i\circ\hat p_i,\; i=1,2,
$$
where $\widehat E_i=\widetilde E_i$, $\widehat W_i=\widetilde W_i$
and $\hat p_i=\tilde p_i:\widehat E_i\to\widehat W_i$, $i=1,2$. The
Banach bundle system above is called a {\bf desingularization} of
the orbibundle $E$ (associated with the above  desingularization of
the orbifold $X$). Moreover, if $S_i:\widetilde W_i\to\widetilde
E_i$, $i=1,2$, are local representations of a section $S:X\to E$,
then we obtain a natural section $\hat S_{12}:\widehat
W_{12}\to\widehat E_{12}$ such that it together $\hat S_i=\tilde
S_i$, $i=1,2$, forms a collection $\hat S=\{\hat S_{12}, \hat S_1,
\hat S_2\}$ that is compatible in the sense that $\hat S_{12}$ is
$\Gamma_{12}$-equivariant and $\hat S_i=(\hat\Pi^{12}_i)^\ast\hat
S_i:=(\hat\Pi^{12}_i)^{-1}\circ\hat S_i\circ\hat\pi^{12}_i$,
$i=1,2$. This collection $\hat S$ is called a desingularization of
the section $S:X\to E$ associated with the above desingularization
of the Banach orbibundle. Any collection with these properties is
called a {\bf global section} of the bundle system.

Now assume that the zero locus $Z(S)$ is compact. Choose open
subsets $W^1_i\subset Cl(W^1_i)\subset U^1_i\subset Cl(U^1_i)\subset
W_i$, $i=1,2$ such that $Z(S)\subset W^1_1\cup W^1_2$. Set
$$
V_{12}=W_1\cap W_2\quad{\rm and}\quad V_i=W^1_i\setminus
Cl(U^1_1)\cap Cl(U^1_2),\;i=1,2.
$$
It is easily checked that $Z(S)\subset V_1\cup V_2\cup V_{12}$ and
$Cl(V_1)\cap V_2=\emptyset=Cl(V_2)\cap V_1$. Define
\begin{eqnarray*}
&&\widehat V_{12}=(\hat\pi_{12})^{-1}(V_{12})\quad{\rm and}\quad
\widehat V_{i}=(\hat\pi_{i})^{-1}(V_{i}),\;i=1,2,\\
&&\widehat F_{12}=(\hat\Pi_{12})^{-1}(E_{12}|_{V_{12}})\quad{\rm
and}\quad \widehat F_{i}=(\hat\Pi_{i})^{-1}(E_{i}|_{V_i}),\;i=1,2.
\end{eqnarray*}
We get a system of Banach bundles
$$
(\widehat{\mathcal F},\widehat{\mathcal V}):=\bigr\{(\widehat F_i,
\Gamma_i,\hat\Pi_i), (\widehat F_{12}, \Gamma_{12}, \hat\Pi_{12}),
\hat p_{12}, \hat p_i, \hat\Pi^{12}_i, \lambda^{12}_i\;\bigm|\;
i=1,2\bigl\}
$$
over the system of Banach manifolds
$$
\bigr\{(\widehat V_i, \Gamma_i,\hat\pi_i), (\widehat V_{12},
\Gamma_{12}, \hat\pi_{12}), \hat\pi^{12}_i, \lambda^{12}_i\;\bigm|\;
i=1,2\bigl\}.
$$
The notion of a global section of this system is defined in the same
way as above. One can prove that each smooth section
$\tilde\sigma:\widetilde W_i\to\widetilde E_i$ with support in
$\widetilde W^1_i:=\pi^{-1}_i(W^1_i)$ may yield a global section
 of the system above, denoted by $\hat\sigma=\{(\hat\sigma)_1, (\hat\sigma)_2, (\hat\sigma)_{12}\}$.

Furthermore, choose nonempty open subsets
$$
 U^0_i\subset Cl(U^0_i)\subset W^1_i,\;i=1,2
$$
such that $Z(S)\subset U^0_1\cup U^0_2$. {\bf We also assume that
$S$ is a Fredholm section of index $r$}, and that $\tilde
p_i:\widetilde W_i\to \widetilde E_i$ has a trivialization
${\mathcal T}_i:\widetilde W_i\times(\widetilde E_i)_{\tilde
x_i}\to\widetilde E_i$, $i=1,2$. Let $\widetilde S^T_i:\widetilde
W_i\to (\widetilde E_i)_{\tilde x_i}$ be the representation of
$\widetilde S_i$ under the trivialization ${\mathcal T}_i$, $i=1,2$.
Since $S$ is Fredholm we can take finitely many sections $\tilde
\sigma_i^{j}:\widetilde W_i\to\widetilde E_i$ with supports in
$\widetilde W^1_i$ such that their trivialization representations
$\tilde \sigma_i^{jT}:\widetilde W_i\to(\widetilde E_i)_{\tilde
x_i}$, $j=1,\cdots, m_i$ and $i=1,2$, satisfy
$$
d\widetilde S_i^{T}(\tilde x)(T_{\tilde x}\widetilde W_i)+ {\rm
span}(\bigl\{\tilde \sigma_i^{1T}(\tilde x),\cdots,
\tilde\sigma_i^{m_iT}(\tilde x)\bigr\})=(\widetilde E_i)_{\tilde
x_i}
$$
for any $\tilde x\in\widetilde U^0_i$ and $i=1,2$. Let
$\hat\sigma_i^{j}=\{(\hat\sigma_{i}^{j})_{12},
(\hat\sigma_{i}^{j})_1, (\hat\sigma_{i}^{j})_2\}$ be the global
sections of the bundle system $(\widehat{\mathcal
F},\widehat{\mathcal V})$, $j=1,\cdots, m_i$ and $i=1,2$, and ${\bf
t}$ denote points $(t_{11},\cdots,t_{1m_1}, t_{21},\cdots,
t_{2m_2})$ in $\R^{m_1+m_2}$. Then it may be proved that for generic
small ${\bf t}$ the section $\Psi^{\bf t}=\{(\Psi^{\bf t})_1,
(\Psi^{\bf t})_2, (\Psi^{\bf t})_{12}\}$ is transversal to the zero
section. That is, the following three sections are transversal to
the zero section,
 \begin{eqnarray*}
 &&\Psi^{\bf t}_{12}:\widehat V_{12}\to\widehat F_{12},\;\hat
 x\mapsto\widehat S_{12}(\hat x)+
 \sum^{m_1}_{j=1}t_{1j}(\hat\sigma_{1}^{j})_{12}(\hat x)+
 \sum^{m_2}_{j=1}t_{2j}(\hat\sigma_{2}^{j})_{12}(\hat x)\\
&&\Psi^{\bf t}_1:\widehat V_{1}\to\widehat F_{1},\;\hat
 x\mapsto\widehat S_{1}(\hat x)+
 \sum^{m_1}_{j=1}t_{1j}(\hat\sigma_{1}^{j})_{1}(\hat x)+
 \sum^{m_2}_{j=1}t_{2j}(\hat\sigma_{2}^{j})_{1}(\hat x)\\
&&\Psi^{\bf t}_2:\widehat V_{2}\to\widehat F_{2},\;\hat
 x\mapsto\widehat S_{2}(\hat x)+
 \sum^{m_1}_{j=1}t_{1j}(\hat\sigma_{1}^{j})_{2}(\hat x)+
 \sum^{m_2}_{j=1}t_{2j}(\hat\sigma_{2}^{j})_{2}(\hat x).
\end{eqnarray*}
Let $\widehat{\mathcal M}^{\bf t}_1=(\Psi^{\bf t}_1)^{-1}(0)$,
$\widehat{\mathcal M}^{\bf t}_2=(\Psi^{\bf t}_2)^{-1}(0)$,
$\widehat{\mathcal M}^{\bf t}_{12}=(\Psi^{\bf t}_{12})^{-1}(0)$.
They are all smooth manifolds of dimension $r$. Denote by
$$
\widehat{\mathcal M}^{{\bf t}\circ}_{12}:=\widehat{\mathcal M}^{\bf
t}_{12}\cap\widehat V_{12}^\circ\quad\hbox{and}\quad
\widehat{\mathcal M}^{{\bf t}\circ}_i:=\widehat{\mathcal M}^{\bf
t}_i\cap\widehat V_i^\circ,\;i=1,2.
$$
Here $\widehat V_{12}^\circ=\widehat V_{12}\cap\widehat
W_{12}^\circ$ and $\widehat V_i^\circ=\widehat V_i\cap\widehat
W_i^\circ$, $i=1,2$. Then the maps
$$
\hat\pi_i^{12}: (\hat\pi_i^{12})^{-1}\bigl(\widehat{\mathcal
M}^{{\bf t}\circ}_i\bigr)\to {\rm Im}(\hat\pi_i^{12})\subset
\widehat{\mathcal M}^{{\bf t}\circ}_i
$$
are $$\frac{\prod^2_{i=1}(\sharp\Gamma_i-\sharp{\rm
Ker}(\Gamma_i,\widetilde W_i)+ 1)}{\sharp\Gamma_i-\sharp{\rm
Ker}(\Gamma_i,\widetilde W_i)+ 1}-\hbox{fold (regular) coverings},\;
i=1,2.
$$
{\bf Assume that $X$ has at least codimension $2$ singularity near
$Z(S)$ and $Z(S)\subset{\rm Int}(X)$.} Then the formal sum
\begin{eqnarray*}
e(E,S)^{\bf
t}:\!\!\!\!\!\!&&\!\!\!=\sum_{i=1}^2\frac{1}{\sharp\Gamma_i-\sharp{\rm
Ker}(\Gamma_i,\widetilde W_i)+ 1)}\bigl\{\hat\pi_i:\widehat{\mathcal
M}^{{\bf t}\circ}_i\to X\bigr\}\\
&&+\frac{1}{\prod^2_{i=1}(\sharp\Gamma_i-\sharp{\rm
Ker}(\Gamma_i,\widetilde W_i)+ 1)}\bigl\{\hat\pi_{12}:
\widehat{\mathcal M}^{{\bf t}\circ}_{12}\to X\bigr\}
\end{eqnarray*}
represents a rational singular cycle in $X$ of dimension $r$, where
we only count once on the overlaps in the summation of the singular
maps above. The homology class of $e(E,S)^{\bf t}$,
$e(E,S)=[e(E,S)^{\bf t}]\in H_r(X,\Q)$ is called the {\bf virtual
Euler class} of the Banach Fredholm orbibundle $(E,X,S)$. We can
prove that this class is independent all related choices. If $r=0$
one may get a well-defined rational number $e(E,S)^{\bf t}$, called
{\bf virtual Euler characteristic} of $(E,X,S)$. In particular, when
$X$ is an orbifold of finite dimension and $E$ is its tangent bundle
it is simply called the {\bf orbifold Euler characteristic of $X$}
(i.e. the Euler characteristic of $X$ as an orbifold).

\begin{example}\label{ex:2.15}
Calculating the orbifold Euler characteristic of a football
orbifold: {\rm It is a weighted projective space. Let ${\bf
q}=(q_1, q_2)$ be a pair of positive integers with highest common
divisor $1$, and $q_1, q_2>1$. The weighted (or twisted)
projective space of ${\bf q}$ is defined by
$$
\CP^1({\bf q})=(\C^2\setminus\{0\})/\C^\ast,
$$
where $\C^\ast=\C\setminus\{0\}$ acts by $\alpha\cdot {\bf
z}=(\alpha^{q_1}z_1, \alpha^{q_2}z_2)$
 for ${\bf z}=(z_1, z_2)\in\C^2\setminus\{0\}$ and $\alpha\in\C^\ast$.
 The  $\C^\ast$-action above is free iff $q_1=q_2=1$.   $\CP^1({\bf q})$ has only isolated orbifold
singularities. Let $[{\bf z}]_{\bf q}$ denote the orbit of ${\bf
z}\in\C^2\setminus\{0\}$ under the above $\C^\ast$-action. It is a
point in $\CP^1({\bf q})$.  Clearly, $X=\CP^1({\bf q})$  has two
cone singularities, one of type $q_1$ (the north pole $p_N$) and
another of type $q_2$ (the south pole $p_S$). The orbifold structure
group $\Gamma_{q_1}$ at $p_N=[1,0]_{\bf q}$ (resp. $\Gamma_{q_2}$ at
$p_S=[0,1]_{\bf q}$) is isomorphic to $\Z/q_1\Z$ (resp. $\Z/q_2\Z$).
 Other points are all smooth
points. Thus $X$ is an effective orbifold and $2$-regular.
(Topologically, it is a $2$-sphere.)  It can be covered by two open
disks, $W_1$ containing $p_N$ and $W_2$ containing $p_S$, that
intersect in an open annulus $W_{12}$. Then $\pi_1:\widetilde W_1\to
W_1$ (resp. $\pi_2:\widetilde W_2\to W_2$) is given by in polar
coordinates $(r,\theta)$ by $(r,\theta)\mapsto (r, q_1\theta)$
(resp. $(s,\phi)$ by $(s,\phi)\mapsto (s, q_2\phi)$). $\Z/q_1\Z$
(resp. $\Z/q_2\Z$) acts on $\widetilde W_1$ (resp. $\widetilde W_2$)
by $(r,\theta)\mapsto (r, \theta+ i/q_1)$, $i=0,\cdots, q_1-1$
(resp. $(s,\phi)\mapsto (s, \phi+ j/q_2)$, $j=0,\cdots, q_2-1$).
Here $\theta,\phi\in S^1=\R/\Z$. Choose smaller discs $V_1, V_2$ in
$W_1, W_2$ respectively. Assume that the annulus $V_{12}$ is given
by $\{(t,\theta): t\in (1,4), \theta\in S^1=\R/\Z\}$ in the polar
coordinate in $W_1$, and that
$$
V_1\cap V_{12}=(1,2)\times S^1\quad{\rm and}\quad V_2\cap
V_{12}=(3,4)\times S^1.
$$
The tangent bundle $\widetilde E_1=T\widetilde W_1\to\widetilde W_1$
(resp. $\widetilde E_2=T\widetilde W_2\to\widetilde W_2$) can be
trivialised away from $p_N$ (resp. $p_S$) by the vector fields
$\tilde\partial_r,\tilde\partial_\theta$ (resp.
$\tilde\partial_s,\tilde\partial_\phi$) which are $\Gamma_{q_1}$
(resp. $\Gamma_{q_2}$) invariant. So
$\tilde\partial_r,\tilde\partial_\theta$ (resp.
$\tilde\partial_s,\tilde\partial_\phi$) descend to sections
$\partial_r, \partial_\theta$ of $E_1\to W_1\setminus\{p_N\}$ (resp.
$\partial_s, \partial_\phi$ of $E_2\to W_2\setminus\{p_S\}$). On the
overlap, $\partial_r=-\partial_\theta$ and
$\partial_s=-\partial_\phi$. Let $\tilde\sigma_1:\widetilde
W_1\to\widetilde E_1$ (resp. $\tilde\sigma_2:\widetilde
W_2\to\widetilde E_2$) be a section of the form
$\beta_1(r)\tilde\partial_r$ (resp. $\beta_2(s)\tilde\partial_s$)
where the function $\beta_1\le 0$ (resp. $\beta_2\ge 0$) has an
isolated zero at $r=0$ (resp. $s=0$). Then the section
$$
\hat\sigma=\{(\hat\sigma_1)_1+ (\hat\sigma_2)_1, (\hat\sigma_1)_2+
(\hat\sigma_2)_2, (\hat\sigma_1)_{12}+ (\hat\sigma_2)_{12}\}
$$
has only zeros $p_N$ and $p_S$ (we here have identified $p_N$ (resp.
$p_S$) with the unique point in $\hat\pi_1^{-1}(p_N)$ (resp.
$\hat\pi_2^{-1}(p_S)$)). Note that both $p_N$ and $p_S$ are regular.
However, by slightly perturbing $\hat\sigma$ near $p_N$ and $p_S$ we
can achieve transversality, and in this case there is a zero near
$p_N$ (resp. near $p_S$) that counts with weight $1/q_1$ (resp.
$1/q_2$). Hence the orbifold Euler characteristic of this orbifold
is $\frac{1}{q_1}+\frac{1}{q_2}$. The reader may refer to \cite{Mc1}
for more details and examples. }
\end{example}

\subsection{ The resolutions of Banach
orbibundles near a compact subset}\label{sec:2.3}

 It is a key step of Liu-Tian's
method in [LiuT1]-[LiuT3] to resolve a Banach orbibundle  near a
compact subset of its base space. We shall presently give a detailed
exposition of this construction; see also \cite{Lu3}. For the sake
of simplicity we {\bf always  assume that $X$ is $1$-regular and
effective} without special statements.  (Actually, for a Banach
Fredholm orbibundle $(X,E,S)$ with compact zero locus $Z(S)$ the
condition that $X$ is $1$-regular can be weaken to the requirement
that $X$ is $1$-regular near any $x\in Z(S)$ because we can deduce
that an open  neighborhood of $Z(S)$ in $X$ is a $1$-regular Banach
orbifold.) In Remark~\ref{rm:2.49} we shall point out how to change
the arguments and results when the effectiveness assumption of $X$
is moved.

\subsubsection{A general resolution}\label{sec:2.2.1}

Let $p:E\to X$ be a  Banach orbibundle and $K\subset X$ be a compact
subset. For each $x\in K$ let $(\widetilde W_x, \Gamma_x,\pi_x)$ be
a  Banach orbifold chart  around it. ({\bf This chart can be
required to be centered at $x$ if necessary}). Since $K$ is compact
we may choose finitely many points $x_{i}\in K$, $i=1,\cdots,n$ (we
assume $n>2$ because of the arguments in \S\ref{sec:2.2}) and
corresponding Banach orbifold charts $(\widetilde W_i,
\Gamma_i,\pi_i)$ around them  on $X$, $i=1,\cdots, n$ such that
their supports $\{W_i\}^n_{i=1}$ satisfy:
\begin{equation}\label{e:2.8}
K\subset\bigcup^n_{i=1}W_i\quad{\rm
and}\quad\bigcap^n_{i=1}W_i=\emptyset.
\end{equation}
 Let
${\mathcal N}$ be the set of all finite subsets $I=\{i_1,\cdots,
i_k\}$ of $\{1,\cdots,n\}$ such that the {\bf intersection}
$W_I:=\cap_{i\in I}W_i$ is nonempty. Then each $I\in{\mathcal N}$
has the length $|I|=\sharp(I)<n$.
  For $I=\{i_1,\cdots, i_k\}\in{\mathcal N}$ we may always assume
  $i_1<\cdots<i_k$. Denote by  the group
$\Gamma_I:=\prod^k_{l=1}\Gamma_{i_l}$ and the {\bf fiber product}
\begin{equation}\label{e:2.9}
\widetilde W_I=\Bigl\{(\tilde u_{i_1},\cdots, \tilde
u_{i_k})\in\prod^k_{l=1} \widetilde W_{i_l}\,\bigm|\, \pi_i(\tilde
u_i)=\pi_j(\tilde u_j)\;\forall i, j\in I\Bigr\}
\end{equation}
equipped with the induced topology from $\prod^k_{l=1} \widetilde
W_{i_l}$.
 Then the obvious {\bf projection}
\begin{equation}\label{e:2.10}
\pi_I: \widetilde W_I\to W_I,\;(\tilde u_{i_1},\cdots,\tilde
u_{i_k})\mapsto \pi_{i_1}(\tilde u_{i_1})=\cdots =\pi_{i_k}(\tilde
u_{i_k}),
 \end{equation}
 has {\bf covering group} $\Gamma_I$ whose action on $\widetilde W_I$ is given by
\begin{equation}\label{e:2.11}
\phi_I\cdot (\tilde u_i)_{i\in I}=(\phi_i\cdot \tilde u_i)_{i\in
I},\quad\forall\phi_I=(\phi_i)_{i\in I}\in \Gamma_I.
\end{equation}
(Hereafter we often write $(\tilde v_{i_1},\cdots, \tilde
v_{i_k})$ as $(\tilde v_i)_{i\in I}$.)  Namely, $\pi_I$ induces a
{\bf homeomorphism} from $\widetilde W_I/\Gamma_I$ onto $W_I$.
Note that $\pi_I$ is always proper.

\begin{claim}\label{cl:2.16}
The action of $\Gamma_I$ on $\widetilde W_I$ in (\ref{e:2.11}) is
effective.
\end{claim}

 Indeed, since each $(\widetilde W_i, \Gamma_i, \pi_i)$
is $1$-regular it follows from Claim~\ref{cl:2.5} that the
nonsingular set $\widetilde W_i^{\circ}$ is a dense open subset in
$\widetilde W_i$. Let
\begin{equation}\label{e:2.12}
\widetilde W_I^\circ:=\{\tilde u_I\in\widetilde
W_I\,|\,\Gamma_I(\tilde u_I)=\{\1\}\}\quad{\rm and}\quad
W_I^\circ:=\pi_I(\widetilde W_I^\circ),
\end{equation}
where $\Gamma_I(\tilde u_I)$ is the isotopy subgroup of $\Gamma_I$
 at $\tilde u_I\in \widetilde W_I$.
Then $\widetilde W_I^\circ=\pi_{i_1}^{-1}(
W_I^\circ)\times\cdots\times\pi_{i_k}^{-1}(W_I^\circ)$ is a dense
open subset in $\widetilde W_I$. Moreover, $W_I^\circ=\cap_{i\in
I}\pi_i(\widetilde W_i^\circ)\cap W_I$ is also a dense open subset
in $W_I$. Note that
$$
\widetilde
W_I=\pi_{i_1}^{-1}(W_I)\times\cdots\times\pi_{i_k}^{-1}(W_I)
$$
and that $g_s\tilde y\ne\tilde y$ for any
$y\in\pi_{i_s}^{-1}(W_I^\circ)$ and
$g_s\in\Gamma_{i_s}\setminus\{\1\}$, $s=1,\cdots, k$. Then for any
$\phi\in\Gamma_I\setminus\{\1\}$ and $(\tilde u_i)_{i\in
I}\in\widetilde W_I^\circ$ it holds that $\phi\cdot(\tilde
u_i)_{i\in I}\ne(\tilde u_i)_{i\in I}$. Claim~\ref{cl:2.16} is
proved.

\begin{remark}\label{rm:2.17}
{\rm As pointed out at the beginning of this section, the regular
set $X^\circ$ in an effective orbifold of finite dimension is always
an open, dense subset in $X$. So the proof of Claim~\ref{cl:2.16}
and thus all arguments till the end of Subsection~\ref{sec:2.5} are
still true even if $X$ is not $1$-regular. However, for an effective
orbifold of infinite dimension we need the $1$-regularity of $X$ to
complete the proof of Claim~\ref{cl:2.16}. Moreover, {\bf if we do
not assume that $X$ is effective and do not make the arguments in
the effective category, Claim~\ref{cl:2.16} is not needed.} See
Remark~\ref{rm:2.49}.}
\end{remark}

Actually, we can give a precise description for singularities of
$\widetilde W_I$ locally. For each $I\in{\cal N}$ we also denote by
\begin{equation}\label{e:2.13}
\widetilde W_I^{sing}:=\{\tilde u_I\in\widetilde
W_I\,|\,\Gamma_I(\tilde u_I)\ne\{\1\}\}\;\;{\rm and}\;\;
W_I^{sing}:=\pi_I(\widetilde W_I^{sing}).
\end{equation}
Clearly, they are relatively closed subsets in $\widetilde W_I$
and $W_I$ respectively.
 For any $u\in W_I^{sing}$ and $\tilde u_I=(\tilde u_{i_1},\cdots, \tilde
u_{i_k})\in\pi_I^{-1}(u)\subset\widetilde W_I$, as in the
arguments above (\ref{e:2.1}) there exists a
 $\Gamma_{i_1}(\tilde u_{i_1})$-equivariant diffeomorphism $F$ from a small neighborhood
 $N(0_{\tilde u_{i_1}})$
 of the origin in $T_{\tilde u_{i_1}}\widetilde W_{i_1}$ onto a neighborhood
 of $\tilde u_{i_1}$ in $\widetilde W_{i_1}$, ${\cal N}(\tilde u_{i_1})=F(N(0_{\tilde
 u_{i_1}}))$, such that
 $$
 {\cal N}(\tilde
 u_{i_1})^{sing}:=\pi_{i_1}^{-1}(W_I^{sing})\cap {\cal N}(\tilde
 u_{i_1})=\!\!\!\!\bigcup_{g\in\Gamma_{i_1}(\tilde u_{i_1})
 \setminus\{\1\}}F\bigl(N(0_{\tilde u_{i_1}})\cap{\nu}(dg(\tilde u_{i_1}))\bigr).
 $$
By Lemma~\ref{lem:2.4}, by shrinking $N(0_{\tilde u_{i_1}})$ if
necessary we can find group isomorphisms ${\mathcal
A}_{i_1i_s}:\Gamma_{i_1}(\tilde u_{i_1})\to\Gamma_{i_s}(\tilde
u_{i_s})$ and  diffeomorphisms $\lambda_{i_1i_s}: {\cal N}(\tilde
u_{i_1})\to \widetilde {\cal N}(\tilde u_{i_s})$ that map $\tilde
u_{i_1}$ to $\tilde u_{i_s}$, and  such that
$\lambda_{i_1i_s}\circ\phi={\mathcal
A}_{i_1i_s}(\phi)\circ\lambda_{i_1i_s}$ for any $\phi\in
\Gamma_{i_1}(\tilde u_{i_1})$ and $1<s\le k$. It follows that

$$
\widetilde W_I^{sing}\cap\bigl(\prod^k_{i=1}{\cal N}(\tilde
u_{i_s})\bigr)=\Bigl\{(\tilde x, \lambda_{i_1i_2}(\tilde
x),\cdots, \lambda_{i_1i_k}(\tilde x))\,\Bigm|\,\tilde x\in{\cal
N}(\tilde
 u_{i_1})^{sing}\Bigr\}.
$$

For any $J\subset I\in {\mathcal  N}$ there are also projections
\begin{equation}\label{e:2.14}
\left.\begin{array}{ll}
 \pi^I_J: \widetilde W_I\to \widetilde W_J,
\;\;(\tilde u_i)_{i\in I}\mapsto (\tilde u_j)_{j\in J},\\
 \lambda^I_J:\Gamma_I\to\Gamma_J,\quad\, (\phi_i)_{i\in
I}\mapsto (\phi_j)_{j\in J}.
\end{array}\right\}
\end{equation}
 Later we also write $\lambda^I_J$ as $\lambda^i_J$ if $I=\{i\}$.
Note that $\pi^I_J$ {\bf is not surjective in general.} To see
this point let us consider the case $I=\{1,2\}$ and $J=\{1\}$.
Then $\widetilde W_I=\{(\tilde u_1,\tilde u_2)\in\widetilde
W_1\times\widetilde W_2\,|\,\pi_1(\tilde u_1)=\pi_2(\tilde
u_2)\}$. Assume that $\pi^I_J:\widetilde W_I\to\widetilde W_J$ is
surjective. For each $\tilde u_1\in\widetilde W_J=\widetilde W_1$
there exists $\tilde u_2\in\widetilde W_2$ such that $(\tilde
u_1,\tilde u_2)\in\widetilde W_I$ and $\pi^I_J((\tilde u_1,\tilde
u_2))=\tilde u_1$. It follows that $u_1=\pi_1(\tilde
u_1)=\pi_2(\tilde u_2)=u_2\in W_2$. That is,  $W_1\subset W_2$.
Clearly, it must not hold in general.

For any $J\subset I\in{\cal N}$, by the definition one easily
checks that $\widetilde W_I^\circ=(\pi^I_J)^{-1}\bigl({\rm
Im}(\pi^I_J)\cap\widetilde W_J^\circ\bigr)$. These imply that for
every $u_I\in W_I^\circ=\pi_I(\widetilde W_I^\circ)\subset
W_J^\circ$  the inverse image $\pi_I^{-1}(u_I)$ (resp.
$\pi_J^{-1}(u_I)$) exactly contains $|\Gamma_I|$ (resp.
$|\Gamma_J|$) points. So for each $\tilde u_J\in\widetilde
W_J^\circ\cap{\rm Im}(\pi^I_J)$ the inverse image
$(\pi^I_J)^{-1}(\tilde u_J)$ has exactly $|\Gamma_I|/|\Gamma_J|$
points. If $(\widetilde E_i, \Gamma_i,\Pi_i) (=(\widetilde E_i,
\widetilde\Gamma_i,\Pi_i)$) are the  Banach orbibundle charts on
$E$ corresponding with $(\widetilde W_i, \Gamma_i,\pi_i)$,
$i=1,\cdots, n$. Repeating the same construction from $\widetilde
E_i$ one obtains the similar projections $\Pi_I:\widetilde E_I\to
E _I$ and $\Pi^I_J:\widetilde E_I\to\widetilde E_J$ for $J\subset
I$. In summary we get:

\begin{proposition}\label{prop:2.18} For any $J\subset I\in{\cal N}$
the projections $\pi_I$, $\lambda^I_J$ and $\pi^I_J$, and $\Pi_I$
and $\Pi^I_J$ satisfy:
 \begin{description}
 \item[(i)] $\pi_J\circ\pi^I_J=\iota^W_{IJ}\circ\pi_I$ (on
 $(\pi^I_J)^{-1}(\widetilde W_J)$) for the inclusion $\iota_{IJ}^W:W_I\hookrightarrow W_J$.
 \item[(ii)]  $\pi^I_J\circ\phi_I=\lambda^I_J(\phi_I)\circ\pi^I_J$ for any $\phi_I\in\Gamma_I$.
 \item[(iii)] $\pi_I$ induces an homeomorphism from the quotient
 $\widetilde W_I/\Gamma_I$  to $W_I$, and the restriction of $\pi_I$
 to $\widetilde W_I^\circ$  is a $|\Gamma_I|$-fold
 covering of $W_I^\circ$.
 \item[(iv)] The generic
fiber of $\pi^I_J$ contains $|\Gamma_I|/|\Gamma_J|=|\prod_{i\in
I\setminus J}\Gamma_i|$ points. Precisely, for each $\tilde u_J\in
{\rm Im}(\pi^I_J)\cap\widetilde W_J^\circ$ the inverse image
$(\pi^I_J)^{-1}(\tilde u_J)\subset\widetilde W_J^\circ$ exactly
contains $|\Gamma_I|/|\Gamma_J|$ points.
\item[(v)]
$\Pi_J\circ\Pi^I_J=\iota^E_{IJ}\circ\Pi_I$
 for  the inclusion $\iota_{IJ}^E:E_I\hookrightarrow E_J$.
 \item[(vi)]   $\Pi^I_J\circ\phi_I=\lambda^I_J(\phi_I)\circ\Pi^I_J$
for any $J\subset I\in{\cal N}$ and any  $\phi_I\in\Gamma_I$.
\item[(vii)] The obvious projection $\tilde p_I:\widetilde
E_I\to\widetilde W_I$ is equivariant with respect to the induced
actions of $\Gamma_I$ on them, i.e. $\tilde
p_I\circ\psi_I=\psi_I\circ\tilde p_I$ for any $\psi_I\in\Gamma_I$.
\item[(viii)] $\pi_I\circ\tilde p_I=p\circ\Pi_I$ for any
 $\phi_I\in\Gamma_I$, and $\pi^I_J\circ\tilde p_I=\tilde
 p_J\circ\Pi^I_J$ for any $J\subset I\in{\cal N}$.
 \item[(ix)] The generic
fiber of $\Pi^I_J$ contains $|\Gamma_I|/|\Gamma_J|$ points.
\end{description}
\end{proposition}

 However, as showed in Lemma~\ref{lem:2.19}
below the projection $\tilde p_I:\widetilde E_I\to\widetilde W_I$ is
not Banach vector bundle if $|I|>1$; accordingly we instead
temporarily call it the {\bf the virtual Banach vector bundle}.
Therefore we get a system of virtual Banach bundles
\begin{equation}\label{e:2.15}
\bigl(\widetilde{\mathcal E}(K), \widetilde W(K)\bigr)=
\bigl\{(\widetilde E_I, \widetilde W_I), \pi_I, \Pi_I, \Gamma_I,
\pi^I_J, \Pi^I_J,\,\lambda^I_J\bigm| J\subset I\in{\mathcal
N}\bigr\}.
\end{equation}

\begin{lemma}\label{lem:2.19}
For each $I\in{\mathcal N}$ with $|I|>1$,  $\tilde p_I:\widetilde
E_I\to\widetilde W_I$ is a Banach bundle variety in the sense that
locally $\widetilde W_I$ is a finite union of Banach manifolds and
the restriction of $\widetilde E_I$ on each of these finitely many
Banach manifolds is a Banach vector bundle.
\end{lemma}

\noindent{\bf Proof.}\quad Let $I=\{i_1,\cdots,i_k\}\in{\mathcal N}$
with $k>1$ and $i_1<i_2<\cdots<i_k$. For a given  $\tilde
u_I\in\widetilde W_I$ and $u_I=\pi_I(\tilde u_I)$ we choose a small
connected open neighborhood  of $u_I$ in $W_I$, $O$ (that is a
support of a Banach orbifold chart in ${\mathcal A}$) and consider
the inverse image $\widetilde O_{i_l}=\pi_{i_l}^{-1}(O)$ of $O$ in
$\widetilde W_{i_l}$, $l=1,\cdots,k$. Then each $\widetilde O_{i_l}$
is an open neighborhood of $\tilde u_{i_l}$ in $\widetilde W_{i_l}$,
and the isotropy subgroup $\Gamma(\tilde u_{i_l})$  of
$\Gamma_{i_l}$ at $\tilde u_{i_l}$ acts on $\widetilde O_{i_l}$. (
To see the second claim, for any $\phi\in \Gamma_{i_l}$ and $\tilde
x\in \widetilde O_{i_l}$ we have $\phi(\tilde x)\in\widetilde
W_{i_l}$. So it follows from $\pi_{i_l}(\phi(\tilde
x))=\pi_{i_l}(\tilde x)\in O$ that $\phi(\tilde x)\in \widetilde
O_{i_l}$. )
 If $O$ is small enough, for any fix $s\in \{1,\cdots,k\}$ it
follows from Lemma~\ref{lem:2.4} that there exist diffeomorphisms
$\lambda_{i_si_l}: \widetilde O_{i_s}\to \widetilde O_{i_l}$ that
map $\tilde u_{i_s}$ to $\tilde u_{i_l}$, and group isomorphisms
${\mathcal A}_{i_si_l}:\Gamma(\tilde u_{i_s})\to\Gamma(\tilde
u_{i_l})$  such that
\begin{equation}\label{e:2.16}
\lambda_{i_si_l}\circ\phi={\mathcal
A}_{i_si_l}(\phi)\circ\lambda_{i_si_l}
\end{equation}
 for any $\phi\in\Gamma(\tilde u_{i_s})$ and $l\in\{1,\cdots, k\}\setminus\{s\}$.
  Those  diffeomorphisms
$\lambda_{i_si_l}$ ($s\ne l$) in (\ref{e:2.16}) are unique up to
composition with elements in $\Gamma(\tilde u_{i_s})$ and
$\Gamma(\tilde u_{i_l})$. Namely, if we have another
diffeomorphism $\lambda'_{i_si_l}: \widetilde O_{i_s}\to
\widetilde O_{i_l}$ that maps $\tilde u_{i_s}$ to $\tilde
u_{i_l}$, and the group isomorphism ${\mathcal
A}'_{i_si_l}:\Gamma(\tilde u_{i_s})\to\Gamma(\tilde u_{i_l})$ such
that
\begin{equation}\label{e:2.17}
\lambda'_{i_si_l}\circ\phi={\mathcal
A}'_{i_si_l}(\phi)\circ\lambda'_{i_si_l}
\end{equation}
 for any $\phi\in\Gamma(\tilde u_{i_s})$, then there exist
 corresponding $\phi_{i_s}\in\Gamma(\tilde u_{i_s})$ and
 $\phi_{i_l}\in\Gamma(\tilde u_{i_l})$ such that
\begin{equation}\label{e:2.18}
\lambda_{i_si_l}=\phi_{i_l}\circ\lambda'_{i_si_l}\quad{\rm
and}\quad \lambda_{i_si_l}=\lambda'_{i_si_l}\circ\phi_{i_s}.
\end{equation}
In particular these imply that
\begin{equation}\label{e:2.19}
\lambda_{i_li_s}\circ\lambda_{i_si_l}=\phi_{i_s}\quad{\rm
and}\quad
\lambda_{i_ti_s}\circ\lambda_{i_li_t}\circ\lambda_{i_si_l}=\phi'_{i_s}.
\end{equation}
for some $\phi_{i_s},\phi'_{i_s}\in\Gamma(\tilde u_{i_s})$, $s\ne
l, t$ and $l\ne t$. For convenience we make the {\bf convention}:
 $\lambda_{i_si_s}=id_{\widetilde O_{i_s}}$ and ${\mathcal
A}_{i_si_s}=1_{\Gamma(\tilde u_{i_s})}$. Then we have
$|\Gamma(\tilde u_{i_s})|^{k-1}$ smooth open embeddings from
$\widetilde O_{i_s}$ into $\widetilde W_I$ given by
\begin{equation}\label{e:2.20}
\phi_I\circ\lambda_I^s:\tilde x_s\mapsto
\bigl(\phi_{i}\circ\lambda_{i_si}(\tilde x_s)\bigr)_{i\in I},
\end{equation}
 where $\lambda_I^s=(\lambda_{i_si})_{i\in I}$ with
 $\lambda_{i_si_s}=id_{\widetilde O_{i_s}}$,
 and $\phi_I=(\phi_{i})_{i\in I}$ belongs to
\begin{equation}\label{e:2.21}
\Gamma(\tilde u_I)_s:=\bigl\{(\phi_{i})_{i\in
I}\in\prod^k_{l=1}\Gamma(\tilde u_{i_l})\,\bigm|\,
\phi_{i_s}=1\bigr\}.
\end{equation}
Hereafter saying $\phi_I\circ\lambda_I^s$ to be a smooth open
embedding from $\widetilde O_{i_s}$ into $\widetilde W_I$ means
that $\phi_I\circ\lambda_I^s$ is a smooth embedding from
$\widetilde O_{i_s}$ into the Banach manifold $\widetilde
W_{i_1}\times\cdots\times\widetilde W_{i_k}$
 and that the image
$\phi_I\circ\lambda_I^s(\widetilde O_{i_s})$ is an open subset in
$\widetilde W_I$.

 Clearly,  $(\tilde
x_{i_1},\cdots,\tilde x_{i_k})\in\prod^k_{l=1}\widetilde
O_{i_l}\subset \prod^k_{l=1}\widetilde W_{i_l}$ is contained in
$\widetilde W_I$ if and only if
$$
(\tilde x_{i_1},\cdots, \tilde
x_{i_k})=\phi_I\circ\lambda_I^s(\tilde x_{i_s})$$
 for some
$\phi_I\in\Gamma(\tilde u_I)_s$. So the neighborhood
\begin{equation}\label{e:2.22}
\widetilde O(\tilde
u_I):=\pi_I^{-1}(O)=\Bigl(\prod^k_{l=1}\widetilde
O_{i_l}\Bigr)\cap\widetilde W_I
\end{equation}
of  $\tilde u_I$ in $\widetilde W_I$ can be identified with an
union of $|\Gamma(\tilde u_{i_s})|^{k-1}$ copies of $\widetilde
O_{i_s}$,
\begin{equation}\label{e:2.23}
\bigcup_{\phi_I\in\Gamma(\tilde
u_I)_s}\phi_I\circ\lambda_I^s(\widetilde O_{i_s}).
\end{equation}
Furthermore, for any two different $\phi_I,\phi'_I$ in
$\Gamma(\tilde u_I)_s$ we have
\begin{eqnarray}
\hspace{-10mm}&&\phi_I\circ\lambda_I^s(\widetilde
O_{i_s})\cap\phi'_I\circ\lambda_I^s(\widetilde O_{i_s})=
\bigl\{(\phi_i\circ\lambda_{i_si}(\tilde v))_{i\in I}\,\bigm|\,
\tilde v\in\widetilde O_{i_s}\,\&\,\nonumber\\
\hspace{-10mm}&&\hspace{20mm} (\phi_i\circ\lambda_{i_si}(\tilde
v))_{i\in I}=(\phi'_i\circ\lambda_{i_si}(\tilde v))_{i\in
I}\bigr\}\supseteq\{\tilde u_I\}.\label{e:2.24}
\end{eqnarray}

 In order to show that $\widetilde W_I$ is a Banach variety, we
need to prove that those $|\Gamma(\tilde u_{i_s})|^{k-1}$ sets,
$\phi_I\circ\lambda_I^s(\widetilde
O_{i_s}),\,\phi_I\in\Gamma(\tilde u_I)_s$, are intrinsic in the
following sense that
\begin{equation}\label{e:2.25}
\bigl\{\phi_I\circ\lambda_I^s(\widetilde O_{i_s})\bigm|
\phi_I\in\Gamma(\tilde
u_I)_s\bigr\}=\bigl\{\phi_I\circ\lambda_I^t(\widetilde
O_{i_t})\bigm| \phi_I\in\Gamma(\tilde u_I)_t\bigr\}
\end{equation}
for any two $s, t$ in $\{1,\cdots,k\}$.  For simplicity we only
prove it for the case that $k=3$ and $s=1, t=2$. Note that we may
identify $\Gamma(\tilde u_I)_1$ (resp. $\Gamma(\tilde u_I)_2$)
with $\Gamma(\tilde u_{i_2})\times\Gamma(\tilde u_{i_3})$ (resp.
$\Gamma(\tilde u_{i_1})\times\Gamma(\tilde u_{i_3})$). Then each
set in $\{\phi_I\circ\lambda_I^1(\widetilde O_{i_1})\,|\,
\phi_I\in\Gamma(\tilde u_I)_1\}$ has the form
\begin{eqnarray*}
&&(\phi_{i_2}\times\phi_{i_3})\circ\lambda_I^1(\widetilde
O_{i_1}):=\\
&&\bigl\{\bigl(\tilde x_{i_1},
\phi_{i_2}\circ\lambda_{i_1i_2}(\tilde x_{i_1}),
\phi_{i_3}\circ\lambda_{i_1i_3}(\tilde x_{i_1})\bigr)\,\bigm|\,
\tilde x_{i_1}\in\widetilde O_{i_1}\bigr\}\nonumber
\end{eqnarray*}
for some $\phi_{i_2}\times\phi_{i_3}\in\Gamma(\tilde
u_{i_2})\times\Gamma(\tilde u_{i_3})$. Similarly, each set in
$\{\phi_I\circ\lambda_I^2(\widetilde O_{i_2})\,|\,
\phi_I\in\Gamma(\tilde u_I)_2\}$ has the form
\begin{eqnarray*}
&&(\phi'_{i_1}\times\phi'_{i_3})\circ\lambda_I^2(\widetilde
O_{i_2}):=\\
&&\bigl\{\bigl(\phi'_{i_1}\circ\lambda_{i_2i_1}(\tilde x_{i_2}),
\tilde x_{i_2}, \phi'_{i_3}\circ\lambda_{i_2i_3}(\tilde
x_{i_2})\bigr)\,\bigm|\, \tilde x_{i_2}\in\widetilde
O_{i_2}\bigr\}\nonumber
\end{eqnarray*}
for some $\phi'_{i_1}\times\phi'_{i_3}\in\Gamma(\tilde
u_{i_1})\times\Gamma(\tilde u_{i_3})$. Thus in the present case
the proof of (\ref{e:2.25}) is reduced to finding
$\phi'_{i_1}\times\phi'_{i_3}\in\Gamma(\tilde
u_{i_1})\times\Gamma(\tilde u_{i_3})$  for given
$\phi_{i_2}\times\phi_{i_3}\in\Gamma(\tilde
u_{i_2})\times\Gamma(\tilde u_{i_3})$ such that
\begin{equation}\label{e:2.26}
(\phi_{i_2}\times\phi_{i_3})\circ\lambda_I^1(\widetilde O_{i_1})
=(\phi'_{i_1}\times\phi'_{i_3})\circ\lambda_I^2(\widetilde
O_{i_2}).
\end{equation}
For any $\phi'_{i_1}\in\Gamma(\tilde u_{i_1})$, since
$\phi'_{i_1}\circ\lambda_{i_2i_1}:\widetilde O_{i_2}\to\widetilde
O_{i_1}$ is a diffeomorphism it is easily checked that the set
\begin{eqnarray*}
&&(\phi_{i_2}\times\phi_{i_3})\circ\lambda^1_I(\widetilde
O_{i_1})\\
&&=\bigl\{\bigl(\tilde x_{i_1},
\phi_{i_2}\circ\lambda_{i_1i_2}(\tilde x_{i_1}),
\phi_{i_3}\circ\lambda_{i_1i_3}(\tilde x_{i_1})\bigr)\,\bigm|\,
\tilde x_{i_1}\in\widetilde O_{i_1}\bigr\}
\end{eqnarray*}
is equal to that of all triples
\begin{eqnarray*}
\bigl(\phi'_{i_1}\circ\lambda_{i_2i_1}(\tilde x_{i_2}),
\phi_{i_2}\circ\lambda_{i_1i_2}(\phi'_{i_1}\circ\lambda_{i_2i_1}(\tilde
x_{i_2})),
\phi_{i_3}\circ\lambda_{i_1i_3}(\phi'_{i_1}\circ\lambda_{i_2i_1}(\tilde
x_{i_2}))\bigr)\\
= \bigl(\phi'_{i_1}\circ\lambda_{i_2i_1}(\tilde x_{i_2}),
\phi_{i_2}\circ\lambda_{i_1i_2}\circ\phi'_{i_1}\circ\lambda_{i_2i_1}(\tilde
x_{i_2}),
\phi_{i_3}\circ\lambda_{i_1i_3}\circ\phi'_{i_1}\circ\lambda_{i_2i_1}(\tilde
x_{i_2})\bigr),
\end{eqnarray*}
where $\tilde x_{i_2}$ takes over $\widetilde O_{i_2}$.
 By (\ref{e:2.17}),
 (\ref{e:2.18}) and (\ref{e:2.19}) we may obtain that
\begin{eqnarray*}
\phi_{i_2}\circ\lambda_{i_1i_2}\circ\phi'_{i_1}\circ\lambda_{i_2i_1}
\!\!\!\!\!\!\!\!&&=\phi_{i_2}\circ{\mathcal
A}_{i_1i_2}(\phi'_{i_1})\circ\lambda_{i_1i_2}\circ\lambda_{i_2i_1}\\
&&=\phi_{i_2}\circ{\mathcal
A}_{i_1i_2}(\phi'_{i_1})\circ\phi''_{i_2}\nonumber
\end{eqnarray*}
for some $\phi''_{i_2}\in\Gamma(\tilde u_{i_2})$, and that
\begin{eqnarray*}
\phi_{i_3}\circ\lambda_{i_1i_3}\circ\phi'_{i_1}\circ\lambda_{i_2i_1}
\!\!\!\!\!\!\!\!&&=\phi_{i_3}\circ{\mathcal
A}_{i_1i_3}(\phi'_{i_1})\circ\lambda_{i_1i_3}\circ\lambda_{i_2i_1}\\
&&=\phi_{i_2}\circ{\mathcal
A}_{i_1i_2}(\phi'_{i_1})\circ\phi'''_{i_2}\circ\lambda_{i_2i_3}\nonumber
\end{eqnarray*}
for some $\phi'''_{i_2}\in\Gamma(\tilde u_{i_2})$. Hence
(\ref{e:2.26}) holds  if we take
\begin{equation*}
\phi'_{i_3}=\phi_{i_2}\circ{\mathcal
A}_{i_1i_2}(\phi'_{i_1})\circ\phi'''_{i_2}\quad{\rm and}\quad
\phi'_{i_1}= {\mathcal
A}_{i_1i_2}^{-1}(\phi^{-1}_{i_2}\circ(\phi''_{i_2})^{-1}).
\end{equation*}

 Similarly, (by shrinking
$O$ if necessary) using the properties of orbibundles one has also
smooth bundle isomorphisms
$$
\Lambda_{i_si_l}:\widetilde E_{i_s}|_{\widetilde
O_{i_s}}\to\widetilde E_{i_l}|_{\widetilde O_{i_l}}
$$
 covering $\lambda_{i_si_l}$, i.e. $\lambda_{i_si_l}\cdot\tilde
p_{i_s}=\tilde p_{i_l}\cdot\Lambda_{i_si_l}$, such that
$$
\Lambda_{i_si_l}\circ\Phi=\boxed{{\mathcal
A}_{i_si_l}(\phi)}\circ\Lambda_{i_si_l}
$$
 for any $\phi\in\Gamma(\tilde u_{i_s})$ and $l\in\{1,\cdots, k\}\setminus\{s\}$.
Hereafter
\begin{equation}\label{e:2.27}
\Phi: \widetilde
E_{i_s}|_{\widetilde O_{i_s}}\to\widetilde E_{i_s}|_{\widetilde
O_{i_s}}\; \bigl({\rm resp.}\;  \boxed{{\cal
 A}_{i_si_l}(\phi)}: \widetilde
E_{i_l}|_{\widetilde O_{i_l}}\to\widetilde E_{i_l}|_{\widetilde
O_{i_l}}\bigr)
\end{equation}
is the  bundle isomorphism lifting of $\phi$ (resp. ${\cal
A}_{i_si_l}(\phi)$) produced in the definition of orbifold bundles.
({\bf As stated in Remark~\ref{rem:2.8} we hereafter write
$(\widetilde E_{i_l}, \widetilde\Gamma_{i_l},\Pi_{i_l})$ as
$(\widetilde E_{i_l}, \Gamma_{i_l},\Pi_{i_l})$ and understand
$\phi_{i_l}\in\Gamma_{i_l}$ to acts on $\widetilde E_{i_l}$ via
$\phi_{i_l}\cdot\tilde\xi=\Phi_{i_l}(\tilde\xi)$ with
$\Phi_{i_l}=\lambda_{W_{i_l}}(\phi_{i_l})$ and
$\tilde\xi\in\widetilde E_{i_l}$}.)
 It is easy to see that $\Phi_{i_l}$ maps $\widetilde
E_{i_l}|_{\widetilde O_{i_l}}$ to $\widetilde E_{i_l}|_{\widetilde
O_{i_l}}$. In particular the restriction of it to $\widetilde
E_{i_l}|_{\widetilde O_{i_l}}$ is also a smooth bundle
automorphism of $\widetilde E_{i_l}|_{\widetilde O_{i_l}}$
covering $\phi_{i_l}|_{\widetilde O_{i_l}}$, and is still denoted
by $\Phi_{i_l}$. So for the above
$\lambda_I^s=(\lambda_{i_si})_{i\in I}$ with
 $\lambda_{i_si_s}=id_{\widetilde O_{i_s}}$,
 and $\phi_I=(\phi_{i})_{i\in I}\in\Gamma(\tilde u_I)_s$ we get
 the corresponding
\begin{equation}\label{e:2.28}
 \Phi_I:=(\Phi_{i})_{i\in
I}\quad{\rm and}\quad\Lambda_I^s:=(\Lambda_{i_si})_{i\in I},
\end{equation}
 where we have made the {\bf convention}:
$$
 \Phi_{i_s}=1\quad{\rm and}\quad
\Lambda_{i_si_s}=id_{\widetilde E_{i_s}|_{\widetilde O_{i_s}}}.
$$
 Corresponding with (\ref{e:2.20}) we obtain $|\Gamma(\tilde
u_{i_s})|^{k-1}$ smooth bundle open embeddings
\begin{equation}\label{e:2.29}
\Phi_I\circ\Lambda_I^s:\widetilde E_{i_s}|_{\widetilde
O_{i_s}}\to\widetilde E_I,\; \tilde \xi\mapsto
\bigl(\Phi_i\circ\Lambda_{i_si}(\tilde \xi)\bigr)_{i\in I}
\end{equation}
for $\phi_I\in\Gamma(\tilde u_I)_s$. Here ``smooth bundle open
embedding'' is to be understood in a similar way as is explained
below (\ref{e:2.21}).
  Clearly, each map
\begin{eqnarray}
 && \Phi_I\circ\Lambda_I^s(\widetilde E_{i_s}|_{\widetilde
O_{i_s})})\to \phi_I\circ\lambda_I^s(\widetilde
O_{i_s}),\label{e:2.30}\\
&& \bigl(\Phi_i\circ\Lambda_{i_si}(\tilde \xi)\bigr)_{i\in
I}\mapsto \bigl(\phi_i\circ\lambda_{i_si}(\tilde p_{i_s}(\tilde
\xi))\bigr)_{i\in I}\nonumber
\end{eqnarray}
 is a  Banach vector bundle over the  Banach
manifold $\phi_I\circ\lambda_I^s(\widetilde O_{i_s})$, and
 $$
 \bigcup_{\phi_I\in\Gamma(\tilde u_I)_s}
 \Phi_I\circ\Lambda_I^s(\widetilde E_{i_s}|_{\widetilde
O_{i_s}})=\widetilde E_I\Bigm|_{\bigcup_{\phi_I\in\Gamma(\tilde
u_I)_s}\phi_I\circ\lambda_I^s(\widetilde O_{i_s})}.
$$
   Corresponding with (\ref{e:2.23}) and
 (\ref{e:2.24}) we have also:
\begin{eqnarray*}
&&\Phi_I\circ\Lambda_I^s(\widetilde E_{i_s}|_{\widetilde
O_{i_s}})\bigcap\Phi'_I\circ\Lambda_I^s(\widetilde
E_{i_s}|_{\widetilde
O_{i_s}})=\\
&&\bigl\{\bigl(\Phi_i\circ\Lambda_{i_si}(\tilde \xi)\bigr)_{i\in
I}\bigm|\tilde \xi\in\widetilde E_{i_s}|_{\widetilde
O_{i_s}}\,\&\, \bigl(\Phi_i\circ\Lambda_{i_si}(\tilde
\xi)\bigr)_{i\in I}=\bigl(\Phi'_i\circ\Lambda_{i_si}(\tilde
\xi)\bigr)_{i\in I}\bigr\}.
\end{eqnarray*}
for any two different $\phi_I,\psi_I$ in $\Gamma(\tilde u_I)_s$,
and
\begin{equation}\label{e:2.31}
\left.\begin{array}{ll}
 \bigl\{\Phi_I\circ\Lambda_I^s(\widetilde
E_{i_s}|_{\widetilde O_{i_s}})\,\bigm|\, \phi_I\in\Gamma(\tilde
u_I)_s\bigr\}\\
=\bigl\{\Phi_I\circ\Lambda_I^t(\widetilde
E_{i_t}|_{\widetilde O_{i_t}})\,\bigm|\, \phi_I\in\Gamma(\tilde
u_I)_t\bigr\}
\end{array}\right\}
\end{equation}
for any  $s, t\in\{1,\cdots,k\}$. Hence $\widetilde E_I\to
\widetilde W_I$ is a Banach bundle variety.
 \hfill$\Box$\vspace{2mm}

\begin{remark}\label{rem:2.20}
{\rm  By shrinking $O$ we may also assume that  each
$\lambda_{i_si_l}$ in (\ref{e:2.16}) is not only  a diffeomorphism
from $\widetilde O_{i_s}$ to $\widetilde O_{i_l}$, but also one from
a neighborhood of the closure of $\widetilde O_{i_s}$ to that of
$\widetilde O_{i_l}$. In this case the union in (\ref{e:2.23}) may
be required to satisfy:
\begin{equation}\label{e:2.32}
\left.\begin{array}{ll}
 {\rm If}\; \tilde
v\in\bigcup_{\phi_I\in\Gamma(\tilde
u_I)_s}\phi_I\circ\lambda_I^s(\widetilde O_{i_s})\;{\rm
and}\;\tilde v\notin \phi_I\circ\lambda_I^t(\widetilde
O_{i_t})\\{\rm for}\;{\rm some}\,t,\;{\rm then}\;  \tilde v\notin
Cl\bigr(\phi_I\circ\lambda_I^t(\widetilde O_{i_t})\bigl).
\end{array}\right\}
\end{equation}
It means that each $\phi_I\circ\lambda_I^t(\widetilde O_{i_t})$ is
relatively closed in $\bigcup_{\phi_I\in\Gamma(\tilde
u_I)_s}\phi_I\circ\lambda_I^s(\widetilde O_{i_s})$. Let us prove
(\ref{e:2.32}) by contradiction arguments. Suppose that $\tilde
v\in Cl\bigr(\phi_I\circ\lambda_I^t(\widetilde
O_{i_t})\bigl)\setminus\phi_I\circ\lambda_I^t(\widetilde
O_{i_t})$. Since $Cl\bigr(\phi_I\circ\lambda_I^t(\widetilde
O_{i_t})\bigl)=\phi_I\circ\lambda_I^t(Cl(\widetilde O_{i_t}))$
there exists a unique $\tilde x\in Cl(\widetilde O_{i_t})$, which
must sit in $Cl(\widetilde O_{i_t})\setminus\widetilde O_{i_t}$,
such that $\phi_I\circ\lambda_I^t(\tilde x)=\tilde v$. Moreover,
since $\tilde v$ belongs to some
$\phi'_I\circ\lambda_I^s(\widetilde O_{i_s})$ one has a unique
$\tilde y\in \widetilde O_{i_s}$ such that
$\phi'_i\circ\lambda_I^s(\tilde y)=\tilde v$. So
$\phi_I\circ\lambda_I^t(\tilde x)=\phi'_I\circ\lambda_I^s(\tilde
y)$. For the sake of clearness we assume $t=1$ and $s=2$ then
\begin{eqnarray*}
&&\bigl(\tilde x, \phi_{i_2}\circ\lambda_{i_1i_2}(\tilde
x),\cdots,
\phi_{i_k}\circ\lambda_{i_1i_k}(\tilde x)\bigr)\\
&&=\bigl(\phi'_{i_1}\circ\lambda_{i_2i_1}(\tilde y), \tilde
y,\phi'_{i_3}\circ\lambda_{i_3i_1}(\tilde y),\cdots,
\phi'_{i_k}\circ\lambda_{i_2i_k}(\tilde x)\bigr).
\end{eqnarray*}
Therefore $\tilde x=\phi'_{i_1}\circ\lambda_{i_2i_1}(\tilde y)$.
Note that  $\lambda_{i_2i_1}:\widetilde O_{i_2}\to\widetilde
O_{i_1}$ and $\phi'_{i_1}:\widetilde O_{i_1}\to\widetilde O_{i_1}$
are both  diffeomorphisms. We deduce that $\tilde
x=\phi'_{i_1}\circ\lambda_{i_2i_1}(\tilde y)\in \widetilde
O_{i_1}$. This contradicts that $\tilde x\in Cl(\widetilde
O_{i_1})\setminus\widetilde O_{i_1}$.  (\ref{e:2.32}) is proved.

Similarly, we also require that
\begin{equation}\label{e:2.33}
\left.\begin{array}{ll} {\rm If}\; \tilde
\xi\in\bigcup_{\Phi_I\in\Gamma(\tilde
u_I)_s}\Phi_I\circ\Lambda_I^s(\widetilde O_{i_s})\;{\rm
and}\;\tilde \xi\notin \Phi_I\circ\Lambda_I^t(\widetilde
O_{i_t})\\
{\rm for}\;{\rm some}\,t,\; {\rm then}\; \tilde \xi\notin
Cl\bigr(\Phi_I\circ\Lambda_I^t(\widetilde O_{i_t})\bigl).
\end{array}\right\}
\end{equation}
Later we always assume that (\ref{e:2.32}) and (\ref{e:2.33}) are
satisfied without special statements.}
\end{remark}

In terms of [LiuT3] we introduce:

\begin{definition}\label{def:2.21}{\rm
The family  of the smooth embeddings
 given by (\ref{e:2.20}),
 $\bigl\{\phi_I\circ\lambda_I^s\,|\, \phi_I\in\Gamma(\tilde u_I)_s\bigr\}$,
  is called a
{\bf local coordinate chart of $\widetilde W_I$ over a neighborhood
$\widetilde O(\tilde u_I)$ of $\tilde u_I$}, and each
$\phi_I\circ\lambda_I^s$ is called a {\bf component} of this chart.
Similarly, we call the family of the  bundle embeddings given by
(\ref{e:2.29}), $\bigl\{\Phi_I\circ\Lambda_I^s\,|\,
\phi_I\in\Gamma(\tilde u_I)_s\bigr\}$,  a {\bf local bundle
coordinate chart of $\widetilde E_I$ over $\widetilde O(\tilde
u_I)$}, and each $\Phi_I\circ\Lambda_I^s$ a {\bf component} of it.
For $\tilde x_I\in\phi_I\circ\lambda_I^s(\widetilde O_{i_s})$, a
connected relative open subset
$W\subset\phi_I\circ\lambda_I^s(\widetilde O_{i_s})$ containing
$\tilde x_I$ is called a {\bf local component of $\widetilde W_I$
near} $\tilde x_I$. In particular,
$\phi_I\circ\lambda_I^s(\widetilde O_{i_s})$ is a local component of
$\widetilde W_I$ near $\tilde u_I$. The restriction of
$\Phi_I\circ\Lambda_I^s(\widetilde E_{i_s}|_{\widetilde O_{i_s}})$
to a  local component near $\tilde u_I$ is called a {\bf local
component of $\widetilde E_I$ near $\tilde u_I$}. Two local
components of $\widetilde W_I$ (or $\widetilde E_I$) near a point
$\tilde u$ of $\widetilde W_I$ are said to be of  {\bf different
kind} if the intersection of both is not a local component of
$\widetilde W_I$ near $\tilde u$. Let ${\mathcal A}$ be a family of
local components of different kind near a point $\tilde u\in
\widetilde W_I$. If the union of sets in ${\mathcal A}$ forms an
open neighborhood ${\mathcal O}(\tilde u)$ of $\tilde u$ in
$\widetilde W_I$ we call ${\mathcal A}$ {\bf a complete family of
local components of $\widetilde W_I$ over ${\mathcal O}(\tilde u)$}.
Clearly, we have also a notion of germ of complete families of local
components.}
 \end{definition}

Note that (\ref{e:2.25}) and (\ref{e:2.31}) show that one can
construct the same complete families of local components starting
from two different $s$ and $t$. Moreover, suppose that $Q\subset O$
is another small open neighborhood of $u_I$ in $W_I$. Then all
constructions in the proof of Lemma~\ref{lem:2.19} work as we
replace $O$ with $Q$, and the corresponding local components of
$\widetilde W_I$ near $\tilde u_I$ and those of $\widetilde E_I$
near $(\widetilde E_I)_{\tilde u_I}$ are, respectively, given by
\begin{eqnarray*}
&&\bigl\{\phi_I\circ\lambda_I^s(\widetilde Q_{i_s})\,\bigm|\,
\phi_I\in\Gamma(\tilde u_I)_s\bigr\}\quad{\rm and}\\
&&\bigl\{\Phi_I\circ\Lambda_I^s(\widetilde E_{i_s}|_{\widetilde
Q_{i_s}})\,\bigm|\, \phi_I\in\Gamma(\tilde u_I)_s\bigr\}.
\end{eqnarray*}
Here  $\widetilde Q_{i_l}=\pi_{i_l}^{-1}(Q)$ are the inverse
images of $Q$ in $\widetilde W_{i_l}$, $l=1,\cdots,k$, and
$\phi_I\circ\lambda_I^s$ and $\Phi_I\circ\Lambda_I^s$ as in
(\ref{e:2.25}) and (\ref{e:2.31}).  More general conclusions,
stated in Lemma 4.3 of [LiuT3], will be:

\begin{proposition}\label{prop:2.22}
The notion of local component is functorial with respect to
restriction and projection, that is:
\begin{description}
\item[(i)] If $\{\phi_I\circ\lambda_I^s\,|\,\phi_I\in\Gamma(\tilde
u_I)_s\}$ is a local coordinate chart near $\tilde
u_I\in\widetilde W_I$, then for each point $\tilde
v_I\in\cup_{\phi_I\in\Gamma(\tilde
u_I)_s}\phi_I\circ\lambda_I^s(\widetilde O_{i_s})$ there exists a
neighborhood ${\mathcal Q}$ of $\tilde v_I$ in $\widetilde W_I$
such that $\{{\mathcal Q}\cap \phi_I\circ\lambda_I^s(\widetilde
O_{i_s})\,|\,\phi_I\in\Gamma(\tilde u_I)_s\}$ is a complete family
of local components of $\widetilde W_I$ over ${\mathcal Q}$.

\item[(ii)] If
$\{\phi_I\circ\lambda_I^s\,|\,\phi_I\in\Gamma(\tilde u_I)_s\}$ is
a local coordinate chart near $\tilde u_I\in\widetilde W_I$ and
$i_s\in J\subset I$, then
 $\{\pi^I_J\circ\phi_I\circ\lambda_I^s\,|\,\phi_I\in\Gamma(\tilde
u_I)_s\}$ is a local coordinate chart near $\pi^I_J(\tilde
u_I)\in\widetilde W_J$ (after  the repeating maps
$\pi^I_J\circ\phi_I\circ\lambda_I^s$ is only taken one at a time).
In particular, not only the projection $\pi^I_J:\widetilde
W_I\to\widetilde W_J$ maps  the local component
$\phi_I\circ\lambda_I^s(\widetilde O_{i_s})$ near $\tilde u_I$ to
the local component $\pi^I_J(\phi_I\circ\lambda_I^s(\widetilde
O_{i_s}))$ near $\pi^I_J(\tilde u_I)$ in $\widetilde W_J$ but also
the restriction
$$\pi^I_J|_{\phi_I\circ\lambda_I^s(\widetilde O_{i_s})}:
\phi_I\circ\lambda_I^s(\widetilde O_{i_s})\to
\pi^I_J(\phi_I\circ\lambda_I^s(\widetilde O_{i_s}))
$$
is a smooth map. Actually one has
$$\pi^I_J\circ\phi_I\circ\lambda_I^s=\lambda^I_J(\phi_I)\circ\lambda_J^s
=\phi_J\circ\lambda_J^s,$$
 where $\lambda_J^s=(\lambda_{i_sl})_{l\in J}$.

\item[(iii)] Correspondingly,  the projection $\Pi^I_J:\widetilde
E_I\to\widetilde E_J$ restricts to a smooth bundle map from
$\widetilde E_I|_{\phi_I\circ\lambda_I^s(\widetilde O_{i_s})}$ to
$\widetilde E_J|_{\pi^I_J(\phi_I\circ\lambda_I^s(\widetilde
O_{i_s}))}$.
\end{description}
\end{proposition}

\noindent{\bf Proof.} (i) By (\ref{e:2.32}) we may choose a
connected open neighborhood ${\mathcal Q}\subset\pi_I^{-1}(O)$ of
$\tilde v$ in $\widetilde W_I$ such that for any
$\phi_I\in\Gamma(\tilde u_I)_s$,
$${\mathcal Q}\cap Cl(\phi_I\circ\lambda_I^s(\widetilde
O_{i_s}))=\emptyset\quad {\rm as}\quad \tilde
v_I\notin\phi_I\circ\lambda_I^s(\widetilde O_{i_s}).$$
 If $v_I\in\phi_I\circ\lambda_I^s(\widetilde O_{i_s})$ then the connected relative
 open subset ${\mathcal Q}\cap \phi_I\circ\lambda_I^s(\widetilde
O_{i_s})$ in $\phi_i\circ\lambda_I^s(\widetilde O_{i_s})$ is a
local component of $\widetilde W_I$ near $\tilde v_I$. It is easy
to see that the union of sets in  $\{{\mathcal Q}\cap
\phi_I\circ\lambda_I^s(\widetilde
O_{i_s})\,|\,\phi_I\in\Gamma(\tilde u_I)_s\}$ is equal to
${\mathcal Q}$. So (i) holds.

(ii) We may assume that $I=\{i_1,\cdots,i_k\}$ and
$J=\{i_1,\cdots,i_{k'}\}$ with $k'<k$. By (\ref{e:2.24}) we may
also take $s=1$. Then by (\ref{e:2.21}) we have
\begin{eqnarray*}
&&\Gamma(\tilde
u_I)_1=\bigl\{(1,\phi_{i_2},\cdots,\phi_{i_k})\bigm|
\phi_{i_l}\in\Gamma(\tilde u_{i_l}), l=2,\cdots,k\bigr\},\\
&&\Gamma(\pi^I_J(\tilde
u_I))_1=\bigl\{(1,\phi_{i_2},\cdots,\phi_{i_{k'}})\bigm|
\phi_{i_l}\in\Gamma(\tilde u_{i_l}), l=2,\cdots,k'\bigr\}.
\end{eqnarray*}
That is, $\lambda^I_J(\Gamma(\tilde u_I)_1)=\Gamma(\pi^I_J(\tilde
u_I))_1$.  So for a given local component
$$\phi_I\circ\lambda_I^1(\widetilde O_{i_1})=\bigl\{
\bigl(\tilde x_1, \phi_{i_2}\circ\lambda_{i_1i_2}(\tilde
x_1),\cdots, \phi_{i_k}\circ\lambda_{i_1i_k}(\tilde
x_1)\bigr)\bigm|\tilde x_1\in\widetilde O_{i_1}\bigr\}
$$
of $\widetilde W_I$ near $\tilde u_I$, we have a corresponding
local component of $\widetilde W_J$ near $\pi^I_J(\tilde u_I)$ as
follow:
$$
\phi_J\circ\lambda_J^1(\widetilde O_{i_1})=\bigl\{ \bigl(\tilde
x_1, \phi_{i_2}\circ\lambda_{i_1i_2}(\tilde x_1),\cdots,
\phi_{i_k'}\circ\lambda_{i_1i_{k'}}(\tilde x_1)\bigr)\bigm|\tilde
x_1\in\widetilde O_{i_1}\bigr\}.
$$
Here
$\phi_J=(1,\phi_{i_2},\cdots,\phi_{i_{k'}})=\lambda^I_J(\phi_I)\in\Gamma(\pi^I_J(\tilde
u_I))_1$ and $\lambda_J^1=(\phi_{i_1l})_{l\in J}$.
 Clearly,
 $\pi^I_J\circ\phi_I\circ\lambda_I^1=\lambda^I_J(\phi_I)\circ\lambda_J^1=\phi_J\circ\lambda_J^1$,
 and so   $\pi^I_J\bigl(\phi_I\circ\lambda_I^1(\widetilde
O_{i_1})\bigr)=\phi_J\circ\lambda_J^1(\widetilde O_{i_1})$.
Moreover the restriction
$$
\pi^I_J|_{\phi_I\circ\lambda_I^1(\widetilde
O_{i_1})}:\phi_I\circ\lambda_I^1(\widetilde O_{i_1})\to
\phi_J\circ\lambda_J^1(\widetilde O_{i_1}),
$$
which maps $\bigl(\tilde x_1,
\phi_{i_2}\circ\lambda_{i_1i_2}(\tilde x_1),\cdots,
\phi_{i_k}\circ\lambda_{i_1i_k}(\tilde x_1)\bigr)$ to
$$\bigl(\tilde x_1, \phi_{i_2}\circ\lambda_{i_1i_2}(\tilde
x_1),\cdots, \phi_{i_{k'}}\circ\lambda_{i_1i_{k'}}(\tilde
x_1)\bigr),$$
  is explicitly  a smooth map. (ii) is proved. (iii) is also
  proved easily.
\hfill$\Box$\vspace{2mm}

The arguments below Definition~\ref{def:2.21} and
Proposition~\ref{prop:2.22} show that the {\bf notion of local
component is intrinsic}. Using this, for each
$I=\{i_1,\cdots,i_k\}\in{\mathcal N}$ with $|I|>1$ we can
desingularize $\widetilde W_I$ to get a true   Banach manifold
$\widehat W_I$. Consider the disjoint union of Banach manifolds
\begin{equation}\label{e:2.34}
\boxed{\widetilde W_I\!}:= \coprod_{\tilde u_I\in\widetilde
W_I}\coprod_{\phi_I\in\Gamma(\tilde u_I)_1}\coprod_{O\in{\mathcal
U}(u_I)}\{(\tilde
u_I,\phi_I)\}\times\phi_I\circ\lambda_I^1(\widetilde O_{i_1}),
\end{equation}
where ${\mathcal U}(u_I)$ is the germ of small connected open
neighborhoods of $u_I=\pi_I(\tilde u_I)$ in $W_I$ (that are supports
of charts in the atlas ${\mathcal A}$) and $\widetilde
O_{i_1}=\pi^{-1}_{i_1}(O)$. In $\boxed{\widetilde W_I\!}$ we define
a relation $\sim$ as follows:
 For $\tilde u_I,\tilde u'_I\in\widetilde W_I$,
 and $\phi_I\in\Gamma(\tilde u_I)_1$, $\phi'_I\in\Gamma(\tilde u'_I)_1$,
 and $\tilde y_I\in\phi_I\circ\lambda_I^1(\widetilde O_{i_1})$
 and $\tilde y_I^\prime\in\phi_I^\prime\circ\lambda_I^{\prime 1}(\widetilde O'_{i_1})$ we define
\begin{equation}\label{e:2.35}
(\tilde u_I,\phi_I, \tilde y_I)\sim(\tilde u'_I,\phi'_I,
 \tilde y'_I)
\end{equation}
if and only if the following two conditions are satisfied:
\begin{description}
\item[(i)] $\tilde y_I=\tilde y'_I$ as points in $ \widetilde
W_I$;

 \item[(ii)] $\phi_I\circ\lambda_I^1(\widetilde O_{i_1})
\bigcap\phi'_I\circ\lambda_I^{\prime 1}(\widetilde O'_{i_1})$ is a
local component of $\widetilde W_I$ near $\tilde y_I=\tilde y'_I$.
\end{description}

By Proposition~\ref{prop:2.22} we immediately obtain:

\begin{claim}\label{cl:2.23}
The relation $\sim$ in (\ref{e:2.35}) is a {\rm regular equivalence}
 in $\boxed{\widetilde W_I\!}\,$. Moreover, any two different
points in $\{(\tilde
u_I,\phi_I)\}\times\phi_I\circ\lambda_I^1(\widetilde O_{i_1})$ are
not equivalent with respect to $\sim$.
\end{claim}

 In \cite[Def.4.3.1, Def.4.3.3]{MaOu}) an equivalence relation
${\mathcal R}\subset X\times X$ on a $C^k$-manifold $X$ is said to
be {\bf regular} if there is a $C^k$-differentiable structure on
$X/{\mathcal R}$ such that the projection map $q:X\to X/{\mathcal
R}$ is a $C^k$-submersion. This unique $C^k$-differentiable
structure on $X/{\mathcal R}$ is called {\bf quotient differentiable
structure} and the associated manifold is called {\bf quotient
manifold} given by ${\mathcal R}$ on $X$. An equivalence relation
${\mathcal R}$ on $X$ is said to preserves the boundary if
${\mathcal R}(\partial X)=\partial X$. For the precise relations
between $\partial X$ and $X/{\mathcal R}$ the reader may refer to
\cite[Th.4.3.14]{MaOu}) (and in particular \cite[Cor.4.3.15]{MaOu}
and \cite[Th.4.3.17]{MaOu} in two special cases that $\partial
X=\emptyset$ and $\partial(X/{\mathcal R})=\emptyset$).\vspace{2mm}

\noindent{\bf Proof of Claim~\ref{cl:2.23}.} We only need to prove
the regularity. For the sake of simplicity we assume that $X$ and
thus $\boxed{\widetilde W_I\!}$ have no boundary. Let ${\cal R}$ be
a subset in $\boxed{\widetilde W_I\!}\times\boxed{\widetilde W_I\!}$
determined by the equivalence relation $\sim$. By Corollary 4.3.15
in \cite{MaOu} we need to show that ${\cal R}$ is a submanifold of
$\boxed{\widetilde W_I\!}\times\boxed{\widetilde W_I\!}$ without
boundary and that $P_1|_{\cal R}:{\cal R}\to\boxed{\widetilde
W_I\!}$ is a submersion. Here $P_1$ is the projection of
$\boxed{\widetilde W_I\!}\times\boxed{\widetilde W_I\!}$ to the
first factor. Let $(\tilde u_I,\phi_I, \tilde y_I)\sim(\tilde
u'_I,\phi'_I,
 \tilde y'_I)$ be as in (\ref{e:2.35}). Then
$\bigl((\tilde u_I,\phi_I, \tilde y_I), (\tilde u'_I,\phi'_I,
 \tilde y'_I)\bigr)$ belongs to ${\cal R}$. Note that for any
$\tilde z_I\in\phi_I\circ\lambda_I^1(\widetilde O_{i_1})
\bigcap\phi'_I\circ\lambda_I^{\prime 1}(\widetilde O'_{i_1})$ we
have
$$
\bigl((\tilde u_I,\phi_I, \tilde z_I), (\tilde u'_I,\phi'_I,
 \tilde z_I)\bigr)\in {\cal R}.
 $$
 So a neighborhood of $\bigl((\tilde u_I,\phi_I, \tilde y_I), (\tilde u'_I,\phi'_I,
 \tilde y'_I)\bigr)$  in ${\cal R}$ can be identified with
$\phi_I\circ\lambda_I^1(\widetilde O_{i_1})
\bigcap\phi'_I\circ\lambda_I^{\prime 1}(\widetilde O'_{i_1})$, and
in this situation the restriction of $P_1|_{\cal R}$ to this
neighborhood is the identity map. The desired conclusions follow
immediately. \hfill$\Box$\vspace{2mm}

 Denote by $[\tilde
u_I,\phi_I, \tilde y_I]$ the equivalence class of $(\tilde
u_I,\phi_I, \tilde y_I)$, and by
\begin{equation}\label{e:2.36}
\widehat W_I=\boxed{\widetilde W_I\!}\bigg/\sim
\end{equation}
 if $|I|>1$, and by $\widehat W_I=\widetilde W_i$ if
 $I=\{i\}\in{\mathcal N}$. Then $\widehat
 W_I$ is a Banach manifold, and each
\begin{equation}\label{e:2.37}
\widehat{\phi_I\circ\lambda_I^1}:\widetilde O_{i_1}\to \widehat
W_I,\;\tilde x\mapsto [\tilde u_I,\phi_I,
\phi_I\circ\lambda_I^1(\tilde x)],
\end{equation}
gives a smooth open embedding (which is called {\bf a local
coordinate chart of $\widehat W_I$ near $[\tilde u_I,\phi_I,\tilde
u_I]$}). The projections $\pi_I$ and $\pi^I_J$ induce natural ones
\begin{equation}\label{e:2.38}
\left.\begin{array}{ll}
\hat\pi_I: \widehat W_I\to W_I,\quad
[\tilde
u_I,\phi_I, \tilde y_I]\mapsto \pi_I(\tilde y_I),\\
 \hat\pi^I_J:\widehat W_I\to \widehat W_J,\;[\tilde u_I,\phi_I, \tilde
y_I]\mapsto [\tilde u_J,\phi_J, \tilde y_J]
\end{array}\right\}
\end{equation}
Here $\phi_J=\lambda^I_J(\phi_I)=(\phi_l)_{l\in J}$, $\tilde
u_J=\pi^I_J(\tilde u_I)=(\tilde u_l)_{l\in J}$ and $\tilde y_J=
\pi^I_J(\tilde y_I)=(\tilde y_l)_{l\in J}$.
 By Proposition~\ref{prop:2.22}(ii) the map $\hat\pi^I_J$ is
well-defined, and  as in the arguments below (\ref{e:2.14}) we can
show that $\hat\pi^I_J$ is not surjective in general if $I\ne
 J$. Note that in the charts $\widehat{\phi_I\circ\lambda_I^1}$ and
$$
\widehat{\phi_J\circ\lambda_J^1}:\widetilde O_{j_1}\to \widehat
W_J,\;\tilde x\mapsto [\tilde u_J,\phi_J,
\phi_J\circ\lambda_J^1(\tilde x)],
$$
 the projection $\hat\pi^I_J$ may be represented by
\begin{equation}\label{e:2.39}
 \widetilde O_{i_1}\to\widetilde O_{j_1},\; \tilde x\mapsto
\lambda_{i_1j_1}(\tilde x)
\end{equation}
because $j_1\in J\subset I$. Here $\lambda_{i_1j_1}=1_{\widetilde
O_{i_1}}$ if $j_1=i_1$.
 So $\hat\pi^I_J$ {\bf is a smooth map}.

Recall that the action of any $\psi_I=(\psi_l)_{l\in
I}\in\Gamma_I$ on
\begin{eqnarray*}
\tilde y_I=\phi_I\circ\lambda_I^1(\tilde x)\!\!\!\!\!\!\!\!\!\!&&=
\bigl(\tilde x, \phi_{i_2}\circ\lambda_{i_1i_2}(\tilde x),\cdots,
\phi_{i_k}\circ\lambda_{i_1i_k}(\tilde x)\bigr)\\
 &&\in \phi_I\circ\lambda_I^1(\widetilde O_{i_1})\subset\widetilde
 W_I
\end{eqnarray*}
is given by
\begin{eqnarray}
\psi_I(\tilde y_I)\!\!\!\!\!\!\!\!\!\!&&=\bigl(\psi_{i_1}(\tilde
x), \psi_{i_2}\circ\phi_{i_2}\circ\lambda_{i_1i_2}(\tilde
x),\cdots, \psi_{i_k}\circ\phi_{i_k}\circ\lambda_{i_1i_k}(\tilde
x)\bigr)\nonumber\\
&&=\bigl(\psi_{i_1}(\tilde x),
\phi^\ast_{i_2}\circ\lambda^\ast_{i_1i_2}(\psi_{i_1}(\tilde
x)),\cdots,
\phi^\ast_{i_k}\circ\lambda^\ast_{i_1i_k}(\psi_{i_1}(\tilde
x)\bigr)\nonumber\\
&&=\phi^\ast_I\circ\lambda_I^{1\ast}(\psi_{i_1}(\tilde x))\in
\phi^\ast_I\circ\lambda_I^{1\ast}(\psi_{i_1}(\widetilde
O_{i_1})).\label{e:2.40}
\end{eqnarray}
Here $\phi^\ast_I=(\phi^\ast_l)_{l\in I}\in\Gamma(\psi_I(\tilde
u_I))_1 $, $\phi^\ast_l= \psi_l\circ\phi_l\circ\psi_l^{-1}$ for
$l\in I$,  and $\lambda^{1\ast}_I=(\lambda^\ast_{i_1l})_{l\in I}$,
$\lambda^\ast_{i_1l}=\psi_l\circ\lambda_{i_1l}\circ\psi_{i_1}^{-1}$
for any $l\in I$. Let us define
\begin{equation}\label{e:2.41}
\psi_I\cdot (\tilde u_I,\phi_I, \tilde y_I)=\bigl(\psi_I(\tilde
u_I), \phi^\ast_I, \psi_I(\tilde y_I)\bigr).
\end{equation}
Then it is easily checked that (\ref{e:2.41}) gives a {\bf smooth
effective action} of $\Gamma_I$ on the space $\boxed{\widetilde
W_I\!}\,$. Observe that for each $s\in\{1,\cdots,k\}$ and
$\phi_I\in\Gamma(\tilde u_I)_s$ the above $\phi^\ast_I$ belongs to
$\Gamma(\psi_I(\tilde u_I))_s$. Denote by
$\lambda^{s\ast}_I=(\lambda^\ast_{i_sl})_{l\in I}$,
$\lambda^\ast_{i_sl}=\psi_l\circ\lambda_{i_sl}\circ\psi_{i_s}^{-1}$
for any $l\in I$. Then
$$
\phi^\ast_I\circ\lambda_I^{s\ast}:\psi_{i_s}(\widetilde
O_{i_s})\to \widetilde W_I,\; \tilde z_s\mapsto
\bigl(\phi_l^\ast\circ\lambda^\ast_{i_sl}(\tilde z_s)\bigr)_{l\in
I}
$$
is a smooth open embedding, and
$\{\phi_I^\ast\circ\lambda_I^{s\ast}\,|\,\phi_I\in\Gamma(\tilde
u_I)_s\}$ is a local coordinate chart of $\widetilde W_I$ over the
neighborhood $\psi_I\bigl(\widetilde O(\tilde u_I)\bigr)$ of
$\psi_I(\tilde u_I)$ in the sense of Definition~\ref{def:2.21}. In
particular, we get a local coordinate chart of $\widehat W_I$ around
$[\psi_I(\tilde u_I),\phi^\ast_I,\psi_I(\tilde u_I)]$,
\begin{equation}\label{e:2.42}
\widehat{\phi^\ast_I\circ\lambda_I^{\ast 1}}:\psi_{i_1}(\widetilde
O_{i_1})\to \widehat W_I,\;\tilde x\mapsto [\psi_I(\tilde
u_I),\phi^\ast_I, \phi^\ast_I\circ\lambda_I^{\ast 1}(\tilde x)].
\end{equation}
Let $(\tilde u_I,\phi_I, \tilde y_I)\sim(\tilde u'_I,\phi'_I,
 \tilde y'_I)$ be as in  (\ref{e:2.35}). Clearly,
 $\psi_I(\tilde y_I)=\psi_I(\tilde y_I^\prime)$ and
 $$
 \phi^\ast_I\circ\lambda_I^{1\ast}\bigl(\psi_{i_1}(\widetilde O_{i_1})\bigr)
\bigcap\phi'^\ast_I\circ\lambda_I^{\prime 1\ast}\bigl(\psi_{i_1}(
\widetilde O'_{i_1})\bigr)
$$
is a local component of $\widetilde W_I$ near $\psi_I(\tilde
y_I)=\psi_I(\tilde y'_I)$. So we get:

\begin{lemma}\label{lem:2.24} The action in (\ref{e:2.41})
preserves the equivalence relation $\sim$, i.e. $(\tilde
u_I,\phi_I, \tilde y_I)\sim(\tilde u'_I,\phi'_I,
 \tilde y'_I)$ if and only if $\psi_I\cdot (\tilde u_I,\phi_I, \tilde y_I)\sim\psi_I\cdot(\tilde
u'_I,\phi'_I,  \tilde y'_I)$. So it induces a natural smooth
effective action of $\Gamma_I$ on $\widehat W_I$:
\begin{equation}\label{e:2.43}
\psi_I\cdot[\tilde u_I,\phi_I, \tilde y_I]=[\psi_I\cdot(\tilde
u_I,\phi_I, \tilde y_I)]
\end{equation}
for any $\psi_I\in\Gamma_I$ and $[\tilde u_I,\phi_I, \tilde
y_I]\in\widehat W_I$.
\end{lemma}

 \noindent{\bf Proof.} We only need to prove that the action is  smooth and
effective. Indeed, in the charts in (\ref{e:2.37}) and
(\ref{e:2.42})
$$[\tilde u_I,\phi_I, \tilde y_I]=[\tilde u_I,\phi_I, \phi_I\circ\lambda_I^1(\tilde
x)]\mapsto\psi_I\cdot[\tilde u_I,\phi_I, \tilde y_I]$$
  is given by the smooth map
$$
\widetilde O_{i_1}\to\psi_{i_1}(\widetilde O_{i_1}),\quad\tilde
x\mapsto \psi_{i_1}(\tilde x).
$$
The effectiveness follows from the fact that the action of
$\Gamma_I$ on $\boxed{\widetilde W_I\!}$ is effective.
\hfill$\Box$\vspace{2mm}

Clearly,   a continuous surjective map
\begin{equation}\label{e:2.44}
 q_I: \widehat W_I\to \widetilde W_I,\;[\tilde u_I,\phi_I, \tilde
 y_I]\mapsto \tilde y_I,
\end{equation}
  is smooth as a map from $\widehat W_I$ to
$\prod_{i\in I}\widetilde W_i$ and also satisfies  $\pi_I\circ
q_I=\hat\pi_I$. It follows from (\ref{e:2.41}) and (\ref{e:2.43})
that $q_I$ commutates with the actions on $\widehat W_I$ and
$\widetilde W_I$, i.e. $q_I(\psi_I\cdot[\tilde u_I,\phi_I, \tilde
y_I])=\psi_I\cdot (q_I([\tilde u_I,\phi_I, \tilde  y_I]))$ for
each $[\tilde u_I,\phi_I, \tilde  y_I]\in \widehat W_I$ and
$\psi_I\in\Gamma_I$. These imply that $q_I$ induces a continuous
surjective map $\hat q_I$ from $\widehat
 W_I/\Gamma_I$ to $\widetilde W_I/\Gamma_I$. In fact,
$\hat q_I$ {\bf is a homeomorphism, and hence $\hat\pi_I$ induces
a homeomorphism from $\widehat W_I/\Gamma_I$ to $W_I$}. To see
this let $[\tilde u_I,\phi_I, \tilde  y_I]\in \widehat W_I,\,
[\tilde u'_I,\phi'_I, \tilde  y'_I]\in \widehat W_I$ be such that
$\pi_I(\tilde y_I)=\pi_I(\tilde y'_I)$. We need to prove that
there exists $\psi_I\in\Gamma_I$ such that $\psi_I\cdot[\tilde
u_I,\phi_I, \tilde  y_I]=[\tilde u'_I,\phi'_I, \tilde  y'_I]$.
Since $\pi_I(\tilde y_I)=\pi_I(\tilde y'_I)$ there is
$\psi_I\in\Gamma_I$ such that $\psi_I(\tilde y_I)=\tilde y'_I$.
Recall that $\tilde y_I\in\phi_I\circ\lambda_I^1(\widetilde
O_{i_1})$  and $\tilde
y_I^\prime\in\phi_I^\prime\circ\lambda_I^{\prime 1}(\widetilde
O'_{i_1})$.
 By (\ref{e:2.40}), $\psi_I(\tilde y_I)\in\phi_I^\ast\circ\lambda_I^{1\ast}(\psi_{i_1}(\widetilde
 O_{i_1}))$. That is, $\psi_I\cdot[\tilde
u_I,\phi_I, \tilde  y_I]=[\psi_I(\tilde u_I), \phi_I^\ast,
\psi_I(\tilde y_I)]$. Note that the local structure of $\widetilde
W_I$ near $\tilde y_I$ is essentially the same as that of
$\widetilde W_I$ near $\tilde y'_I$. (By composing with a suitable
element $\varphi_I\in (\Gamma_I)_{\tilde y'_I}$) $\psi_I$ can be
chosen so that the intersection $\phi'_I\circ\lambda_I^{\prime
1}(\widetilde O'_{i_1})\cap
\phi_I^\ast\circ\lambda_I^{1\ast}(\psi_{i_1}(\widetilde
 O_{i_1}))$ is a local component of $\widetilde W_I$ near
 $\tilde y'_I=\psi_I(\tilde y_I)$. Let
\begin{equation}\label{e:2.45}
\left.\begin{array}{ll}
 \widehat W_I^{\circ}:=\{\hat u_I\in\widehat
W_I\,|\,\Gamma_I(\hat
u_I)=\{1\}\}\quad{\rm and}\\
\widehat W_I^{sing}:=\{\hat u_I\in\widehat W_I\,|\,\Gamma_I(\hat
u_I)\ne\{1\}\}
\end{array}
\right\}
\end{equation}
where $\Gamma_I(\hat u_I)$ is the isotropy subgroup of $\Gamma_I$ at
$\hat u_I\in \widehat W_I$. Then $\widehat W_I^{\circ}$ is an open
and dense subset in $\widehat W_I$, and $\widehat W_I^{sing}$ is
relatively closed in $\widehat W_I$. Moreover, $\hat\pi_I(\widehat
W_I^\circ)$ and $\hat\pi_I(\widehat W_I^{sing})$ are equal to
$W_I^\circ$ in (\ref{e:2.12}) and $W_I^{sing}$ in (\ref{e:2.13})
respectively. Summarizing the above arguments we get:

\begin{proposition}\label{prop:2.25}
\begin{description}
\item[(i)]  $\widehat
 W_I$ is a Banach manifold and $\hat\pi^I_J$
is a smooth map.

\item[(ii)] There exists a smooth effective action of $\Gamma_I$
on $\widehat  W_I$ (given by (\ref{e:2.43})) such that $\hat\pi_I$
is invariant under the action and induces a homeomorphism from
$\widehat W_I/\Gamma_I$ to $W_I$.\footnote{This homeomorphism fact
is not used in the arguments of this paper; we only need it to be
a proper, continuous surjective map, which is obvious.}

\item[(iii)]
$\hat\pi_J\circ\hat\pi_J^I=\iota^W_{IJ}\circ\hat\pi_I$ for any
$J\subset I\subset{\cal N}$ and the inclusion $\iota^W_{IJ}:
W_I\hookrightarrow W_J$.

 \item[(iv)] $\hat\pi^I_J\circ\psi_I=\lambda^I_J(\psi_I)\circ\hat\pi^I_J$ for
any $J\subset I\subset{\cal N}$ and  $\psi_I\in\Gamma_I$.

\item[(v)] For any $L, J, I\in{\cal N}$ with $L\subset J\subset
 I$ it holds that $\hat\pi^J_L\circ\hat\pi^I_J=\hat\pi^I_L$.

\item[(vi)]  For any $I=\{i_1,\cdots, i_k\}\in{\cal N}$  the
isotropy group of $\Gamma_I$ at $\hat u_I\in\widehat W_I$ is given
by
$$
\Gamma_I({\hat u_I})=\Gamma_I({\tilde u_I})=\Gamma(\tilde
u_{i_1})\times\cdots\times\Gamma(\tilde u_{i_k}),
$$
where $\tilde u_I=q_I(\hat u_I)$ is given by (\ref{e:2.44}) and
$\Gamma(\tilde u_{i_l})$ is  the isotropy group of $\Gamma_{i_l}$
at $\tilde u_{i_l}$ as before. For $l=1,\cdots, |I|$,
$|(\Gamma_I)_{\hat u_I}|=(|\Gamma(\tilde u_{i_l})|)^{|I|}$ because
these $\Gamma(\tilde u_{i_l})$ are isomorphic. Moreover,  for
$u_I\in W_I$ the inverse image $(\hat \pi_I)^{-1}(u_I)$ contains
$$|\Gamma_I|-|\Gamma_I({\hat
u_I})|+1=\prod^{|I|}_{l}(|\Gamma_{i_l}|-|\Gamma(\tilde
u_{i_l})|+1)
$$
elements exactly.

\item[(vii)]  For any $J\subset I\in{\cal N}$ and $\hat u_J\in{\rm
Im}(\hat\pi^I_J)$ the inverse image $(\hat\pi^I_J)^{-1}(\hat u_J)$
exactly contains
$$\frac{|\Gamma_I|-|\Gamma_I({\hat
u_I})|+1}{|\Gamma_J|-|\Gamma_J({\hat u_J})|+1}=\prod_{i\in
I\setminus J}(|\Gamma_i|-|\Gamma_J({\hat u_J})|^{1/|J|}+1)$$
elements, where $\hat u_I\in(\hat\pi^I_J)^{-1}(\hat u_J)$. In
particular, if $\hat u_J\in{\rm Im}(\hat\pi^I_J)\cap\widehat
W_J^\circ$ then $(\hat\pi^I_J)^{-1}(\hat u_J)$ contains
$|\Gamma_I|/|\Gamma_J|$ elements exactly.

\item[(viii)] For any $J\subset I\in{\cal N}$ it holds that
$(\hat\pi^I_J)^{-1}\bigl({\rm Im}(\hat\pi^I_J)\cap\widehat
W_J^\circ\bigr)=\widehat W_I^\circ$ and
$(\hat\pi^I_J)^{-1}\bigl({\rm Im}(\hat\pi^I_J)\cap\widehat
W_J^{sing}\bigr)=\widehat W_I^{sing}$. Moreover the restriction of
$\hat\pi^I_J$ to $\widehat W_I^\circ$ is a
$|\Gamma_I|/|\Gamma_J|$-fold (regular) covering to ${\rm
Im}(\hat\pi^I_J)\cap\widehat W_J^\circ$.

\item[(ix)] For any $J\subset I\in{\cal N}$ and $\hat u_I\in
\widehat W_I$, $\lambda^I_J\bigl(\Gamma_I({\hat
u_I})\bigr)=\Gamma_J({\hat u_J})$, where $\hat
u_J=\hat\pi^I_J(\hat u_I)$.
\end{description}
\end{proposition}

 (v) may easily follow from (\ref{e:2.38}),  and
 the final (vi)-(ix) are also easily proved by the
above construction. But (vii) is not derived from (ix).

\begin{remark}\label{rem:2.26}
{\rm By (\ref{e:2.25}) it is not hard to check that $\widehat W_I$
is intrinsic in the following sense. Assume that the space
$\boxed{\widetilde W_I\!}$ in (\ref{e:2.34}) is replaced by
$$
\coprod_{\tilde u_I\in\widetilde W_I}\coprod_{O\in{\mathcal
U}(u_I)}\coprod_{1\le s\le |I|}\coprod_{\phi_I\in\Gamma(\tilde
u_I)_s}\phi_I\circ\lambda_I^s(\widetilde O_{i_s}(\tilde u_{i_s}))
$$
 and that the equivalence relation in it is defined by
$$
(\tilde u_I,\phi_I, \tilde y_I)\sim(\tilde u'_I,\phi'_I,
 \tilde y'_I)
$$
if and only if the following two conditions are satisfied:
\begin{description}
\item[(i)] $\tilde y_I=\tilde y'_I$ as points in $ \widetilde
W_I$;

 \item[(ii)] $\phi_I\circ\lambda_I^s(\widetilde O_{i_s})
\bigcap\phi'_I\circ\lambda_I^{\prime t}(\widetilde O'_{i_t})$ is a
local component of $\widetilde W_I$ near $\tilde y_I=\tilde y'_I$.
\end{description}
Here $\tilde u_I,\tilde u'_I\in\widetilde W_I$,
  $\phi_I\in\Gamma(\tilde u_I)_s$, and $\phi'_I\in\Gamma(\tilde u'_I)_t$,
 and $\tilde y_I\in\phi_I\circ\lambda_I^s(\widetilde O_{i_s})$,
  $\tilde y_I^\prime\in\phi_I^\prime\circ\lambda_I^{t\prime}(\widetilde
  O'_{i_t})$, and ${\mathcal U}(u_I)$ is the germ of small
connected open neighborhood of $u_I$ in $W_I$, and $\widetilde
O_{i_s}(\tilde u_{i_s})=\pi^{-1}_{i_s}(O)$, $s=1,\cdots,|I|$.
 Then we get the same $\widehat W_I$, $\hat\pi_I$ and
$\hat\pi^I_J$.  In particular $\widehat
 W_I$ only depends on $\widetilde W_I$.}
 \end{remark}

  Similarly, we can define a desingularization  $\widehat E_I$ of
  $\widetilde E_I$ and the bundle projection $\hat p_I:\widehat E_I
  \to\widehat W_I$ so that it
   is a Banach vector  bundle. For later convenience
we shall give necessary details. Corresponding with (\ref{e:2.34})
let us consider the Banach manifold
\begin{equation}\label{e:2.46}
\boxed{\widetilde E_I\!}:= \coprod_{\tilde u_I\in\widetilde
W_I}\coprod_{\phi_I\in\Gamma(\tilde u_I)_1}\coprod_{O\in{\mathcal
U}(u_I)}\{(\tilde
u_I,\phi_I)\}\times\Phi_I\circ\Lambda_I^1(\widetilde
E_{i_1}|_{\widetilde O_{i_1}}),
\end{equation}
which is  a Banach bundle over $\boxed{\widetilde W_I\!}$ clearly.
Here ${\mathcal U}(u_I)$ and $\widetilde
O_{i_1}=\pi^{-1}_{i_1}(O)$ are as in (\ref{e:2.34}), and $\Phi_I$
and $\Lambda_I^s$ as in (\ref{e:2.28}). Let us define a relation
$\overset{e}{\sim}$ in $\boxed{\widetilde E_I\!}$ as follows:
 For $\tilde u_I,\tilde u'_I\in\widetilde W_I$,
 and $\phi_I\in\Gamma(\tilde u_I)_1$, $\phi'_I\in\Gamma(\tilde u'_I)_1$,
 and  $\tilde\xi_I\in\Phi_I\circ\Lambda_I^1(\widetilde E_{i_1}|_{\widetilde O_{i_1}})$
 and $\tilde\xi_I^\prime\in\Phi_I^\prime\circ\Lambda_I^{\prime 1}(\widetilde E_{i_1}|_{\widetilde O'_{i_1}})$
  we define
\begin{equation}\label{e:2.47}
(\tilde u_I,\phi_I, \tilde \xi_I)\overset{e}{\sim}(\tilde
u'_I,\phi'_I,
 \tilde\xi'_I)
 \end{equation}
if and only the following two conditions hold:
\begin{description}
\item[($1^\circ$)] $\tilde \xi_I=\tilde \xi'_I$ as points in
$\widetilde W_I$;

\item[($2^\circ$)] $\Phi_I\circ\Lambda_I^1(\widetilde
E_{i_1}|_{\widetilde O_{i_1}})
\bigcap\Phi'_I\circ\Lambda_I^{\prime 1}(\widetilde
E_{i_1}|_{\widetilde O'_{i_1}})$ is a local component of
$\widetilde E_I$ near $\tilde p_I(\tilde \xi_I)=\tilde p_I(\tilde
\xi'_I)$, where $\tilde p_I:\widetilde E_I\to\widetilde W_I$ is
the obvious projection.
\end{description}
As in Claim~\ref{cl:2.23} we can prove that $\overset{e}{\sim}$ is a
regular equivalence.
 It is easy to see that
the relation $\overset{e}{\sim}$  is compatible with $\sim$ in
(\ref{e:2.35}). That is,
\begin{eqnarray*}
&&(\tilde u_I,\phi_I, \tilde \xi_I)\overset{e}{\sim}(\tilde
u'_I,\phi'_I,
 \tilde\xi'_I)\Longrightarrow (\tilde u_I,\phi_I,
\tilde p_I(\tilde\xi_I))\sim(\tilde u'_I, \phi'_I,
 \tilde p_I(\tilde\xi'_I))\quad{\rm and}\\
&&(\tilde u_I,\phi_I, \tilde y_I)\sim(\tilde u'_I, \phi'_I,
 \tilde y'_I)\Longrightarrow
 (\tilde u_I,\phi_I,
\tilde o_I(\tilde y_I))\overset{e}{\sim}(\tilde u'_I,\phi'_I,
 \tilde o_I(\tilde y'_I)).
 \end{eqnarray*}
 Here $\tilde o_I:\widetilde W_I\to \widetilde
 E_I$ is the zero section, which is well-defined though
 $\tilde p_I:\widetilde E_I\to\widetilde W_I$ is not a bundle projection.

Let $\langle\tilde u_I,\phi_I, \tilde \xi_I\rangle$ denote the
equivalence class of $(\tilde u_I,\phi_I, \tilde \xi_I)$ with
respect to $\overset{e}{\sim}$. Denote by
\begin{equation}\label{e:2.48}
\widehat E_I=\boxed{\widetilde E_I\!}\bigg/\overset{e}{\sim}
\end{equation}
 if $|I|>1$, and by $\widehat E_I=\widetilde E_i$ if
 $I=\{i\}\in{\mathcal N}$. Then
\begin{equation}\label{e:2.49}
 \hat p_I: \widehat
 E_I\to\widehat
 W_I,\;\langle\tilde u_I,\phi_I, \tilde
 \xi_I\rangle\mapsto [\tilde u_I,\phi_I, \tilde p_I(\tilde \xi_I)]
\end{equation}
  is a Banach bundle, and each
\begin{equation}\label{e:2.50}
\widehat{\Phi_I\circ\Lambda_I^1}:\widetilde E_{i_1}|_{\widetilde
O_{i_1}}\to \widehat E_I,\;\tilde \xi\mapsto \langle\tilde
u_I,\phi_I, \Phi_I\circ\Lambda_I^1(\tilde \xi)\rangle,
\end{equation}
gives a smooth bundle open embedding, called {\bf a local bundle
coordinate chart of $\widehat E_I$ near $\langle\tilde u_I,\phi_I,
\tilde o_I(\tilde u_I)\rangle$}. The projections $\Pi_I$ and
$\Pi^I_J$ induce natural ones
\begin{equation}\label{e:2.51}
\left.\begin{array}{ll}
 \hat\Pi_I: \widehat E_I\to
E_I,\;\langle\tilde
u_I,\phi_I, \tilde \xi_I\rangle\mapsto \Pi_I(\tilde \xi_I),\\
 \hat\Pi^I_J:\widehat E_I\to \widehat E_J,\;\langle\tilde u_I,\phi_I, \tilde
\xi_I\rangle\mapsto \langle\tilde u_J,\phi_J, \tilde \xi_J\rangle
\end{array}\right\}
\end{equation}
 which are well-defined because of
Proposition~\ref{prop:2.22}(iii). Here
$\phi_J=\lambda^I_J(\phi_I)=(\phi_l)_{l\in J}$, $\tilde
u_J=\pi^I_J(\tilde u_I)=(\tilde u_l)_{l\in J}$ and $\tilde \xi_J=
\Pi^I_J(\tilde \xi_I)=(\tilde \xi_l)_{l\in J}$. Moreover, in the
bundle charts $\widehat{\Phi_I\circ\Lambda_I^1}$ and
$$
 \widehat{\Phi_J\circ\Lambda_J^1}:
\widetilde E_{j_1}|_{\widetilde O_{j_1}}\to \widehat
E_J,\;\tilde\eta\mapsto \langle\tilde u_J,\phi_J,
\Phi_J\circ\Lambda_J^1(\tilde \eta)\rangle,
 $$
 the projection $\hat\Pi^I_J$ may be represented by
\begin{equation}\label{e:2.52}
 \widetilde E_{i_1}|_{\widetilde O_{i_1}}\to \widetilde E_{j_1}|_{\widetilde O_{j_1}},
 \; \tilde \xi\mapsto
\Lambda_{i_1j_1}(\tilde \xi)
\end{equation}
because $j_1\in J\subset I$. Here $\Lambda_{i_1j_1}=1_{\widetilde
E_{i_1}|_{\widetilde O_{i_1}}}$ if $j_1=i_1$. It follows that
$\hat\Pi^I_J$ is a {\bf smooth bundle map} and also {\bf
isomorphically maps} the fibre of $\widehat E_I$ at $[\tilde u_I,
\phi_I,\tilde p_I(\tilde\xi_I)]\in\widehat W_I$ to that of
$\widehat E_J$ at $[\tilde u_J, \phi_J,\tilde
p_J(\tilde\xi_J)]\in\widehat W_J$.

In order to see how the action of $\Gamma_I$ on $\widetilde E_I$
induces a natural one on $\boxed{\widetilde E_I\!}\,$, note that
$\psi_I=(\psi_l)_{l\in I}\in\Gamma_I$ acts on $\widetilde E_I$ via
\begin{equation}\label{e:2.53}
\psi_I\cdot\tilde\xi_I=\Psi_I(\tilde\xi_I)=(\Psi_l(\tilde\xi_l))_{l\in
I}.
\end{equation}
So for
\begin{eqnarray*}
\tilde \xi_I=\Phi_I\circ\Lambda_I^1(\tilde
\eta)\!\!\!\!\!\!\!\!&&= \bigl(\tilde \eta,
\Phi_{i_2}\circ\Lambda_{i_1i_2}(\tilde \eta),\cdots,
\Phi_{i_k}\circ\Lambda_{i_1i_k}(\tilde \eta)\bigr)\\
&&\in \Phi_I\circ\Lambda_I^1(\widetilde E_{i_1}|_{\widetilde
O_{i_1}})\subset\widetilde E_I
\end{eqnarray*}
it holds that
\begin{eqnarray*}
\psi_I\cdot\tilde
\xi_I\!\!\!\!\!\!\!\!\!\!&&=\bigl(\Psi_{i_1}(\tilde \eta),
\Psi_{i_2}\circ\Phi_{i_2}\circ\Lambda_{i_1i_2}(\tilde
\eta),\cdots,
\Psi_{i_k}\circ\Phi_{i_k}\circ\Lambda_{i_1i_k}(\tilde
\eta)\bigr)\\
&&=\bigl(\Psi_{i_1}(\tilde\eta),
\Phi^\ast_{i_2}\circ\Lambda^\ast_{i_1i_2}(\Psi_{i_1}(\tilde
\eta)),\cdots,
\Phi^\ast_{i_k}\circ\Lambda^\ast_{i_1i_k}(\Psi_{i_1}(\tilde
\eta))\bigr)\nonumber\\
&&=\Phi^\ast_I\circ\Lambda_I^{1\ast}(\Psi_{i_1}(\tilde \eta))\in
\Phi^\ast_I\circ\Lambda_I^{1\ast}\bigr(\Psi_{i_1}(\widetilde
E_{i_1}|_{\widetilde O_{i_1}})\bigl),\nonumber
\end{eqnarray*}
where $\Phi^\ast_I=(\Phi^\ast_l)_{l\in I}=
(\Psi_l\circ\Phi_l\circ\Psi_l^{-1})_{l\in I}$,  and
$\Lambda^{1\ast}_I=(\Lambda^\ast_{i_1l})_{l\in I}$,
$\Lambda^\ast_{i_1l}=\Psi_l\circ\Lambda_{i_1l}\circ\Psi_{i_1}^{-1}$
for any $l\in I$. These motivate us to define
\begin{equation}\label{e:2.54}
\psi_I\cdot (\tilde u_I,\phi_I, \tilde \xi_I)=\bigl(\psi_I(\tilde
u_I), \phi^\ast_I, \psi_I\cdot\tilde \xi_I\bigr).
\end{equation}
Here $\phi^\ast_I=(\phi^\ast_l)_{l\in
I}=(\psi_l\circ\phi_l\circ\psi_l^{-1})_{l\in I}
\in\Gamma(\psi_I(\tilde u_I))_1 $ as in (\ref{e:2.40}). It gives a
smooth effective action of $\Gamma_I$ on $\boxed{\widetilde
E_I\!}\,$. As above $\phi^\ast_I$ belongs to $\Gamma(\psi_I(\tilde
u_I))_s$ as $\phi_I\in\Gamma(\tilde u_I)_s$. And  for each
$s\in\{1,\cdots,k\}$, if we define
$\Lambda^{s\ast}_I=(\Lambda^\ast_{i_sl})_{l\in I}$,
$\Lambda^\ast_{i_sl}=\Psi_l\circ\Lambda_{i_sl}\circ\Psi_{i_s}^{-1}$
for any $l\in I$, then
\begin{equation}\label{e:2.55}
\Phi^\ast_I\circ\Lambda_I^{s\ast}:\Psi_{i_s}(\widetilde
E_{i_s}|_{\widetilde O_{i_s}})\to \widetilde E_I,\; \tilde
\xi_s\mapsto \bigl(\Phi_l^\ast\circ\Lambda^\ast_{i_sl}(\tilde
\xi_s)\bigr)_{l\in I}
\end{equation}
is a smooth bundle open embedding, and
$\{\Phi_I^\ast\circ\Lambda_I^{s\ast}\,|\,\phi_I\in\Gamma(\tilde
u_I)_s\}$ is a {\bf local bundle coordinate chart} of $\widetilde
E_I$ over the neighborhood $\psi_I\bigl(\widetilde O(\tilde
u_I)\bigr)$ of $\psi_I(\tilde u_I)$.  In particular, we get a
local bundle coordinate chart of $\widehat E_I$ around
$\langle\psi_I(\tilde u_I),\phi^\ast_I,\tilde o_I(\psi_I(\tilde
u_I))\rangle$,
\begin{equation}\label{e:2.56}
\widehat{\Phi^\ast_I\circ\Lambda_I^{\ast 1}}:\widetilde
E_{i_1}|_{\psi_{i_1}(\widetilde O_{i_1})}\to \widehat E_I,\;\tilde
\xi\mapsto \langle\psi_I(\tilde u_I),\phi^\ast_I,
\Phi^\ast_I\circ\Lambda_I^{\ast 1}(\tilde \xi)\rangle.
\end{equation}
Let $(\tilde u_I,\phi_I, \tilde\xi_I)\overset{e}{\sim}(\tilde
u'_I,\phi'_I,
 \tilde\xi'_I)$ be as in  (\ref{e:2.47}). Then
 $\Psi_I(\tilde\xi_I)=\Psi_I(\tilde\xi_I^\prime)$ and
 $$
 \Phi^\ast_I\circ\Lambda_I^{1\ast}\bigl(\widetilde E_{i_1}|_{\psi_{i_1}(\widetilde O_{i_1})}\bigr)
\bigcap\Phi'^\ast_I\circ\Lambda_I^{\prime 1\ast}\bigl(\widetilde
E_{i_1}|_{\psi_{i_1}( \widetilde O'_{i_1}})\bigr)
$$
is a local component of $\widetilde E_I$ near $\tilde
p_I(\Psi_I(\tilde \xi_I))=\tilde p_I(\Psi_I(\tilde\xi'_I))$. So
the action in (\ref{e:2.54}) preserves the equivalence relation
$\overset{e}{\sim}$ in $\boxed{\widetilde E_I\!}\,$, i.e.
$$
 (\tilde u_I,\phi_I, \tilde \xi_I)\overset{e}{\sim}(\tilde
u'_I,\phi'_I,
 \tilde\xi'_I)\Longleftrightarrow
\psi_I\cdot (\tilde u_I,\phi_I, \tilde
\xi_I)\overset{e}{\sim}\psi_I\cdot(\tilde u'_I,\phi'_I,
\tilde\xi'_I).
$$
It follows that (\ref{e:2.54}) induces a natural smooth effective
action of $\Gamma_I$ on $\widehat E_I$:
\begin{equation}\label{e:2.57}
\psi_I\cdot\langle\tilde u_I,\phi_I, \tilde
\xi_I\rangle=\langle\psi_I\cdot(\tilde u_I,\phi_I,
\tilde\xi_I)\rangle.
\end{equation}
 In fact, in the charts in (\ref{e:2.50})
and (\ref{e:2.56}) the action
$$\langle\tilde u_I,\phi_I, \tilde\xi_I\rangle=\langle\tilde u_I,\phi_I,
\Phi_I\circ\Lambda_I^1(\tilde\eta)\rangle\mapsto\psi_I\cdot\langle\tilde
u_I,\phi_I, \tilde\xi_I\rangle$$
 is given by the smooth map
$$
 (\widetilde E_{i_1}|_{\widetilde
O_{i_1}}\to\widetilde E_{i_1}|_{\psi_{i_1}(\widetilde
O_{i_1})})\quad\tilde \eta\mapsto \Psi_{i_1}(\tilde\eta).
$$
The effectiveness is also easily proved as before.
 Moreover, the clear continuous surjective map
$$
 Q_I: \widehat E_I\to \widetilde E_I,\;\langle\tilde u_I,\phi_I, \tilde
 \xi_I\rangle\mapsto \tilde\xi_I
$$
 commutates with the above actions of $\Gamma_I$ on
$\widehat E_I$ and $\widetilde E_I$ and also satisfies $\Pi_I\circ
Q_I=\hat\Pi_I$. As in the arguments above (\ref{e:2.45}) we can
derive that the projection $\hat\Pi_I$ is invariant under the
 $\Gamma_I$-action and induces a homeomorphism
$$
  \widehat E_I/\Gamma_I\to  E_I,\;\langle\tilde u_I,\phi_I, \tilde
 \xi_I\rangle\mapsto \Pi_I(\tilde\xi_I).
$$
 These arguments yield:

\begin{proposition}\label{prop:2.27}
\begin{description}
\item[(i)] Each $\hat p_I:\widehat
 E_I\to\widehat
 W_I$ is a  Banach bundle.

\item[(ii)] There exists a smooth effective action of $\Gamma_I$
on $\widehat
 E_I$ (given by
(\ref{e:2.57})), and $\hat\Pi_I$ is invariant under the action and
induces a homeomorphism from $\widehat
 E_I/\Gamma_I$ to $E_I$.

\item[(iii)] The projection $\hat p_I:\widehat E_I\to\widehat W_I$
is equivariant with respect to the actions of $\Gamma_I$ on them,
i.e. $\hat p_I\circ\psi_I=\psi_I\circ\hat p_I$ for any
$\psi_I\in\Gamma_I$.

 \item[(iv)]  $\hat\Pi_J\circ\hat\Pi_J^I=\iota^E_{IJ}\circ\hat\Pi_I$
   for any $J\subset I\in{\cal N}$ and the inclusion $\iota^E_{IJ}:
E_I\hookrightarrow E_J$.

 \item[(v)] For any $J\subset I\in{\cal N}$, $\hat\Pi^I_J$
is a smooth bundle map covering $\hat\pi^I_J$ and also
isomorphically maps the fibre of $\widehat E_I$ at $\hat u_I\in
(\pi^I_J)^{-1}(\widehat W_J)$ to that of $\widehat E_J$ at
$\hat\pi^I_J(\hat u_I)$. Moreover, it is $\Gamma_I$-equivariant,
i.e., $\hat\Pi^I_J\circ\psi_I=\lambda^I_J(\psi_I)\circ\hat\Pi^I_J$
for  $\psi_I\in\Gamma_I$.

\item[(vi)] $\hat\pi_I\circ\hat p_I=p\circ\hat\Pi_I$ for any
 $\phi_I\in\Gamma_I$.

\item[(vii)] $\hat p_J\circ\hat\Pi^I_J=\hat\pi^I_J\circ\hat p_I$
for any $J\subset I\in{\cal
 N}$.

\item[(viii)] For any $L, J, I\in{\cal N}$ with $L\subset J\subset
 I$ it holds that $\hat\Pi^J_L\circ\hat\Pi^I_J=\hat\Pi^I_L$.
 \end{description}
\end{proposition}

 (vii) may easily follow from (\ref{e:2.51}).
 We get a system of Banach bundles
\begin{equation}\label{e:2.58}
\bigl(\widehat{\mathcal E}(K), \widehat W(K)\bigr)=
\bigl\{\bigl(\widehat E_I, \widehat W_I\bigr), \hat\pi_I,
\hat\Pi_I, \Gamma_I, \hat\pi^I_J,\hat\Pi^I_J, \lambda^I_J\bigm|
J\subset I\in{\mathcal N}\bigr\},
\end{equation}
which is called a {\bf desingularization} of
$\bigl(\widetilde{\mathcal E}(K), \widetilde W(K)\bigr)$ in
(\ref{e:2.15}).

\begin{definition}\label{def:2.28} {\rm A {\bf global section} of
the bundle system $(\widehat {\mathcal E}(K), \widehat W(K))$ is a
compatible collection $\widehat S=\{\widehat S_I\,|\,
I\in{\mathcal N}\}$ of sections $\widehat S_I$ of $\widehat
E_I\to\widehat W_I$ in the sense that
\begin{equation}\label{e:2.59}
\widehat S_I=(\hat\Pi^I_J)^\ast \widehat
S_J:=(\hat\Pi^I_J)^{-1}\circ \widehat S_J\circ\hat\pi^I_J,
\end{equation}
i.e. $\widehat S_J(\hat\pi^I_J(x))=\hat\Pi^I_J(\widehat S_I(x))$
for any $x\in \widehat W_I$ and $J\subset I\in{\cal N}$.
  $\widehat S$ is said to be {\bf smooth} if each $\widehat S_I$ is. $\widehat S$ is called {\bf
transversal} to the zero section if each section $\widehat
S_I:\widehat W_I\to\widehat E_I$ is so. Furthermore, a global
smooth section $\widehat S=\{\widehat S_I\,|\, I\in{\mathcal N}\}$
is called {\bf Fredholm} if each $\widehat S_I$ is a Fredholm
section. If each section $\widehat S_I$ is $\Gamma_I$-equivariant
the global section $\widehat S=\{\widehat S_I\}$ is said to be
{\bf equivariant}.}
\end{definition}

\begin{lemma}\label{lem:2.29}
 For each $I\in{\cal N}$ with $|I|>1$, and $l\in I$ we denote by $I_l=\{l\}$ for $1\le l\le n$.
  Then each smooth section
  $\tilde\sigma_l:\widetilde W_l\to\widetilde E_l$ can define one
  $\hat\sigma_{lI}:\widehat W_I\to\widehat   E_I$ by
 \begin{equation}\label{e:2.60}
 \widehat W_I\to\widehat E_I,\; \hat x\mapsto
\bigl((\hat\Pi^I_l)^\ast\tilde\sigma_l\bigr)(\hat x)=
   (\hat\Pi^I_l)^{-1}\bigl(\tilde\sigma_l(\hat \pi^I_{I_l}(\hat x))\bigr)
\end{equation}
which is equal to $\tilde\sigma_l$ if $I=I_{l}$. Moreover, if
$\tilde\sigma_l$ is Fredholm then so is $\hat\sigma_{lI}$.
\end{lemma}

\noindent{\bf Proof.}\quad We first prove that the section
$\hat\sigma_{lI}$  is smooth. To this end let us write it in the
charts in (\ref{e:2.37}) and (\ref{e:2.50}). Let $\hat
x=\widehat{\phi_I\circ\lambda_I^1}(\tilde x)=[\tilde u_I,\phi_I,
\phi_I\circ\lambda^1_I(\tilde x)]$ for $\tilde x\in \widetilde
O_{i_1}$. Then $\pi^I_{I_l}(\hat x)=[\tilde
u_l,\phi_l,\phi_l\circ\lambda_{i_1l}(\tilde
x)]=\phi_l\circ\lambda_{i_1l}(\tilde x)$, and specially equals
$\tilde x$ if $l=i_1$. So
$\tilde\sigma_l(\phi_l\circ\lambda_{i_1l}(\tilde x))=\langle\tilde
u_l,\phi_l, \tilde\sigma_l(\phi_l\circ\lambda_{i_1l}(\tilde
x))\rangle$. As in (\ref{e:2.52}), under the bundle charts
$\widehat{\Phi_I\circ\Lambda^1_I}$ in (\ref{e:2.50}) and
\begin{eqnarray*}
&&\widehat{\Phi_{I_l}\circ\Lambda^1_{I_l}}: \widetilde
E_{l}|_{\widetilde O_{l}}\to \widehat E_{I_l},\\
&&\qquad\tilde\eta\mapsto \langle\tilde u_{I_l},\phi_{I_l},
\Phi_{I_l}\circ\Lambda^1_{I_l}(\tilde \eta)\rangle=\langle\tilde
u_l,\phi_l, \tilde \eta\rangle=\tilde\eta,
\end{eqnarray*}
the projection $\hat\Pi^I_{I_l}$ may be represented by
$$
 \widetilde E_{i_1}|_{\widetilde O_{i_1}}\to \widetilde E_{l}|_{\widetilde O_{l}},
 \; \tilde \xi\mapsto
\Phi_l\circ\Lambda_{i_1l}(\tilde \xi).
$$
It follows that
$\hat\sigma_{lI}$ may locally be represented as
\begin{equation}\label{e:2.61}
\widetilde O_{i_1}\to \widetilde E_{i_1}|_{\widetilde O_{i_1}},
 \; \tilde x\mapsto (\Phi_l\circ\Lambda_{i_1l})^{-1}\bigl(\tilde\sigma_l(\phi_l\circ\lambda_{i_1l}(\tilde x))\bigr).
\end{equation}
This shows that $\hat\sigma_{lI}$ is smooth. The other claims
follows from this clearly. \hfill$\Box$\vspace{2mm}

 A class of special sections of the bundle system
$(\widehat {\mathcal E}(K), \widehat W(K))$ may be derived from
those of the Banach orbibundle $E\to X$.

\begin{lemma}\label{lem:2.30}
Each (smooth) section $S$ of the Banach orbibundle $E\to X$ may
yield a global (smooth) section of the bundle system $(\widehat
{\mathcal E}(K), \widehat W(K))$, denoted by $\widehat
S=\{\widehat S_I\,|\, I\in{\mathcal N}\}$. The section $\widehat
S$ is equivariant in the sense that each $\widehat S_I$ is
$\Gamma_I$-equivariant. If $S$ is Fredholm so is $\widehat S$ and
they has the same index.
\end{lemma}

\noindent{\bf Proof.}\quad Let  $\tilde S_i:\widetilde
W_i\to\widetilde E_i$ unique $\Gamma_i$-equivariant lifts of
$S|_{W_i}$, $i=1,\cdots,n$. They are compatible in the sense that
$$
\tilde S_j\circ\lambda_{ij}=\Lambda_{ij}\circ\tilde S_i\quad{\rm
on}\quad\widetilde O_i,\;i, j=1,\cdots,n,
$$
where $\lambda_{ij}$,
$\Lambda_{ij}$ and $\widetilde O_i$ are as above.
 Since $\Pi_i\circ\tilde S_i=S\circ\pi_i$, $i=1,\cdots,n$,
for each $I\in{\cal N}$ we have a map
$$
\tilde S_I: \widetilde W_I\to \widetilde E_I,\;\tilde x_I=(\tilde
x_i)_{i\in I}\mapsto (\tilde S_i(\tilde x_i))_{i\in I}.
$$
 This map determines a (smooth) section of the
Banach bundle $\hat p_I:\widehat W_I\to\widehat E_I$. In fact, let
$I=\{i_1,\cdots, i_k\}$, and $(\tilde
u_I,\phi_I,\phi_I\circ\lambda_I^1(\tilde x))\sim(\tilde
u'_I,\phi'_I,\phi'_I\circ\lambda'^1_I(\tilde x'))$ in
$\boxed{\widehat W_I}$ for $\tilde x\in\widetilde O_{i_1}$ and
$\tilde x'\in\widetilde O'_{i_1}$. Using the local chart
$\phi_I\circ\lambda_I^1$ in the proof of Proposition~\ref{prop:2.22}
we have
\begin{eqnarray*}
\tilde S_I\bigl(\phi_I\circ\lambda_I^1(\tilde
x)\bigr)\!\!\!\!\!\!\!\!\!\!\!&&= \Bigl(\tilde S_{i_1}(\tilde x),
\tilde S_{i_2}\bigl(\phi_{i_2}\circ\lambda_{i_1i_2}(\tilde
x)\bigr),\cdots, \tilde
S_{i_k}\bigl(\phi_{i_k}\circ\lambda_{i_1i_k}(\tilde
x)\bigr)\Bigr)\\
&&= \Bigl(\tilde S_{i_1}(\tilde x), \phi_{i_2}\cdot\tilde
S_{i_2}\bigl(\lambda_{i_1i_2}(\tilde x)\bigr),\cdots,
\phi_{i_k}\cdot\tilde S_{i_k}\bigl(\lambda_{i_1i_k}(\tilde
x)\bigr)\Bigr)\\
&&= \Bigl(\tilde S_{i_1}(\tilde x), \Phi_{i_2}\bigl(\tilde
S_{i_2}\bigl(\lambda_{i_1i_2}(\tilde x)\bigr)\bigr),\cdots,
\Phi_{i_k}\bigl(\tilde S_{i_k}\bigl(\lambda_{i_1i_k}(\tilde
x)\bigr)\bigr)\Bigr)\\
&&= \Bigl(\tilde S_{i_1}(\tilde x),
\Phi_{i_2}\bigl(\Lambda_{i_1i_2}(\tilde S_{i_1}(\tilde
x))\bigr),\cdots, \Phi_{i_k}\bigl(\Lambda_{i_1i_k}(\tilde
S_{i_1}(\tilde
x))\bigr)\Bigr)\\
&&=\Phi_I\circ\Lambda_I^1(\tilde S_{i_1}(\tilde x)),
\end{eqnarray*}
\begin{eqnarray*}
\tilde S_I\bigl(\phi'_I\circ\lambda'^1_I(\tilde
x')\bigr)\!\!\!\!\!\!\!\!\!\!&&= \Bigl(\tilde S_{i_1}(\tilde x'),
\tilde S_{i_2}\bigl(\phi'_{i_2}\circ\lambda'_{i_1i_2}(\tilde
x')\bigr),\cdots, \tilde
S_{i_k}\bigl(\phi'_{i_k}\circ\lambda'_{i_1i_k}(\tilde
x')\bigr)\Bigr)\\
&&= \Bigl(\tilde S_{i_1}(\tilde x'), \phi_{i_2}\cdot\tilde
S_{i_2}\bigl(\lambda'_{i_1i_2}(\tilde x')\bigr),\cdots,
\phi'_{i_k}\cdot\tilde S_{i_k}\bigl(\lambda'_{i_1i_k}(\tilde
x')\bigr)\Bigr)\\
&&= \Bigl(\tilde S_{i_1}(\tilde x'), \Phi'_{i_2}\bigl(\tilde
S_{i_2}\bigl(\lambda'_{i_1i_2}(\tilde x')\bigr)\bigr),\cdots,
\Phi'_{i_k}\bigl(\tilde S_{i_k}\bigl(\lambda'_{i_1i_k}(\tilde
x')\bigr)\bigr)\Bigr)\\
&&= \Bigl(\tilde S_{i_1}(\tilde x'),
\Phi'_{i_2}\bigl(\Lambda'_{i_1i_2}(\tilde S_{i_1}(\tilde
x'))\bigr),\cdots, \Phi'_{i_k}\bigl(\Lambda'_{i_1i_k}(\tilde
S_{i_1}(\tilde
x'))\bigr)\Bigr)\\
&&=\Phi'_I\circ\Lambda'^1_I(\tilde S_{i_1}(\tilde x')).
\end{eqnarray*}
Since $\phi_I\circ\lambda_I^1(\tilde
x)=\phi'_I\circ\lambda'^1_I(\tilde x')$($\Longrightarrow \tilde
x=\tilde x'$) we get that
\begin{eqnarray*}
\tilde S_I\bigl(\phi_I\circ\lambda_I^1(\tilde
x)\bigr)\!\!\!\!\!\!\!\!\!&&=\Phi_I\circ\Lambda_I^1(\tilde S_{i_1}(\tilde x))\\
&&=\Phi'_I\circ\Lambda'^1_I(\tilde S_{i_1}(\tilde x'))=\tilde
S_I\bigl(\phi'_I\circ\lambda'^1_I(\tilde x')\bigr).
\end{eqnarray*}
 Moreover, $\phi_I\circ\lambda_I^1(\widetilde
O_{i_1} )\cap\phi'_I\circ\lambda'^1_I(\widetilde O'_{i_1})$ is a
local component of $\widehat W_I$ near
$\phi_I\circ\lambda_I^1(\tilde x)=\phi'_I\circ\lambda'^1_I(\tilde
x')$. It follows that
$$
\Phi_I\circ\Lambda_I^1(\widetilde E_{i_1}|_{\widetilde O_{i_1}}
)\cap\Phi'_I\circ\Lambda'^1_I(\widetilde E_{i_1}|_{\widetilde
O'_{i_1}})
$$
is a local component of $\widehat E_I$ near
$\phi_I\circ\lambda_I^1(\tilde x)=\phi'_I\circ\lambda'^1_I(\tilde
x')$. This shows that
\begin{eqnarray*}
 &&\Bigl(\tilde
u_I,\tilde\phi_I,\tilde S_I\bigl(\phi_I\circ\lambda^1_I(\tilde
x)\bigr)\Bigr)\overset{e}{\sim} \Bigl(\tilde
u'_I,\tilde\phi'_I,\tilde S_I\bigl(\phi'_I\circ\lambda'^1_I(\tilde
x')\bigr)\Bigr)\quad{\rm or}\\
&&\widehat{\Phi'_I\circ\Lambda'^1_I}(\tilde S_{i_1}(\tilde x'))=
\widehat{\Phi_I\circ\Lambda_I^1}(\tilde S_{i_1}(\tilde x)).
\end{eqnarray*}
So we may define $\widehat S_I:\widehat W_I\to \widehat E_I$ by
\begin{equation}\label{e:2.62}
[\tilde u_I,\phi_I, \tilde x_I]\mapsto \langle\tilde u_I,\phi_I,
\tilde s_I(\tilde x_I)\rangle.
\end{equation}
It is easy to see that in the charts of (\ref{e:2.37}) and
(\ref{e:2.50}) the section may be expressed as
\begin{equation}\label{e:2.63}
\widetilde O_{i_1}\to\widetilde E_{i_1}|_{\widetilde
O_{i_1}},\;\tilde x\mapsto\tilde S_{i_1}(\tilde x).
\end{equation}
Thus $\widehat S_I$ is a {\bf smooth section}, and its zero set is
given by
$$
Z(\widehat S_I)=(\hat \pi_I)^{-1}\bigl(Z(S)\cap W_I\bigr).
$$
 From (\ref{e:2.38}) and (\ref{e:2.51}) one can
immediately derive
\begin{equation}\label{e:2.64}
\widehat S_J\circ\hat\pi^I_J=\hat\Pi^I_J\circ\widehat S_I
\end{equation}
 for any $J\subset I\in{\cal N}$. Finally, one easily checks that each $\widehat S_I$ is
 $\Gamma_I$-equivariant. From (\ref{e:2.63}) it easily follows that
 each $\widehat S_I$ and thus $\widehat S$ are Fredholm if $S$ is so. \hfill$\Box$\vspace{2mm}

As done in Section~\ref{sec:1} our aim is to give a small
perturbation of the Fredholm section $\widehat S$ so that it becomes
transversal to the zero section. The arguments at the beginning of
this section show that such a purpose can only be realized by
perturbations of the nonequivariant global sections. However, even
if $\tilde\sigma_l:\widetilde W_l\to\widetilde E_l$ has compact
support contained in $\widetilde W_l$ one mayt not be able to extend
the collection $\{\hat\sigma_{lI}: I\in{\cal N},\;l\in I\}$ of
sections given by Lemma~\ref{lem:2.29} into a global smooth section
of $(\widehat {\mathcal E}(K), \widehat W(K))$ because the sets
$W_i$ overlap too much to satisfy the compatibility conditions of
the definition. (cf. Ex.4.7 in \cite{Mc1}.) So we need to make some
improvements on it.

 \subsubsection{Renormalization of the
 resolution in (\ref{e:2.58})}\label{sec:2.2.2}

In this subsection we shall  improve the resolution in last
subsection. Note that each $I\in{\mathcal N}$ has length $|I|<n$.
For the above open sets $W_i$, $i=1, 2, \cdots, n$, one may take
 pairs of open sets
 $$
 W^j_i\subset\subset U^j_i,\; j=1, 2, \cdots,n-1,
$$
 such that
$$
 U^1_i\subset\subset
W^2_i\subset\subset U^2_i\cdots \subset\subset
W^{n-1}_i\subset\subset U^{n-1}_i\subset\subset W_i.
$$
Then for each $I\in{\cal N}$ with $|I|=k$ define
\begin{equation}\label{e:2.65}
 V_I:=\Bigl(\bigcap_{i\in I}W_i^k\Bigr)\setminus
\Bigl(\bigcup_{J:|J|>k}\bigl(\bigcap_{j\in
J}Cl(U_j^k)\bigr)\Bigr).
\end{equation}
Then the second condition in (\ref{e:2.8}) implies
$$
V_I=W_I^{n-1}:=\cap_{i\in I}W_i^{n-1}\quad\forall I\in{\cal
N}\;\;{\rm with}\; |I|=n-1.
$$
 For $\{V_I\,|\, I\in{\cal
N}\}$, unlike $\{W_I\,|\, I\in{\cal N}\}$ we cannot guarantee that
$V_I\subset V_J$ even if $J\subset I$. However, Lemma 4.3 in
[LiuT1] showed:

\begin{lemma}\label{lem:2.31}
 $\{V_I\,|\,I\in{\mathcal
N}\}$ is an open covering of $\cup^n_{i=1}W_i^1$  and satisfies:
\begin{description}
\item[(i)] $V_I\subset W_I$  for any $I\in {\mathcal N}$.

\item[(ii)] $Cl(V_I)\cap V_J\ne\emptyset$  only if $I\subset J$ or
$J\subset I$. (Actually  $Cl(V_I)\cap V_J\ne\emptyset$ implies
$V_I\cap V_J\ne\emptyset$, and thus $Cl(V_I)\cap V_J=\emptyset$
for any $I, J\in{\cal N}$ with $|I|=|J|=n-1$.)

\item[(iii)] For any nonempty open subset ${\cal
W}^\ast\subset\cup^n_{i=1}W_i^1$ and $I\in{\cal N}$ let
$W_I^\ast=W_I\cap{\cal W}^\ast$ and $V_I^\ast=V_I\cap{\cal
W}^\ast$. Then $\{V_I^\ast\,|\, I\in{\cal N}\}$ is an open
covering of ${\cal W}^\ast$, and also satisfies the above
corresponding properties (i)-(ii).
\end{description}
\end{lemma}

\noindent{\bf Proof.}\quad {\it Step 1}. Let ${\cal
N}_k=\{I\in{\cal N}\,|\, |I|=k\}$, $k=1,\cdots,n-1$. We first
prove that
 $\{V_I\,|\,I\in{\mathcal N}\}$  is an open covering of
$\cup^n_{i=1}W_i^1$. It is easy to see that for any $I\in{\mathcal
N}$ with $|I|=k$,
$$
V_I=\Bigl(\bigcap_{i\in I}W_i^k\Bigr)\setminus
\Bigl(\bigcup_{J\in{\cal N}_{k+1}}\bigl(\bigcap_{j\in
J}Cl(U_j^k)\bigr)\Bigr).
 $$
 For $x\in W_i$, set $I_1=\{i\}$.
 If $x\in V_{I_1}$ nothing is done. Otherwise, there is
 $J_1\in{\mathcal N}$ with $|J_1|=2$  such that
 $x\in\cap_{j\in J_1}Cl(U^1_j)\subset W^2_{J_1}$. Set $I_2=I_1\cup J_1$ then
 $x\in W^1_{I_1}\cap W^2_{J_1}\subset W^2_{I_1}\cap W^2_{J_1} $ and
 thus
 $$x\in W^2_{I_2}\subset W^{|I_2|}_{I_2},\;  |I_2|\ge |J_1|=|I_1|+1=2.$$
If $x\in V_{I_2}$ then nothing is done. Otherwise, because of
(\ref{e:2.65}) there is $J_2\in{\mathcal N}$ with $|J_2|=|I_2|+1$
such that $x\in \cap_{j\in J_2}Cl(U_j^{|I_2|})\subset
W_{J_2}^{|J_2|}$. Set $I_3=I_2\cup J_2$ then
 $x\in  W^2_{I_2}\cap W_{J_2}^{|J_2|}\subset W^{|I_3|}_{I_2}\cap W_{J_2}^{|I_3|}$ and
 thus
$$
x\in W_{I_3}^{|I_3|},\;  |I_3|\ge |J_2|=|I_2|+1\ge 3.
$$
After repeating finite times this process there must exist some
$I_k\in{\mathcal N}$ such that $x\in V_{I_k}$ because $
V_I=W_I^{n-1}$ for any $I\in{\cal N}$ with $|I|=n-1$.
 The desired conclusion is proved.

{\it Step 2.}  We prove (i) and (ii). (i) is obvious.  We only
need to prove (ii). Let $ Cl(V_I)\cap V_J\ne\emptyset$ for two
different  $I, J\in{\mathcal N}$ with $|I|=k$ and $|J|=l$. Since
$I=\{1,\cdots, n\}\notin{\cal N}$ we can assume that $n>k\ge l$
below. Suppose that there exists a $r\in J\setminus I$. Let $x\in
Cl(V_I)\cap V_J$. Take a sequence $\{x_k\}\subset V_I$ such that
$x_k\to x$. Since $V_J$ is open then $x_k\in V_J$ for $k$
sufficiently large. So  $V_I\cap V_J\ne\emptyset$ and we may
assume that $x\in V_I\cap V_J$. Now one hand $x\in V_J$ implies
that
$$x\in Cl\bigl(\cap_{i\in J}W^l_i\bigr)\subset Cl\bigl(\cap_{i\in J}W^k_i\bigr)\subset Cl(W^k_r).$$
 On the other hand $x\in V_I$
 implies that $x$ does not belong  to
 \begin{eqnarray*}
 \cap_{i\in I}Cl(U^k_i)\cap Cl(U^k_r)\!\!\!\!\!\!\!&&\supseteq\cap_{i\in I}Cl(W^k_i)\cap
 Cl(W^k_r)\\
&&\supseteq Cl(\cap_{i\in I}W^k_i)\cap Cl(W^k_r).
\end{eqnarray*}
So $x\notin Cl(W^k_r)$ because $x\in V_I\subset Cl\bigl(\cap_{i\in
I}W^k_i\bigr)$. This contradiction shows that $J\subset I$.

{\it Step 3.} Note that $Cl(V_I^\ast)\cap V_J^\ast\subset
Cl(V_I)\cap V_J\cap{\cal W}^\ast$ for any $I, J\in{\cal N}$. The
desired conclusions follow immediately.
   \hfill$\Box$\vspace{2mm}

\begin{remark}\label{rm:2.32}
{\rm From the proof of Lemma~\ref{lem:2.31} one easily sees: if
there exists a positive integer $1<n_0<n$ such that $W_I=\emptyset$
for any $I\in{\cal N}$ with $|I|\ge n_0$ then one only needs to take
$n(n_0-1)$ pairs of open sets,
 $W^j_i\subset\subset U^j_i$, $j=1, 2, \cdots,n_0-1$ such that
$$
 U^1_i\subset\subset
W^2_i\subset\subset U^2_i\cdots \subset\subset
W^{n_0-1}_i\subset\subset U^{n_0-1}_i\subset\subset W_i.
$$
In this case for each $I\in{\cal N}$ with $|I|=k<n_0$ we still use
(\ref{e:2.65}) to define $V_I$. Then $V_I=W_I^{n_0-1}$ for any
$I\in{\cal N}$ with $|I|=n_0-1$, and $\{V_I\,|\, I\in{\cal N}\}$
still satisfies Lemma~\ref{lem:2.31}. Our later arguments are also
effective actually for this general case. } \end{remark}

  Set
\begin{equation}\label{e:2.66}
 \widehat V_I=(\hat\pi_I)^{-1}(V_I)\quad{\rm
and}\quad\widehat F_I=(\widehat\Pi_I)^{-1}(E_I|_{V_I}),
\end{equation}
For $J\subset I\in{\cal N}$, though $\hat\pi^I_J(\widehat
W_I)\subset\widehat W_J$ we cannot guarantee that
$\hat\pi^I_J(\widehat V_I)\subset\widehat V_J$. So $\hat\pi^I_J$
only defines a smooth map from the open subset
$(\hat\pi^I_J)^{-1}(\widehat V_J)\cap\widehat
V_I=(\hat\pi_I)^{-1}(V_I\cap V_J)$
 to $\widehat V_J$. However, {\bf we still denote by}
$\hat\pi^I_J$ the restriction of $\hat\pi^I_J$ to
$(\hat\pi^I_J)^{-1}(\widehat V_J)\cap\widehat V_I$. In this case
it holds that
$$
\hat\pi_J\circ\hat\pi^I_J=\iota^V_{IJ}\circ\hat\pi_I
$$
for any $J\subset I\in {\cal N}$, where
$\iota^V_{IJ}:V_I\cap V_J\hookrightarrow V_J$ is the inclusion.
Similarly, $\widehat\Pi^I_J$ denote the restriction of
$\widehat\Pi^I_J$ to $(\widehat\Pi^I_J)^{-1}(\widehat
F_J)\cap\widehat F_I$.
   We still use
$\hat p_I$ to denote the bundle projection $\widehat
F_I\to\widehat V_I$ induced by one $\hat p_I:\widehat
E_I\to\widehat W_I$ in (\ref{e:2.49}). The system of Banach
bundles
\begin{equation}\label{e:2.67}
(\widehat {\mathcal F}(K), {\widehat V}(K))=\bigl \{(\widehat F_I,
\widehat V_I), \hat\pi_I,
 \hat\pi^I_J, \widehat\Pi_I, \widehat\Pi^I_J, \hat p_I,\Gamma_I\bigm| J\subset I\in{\mathcal
 N}\bigr\},
\end{equation}
 is called {\bf renormalization} of
$(\widehat {\mathcal E}(K), \widehat W(K))$.

 We can define a {\bf global section} of the bundle system  $(\widehat {\mathcal F}(K),
{\widehat V}(K))$ as in Definition~\ref{def:2.28}, but the
compatibility condition (\ref{e:2.59}) is replaced by
\begin{equation}\label{e:2.68}
 S_I=(\hat\Pi^I_J)^\ast S_J:=(\hat\Pi^I_J)^{-1}\circ
S_J\circ\hat\pi^I_J,\quad{\rm on}\; (\hat\pi^I_J)^{-1}(\widehat
V_J).
\end{equation}
Clearly, each section of the bundle system $(\widehat {\mathcal
E}(K), \widehat W(K))$ restricts to one of $(\widehat {\mathcal
F}(K), \widehat V(K))$, still denoted by $\widehat S=\{\widehat
S_I\,|\, I\in{\mathcal N}\}$. The other notions are defined with
the same way.

\begin{lemma}\label{lem:2.33}
 For $1\le l\le n$, if a smooth section $\tilde\sigma_l:\widetilde W_l\to\widetilde E_l$
  has a support contained in $\widetilde W_l^1:=\pi^{-1}_l(W_l^1)$ then
   it may determine a smooth global section $\hat
\sigma_l=\{(\hat\sigma_l)_I: I\in{\cal N}\}$ of the bundle system
$(\widehat {\mathcal F}(K), \widehat V(K))$ and for each
$I\in{\cal N}$ with $l\in I$ the section $(\hat\sigma_l)_I$ is
exactly the restriction of the section given by (\ref{e:2.60}) to
$\widehat V_I$.
 \end{lemma}

\noindent{\bf Proof.}\quad For $I\in {\mathcal N}$ with $|I|=k$
and $l\notin I$, (\ref{e:2.65}) implies that

\begin{eqnarray*}
V_I\!\!\!\!\!&&=\Bigl(\bigcap_{i\in I}W_i^k\Bigr)\setminus
\Bigl(\bigcup_{J:|J|=k+1}\bigl(\bigcap_{j\in
J}Cl(U^k_j)\bigr)\Bigr)\\
&&\subset\Bigl(\bigcap_{i\in I}W_i^k\Bigr)\setminus
\Bigl(\bigl(\bigcap_{i\in I}Cl(U^k_i)\bigr)\cap
Cl(U^k_l)\Bigr)\\
&&\subset\Bigl(\bigcap_{i\in I}W_i^k\Bigr)\setminus
\Bigl(\bigl(\bigcap_{i\in I}Cl(U^k_i)\bigr)\cap
Cl(W^1_l)\Bigr)\\
&&=\Bigl(\bigcap_{i\in I}W_i^k\Bigr)\setminus
 Cl(W^1_l)\\
&&\subset\Bigl(\bigcap_{i\in I}Cl(W_i^k)\Bigr)\setminus
 W^1_l
\end{eqnarray*}
where the second inclusion is because $W^1_l\subset U^1_l\subset
U^k_l$, and the second equality comes from the fact that
$\cap_{i\in I}W_i^k\subset\cap_{i\in I}U^k_i$. It follows that
$$
Cl(V_I)\subset\Bigl(\bigcap_{i\in I}Cl(W_i^k)\Bigr)\setminus
 W^1_l
$$
since the left side is closed.   But $\pi_l({\rm
supp}(\tilde\sigma_l))\subset W_l^1$.
 Hence
\begin{equation}\label{e:2.69}
\pi_l({\rm supp}(\tilde\sigma_l))\cap
Cl(V_I)=\emptyset\quad\forall\, l\notin I,\,1\le l\le n.
 \end{equation}
 We define
\begin{equation}\label{e:2.70}
(\hat\sigma_l)_I=\left\{\begin{array}{lc}
 & 0\;\qquad {\rm if}\; l\notin
I,\\
&\tilde\sigma_{lI}|_{\widehat V_I}\;{\rm if}\; l\in I.
\end{array}\right.
 \end{equation}
Here $\sigma_{lI}$ is given by (\ref{e:2.60}). It is
 easily checked that  the
collection $\{(\hat\sigma_l)_I: I\in{\cal N}\}$ is compatible in
the sense of (\ref{e:2.68}). Hence $\hat
\sigma_l=\{(\hat\sigma_l)_I: I\in{\cal N}\}$ is a smooth global
section of the bundle system $(\widehat{\cal F}(K),\widehat
V(K))$.
 \hfill$\Box$\vspace{2mm}

For an open neighborhood ${\cal W}^\ast$ of $K$ in
$\cup^n_{i=1}W_i^1$ we set $V^\ast_I:=V_I\cap {\cal W}^\ast$ and
\begin{equation}\label{e:2.71}
\widehat V^\ast_I=(\hat\pi_I)^{-1}(V^\ast_I),\qquad\widehat
F^\ast_I=(\hat\Pi_I)^{-1}(F_I|_{V^\ast_I}).
\end{equation}
As above we get a  system of Banach bundles
\begin{equation}\label{e:2.72}
(\widehat {\mathcal F}^\ast, \widehat V^\ast)=\bigl \{(\widehat
F_I^\ast, \widehat V_I^\ast), \hat\pi_I, \hat\pi^I_J,
\widehat\Pi_I, \widehat\Pi^I_J, \hat p_I,\Gamma_I\bigm| J\subset
I\in{\mathcal N}\bigr\},
\end{equation}
which is called the {\bf restriction} of $(\widehat {\mathcal F},
\widehat V)$ to the open subset ${\cal W}^\ast$.
 Similarly, we can
define its global section. Clearly, each global section $\hat
\sigma =\{(\hat\sigma)_I\,|\,I\in{\cal N}\}$ of $(\widehat{\cal
F},\widehat V)$ restricts to a global section of $(\widehat
{\mathcal F}^\ast, \widehat V^\ast)$, still denoted by
$\hat\sigma$ without confusions.

\begin{remark}\label{rem:2.34}{\rm  For each $i=1,\cdots,n$ let us
take pairs of open subsets $W^{+j}_i\subset\subset U^{+j}_i$,
$j=1,\cdots, n-1$ such that
$$
W^j_i\subset\subset W^{+j}_i\subset\subset U^{+j}_i\subset\subset
U^j_i,\;j=1,\cdots, n-1.
$$
Then for each $I\in{\cal N}$ with $|I|=k$ we follow (\ref{e:2.65})
to define
$$
 V^+_I:=\Bigl(\bigcap_{i\in I}W^{+k}_i\Bigr)\setminus
\Bigl(\bigcup_{J:|J|>k}\bigl(\bigcap_{j\in
J}Cl(U^{+k}_j)\bigr)\Bigr)
$$
 and  get another renormalization  $(\widehat {\mathcal F}^+(K), \widehat V^+(K))$ of
 the system of Banach bundles
$(\widehat {\mathcal E}(K), \widehat W(K))$.
 It is easily seen that for any section
$\tilde\sigma_l:\widetilde W_l\to\widetilde E_l$ with  support in
$\widetilde W_i^1$, the same reason as above may yield a global
smooth section $\hat\sigma_l^+=\{(\hat\sigma_l)^+_I: I\in{\cal
N}\}$ of the bundle system $(\widehat {\mathcal F}^+(K), \widehat
V^+(K))$ which  restricts to $\hat\sigma_l$ on $(\widehat
{\mathcal F}(K), {\widehat V}(K))$. Note that for any $I\in{\cal
N}$,
\begin{eqnarray*}
Cl(V_I)\!\!\!\!\!&&\subset\Bigl(\bigcap_{i\in
I}Cl(W^{k}_i)\Bigr)\setminus
\Bigl(\bigcup_{J:|J|>k}\bigl(\bigcap_{j\in J}U^{k}_j\bigr)\Bigr)\\
&&\subset\Bigl(\bigcap_{i\in I}W^{+k}_i\Bigr)\setminus
\Bigl(\bigcup_{J:|J|>k}\bigl(\bigcap_{j\in
J}Cl(U^{+k}_j)\bigr)\Bigr)=V_I^+
\end{eqnarray*}
because $Cl(W^{k}_i)\subset W^{+k}_i$ and $Cl(U^{+k}_i)\subset
U^{k}_i$ for any $i, k$.
 Let $Cl(\widehat
V_I^\ast)$ (resp. $Cl(\widehat F_I^\ast)$) be the closure
$\widehat V_I^\ast$ (resp. $\widehat F_I^\ast$) in $\widehat
V_I^+$ (resp. $\widehat F_I^+$). Then we may also consider the
system of Banach bundles
$$
(Cl(\widehat {\mathcal F}^\ast), Cl(\widehat V^\ast))=\bigl
\{(Cl(\widehat F_I^\ast), Cl(\widehat V_I^\ast)), \hat\pi_I,
 \hat\pi^I_J, \widehat\Pi_I,\widehat\Pi^I_J, \hat p_I,\Gamma_I\bigm| J\subset I\in{\mathcal N}\bigr\}
$$
and define its global section by requiring that
$S_J(\hat\pi^I_J(x))=\hat\Pi^I_J(S_I(x))$  for $x\in
(\hat\pi^I_J)^{-1}(Cl(\widehat V_J^\ast))\cap Cl(\widehat V_I^\ast)$
and $J\subset I\in{\cal N}$. In this case, the above section
$\hat\sigma_l^+$  naturally restricts to a section $(Cl(\widehat
{\mathcal F}^\ast), Cl(\widehat V^\ast))$, still denoted by
$\hat\sigma_l$ without confusions. Moreover, the section $\widehat
S=\{\widehat S_I\,|\, I\in{\mathcal N}\}$ of the bundle system
$(\widehat {\mathcal E}, \widehat W)$  in Lemma~\ref{lem:2.33}
naturally restricts to a global (smooth) section of $(Cl(\widehat
{\mathcal F}^\ast), Cl(\widehat V^\ast))$. }
\end{remark}

\subsection{Transversality}\label{sec:2.4}

As in Section~\ref{sec:2.3} we also assume that $X$ is $1$-regular
and effective.  Remark~\ref{rm:2.49} is still effective.
 From now on we always make:

\begin{assumption}\label{ass:2.35}{\rm
 $(X, E, S)$  is a Banach Fredholm orbibundle of index $r$ and
 compact zero locus $Z(S)$.}
\end{assumption}

Taking $K=Z(S)$ then (\ref{e:2.58}) and (\ref{e:2.67}) give a
system of Banach bundles
$$
\bigl(\widehat{\mathcal E}, \widehat W\bigr)=
\Bigl\{\bigl(\widehat E_I, \widehat W_I\bigr), \hat\pi_I,
\widehat\Pi_I,  \hat\pi^I_J,\widehat\Pi^I_J, \hat p_I,
\Gamma_I,\lambda^I_J\,\bigm|\, J\subset I\in{\mathcal N}\Bigr\},
$$
 and its renormalization system of Banach bundles
\begin{equation}\label{e:2.73}
(\widehat {\mathcal F}, {\widehat V})=\bigl \{(\widehat F_I,
\widehat V_I), \hat\pi_I, \hat\pi^I_J, \widehat\Pi_I,
\widehat\Pi^I_J, \hat p_I,\Gamma_I\,\lambda^I_J\bigm| J\subset
I\in{\mathcal
 N}\bigr\}.
\end{equation}
Since the Fredholm map is locally proper, by shrinking
 open sets $W_i$ we can make:

\begin{assumption}\label{ass:2.36}
{\rm Each Banach bundle $\widetilde E_i\to \widetilde W_i$ is
trivializable  and the unique $\Gamma_i$-equivariant lift section
$\widetilde S_i:\widetilde W_i\to\widetilde E_i$ of $S|_{W_i}$ has
a proper Fredholm trivialization representative.}
\end{assumption}

Let $\widehat S=\{\widehat S_I\,|\,I\in{\cal N}\}$ be the Fredholm
section of the bundle systems $\bigl(\widehat{\mathcal E}, \widehat
W\bigr)$ produced by Lemma~\ref{lem:2.29} from the section $S$. It
restricts to such a section of the bundle systems
 $(\widehat {\mathcal E}, {\widehat V})$, also denoted by
 $\widehat S=\{\widehat S_I\,|\,I\in{\cal N}\}$
 without confusions.

\subsubsection{Local transversality}

For each $i=1,\cdots, n$ let $\widetilde S_i:\widetilde W_i\to
\widetilde E_i$ be the $\Gamma_i$-equivariant lifts of $S|_{W_i}$.
Fix a trivialization
\begin{equation}\label{e:2.74}
{\cal T}_i: \widetilde W_i\times (\widetilde E_i)_{\tilde x_i}\to
\widetilde E_i,
\end{equation}
 we get a representation of $\widetilde S_i$,
\begin{equation}\label{e:2.75}
\widetilde S_i^T:\widetilde W_i\to(\widetilde E_i)_{\tilde
x_i},\;i.e.,\; \widetilde S_i(\tilde x)={\cal T}_i(\tilde x,
\widetilde S_i^T(\tilde x))\;\forall \tilde x\in\widetilde W_i.
\end{equation}
It is a proper Fredholm map  by Assumption~\ref{ass:2.36}. Since the
differential $d\widetilde S^T_i(\tilde x_i): T_{\tilde
x_i}\widetilde W_i\to (\widetilde E_i)_{\tilde x_i}$ is a linear
Fredholm operator there exist finitely many nonzero elements
$v_{ij}\in (\widetilde E_i)_{\tilde x_i}$, $j=1,\cdots, m_i$, such
that
\begin{equation}\label{e:2.76}
d\widetilde S^T_i(\tilde x_i)(T_{\tilde x_i}\widetilde W_i)+ {\rm
span}(\{v_{i1},\cdots, v_{im_i}\})=(\widetilde E_i)_{\tilde x_i}.
\end{equation}
 Using $v_{ij}$ and the trivialization ${\cal T}_i$  we
 get smooth sections $\tilde s_{ij}:\widetilde W_i\to\widetilde
E_i$ by
$$\tilde s_{ij}(\tilde x)={\cal T}_i(\tilde x, v_{ij})\;\forall \tilde x\in\widetilde W_i,\;j=1,\cdots,m_i.$$
Note that $\tilde s_{ij}(\tilde x)\ne 0$ for any $\tilde
x\in\widetilde W_i$ since under the trivialization ${\cal T}_i$
the section $\tilde s_{ij}$ has exactly a representation,
$$
\widetilde W_i\to (\widetilde E_i)_{\tilde x_i},\;\tilde x\mapsto
v_{ij}.
$$
So using the implicit function theorem  (by furthermore shrinking
$W_i$ and increasing $n$) we can make:

\begin{assumption}\label{ass:2.37}
$d\widetilde S^T_i(\tilde x)(T_{\tilde x}\widetilde W_i)+ {\rm
span}(\{v_{i1},\cdots, v_{im_i}\})=(\widetilde E_i)_{\tilde x_i}$
for any $\tilde x\in\widetilde W_i$, $i=1,\cdots,n$.
\end{assumption}

For each $i=1,\cdots, n$ let us take a $\Gamma_i$-invariant smooth
cut-off function  $\gamma_i:\widetilde W_i^1\to [0,1]$ with
support ${\rm supp}(\gamma)\subset\subset\widetilde W_i^1$ and
denote by
\begin{equation}
\widetilde U^0_i:=\{\tilde x\in\widetilde W_i^1\,|\,
\gamma_i(\tilde x)>0\}\quad{\rm and}\quad U^0_i=\pi_i(\widetilde
U^0_i).\label{e:2.77}
\end{equation}
Then $U^0_i\subset\subset W_i^1$, $i=1,\cdots,n$. We may also
assume that
\begin{equation}\label{e:2.78}
Z(S)\subset\bigcup^n_{i=1} U^0_i.
\end{equation}
This can be done since $Z(S)\subset\cup^n_{i=1}W_i^1$. Setting
\begin{equation}\label{e:2.79}
\tilde\sigma_{ij}:=\gamma_i\cdot \tilde s_{ij},\; j=1,\cdots,m_i,
\end{equation}
they are smooth sections of the Banach bundle $\widetilde
E_i\to\widetilde W_i$ with support in $\widetilde W_i^1$, and have
representations under the trivialization ${\cal T}_i$,
\begin{equation}\label{e:2.80}
\tilde\sigma^T_{ij}:\widetilde W_i\to (\widetilde E_i)_{\tilde
x_i},\;\tilde x\mapsto \gamma(\tilde x)v_{ij},\;j=1,\cdots, m_i.
\end{equation}
 It easily follows from Assumption~\ref{ass:2.37} that
\begin{equation}\label{e:2.81}
\left.\begin{array}{ll}
 d\widetilde S^T_i(\tilde x)(T_{\tilde
x}\widetilde W_i)+ {\rm span}(\{\tilde\sigma^T_{i1}(\tilde
x),\cdots, \tilde\sigma^T_{im_i}(\tilde x)\})=(\widetilde
E_i)_{\tilde x_i}\\
\forall\tilde x\in \widetilde U_i^0,\; i=1,\cdots, n.
\end{array}
\right\}
\end{equation}
   By Lemma~\ref{lem:2.33}
 each $\tilde\sigma_{ij}$  yields a smooth section
$\hat\sigma_{ij}=\{(\hat\sigma_{ij})_I\,|\, I\in{\cal N}\}$ of the
 bundle system $(\widehat {\mathcal E}, {\widehat V})$ in (\ref{e:2.73}),
 $j=1,\cdots, m_i$, $i=1,\cdots,n$.

\subsubsection{Global transversality}

 Let $m=m_1+\cdots+m_n$.  Consider
the obvious pullback  Banach bundle system
\begin{eqnarray}
&&\bigl({\bf P}_1^\ast\widehat {\mathcal F}, \widehat
V\times\R^m\bigr)\label{e:2.82}\\
&&=\bigl \{({\bf P}_1^\ast\widehat F_I, \widehat V_I\times\R^m),
\hat\pi_I,  \hat\pi^I_J, \widehat\Pi_I, \widehat\Pi^I_J, \hat
p_I,\Gamma_I\,\bigm|\, J\subset I\in{\mathcal N}\bigr\}, \nonumber
\end{eqnarray}
where ${\bf P}_1$ are the projections to the first factor, and
$\hat\pi_I, \hat\pi^I_J, \widehat\Pi_I,\widehat\Pi^I_J, \hat p_I$
are naturally pullbacks of those projections in (\ref{e:2.73}). It
has a Fredholm section $\Upsilon=\{\Upsilon_I\,|\, I\in{\cal N}\}$
given by
\begin{eqnarray}
&&\Upsilon_I: \widehat V_I\times\R^m\to {\bf P}_1^\ast\widehat
F_I,\label{e:2.83}\\
&&\qquad \Bigl(\hat x_I, \{t_{ij}\}_{\substack{
 1\le j\le m_i\\
 1\le i\le n}}\Bigr)\mapsto \hat S_I(\hat x_I)+
 \sum^n_{i=1}\sum^{m_i}_{j=1}t_{ij}(\hat\sigma_{ij})_I(\hat
 x_I)\nonumber\\
&&\hspace{35mm} =\hat S_I(\hat x_I)+
 \sum_{i\in I}\sum^{m_i}_{j=1}t_{ij}(\hat\sigma_{ij})_I(\hat
 x_I).\nonumber
 \end{eqnarray}
 The final equality comes from (\ref{e:2.69}).
 Clearly, $\Upsilon_I(\hat x_I, 0)=0$ for any zero $\hat x_I$ of
 $\hat S_I$ in $\widehat V_I$.

\begin{theorem}\label{th:2.38}
There exists an open neighborhood ${\cal W}^\ast$ of $Z(S)$ in
$\cup^n_{i=1}U_i^0$ and $\varepsilon>0$ such that for the
restriction bundle system $(\widehat{\cal F}^\ast,\widehat V^\ast)$
of $(\widehat{\cal F},\widehat V)$ to $W^\ast$ the restriction of
the global section $\Upsilon=\bigl\{\Upsilon_I\,|\, I\in{\cal
N}\bigr\}$ to $({\bf P}_1^\ast\widehat {\mathcal F}^\ast, \widehat
V^\ast\times B_\varepsilon(\R^m))$  is Fredholm and transversal to
the zero section. Consequently, for each $I\in{\cal N}$ the set
\begin{equation}\label{e:2.84}
\widehat\Omega_I(S):=\{(\hat x_I, {\bf t})\in\widehat
V_I^\ast\times B_\varepsilon(\R^m)\,|\, \Upsilon_I(\hat x_I,{\bf
t})=0\}
\end{equation}
is a smooth manifold of dimension ${\rm Ind}(S)+m$ and with
compact closure in $Cl(\widehat V^\ast_I)\times\R^m\subset\widehat
V_I^{+}\times\R^m$. Moreover, the family
$\{\widehat\Omega_I(S)\,|\, I\in{\cal N}\}$ is compatible in the
sense that for any $J\subset I\in{\cal N}$,
$$
\hat\pi^I_J\bigl((\hat \pi^I_J)^{-1}(\widehat V_J^\ast\times
B_\varepsilon(\R^m))\cap\widehat\Omega_I(S)\bigr)={\rm
Im}(\hat\pi^I_J)\cap\widehat\Omega_J(S)\quad{\rm and}
$$
that for any $(\hat x_J, {\bf t})\in{\rm
Im}(\hat\pi^I_J)\cap\widehat\Omega_J(S)$ all
$$\frac{|\Gamma_I|-|\Gamma_I({\hat
x_I})|+1}{|\Gamma_J|-|\Gamma_J({\hat x_J})|+1}=\prod_{i\in
I\setminus J}(|\Gamma_i|-|\Gamma_J({\hat x_J})|^{1/|J|}+1)$$
elements of the inverse image $(\hat\pi^I_J)^{-1}((\hat x_J, {\bf
t}))$ sit in $(\hat \pi^I_J)^{-1}(\widehat V_J^\ast\times
B_\varepsilon(\R^m))\cap\widehat\Omega_I(S)$
 for any $(\hat x_I,{\bf t})\in(\hat\pi^I_J)^{-1}((\hat x_J, {\bf t}))$.
So for
\begin{equation}\label{e:2.85}
\widehat\Omega_J(S)^\circ:=\widehat\Omega_J(S)\cap(\widehat
W_J^\circ\times B_\varepsilon(\R^m)),
\end{equation}
the restriction of the projection
\begin{equation}\label{e:2.86}
\hat\pi^I_J:(\hat \pi^I_J)^{-1}(\widehat V_J^\ast\times
B_\varepsilon(\R^m))\cap\widehat\Omega_I(S)\to {\rm
Im}(\hat\pi^I_J)\cap\widehat\Omega_J(S)
\end{equation}
to $(\hat \pi^I_J)^{-1}(\widehat V_J^\ast\times
B_\varepsilon(\R^m))\cap\widehat\Omega_I(S)^\circ$
 is a $|\Gamma_I|/|\Gamma_J|$-fold (regular) smooth covering to
 ${\rm Im}(\hat\pi^I_J)\cap\widehat\Omega_J(S)^\circ$.
\end{theorem}

\noindent{\bf Proof.}\quad Let $(\widehat {\mathcal F}^+, \widehat
V^+)$ be another renormalization of  the system of Banach bundles
$(\widehat {\mathcal E}, \widehat W)$  as in Remark~\ref{rem:2.34}.
Then $Cl(\widehat V_I)\subset \widehat V^+_I$ for each $I\in{\cal
N}$. Moreover, $\Upsilon=\{\Upsilon_I\,|\, I\in{\cal N}\}$ can
naturally extend to a Fredholm section $({\bf P}_1^\ast\widehat
{\mathcal F}^+, \widehat V^+\times\R^m)$, still denoted
$\Upsilon=\{\Upsilon_I\,|\, I\in{\cal N}\}$. Let $\hat u_I\in
Cl(\widehat V_I)\cap Z(\widehat S_I)$. Then $\Upsilon_I(\hat u_I,
0)=0$. Note that $u_I=\hat\pi_I(\hat u_I)\in\cup^n_{i=1}U_i^0$ and
thus  sits in $U^0_{i_q}\subset W_{i_q}^1$ for some $i_q\in I$
because of (\ref{e:2.69}). By Remark~\ref{rem:2.26} and as in
(\ref{e:2.37}) we have the local chart of the Banach manifold
$\widehat W_I$ near $\hat u_I$,
\begin{equation}\label{e:2.87}
\widehat{\phi_I\circ\lambda_I^q}:\widetilde O_{i_q}\to \widehat
W_I,\;\tilde x\mapsto [\tilde u_I,\phi_I,
\phi_I\circ\lambda_I^q(\tilde x)],
\end{equation}
such that $\hat u_I=\widehat{\phi_I\circ\lambda_I^q}(\tilde
u_{i_q})$ for $\tilde u_{i_q}\in\widetilde O_{i_q}$.  We also have
the corresponding Banach bundle chart as in (\ref{e:2.50}),
\begin{equation}\label{e:2.88}
\widehat{\Phi_I\circ\Lambda_I^q}:\widetilde E_{i_q}|_{\widetilde
O_{i_q}}\to \widehat E_I,\;\tilde \xi\mapsto \langle\tilde
u_I,\phi_I, \Phi_I\circ\Lambda_I^q(\tilde \xi)\rangle.
\end{equation}
In these charts, as in (\ref{e:2.63}) the section $\hat S_I$ has
the local expression
$$
 \widetilde O_{i_q}\to\widetilde
E_{i_q}|_{\widetilde O_{i_q}},\;\tilde x\mapsto\tilde
S_{i_q}(\tilde x),
 $$
 and as in (\ref{e:2.61}) the section
$(\hat\sigma_{ij})_I$ can be represented as
$$
 \widetilde O_{i_q}\to \widetilde
E_{i_q}|_{\widetilde O_{i_q}},
 \; \tilde x\mapsto (\Phi_i\circ\Lambda_{i_qi})^{-1}
 \bigl(\tilde\sigma_l(\phi_i\circ\lambda_{i_qi}(\tilde x))\bigr)
$$
 if $i\in I$. Using these we can get the following
local expression of $\Upsilon_I$  in the natural pullback charts
of those in (\ref{e:2.87}) and (\ref{e:2.88}),
\begin{eqnarray}
&&\Upsilon_{Iq}:\widetilde O_{i_q}\times\R^m\to{\bf
P}_1^\ast(\widetilde E_{i_q}|_{\widetilde
O_{i_q}}),\label{e:2.89}\\
&&\Bigl(\tilde x, \{t_{ij}\}_{\substack{
 1\le j\le m_i\\
 1\le i\le n}}\Bigr) \mapsto\tilde
S_{i_q}(\tilde x) + \sum_{i\in
I}\sum^{m_i}_{j=1}t_{ij}(\Phi_i\circ\Lambda_{i_qi})^{-1}
\bigl(\tilde\sigma_{ij}(\phi_i\circ\lambda_{i_qi}(\tilde
x))\bigr).\nonumber
\end{eqnarray}
(Actually we here should replace $\widetilde O_{i_q}$ by its open
subset $(\widehat{\phi_I\circ\lambda_I^q})^{-1}(\widehat
V_I)\subset\widetilde O_{i_q}$. To save notations we still use
$\widetilde O_{i_q}$.) Let
$\widehat{\phi_I\circ\lambda_I^q}(\tilde u_{i_q})=\hat u_I$. Then
$\tilde u_{i_q}\in\widetilde U^0_{i_q}$ because $u_I\in
U^0_{i_q}$. We need to prove that  the section $\Upsilon_{Iq}$ in
(\ref{e:2.89}) is transversal to the zero section at $(\tilde
u_{i_q}, 0)$.
 Notice that
\begin{eqnarray*}
&&\tilde S_{i_q}(\tilde x) + \sum_{i\in
I}\sum^{m_i}_{j=1}t_{ij}(\Phi_i\circ\Lambda_{i_qi})^{-1}
\bigl(\tilde\sigma_{ij}(\phi_i\circ\lambda_{i_qi}(\tilde
x))\bigr)\\
&&=\tilde S_{i_q}(\tilde x) + \sum^{m_{i_q}}_{j=1}t_{i_qj}
\tilde\sigma_{i_qj}(\tilde x)\\
&&+ \sum_{i\in I\setminus\{
i_q\}}\sum^{m_i}_{j=1}t_{ij}(\Phi_i\circ\Lambda_{i_qi})^{-1}
\bigl(\tilde\sigma_{ij}(\phi_i\circ\lambda_{i_qi}(\tilde
x))\bigr).\nonumber
\end{eqnarray*}
Let us denote by
\begin{equation}\label{e:2.90}
\tilde\tau^q_{ij}(\tilde x) =(\Phi_i\circ\Lambda_{i_qi})^{-1}
\bigl(\tilde\sigma_{ij}(\phi_i\circ\lambda_{i_qi}(\tilde
x))\bigr).
\end{equation}
Then in view of the trivialization representations in
(\ref{e:2.75}) and (\ref{e:2.80}) the section
$$\Upsilon_{Iq}=\tilde S_{i_q} +
\sum^{m_{i_q}}_{j=1}t_{i_qj} \tilde\sigma_{i_qj}+ \sum_{i\in
I\setminus\{i_q\}}\sum^{m_i}_{j=1}t_{ij}\tilde\tau^q_{ij}$$
 has the trivialization representation:
\begin{eqnarray}
&&\Upsilon_{Iq}^T: \widetilde O_{i_q}\times\R^m \to (\widetilde
E_{i_q})_{\tilde x_{i_q}},\label{e:2.91}\\
&&(\tilde x, {\bf t})\mapsto S^T_{i_q}(\tilde x)+
\sum^{m_{i_q}}_{j=1}t_{i_qj} \sigma^T_{i_qj}(\tilde x)+ \sum_{i\in
I\setminus\{i_q\}}\sum^{m_i}_{j=1}t_{ij}\tau^{qT}_{ij}(\tilde
x),\nonumber
\end{eqnarray}
where $\tau^{qT}_{ij}$,   as in (\ref{e:2.75}), is defined by
$$\widetilde \tau^q_{ij}(\tilde x)={\cal T}_i(\tilde x, \widetilde
\tau^{qT}_{ij}(\tilde x))\;\forall \tilde x\in\widetilde W_i.$$
  Note that the tangent map of $\Upsilon_{Iq}^T$ at $(\tilde u_{i_q}, 0)$ is given by
\begin{eqnarray}
&&d\Upsilon_{Iq}^T(\tilde u_{i_q}, 0):T_{\tilde u_{i_q}}\widetilde
W_{i_q}\times\R^m\to (\widetilde
E_{i_q})_{\tilde x_{i_q}},\label{e:2.92}\\
  && (\tilde \xi,{\bf v})\mapsto dS^T_{i_q}(\tilde u_{i_q})(\tilde\xi)+
\sum^{m_{i_q}}_{j=1}v_{i_qj} \sigma^T_{i_qj}(\tilde u_{i_q})+
\sum_{i\in
I\setminus\{i_q\}}\sum^{m_i}_{j=1}v_{ij}\tau^{qT}_{ij}(\tilde
u_{i_q}).\nonumber
\end{eqnarray}
It immediately follows from (\ref{e:2.81}) that the map in
(\ref{e:2.92}) is a surjective continuous linear Fredholm
operator. This shows that
  the section $\Upsilon_{Iq}$ in (\ref{e:2.89})  is transversal to the zero
 section at $(\tilde u_{i_q}, 0)$. Consequently, $\Upsilon_I$ is
transversal to the zero section at $(\hat u_I, 0)$. Note that
$Cl(\widehat V_I)\cap Z(\widehat S_I)$ is a {\bf compact} subset
in $\widehat V^{+}_I$. Using the implicit function theorem we may
get an open neighborhood $U\bigl(Cl(\widehat V_I)\cap Z(\widehat
S_I)\bigr)\subset\subset\widehat V_I^{+}$ of $Cl(\widehat V_I)\cap
Z(\widehat S_I)$ in $\widehat V_I^{+}$ and $\varepsilon>0$ such
that:
\begin{description}
\item[(A)] For any $(\hat u_I, {\bf t})\in U\bigl(Cl(\widehat
V_I)\cap Z(\widehat S_I)\bigr)\times B_\varepsilon(\R^m)$ and the
local trivialization representative $\Upsilon_{Iq}^T$ as above the
tangent map
$$
d\Upsilon_{Iq}^T(\tilde u_{i_q}, {\bf t}):T_{\tilde
u_{i_q}}\widetilde W_{i_q}\times\R^m\to (\widetilde
E_{i_q})_{\tilde x_{i_q}}
$$
is a surjective continuous linear Fredholm operator. In
particular, the restriction of $\Upsilon_I$  to the open subset
$U\bigl(Cl(\widehat V_I)\cap Z(\widehat S_I)\bigr)\times
B_\varepsilon(\R^m)$ is Fredholm and also transversal to the zero
section.

\item[(B)] The restriction of $\Upsilon_I$ to the closure of
$U\bigl(Cl(\widehat V_I)\cap Z(\widehat S_I)\bigr)\times
B_\varepsilon(\R^m)$ (in $\widehat V_I^+\times\R^m$) has compact
zero set.
\end{description}
 Let us take an open
neighborhood ${\cal W}^\ast$ of $Z(S)$ in $\cup^n_{i=1}W_i$ so
that
$$
\hat\pi_I^{-1}({\cal W}^\ast\cap V_I)\subset U\bigl(Cl(\widehat
V_I)\cap Z(\widehat S_I)\bigr)\;\forall I\in{\cal N}.
$$
Then for the bundle $\widehat F_I^\ast\to \widehat V_I^\ast$
obtained from this ${\cal W}^\ast$ and as in (\ref{e:2.71}) it
holds that:
\begin{description}
\item[(C)] For any $(\hat u_I, {\bf t})\in\widehat V^\ast_I\times
B_\varepsilon(\R^m)$ and the local trivialization representative
$\Upsilon_{Iq}^T$ as above the tangent map
$$
d\Upsilon_{Iq}^T(\tilde u_{i_q}, {\bf t}):T_{\tilde
u_{i_q}}\widetilde W_{i_q}\times\R^m\to (\widetilde
E_{i_q})_{\tilde x_{i_q}}
$$
is a {\bf surjective} continuous linear Fredholm operator. In
particular, the section $\Upsilon_I: \widehat V^\ast_I\times
B_\varepsilon(\R^m)\to {\bf P}_1^\ast\widehat F^\ast_I$ is
Fredholm and also transversal to the zero section.

\item[(D)] the zero set of $\Upsilon_I$ has {\bf compact closure}
in $\widehat V_I^{+}\times\R^m$.
\end{description}
 Since ${\cal N}$ is finite, by shrinking $\varepsilon>0$ it follows that the
restriction of the  section $\Upsilon=\bigl\{\Upsilon_I\,|\,
I\in{\cal N}\bigr\}$ to $({\bf P}_1^\ast\widehat {\mathcal
F}^\ast, \widehat V^\ast\times B_\varepsilon(\R^m))$  is Fredholm
and also transversal to the zero section. In particular, from the
implicit function theorem we derive that each
$\widehat\Omega_I(S)$ is a smooth manifold of dimension ${\rm
Ind}(S)+m$ and with compact closure in $Cl(\widehat
V^\ast_I)\times\R^m\subset\widehat V_I^{+}\times\R^m$.

Since $\Upsilon$ restricts to a global section of  $\bigl({\bf
P}_1^\ast\widehat {\mathcal F}^\ast, \widehat V^\ast\times
B_\varepsilon(\R^m)\bigr)$ it is clear that
$\hat\pi^I_J\bigl((\hat \pi^I_J)^{-1}(\widehat V_J^\ast\times
B_\varepsilon(\R^m))\cap\widehat\Omega_I(S)\bigr)={\rm
Im}(\hat\pi^I_J)\cap\widehat\Omega_J(S)$ for any $J\subset
I\in{\cal N}$. Next note that the section $\hat S_I$ (resp.
$(\hat\sigma_{ij})_I$) is defined by the pullback of $\hat S_J$
(resp. $(\hat\sigma_{ij})_J$) under the projection $\hat\pi^I_J$.
It follows that for any $(\hat x_J, {\bf t})\in{\rm
Im}(\hat\pi^I_J)\cap\widehat\Omega_J(S)$ the inverse image
$(\hat\pi^I_J)^{-1}((\hat x_J, {\bf t}))$ must be contained in
$(\hat \pi^I_J)^{-1}(\widehat V_J^\ast\times
B_\varepsilon(\R^m))\cap\widehat\Omega_I(S)$. The other
conclusions are easily derived.
 \hfill$\Box$\vspace{2mm}

{\bf From now on} assume that $X$ is {\bf separable} so that the
Sard-Smale theorem may be used.  Since ${\cal N}$ is a finite set
Theorem~\ref{th:2.38} and the Sard-Smale theorem immediately give:

\begin{corollary}\label{cor:2.39}
There exists a residual subset $B_\varepsilon(\R^m)_{res}\subset
B_\varepsilon(\R^m)$ such that for each ${\bf t}\in
B_\varepsilon(\R^m)_{res}$ the global section
 $\Upsilon^{({\bf t})}=\{\Upsilon^{({\bf t})}_I\,|\, I\in{\cal
N}\}$ of the Banach bundle system $\bigl(\widehat {\mathcal
F}^\ast, \widehat V^\ast\bigr)$ is Fredholm and transversal to the
zero section.
Here ${\bf t}=\{t_{ij}\}_{\substack{ 1\le j\le m_i\\
 1\le i\le n}}$ and
\begin{eqnarray}
&&\Upsilon^{({\bf t})}_I: \widehat V^\ast_I\to \widehat F^\ast_I,\label{e:2.93}\\
&&\qquad\hat x_I \mapsto \Upsilon_{I}(\tilde x, {\bf t})=\hat
S_I(\hat x_I)+
 \sum^n_{i=1}\sum^{m_i}_{j=1}t_{ij}(\hat\sigma_{ij})_I(\hat
 x_I).\nonumber
 \end{eqnarray}
 So the set
$\widehat{\cal M}^{\bf t}_I(S):= (\Upsilon^{({\bf t})}_I)^{-1}(0)$
is a manifold of dimension $r={\rm Ind}(S)$ and has a compact
closure in $Cl(\widehat V^\ast_I)\subset\widehat V_I^{+}$.
Moreover, the family
\begin{equation}\label{e:2.94}
\widehat{\cal M}^{\bf t}(S)=\bigl\{\widehat{\cal M}^{\bf
t}_I(S)\,:\, I\in{\cal N}\bigr\}
\end{equation}
is compatible in the sense that for any $J\subset I\in{\cal N}$,
$$
 \hat\pi^I_J\bigl((\hat
\pi^I_J)^{-1}(\widehat V_J^\ast)\cap\widehat{\cal M}^{\bf
t}_I(S)\bigr)={\rm Im}(\hat\pi^I_J)\cap\widehat{\cal M}^{\bf
t}_J(S)\quad{\rm and}
$$
that for any $\hat x_J\in{\rm Im}(\hat\pi^I_J)\cap\widehat{\cal
M}^{\bf t}_J(S)$ the inverse image $(\hat\pi^I_J)^{-1}(\hat x_J)$
contains
$$\frac{|\Gamma_I|-|\Gamma_I({\hat
x_I})|+1}{|\Gamma_J|-|\Gamma_J({\hat x_J})|+1}=\prod_{i\in
I\setminus J}(|\Gamma_i|-|\Gamma_J({\hat x_J})|^{1/|J|}+1)$$
elements and all sit in $(\hat \pi^I_J)^{-1}(\widehat
V_J^\ast)\cap\widehat{\cal M}^{\bf t}_I(S)$
 for any $\hat x_I\in(\hat\pi^I_J)^{-1}(\hat x_J)$.
So for
\begin{equation}\label{e:2.95}
\left.\begin{array}{ll}
\widehat{\cal M}^{\bf
t}_J(S)^\circ:=\widehat{\cal M}^{\bf t}_J(S)\cap\widehat
W_J^\circ=\widehat{\cal M}^{\bf
t}_J(S)\cap\widehat V_J^{\ast\circ}\quad{\rm and}\\
\widehat{\cal M}^{\bf t}_J(S)^{sing}:=\widehat{\cal M}^{\bf
t}_J(S)\cap\widehat W_J^{sing}=\widehat{\cal M}^{\bf
t}_J(S)\cap\widehat V_J^{\ast sing}
\end{array}\right\}
\end{equation}
it also holds that
\begin{eqnarray*}
&&\hat\pi^I_J\bigl((\hat \pi^I_J)^{-1}(\widehat
V_J^\ast)\cap\widehat{\cal M}^{\bf t}_I(S)^\circ\bigr)={\rm
Im}(\hat\pi^I_J)\cap\widehat{\cal M}^{\bf t}_J(S)^\circ,\\
&&\hat\pi^I_J\bigl((\hat \pi^I_J)^{-1}(\widehat
V_J^\ast)\cap\widehat{\cal M}^{\bf t}_I(S)^{sing}\bigr)={\rm
Im}(\hat\pi^I_J)\cap\widehat{\cal M}^{\bf t}_J(S)^{sing}
\end{eqnarray*}
and the restriction of the projection
\begin{equation}\label{e:2.96}
 \hat\pi^I_J: (\hat \pi^I_J)^{-1}(\widehat
V_J^\ast)\cap\widehat{\cal M}^{\bf t}_I(S)\to {\rm
Im}(\hat\pi^I_J)\cap\widehat{\cal M}^{\bf t}_J(S)
\end{equation}
to $(\hat \pi^I_J)^{-1}(\widehat V_J^\ast)\cap\widehat{\cal
M}^{\bf t}_I(S)^\circ$
 is a $|\Gamma_I|/|\Gamma_J|$-fold
(regular) smooth covering to ${\rm
Im}(\hat\pi^I_J)\cap\widehat{\cal M}^{\bf t}_J(S)^\circ$. Also
note that each $\hat\pi_I: \widehat{\cal M}^{\bf t}_I(S) ({\rm
resp.}\, \widehat{\cal M}^{\bf t}_I(S)^\circ )\to X$ is  proper.
\end{corollary}

Furthermore, we also have:

\begin{proposition}\label{prop:2.40}
\begin{description}
\item[(i)]  For any given small open neighborhood ${\cal U}$ of
$Z(S)$ in $\cup^n_{i=1}U_i^0$ there exists $\epsilon\in
(0,\varepsilon)$ such that $\hat\pi_I(\widehat{\cal M}^{\bf
t}_I(S))\subset {\cal U}$ and ${\bf t}\in B_\epsilon(\R^m)$.

\item[(ii)] For any two ${\bf t},\, {\bf t}'\in
B_\varepsilon(\R^m)_{res}$ there exist generic path
$\gamma:[0,1]\to B_\varepsilon(\R^m)$ with $\gamma(0)={\bf t}$ and
$\gamma(1)={\bf t}'$ such that for each $I\in{\cal N}$    the set
$$(\Psi_I^\gamma)^{-1}(0):=\cup_{t\in
[0,1]}\{t\}\times(\Upsilon_I^{\gamma(t)})^{-1}(0)
$$
is a smooth manifold of dimension ${\rm Ind}(S)+1$ and with boundary
\footnote{We here have assumed that $(X, E, S)$ is oriented and the
moduli spaces have been oriented according to
Proposition~\ref{prop:2.44}. Otherwise the following negative sign
should be removed.} $\{0\}\times\widehat{\mathcal M}_I^{{\bf
t}}(S)\cup(-\{1\}\times\widehat{\mathcal M}_I^{{\bf t}'}(S))$;
Moreover, the family
$(\Psi^\gamma)^{-1}(0)=\{(\Psi_I^\gamma)^{-1}(0)\,|\, I\in{\cal
N}\}$ is compatible in the similar sense to (\ref{e:2.96}).

\item[(iii)] Assume $\delta\in (0,\varepsilon]$ so that
 the closure of $\cup_{I\in{\mathcal N}}
\hat\pi_I(\widehat{\mathcal M}_I^{{\bf t}}(S))$ in $\cup^n_{i=1}W_i$
may be contained in $\cup_{I\in{\cal N}}V_I^\ast$ for any ${\bf
t}\in B_\delta(\R^m)$. (This is always possible by (i) and
Theorem~\ref{th:2.38}.) Then for any ${\bf t}\in B_\delta(\R^m)$ the
set $\cup_{I\in{\mathcal N}} \hat\pi_I(\widehat{\mathcal M}_I^{{\bf
t}}(S))$ is compact.
 Consequently, the
family $\widehat{\cal M}^{\bf t}(S)=\bigl\{\widehat{\cal M}^{\bf
t}_I(S)\,|\, I\in{\cal N}\bigr\}$ is ``like'' an open cover of a
compact manifold.  Specially, each $\widehat{\mathcal M}_I^{{\bf
t}}(S)$ is a finite set provided that ${\rm Ind}(S)=0$.

\item[(iv)] For any $z\in Z(S)$ there exists a residual subset
$B_\varepsilon(\R^m)_{res}^z\subset B_\varepsilon(\R^m)_{res}$ such
that for each ${\bf t}\in B_\varepsilon(\R^m)_{res}^z$ the union
$\cup_{I\in{\cal N}}\widehat{\cal M}^{\bf t}(S)$ cannot contain $z$.
Furthermore, there exist a small neighborhood $O_z$ of $z$ in $X$
and a subset $H_z\subset\R^m$, which is the union of finitely many
subspaces in $\R^m$ of codimension at least one,  such that for a
small open neighborhood ${\cal U}(S^m\cap H_z)$ of $S^m\cap H_z$ in
the unit sphere $S^m$ of $\R^m$ and the open cone $\angle{\cal
U}(S^m\cap H_z)$ spanned by ${\cal U}(S^m\cap H_z)$ it holds that
$$\bigl(\cup_{I\in{\cal N}}\widehat{\cal M}^{\bf t}(S)\bigr)\cap(O_x\cap Z(S))=\emptyset$$
for any ${\bf t}\in B_\varepsilon(\R^m)_{res}\setminus\angle{\cal
U}(S^m\cap H_z)$.
\end{description}
\end{proposition}

\noindent{\bf Proof.}\quad (i) Assume that the conclusion does not
hold. Since ${\cal N}$ is finite there exists a small open
neighborhood ${\cal U}$ of $Z(S)$ in $\cup^n_{i=1}U_i^0$, some
$I\in{\cal N}$, ${\bf t}_{k}\in\R^m$ and $\hat u_{Ik}\in
\widehat{\cal M}^{{\bf t}_k}_I(S)$, $k=1,2,\cdots,$ such that ${\bf
t}_k\to 0$ and $u_{Ik}=\hat\pi_I(\hat u_{Ik})\notin {\cal U}$ for
$k=1,2,\cdots$. Note that $(u_{Ik},{\bf t}_k)\in
\widehat\Omega_I(S)$ for each $k$. By Theorem~\ref{th:2.38},
$\widehat\Omega_I(S)$  has a compact closure in $Cl(\widehat
V^\ast_I)\times\R^m\subset\widehat V_I^+\times\R^m$. We may assume
that $\hat u_{Ik}\to \hat u_I\in Cl(\widehat
V^\ast_I)\subset\widehat V_I^{+}$. It follows that $\Upsilon_I(\hat
u_I, 0)=0$ and thus $u_I=\hat\pi_I(\hat u_I)\in Z(S)\subset {\cal
U}$. This contradicts to the fact that $u_{Ik}\notin{\cal U}$ for
each $k$.

 (ii) Let ${\cal P}^l({\bf t},{\bf t}')$
denote the space of all $C^l$-smooth paths $\gamma:[0,1]\to
B_\varepsilon(\R^m)$ connecting ${\bf t}$ to ${\bf t}'$. It is a
Banach manifold. Consider the pullback Banach bundle system
\begin{eqnarray*}
&&\bigl({\bf P}_1^\ast\widehat {\mathcal F}^\ast, \widehat
V^\ast\times{\cal P}^l({\bf t},{\bf t}')\times [0, 1]\bigr)\\
&&= \bigl\{\bigl({\bf P}_1^\ast\widehat F_I^\ast, \widehat
V_I^\ast\times{\cal P}^l({\bf t},{\bf t}')\times [0,1]\bigr),
\hat\pi_I, \hat\pi^I_J, \hat\Pi_I,\hat\Pi^I_J, \hat
p_I,\Gamma_I\,\bigm|\, J\subset I\in{\mathcal N}\bigr\}\nonumber
\end{eqnarray*}
and its global section $\Psi=\{\Psi_I\,|\, I\in{\cal N}\}$ given
by
\begin{eqnarray}
&&\Psi_I: \widehat V_I^\ast\times{\cal P}^l({\bf t},{\bf
t}')\times [0,1]\to {\bf P}_1^\ast\widehat F_I^\ast,\label{e:2.97}\\
&& (\hat x_I,\gamma, s)\mapsto \hat S_I(\hat x_I)+
\sum^n_{i=1}\sum^{m_i}_{j=1}\gamma(s)_{ij}(\hat\sigma_{ij})_I(\hat
x_I).\nonumber
\end{eqnarray}
Here $\gamma(s)=\{\gamma(s)_{ij}\}_{\substack{
 1\le j\le m_i\\
 1\le i\le n}}$, ${\bf P}_1$ is the projection to the first factor, and
$\hat\pi_I, \hat\pi^I_J, \hat\Pi_I,\hat\Pi^I_J, \hat p_I$ are
naturally pullbacks of those projections in (\ref{e:2.72}). We want
to prove that $\Psi_I$ is transversal to the zero section. Let
$(\hat u_I, \gamma_0,  \tau)\in\widehat V_I^\ast\times{\cal
P}^l({\bf t},{\bf t}')\times [0, 1]$ be the zero of $\Psi_I$. As in
proof of Theorem~\ref{th:2.38} we have $u_I=\hat\pi_I(\hat
u_I)\in\cup^n_{i=1}U_i^0$ and thus may assume that
 $u_I\in U^0_{i_q}\subset W_{i_q}$ for some $i_q\in I$.
Similar to (\ref{e:2.89}), using the pullback charts of those in
(\ref{e:2.87}) and (\ref{e:2.88}) we can get the local expression
of $\Psi_I$,
\begin{eqnarray*}
&&\widetilde O_{i_q}\times{\cal P}^l({\bf t},{\bf t}')\times
[0,1]\to{\bf P}_1^\ast(\widetilde E_{i_q}|_{\widetilde
O_{i_q}}),\\
&&(\tilde x,\gamma, s) \mapsto\tilde S_{i_q}(\tilde x) +
\sum^{m_{i_q}}_{j=1}\gamma(s)_{ij} \tilde\sigma_{i_qj}+ \sum_{i\in
I\setminus\{i_q\}}\sum^{m_i}_{j=1}\gamma(s)_{ij}\tilde\tau^q_{ij},
\end{eqnarray*}
where $\tilde\tau^q_{ij}$ is defined by (\ref{e:2.90}). As in
(\ref{e:2.91}) it has the following trivialization representation:
\begin{eqnarray}
&&\Psi_{Iq}^T:\widetilde O_{i_q}\times{\cal P}^l({\bf t},{\bf
t}')\times [0,1]\to (\widetilde
E_{i_q})_{\tilde x_{i_q}},\label{e:2.98}\\
&&(\tilde x, \gamma, s)\mapsto S^T_{i_q}(\tilde x) +
\sum^{m_{i_q}}_{j=1}\gamma(s)_{i_qj} \sigma^T_{i_qj}(\tilde x)+
\sum_{i\in
I\setminus\{i_q\}}\sum^{m_i}_{j=1}\gamma(s)_{ij}\tau^{qT}_{ij}(\tilde
x).\nonumber
\end{eqnarray}
Here $\tau^{qT}_{ij}$ is as in (\ref{e:2.91}). Then for
$\Upsilon^T_{Iq}$ in (\ref{e:2.91}) it holds that
$\Psi_{Iq}^T(\tilde x, \gamma, s)=\Upsilon^T_{Iq}(\tilde x,
\gamma(s))$. So for $\tau\ne 0,1$ and any $(\tilde\xi,\alpha)\in
T_{(\tilde u_{i_q}, \gamma_0)}\bigl(\widetilde O_{i_q}\times{\cal
P}^l({\bf t},{\bf t}')\bigr)$ it holds that
$$d\Psi_{Iq}^T(\tilde u_{i_q}, \gamma_0, \tau)(\tilde\xi, \alpha, 0)
=d\Upsilon^T_{Iq}(\tilde u_{i_q}, \gamma(\tau))(\tilde\xi,
\alpha(\tau)).
$$
This and (C) in the proof of Theorem~\ref{th:2.38} together imply
that the tangent map $d\Psi_{Iq}^T(\tilde u_{i_q}, \gamma_0, \tau)$
is a surjective continuous linear Fredholm operator.

If $\tau=0$ or $1$ then $\gamma(0)={\bf t}$ and $\gamma(1)={\bf
t}'$. Note that the sections $\Upsilon^{({\bf t})}$ and
$\Upsilon^{({\bf t}')}$ are transversal to the zero section. That
is,
\begin{eqnarray*}
&&(\widehat{\phi_I\circ\lambda_I^1})^{-1}(\widehat
V_I)\subset\widetilde O_{i_q}\to\widetilde E_{i_q}|_{\widetilde
O_{i_q}},\nonumber\\
&&\tilde x \mapsto\tilde S_{i_q}(\tilde x) + \sum_{i\in
I\setminus\{i_q\}}\sum^{m_i}_{j=1}t_{ij}\tilde\tau_{ij}(\tilde x)\quad{\rm and}\\
&&\tilde x\mapsto\tilde S_{i_q}(\tilde x) + \sum_{i\in
I\setminus\{i_q\}}\sum^{m_i}_{j=1}t'_{ij}\tilde\tau_{ij}(\tilde x)
\end{eqnarray*}
 are transversal to the zero section, or their trivialization
representations
\begin{eqnarray*}
&& (\widehat{\phi_I\circ\lambda_I^q})^{-1}(\widehat
V_I)\cap\widetilde O_{i_q}^\alpha\to (\widetilde
E_{i_q\alpha})_{\tilde x_{i_q\alpha}},\nonumber\\
&&\tilde x\mapsto S^T_{i_q}(\tilde x)+
\sum^{m_{i_q}}_{j=1}t_{i_qj} \sigma^T_{i_qj}(\tilde x)+ \sum_{i\in
I\setminus\{i_q\}}\sum^{m_i}_{j=1}t_{ij}\tau^{qT}_{ij}(\tilde
x)\quad{\rm and}\\
&&\tilde x\mapsto S^T_{i_q}(\tilde x)+
\sum^{m_{i_q}}_{j=1}t'_{i_qj} \sigma^T_{i_qj}(\tilde x)+
\sum_{i\in
I\setminus\{i_q\}}\sum^{m_i}_{j=1}t'_{ij}\tau^{qT}_{ij}(\tilde x)
\end{eqnarray*}
have the surjective tangent map at $\tilde u_{i_q}$ respectively.
These imply that if $\tau=0$ (resp. $\tau=1$) the map
$\Psi^T_{Iq}$ in (\ref{e:2.98}) has the surjective tangent map at
$(\tilde u_{i_q},\gamma_0, 0)$ (resp. $(\tilde u_{i_q}, \gamma_0,
1)$).

In summary, we have proved that the tangent map
$d\Psi^T_{Iq}(\tilde u_{i_q}, \gamma_0, \tau)$ is a surjective
continuous linear Fredholm operator. So the section $\Psi_I$ is
transversal to the zero section at $(\hat u_I, \gamma_0, \tau)$.
Then the standard arguments show that the section
\begin{eqnarray}
&& \Psi_I^\gamma: \widehat
V_I^\ast\times [0,1]\to {\bf P}_1^\ast\widehat F_I,\label{e:2.99}\\
&&\qquad (\hat x_I, s)\mapsto \hat S_I(\hat x_I)+
\sum^n_{i=1}\sum^{m_i}_{j=1}\gamma(s)_{ij}(\hat\sigma_{ij})_I(\hat
x_I)\nonumber
 \end{eqnarray}
is transversal to the zero section for a generic $\gamma\in{\cal
P}^l({\bf t},{\bf t}')$. Since ${\cal N}$ is finite we get that the
global section $\Psi^\gamma=\{\Psi^\gamma_I\,|\, I\in{\cal N}\}$ of
the pullback  Banach bundle system
\begin{eqnarray*}
&&\bigl({\bf P}_1^\ast\widehat {\mathcal F}^\ast, \widehat
V^\ast\times [0, 1]\bigr)\\
&&= \bigl\{({\bf P}_1^\ast\widehat F_I^\ast, \widehat
V_I^\ast\times [0,1]), \hat\pi_I, \hat\pi^I_J,
\hat\Pi_I,\hat\Pi^I_J, \hat p_I,\Gamma_I\,\bigm|\, J\subset
I\in{\mathcal N}\bigr\},\nonumber
\end{eqnarray*}
 is transversal to
the zero section for a generic $\gamma\in{\cal P}^l({\bf t},{\bf
t}')$. Fix a generic $\gamma\in{\cal P}^l({\bf t},{\bf t}')$ then
the (smooth) manifold with boundary (and corner)
\begin{equation}\label{e:2.100}
(\Psi^\gamma_I)^{-1}(0)=\bigl\{(\hat x_I, s)\in \widehat
V_I^\ast\times [0,1]\,|\, \Psi_I(\hat x_I,\gamma,s)=0\bigr\}
\end{equation}
forms a cobordism between $\widehat{\mathcal M}_I^{{\bf t}}(S)$ and
$\widehat{\mathcal M}_I^{{\bf t}'}(S)$ for each $I\in{\cal N}$. Note
that $\Upsilon_I^{\gamma(t)}(x)=\Psi_I^{\gamma}(x,t)$ for
$(x,t)\in\widehat V_I^\ast\times [0, 1]$. (ii) is proved.

(iii) By the choice of $\delta$, for any ${\bf t}\in
B_\delta(\R^m)$, the closure of $\cup_{I\in{\mathcal N}}
\hat\pi_I(\widehat{\mathcal M}_I^{{\bf t}}(S))$ in
$\cup^n_{i=1}W_i$ is contained in $\cup_{I\in{\cal N}}V_I^\ast$.
We want to prove that $\cup_{I\in{\mathcal N}}
\hat\pi_I(\widehat{\mathcal M}_I^{{\bf t}}(S))$ is compact.
 By Corollary~\ref{cor:2.39} each
$\hat\pi_I(\widehat{\mathcal M}_I^{{\bf t}}(S))$ has the compact
closure in $\cup_{I\in{\cal N}}\widehat
V_I^+\subset\cup^n_{i=1}W_i$, and the compact subset
$$Cl\bigl(\cup_{I\in{\mathcal N}} \hat\pi_I(\widehat{\mathcal
M}_I^{{\bf t}}(S))\bigr)=\cup_{I\in{\mathcal N}}
Cl\bigl(\hat\pi_I(\widehat{\mathcal M}_I^{{\bf
t}}(S))\bigr)\subset\cup_{I\in{\cal N}}V_I^\ast.
$$
It follows that any  $x\in Cl\bigl(\hat\pi_I(\widehat{\mathcal
M}_I^{{\bf t}}(S))\bigr)\setminus \hat\pi_I(\widehat{\mathcal
M}_I^{{\bf t}}(S))\subset Cl(V_I^\ast)$ may be contained in
$V_L^\ast$ for some $L\in{\cal N}$. So $Cl(V_I)\cap
Cl(V_L)\ne\emptyset$. By Lemma~\ref{lem:2.31} it holds that
$I\subset L$ or $L\subset I$. Let $I\subset L$ and $x=\hat\pi_L(\hat
x_L)$ for some $\hat x_L\in\widehat V_L^\ast$. By
Remark~\ref{rem:2.34} we may assume $x=\hat\pi_I(\hat x_I)$ for some
$\hat x_I\in Cl(\widehat V_I^\ast)$. Since $\Upsilon^{({\bf
t})}=\{\Upsilon^{({\bf t})}_I\,|\, I\in{\cal N}\}$ may be extended
into a global section of the bundle system $\bigl(Cl(\widehat
{\mathcal F}^\ast), Cl(\widehat V^\ast)\bigr)$ naturally. We get
that $\Upsilon^{({\bf t})}_I(\hat x_I)=0$. Since $\hat\pi_I^L(\hat
x_L)=\hat x_I$ we get that $\Upsilon^{({\bf t})}_L(\hat x_L)=0$ and
thus $x\in\hat\pi_L(\widehat{\mathcal M}_L^{{\bf t}}(S))$ because of
$x\in V_L^\ast$. This shows that
$Cl\bigl(\hat\pi_I(\widehat{\mathcal M}_I^{{\bf
t}}(S))\bigr)\subset\cup_{J\in{\mathcal N}}
\hat\pi_J(\widehat{\mathcal M}_J^{{\bf t}}(S)$. For the case
$L\subset I$ we have $\hat\pi^I_L(\hat x_I)=\hat x_L$. It follows
from $\Upsilon^{({\bf t})}_I(\hat x_I)=0$ that $\Upsilon^{({\bf
t})}_L(\hat x_L)=0$. Hence it also holds that $x=\hat\pi_L(\hat
x_L)\in \hat\pi_L(\widehat{\mathcal M}_L^{{\bf t}}(S))$.
 In summary we get that $\cup_{I\in{\mathcal
N}}Cl\bigr(\hat\pi_I(\widehat{\mathcal M}_I^{{\bf
t}}(S))\bigr)\subset\cup_{I\in{\mathcal N}}
\hat\pi_I(\widehat{\mathcal M}_I^{{\bf t}}(S))$ and hence
 $\cup_{I\in{\mathcal N}} \hat\pi_I(\widehat{\mathcal M}_I^{{\bf
t}}(S))=\cup_{I\in{\mathcal N}}Cl\bigr(\hat\pi_I(\widehat{\mathcal
M}_I^{{\bf t}}(S))\bigr)$
 is compact.

The final claim is actually a direct consequence of the above
arguments. For the sake of clearness we prove it as follows.
 If ${\rm Ind}(S)=0$ then each $\widehat{\mathcal M}_I^{{\bf t}}(S)$
is a manifold of dimension zero. Assume that it contains infinitely
many points $\hat x^{(k)}$, $k=1, 2,\cdots$. We may assume that
$\{\hat\pi_I(\hat x^{(k)})\}$ converges to some
$x\in\hat\pi_L\bigl(\widehat{\mathcal M}_L^{{\bf t}}(S)\bigr)$
because $\cup_{I\in{\mathcal N}} \hat\pi_I(\widehat{\mathcal
M}_I^{{\bf t}}(S))$ is compact. Note that the inverse image of each
point by $\hat\pi_I$ contains at most $|\Gamma_I|$ points and that
$\{\hat\pi_I(\hat x^{(k)})\}$ are contained in the closed subset
$Cl(\widehat V_I^\ast)\subset\widehat V_I^{+}$. After passing to a
subsequence we may assume that $\{\hat x^{(k)}\}$ converges to some
$\hat x_I\in Cl\bigl(\widehat{\mathcal M}_I^{{\bf t}}(S)\bigr)$.
Thus $\Upsilon^{({\bf t})}_I(\hat x_I)=0$. However, the section
$\Upsilon^{({\bf t})}_I$ is still transversal to the zero section at
$\hat x_I$ because ${\bf t}\in B_\varepsilon(\R^m)_{res}$. This
destroys the manifold structure of $(\Upsilon^{({\bf
t})}_I)^{-1}(0)$ near $\hat x_I$. The desired conclusion is proved.

(iv) By (\ref{e:2.78}) we may assume $z\in U_1^0$ without loss of
generality. It follows from the constructions of $\tilde s_{1j}$ and
$\tilde\sigma_{1j}$ in (\ref{e:2.79}) that $\tilde\sigma_{1j}(\tilde
z)\ne 0$ for any $\tilde z\in\pi_1^{-1}(z)$. By Lemma~\ref{lem:2.33}
one easily sees that for $I\in{\cal N}$ with $1\notin I$,
$(\hat\sigma_{1j})_I(\hat z)\ne 0$ for any $\hat
z_I\in\hat\pi_I^{-1}(z)$. For each such $I$ and $\hat z_I$ the
kernel of the nonzero linear map
$$
\R^m\to \hat F_{I\hat z},\;{\bf t}\to
\sum^n_{i=1}\sum^{m_i}_{j=1}t_{ij}(\hat\sigma_{ij})_I(\hat z_I)
$$
has at lease codimension one. Denote the kernel by $H_{\hat z_I}$.
Let ${\cal N}_1=\{I\in{\cal N}\,|\, 1\notin I\}$ and
$H_z=\cup_{I\in{\cal N}_1}\cup_{\hat
z_I\in\hat\pi_I^{-1}(z)}H_{\hat z_I}$. Then
$B_\varepsilon(\R^m)^z_{res}:=(\R^m\setminus H_z)\cap
B_\varepsilon(\R^m)_{res}$ is also a residual subset in
$B_\varepsilon(\R^m)$ and for each ${\bf t}\in
B_\varepsilon(\R^m)^z_{res}$, $I\in{\cal N}_1$ and $\hat
z_I\in\hat\pi_I^{-1}(z)$ we have
$$\sum^n_{i=1}\sum^{m_i}_{j=1}t_{ij}(\hat\sigma_{ij})_I(\hat
z_I)\ne 0.
$$
This implies the conclusion.

Next the proof of the second claim is similar to
Remark~\ref{rem:1.8}. For  a given small open neighborhood ${\cal
U}(S^m\cap H_z)$ of $S^m\cap H_z$ in the unit sphere $S^m$ in
$\R^m$ it is easily proved that there exists a small neighborhood
$O_z$ of $z$ in $X$ such that for each $I\in{\cal N}_1$, $\hat
x_I\in\hat\pi_I^{-1}(O_z)$ and ${\bf t}\in S^m\setminus{\cal
U}(S^m\cap H_z)$ one has
$$\sum^n_{i=1}\sum^{m_i}_{j=1}t_{ij}(\hat\sigma_{ij})_I(\hat
x_I)\ne 0.
$$
So for the open cone $\angle{\cal U}(S^m\cap H_z)$ spanned by
${\cal U}(S^m\cap H_z)$ and any ${\bf t}\in
B_\varepsilon(\R^m)_{res}\setminus\angle{\cal U}(S^m\cap H_z)$ it
holds that
$$\bigl(\cup_{I\in{\cal N}}\widehat{\cal M}^{\bf t}(S)\bigr)\cap(O_z\cap Z(S))=\emptyset.$$
\hfill$\Box$\vspace{2mm}

\begin{remark}\label{rem:2.41}
{\rm Now $\widehat{\cal M}^{\bf t}_J(S)^\circ$ in (\ref{e:2.95})
is not necessarily dense in $\widehat{\cal M}^{\bf t}_J(S)$. It
might even be empty. The set $\widehat{\cal M}^{\bf
t}_J(S)^{sing}$ in (\ref{e:2.95}) is not necessarily a subvariety
in $\widehat{\cal M}^{\bf t}_J(S)$ of lower dimension. }
\end{remark}

Because of the shortcoming stated in the above remark one cannot use
the family $\{\widehat{\cal M}^{\bf t}_J(S)\,|\, I\in{\cal N}\}$  to
construct the desired Euler class. We need to  suitably refine  the
above arguments so as to improve our present results.

Firstly, we require that the Banach orbifold charts $(\widetilde
W_i, \Gamma_i,\pi_i)$ at the beginning of Section~\ref{sec:2.2.1}
are all centred at some $x_i\in X$, $i=1,\cdots, n$. This means that
$\pi_i^{-1}(x_i)=\{\tilde x_i\}$ and $\Gamma_i=\Gamma_i(\tilde
x_i)$, $i=1,\cdots, n$. If $\Gamma_i\ne\{\1\}$, i.e., $x_i\in
X^{sing}$, by the proof of Claim~\ref{cl:2.5} we may shrink $W_i$ to
require that
\begin{equation}\label{e:2.101}
 \pi_i^{-1}(X^{sing}\cap W_i)=\bigcup_{g\in\Gamma_i
 \setminus\{\1\}}F_i\bigl(N(0_{\tilde x_i})\cap{\nu}(dg(\tilde x_i))\bigr),
 \end{equation}
where ${\nu}(dg(\tilde x_i))={\rm Ker}({\1}-dg(\tilde x_i))$ and
$F_i$ is a $\Gamma_i$-equivariant diffeomorphism  from a
neighborhood $N(0_{\tilde x_i})$
 of the origin in $T_{\tilde x_i}\widetilde W_i$ onto $\widetilde
 W_i$. Note that $\Gamma_i$ is finite. There exist two cases:

 $\bullet$ Each subspace  in $\{\nu(dg(\tilde x_i)-\1)| g\in\Gamma_i\setminus\{\1\}\}$
is the zero space $\{0\}$ and thus $x_i$ must be an isolated
singular point.

$\bullet$ There exist at least a nonzero space  in
$\{\nu(dg(\tilde x_i)-\1)| g\in\Gamma_i\setminus\{\1\}\}$. In this
case we may assume that
 $H_{i1},\cdots, H_{il_i}$ are all nontrivial subspaces consisting
 of the subspaces  in $\{\nu(dg(\tilde x_i)-\1)| g\in\Gamma_i\}$
and their intersections. If some $H_{is}$ has codimension more
than $r+1={\rm Ind}(S)+1$, by Proposition~\ref{prop:C.2} we can
choose a closed subspace of codimension $r+1$ containing it to
replace $H_{is}$. Let us also denote this subspace by $H_{is}$.
Then we can assume that
\begin{equation}\label{e:2.102}
 {\rm codim}(H_{i1})\le {\rm codim}(H_{i2})\le\cdots\le {\rm
 codim}(H_{il_i})\le r+1.
 \end{equation}
Clearly, $\cap^{l_i}_{s=1}H_{is}=\{0\}$. Set $\widetilde P_{is}=
F_i\bigl(N(0_{\tilde x_i})\cap H_{is}\bigr)$, $s=1,\cdots, l_i$.
Then each $\widetilde P_{is}$ is a (relatively closed connected)
Banach submanifold of $\widetilde W_i$ which contains $\tilde x_i$
as a (relative) interior point, and has also the same codimension
as $H_{is}$, $s=1,\cdots, l_i$.
 Moreover it also holds that
\begin{equation}\label{e:2.103}
\left.\begin{array}{ll}
 \pi_i^{-1}(X^{sing}\cap W_i)\subset\cup^{l_i}_{s=1}\widetilde
 P_{is}\quad{\rm and}\quad\cap^{l_i}_{s=1}\widetilde P_{is}=\{\tilde x_i\}\\
 T_{\tilde x_i}\widetilde
P_{is}=H_{is},\;s=1,\cdots, l_i.
 \end{array}\right\}
 \end{equation}

 It follows from
Proposition~\ref{prop:C.1} that the restriction of the section
$\widetilde S_i:\widetilde W_i\to\widetilde E_i$ in (\ref{e:2.75})
to each $\widetilde P_{is}$ is also Fredholm and has  index
$r-{\rm codim}\widetilde P_{is}$. Moreover under the
trivialization naturally induced by ${\cal T}_i$ in (\ref{e:2.74})
these restrictions have trivialization representatives
$$
\widetilde S_i^T|_{\widetilde P_{is}}: \widetilde P_{is}\to
(\widetilde E_i)_{\tilde x_i},\;s=1,\cdots, n_i.
$$
These are Fredholm. By increasing $m_i$ in (\ref{e:2.76}) and
shrinking $W_i$ we may also require that
$$
d\widetilde S_i^T|_{\widetilde P_{is}}(\tilde x_i)(T_{\tilde
x_i}\widetilde P_{is})+ {\rm span}(\{v_{i1},\cdots,
v_{im_i}\})=(\widetilde E_i)_{\tilde x_i}
$$
and furthermore that for $i=1,\cdots, n$ and $s=1,\cdots, n_i$
\begin{equation}\label{e:2.104}
\left.\begin{array}{ll}
 d\widetilde S_i^T|_{\widetilde
P_{is}}(\tilde y)(T_{\tilde y}\widetilde P_{is})+ {\rm
span}(\{v_{i1},\cdots,
v_{im_i}\})=(\widetilde E_i)_{\tilde x_i},\\
d\widetilde S_i^T(\tilde x)(T_{\tilde x}\widetilde P_{is})+ {\rm
span}(\{v_{i1},\cdots, v_{im_i}\})=(\widetilde E_i)_{\tilde
x_i}\\
\qquad\forall \tilde y\in\widetilde W_i\cap \widetilde
P_{is}\;{\rm and}\;\forall \tilde x\in\widetilde W_i.
\end{array}\right\}
\end{equation}
In this case the sections $\tilde\sigma_{ij}$ (and $\tilde
s_{ij}$) obtained as in (\ref{e:2.79}) are called {\bf optimal}
(for $(X, E, S)$). From these it follows that if $\tilde u_{i_q}$
in (\ref{e:2.92}) belongs to some $\widetilde P_{i_qs}$ then the
tangent map of the restriction of $\Upsilon^T_{Iq}$ at $(\tilde
u_{i_q}, 0)$,
\begin{eqnarray*}
&&d\Upsilon_{Iq}^T(\tilde u_{i_q}, 0):T_{\tilde u_{i_q}}\widetilde
H_{i_qs}\times\R^m\to (\widetilde
E_{i_q})_{\tilde x_{i_q}},\nonumber\\
&&\hspace{-5mm}(\tilde \xi,{\bf v})\mapsto
d(S^T_{i_q}|_{\widetilde H_{i_qs}})(\tilde u_{i_q})(\tilde\xi)+
\sum^{m_{i_q}}_{j=1}v_{i_qj} \sigma^T_{i_qj}(\tilde u_{i_q})+
\sum_{i\in
I\setminus\{i_q\}}\sum^{m_i}_{j=1}v_{ij}\tau^{qT}_{ij}(\tilde
u_{i_q})\nonumber
\end{eqnarray*}
is also surjective. So by shrinking $W^\ast$ and $\varepsilon>0$ in
Theorem~\ref{th:2.38} we get that for each $I\in{\cal N}$, $i\in I$
and $s=1,\cdots, l_i$ the restriction of the section $\Upsilon_I$ to
each  Banach submanifold
$$
\hat\pi_I^{-1}(\widehat V_I^{\ast}\cap P_{is})\times
B_\varepsilon(\R^m)\subset\widehat V_I^\ast\times
B_\varepsilon(\R^m)
$$
is also transversal to the zero section, where
$P_{is}=\pi_i(\widetilde P_{is})$. In this case we say the section
$\Upsilon_I$ to be {\bf strongly transversal} to the zero section.
It follows that for generic ${\bf t}\in B_\varepsilon(\R^m)$ the
restriction of the section $\Upsilon_I^{({\bf t})}$ in
Corollary~\ref{cor:2.39} to each Banach submanifold
$\hat\pi_I^{-1}(\widehat V_I^{\ast}\cap P_{is})\subset\widehat
V_I^{\ast}$ is also transversal to the zero section. Thus the
submanifold
$$\hat\pi_I^{-1}(P_{is})\cap\widehat{\cal
M}^{\bf t}_I(S)=\hat\pi_I^{-1}(\widehat V_I^{\ast}\cap
P_{is})\cap\widehat{\cal M}^{\bf t}_I(S)\subset\widehat{\cal
M}_I(S)
$$ has codimension ${\rm
codim}P_{is}$. Set
\begin{equation}\label{e:2.105}
 \widehat{\cal M}^{\bf
t}_I(S)^{v}=\cup_{i\in I}\cup^{l_i}_{s=1}\hat\pi_I^{-1}(\widehat
V_I^{\ast}\cap P_{is})\cap\widehat{\cal M}^{\bf t}_I(S).
\end{equation}
It is a subvariety in $\widehat{\cal M}^{\bf
t}_I(S)$ and also contains $\widehat{\cal M}^{\bf t}_I(S)^{sing}$.

Furthermore, carefully checking the proof of
Proposition~\ref{prop:2.40}(ii) one easily sees that for each
$I\in{\cal N}$, $i\in I$ and $s=1,\cdots, l_i$, the restriction of
the section $\Psi_I$ in (\ref{e:2.97}) to $\hat\pi_I^{-1}(\widehat
V_I^{\ast}\cap P_{is})\times{\cal P}^l({\bf t},{\bf t}')\times
[0,1]$ is not only Fredholm but also transversal to the zero
section. Since $Z(S)$ is compact and ${\cal N}$ is finite it follows
that for a generic path $\gamma\in{\cal P}^l({\bf t},{\bf t}')$ the
restriction of the section $\Psi_I^\gamma$ in (\ref{e:2.99}) to each
$\hat\pi_I^{-1}(\widehat V_I^{\ast}\cap P_{is})\times [0,1]$ is
Fredholm and transversal to the zero section. Denote this
restriction by $\Psi_{Iis}^\gamma$. Then the smooth manifold
$(\Psi_{Iis}^\gamma)^{-1}(0)$ has codimension ${\rm codim}P_{is}$ in
$(\Psi_{I}^\gamma)^{-1}(0)$ and also boundary
$$
\bigl(\hat\pi_I^{-1}(\widehat V_I^{\ast}\cap
P_{is})\cap\widehat{\cal M}^{\bf t}_I(S)\bigr)\bigcup
\bigl(-\hat\pi_I^{-1}(\widehat V_I^{\ast}\cap
P_{is})\cap\widehat{\cal M}^{\bf t'}_I(S)\bigr).
$$
\footnote{The negative signs here and in the following
(\ref{e:2.108}) are to be understood as in the footnote to
Proposition~\ref{prop:2.40}(ii).} Here $(\Psi^\gamma_I)^{-1}(0)$ is
as in (\ref{e:2.100}). Denote by
\begin{equation}\label{e:2.106}
(\Psi^\gamma_I)^{-1}(0)^v=\bigl(\cup_{i\in
I}\cup^{l_i}_{s=1}\hat\pi_I^{-1}(\widehat V_I^{\ast}\cap
P_{is})\times [0, 1]\bigr)\cap(\Psi^\gamma_I)^{-1}(0).
\end{equation}
It is a subvariety in $(\Psi^\gamma_I)^{-1}(0)$ and contains all
singular points of  $(\Psi^\gamma_I)^{-1}(0)$. Moreover, if we
understand the boundary of $(\Psi^\gamma_I)^{-1}(0)^v$ to be the
union of all branch submanifolds, i.e.,
\begin{equation}\label{e:2.107}
\partial(\Psi^\gamma_I)^{-1}(0)^v=\bigcup_{i\in
I}\bigcup^{l_i}_{s=1}\partial\bigl((\hat\pi_I^{-1}(\widehat
V_I^{\ast}\cap P_{is})\times [0,
1])\cap(\Psi^\gamma_I)^{-1}(0)\bigl).
\end{equation}
Then it holds that
\begin{equation}\label{e:2.108}
\partial(\Psi^\gamma_I)^{-1}(0)^v=\widehat{\cal M}^{\bf t}_I(S)^{v}\cup
(-\widehat{\cal M}^{\bf t'}_I(S)^{v}).
\end{equation}
Summarizing the above arguments we get:

\begin{proposition}\label{prop:2.42}
There exists a residual subset $B_\varepsilon(\R^m)_{res}\subset
B_\varepsilon(\R^m)$ such that for each ${\bf t}\in
B_\varepsilon(\R^m)_{res}$ the family $\widehat{\cal M}^{\bf t}(S)$
in Corollary~\ref{cor:2.39} can be chosen to have also properties:
For each $I\in{\cal N}$ there exists a subvariety $\widehat{\cal
M}^{\bf t}_I(S)^{v}\subset\widehat{\cal M}^{\bf t}_I(S)$ such that
\begin{description}
\item[(i)] $\widehat{\cal M}^{\bf
t}_I(S)^{sing}\subset\widehat{\cal M}^{\bf t}_I(S)^{v}$,

 \item[(ii)] each branch of $\widehat{\cal M}^{\bf
t}_I(S)^{v}$ has dimension lower than $r$ and thus $\widehat{\cal
M}^{\bf t}_I(S)^\circ=\widehat{\cal M}^{\bf
t}_I(S)\setminus\widehat{\cal M}^{\bf t}_I(S)^{sing}\supseteq
\widehat{\cal M}^{\bf t}_I(S)\setminus\widehat{\cal M}^{\bf
t}_I(S)^{v}$ is an open and dense submanifold in $\widehat{\cal
M}^{\bf t}_I(S)$,

\item[(iii)] if each $x\in X^{sing}\cap Z(S)$  has singular
codimension $r'\ge 1$  then $\widehat{\cal M}^{\bf t}_I(S)^{v}$
has no branch submanifolds of codimension $r'\ge 1$,

\item[(iv)] the family $\widehat{\cal M}^{\bf
t}(S)^v=\{\widehat{\cal M}^{\bf t}_I(S)^v\,|\, I\in{\cal N}\}$ is
compatible in the sense that that for any $J\subset I\in{\cal N}$,
$$
 \hat\pi^I_J\bigl((\hat
\pi^I_J)^{-1}(\widehat V_J^\ast)\cap\widehat{\cal M}^{\bf
t}_I(S)^v\bigr)={\rm Im}(\hat\pi^I_J)\cap\widehat{\cal M}^{\bf
t}_J(S)^v,
$$

\item[(v)] $\cup_{I\in{\cal N}}\hat\pi_I\bigr(\widehat{\cal
M}^{\bf t}_I(S)^v\bigr)$ is a compact subset in $\cup_{I\in{\cal
N}}\hat\pi_I\bigr(\widehat{\cal M}^{\bf t}_I(S)\bigr)$,

\item[(vi)] for any two ${\bf t},\, {\bf t}'\in
B_\varepsilon(\R^m)_{res}$  there exist generic paths
$\gamma:[0,1]\to B_\varepsilon(\R^m)$ with $\gamma(0)={\bf t}$ and
$\gamma(1)={\bf t}'$ such that for each $I\in{\cal N}$ there exists
a subvariety $\partial(\Psi^\gamma_I)^{-1}(0)^v$ of the manifold
$\partial(\Psi^\gamma_I)^{-1}(0)$ in Proposition~\ref{prop:2.40}(ii)
such that:

$\bullet$ containing all singular points of
$\partial(\Psi^\gamma_I)^{-1}(0)$,

$\bullet$ (\ref{e:2.108}) holding, i.e.,
$\partial(\Psi^\gamma_I)^{-1}(0)^v=\widehat{\cal M}^{\bf
t}_I(S)^{v}\cup (-\widehat{\cal M}^{\bf t'}_I(S)^{v})$,

$\bullet$ satisfying the obvious properties corresponding  with
(ii)-(v) above.
\end{description}
\end{proposition}

Later we call  $(\Psi^\gamma)^{-1}(0)=\{(\Psi^\gamma_I)^{-1}(0)|
I\in{\cal N}\}$ as in (vi) a  {\bf strong cobordsim} between
$\widehat{\cal M}^{\bf t}(S)=\{\widehat{\cal M}^{\bf t}_I(S)|
I\in{\cal N}\}$ and $\widehat{\cal M}^{\bf t'}(S)=\{\widehat{\cal
M}^{\bf t'}_I(S)| I\in{\cal N}\}$.

\begin{remark}\label{rem:2.43}
{\rm (i) If $X$ has no boundary (or more generally
$Z(S)\subset{\rm Int}(X)$) then (smooth) manifold
$(\Psi^\gamma_I)^{-1}(0)$ in (\ref{e:2.100}) has only boundary but
no corner for $\varepsilon>$ sufficiently small. In this case the
family $\widehat{\cal M}^{\bf t}(S)=\bigl\{\widehat{\cal M}^{\bf
t}_I(S)\,|\, I\in{\cal N}\bigr\}$ is ``like'' an open cover of a
closed manifold. If $Z(S)\cap\partial X\ne\emptyset$ then the
manifold $(\Psi^\gamma_I)^{-1}(0)$ in (\ref{e:2.100}) might has
only boundary but
corner.\\
(ii) If $S$ has been transversal to the zero section so is the
section $\widehat S=\{\widehat S_I:\,I\in{\cal N}\}$. In this case
$S^{-1}(0)$ is a compact orbifold of dimension ${\rm Ind}(S)$ and
${\bf t}={\bf 0}\in B_\varepsilon(\R^m)_{res}$. It follows from
Proposition~\ref{prop:2.42} that $\cup_{I\in{\mathcal N}}
\hat\pi_I(\widehat{\mathcal M}_I^{{\bf 0}}(S))=S^{-1}(0)$ and that
$\widehat{\mathcal M}^{{\bf 0}}(S)=\{\widehat{\mathcal M}_I^{{\bf
0}}(S):\,I\in{\cal N}\}$ is cobordant to $\widehat{\mathcal M}^{{\bf
t}}(S)=\{\widehat{\mathcal M}_I^{{\bf t}}(S):\,I\in{\cal N}\}$ for
any ${\bf t}\in B_\varepsilon(\R^m)_{res}$. Actually,
$\widehat{\mathcal M}^{{\bf t}}(S)$ is exactly a resolution of the
orbifold
$S^{-1}(0)$.\\
(iii) Sometimes it is useful and convenient to write the
compatible family $\widehat{\mathcal M}^{{\bf t}}(S)=\{
\widehat{\mathcal M}_I^{{\bf t}}(S)\,|\, I\in{\cal N}\}$ in the
tight version.  Following \cite{Mc1} consider the disjoint union
$\coprod_{I\in{\cal N}}\widehat{\mathcal M}_I^{{\bf t}}(S)$ and
let $\sim$ be the equivalence relation that is generated by
setting $\hat y_I\sim\hat z_J$ if $J\subset I$ and
$\hat\pi^I_J(\hat y_I)=\hat z_J$. Denote by
$$
\overline{\mathcal M}^{{\bf t}}(S):= \coprod_{I\in{\cal
N}}\widehat{\mathcal M}_I^{{\bf t}}(S)/\sim
$$
and by $\hat q_I:\widehat{\mathcal M}_I^{{\bf t}}(S)\to
\overline{\mathcal M}^{{\bf t}}(S)$ the restriction of the obvious
quotient map to $\widehat{\mathcal M}_I^{{\bf t}}(S)$. Let us
write  $\overline{\mathcal M}_I^{{\bf t}}(S):=\hat
q_I(\widehat{\mathcal M}_I^{{\bf t}}(S))$. In general, $\hat q_I$
is not a homeomorphism from $\widehat{\mathcal M}_I^{{\bf t}}(S)$
onto $\overline{\mathcal M}_I^{{\bf t}}(S)$.  We also write
$$
\overline{\mathcal M}^{{\bf t}}(S)^\circ:= \coprod_{I\in{\cal
N}}\widehat{\mathcal M}_I^{{\bf t}}(S)^\circ/\sim\quad{\rm
and}\quad\overline{\mathcal M}_I^{{\bf t}}(S)^\circ=\hat
q_I(\widehat{\mathcal M}_I^{{\bf t}}(S)^\circ)
$$
for each $I\in{\cal N}$. Set
$$\lambda_I:\overline{\mathcal M}_I^{{\bf t}}(S)\to\Q,\;\bar x\mapsto\frac{|\hat
q_I^{-1}(\bar x)|}{|\Gamma_I|-|(\Gamma_I)_{\hat x_I}|+1},
$$
where $|\hat q_I^{-1}(\bar x)|$ is the number of elements in the
set $\hat q_I^{-1}(\bar x)$ and $\hat x_I$ is any element in $\hat
q_I^{-1}(\bar x)$. If $J\subset I$ and $\bar x\in
\overline{\mathcal M}_I^{{\bf t}}(S)\cap \overline{\mathcal
M}_J^{{\bf t}}(S)$ it is easily checked that $\lambda_I(\bar
x)=\lambda_J(\bar x)$. So all these $\lambda_I$ give a positive
rational function on $\overline{\mathcal M}^{{\bf t}}(S)$, denoted
by $\lambda$, is called the {\bf label function}. In particular,
$\lambda(\bar x)=|\hat q_I^{-1}(\bar x)|/|\Gamma_I|$ for any $\bar
x\in\overline{\mathcal M}_I^{{\bf t}}(S)^\circ$. All $\hat\pi_I$
can be glued into a natural map
$$
\hat\pi:\overline{\cal M}^{\bf t}(S)\to X.
$$
Its restriction on $\overline{\cal M}^{\bf
t}(S)^\circ\subset\overline{\cal M}^{\bf t}(S)$ is also denoted by
$\hat\pi$ when there is no danger of confusion.}
\end{remark}

Later we directly write $(\hat \pi^I_J)^{-1}(\widehat
V_J^\ast)\cap\widehat{\cal M}^{\bf t}_I(S)$ as $(\hat
\pi^I_J)^{-1}(\widehat{\cal M}^{\bf t}_J(S))$ when there is no
danger of confusion. Moreover, we always assume that the family
$\widehat{\cal M}^{\bf t}(S)=\bigl\{\widehat{\cal M}^{\bf
t}_I(S)\,|\, I\in{\cal N}\bigr\}$ also satisfies
Proposition~\ref{prop:2.42}.

\subsection{Orientation}\label{sec:2.5}

We assume that $X$ is as in Section~\ref{sec:2.4}. Now we begin to
consider the orientation of the family of manifolds $\widehat{\cal
M}^{\bf t}(S)=\bigl\{\widehat{\cal M}^{\bf t}_I(S)\,|\, I\in{\cal
N}\bigr\}$. By Proposition~\ref{prop:2.42} each $\widehat{\cal
M}^{\bf t}_I(S)^\circ$ is an open dense subset of $\widehat{\cal
M}^{\bf t}_I(S)$. Moreover for any $J\subset I\in{\cal N}$, by
(\ref{e:2.96}) the projection $\hat\pi^I_J$ restricts to a
$|\Gamma_I|/|\Gamma_J|$-fold regular smooth covering
$$
(\hat \pi^I_J)^{-1}(\widehat V_J^\ast)\cap\widehat{\cal M}^{\bf
t}_I(S)^\circ\to {\rm Im}(\hat\pi^I_J)\cap\widehat{\cal M}^{\bf
t}_J(S)^\circ.
$$
Thus we may say it to be {\bf orientation preserving} provided that
both $\widehat{\cal M}^{\bf t}_I(S)$ and $\widehat{\cal M}^{\bf
t}_J(S)$ are oriented. Later, the statement that $\hat\pi^I_J$ is
orientation preserving should always be understood in the way. We
call the family $\widehat{\cal M}^{\bf t}(S)=\bigl\{\widehat{\cal
M}^{\bf t}_I(S)\,|\, I\in{\cal N}\bigr\}$ {\bf oriented} if each
$\widehat{\cal M}^{\bf t}_I(S)$ is oriented and each projection
$\hat\pi^I_J$, $J\subset I$, is orientation preserving.

\begin{proposition}\label{prop:2.44}
If the Banach Fredholm orbibundle $(X, E, S)$ is oriented then the
family $\widehat{\cal M}^{\bf t}(S)$ has an induced orientation.
\end{proposition}

\noindent{\bf Proof.}\quad We only need to prove that {\it
$\widehat S=\{\widehat S_I\,|\,I\in{\cal N}\}$ {\bf is oriented as
a Fredholm section of the system of Banach bundles $(\widehat{\cal
E},\widehat W)$}}. This means that each determinant line bundle
${\rm det}(D\widehat S_I)\to\widehat W_I$ is given an orientation,
i.e. a continuous nowhere zero section, and that for any
$I\in{\cal N}$ and $J\subset I$ the bundle map
$$
{\rm det}(D\widehat S_I)|_{(\hat\pi^I_J)^{-1}(\widehat
W_I^\circ)}\to {\rm det}(D\widehat S_J)|_{\widehat W_J^\circ}
$$
induced by the projection $(\hat\pi^I_J,\hat\Pi^I_J)$, which  maps
a fiber at $\hat x_I\in (\hat\pi^I_J)^{-1}(\widehat W_I^\circ)$ to
that of $\hat\pi^I_J(\hat x_I)\in \widehat W_J^\circ$
isomorphically, is orientation preserving.

 By the definition below
Definition~\ref{def:2.11} let $\tilde S_i:\widetilde W_i\to
\widetilde E_i$ be the $\Gamma_i$-equivariant lift of $S|_{W_i}$,
$i=1,\cdots,n$. Then for any $u\in W_i\cap W_j$, $\tilde
u_i=\pi_i^{-1}(u)$ and $\tilde u_j=\pi_j^{-1}(u)$ there exist a
connected open neighborhood $O$ of $u$ in $W_i\cap W_j$, a group
isomorphism ${\cal A}_{ij}: \Gamma(\tilde u_i)\to\Gamma(\tilde
u_j)$, a ${\cal A}_{ij}$-equivariant diffeomorphism
$\lambda_{ij}:\widetilde O_i=\pi_i^{-1}(O)\to\widetilde
O_j=\pi_j^{-1}(O)$ and a bundle isomorphism
 $\Lambda_{ij}:\widetilde E_i|_{\widetilde O_i}\to
\widetilde E_j|_{\widetilde O_j}$ covering $\lambda_{ij}$ such
that
\begin{equation}\label{e:2.109}
\widetilde S_j\circ\lambda_{ij}=\Lambda_{ij}\circ\widetilde
S_i\quad{\rm on}\quad\widetilde O_i.
\end{equation}
That $S$ is oriented means that each $\tilde S_i$ is oriented,
i.e. the determinant bundle ${\rm det}(D\tilde S_i)$ over
$Z(\tilde S_i)$ is given a continuous nowhere zero section
$o(\tilde S_i)$, and that for every pair $i, j=1,\cdots,n$ there
exist positive continuous functions $f_{ij}:\widetilde O_i\to\R$,
 such that for any $\tilde x\in\widetilde O_i$,
$$
 {\rm
det}(\Lambda_{ij})\bigl(o(\widetilde S_i)(\tilde
x)\bigr)=f_{ij}(\lambda_{ij}(\tilde x))\cdot o(\widetilde
S_j)\bigl(\lambda_{ij}(\tilde x)\bigr).
$$

 We first show that they induce such a section of
${\rm det}(D\widehat S_I)$ for each $I\in{\cal N}$ with $|I|>1$.
Here
 $\widehat S_I:\widehat W_I\to \widehat E_I$ is defined by (\ref{e:2.62}). Let
$I=\{i_1,\cdots,i_k\}\in{\cal N}$ with $k>1$. For any $\tilde
u_I\in\widetilde W_I$ and $u_I=\pi_I(\tilde u_I)$ let $\widetilde
O_{i_l}=\pi_{i_l}^{-1}(O)$, $\lambda_{i_qi_l}$ and ${\cal
A}_{i_qi_l}:\Gamma(\tilde u_{i_q})\to\Gamma(\tilde u_{i_l})$ be as
in the proof of Lemma~\ref{lem:2.19}. Let $o(\tilde S_i)$ be the
continuous nowhere zero section of ${\rm det}(D\tilde S_i)$ defining
the orientation of $\tilde S_i$. For any two $i_q, i_l\in I$, by the
bundle isomorphism $\Lambda_{i_qi_l}$ above (\ref{e:2.17}) and the
first paragraph there exist positive continuous functions
$f_{i_li_q}:\widetilde O_{i_l}\to\R$ such that
$$
{\rm det}(\Lambda_{i_qi_l})\bigl(o(\tilde S_{i_q})(\tilde
x)\bigr)=f_{i_li_q}(\lambda_{i_qi_l}(\tilde x))\cdot o(\tilde
S_{i_l})\bigl(\lambda_{i_qi_l}(\tilde x)\bigr)
$$
for any $\tilde x\in\widetilde O_{i_q}$. For the section in
(\ref{e:2.63}),
$$
\widehat S_I:\widehat W_I\to \widehat E_I, [\tilde u_I,\phi_I,
\tilde x_I]\mapsto \langle\tilde u_I,\phi_I, \tilde S_I(\tilde
x_I)\rangle,
$$
it may be expressed as
\begin{equation}\label{e:2.110}
\widetilde O_{i_q}\to\widetilde E_{i_q}|_{\widetilde
O_{i_q}},\;\tilde x\mapsto\widetilde S_{i_q}(\tilde x)
\end{equation}
 in the charts $\widehat{\phi_I\circ\lambda_I^q}$ in (\ref{e:2.87}) and
 $\widehat{\Phi_I\circ\Lambda_I^q}$ in (\ref{e:2.88}), $q=1,\cdots,k$.
So the orientation of $\widetilde S_{i_q}$ naturally gives one of
$\widehat S_I|_{\widehat{\phi_I\circ\lambda_I^q}(\widetilde
O_{i_q})}$. What we need to prove is that these orientations agree
on overlaps.

Let $\widehat W_I^{\circ}$ be as in (\ref{e:2.45}).  Then it
suffices to check that those orientations agree on overlaps near any
$\hat u_I\in \widehat W_I^{\circ}$. In this case we may require that
$\widetilde O_{i_q}\subset\widetilde W_{i_q}^\circ$, and thus the
corresponding charts in (\ref{e:2.87}) and (\ref{e:2.88}) may be
written as
\begin{eqnarray*}
&&\widehat{\lambda_I^q}:\widetilde O_{i_q}\to \widehat
W_I^{\circ},\;\tilde x\mapsto [\tilde u_I, 1_I,
\lambda_I^q(\tilde x)],\\
&&\widehat{\Lambda_I^q}:\widetilde E_{i_q}|_{\widetilde
O_{i_q}}\to \widehat E_I,\;\tilde \xi\mapsto \langle\tilde u_I,
1_I, \Lambda_I^q(\tilde \xi)\rangle
\end{eqnarray*}
respectively. Here $\lambda_I^q=(\lambda_{i_qi})_{i\in I}$ with
 $\lambda_{i_qi_q}=id_{\widetilde O_{i_q}}$,
  $\Lambda_I^q=(\Lambda_{i_qi})_{i\in I}$
and each $\Lambda_{i_qi}$ is the lifting of $\lambda_{i_qi}$ as
before and thus   $\Lambda_{i_qi_q}=id_{\widetilde
E_{i_q}|_{\widetilde O_{i_q}}}$. So we have transition functions
\begin{eqnarray*}
&&(\widehat{\lambda_I^l})^{-1}\circ\widehat{\lambda_I^q}:\widetilde
O_{i_q}\to \widetilde O_{i_l},\;\tilde x\mapsto
\lambda_{i_qi_l}(\tilde x),\\
&&(\widehat{\Lambda_I^l})^{-1}\circ\widehat{\Lambda_I^q}:\widetilde
E_{i_q}|_{\widetilde O_{i_q}}\to \widetilde E_{i_l}|_{\widetilde
O_{i_l}},\;\tilde \xi\mapsto \Lambda_{i_qi_l}(\tilde \xi).
\end{eqnarray*}
This, (\ref{e:2.109}) and (\ref{e:2.110}) for $i=i_q, j=i_l$
together imply that the orientation of $\widehat
S_I|_{\widehat{\phi_I\circ\lambda_I^1}(\widetilde O_{i_1})}$ given
by the above one of $\widetilde S_{i_1}$ agree with (on their
overlaps) those of $\widehat
S_I|_{\widehat{\phi_I\circ\lambda_I^q}(\widetilde O_{i_q})}$ given
by the ones of $\widetilde S_{i_q}$, $q=1,\cdots,n$.

 Moreover, for $\hat u_I'\in\hat\pi_I^{-1}(O)$ and a small connected open
neighborhood $O'$ of $u'_I$ in $W_I^{\circ}=\hat\pi_I(\widehat
W_I^{\circ})$ and $\widetilde O'_{i_l}=\pi_{i_l}^{-1}(O')$,
$l=1,\cdots,k$, we have the corresponding charts
\begin{eqnarray*}
&&\widehat{\lambda'^l_I}:\widetilde O'_{i_l}\to \widehat
W_I,\;\tilde x\mapsto [\tilde u_I, 1_I,
\lambda'^l_I(\tilde x)],\\
&&\widehat{\Lambda'^l_I}:\widetilde E_{i_l}|_{\widetilde
O'_{i_l}}\to \widehat E_I,\;\tilde \xi\mapsto \langle\tilde u_I,
1_I, \Lambda'^l_I(\tilde \xi)\rangle.
\end{eqnarray*}
Note that $\pi_{i_q}^{-1}(O'\cap O)=\widetilde
O'_{i_q}\cap\widetilde O_{i_q}$ and that both $\lambda_{i_qi_l}$
both $\lambda'_{i_qi_l}$ are diffeomorphisms from $\widetilde
O'_{i_q}\cap\widetilde O_{i_q}$ to $\widetilde
O'_{i_l}\cap\widetilde O_{i_l}$, $q, l=1,\cdots,k$. Using the fact
that $\pi_{i_l}\circ\lambda_{i_qi_l}=\pi_{i_q}$ and
$\pi_{i_l}\circ\lambda'_{i_qi_l}=\pi_{i_q}$ we derive that
$$
\lambda_{i_qi_l}=\lambda'_{i_qi_l}\;{\rm on}\;\widetilde
O'_{i_q}\cap\widetilde O_{i_q}\quad{\rm and}\quad
\Lambda_{i_qi_l}=\Lambda'_{i_qi_l}\;{\rm on}\;\widetilde
E_{i_q}|_{\widetilde O'_{i_q}\cap\widetilde O_{i_q}}
$$
 for any $q,l=1,\cdots,k$. In particular, we have
\begin{eqnarray*}
&&(\widehat{\lambda'^1_I})^{-1}\circ\widehat{\lambda^1_I}:\widetilde
O'_{i_1}\cap\widetilde O_{i_1}\to \widetilde
O'_{i_1}\cap\widetilde O_{i_1},\;\tilde x\mapsto
\tilde x,\\
&&(\widehat{\Lambda'^1_I})^{-1}\circ\widehat{\Lambda^1_I}:\widetilde
E_{i_1}|_{\widetilde O'_{i_1}\cap\widetilde O_{i_1}}\to \widetilde
E_{i_1}|_{\widetilde O'_{i_1}\cap\widetilde O_{i_1}},\;\tilde
\xi\mapsto \tilde \xi.
\end{eqnarray*}
It follows from these and the above arguments that the orientation
of $\widehat S_I$ defined as above is independent of choices of
the coordinate charts.

Finally we show that for any $I, J\in{\cal N}$ with $J\subset I$ the
covering projection $\hat\Pi^I_J$ preserves the orientation that we
just defined. In fact, using the local expression of $\widehat S_I$
in (\ref{e:2.52}) the relation in (\ref{e:2.64}) may be written as
$$
\tilde S_{j_1}\circ\lambda_{i_1j_1}=\Lambda_{i_1j_1}\circ\tilde
S_{i_1}.
$$
So it follows from  (\ref{e:2.110}) that
 $\hat\Pi^I_J$ is orientation preserving for any $J\subset I\in{\cal N}$.
 That is, there exists a positive continuous function
$f_{IJ}:\widehat W_I\to\R$ such that
$$
{\rm
det}(\hat\Pi^I_J)\bigl(o(\widehat S_I)(\hat
x)\bigr)=f_{IJ}(\hat\pi^I_J(\hat
 x))\cdot o(\widehat S_J)\bigl(\pi^I_J(\hat x)\bigr)
$$
for any $\hat x\in (\hat\pi^I_J)^{-1}(\widehat
W_J)$. The proof is completed. \hfill$\Box$\vspace{2mm}

\subsection{Virtual Euler cycle and class}\label{sec:2.6}

Assume that $X$ is as in Section~\ref{sec:2.5}. Let $(X, E, S)$ be
an oriented Banach Fredholm bundle of index $r$ and with compact
zero locus $Z(S)$. Then we have a family of compatible oriented
manifolds $\widehat{\cal M}^{\bf t}(S)=\bigl\{\widehat{\cal M}^{\bf
t}_I(S)\,:\, I\in{\cal N}\bigr\}$, ${\bf t}\in
B_\varepsilon(\R^m)_{reg}$ to satisfy Corollary~\ref{cor:2.39},
Proposition~\ref{prop:2.40}, Proposition~\ref{prop:2.42} and
Proposition~\ref{prop:2.44}. However, in order to use  the family
$\widehat{\cal M}^{\bf t}(S)$ to construct a cycle we also assume
that $(X, E, S)$  satisfies the following:

\begin{assumption}\label{ass:2.45}
{\rm  Each $x\in X^{sing}\cap Z(S)$ has singularity of codimension
at least two. More precisely, for some Banach orbifold chart
$(\widetilde W, \Gamma,\pi)$ at centred at $x$, i.e.
$\pi^{-1}(x)=\{\tilde x\}$ and $\Gamma=\Gamma_{\tilde x}$, if $1$
is an eigenvalue of the linearization of any
$g\in\Gamma\setminus\{\1\}$, $dg(\tilde x): T_{\tilde x}\widetilde
W\to T_{\tilde x}\widetilde W$ then the corresponding eigenspace
${\rm Ker}(dg(\tilde x)-\1)$ has codimension at least two.}
\end{assumption}

Clearly this assumption is weaker than the condition that $X$ is
$2$-regular near $Z(S)$. The Assumption~\ref{ass:2.45} is
automatically satisfied if $X$ is either an oriented orbifold of
finite dimension or a
 complex Banach orbifold. On the other hand there exist
 examples of oriented Banach Fredholm bundles not to satisfy
 Assumption~\ref{ass:2.45}. For example, let $\widetilde X=\R\times
 S^n$, $\widetilde E=\widetilde X\times\R$ and $\Gamma={\Z}_2=\{\1, e\}$
 acts on them by
$$e\cdot (t, x)=(-t, x)\quad{\rm and}\quad e\cdot (t, x; t')=(-t,
x; -t').
$$
Let the section $\widetilde S:\widetilde X\to\widetilde E$ be
given by $(t, x)\mapsto (t, x; t)$. Then $Z(\widetilde
S)=\{0\}\times S^n$ and $X^{sing}=\{0\}\times S^n$. Let
$X=\widetilde X/\Gamma$, $E=\widetilde E/\Gamma$ and $S:X\to E$ be
the natural descent of $\widetilde S$. One can check that $(X, E,
S)$ is an orientable Banach Fredholm bundle with compact zero
locus. However, each $x\in X^{sing}\cap Z(S)$ has singularity of
codimension one.

Now under Assumption~\ref{ass:2.45} we can require that each
$\widehat{\cal M}^{\bf t}_I(S)^\circ$ is an open and dense in
$\widehat{\cal M}^{\bf t}_I(S)$ and that the singular set
$\widehat{\cal M}^{\bf t}_I(S)^{sing}$ is contained in a subvariety
of $\widehat{\cal M}^{\bf t}_I(S)$ of codimension two. It follows
that for any $J\subset I\in{\cal N}$ the map
$$
\hat\pi^I_J: (\hat \pi^I_J)^{-1}\bigl(\widehat{\cal M}^{\bf
t}_J(S)^\circ\bigr)\to {\rm Im}(\hat\pi^I_J)\subset\widehat{\cal
M}^{\bf t}_J(S)^\circ
$$
is a $|\Gamma_I|/|\Gamma_J|$-fold regular smooth covering
preserving orientation. As in \cite{LiuT1}, for each $I\in{\cal
N}$ we take an increasing sequence of compact submanifolds with
boundary of $\widehat{\cal M}^{\bf t}_I(S)^\circ$, $\{\widehat
K_I^p\}^\infty_{p=1}$, and choose a triangulation of each
$\widehat K_I^p$ such that:
\begin{description}
\item[(i)] $\widehat{\cal M}^{\bf
t}_I(S)^\circ=\cup^\infty_{p=1}\widehat K_I^p$.

\item[(ii)] For each simplex $c^p_I\in\widehat K^p_I$,
$\hat\pi^I_J(c^p_I)$ is a simplex of $\widehat K^p_J$ provided that
$J\subset I$ and $\widehat K^p_J\cap \hat\pi^I_J(\widehat
K^p_I)\ne\emptyset$.

\item[(iii)] If $q>p$, $\widehat K^q_I|_{\widehat K^p_I}$ is
obtained from $\widehat K^p_I$ by some divisions as simplicial
complex.
\end{description}
For each simplex in $\widehat K^p_I$ of dimension $r={\rm Ind}(S)$
we give it an orientation induced from that of $\widehat{\cal
M}^{\bf t}_I(S)$, and then take the summation of them to get an
integral chain in $\widehat{\cal M}^{\bf t}_I(S)$, denoted by
$S(\widehat K^p_I)$. Let
$$\widehat S(\widehat K^p_I):=\frac{1}{|\Gamma_I|}\hat\pi_I\circ {
S(\widehat K^p_I)}
$$
be the corresponding  smooth singular rational chain in $X$. All
these $\widehat S(\widehat K^p_I)$, $I\in{\cal N}$, form a smooth
singular rational chain in $X$,
$$
C(\widehat K^p):=``\sum_{I\in{\cal N}}"\widehat S(\widehat K^p_I).
$$
  Here the summation  is put in
double quotation marks because we only count them once on those
overlaps where more than one singular maps appear. Let $e(E,
S)^{\bf t}$ denote the inverse limit of $\{C(\widehat K^p)\}_p$.
It is a rational singular chain in $X$ of dimension $r$. We also
formally write it as
\begin{equation}\label{e:2.111}
 e(E, S)^{\bf t}:=``\sum_{I\in{\cal
N}}"\frac{1}{|\Gamma_I|}\bigl\{\hat\pi_I:\widehat{\cal M}^{\bf
t}_I(S)^\circ\to X\bigr\}
\end{equation}
for any ${\bf t}\in B_\varepsilon(\R^m)_{res}$ if no confusion
occurs, where $\sum_{I\in{\cal N}}$ is put in double quotation
marks in (\ref{e:2.111})  because of the same reason as above. We
call $e(E, S)^{\bf t}$ a {\bf virtual Euler chain} of the triple
$(X, E, S)$.

If $X$ has boundary and $\partial X\cap Z(S)\ne\emptyset$, then it
is possible that
$$
\partial\widehat{\cal M}^{\bf t}_I(S)^\circ:=\widehat{\cal
M}^{\bf t}_I(S)^\circ\cap\hat\pi_I^{-1}(\partial X)\ne\emptyset.
$$
In this case,  for any $J\subset I$ the covering $\hat\pi^I_J$
also restricts to a regular one
$$
\hat\pi^I_J:(\hat \pi^I_J)^{-1}\bigl(\partial\widehat{\cal M}^{\bf
t}_J(S)^\circ\bigr)\to {\rm
Im}(\hat\pi^I_J)\cap\partial\widehat{\cal M}^{\bf t}_J(S)^\circ.
$$
We can require that the above triangulation is compatible with
these boundaries. Namely, the triangulation is an extension of a
given  triangulation of the family $\{\partial\widehat{\cal
M}^{\bf t}_J(S)^\circ\,|\, I\in{\cal N}\}$. Then the corresponding
 singular rational chain in $\partial X$ constructed from the
latter is exactly the boundary of $C(\widehat K^p_I)$,
$$
\partial C(\widehat K^p)=``\sum_{I\in{\cal N}}"\partial\widehat S(\widehat
K^p_I)=``\sum_{I\in{\cal
N}}"\frac{1}{|\Gamma_I|}\hat\pi_I\circ{\partial S(\widehat
K^p_I)}.
$$
We later write it as
$$
\partial e(E, S)^{\bf t}=``\sum_{I\in{\cal
N}}"\frac{1}{|\Gamma_I|}\bigl\{\hat\pi_I: \partial\widehat{\cal
M}^{\bf t}_I(S)^\circ\to \partial X\bigr\}
$$
without special statements. Consequently, we get:

\begin{claim}\label{cl:2.46}
If $X$ has no $2$-boundary $\partial_2 X$ then the class
$[\partial e(E, S)^{\bf t}]$ is zero in $H_{r-1}(X,\Q)$.
\end{claim}

Actually, for any smooth map $f$ from $X$ to a finite dimensional
smooth orbifold $M$ and any closed $(r-1)$-form $\omega$ on $M$
with compact support it directly follows from the Stokes theorem
that
\begin{eqnarray*}
\int_{\partial e(E, S)^{\bf
t}}f^\ast\omega\!\!\!\!\!\!&&:=``\sum_{I\in{\cal
N}}"\frac{1}{|\Gamma_I|}\int_{
\partial\widehat{\cal M}^{\bf t}_I(S)^\circ}(f\circ\hat\pi_I)^\ast\omega\\
&&=``\sum_{I\in{\cal N}}"\frac{1}{|\Gamma_I|}\int_{ \widehat{\cal
M}^{\bf t}_I(S)^\circ}(f\circ\hat\pi_I)^\ast(d\omega)=0.
\end{eqnarray*}
Here the properness of $\hat\pi_I$ is used.

Claim~\ref{e:2.43} is an analogue of the fact that for an
$n$-dimensional compact oriented manifold (resp. orbifold) $M$
with boundary the boundary $\partial M$ is an $(n-1)$-dimensional
manifold (resp. orbifold) with a natural orientation induced from
the orientation of $M$ and the image of the fundamental class
$[\partial M]$ is zero in $H_{n-1}(M,\Z)$ (resp. $H_{n-1}(M,\Q)$).

If $Z(S)\subset{\rm Int}(X)$ it follows from the above arguments
that for ${\bf t}$ sufficiently small, $\partial e(E, S)^{\bf
t}=\emptyset$ and thus  $e(E, S)^{\bf t}$ is a rational singular
cycle in $X$ of dimension $r$. In this case
 $e(E, S)^{\bf t}$ is called  a {\bf virtual Euler cycle}
of the triple $(X, E, S)$.  By Proposition~\ref{prop:2.40}(i)-(ii)
and Proposition~\ref{prop:2.42}(v)-(vi) it is easily seen that the
homology class of this cycle is independent of choices of generic
small ${\bf t}$.  Thus the homology class $[e(E, S)^{\bf t}]\in
H_r(X,\Q)$ is well-defined. We call the class
\begin{equation}\label{e:2.112}
e(E, S)=[e(E, S)^{\bf t}]\in H_r(X,\Q),
 \end{equation}
  the {\bf virtual Euler class} of the triple because we shall
  prove in next two subsections that it is independent of the different choices above.

In Remark~\ref{rem:2.43}, as singular chains $e(E, S)^{\bf t}$ can
be identified with $\hat\pi:\overline{\cal M}_I^{\bf t}(S)^\circ\to
X$. Later we shall also write $e(E, S)^{\bf
t}=\{\hat\pi:\overline{\cal M}_I^{\bf t}(S)^\circ\to X\}$.

 An important example of the
Banach Fredholm orbibundles is a triple consisting a finitely
dimensional orbibundle $p:E\to X$ and a smooth section $S:X\to E$.
As pointed out below (\ref{e:2.1}), if $X$ is also orientable it has
no singularities of codimension one naturally. In particular, if $X$
is a closed oriented orbifold we may take $E$ to be the tangent
orbibundle $TX$ and $S$ to be the zero section. The corresponding
Euler class $e(TX, 0)$ becomes a rational number, called the {\bf
orbifold Euler characteristic (number)} of $X$ and simply denoted by
$e(X)$.

 Let  $P$  be a smooth orbifold of finite
dimension and $f:X\to P$ be a smooth map. Then $f_\ast(e(E, S))\in
H_r(P,\Q)$, and $\langle f_\ast(e(E, S)),\alpha\rangle\in\Q$ for
any $\alpha\in H^r(P,\Q)$. Let $\alpha^\ast$ be a closed
differential form representative of it. As usual we have
$$
\langle \alpha, f_\ast(e(E,S))\rangle=\int_{f\circ e(E, S)^{\bf
t}}\alpha^\ast=\int_{e(E, S)^{\bf t}}f^\ast\alpha^\ast
$$
for any ${\bf t}\in B_\varepsilon(\R^m)_{res}$. By the definition of
integration over chains it is not hard to check that
$$
\int_{e(E, S)^{\bf t}}f^\ast\alpha^\ast=``\sum_{I\in{\cal
N}}"\frac{1}{|\Gamma_I|}\int_{\widehat{\cal M}^{\bf t}_I
(S)^\circ}(f\circ\hat\pi_I)^\ast\alpha^\ast.
$$
So we get that for any ${\bf t}\in B_\varepsilon(\R^m)_{res}$,
\begin{eqnarray}\label{e:2.113}
\langle f^\ast\alpha,
e(E,S)\rangle\!\!\!\!\!\!\!\!\!&&=``\sum_{I\in{\cal
N}}"\frac{1}{|\Gamma_I|}\int_{\widehat{\cal M}^{\bf t}_I
(S)^\circ}(f\circ\hat\pi_I)^\ast\alpha^\ast\\
&&=\int_{\overline{\cal M}_I^{\bf
t}(S)^\circ}(f\circ\hat\pi)^\ast\alpha^\ast.\nonumber
\end{eqnarray}

\begin{remark}\label{rem:2.47}
{\rm If $E\to X$ is a Banach bundle then all groups $\Gamma_i$ and
$\Gamma_I$ in the arguments above are trivial. In this case each
$\hat\pi_I$ (resp. $\hat\Pi^I_J$) is an open embedding (resp.
bundle embedding) from Banach manifolds (resp. bundles) $\widehat
W_I$ (resp. $\widehat E_I$) into $X$ (resp. $E$), and for any
$J\subset I\in{\cal N}$ the projection $\hat\pi^I_J$ (resp.
$\hat\Pi^I_J$) is an open embedding (resp. bundle embedding) from
Banach manifolds $(\hat \pi^I_J)^{-1}\bigl(\widehat W_I)\bigr)$
(resp. bundles $\widehat E_I|_{(\hat \pi^I_J)^{-1}(\widehat
W_J)}$) to ${\it Im}(\hat\pi^I_J)\subset\widehat W_J$ (resp. ${\it
Im}(\hat\Pi^I_J)\subset\widehat E_J$). So $\widehat W_I$ (resp.
$\widehat E_I$) may be identified an open subset $W_I$ of $X$
(resp. $E|_{W_I}=E_I$), and $\hat\pi^I_J$ (resp. $\hat\Pi^I_J$) is
identified with the inclusion $W_I\hookrightarrow W_J$ (resp.
$E_I\hookrightarrow E_J$). The systems of Banach bundles
$(\widehat{\cal E}, \widehat W)$ and $(\widehat{\cal F},\widehat
V)$ actually becomes Banach bundles $E|_W\to W$ and $E|_V\to V$
respectively. Moreover, the global section $\widehat S=\{\widehat
S_I\,|\, I\in{\cal N}\}$ of $(\widehat{\cal E}, \widehat W)$
(resp. $(\widehat {\cal F}, \widehat V)$) is actually the
restriction of the section $S:X\to E$ to $W$ (resp. $V$). As the
global sections of $(\widehat{\cal F}, \widehat V)$,
$\hat\sigma_{ij}$ become ones of the Banach bundle $E\to M$ with
support in $W_i$, $j=1,\cdots,m_i$. Finally the virtual Euler
cycle $e(E, S)^{\bf t}$ in (\ref{e:2.111}) becomes the Euler cycle
in Section~\ref{sec:1}. Hence our virtual Euler cycle is a
generalization of the Euler cycle in Section~\ref{sec:1}.}
\end{remark}

As expected  we need to show that the class
 $e(E,S)$ is independent of the different choices of the sections
$\sigma_{ij}$ and the open sets $W_i$.  It should be first noted
that the dependence on choices of the smooth cut-off function
$\gamma_i$ may be reduced to one on choices of sections
$\sigma_{ij}$. This may be reduced to the cases in the following
two subsections. We follow the arguments in \cite{Lu3}.

\subsubsection{ The dependence on choices of sections
$\tilde\sigma_{ij}$}\label{sec:2.5.1}

We only need to consider the case  in which a section
$\sigma_{n(m_n+1)}$ is added to
$\{\sigma_{n1},\cdots,\sigma_{nm_n}\}$. Let $m'_n=m_n+1$ and
$m'_i=m_i$, $i=1,\cdots,n-1$. Then $m'_1+\cdots+m'_n=m+1$.
Replacing the system of Banach bundles in (\ref{e:2.82}) and its
section in (\ref{e:2.83}) we use the system of Banach bundles
\begin{eqnarray*}
&&\bigl({\bf P}_1^\ast\widehat {\mathcal F}, \widehat
V\times\R^{m+1}\bigr)\\
&&=\bigl \{({\bf P}_1^\ast\widehat F_I, \widehat
V_I\times\R^{m+1}), \hat\pi_I,  \hat\pi^I_J,
\hat\Pi_I,\hat\Pi^I_J, \hat p_I,\Gamma_I\,\bigm|\, J\subset
I\in{\mathcal N}\bigr\},
\end{eqnarray*}
and its corresponding  global Fredholm section
$\Upsilon'=\{\Upsilon'_I\,|\, I\in{\cal N}\}$,
\begin{eqnarray*}
&&\Upsilon'_I: \widehat V_I\times\R^{m+1}\to {\bf
P}_1^\ast\widehat F_I,\\
 &&\qquad\Bigl(\hat x_I, \{t_{ij}\}_{\substack{
 1\le j\le m'_i\\
 1\le i\le n}}\Bigr)\mapsto \hat S_I(\hat x_I)+
 \sum^n_{i=1}\sum^{m'_i}_{j=1}t_{ij}(\hat\sigma_{ij})_I(\hat
 x_I)\\
&&\hspace{35mm} =\hat S_I(\hat x_I)+
 \sum_{i\in I}\sum^{m'_i}_{j=1}t_{ij}(\hat\sigma_{ij})_I(\hat
 x_I).\nonumber
 \end{eqnarray*}
Here ${\bf P}_1$ is still the projection to the first factor, and
$\hat\pi_I, \hat\pi^I_J, \hat\Pi_I,\hat\Pi^I_J, \hat p_I$ are
naturally pullbacks of those projections in (\ref{e:2.72}). (Later
we shall use these notations without confusion occurs.) By shrinking
${\cal W}^\ast$ and $\varepsilon>0$ in Theorem~\ref{th:2.38} we
assume that  $(\widehat{\cal F}^\ast,\widehat V^\ast)$ and $({\bf
P}_1^\ast\widehat{\cal F}^\ast,\widehat V^\ast\times
B_\varepsilon(\R^{m+1}))$ therein are the corresponding systems of
Banach bundles and that $\Upsilon'=\{\Upsilon'_I\,|\,I\in{\cal N}\}$
is the corresponding section which is {\bf strongly} transversal to
the zero section. Then
 corresponding with the manifold $\widehat\Omega_I(s)$
of dimension ${\rm Ind}(S)+m$ in (\ref{e:2.84}) and the projection
in (\ref{e:2.86}) we have the manifold of dimension ${\rm
Ind}(S)+m+1$,
$$
\widehat\Omega'_I(S):=\{(\hat x_I, {\bf t})\in\widehat
V_I^\ast\times B_\varepsilon(\R^{m+1})\,|\, \Upsilon'_I(\hat
x_I,{\bf t})=0\}
$$
and the  projection
$$
\hat\pi^I_J: (\hat
\pi^I_J)^{-1}\bigl(\widehat\Omega'_J(S)\bigr)\to {\rm
Im}(\hat\pi^I_J)\subset\widehat\Omega'_J(S)
$$
for any $J\subset I\in{\cal N}$. Note that for any $I\in{\cal N}$
and $(\hat x_I, {\bf t})\in\widehat V_I^\ast\times
B_\varepsilon(\R^{m+1})$,
$$
\Upsilon'_I\Bigl(\hat x_I,
\{t_{ij}\}_{\substack{
 1\le j\le m'_i\\
 1\le i\le n}}\Bigr)=\Upsilon_I\Bigl(\hat x_I, \{t_{ij}\}_{\substack{
 1\le j\le m_i\\
 1\le i\le n}}\Bigr)+t_{n(m_n+1)}(\hat\sigma_{n(m_n+1)})_I(\hat x_I)
$$
if $n\in I$, and
\begin{equation}\label{e:2.114}
\Upsilon'_I\Bigl(\hat x_I, \{t_{ij}\}_{\substack{
 1\le j\le m'_i\\
 1\le i\le n}}\Bigr)=\Upsilon_I\Bigl(\hat x_I, \{t_{ij}\}_{\substack{
 1\le j\le m_i\\
 1\le i\le n}}\Bigr)
\end{equation}
if $n\notin I$. Let $t_{n(m_n+1)}$ denote the final coordinate in
$\R^{m+1}$ and ${\bf t}^{(1)}$ the first $m$ coordinates. It
follows from (\ref{e:2.114}) that
\begin{eqnarray*}
&&(\Upsilon'_I)^{-1}(0)=\bigl\{(\hat x_I, {\bf t})\in\widehat
V_I^\ast\times B_\varepsilon(\R^{m+1})\,|\,\Upsilon'_I(\hat x_I,
{\bf t})=0\bigr\}\\
&&\hspace{16mm}=\bigl\{(\hat x_I, {\bf t})\in\widehat
V_I^\ast\times B_\varepsilon(\R^{m+1})\,|\,\Upsilon_I(\hat x_I,
{\bf
t}^{(1)})=0\bigr\}\\
&&\hspace{16mm}=(\Upsilon_I)^{-1}(0)\times (-\varepsilon,
\varepsilon)
\end{eqnarray*}
if $n\notin I$. Here $\Upsilon_I$ is as in Theorem~\ref{th:2.38}. In
this case, by Lemma~\ref{lem:1.7} we may choose ${\bf t}=( {\bf
t}^{(1)}, t_{n(m_n+1)})\in B_\varepsilon(\R^{m+1})_{res}$ such that
${\bf t}^{(1)}\in B_\varepsilon(\R^m)_{res}$. So for such a ${\bf
t}$,
\begin{equation}\label{e:2.115}
\widehat{\cal M}'^{\bf t}_I(S)=\widehat{\cal M}_I^{{\bf
t}^{(1)}}(s)\times\bigl\{t_{n(m_n+1)}\bigr\}\;{\rm if}\;n\notin I.
\end{equation}
Here $\widehat{\cal M}'^{\bf t}_I(S):= (\Upsilon'^{({\bf
t})}_I)^{-1}(0)$  and    $\widehat{\cal M}^{{\bf t}^{(1)}}_I(S)=
(\Upsilon^{({\bf t}^{(1)})}_I)^{-1}(0)$ is as in
Corollary~\ref{cor:2.39}.

Now let $n\in I$.  In order to relate $\widehat\Omega'_I(S)$ to
$\widehat\Omega_I(S)$  we need to consider the system of Banach
bundles
\begin{eqnarray*}
&&\bigl({\bf P}_1^\ast\widehat {\mathcal F}, \widehat
V\times\R^{m+1}\times[0,1]\bigr)\\
&&=\bigl \{({\bf P}_1^\ast\widehat F_I, \widehat
V_I\times\R^{m+1}\times[0,1]), \hat\pi_I, \hat\pi^I_J,
\hat\Pi_I,\hat\Pi^I_J, \hat p_I,\Gamma_I\,\bigm|\, J\subset
I\in{\mathcal N}\bigr\},
\end{eqnarray*}
and the Fredholm section $\Phi=\{\Phi_I\,|\, I\in{\cal N}\}$,
\begin{eqnarray*}
&&\Phi_I: \widehat V_I\times\R^{m+1}\times[0,1]\to
{\bf P}_1^\ast\widehat F_I,\\
&& \Bigl(\hat x_I, \{t_{ij}\}_{\substack{ 1\le j\le m'_i\\
 1\le i\le n}}, s\Bigr)\mapsto \hat S_I(\hat x_I) +
 \sum^n_{i=1}\sum^{m_i}_{j=1}t_{ij}(\hat\sigma_{ij})_I(\hat
 x_I)\\
 &&\hspace{39mm}+ s\cdot t_{n(m_n+1)}(\hat\sigma_{n(m_n+1)})_I(\hat
 x_I) \\
&&\hspace{39mm} =\hat S_I(\hat x_I)+ \sum_{i\in
I}\sum^{m_i}_{j=1}t_{ij}(\hat\sigma_{ij})_I(\hat
 x_I)\\
&&\hspace{39mm} + s\cdot
t_{n(m_n+1)}(\hat\sigma_{n(m_n+1)})_I(\hat
 x_I).
 \end{eqnarray*}
As in (\ref{e:2.92}) we have its trivialization representation:
\begin{eqnarray*}
&& \widetilde O_{i_q}\times\R^{m+1}\times[0,1] \to (\widetilde
E_{i_q})_{\tilde x_{i_q}},\\
&&(\tilde x, {\bf t}, s)\mapsto S^T_{i_q}(\tilde x)+
\sum^{m_{i_q}}_{j=1}t_{i_qj} \sigma^T_{i_qj}(\tilde x)+ \sum_{i\in
I\setminus\{i_q\}}\sum^{m_i}_{j=1}t_{ij}\tau^{qT}_{ij}(\tilde
x)\\
&&\hspace{14mm}+ s\cdot t_{n(m_n+1)} \sigma^T_{n(m_n+1)}(\tilde
x).
\end{eqnarray*}
As in the proofs of Theorem~\ref{th:2.38} and
Proposition~\ref{prop:2.42}, by shrinking ${\cal W}^\ast$ and
$\varepsilon>0$ therein we can prove that the restriction of the
section $\Phi=\{\Phi_I\,|\, I\in{\cal N}\}$ to the system of bundles
$\bigl({\bf P}_1^\ast\widehat {\mathcal F}^\ast, \widehat
V^\ast\times B_\varepsilon(\R^{m+1})\times[0,1]\bigr)$ is {\bf
strongly} transversal to the zero section.  By Lemma~\ref{lem:1.7}
we can choose a regular value ${\bf t}$
 of the
projection
$${\cal P}_I: (\Phi_I)^{-1}(0)\to \R^{m+1},\;(\hat x, {\bf t},
s)\mapsto  {\bf t},$$
 i.e. ${\bf t}\in B_\varepsilon(\R^{m+1})_{res}$ such that ${\bf
t}^{(1)}\in B_\varepsilon(\R^m)_{res}$. Then the manifold
$${\cal P}_I^{-1}({\bf t}):=\{(\hat x, {\bf t}, s)\in \widehat
V^\ast\times B_\varepsilon(\R^{m+1})\times[0,1]\,|\, \Phi(\hat x,
{\bf t}, s)=0\}$$
 forms a strong cobordsim between $\widehat{\cal M}'^{\bf t}_I(S)$ and
$$
 \bigl\{(\hat x_I, {\bf t}, 0)\in {\cal P}_I^{-1}({\bf t})\bigr\}=\widehat{\cal M}^{{\bf
 t}^{(1)}}_I(S)\times\{t_{n(m_n+1)}\}
 \approx\widehat{\cal M}^{{\bf t}^{(1)}}_I(S).
 $$
 Combing this and (\ref{e:2.115}) we get that
$\widehat{\cal M}'^{\bf t}_I(S)$ and $\widehat{\cal M}^{{\bf
 t}^{(1)}}_I(S)$ are always strongly cobordant.
Clearly, for any $J\subset I\in{\cal N}$ the corresponding
 projection
$$
 \hat\pi^I_J: (\hat
\pi^I_J)^{-1}\bigl({\cal P}_J^{-1}({\bf t})\bigr)\to {\rm
Im}(\hat\pi^I_J)\subset {\cal P}_J^{-1}({\bf t})
$$
also forms a cobordsim between the projection
$$
 \hat\pi^I_J: (\hat
\pi^I_J)^{-1}\bigl(\widehat{\cal M}^{{\bf t}^{(1)}}_J(S)\bigr)\to
{\rm Im}(\hat\pi^I_J)\subset\widehat{\cal M}^{{\bf t}^{(1)}}_J(S)
$$
and
$$
 \hat\pi^I_J: (\hat
\pi^I_J)^{-1}\bigl(\widehat{\cal M}'^{\bf t}_J(S)\bigr)\to {\rm
Im}(\hat\pi^I_J)\subset\widehat{\cal M}'^{\bf t}_J(S).
$$
Consider the orientation given in Proposition~\ref{prop:2.44}. The
above cobordsims are also oriented ones.  In summary we arrive at:

\begin{proposition}\label{prop:2.48}
For different choices of optimal sections $\tilde\sigma_{ij}$, the
corresponding families such as $\widehat{\cal M}^{\bf t}(S)$ are
oriented strongly cobordant, and thus the corresponding singular
chains as in (\ref{e:2.110}), $e(E, S)^{\bf t}$, are {\rm oriented
cobordant}. In particular, if $Z(S)\cap\partial X=\emptyset$ then
the class $e(E, S)$ is well-defined and is independent of different
choices of optimal sections $\tilde\sigma_{ij}$.
\end{proposition}

\subsubsection{The dependence on choices of
$W_i$}\label{sec:2.5.2}

In this subsection we assume that  $Z(S)\subset{\rm Int}(X)$ and
Assumption~\ref{ass:2.45} is satisfied (so that $e(E, S)$ is
well-defined). We only need to consider the case that {\bf an open
subset $W_{n+1}$ and the corresponding optimal sections
$\tilde\sigma_{(n+1)j}=\gamma_{(n+1)j}\cdot\tilde s_{(n+1)j}$,
$j=1,\cdots, m_{n+1}$, are added.}

The proof can be done as in \S 4.3 in \cite{Lu3}. Here we shall
prove it as a consequence of Claim 2.43 and the restriction
principle in Proposition~\ref{prop:2.51}. By the above we may assume
that $X$ has no boundary. Consider the Banach orbifold $X\times [0,
1]$ with boundary $\partial(X\times [0, 1])=(X\times\{0\})\cup
(X\times\{1\})$ and natural pullback Fredholm orbibundle $(X\times
[0, 1], {\rm pr}^\ast E, {\rm pr}^\ast S)$ via the projection to the
first factor ${\rm pr}:X\times [0, 1]\to X$. $(X\times [0, 1], {\rm
pr}^\ast E, {\rm pr}^\ast S)$ has an orientation induced from one of
$(X, E, S)$. Let $Y_i=X\times\{i\}$, $i=0,1$. They are Banach
suborbifolds of $X\times [0, 1]$, and also closed subsets of the
latter. Obviously, $(Y_i, ({\rm pr}^\ast E)|_{Y_i}, ({\rm pr}^\ast
S)|_{Y_i})$, $i=0,1$ can be identified with $(X, E, S)$. Let
$(\widetilde W_i, \Gamma_i,\pi_i)$ be the chosen Banach orbifold
charts centred at $x_i\in Z(S)$  on $X$,
  $i=1,\cdots,n+1$  such that $Z(S)\subset\cup^n_{i=1}W_i$ and
  $\cap^n_{i=1}W_i=\emptyset$.
Suppose also that $(\widetilde E_i, \Gamma_i,\Pi_i)$, $i=1,\cdots,
n+1$ are the corresponding orbibundle charts on them. As before
let pairs of open sets
 $W^j_i\subset\subset U^j_i$, $i=1,\cdots, n+1$, $j=1, 2, \cdots,n$, be such that
 $Z(S)\subset\cup^n_{i=1}W_i^1$ and
$$
 U^1_i\subset\subset
W^2_i\subset\subset U^2_i\cdots \subset\subset
W^{n}_i\subset\subset U^{n}_i\subset\subset W_i.
$$
For $i=1,\cdots, n$ and $j=1,\cdots, n+1$ denote by
$$
\begin{array}{lcccccr}
 \widetilde W_i^\prime=\widetilde W_i\times[0, \frac{3}{4}),
 \;\!\!\! &W_i^\prime=W_i\times[0, \frac{3}{4}),\;
 &\pi_i^\prime=\pi_i\times{\bf 1},\\
\widetilde W_{n+j}^\prime=\widetilde W_j\times(\frac{1}{4}, 1],\;
&W_{n+j}^\prime=W_j\times(\frac{1}{4}, 1],\;
 &\pi_{n+j}^\prime=\pi_j\times{\bf 1},\\
\Gamma_i^\prime=\{g\times{\bf 1}\,|\,g\in\Gamma_i\},
&\Gamma_{n+j}^\prime=\{g\times{\bf 1}\,|\,
g\in\Gamma_j\},\\
\widetilde E_i^\prime=p_i^\ast\widetilde E_i,\;\!\!\!
&E_i^\prime=({\rm pr}^\ast{E})|_{W_i^\prime},\;
 &p_i: \widetilde W_i\times[0, \frac{3}{4})\to\widetilde W_i,\\
\widetilde E_{n+j}^\prime=q_j^\ast\widetilde E_j,\;\!\!\!
&E_{n+j}^\prime=({\rm pr}^\ast{E})|_{W_{n+j}^\prime},\;
 &q_j: \widetilde W_j\times(\frac{1}{4}, 1]\to\widetilde W_j,
\end{array}$$
where all $p_i, q_j$ are the projections to the first factors. Let
\begin{eqnarray*}
&&{\mathcal N}=\{I\subset\{1,\cdots, n\}\,|\, W_I=\cap_{i\in
I}W_i\ne\emptyset\},\\
&&{\mathcal N}''=\{I\subset\{1,\cdots, n+1\}\,|\, W_I=\cap_{i\in
I}W_i\ne\emptyset\},\\
&&{\mathcal N}^{\prime}=\{I\subset\{1,\cdots, 2n+1\}\,|\,
W'_I=\cap_{i\in I}W'_i\ne\emptyset\}.
\end{eqnarray*}
Then each $I\in{\cal N}^{\prime}$ contains at most $2n-2$ elements.
By  Remark~\ref{rm:2.32} we only need to take pairs of open sets
 $$
 W'^j_i\subset\subset U'^j_i,\; i=1,\cdots, 2n+1, j=1, 2, \cdots,2n-2,
 $$
 such that
 $Z(S)\times [0, 1]\subset\cup^{2n+1}_{i=1}W_i'^1$ and
$$
 U'^1_i\subset\subset
W'^2_i\subset\subset U'^2_i\cdots \subset\subset
W'^{2n-2}_i\subset\subset U'^{2n-2}_i\subset\subset W'_i.
$$
Clearly, for a very small positive $\eta<\frac{1}{100}$ one can
require that
\begin{eqnarray*}
&&W^j_i\times [0, \eta)=(X\times [0, \eta])\cap W'^j_i,\\
&&U^j_i\times [0,\eta)=(X\times [0, \eta])\cap U'^j_i,
\end{eqnarray*}
for $i=1,\cdots, n$, $j=1,\cdots, n-1$, and that
\begin{eqnarray*}
&&W^j_i\times (1-\eta, 1]=(X\times [1-\eta, 1])\cap W'^j_{n+i},\\
&&U^j_i\times (1-\eta, 1]=(X\times [1-\eta, 1])\cap U'^j_{n+i}
\end{eqnarray*}
for $i=1,\cdots, n+1$, $j=1,\cdots, n$. These imply that
\begin{equation}\label{e:2.116}
\left.\begin{array}{ll}
\widetilde W'^j_i\cap(\pi'_i)^{-1}(W'_i)=W^j_i\times [0,\eta),\\
\widetilde U'^j_i\cap(\pi'_i)^{-1}(W'_i)=U^j_i\times [0,\eta),\\
Cl(U^j_i)\times\{0\}=(X\times\{0\})\cap Cl(U'^j_i)
\end{array}\right\}
\end{equation}
for $i=1,\cdots, n$, $j=1,\cdots, n-1$, and that
\begin{equation}\label{e:2.117}
\left.\begin{array}{ll}
\widetilde W'^j_{i+n}\cap(\pi'_{i+n})^{-1}(W'_{i+n})=W^j_i\times (1-\eta, 1],\\
\widetilde U'^j_{i+n}\cap(\pi'_{i+n})^{-1}(W'_{i+n})=U^j_i\times (1-\eta, 1],\\
Cl(U^j_i)\times\{1\}=(X\times\{1\})\cap Cl(U'^j_{n+i})
\end{array}\right\}
\end{equation}
for $i=1,\cdots, n+1$, $j=1,\cdots, n$.

Let $\gamma_k$, $k=1,\cdots,n+1$, are the corresponding smooth
cut-off functions on $\widetilde W_k$ as given by (\ref{e:2.77}).
One can take $\Gamma_k^\prime$-invariant smooth cut-off functions
$\gamma_k^\prime$ on $\widetilde W_k^\prime$, $k=1,\cdots, 2n+1$
such that:
\begin{itemize}
\item $\gamma_k^\prime(x,0)=\gamma_k(x)$, $\forall (x,
0)\in\widetilde W_k^\prime$ and $k=1,\cdots, n$;

\item $\gamma_k^\prime(x,1)=\gamma_{k-n}(x)$, $\forall (x,
1)\in\widetilde W_k^\prime$ and $k=n+1,\cdots, 2n+1$;

\item $\widetilde U_k^{\prime 0}:=\{e\in\widetilde W_k^\prime\;|\;
\gamma_k^\prime(e)>0\}\subset\subset\widetilde W_k^\prime$,
$k=1,\cdots, 2n+1$;

\item $\cup_{k=1}^{2n+1} U_k^{\prime 0}\supset Z(S)\times[0,1]$,
where $U_k^{\prime 0}:=\pi_k'(\widetilde U_k^{\prime 0})$,
$k=1,\cdots, 2n+1$.
\end{itemize}
Then we have:
\begin{eqnarray*}
 && U_k^{\prime 0}\cap
 (W_k\times\{0\})=U_k^0\times\{0\},\quad k=1,\cdots,n,\\
&& U_{n+ k}^{\prime 0}\cap
(W_k\times\{1\})=U_k^0\times\{1\},\;k=1,\cdots,n+1.
 \end{eqnarray*}

Let $\tilde \sigma_{ij}=\gamma_i\cdot\tilde s_{ij}$ be the
(optimal) sections of the bundle $\widetilde E_i\to\widetilde W_i$
as in (\ref{e:2.79}), $i=1,\cdots, n+1$ and $j=1,\cdots, m_i$.
Then the  sections of the bundle $\widetilde
E_k^\prime\to\widetilde W_k^\prime$,
$$
\tilde \sigma^\prime_{kj}(x,t)=\gamma_k^\prime(x,t)\cdot\tilde
s_{kj}(x)\quad\forall\;(x,t)\in\widetilde W_k^\prime,
$$
 $k=1,\cdots,
n, j=1,\cdots, m_k$, and  ones of the bundle $\widetilde
E_{n+k}^\prime\to\widetilde W_{n+k}^\prime$,
$$
\tilde
\sigma^\prime_{(n+k)j}(x,t)=\gamma_{n+k}^\prime(x,t)\cdot\tilde
s_{kj}(x)\quad\forall (x,t)\in\widetilde W_{n+k}^\prime,
$$
$k=1,\cdots,n+1, j=1,\cdots, m_k$, form a  optimal family, and
these $m':=m+ m''$   sections satisfy the corresponding
 properties with  (\ref{e:2.104}). Here
 $$
 m=\sum^n_{i=1}m_i\quad{\rm and}\quad m''=\sum^{n+1}_{i=1}m_i.
 $$
Let ${\cal W}'^\ast$ be a small open neighborhood $Z(p_1^\ast
S)=Z(S)\times [0, 1]$ in $X\times [0, 1]$ such that $Z(S)\times [0,
1]\subset\cup_{k=1}^{2n+1} U_k^{\prime 0}$. As usual we can use all
these data to construct a family of strongly cobordant virtual Euler
chains of $(X\times [0, 1], {\rm pr}^\ast E, {\rm pr}^\ast S)$,
$$
 e({\rm pr}^\ast E, {\rm pr}^\ast S)^{{\bf
t}'}=``\sum_{I\in{\cal N}'}"
\frac{1}{|\Gamma'_I|}\{\hat\pi'_I:\widehat{\cal M}^{{\bf t}'}_I
({\rm pr}^\ast S)^\circ\to X\times [0,1]\}
$$
for ${\bf t}'=\{t_{ij}\}_{\substack{ 1\le j\le m_i\\
 1\le i\le 2n+1}}\in B_\varepsilon(\R^{m'})_{res}$.
Clearly, we can write ${\bf t}'=({\bf t}, {\bf t}'')$,
$$
{\bf t}:=\{t_{ij}\}_{\substack{ 1\le j\le m_i\\
 1\le i\le n}}\in B_\varepsilon(\R^{m})\quad{\rm and}\quad{\bf t}'':=\{t_{ij}\}_{\substack{ 1\le j\le m_i\\
 n+1\le i\le 2n+1}}\in B_\varepsilon(\R^{m''}).
$$
Note that the boundary of $e({\rm pr}^\ast E, {\rm pr}^\ast
S)^{{\bf t}'}$,
$$
 \partial e({\rm pr}^\ast E, {\rm pr}^\ast S)^{{\bf
t}'}=``\sum_{I\in{\cal N}'}"
\frac{1}{|\Gamma'_I|}\{\hat\pi'_I:\partial\widehat{\cal M}^{{\bf
t}'}_I ({\rm pr}^\ast S)^\circ\to \partial(X\times [0,1])\},
$$
is equal to $\triangle(0, {\bf t})\cup (-\triangle(1, {\bf
t}''))$. Here
\begin{eqnarray*}
\triangle(0, {\bf t})=``\sum_{I\in{\cal N}'}"
\frac{1}{|\Gamma'_I|}\{\hat\pi'_I:\partial\widehat{\cal M}^{{\bf
t}'}_I ({\rm pr}^\ast S)^\circ\cap (\hat\pi'_I)^{-1}(X\times\{0\})
\to X\times\{0\}\}\\
\qquad\qquad=``\sum_{I\in{\cal N}}"
\frac{1}{|\Gamma'_I|}\{\hat\pi'_I:\partial\widehat{\cal M}^{{\bf
t}'}_I ({\rm pr}^\ast S)^\circ\cap (\hat\pi'_I)^{-1}(X\times\{0\})
\to X\times\{0\}\},\\
\triangle(1, {\bf t}'')=``\sum_{I\in{\cal N}'}"
\frac{1}{|\Gamma'_I|}\{\hat\pi'_I:\partial\widehat{\cal M}^{{\bf
t}'}_I ({\rm pr}^\ast S)^\circ\cap (\hat\pi'_I)^{-1}(X\times\{1\})
\to X\times\{1\}\}\\
\qquad\qquad=``\sum_{I\in{\cal N}^\ast}"
\frac{1}{|\Gamma'_I|}\{\hat\pi'_I:\partial\widehat{\cal M}^{{\bf
t}'}_I ({\rm pr}^\ast S)^\circ\cap (\hat\pi'_I)^{-1}(X\times\{1\})
\to X\times\{1\}\},
\end{eqnarray*}
where (\ref{e:2.116}) and (\ref{e:2.117}) are used, and
$$
{\cal N}^\ast=\{I\in{\cal N}'\,|\, \max I>n\}=\{J+n\,|\, J\in{\cal
N}''\}.
$$
Let $e(E,S)^{\bf t}$ and $e(E,S)^{{\bf t}''}$ be respectively the
virtual Euler cycles of $(X, E, S)$ constructed from the following
data of two groups
$$
 \left\{ \begin{array}{ll}
 {\cal
W}^\ast:=\{x\in X\,|\, (x, 0)\in {\cal
W}'^\ast\cap(X\times\{0\})\},\;\{W_i\}^n_{i=1},\\
 \{\tilde\sigma_{ij}\}_{\substack{ 1\le j\le m_i\\
 1\le i\le n}},\;\{W^j_i\subset\subset U^j_i\}_{\substack{ 1\le i\le n\\
 1\le j\le n-1}},\;{\bf t}:=\{t_{ij}\}_{\substack{ 1\le j\le m_i\\
 1\le i\le n}}\in B_\varepsilon(\R^{m})
\end{array}\right.$$
$$
\left\{ \begin{array}{ll}
{\cal W}''^\ast:=\{x\in X\,|\, (x, 1)\in
{\cal
W}'^\ast\cap(X\times\{0\})\},\; \{W_i\}^{n+1}_{i=1},\;\\
 \{\tilde\sigma_{ij}\}_{\substack{ 1\le j\le m_i\\
 1\le i\le n+1}},\;\{W^j_i\subset\subset U^j_i\}_{\substack{ 1\le i\le n+1\\
 1\le j\le n}}),\;{\bf t}:=\{t_{ij}\}_{\substack{ 1\le j\le m_i\\
 1\le i\le n+1}}\in B_\varepsilon(\R^{m''}).
\end{array}\right.
$$
Then by the above choices and the restriction principle in
Proposition~\ref{prop:2.51} one easily checks that $\triangle(0,
{\bf t})$ (resp.  $\triangle(1, {\bf t}'')$) can be exactly
identified with the virtual Euler cycle
 $e(E,S)^{\bf t}$ (resp. $e(E,S)^{{\bf t}''}$). So
  Claim 2.43 yields
\begin{eqnarray*}
[e(E,S)^{\bf t}]-[e(E,S)^{{\bf t}''}]\!\!\!\!\!&&=[e(E,S)^{\bf
t}]+[-e(E,S)^{{\bf t}''}]\\
&&=[\triangle(0, {\bf t})]+[-\triangle(1, {\bf
t}'')]\\
&&=[\triangle(0, {\bf t})\cup(-\triangle(1, {\bf
t}''))]\\
&&=[\partial e({\rm pr}^\ast E, {\rm pr}^\ast S)^{{\bf t}'}]=0.
\end{eqnarray*}
The desired result is proved.

From our above result one easily derives that $e(E,S)$ is
independent of concrete choices of the cut-off functions $\gamma_k$
and pairs of open sets $W_i^j\subset\subset U_i^j$.

\begin{remark}\label{rm:2.49}
{\rm If we only assume that $X$ is $1$-regular, the affectivity
arguments in Claim~\ref{cl:2.16} is not needed, and
Claim~\ref{cl:2.23} still holds true. In this case, after some
numbers such as $|\Gamma_I|$, $|\Gamma(\tilde u_i)|$ are replaced by
\begin{eqnarray*}
&&\left\{\begin{array}{ll}
 |\Gamma_i|-|{\rm Ker}(\Gamma_i,\widetilde
W_i)|+ 1\quad&\hbox{if}\;
I=\{i\},\\
\prod^k_{s=1}(|\Gamma_{i_s}|-|{\rm Ker}(\Gamma_{i_s},\widetilde
W_{i_s})|+ 1)\quad&\hbox{if}\; I=\{i_1,\cdots, i_k\}\in{\mathcal N},
\end{array}\right.\\
&&|\Gamma(\tilde u_i)|-|{\rm Ker}(\Gamma(\tilde u_i), \widetilde
O(\tilde u_i))|+1,
\end{eqnarray*}
respectively,  all arguments in this subsection and next two
subsections are all true. We leave them to reader for checking. }
\end{remark}

\subsection{Localization formula and restriction principle}\label{sec:2.7}

Motivated by the localization formula  in \S5 of \cite{Lu3} we here
shall give two abstract versions of them. The second one is given as
a consequence of a restriction principle. Without special
statements, {\bf all Banach orbifolds are also supposed to be
$1$-regular and effective, and all  Banach Fredholm orbibundles are
oriented, have compact zero loci and also satisfy
Assumption~\ref{ass:2.45}.} (Remark~\ref{rm:2.49} is still
effective.)  In particular, $(X,E,S)$ denote  an oriented Banach
Fredholm orbibundle of index $r$ and with compact zero set $Z(S)$,
and satisfy $Z(S)\subset{\rm Int}(X)$ and Assumption~\ref{ass:2.45}.
Another approach to localization in symplectic geometry was
presented in \cite{Chen}.

 Let
$\Lambda\subset Z(S)$ be a nonempty compact subset. In the
construction of the above virtual Euler cycle of
$(X, E, S)$ if one only requires: \\
(i) the finitely many points $x_i$, $i=1,\cdots, n$ belong to
$\Lambda$,\\
 (ii) $\{W_i\}^n_{i=1}$ is an open cover of $\Lambda$,\\
(iii) ${\cal W}^\ast$ is an open neighborhood of $\Lambda$ contained in $\cup^n_{i=1}U_i^0$,
denoted by ${\cal W}^\ast(\Lambda)$,\\
then we get the corresponding systems of Banach bundles,
\begin{eqnarray*}
&&\bigl(\widehat{\mathcal E}(\Lambda), \widehat W(\Lambda)\bigr):=
\bigl\{\bigl(\widehat E_I, \widehat W_I\bigr), \hat\pi_I,
\hat\Pi_I, \Gamma_I, \hat\pi^I_J,\hat\Pi^I_J,
\lambda^I_J\bigm| J\subset I\in{\mathcal N}\bigr\},\\
&&\bigl(\widehat{\mathcal F}(\Lambda), \widehat V(\Lambda)\bigr):=
\bigl\{\bigl(\widehat F_I, \widehat V_I\bigr), \hat\pi_I,
\hat\Pi_I, \Gamma_I, \hat\pi^I_J,\hat\Pi^I_J, \lambda^I_J\bigm|
J\subset I\in{\mathcal
N}\bigr\},\\
&&\bigl(\widehat{\mathcal F}^\ast(\Lambda), \widehat
V^\ast(\Lambda)\bigr):= \bigl\{\bigl(\widehat F_I^\ast, \widehat
V_I^\ast\bigr), \hat\pi_I, \hat\Pi_I, \Gamma_I,
\hat\pi^I_J,\hat\Pi^I_J, \lambda^I_J\bigm| J\subset I\in{\mathcal
N}\bigr\},
\end{eqnarray*}
where $V_I$, $\widehat V_I$, $\widehat F_I$ (resp. $V_I^\ast$,
$\widehat V_I^\ast$, $\widehat F_I^\ast$) are defined by
(\ref{e:2.65}) and (\ref{e:2.66}) (resp.  (\ref{e:2.71})). These
are called the {\bf Banach bundle systems of $(X, E, S)$ relative
to $\Lambda$}.  Corresponding with (\ref{e:2.93}) we have a family
of Fredholm sections of $\bigl(\widehat{\mathcal F}^\ast(\Lambda),
\widehat V^\ast(\Lambda)\bigr)$, $\Upsilon^{({\bf
t})}(\Lambda)=\{\Upsilon^{({\bf t})}_I(\Lambda)\,|\, I\in{\cal
N}\}$, where ${\bf t}=\{t_{ij}\}_{\substack{ 1\le j\le m_i\\
 1\le i\le n}}\in B_\varepsilon(\R^m)$ and
 the section $\Upsilon^{({\bf t})}_I(\Lambda): \widehat V^\ast_I\to \widehat
F^\ast_I$ is given by
$$
 \hat x_I \mapsto \Upsilon_{I}(\tilde x, {\bf t})=\hat
S_I(\hat x_I)+
 \sum^n_{i=1}\sum^{m_i}_{j=1}t_{ij}(\hat\sigma_{ij})_I(\hat x_I).
$$
 Moreover, there exists a residual subset
$B_\varepsilon(\R^m)_{res}^\Lambda\subset B_\varepsilon(\R^m)$
such that for each ${\bf t}\in B_\varepsilon(\R^m)_{res}^\Lambda$
the global section
 $\Upsilon^{({\bf t})}(\Lambda)$ is strongly transversal to the
zero section. So we get a family of strongly cobordant rational
chains near $\Lambda$ of dimension $r$, denoted by
$$
 e(E, S)^{\bf
t}_\Lambda:=``\sum_{I\in{\cal N}}"
\frac{1}{|\Gamma_I|}\{\hat\pi_I:\widehat{\cal M}^{\bf t}_I
(S,\Lambda)^\circ\to X\},\forall{\bf t}\in
B_\varepsilon(\R^m)_{res}^\Lambda.
$$
They are called {\bf virtual Euler chains relative to} $\Lambda$ of
the triple $(X, E, S)$. Here  $\widehat{\cal M}^{\bf t}_I(S,
\Lambda):= (\Upsilon^{({\bf t})}_I(\Lambda))^{-1}(0)$ is a manifold
of dimension ${\rm Ind}(S)$ and with compact closure in $Cl(\widehat
V^\ast_I)\subset\widehat V_I^{+}$, and $\widehat{\cal M}^{\bf t}_I
(S,\Lambda)^\circ$ is the regular part of $\widehat{\cal M}^{\bf
t}_I (S,\Lambda)$. Clearly, if $Z(S)$ is compact then for any
compact subset $\Lambda\subset Z(S)$ one may view $e(E, S)^{\bf t}$
as virtual Euler chains  of $(X, E, S)$ relative to $\Lambda$. As
below (\ref{e:2.112}) we can also define $\overline{\cal M}^{\bf
t}(S,\Lambda)$, $\overline{\cal M}^{\bf t}(S,\Lambda)^\circ$ and
identify $e(E, S)^{\bf t}_\Lambda$ with $\hat\pi:\overline{\cal
M}^{\bf t}(S,\Lambda)^\circ\to X$.

\begin{proposition}\label{prop:2.50}({\bf First localization formula}).
Let  $P$  be a smooth orbifold of finite dimension and $f:X\to P$
be a smooth map. If $\alpha\in H^r(P, \R)$ has a closed
differential form representative $\alpha^\ast$ such that the
support ${\rm supp}(f^\ast\alpha^\ast)$ of $f^\ast\alpha^\ast$ is
contained in $\Lambda$, then
\begin{eqnarray*}
\langle f^\ast\alpha, e(E,S)\rangle\!\!\!\!\!\!\!\!\!&&=\int_{e(E,
S)^{\bf
t}_\Lambda}f^\ast\alpha^\ast\\
&&:=``\sum_{I\in{\cal N}}"\frac{1}{|\Gamma_I|}\int_{\widehat{\cal
M}^{\bf t}_I
(S,\Lambda)^\circ}(f\circ\hat\pi_I)^\ast\alpha^\ast\nonumber\\
&&=\int_{\overline{\cal M}^{\bf
t}(S,\Lambda)^\circ}(f\circ\hat\pi)^\ast\alpha^\ast\nonumber
\end{eqnarray*}
for any ${\bf t}\in B_\varepsilon(\R^m)^\star_{res}\subset
B_\varepsilon(\R^m)_{res}^\Lambda$. Here
$B_\varepsilon(\R^m)^\star_{res}$ is also a residual subset in
$B_\varepsilon(\R^m)$, and $\sum_{I\in{\cal N}}$ is put in double
quotation marks  because of the same reason as before, i.e. the
repeating part is only counted once.
\end{proposition}

\noindent{\bf Proof.}\quad Recall the construction of $e(E,
S)^{\bf t}_\Lambda$. We have chosen  pairs of open sets
$W^j_i\subset\subset U^j_i$ such that
$$
 W^1_i\subset\subset U^1_i\subset\subset
W^2_i\subset\subset U^2_i\cdots \subset\subset
W^{n-1}_i\subset\subset U^{n-1}_i\subset\subset W_i,
$$
 $i=1,\cdots, n$, $j=1,\cdots, n-1$.  For each $i=1,\cdots, n$, let
  $\tilde\sigma_{il}=\gamma_i\cdot \tilde s_{il}$ be the chosen optimal
 sections with support in $\widetilde W_i^1$, $l=1,\cdots, m_i$,
 and $\widetilde U_i^0=\{\tilde x\in \widetilde W_i\,|\, \gamma_i(\tilde
 x)>0\}$ and $U_i^0=\pi_i(\widetilde U_i)$.

Fix a small open neighborhood ${\cal W}(\Lambda)_0$ of $\Lambda$
in $X$ so that
$$
{\cal W}(\Lambda)_0\subset\subset{\cal W}^\ast(\Lambda).
$$
Here ${\cal W}^\ast(\Lambda)$ is a small open neighborhood of
$\Lambda$ contained in $\cup^n_{i=1}U_i^0$.  Since
$Z(S)\setminus{\cal W}^\ast(\Lambda)$ is compact we may
furthermore take finite many points $x_i\in Z(S)\setminus{\cal
W}^\ast(\Lambda)$ and the corresponding orbifold charts
$(\widetilde W_i, \Gamma_i,\pi_i)$ centred at $x_i$ in $X$,
$i=n+1,\cdots, n'$ such that
\begin{equation}\label{e:2.118}
Z(S)\setminus{\cal W}^\ast(\Lambda)\subset{\cal
W}(\Lambda)_2\quad{\rm and}\quad {\cal W}(\Lambda)_2\cap{\cal
W}(\Lambda)_0=\emptyset
\end{equation}
for ${\cal
W}(\Lambda)_2:=\cup^{n'}_{i=n+1}W_i=\cup^{n'}_{i=n+1}\pi_i(\widetilde
W_i)$. Of course we  have also the corresponding  orbibundle
charts $(\widetilde E_i, \Gamma_i,\Pi_i)$, $i=n+1,\cdots, n'$.
Furthermore we choose pairs of open sets $W^j_i\subset\subset
U^j_i$, $i=1,\cdots, n$, $j=n,\cdots, n'-1$, and
$W^l_k\subset\subset U^l_k$, $k=n+1,\cdots, n'$, $l=1,\cdots,
n'-1$  such that
$$
 W^1_i\subset\subset U^1_i\subset\subset
W^2_i\subset\subset U^2_i\cdots \subset\subset
W^{n'-1}_i\subset\subset U^{n'-1}_i\subset\subset W_i,
$$
 $i=1,\cdots, n'$, $j=1,\cdots, n'-1$. By (\ref{e:2.118}) we can
 require that
$$
 Z(S)\setminus{\cal
W}^\ast(\Lambda)\subset{\cal W}(\Lambda)'_2\quad{\rm for}\quad
{\cal W}(\Lambda)'_2:=\cup^{n'}_{i=n+1}W^1_i.
$$
 So for each $i=n+1,\cdots, n'$, we can take $\Gamma_i$-invariant smooth cut-off
 function $\gamma_i:\widetilde W_i\to [0, 1]$ with support in $\widetilde
 W_i^1$, and smooth sections   $\tilde\sigma_{il}=\gamma_i\cdot \tilde s_{il}$,
 $l=1,\cdots, m_i$, such that
\begin{equation}\label{e:2.119}
Z(S)\setminus{\cal
W}^\ast(\Lambda)\subset\cup^{n'}_{i=n+1}U_i^0\quad{\rm for}\quad
\widetilde U_i^0=\{\tilde x\in \widetilde W_i\,|\, \gamma_i(\tilde
 x)>0\},
\end{equation}
and that all sections $\tilde\sigma_{ij}$, $j=1,\cdots, m_i$ and
$i=1,\cdots, n'$, are optimal for $(X, E, S)$.  Denote by ${\cal
N}'$ the collection of all subsets $I$ of $\{1,\cdots, n'\}$ such
that $W_I:=\cap_{i\in I}W_i\ne\emptyset$. Then we have a system of
Banach bundles,
$$
\bigl(\widehat{\mathcal E}', \widehat W'\bigr):=
\bigl\{\bigl(\widehat E'_I, \widehat W'_I\bigr), \hat\pi'_I,
\hat\Pi'_I, \Gamma'_I, \hat\pi'^I_J,\hat\Pi'^I_J,
\lambda'^I_J\,\bigm|\, J\subset I\in{\mathcal N}'\bigr\}.
$$
It is easily checked that for any $J\subset I\in{\cal N}$,
\begin{equation}\label{e:2.120}
\left.\begin{array}{ll} (\widehat E'_I, \widehat
W'_I\bigr)=(\widehat E_I, \widehat W_I\bigr),\;
 \hat\pi'_I=\hat\pi_I,\;\hat\Pi'_I=\hat\Pi_I,\\
 \Gamma'_I=\Gamma_I,\; \hat\pi'^I_J=\hat\pi^I_J,\;\hat\Pi'^I_J=\hat\Pi^I_J,\;
  \lambda'^I_J=\lambda^I_J.
  \end{array}\right\}
\end{equation}
 Following (\ref{e:2.65}), for each $I\in{\cal N}'$ with $|I|=k$ we define
$$
 V'_I:= \Bigl(\bigcap_{i\in I}W^k_i\Bigr)\setminus
\Bigl(\bigcup_{J:|J|>k}\bigl(\bigcap_{j\in
J}Cl(U^k_j)\bigr)\Bigr).
$$
 Since $\cap^n_{i=1}W_i=\emptyset$, then
$V'_I=W^{|I|}_I$ for any $I\in{\cal N}'$ containing
$\{1,\cdots,n\}$.
 Note that $V_I'\subset
V_I$ for any $I\in{\cal N}$ in general. Setting $\widehat
V_I'=\hat\pi_I^{-1}(V_I')$ and $\widehat F_I'=\widehat
E'_I|_{\widehat V_I'}$ we get the renormalization of
$\bigl(\widehat{\mathcal E}', \widehat W'\bigr)$,
$$
 \bigl(\widehat{\mathcal F}',
\widehat V'\bigr):= \bigl\{\bigl(\widehat F'_I, \widehat
V'_I\bigr), \hat\pi'_I, \hat\Pi'_I, \Gamma'_I,
\hat\pi'^I_J,\hat\Pi'^I_J, \lambda'^I_J\,\bigm|\, J\subset
I\in{\mathcal N}'\bigr\}.
$$
By (\ref{e:2.119}) one can take an open neighborhood ${\cal W}^\ast$
of $Z(S)$ in $X$ so small that
$$
 {\cal W}^\ast\subset {\cal
 W}^\ast(\Lambda)\cup(\cup^{n'}_{i=n+1}U_i^0).
$$
Setting $V'^\ast_I=V'_I\cap{\cal W}^\ast$, $\widehat
V'^\ast_I=(\hat\pi'_I)^{-1}(V'^\ast_I)$ and $\widehat
F'^\ast_I=\widehat F'_I|_{\widehat V'^\ast_I}$ we get the
restriction of $\bigl(\widehat{\mathcal F}', \widehat V'\bigr)$ to
${\cal W}^\ast$,
$$
 \bigl(\widehat{\mathcal F}'^\ast,
\widehat V'^\ast\bigr):= \bigl\{\bigl(\widehat F'^\ast_I, \widehat
V'^\ast_I\bigr), \hat\pi'_I, \hat\Pi'_I, \Gamma_I,
\hat\pi'^I_J,\hat\Pi'^I_J, \lambda'^I_J\,\bigm|\, J\subset
I\in{\mathcal N}'\bigr\}.
$$
By shrinking ${\cal W}^\ast$ and $\varepsilon>0$ we have a family of
sections of $\bigl(\widehat{\mathcal F}'^\ast, \widehat
V'^\ast\bigr)$ that are strongly transversal to the zero section,
$\Upsilon^{{\bf t}'}=\{\Upsilon^{({\bf t}')}_I\,|\, I\in{\cal
N}'\}$,
$$
 \Upsilon^{({\bf t}')}_I: \widehat
V'^\ast_I\to \widehat F'^\ast_I,\; \hat x_I \mapsto
\Upsilon_{I}(\tilde x, {\bf t}')=\hat S'_I(\hat x_I)+
 \sum^k_{i=1}\sum^{m_i}_{j=1}t_{ij}(\hat\sigma'_{ij})_I(\hat x_I),
$$
where ${\bf t}'=\{t_{ij}\}_{\substack{ 1\le j\le m_i\\
 1\le i\le l}}\in B_\varepsilon(\R^{m'})_{res}$ and $m'=\sum^{n'}_{i=1}m_i$. These yield
  a family of cobordant virtual Euler cycles of
dimension $r$,
$$
e(E, S)^{{\bf t}'}=``\sum_{I\in{\cal N}'}"
\frac{1}{|\Gamma'_I|}\{\hat\pi'_I:\widehat{\cal M}^{{\bf t}'}_I
(S)^\circ\to X\}\;\forall{\bf t}'\in B_\varepsilon(\R^{m'})_{res},
$$
 where $\widehat{\cal M}^{{\bf t}'}_I
(S)=(\Upsilon^{({\bf t}')}_I)^{-1}(0)$. Note that
$$
 \langle f^\ast\alpha, e(E,S)\rangle=``\sum_{I\in{\cal
N}'}"\frac{1}{|\Gamma'_I|}\int_{\widehat{\cal M}^{{\bf t}'}_I
(S)^\circ}(\hat\pi_I'\circ f)^\ast\alpha^\ast
$$
for any ${\bf t}'\in B_\varepsilon(\R^{m'})_{res}$. By
Lemma~\ref{lem:1.7} the set
$$
\Bigl\{{\bf t}=\{t_{ij}\}_{\substack{ 1\le j\le m_i\\
 1\le i\le n}}\Bigm| {\bf t}'=\{t_{ij}\}_{\substack{ 1\le j\le m_i\\
 1\le i\le n'}}\in B_\varepsilon(\R^{m'})_{res}\Bigr\}
 $$
 is a residual subset in $B_\varepsilon(\R^m)$. Let
$B_\varepsilon(\R^m)^\star_{res}$ be the intersection of it and
$B_\varepsilon(\R^m)_{res}^\Lambda$, which is also a residual
subset in $B_\varepsilon(\R^m)$.
Then for each ${\bf t}=\{t_{ij}\}_{\substack{ 1\le j\le m_i\\
 1\le i\le n}}\in B_\varepsilon(\R^m)^\star_{res}$ one may take
 a ${\bf t}'\in B_\varepsilon(\R^{m'})_{res}$ of form $\{t_{ij}\}_{\substack{ 1\le j\le m_i\\
 1\le i\le n'}}\in B_\varepsilon(\R^{m'})_{res}$. So it suffice to
 prove  that
$$
 ``\sum_{I\in{\cal
N}'}"\frac{1}{|\Gamma'_I|}\int_{\widehat{\cal M}^{{\bf t}'}_I
(S)^\circ}(\hat\pi_I'\circ f)^\ast\alpha^\ast=``\sum_{I\in{\cal
N}}"\frac{1}{|\Gamma_I|}\int_{\widehat{\cal M}^{\bf t}_I
(S,\Lambda)^\circ}(f\circ\hat\pi_I)^\ast\alpha^\ast
$$
for such ${\bf t}$ and ${\bf t}'$ very small. To this end note
that the second condition in (\ref{e:2.118}) implies
$$\int_{\widehat{\cal M}^{{\bf t}'}_I (S)^\circ}(\hat\pi_I'\circ
f)^\ast\alpha^\ast=0
$$
if $I\in{\cal N}'$ contains a number in $\{n+1,\cdots, n'\}$
because the support ${\rm supp}(f^\ast\alpha^\ast)$ of
$f^\ast\alpha^\ast$ is contained in $\Lambda$. So we only need to
prove that
\begin{equation}\label{e:2.121}
``\sum_{I\in{\cal N}}"\frac{1}{|\Gamma'_I|}\int_{\widehat{\cal
M}^{{\bf t}'}_I (S)^\circ}(\hat\pi_I'\circ
f)^\ast\alpha^\ast=``\sum_{I\in{\cal
N}}"\frac{1}{|\Gamma_I|}\int_{\widehat{\cal M}^{\bf t}_I
(S,\Lambda)^\circ}(f\circ\hat\pi_I)^\ast\alpha^\ast.
\end{equation}
For any $I\in{\cal N}$, it is clear that $\widehat
V_I'\subset\widehat V_I$ and $\widehat F'^\ast_I=\widehat
F^\ast_I|_{\widehat V_I'}$. Moreover,  by (\ref{e:2.70}),
$(\hat\sigma'_{ij})_I=0$  for any $i>n$. So $\Upsilon^{({\bf
t}')}_I$ is equal to the restriction of $\Upsilon^{({\bf t})}_I$
to $\widehat V'^\ast_I$. From these  it follows that
\begin{equation}\label{e:2.122}
\widehat{\cal M}^{{\bf t}'}_I (S)\subset\widehat{\cal M}^{\bf t}_I
(S,\Lambda)\quad{\rm and}\quad\hat\pi_I|_{\widehat{\cal M}^{{\bf
t}'}_I (S)}=\hat\pi'_I.
\end{equation}
Moreover, for any $J\subset I\in{\cal N}$ it is also easily seen
that the inclusion and restriction in (\ref{e:2.122}) are compatible
with the projections
\begin{eqnarray*}
&&\hat\pi^I_J: (\hat \pi^I_J)^{-1}\bigl(\widehat{\cal M}^{\bf
t}_J(S,\Lambda)\bigr)\to {\rm Im}(\hat\pi^I_J)\subset\widehat{\cal
M}^{\bf
t}_J(S, \Lambda)\quad{\rm and}\\
&&\hat\pi'^I_J: (\hat \pi'^I_J)^{-1}\bigl(\widehat{\cal M}^{{\bf
t}'}_J(S)\bigr)\to {\rm Im}(\hat\pi'^I_J)\subset\widehat{\cal
M}^{{\bf t}'}_J(S).
\end{eqnarray*}
Note also that (\ref{e:2.119}) and the second condition in
(\ref{e:2.118}) implies
$$\hat\pi_I\Bigl(\widehat{\cal M}^{\bf t}_I
(S,\Lambda)\setminus\widehat{\cal M}^{{\bf t}'}_I
(S)\Bigr)\cap{\cal W}(\Lambda)_0\subset \hat\pi_I(\widehat
V_I\setminus\widehat V_I')\cap{\cal W}(\Lambda)_0=\emptyset
$$
for each $I\in{\cal N}$. Then the desired result (\ref{e:2.121})
follows from these, (\ref{e:2.120}) and the fact that ${\rm
supp}(f^\ast\alpha^\ast)\subset\Lambda\subset{\cal W}(\Lambda)_0$.
\hfill$\Box$\vspace{2mm}

Similar to the arguments above Proposition~\ref{prop:1.12}, for a
connected component $\Lambda$ of $Z(S)$ we choose the above
neighborhood ${\cal W}(\Lambda)$ of $\Lambda$ such that ${\cal
W}(\Lambda)\cap (Z(S)\setminus\Lambda)=\emptyset$, and then get a
family of cobordant cycles
 $\{e(E,S)_\Lambda^{\bf t}\,|\,{\bf t}\in
B_\epsilon(\R^l)_{res}\}$ and corresponding homology class in
$e(E,S)_\Lambda\in H_r(X,\Q)$. If $\Lambda_i$, $i=1,\cdots, p$,
are all connected components of $\Lambda$, then
$$
e(E,S)=\sum^p_{i=1}e(E,S)_{\Lambda_i}.
$$
This easily follows from the following result.

\begin{proposition}\label{prop:2.51}
 ({\bf The restriction principle}). Let $(X, E, S)$  be an oriented Banach Fredholm
 orbibundle with compact zero locus. Assume that a closed subset $Y\subset X$ is a
 Banach suborbifold of finite positive
 codimension. Then $(Y, E|_Y, S|_Y)$ is also an oriented
Banach Fredholm  orbibundle with compact zero locus, and there exist
a family of strongly cobordant virtual Euler chains  of the triple
$(X, E, S)$,
$$
 e(E, S)^{{\bf
t}'}=``\sum_{I\in{\cal N}'}"
\frac{1}{|\Gamma_I|}\{\hat\pi_I:\widehat{\cal M}^{{\bf t}'}_I
(S)^\circ\to X\}\;\forall{\bf t}'\in B_\varepsilon(\R^{m'})_{res},
$$
and that of strongly cobordant virtual Euler chains of the triple
$(Y, E|_Y, S|_Y)$,
$$
 e(E|_Y, S|_Y)^{\bf
t}=``\sum_{I\in{\cal N}}"
\frac{1}{|\Gamma_I|}\{\hat\pi_I^y:\widehat{\cal M}^{\bf t}_I
(S|_Y)^\circ\to Y\}\;\forall{\bf t}\in B_\varepsilon(\R^m)_{res},
$$
such that $m<m'$, ${\cal N}\subset{\cal N}'$,  $\widehat{\cal
M}^{{\bf t}'}_I(S)\cap \hat\pi_I^{-1}(Y)=\emptyset$ for any
$I\in{\cal N}'\setminus{\cal N}$, and that
\begin{equation}\label{e:2.123}
\left.\begin{array}{ll}
 \widehat{\cal M}^{\bf t}_I
(S|_Y)=\widehat{\cal M}^{{\bf t}'}_I(S)\cap
\hat\pi_I^{-1}(Y),\\
\widehat{\cal M}^{\bf t}_I (S|_Y)^\circ=\widehat{\cal M}^{{\bf
t}'}_I(S)^\circ\cap \hat\pi_I^{-1}(Y),\\
\widehat{\cal M}^{\bf t}_I (S|_Y)^{sing}=\widehat{\cal M}^{{\bf
t}'}_I(S)^{sing}\cap \hat\pi_I^{-1}(Y),\\
 \hat\pi^{yI}_J=\hat\pi^I_J|_{\widehat{\cal M}^{\bf t}_I (S|_Y)},
\end{array}
\right\}
\end{equation}
for any $I\in{\cal N}$ and $J\subset I$, and
$$
{\bf t}'=\{t_{ij}\}_{\substack{ 1\le j\le m_i\\
 1\le i\le n'}}\in B_\varepsilon(\R^{m'})_{res}\quad{\rm and}\quad
{\bf t}:=\{t_{ij}\}_{\substack{ 1\le j\le m_i\\
 1\le i\le n}}.
$$
 Consequently, for such ${\bf t}'$ and ${\bf t}$ it holds that
\begin{equation}\label{e:2.124}
\overline{\cal M}^{\bf t}(S|_Y)^\circ=\overline{\cal M}^{{\bf
t}'}(S)^\circ\cap \hat\pi^{-1}(Y)\quad{\rm and}\quad
\hat\pi^y=\hat\pi|_{\overline{\cal M}^{\bf t}(S|_Y)^\circ}.
\end{equation}
In this case we say $e(E|_Y, S|_Y)^{\bf t}$ to be the intersection
 of $e(E, S)^{{\bf t}'}$ with $Y$, or the restriction of
$e(E, S)^{{\bf t}'}$ to $Y$.
\end{proposition}

\noindent{\bf Proof.}\quad Since $Z(S)$ is compact and $Y$ is closed
in $X$, the zero locus $Z(S|_Y)=Z(S)\cap Y$ is also a compact subset
in $Y$. One can check that Assumption~\ref{ass:2.45} also holds for
the suborbifold $Y$. By the definition of the suborbifold in
Section~\ref{sec:2.1} we can choose finitely many points
\begin{eqnarray*}
x_i\in Z(S|_Y), i=1,\cdots, n,\\
x_i\in Z(S)\setminus Z(S|_Y), i=n+1,\cdots, n'
 \end{eqnarray*}
 and orbifold charts $(\widetilde W_i, \Gamma_i, \pi_i)$ on $X$
centered at $x_i$, $i=1,\cdots, n'$ such that:
\begin{description}
\item[($1^\circ$)] For each $i=1,\cdots, n$ there exists a
submanifold $\widetilde W_i^y\subset\widetilde W_i$ that is stable
under $\Gamma_i$ and compatible with the inclusion map, such that
the restriction $(\widetilde W_i^y, \Gamma_i^y,
\pi_i^y)=(\widetilde W_i^y, \Gamma_i|_{\widetilde W_i^y},
\pi_i|_{\widetilde W_i^y})$ is a Banach orbifold chart on $Y$, and
thus $W_i^y:=\pi_i(\widetilde W_i^y)=W_i\cap Y$. ({\it Warning}:
$\Gamma_i^y=\Gamma_i$ as abstract groups. Thus later we also write
$\Gamma_i^y$ as $\Gamma_i$ without danger of confusions.)

\item[($2^\circ$)] $\{W_i^y=\pi_i(\widetilde W_i^y)\}^n_{i=1}$ is
an open cover of $Z(S|_Y)$ with $\cap^n_{i=1}W_i^y=\emptyset$, and
 $\{W_i=\pi_i(\widetilde W_i)\}^{n'}_{i=1}$ is an open cover of
$Z(S)$ in $X$ with $\cap^{n'}_{i=1}W_i=\emptyset$.

\item[($3^\circ$)] For $i=n+1,\cdots, n'$, $W_i=\pi_i(\widetilde
W_i)$ are disjoint with $Y$.
\end{description}
 Similarly,  let $(\widetilde
E_i,\Gamma_i,\Pi_i)$ be the orbifold bundle charts on $E$
corresponding with $(\widetilde W_i, \Gamma_i, \pi_i)$ for
$i=1,\cdots, n'$, then we get such a chart  on $E|_Y$,
$$
(\widetilde E_i^y,\Gamma_i,\Pi_i^y)=(\widetilde E_i|_{W_i^y}
,\Gamma_i,\Pi_i|_{\widetilde E_i^y})
$$
which is corresponding with $(\widetilde W_i^y, \Gamma_i,
\pi_i^y)$ for $i=1,\cdots, n$. As before take
 pairs of open sets
 $W^j_i\subset\subset U^j_i$, $i=1,\cdots, n',\; j=1, 2, \cdots,n'-1$,
 such that
$$
 U^1_i\subset\subset
W^2_i\subset\subset U^2_i\cdots \subset\subset
W^{n'-1}_i\subset\subset U^{n'-1}_i\subset\subset W_i.
$$
We can also require that
\begin{equation}\label{e:2.125}
Cl(U^{j}_i\cap Y)=Cl(U^j_i)\cap Y,\;i=1,\cdots,n,\;j=1,\cdots,
n-1.
\end{equation}
(These are easily satisfied by using local coordinate charts and
shrinking $W_i$ if necessary.)  As before let ${\cal N}'$ be the
set of all finite subsets $I\in\{1,\cdots, n'\}$ with
$W_I=\cap_{i\in I}W_i\ne\emptyset$, and ${\cal N}=\{I\in{\cal
N}'\,|\, \max I\le n\}$.
 For each $I\in{\cal N}'$ with $|I|=k$ one uses
(\ref{e:2.65}) to define $V_I$. Then for
$i=1,\cdots,n,\;j=1,\cdots,n-1$ set $W^{yj}_i=W^j_i\cap Y$ and
$U^{yj}_I=U^j_i\cap Y$. Then $W^{yj}_i\subset\subset U_i^{yj}$ and
for $i=1,\cdots, n$,
$$
 U^{y1}_i\subset\subset
W^{y2}_i\subset\subset U^{y2}_i\cdots \subset\subset
W^{y(n-1)}_i\subset\subset U^{y(n-1)}_i\subset\subset W_i^y.
$$
Such choices are reasonable because of $\cap^n_{i=1}W_i^y=\emptyset$
and Remark~\ref{rm:2.32}. As before, for each $I\in{\cal N}$ with
$|I|=k$ one uses (\ref{e:2.65}) to define
$$
 V_I^y:=\Bigl(\bigcap_{i\in I}W_i^{yk}\Bigr)\setminus
\Bigl(\bigcup_{J:|J|>k}\bigl(\bigcap_{j\in
J}Cl(U_j^{yk})\bigr)\Bigr).
$$
 Using these we can, as before, get
the  system of Banach bundles of $(X,E,S)$,
\begin{eqnarray*}
&&\bigl(\widehat{\mathcal E}, \widehat W\bigr):=
\bigl\{\bigl(\widehat E_I, \widehat W_I\bigr), \hat\pi_I,
\hat\Pi_I, \Gamma_I, \hat\pi^{I}_J,\hat\Pi^{I}_J,
\lambda^I_J\,\bigm|\, J\subset I\in{\mathcal N}'\bigr\},\\
&&\bigl(\widehat{\mathcal F}, \widehat V\bigr):=
\bigl\{\bigl(\widehat F_I, \widehat V_I\bigr), \hat\pi_I,
\hat\Pi_I, \Gamma_I, \hat\pi^{I}_J,\hat\Pi^{I}_J,
\lambda^I_J\,\bigm|\, J\subset I\in{\mathcal
N}'\bigr\},\\
&&\bigl(\widehat{\mathcal F}^{\ast}, \widehat V^{\ast}\bigr):=
\bigl\{\bigl(\widehat F_I^{\ast}, \widehat V_I^{\ast}\bigr),
\hat\pi_I, \hat\Pi_I, \Gamma_I, \hat\pi^{I}_J,\hat\Pi^{I}_J,
\lambda^I_J\,\bigm|\, J\subset I\in{\mathcal N}'\bigr\},
\end{eqnarray*}
and those corresponding with $(Y, E|_Y, S|_Y)$,
\begin{eqnarray*}
&&\bigl(\widehat{\mathcal E}^y, \widehat W^y\bigr):=
\bigl\{\bigl(\widehat E_I^y, \widehat W_I^y\bigr), \hat\pi_I^y,
\hat\Pi_I^y, \Gamma_I, \hat\pi^{yI}_J,\hat\Pi^{yI}_J,
\lambda^I_J\,\bigm|\, J\subset I\in{\mathcal N}\bigr\},\\
&&\bigl(\widehat{\mathcal F}^y, \widehat V^y\bigr):=
\bigl\{\bigl(\widehat F_I^y, \widehat V_I^y\bigr), \hat\pi_I^y,
\hat\Pi_I^y, \Gamma_I, \hat\pi^{yI}_J,\hat\Pi^{yI}_J,
\lambda^I_J\,\bigm|\, J\subset I\in{\mathcal
N}\bigr\},\\
&&\bigl(\widehat{\mathcal F}^{y\ast}, \widehat V^{y\ast}\bigr):=
\bigl\{\bigl(\widehat F_I^{y\ast}, \widehat V_I^{y\ast}\bigr),
\hat\pi_I^y, \hat\Pi_I^y, \Gamma_I, \hat\pi^{yI}_J,\hat\Pi^{yI}_J,
\lambda^I_J\,\bigm|\, J\subset I\in{\mathcal N}\bigr\},
\end{eqnarray*}
where the small neighborhood ${\cal W}^{y\ast}$ of $Z(S|_Y)$ in
$Y$ is taken as ${\cal W}^\ast\cap Y$.

Moreover, by ($3^\circ$) and (\ref{e:2.125}) it is easy to see
that
\begin{equation}\label{e:2.126}
\left.\begin{array}{ll}
 W_I\cap Y=\emptyset\quad\forall I\in{\cal N}'\setminus{\cal N},\\
 V_I\cap Y=V_I^y\quad\forall
I\in{\cal N}.
\end{array}\right\}
\end{equation}
It follows that for any $I\in{\cal N}$,
\begin{equation}\label{e:2.127}
 \left.\begin{array}{ll}
 \bigl(\widehat E_I^y, \widehat
W_I^y\bigr)=\bigl(\widehat
E_I|_{\widehat W_I^y}, \widehat W_I\cap(\hat\pi_I)^{-1}(Y)\bigr)\\
\bigl(\widehat F_I^y, \widehat V_I^y\bigr)=\bigl(\widehat
F_I|_{\widehat V_I^y}, \widehat V_I\cap(\hat\pi_I)^{-1}(Y)\bigr)\\
\bigl(\widehat F_I^{y\ast}, \widehat
V_I^{y\ast}\bigr)=\bigl(\widehat F^\ast_I|_{\widehat V_I^{y\ast}},
\widehat V^\ast_I\cap(\hat\pi_I)^{-1}(Y)\bigr).
\end{array}\right\}
\end{equation}

As before we assume that each bundle $\widetilde E_i\to\widetilde
W_i$ is trivial. Then the  trivializations ${\cal T}_i: \widetilde
W_i\times (\widetilde E_i)_{\tilde x_i}\to \widetilde E_i$ in
(\ref{e:2.74}) restrict to  natural ones
$$
 {\cal T}_i^y: \widetilde
W_i^y\times (\widetilde E_i)_{\tilde x_i}\to \widetilde
E_i^y,\quad i=1,\cdots,n.
$$
 Note that the lift section $\widetilde
S_i:\widetilde W_i\to\widetilde E_i$ of $S|_{W_i}$ also restricts
to such a section
$$\widetilde S_i^y=\widetilde S_i|_{\widetilde
W_i^y}:\widetilde W_i^y\to\widetilde E_i^y$$
 of $S|_{W_i^y}$. Let $\widetilde S_i^T$ be the  trivialization
representation of $\widetilde S_i$ under ${\cal T}_i$.  It is
easily checked that under ${\cal T}_i^y$ the  trivialization
representation of $\widetilde S_i^y$,
$$\widetilde S_i^{yT}:\widetilde W_i^y\to
(\widetilde E_i^y)_{\tilde x_i},$$
  is exactly the restriction of $\widetilde S_i^T$ to $\widetilde
W_i^y$, i.e., $\widetilde S_i^{yT}=\widetilde S_i^T|_{\widetilde
W_i^y}$.

Moreover, let  $\gamma_i$ be the smooth cut-off functions in
(\ref{e:2.77}), $i=1,\cdots, n'$. Then $\gamma_i$, $i=1,\cdots,
n$, naturally restricts to such a one $\gamma_i^y:\widetilde
W_i^{y1}\to [0,1]$ with support in $\widetilde W_i^{y1}$,  and
such that
\begin{eqnarray*}
&&\widetilde U^{y0}_i:=\{\tilde x\in\widetilde W_i^y\,|\,
\gamma_i(\tilde x)>0\}=\widetilde U^{0}_i\cap\widetilde
W_i^y\subset\subset\widetilde W_i^{y1} \quad{\rm and}\\
&&Z(S|_Y)\subset\cup^n_{i=1}U^{y0}_i=\cup^n_{i=1}\pi_i(\widetilde
U^{y0}_i).
\end{eqnarray*}
Carefully checking (\ref{e:2.79}) and (\ref{e:2.80}) we see that
the sections $\tilde\sigma_{ij}^y:=\gamma_i\cdot \tilde s_{ij}^y:
\widetilde W_i^y\to\widetilde E_i^y$ and their trivialization
representations
$$\tilde\sigma^{yT}_{ij}:\widetilde W_i^y\to
(\widetilde E_i)_{x_i},\;\tilde x\mapsto \gamma(\tilde x)v_{ij}
$$
are exactly the restrictions of $\tilde\sigma_{ij}$ in
(\ref{e:2.79}) and $\tilde\sigma^T_{ij}$ in (\ref{e:2.80}) to
$\widetilde W_i^y$, $j=1,\cdots, m_i$. As before we use these to
get a global section  of $\bigl(\widehat{\mathcal F}^{\ast},
\widehat V^{\ast}\bigr)$, $\Upsilon^{({\bf t}')}=\{\Upsilon^{({\bf
t}')}_I\,|\, I\in{\cal N}'\}$ given by
$$
 \Upsilon^{({\bf t}')}_I: \widehat
V^{\ast}_I\to \widehat F^{\ast}_I,\; \hat x_I \mapsto\hat S_I(\hat
x_I)+
 \sum^{n'}_{i=1}\sum^{m_i}_{j=1}t_{ij}(\hat\sigma_{ij})_I(\hat x_I),
$$
and  that of $\bigl(\widehat{\mathcal F}^{y\ast}, \widehat
V^{y\ast}\bigr)$, $\Upsilon^{y({\bf t})}=\{\Upsilon^{y({\bf
t})}_I\,|\, I\in{\cal N}\}$ given by
$$
\Upsilon^{y({\bf t})}_I: \widehat V^{y\ast}_I\to \widehat
F^{y\ast}_I,\; \hat x_I \mapsto\hat S_I^y(\hat x_I)+
 \sum^n_{i=1}\sum^{m_i}_{j=1}t_{ij}(\hat\sigma_{ij}^y)_I(\hat x_I).
 $$
Here ${\cal N}=\{I\in{\cal N}'\,|\, \max I\le n\}$. Furthermore we
can find a small $\varepsilon>0$ and a residual subset
$B_\varepsilon(\R^{m'})_{res}\subset B_\varepsilon(\R^{m'})$,
$m'=\sum^{n'}_{i=1}m_i$ such that for any
\begin{equation}\label{e:2.128}
{\bf t}'=\{t_{ij}\}_{\substack{ 1\le j\le m_i\\
 1\le i\le n'}}\in B_\varepsilon(\R^{m'})_{res}\quad{\rm and}\quad
{\bf t}:=\{t_{ij}\}_{\substack{ 1\le j\le m_i\\
 1\le i\le n}}
\end{equation}
({\it Note}: when ${\bf t}'$ takes over
$B_\varepsilon(\R^{m'})_{res}$ all corresponding ${\bf t}$ form a
residual subset $B_\varepsilon(\R^{m})_{res}\subset
B_\varepsilon(\R^{m})$), the above sections are {\bf strongly
transversal} to the zero section (by shrinking ${\cal W}^\ast$ and
increasing $n$ and $n'$ if necessary). As before we set
\begin{eqnarray*}
&& \widehat{\cal M}^{{\bf t}'}(S)=\bigl\{\widehat{\cal M}^{{\bf
t}'}_I(S)=(\Upsilon^{({{\bf t}'})}_I)^{-1}(0) \,|\,I\in{\cal N}'\bigr\}\quad{\rm and}\\
 &&\widehat{\cal M}^{\bf t}(S|_Y)=\bigl\{\widehat{\cal M}^{\bf
t}_{I}(S|_Y)=(\Upsilon^{({\bf t})}_{I})^{-1}(0) \,|\,I\in{\cal
N}\bigr\}.
\end{eqnarray*}
Note that $\hat\pi_I$ (resp. $\hat\Pi_I$) restricts to $\hat\pi_I^y$
(resp. $\hat\Pi_I^y$) on $\widehat W_I^y$, $\widehat V_I^y$ and
$\widehat V_I^{y\ast}$ (resp. $\widehat E_I^y$, $\widehat F_I^y$ and
$\widehat F_I^{y\ast}$) and these restrictions are also compatible
with $\hat\pi^I_J$ and $\hat\pi^{yI}_J$ for any $J\subset I\in{\cal
N}$. One can checks that
$$
(\hat\sigma_{ij}^y)_I=\left\{\begin{array}{ll}
(\hat\sigma_{ij})_I|_{\widehat V_I^y}\;&{\rm if}\; I\in{\cal N},\\
0\;&{\rm if}\;I\in{\cal N}'\setminus{\cal N}.
\end{array}\right.
$$
It follows from these, (\ref{e:2.126}) and (\ref{e:2.127}) that
for ${\bf t}, {\bf t}'$ in (\ref{e:2.128}), $\Upsilon^{y({\bf
t})}_I= \Upsilon^{({\bf t}')}_I|_{\widehat V_I^{y\ast}}$ and thus
\begin{eqnarray*}
&&\widehat{\cal M}^{{\bf
t}'}_I(S)\cap(\hat\pi_I)^{-1}(Y)=\widehat{\cal M}^{\bf t}_{I}(S|_Y),\quad\forall I\in{\cal N}\quad{\rm and}\\
&&\widehat{\cal M}^{{\bf
t}'}_I(S)\cap(\hat\pi_I)^{-1}(Y)=\emptyset, \quad \forall
I\in{\cal N}'\setminus{\cal N}.
\end{eqnarray*}
(\ref{e:2.123}) and (\ref{e:2.124}) follow immediately.
\hfill$\Box$\vspace{2mm}

\begin{corollary}\label{cor:2.52}
 ({\bf Second localization formula}). Let $f$ be an oriented smooth orbifold of finite
 dimension, $Q\subset P$ be an oriented closed suborbifold of
 dimension $k$  and a smooth map $f:X\to P$ be such that $Y:=f^{-1}(Q)\subset X$
 is a Banach  suborbifold of codimension $\dim P-k$.
 Then $(Y, E|_Y, S|_Y)$ is still an oriented Banach Fredholm
 orbibundle of index ${\rm Ind}(S|_Y)=r+k-\dim P$ and
for any class $\kappa\in H_\ast(Q,\R)$ of dimension $\dim P-r$ it
holds that
$$\langle PD_Q(\kappa), f_\ast(e(E|_Y, S|_Y))\rangle=\langle PD_P(\kappa),
f_\ast(e(E,S))\rangle.
$$
Here $PD_Q(a)$ (resp. $PD_P(a)$) is the Poincar\'e dual of $a$ in
$H^\ast(Q,\R)$ (resp. $H^\ast(P,\R)$).
\end{corollary}

\noindent{\bf Proof.}\quad Let $R_S$ be a retraction from a tubular
neighborhood $U(Q)$ of $Q$ in $P$ onto $Q$ (which is identified with
the projection the normal bundle of $Q$ in $P$.) We can take a
representative form $Q^\ast$ of the Poincar\'e dual $PD_P(Q)$ of $Q$
in $P$ whose support can be required to be contained in $U(Q)$. Then
for a representative form $\kappa_Q^\ast$ of $PD_Q(\kappa)$ the form
$\kappa^\ast:=R_S^\ast(\kappa^\ast_Q)\cup Q^\ast$ is  a
representative form of $PD_P(\kappa)$ whose support is contained in
$U(Q)$.

Let $e(E, S)^{{\bf t}'}$ and $e(E|_Y, S|_Y)^{\bf t}$ be  as
constructed in Proposition~\ref{prop:2.51}. As in the proof of
Proposition~\ref{prop:1.12}, it follows from (\ref{e:2.113}) and
(\ref{e:2.124}) that

\begin{eqnarray*}
\langle PD_Q(\kappa), (f|_Y)_\ast(e(E|_Y,
S|_Y))\rangle\!\!\!\!\!\!\!\!\!\!&&=\langle
(f|_Y)^\ast(PD_Q(\kappa)), e(E|_Y, S|_Y)\rangle\\
&&=\int_{e(E|_Y, S|_Y)^{\bf
t}}(f|_Y)^\ast\kappa^\ast_Q\\
&&=\int_{\overline{\cal M}^{\bf t} (S|_Y)^\circ
}(f|_Y\circ\hat\pi_I^y)^\ast\kappa^\ast_Q\nonumber\\
&&=\int_{(f|_Y\circ\hat\pi_I^y)(\overline{\cal M}^{\bf t}
(S|_Y)^\circ)
}\kappa^\ast_Q\nonumber\\
&&=\int_{(f\circ\hat\pi_I)(\overline{\cal M}^{{\bf t}'}
(S)^\circ)\cap Q}\kappa^\ast_Q\nonumber\\
&&=\int_{R_Q((f\circ\hat\pi_I)(\overline{\cal M}^{{\bf t}'}
(S)^\circ)\cap Q)}\kappa^\ast_Q\nonumber\\
&&=\int_{(f\circ\hat\pi_I)(\overline{\cal M}^{{\bf t}'}
(S)^\circ)\cap Q}R_Q^\ast\kappa^\ast_Q\nonumber\\
&&=\int_{(f\circ\hat\pi_I)(\overline{\cal M}^{{\bf t}'}
(S)^\circ)}(R_Q^\ast\kappa^\ast_Q)\wedge Q^\ast\nonumber\\
&&=\int_{(f\circ\hat\pi_I)(\overline{\cal M}^{{\bf t}'}
(S)^\circ)}\kappa^\ast\nonumber\\
&&=\int_{e(E,S)^{{\bf t}'}}f^\ast\kappa^\ast\nonumber\\
&&=\langle f^\ast(PD_P(\kappa)), e(E,S)\rangle\\
&&=\langle PD_P(\kappa), f_\ast(e(E,S))\rangle.
\end{eqnarray*}
This completes the proof of the proposition.
  \hfill$\Box$\vspace{2mm}

\subsection{Properties}\label{sec:2.8}

In this subsection we shall generalize the other properties for
Euler classes of the Banach Fredholm bundles in \S\ref{sec:1.3} to
virtual Euler classes. Without special statements, {\bf all Banach
orbifolds are supposed to $1$-regular and effective, and all  Banach
Fredholm orbibundles are oriented, have compact zero loci and also
satisfy Assumption~\ref{ass:2.45}.} (Remark~\ref{rm:2.49} is still
effective.)

\begin{lemma}\label{lem:2.53}
Let $p_i:E^{(i)}\to X$, $i=1,2$, be two Banach orbibundles, and
 $E:=E^{(1)}\oplus E^{(2)}\to X$ be their direct sum. For
any compact subset $K\subset X$ let $(\widetilde
W_j,\Gamma_j,\pi_j)$ be orbifold charts centred at $x_j\in K$  on
$X$, $j=1,\cdots, n$, and them satisfy (\ref{e:2.8}).
 For $i=1,2$ and $j=1,\cdots, n$ let  $(\widetilde E^{(i)}_j,
\Gamma_j,\Pi^{(i)}_j)$ be the corresponding Banach orbibundle
charts on $(\widetilde W_j,\Gamma_j,\pi_j)$ of $E^{(i)}$ with
support $E^{(i)}_j=(p^{(i)})^{-1}(W_j)$  such that $(\widetilde
E^{(1)}_j\oplus\widetilde E^{(2)}_j,
\Gamma_j,\Pi^{(1)}_j\oplus\Pi^{(2)}_j)$ are such charts on
$(\widetilde W_j,\Gamma_j,\pi_j)$ of $E^{(1)}\oplus E^{(2)}$ with
support $E_j:=E^{(1)}_j\oplus E^{(2)}_j=(p^{(1)}\oplus
p^{(2)})^{-1}(W_j)$. Denote by the bundle projections $\tilde
p^{(i)}_j: \widetilde E^{(i)}_j\to\widetilde W_j$ and $\tilde
p_j:=\tilde p^{(1)}_j\oplus\tilde p^{(2)}_j: \widetilde
E_j:=\widetilde E^{(1)}_j\oplus\widetilde E^{(2)}_j$, $i=1,2$ and
$j=1,\cdots,n$. Note that the action of $\Gamma_j$ on $\widetilde
E_j$ is given by $\phi\cdot\tilde u=(\phi\cdot\tilde u_1,
\phi\cdot\tilde u_2)$ for $\phi\in\Gamma_j$ and $\tilde u=(\tilde
u_1,\tilde u_2)\in\widetilde E_j$. Assume that
$$
(\widehat{\cal E}^{(i)}(K), \widehat W(K))= \bigl\{\bigl(\widehat
E_I^{(i)}, \widehat W_I\bigr), \hat\pi_I, \hat\Pi^{(i)}_I,
\Gamma_I, \hat\pi^I_J,\hat\Pi^{(i)I}_J, \lambda^I_J\,\bigm|\,
J\subset I\in{\mathcal N}\bigr\}
$$
are resolutions of $(E^{(i)}, X)$ as given in \S 2.2, $i=1,2$.
Then their direct sum
\begin{eqnarray*}
&&(\widehat{\cal E}^{(1)}(K)\oplus\widehat{\cal E}^{(2)}(K),
\widehat W(K)):=\\
&& \bigl\{\bigl(\widehat E_I^{(1)}\oplus\widehat E_I^{(2)},
\widehat W_I\bigr), \hat\pi_I, \hat\Pi_I, \Gamma_I,
\hat\pi^I_J,\hat\Pi^I_J, \lambda^I_J\,\bigm|\, J\subset
I\in{\mathcal N}\bigr\}
\end{eqnarray*}
is such a resolution of $(E, X)$ near $K$, where
$\hat\Pi_I=\hat\Pi_I^{(1)}\oplus\hat\Pi_I^{(2)}$,
 and $\hat\Pi^I_J=\hat\Pi^{(1)I}_J\oplus\hat\Pi^{(2)I}_J$.
Moreover, for the chosen
 pairs of open sets
 $W^j_i\subset\subset U^j_i$, $i=1,\cdots, n,\; j=1, 2, \cdots,n-1$,
 such that
$$
 U^1_i\subset\subset
W^2_i\subset\subset U^2_i\cdots \subset\subset
W^{n-1}_i\subset\subset U^{n-1}_i\subset\subset W_i,
$$
the renormalization of  $(\widehat{\cal
E}^{(1)}(K)\oplus\widehat{\cal E}^{(2)}(K), \widehat W(K))$, denoted
by  $(\widehat{\cal F}, \widehat V(K))$, is  equal to the direct sum
$(\widehat{\cal F}^{(1)}(K)\oplus\widehat{\cal F}^{(2)}(K), \widehat
V(K))$, where $(\widehat{\cal F}^{(i)}(K), \widehat V(K))$ are the
renormalizations of $(\widehat{\cal E}^{(i)}(K), \widehat W(K))$
constructed from the above pairs of open sets, $i=1,2$. Furthermore,
one has also:
\begin{description}
\item[(a)] If $\widehat S^{(i)}=\{\widehat S^{(i)}_I:\,I\in{\cal
N}\}$ are the global sections of $(\widehat{\cal E}^{(i)}(K),
\widehat W(K))$ generated by the sections $S^{(i)}:X\to E^{(i)}$,
$i=1,2$ then $\widehat S^{(1)}\oplus\widehat S^{(2)}=\{\widehat
S^{(1)}_I\oplus\widehat S^{(2)}_I\,:\,I\in{\cal N}\}$ is the
section of $(\widehat{\cal E}^{(1)}(K)\oplus\widehat{\cal
E}^{(2)}(K), \widehat W(K))$ generated by the section
$S^{(1)}\oplus S^{(2)}:X\to E^{(1)}\oplus E^{(2)}$.

 \item[(b)] If $\hat\sigma_l=\{(\hat\sigma_l)_I: I\in{\cal N}\}$ are
 the sections of $(\widehat{\cal F}^{(i)}, \widehat V(K))$
  generated by smooth sections $\tilde\sigma_l^{(i)}:\widetilde W_l\to\widetilde
 E_l^{(i)}$ with supports in $\widetilde W_l^1$ then
$\hat\sigma^{(1)}_l\oplus\hat\sigma^{(2)}_l:=\{(\hat
\sigma^{(1)}_l)_I\oplus(\hat \sigma^{(2)}_l)_I\,:\,I\in{\cal N}\}$
is such a section of $(\widehat{\cal
F}^{(1)}(K)\oplus\widehat{\cal F}^{(2)}(K), \widehat V(K))$
generated by the section $\tilde\sigma_l^{(1)}\oplus
\tilde\sigma_l^{(2)}:\widetilde W_l\to \widetilde
E^{(1)}_l\oplus\widetilde E^{(2)}_l$.
\end{description}
\end{lemma}

\noindent{\bf Proof.}\quad  Carefully checking the constructions and
arguments in \S~\ref{sec:2.3} one easily gets the proof of the
lemma. \hfill$\Box$\vspace{2mm}

As a generalization of Proposition~\ref{prop:1.13} we have:

\begin{proposition}\label{prop:2.54}
({\bf Stability}).  Let the Banach orbibundles $p_i:E^{(i)}\to X$,
$i=1,2$ and $E:=E^{(1)}\oplus E^{(2)}\to X$ be as in
Lemma~\ref{lem:2.53}. For sections $S^{(i)}:X\to E^{(i)}$, $i=1,2$
and their natural sum  $S:=S^{(1)}\oplus S^{(2)}:X\to E$, assume
that $(X, E, S)$ is an oriented Banach Fredholm orbibundle of
index $r$ and with compact zero locus $Z(S)$. Clearly,
$Z(S)=Z(S^{(1)})\cap Z(S^{(2)})$ and
$Z(S)=Z(S^{(1)}|_{Z(S^{(2)}})$. Also assume  that $S^{(2)}$ is
transversal to the zero section at each point $x\in Z(S)$. (So
there exists an open neighborhood ${\cal U}$ of $Z(S)$ in $X$ such
that $Z(S^{(2)})^\star:=Z(S^{(2)})\cap{\cal U}$ is a smooth
orbifold).  Then
$$D(S^{(1)}|_{Z(S^{(2)})}):(TZ(S^{(2)}))|_{Z(S)} \to
E^{(1)}|_{Z(S)}$$
 is Fredholm and has index ${\rm
Ind}(D(S^{(1)}|_{Z(S^{(2)})})={\rm Ind}(DS)$.  Furthermore if both
$(X, E, S)$ and $(Z(S^{(2)})^\star, E^{(1)}|_{Z(S^{(2)})^\star},
S^{(1)}|_{Z(S^{(2)})^\star})$ also satisfy Assumption~\ref{ass:2.45}
then there exists a virtual Euler cycle representative $C$ (resp.
$C_1$) of $e(E,S)$ (resp. $e(E^{(1)}|_{Z(S^{(2)})^\star},
S^{(1)}|_{Z(S^{(2)})^\star})$) such that $C=C_1$ and thus
$$e(E,S) = (i_{Z(S^{(2)})^\star})_\ast e(E^{(1)}|_{Z(S^{(2)})^\star}, S^{(1)}|_{Z(S^{(2)})^\star}).$$
Here $(i_{Z(S^{(2)})^\star})_\ast$ is the homomorphism between
homology groups induced by the inclusion $i_{Z(S^{(2)})^\star}:
Z(S^{(2)})^\star\hookrightarrow X$.
\end{proposition}

\noindent{\bf Proof.}\quad   Since the Euler class of $(X, E, S)$
depends only on the restriction of $(X, E, S)$ to a neighborhood of
$Z(S)$ in $X$ we have assumed that $Z(S^{(2)})$ itself is a smooth
orbifold. As in the proof of Proposition~\ref{prop:2.51} we may use
the definition of the suborbifold to choose finitely many points
$x_i\in Z(S)\subset Z(S^{(2)})$ and orbifold charts $(\widetilde
W_i, \Gamma_i, \pi_i)$  centered at $x_i$ on $X$ and with
$W_i\subset{\cal U}$, $i=1,\cdots, n$ such that:
\begin{description}
\item[(i)] For each $i$ there exists a submanifold $\widetilde
W_i^z\subset\widetilde W_i$ that is stable under $\Gamma_i$ and
compatible with the inclusion map, such that the restriction
$(\widetilde W_i^z, \Gamma_i^z, \pi_i^z)=(\widetilde W_i^z,
\Gamma_i|_{\widetilde W_i^z}, \pi_i|_{\widetilde W_i^z})$ is a
Banach orbifold chart on $Z(S^{(2)})$. (As explained in the proof of
Proposition~\ref{prop:2.51}, $\Gamma_i^z=\Gamma_i$ and they shall be
written as $\Gamma_i$ below.)
 \item[(ii)]
$Z(S)\subset\cup^n_{i=1}W_i^z=\cup^n_{i=1}\pi_i(\widetilde W_i^z)$
($\subset\cup^n_{i=1}W_i=\cup^n_{i=1}\pi_i(\widetilde W_i)$).
\end{description}
For $k=1,2$ and each $i$ let $(\widetilde
E_i^{(k)},\Gamma_i,\Pi_i^{(k)})$ be the orbifold bundle charts on
$E^{(k)}$ corresponding with $(\widetilde W_i, \Gamma_i, \pi_i)$,
then we get such a chart
$$(\widetilde
E_i^{(1)z},\Gamma_i,\Pi_i^{(1)z})=(\widetilde E^{(1)}_i|_{W_i^z}
,\Gamma_i,\Pi^{(1)}_i|_{\widetilde E_i^{(1)z}})$$
 on $E^{(1)}|_{Z(S^{(2)})}$ corresponding with $(\widetilde W_i^z,
\Gamma_i, \pi_i^z)$. As before we assume that all bundles
$\widetilde E_i^{(k)}\to\widetilde W_i$ are trivial,
$i=1,\cdots,n$ and $k=1,2$. For the chosen trivializations
\begin{equation}\label{e:2.129}
{\cal T}_i^{(k)}: \widetilde W_i\times (\widetilde
E^{(k)}_i)_{x_i}\to \widetilde E_i^{(k)},\;i=1,\cdots, n,
\end{equation}
they restrict to
natural ones
\begin{equation}\label{e:2.130}
{\cal T}_i^{(1)z}: \widetilde W_i^z\times (\widetilde
E_i^{(1)})_{\tilde x_i}\to \widetilde E_i^{(1)z},\; i=1,\cdots, n.
\end{equation}
Now the lift section $\widetilde S_i^{(1)}:\widetilde
W_i\to\widetilde E_i^{(1)}$ of $S^{(1)}|_{W_i}$  restricts to such a
section $\widetilde S^{(1)z}_i=\widetilde S^{(1)}_i|_{\widetilde
W_i^z}:\widetilde W_i^z\to\widetilde E_i^{(1)z}$ of
$S^{(1)}|_{W_i^z}$. It is easily checked that the trivialization
representation $\widetilde S_i^{(1)zT}:\widetilde W_i^z\to
(\widetilde E_i^{(1)z})_{\tilde x_i}$ of $\widetilde S_i^{(1)z}$
under ${\cal T}_i^{(1)z}$ is exactly the restriction of $\widetilde
S_i^{(1)T}$ to $\widetilde W_i^z$, i.e. $\widetilde
S_i^{(1)zT}=\widetilde S_i^{(1)T}|_{\widetilde W_i^z}$. As in the
proof of Proposition~\ref{prop:1.13} we can show that the vertical
differential $D\widetilde S^{(1)z}_i(x)$ at each zero point $x$ of
$\widetilde S^{(1)z}_i$ is Fredholm and has the same index as
$D\widetilde S^{(1)}_i(x)$. This yields the first conclusion.

Take pairs of open sets
 $W^j_i\subset\subset U^j_i$, $i=1,\cdots, n,\; j=1, 2, \cdots,n-1$,
 such that $Z(S)\subset\cup^n_{i=1}W_i^1$ and that
$$
 U^1_i\subset\subset
W^2_i\subset\subset U^2_i\cdots \subset\subset
W^{n-1}_i\subset\subset U^{n-1}_i\subset\subset W_i.
$$
Let $(\widehat{\cal F}^{(j)}, \widehat V)$, $j=1,2$ and
$(\widehat{\cal F}, \widehat V)=(\widehat{\cal
F}^{(1)}\oplus\widehat{\cal F}^{(2)}, \widehat V) $ be the
corresponding bundle systems constructed in Lemma~\ref{lem:2.53} (by
taking $K=Z(S)$).
 As the construction of the bundle system  $\bigl(\widehat{\mathcal F}^{y}, \widehat
V^{y}\bigr)$ in the proof of Proposition~\ref{prop:2.51} we now can
get  such bundle systems of $(Z(S^{(2)}), E^{(j)}|_{Z(S^{(2)})},
S^{(j)}|_{Z(S^{(2)})})$,
$$
\bigl(\widehat{\mathcal F}^{(j)z}, \widehat V^{z}\bigr):=
\bigl\{\bigl(\widehat F_I^{(j)z}, \widehat V_I^{z}\bigr),
\hat\pi_I^z, \hat\Pi_I^{(j)z}, \Gamma_I,
\hat\pi^{zI}_J,\hat\Pi^{(j)zI}_J, \lambda^I_J\,\bigm|\, J\subset
I\in{\mathcal N}\bigr\},
$$
which are exactly the restriction of the bundle system
$(\widehat{\cal F}^{(j)}, \widehat V)$, $j=1,2$.
 Let us choose
a small neighborhood ${\cal W}^\ast$ of $Z(S)$ in $X$ and set
${\cal W}^{\ast z}:={\cal W}^\ast\cap Z(S^{(2)})$. Denote by
$\bigl(\widehat{\mathcal F}^{(j)z\ast}, \widehat V^{z\ast}\bigr)$
(resp. $\bigl(\widehat{\mathcal F}^{\ast}, \widehat
V^{\ast}\bigr)$)  the restriction system of bundles of
$\bigl(\widehat{\mathcal F}^{(j)z}, \widehat V^{z}\bigr)$ (resp.
$\bigl(\widehat{\mathcal F}, \widehat V\bigr)$) to ${\cal
W}^{z\ast}$ (resp. ${\cal W}^{\ast }$), $j=1, 2$. Then
$\bigl(\widehat{\mathcal F}^{(1)z\ast}, \widehat V^{z\ast}\bigr)$
is exactly the restriction of $\bigl(\widehat{\mathcal
F}^{(1)\ast}, \widehat V^{\ast}\bigr)$ to ${\cal W}^{z\ast}$.

Moreover, the smooth cut-off function $\gamma_i$ in (\ref{e:2.77})
naturally restricts to such a one $\gamma_i^z:\widetilde W_i^z\to
[0,1]$ and such that
\begin{eqnarray*}
&&\widetilde U^{z0}_i:=\{\tilde x\in\widetilde W_i^z\,|\,
\gamma_i(\tilde x)>0\}=\widetilde U^{0}_i\cap\widetilde
W_i^z\subset\subset\widetilde W_i^{z1},\\
&&Z(S)\subset\cup^n_{i=1}U^{z0}_i=\cup^n_{i=1}\pi_i(\widetilde
U^{z0}_i).
\end{eqnarray*}

 Note that $S^{(2)}$ is transversal to the
zero section. As before, for each $i=1,\cdots,n$ we can choose
$m_i$ vectors $v^{(1)}_{il}\in (\widetilde E^{(1)}_i)_{x_i}$,
$l=1,\cdots,m_i$ and then use $\gamma_i$, $\gamma_i^z$ and the
trivializations in (\ref{e:2.129}) and (\ref{e:2.130}) to
construct the sections $\tilde\sigma_{ij}^{(1)}=\gamma_i\cdot
\tilde s_{ij}^{(1)}: \widetilde W_i\to\widetilde
E_i^{(1)}\subset\widetilde E_i$ and
$$
\tilde\sigma_{ij}^{(1)z}=\gamma_i\cdot \tilde
s_{ij}^{(1)z}: \widetilde W_i^z\to\widetilde E_i^{(1)z}
$$
such that for sufficiently small ${\cal W}^\ast$, $\varepsilon>0$
and each ${\bf t}=\{t_{ij}\}_{\substack{ 1\le j\le m_i\\
 1\le i\le n}}$ in a residual subset
$B_\varepsilon(\R^m)_{res}\subset B_\varepsilon(\R^m)$ the global
section  of the Banach bundle system $\bigl(\widehat {\mathcal
F}^\ast, \widehat V^\ast\bigr)$, $\Upsilon^{({\bf
t})}=\{\Upsilon^{({\bf t})}_I\,|\, I\in{\cal N}\}$,
$$
 \Upsilon^{({\bf t})}_I: \widehat
V^\ast_I\to \widehat F^\ast_I,\; \hat x_I \mapsto\hat S_I(\hat
x_I)+
 \sum^n_{i=1}\sum^{m_i}_{j=1}t_{ij}(\hat\sigma^{(1)}_{ij})_I(\hat x_I),
 $$
and that of $\bigl(\widehat {\mathcal F}^{(1)z\ast}, \widehat
V^{z\ast}\bigr)$, $\Upsilon^{({\bf t})(1)z}=\{\Upsilon^{({\bf
t})(1)z}_I\,|\, I\in{\cal N}\}$,
$$
 \Upsilon^{({\bf t})(1)z}_I: \widehat
V^{z\ast}_I\to \widehat F^{(1)z\ast}_I,\; \hat x_I \mapsto\hat
S_I^{(1)z}(\hat x_I)+
 \sum^n_{i=1}\sum^{m_i}_{j=1}t_{ij}(\hat\sigma_{ij}^{(1)z})_I(\hat x_I),
$$
 are strongly transversal to the zero section. Here
 $\{(\hat\sigma^{(1)}_{ij})_I\,|\,I\in{\cal N}\}$ denotes the
 global section of $\bigl(\widehat {\mathcal
F}^\ast, \widehat V^\ast\bigr)$  generated by
$\tilde\sigma_{ij}^{(1)}$. Clearly, $(\hat\sigma^{(1)}_{ij})_I$
takes values in $\widehat F^{(1)z\ast}_I$. Consider the families
of manifolds of the zero loci of these sections,
\begin{eqnarray*}
&& \widehat{\cal M}^{\bf t}(S)=\{\widehat{\cal M}^{\bf t}_I
(S)\,|\,I\in{\cal N}\}
 \quad{\rm and}\\
 && \widehat{\cal M}^{\bf t}
(S^{(1)}|_{Z(S^{(2)})^\star})=\{\widehat{\cal M}^{\bf t}_I
(S^{(1)}|_{Z(S^{(2)})^\star})\,|\, I\in{\cal N}\}.
\end{eqnarray*}
We claim that
\begin{equation}\label{e:2.131}
\widehat{\cal M}^{\bf t}_I (S)=\widehat{\cal M}^{\bf t}_I
(S^{(1)}|_{Z(S^{(2)})^\star})\quad\forall I\in{\cal N}.
\end{equation}
In fact, for any $\hat x_I\in\widehat V^\ast_I$ one can derive from
Lemma~\ref{lem:2.53} that
$$
\hat S_I(\hat x_I)=\hat S_I^{(1)}(\hat x_I)\oplus\hat
S_I^{(2)}(\hat x_I)\in (\widehat F^\ast_I)_{\hat x_I}=(\widehat
F^{(1)\ast}_I)_{\hat x_I}\oplus(\widehat F^{(2)\ast}_I)_{\hat
x_I}
$$
and $(\hat\sigma^{(1)}_{ij})_I(\hat x_I)\in (\widehat
F^{(1)\ast}_I)_{\hat x_I}\oplus 0\subset (\widehat
F^{(1)\ast}_I)_{\hat x_I}\oplus(\widehat F^{(2)\ast}_I)_{\hat
x_I}=(\widehat F^{(1)\ast}_I)_{\hat x_I}$. So $\hat x_I\in
\widehat{\cal M}^{\bf t}_I (S)_I$, i.e.
$$
\hat S_I(\hat x_I)+
 \sum^n_{i=1}\sum^{m_i}_{j=1}t_{ij}(\hat\sigma_{ij}^{(1)})_I(\hat
 x_I)=0
$$
implies that $\hat S^{(2)}_I(\hat x_I)=0$ and
\begin{equation}\label{e:2.132}
 \hat S^{(1)}_I(\hat x_I)+
 \sum^n_{i=1}\sum^{m_i}_{j=1}t_{ij}(\hat\sigma_{ij}^{(1)})_I(\hat
 x_I)=0.
\end{equation}
It follows from the first equality that $\hat x_I\in\widehat
V^{z\ast}_I$. Note that the restriction of
$(\hat\sigma_{ij}^{(1)})_I$ to $\widehat V^{z\ast}_I$ is equal to
$(\hat\sigma_{ij}^{(1)z})_I$. We derive from (\ref{e:2.132}) that
$$
\hat S^{(1)}_I(\hat x_I)+
 \sum^n_{i=1}\sum^{m_i}_{j=1}t_{ij}(\hat\sigma_{ij}^{(1)z})_I(\hat
 x_I)=0.
$$
That is, $\hat x_I\in\widehat{\cal M}^{\bf t}_I
(S^{(1)}|_{Z(S^{(2)})})$. This shows that $\widehat{\cal M}^{\bf
t}_I\subset\widehat{\cal M}^{\bf t}_I
(S^{(1)}|_{Z(S^{(2)})^\star})$.

On the other hand $\hat S_I(\hat x_I)=\hat S^{(1)}_I(\hat x_I)$
and $(\hat\sigma_{ij}^{(1)z})_I(\hat
x_I)=(\hat\sigma_{ij}^{(1)})_I(\hat x_I)$ for any $\hat
x_I\in\widehat V^{z\ast}_I$. So it is clear that $\widehat{\cal
M}^{\bf t}_I (S^{(1)}|_{Z(S^{(2)})^\star})\subset\widehat{\cal
M}^{\bf t}_I$. (\ref{e:2.131}) is proved.

One can also easily prove that for any $J\subset I\in{\cal N}$ the
projections $\hat\pi^I_J$ and $\hat\pi_J^{zI}$ have the same domain
and image  as maps from $\widehat{\cal M}^{\bf t}_I
(S)=\widehat{\cal M}^{\bf t}_I (S^{(1)}|_{Z(S^{(2)})^\star})$ to
$\widehat{\cal M}^{\bf t}_J (S)=\widehat{\cal M}^{\bf t}_J
(S^{(1)}|_{Z(S^{(2)})^\star})$.  It follows that for each ${\bf
t}\in B_\varepsilon(\R^m)_{res}$ the corresponding virtual Euler
chains (or cycles)
\begin{eqnarray*}
&&e(E, S)^{\bf t}=``\sum_{I\in{\cal N}}"
\frac{1}{|\Gamma_I|}\{\hat\pi_I:\widehat{\cal M}^{\bf t}_I
(S)^\circ\to X\}\quad{\rm and}\\
&& e(E^{(1)}|_{Z(S^{(2)})^\star},
S^{(1)}|_{Z(S^{(2)})^\star})^{\bf
t}\\
&&=``\sum_{I\in{\cal N}}"
\frac{1}{|\Gamma_I|}\{\hat\pi_I^z:\widehat{\cal M}^{\bf t}_I
(S^{(1)}|_{Z(S^{(2)})^\star})^\circ\to Z(S^{(2)})^\star\}\;
\end{eqnarray*}
are same. The desired result is proved. \hfill$\Box$\vspace{2mm}

As a consequence we can get the following slight generalization of
Proposition 2.8 in \cite{R1}.

\begin{corollary}\label{cor:2.55}
 Let $X$ be a separable Banach orbifold satisfying
 Assumption~\ref{ass:2.45} and  $(X, E, S)$ be an oriented
 Banach Fredholm orbibundle of index $r$
  and with compact zero locus $Z(S)$. Assume:
\begin{description}
\item[(i)] $\dim{\rm Coker}DS(y)=k$ on $Z(S)$, $Z(S)$ is a closed
smooth orbifold of dimension $r+ k$ and thus ${\rm Coker}DS$ forms
an obstruction bundle ${\mathcal E}$ over $Z(S)$.

\item[(ii)] There exist an open neighborhood ${\cal W}$ of $Z(S)$
and a decomposition of the direct sum of Banach orbibundles
$$E|_{{\cal W}}=F\oplus F^c$$
such that $F^c|_{Z(S)}$ is isomorphic to ${\cal E}$ and that the
Banach orbisubbundle $F$ of corank $k$  restricts to ${\rm
Im}(DS)$ on $Z(S)$, i.e., $F_y={\rm Im}(DS(y))$ for any $y\in
Z(S)$.
\end{description}
Then there exist a virtual Euler cycle representation $M$ of
$e(E,S)$ and that $N$ of ${\cal E}$ such that $M=N$, and in
particular $ e(E,S)=e({\mathcal E})\cap [Z(S)]$.
 \end{corollary}

Indeed, as in the proof of Corollary~\ref{cor:1.14} we can use
Proposition~\ref{prop:2.54} to derive it.

 Similar to the arguments in \S 2.5.2 we can directly use Claim
2.43 and the restriction principle in Proposition~\ref{prop:2.50} to
prove the following property.

\begin{proposition}\label{prop:2.56}
({\bf Homotopy}) Let $(X, E^{(i)}, S^{(i)})$, $i=0,1$, be two
oriented Banach Fredholm bundles
 of index $r$ and with compact zero locus. Assume that $\partial X=\emptyset$ and
 $(X\times [0, 1], E, S)$ is an oriented homotopy between them.
 Then there exist virtual Euler chains of $(X\times [0, 1], E, S)$, $(X, E^{(0)}, S^{(0)})$
and $(X, E^{(1)}, S^{(1)})$ respectively, $e(E, S)^{\bf t}$,
$e(E^{(0)}, S^{(0)})^{{\bf t}_0}$ and $e(E^{(1)}, S^{(1)})^{{\bf
t}_1}$ such that
$$
\partial e(E, S)^{\bf t}=e(E^{(0)}, S^{(0)})^{{\bf t}_0}\cup
(-e(E^{(1)}, S^{(1)})^{{\bf t}_1}).
$$
In particular if  the virtual Euler classes of $(X, E^{(i)},
S^{(i)})$, $i=0,1$, also exist, then $e(E^{(0)},
 S^{(0)})=e(E^{(1)}, S^{(1)})$.
 \end{proposition}

Actually, in the theorem above,  for any $m\in\N$ if the
$m$-boundary $\partial^m X\ne\emptyset$ then $e(E, S)^{\bf t}$ has
 the $m$-boundary
 (defined in an obvious way),
$$
\partial^m e(E, S)^{\bf t}=\partial^{m-1}e(E^{(0)}, S^{(0)})^{{\bf t}_0}\cup
\partial^{m-1}e(E^{(1)}, S^{(1)})^{{\bf t}_1}\cup e(E, S)^{\bf t}_m,
$$
 where we neglect the orientation and
$$
e(E, S)^{\bf t}_m:=``\sum_{I\in{\cal
N}}"\frac{1}{|\Gamma_I|}\bigl\{\hat\pi_I:\widehat{\cal M}^{\bf
t}_I(S)^\circ\cap(\hat\pi_I)^{-1}(\partial^m X\times (0, 1))\to
X\times [0, 1]\bigr\}.
$$

Next we  study the {\bf functoriality} of the Euler class of the
Banach  Fredholm orbibundles. Let $(X, E, S)$ and $(X', E', S')$
be two oriented Banach Fredholm orbibundles with compact zero
loci.  A {\bf morphism} from $(X, E, S)$ to $(X', E', S')$ is an
orbibundle embedding $(\psi, \Psi):E\to E'$ with following
properties:
\begin{description}
\item[(A)] $S\circ \psi=\Psi\circ S'$ and $Z(S')=\psi(Z(S))$;

\item[(B)] For any $x\in Z(S)$, the differential $d\psi(x):T_xX\to
T_{\psi(x)}X'$ and the above restriction $\Psi_x: E_x\to
E'_{f(x)}$ induce isomorphisms
\begin{eqnarray*}
&&d\psi(x): {\rm Ker}(DS(x))\to {\rm Ker}(DS'(\psi(x)))\quad{\rm and}\\
&&[\Psi_x]: {\rm Coker}(DS(x))\to {\rm Coker}(DS'(\psi(x))),
\end{eqnarray*}
and the resulting isomorphism from ${\rm det}(DS)$ to ${\rm
det}(DS')$ is orientation preserving. (These are all understand on
levels of lifts.)
\end{description}

Clearly, $(X, E, S)$ and $(X', E', S')$ have the same index.
Corresponding with Proposition~\ref{prop:1.17} we have:

\begin{proposition}\label{prop:2.57}
Let $(X, E, S)$ and $(X', E', S')$ be two oriented Banach Fredholm
orbibundles with compact zero loci, and $(\psi, \Psi)$ be a morphism
from $(X, E, S)$ to $(X', E', S')$. Then there exist a family of
strongly cobordant virtual Euler chains of the triple $(X', E',
S')$,
$$
 e(E', S')^{{\bf
t}}=``\sum_{I\in{\cal N}}"
\frac{1}{|\Gamma'_I|}\{\hat\pi'_I:\widehat{\cal M}^{{\bf t}}_I
(S')^\circ\to X'\}\;\forall{\bf t}\in B_\varepsilon(\R^{m})_{res},
$$
and that of strongly cobordant virtual Euler chains of the triple
$(X, E, S)$,
$$
 e(E, S)^{{\bf
t}}=``\sum_{I\in{\cal N}}"
\frac{1}{|\Gamma'_I|}\{\hat\pi_I:\widehat{\cal M}^{{\bf t}}_I
(S)^\circ\to X\}\;\forall{\bf t}\in B_\varepsilon(\R^{m})_{res},
$$
and a family of embeddings $\{\widehat\psi_I\,|\, I\in{\cal N}\}$,
$\widehat\psi_I:\widehat W_I\to\widehat W'_I$, which are also
compatible with the projections $\hat\pi^I_J$ and $\hat\pi'^I_J$,
such that for any $I\in{\cal N}$ and $J\subset I$,
\begin{equation}\label{e:2.133}
\left.\begin{array}{ll}
\widehat\psi_I(\widehat{\cal M}^{\bf
t}_I(S)^\circ)=\widehat{\cal M}^{\bf t}_I(S')^\circ\cap(\hat\pi'_I)^{-1}(\psi(X)),\\
\widehat\psi_I(\widehat{\cal M}^{\bf
t}_I(S)^{sing})=\widehat{\cal M}^{\bf t}_I(S')^{sing}\cap(\hat\pi'_I)^{-1}(\psi(X)),\\
\widehat\psi_J\circ\hat\pi^I_J
 =\hat\pi'^I_J|_{\widehat\psi_I(\widehat{\cal M}^{\bf t}_I (S))}\circ \widehat\psi_I,\\
\Bigl(\hat\pi'^I_J\bigl((\hat \pi'^I_J)^{-1}(\widehat
V'^\ast_J)\cap\widehat{\cal M}^{\bf t}_I(S')\bigr)\Bigr)\cap
(\hat\pi'_J)^{-1}(\psi(X))\\
=\widehat\psi_I\Bigl(\hat\pi^I_J\bigl((\hat \pi^I_J)^{-1}(\widehat
V_J^\ast)\cap\widehat{\cal M}^{\bf
t}_I(S)\bigr)\Bigr),\\
\Bigl({\rm Im}(\hat\pi'^I_J)\cap\widehat{\cal M}^{\bf
t}_J(S')\Bigr)\cap(\hat\pi'_J)^{-1}(\psi(X))\\
=\widehat\psi_I\bigl({\rm Im}(\hat\pi^I_J)\cap\widehat{\cal
M}^{\bf t}_J(S)\bigr)
\end{array}
\right\}
\end{equation}
If $\overline W$ and $\overline W'$ are the spaces constructed
from $\{\widehat W_I\,|\, I\in{\cal N}\}$ and $\{\widehat
W'_I\,|\, I\in{\cal N}\}$ respectively then $\{\widehat\psi_I\,|\,
I\in{\cal N}\}$ induce a natural map $\bar\psi:\overline
W\to\overline W'$ such that
\begin{equation}\label{e:2.134}
\left.\begin{array}{ll}\bar\psi(\overline{\cal M}^{\bf
t}(S)^\circ)=\overline{\cal M}^{\bf t}(S')^\circ\cap
(\hat\pi')^{-1}(\psi(X)),\\
\bar\psi\circ\hat\pi=\hat\pi'|_{\bar\psi(\overline{\cal M}^{\bf
t}(S)^\circ)}\circ\bar\psi
\end{array}
\right\}
\end{equation}
Consequently, if $e(E,S)$ and $e(E', S')$ exist,
$\langle\psi^\ast\alpha, e(E,S)\rangle= \langle\alpha, e(E',
S')\rangle$ for any $\alpha\in H^\ast(X', X'\setminus
\psi(X);\R)$. (In this case we say $\psi_\ast e(E, S)$ to be the
intersection  of $e(E', S')$ with $\psi(X)$.)
\end{proposition}

\noindent{\bf Proof.}\quad By the definition of the morphism we may
assume: {\bf $X$ is a Banach suborbifold of $X'$, $E$ is a
suborbibundle  of $E'|_X$ and $S=S'|_X$}. Of course, both $f$ and
$F$ are inclusion maps. In particular, by (A) we have $Z(S')=Z(S)$.
As in the proof of Proposition~\ref{prop:2.51} we can choose
finitely many points
$$
x_i\in Z(S)=Z(S'), i=1,\cdots, n,
 $$
 and orbifold charts $(\widetilde W'_i, \Gamma'_i, \pi'_i)$ on $X'$
centered at $x_i$, $i=1,\cdots, n$ such that:
\begin{description}
\item[($1^\circ$)] For each $i$ there exists a submanifold
$\widetilde W_i\subset\widetilde W'_i$ that is stable under
$\Gamma'_i$ and compatible with the inclusion map, such that the
restriction $(\widetilde W_i, \Gamma'_i, \pi_i)=(\widetilde W_i,
\Gamma'_i|_{\widetilde W_i}, \pi'_i|_{\widetilde W_i})$ is a Banach
orbifold chart on $X$. (see the explanations in the proof of
Proposition~\ref{prop:2.51} for the reason that
$\Gamma'_i=\Gamma'_i|_{\widetilde W_i}$.)

\item[($2^\circ$)] $\{W_i=\pi_i(\widetilde W_i)\}^n_{i=1}$ is an
open cover of $Z(S)$ in $X$ and thus  $\{W'_i=\pi'_i(\widetilde
W'_i)\}^n_{i=1}$ is such a cover of $Z(S')$ in $X'$.
\end{description}
 Similarly, for each $i$ let $(\widetilde
E'_i,\Gamma'_i,\Pi'_i)$ be the orbifold bundle chart on $E'$
corresponding with $(\widetilde W'_i, \Gamma'_i, \pi'_i)$, then
there exists a suborbibundle $\widetilde E_i\subset\widetilde
E'_i|_{W_i}$ that is stable under $\Gamma'_i$ and compatible with
the inclusion map, such that the restriction $(\widetilde E_i,
\Gamma'_i, \Pi_i)=(\widetilde E_i, \Gamma'_i|_{\widetilde E_i},
\Pi'_i|_{\widetilde E_i})$ is a Banach orbibundle chart on $E$
which is corresponding with $(\widetilde W_i, \Gamma'_i, \pi_i)$.

As usual take pairs of open sets
 $W'^j_i\subset\subset U'^j_i$, $i=1,\cdots, n,\; j=1, 2, \cdots,n-1$,
 such that $Z(S)\subset\cup^n_{i=1}W'^1_i$ and that
$$
 U'^1_i\subset\subset
W'^2_i\subset\subset U'^2_i\cdots \subset\subset
W'^{n-1}_i\subset\subset U'^{n-1}_i\subset\subset W'_i.
$$
Set  $W^j_i:=W'^j_i\cap X$ and $U^j_i:=U'^j_i\cap X$ then they are
 pairs of open sets in $X$ and $Z(S)\subset\cup^n_{i=1}W_i^1$.
Moreover, for $i=1,\cdots, n,\; j=1, 2, \cdots,n-1$ it also holds
that $W^j_i\subset\subset U^j_i\subset\subset W_i$ and
$$
 U^1_i\subset\subset
W^2_i\subset\subset U^2_i\cdots \subset\subset
W^{n-1}_i\subset\subset U^{n-1}_i\subset\subset W_i.
$$
As before we can use these to construct the corresponding bundle
systems $(\widehat{\cal E}, \widehat W)$, $(\widehat{\cal F},
\widehat V)$ $(\widehat{\cal E}', \widehat W')$, $(\widehat{\cal
F}', \widehat V')$ and
\begin{eqnarray*}
&&\bigl(\widehat{\mathcal F}^{\ast}, \widehat V^{\ast}\bigr):=
\bigl\{\bigl(\widehat F_I^{\ast}, \widehat V_I^{\ast}\bigr),
\hat\pi_I, \hat\Pi_I, \Gamma'_I, \hat\pi^{I}_J,\hat\Pi^{I}_J,
\lambda'^I_J\,\bigm|\, J\subset I\in{\mathcal N}\bigr\},\\
&&\bigl(\widehat{\mathcal F}'^{\ast}, \widehat V'^{\ast}\bigr):=
\bigl\{\bigl(\widehat F'^{\ast}_I, \widehat V'^{\ast}_I\bigr),
\hat\pi'_I, \hat\Pi'_I, \Gamma'_I, \hat\pi'^{I}_J,\hat\Pi'^{I}_J,
\lambda'^{I}_J\,\bigm|\, J\subset I\in{\mathcal N}\bigr\},
\end{eqnarray*}
where the small neighborhoods ${\cal W}^{\ast}$ of $Z(S)$ in $X$
is taken as ${\cal W}'^\ast\cap X$.

For each $i=1,\cdots, n$, as before we take $\Gamma'_i$-invariant
smooth cut-off function $\gamma'_i:\widetilde W'_i\to [0,1]$ with
support in $\widetilde W'^1_i$,  such that
$$
Z(S')\subset\cup^n_{i=1}U'^{0}_i=\cup^n_{i=1}\pi'_i(\widetilde
U'^{0}_i)
$$
for $\widetilde U'^{0}_i:=\{\tilde x\in\widetilde W'_i\,|\,
\gamma'_i(\tilde x)>0\}$. Then for $U^0_i:=U'^0_i\cap X$,
$i=1,\cdots,n$, it holds that $Z(S)\subset\cup^n_{i=1}U^{0}_i$.

Let $\widetilde S_i:\widetilde W_i\to\widetilde E_i$ (resp.
$\widetilde S'_i:\widetilde W'_i\to\widetilde E'_i$) be the lift
of $S|_{W_i}$ (resp. $S'|_{W'_i}$). Then by (A) and (B) it also
holds that $Z(\widetilde S'_i)=Z(\widetilde S_i)$, and that for
any $\tilde x\in Z(\widetilde S_i)$, $D\widetilde S'_i(\tilde
x)|_{T_{\tilde x}\widetilde W_i}=D\widetilde S_i(\tilde x)$, ${\rm
Ker}(D\widetilde S_i(\tilde x))={\rm Ker}(D\widetilde S'_i(\tilde
x))$ and the inclusion $(\widetilde E_i)_{\tilde x}\hookrightarrow
(\widetilde E_i')_{\tilde x}$ induces an isomorphism
\begin{eqnarray*}
 &&(\widetilde E_i)_{\tilde x}/D{\widetilde S_i}(\tilde x)(T_{\tilde x}\widetilde W_)
 \to (\widetilde E'_i)_{\tilde x}/D(\widetilde S'_i)(\tilde x)(T_{\tilde x}\widetilde W_i'),\\
&&  \tilde v+ D\widetilde S_i({\tilde x})(T_{\tilde x}\widetilde
W_i)\mapsto \tilde v + D\widetilde S_i'(\tilde x)(T_{\tilde
x}\widetilde W_i').
\end{eqnarray*}
As a consequence we get:

\begin{claim}\label{cl:2.58}
 Let $\tilde x\in Z(\widetilde S'_i)=Z(\widetilde S_i)$ and $\tilde
v_i\in (\widetilde E_i)_{\tilde x}, i=1,\cdots, k$, satisfy
$$
D\widetilde S_i(\tilde x)(T_{\tilde x}\widetilde W_i)+ {\rm
span}\{\tilde v_1,\cdots, \tilde v_k\}=(\widetilde E_i)_{\tilde
x}.
$$
Then $D\widetilde S'_i(\tilde x)(T_{\tilde x}\widetilde W_i')+ {\rm
span}\{\tilde v_1,\cdots, \tilde v_k\}=(\widetilde E'_i)_{\tilde
x}$.
 \end{claim}

For the chosen trivializations
$$
{\cal T}'_i: \widetilde W'_i\times (\widetilde E'_i)_{\tilde
x_i}\to \widetilde E'_i,\;i=1,\cdots, n,
$$
 they restrict to natural ones
$$
 {\cal T}_i: \widetilde
W_i\times (\widetilde E_i)_{\tilde x_i}\to \widetilde E_i,\;
i=1,\cdots, n.
$$
By increasing $n$ and shrinking $W'_i$ one can take $\tilde
v_{ij}\in (\widetilde E_i)_{\tilde x_i}$, $j=1,\cdots, m_i$ such
that for any $\tilde x\in\widetilde W_i$ and $\tilde s_{ij}(\tilde
x):={\cal T}_i(\tilde x, \tilde v_{ij})$, $j=1,\cdots, m_i$,
\begin{equation}\label{e:2.135}
D\widetilde S_i(\tilde x)(T_{\tilde x}\widetilde W_i)+ {\rm
span}\{\tilde s_{i1}(\tilde x),\cdots, \tilde s_{im_i}(\tilde
x)\}=(\widetilde E_i)_{\tilde x}.
\end{equation}
Note that for any $\tilde x\in\widetilde W_i$, $\tilde
s'_{ij}(\tilde x):={\cal T}'_i(\tilde x, \tilde v_{ij})$ are equal
to $\tilde s_{ij}(\tilde x)$ for $i=1,\cdots, m_i$. By
Claim~\ref{cl:2.58} we get that for any $\tilde x\in\widetilde
W'_i$,
\begin{equation}\label{e:2.136}
D\widetilde S'_i(\tilde
x)(T_{\tilde x}\widetilde W_i')+ {\rm span}\{\tilde s'_{i1}(\tilde
x),\cdots, \tilde s'_{im_i}(\tilde x)\}=(\widetilde E'_i)_{\tilde
x}.
\end{equation}

Consider the sections $\sigma'_{ij}:=\gamma'_i\cdot\tilde s'_{ij}$
of $\widetilde E'_i$ with supports in $\widetilde W'^1_i$, and those
$\sigma_{ij}:=\gamma'_i\cdot\tilde s_{ij}$ of $\widetilde E_i$ with
supports in $\widetilde W^1_i$, $j=1,\cdots, m_i$. Clearly,
$\sigma_{i1}=\sigma'_{i1}|_{\widetilde W_i}, \cdots,
\sigma_{im_i}=\sigma'_{im_i}|_{\widetilde W_i}$, $i=1,\cdots, n$.
Assume that the above small neighborhood ${\cal W}'^\ast$ of $Z(S')$
in $X$ is contained in $\cup^n_{i=1}U'^0_i$ (and thus ${\cal
W}^{\ast}={\cal W}'^\ast\cap X\subset\cup^n_{i=1}U^0_i$). By
(\ref{e:2.135}) and (\ref{e:2.136}), as before we can shrink ${\cal
W}'^\ast$  and find a small $\varepsilon>0$ such that for
 each ${\bf t}=\{t_{ij}\}_{\substack{ 1\le j\le m_i\\
 1\le i\le n}}$ in a residual subset
$B_\varepsilon(\R^m)_{res}\subset B_\varepsilon(\R^m)$ with
$m=\sum^n_{i=1}m_i$ the sections sections
\begin{eqnarray*}
 \Upsilon^{({\bf t})}_I: \widehat
V^\ast_I\to \widehat F^\ast_I,\; \hat x_I \mapsto\hat S_I(\hat
x_I)+
 \sum^n_{i=1}\sum^{m_i}_{j=1}t_{ij}(\hat\sigma_{ij})_I(\hat x_I),\\
  \Upsilon'^{({\bf t})}_{I}: \widehat
V'^{\ast}_I\to \widehat F'^{\ast}_I,\; \hat x_I \mapsto\hat
S'_I(\hat x_I)+
 \sum^{n}_{i=1}\sum^{m_i}_{j=1}t_{ij}(\hat\sigma'_{ij})_I(\hat
 x_I),
\end{eqnarray*}
are strongly transversal to the zero section for all $I\in{\cal N}$.

One also easily checks that $\widehat F_I^\ast$ is a Banach
subbundle of $\widehat F'^\ast_I|_{\widehat V_I^\ast}$,
$$
\widehat V_I^\ast=\widehat V'^\ast_I\cap(\hat\pi'_I)^{-1}(X),\quad
\widehat F_I^\ast=\widehat F'^\ast_I\cap(\hat\Pi'_I)^{-1}(E),\quad
\Upsilon^{({\bf t})}_I=\Upsilon'^{({\bf t})}_{I}|_{\widehat
V_I^\ast}
$$
and  these relation are also compatible with $\hat\pi^{I}_J$,
$\hat\Pi^{I}_J$ and $\hat\pi'^{I}_J$, $\hat\Pi'^{I}_J$. As usual
let $\widehat{\cal M}^{\bf t}_I(S)=(\Upsilon^{({\bf
t})}_I)^{-1}(0)$ and $\widehat{\cal M}^{\bf t}_I(S')
=(\Upsilon'^{({\bf t})}_{I})^{-1}(0)$. Then (\ref{e:2.133}) and
(\ref{e:2.134}) easily follows from these.
 \hfill$\Box$\vspace{2mm}

  Note that the family $\widehat{\cal M}^{\bf t}(S)=\{\widehat{\cal M}^{\bf t}_I(S)\,|\, I\in{\cal N}\}$
may be viewed a `subset' of the system
$$
 {\widehat V}=\bigl \{\widehat V_I, \hat\pi_I, \hat\pi^I_J,
 \hat p_I,\Gamma_I,\,\lambda^I_J\,\bigm|\, J\subset I\in{\mathcal
 N}\bigr\}
 $$
( a resolution of $X$ near $Z(S)$) and the Euler chain $e(E,S)^{\bf
t}$ may be viewed as a smooth map from $\widehat{\cal M}^{\bf t}(S)$
to $X$.  One can also consider the Euler class of pull-back Banach
Fredholm orbibundles by Fredholm maps. It is possible for us to get
the following generalization of Proposition~\ref{prop:1.18}.

\begin{proposition}\label{prop:2.59}
Let $X$ and $X'$ be two separable  Banach orbifolds satisfying
Assumption~\ref{ass:2.45}, and $(X, E, S)$ (resp. $(X', E', S')$) be
the Banach Fredholm orbibundle of index $r$ (resp. $r'$). Let
$f:X\to X'$ be a proper Fredholm map of index $d$. Assume that that
$\bar f:(X, E, S)\to (X', E', S')$ is a bundle map covering $f$,
i.e. a bundle map $\bar f: E\to E'$ satisfying $S'\circ f=\bar
f\circ S$, and  that each local lifting of $\bar f$ is a Banach
space isomorphism when restricted to each fiber. (That is, $(X, E,
S)$ is the pull-back of $(X', E', S')$ via $f$). Furthermore assume
that $f$ is also a submersion. Then one has:
\begin{description}
 \item[(i)] $r=r'+d$ and $Z(S')=f^{-1}(Z(S))$
is compact. \item[(ii)] There exist a resolution $\widehat V$ of
$X$ near $Z(S)$ and that $\widehat V'$ of $X'$ near $Z(S')$, and a
natural map $\hat f:\widehat V\to\widehat V'$ induced by $f$,  and
the virtual Euler chain $e(E, S)^{\bf t}$ (resp. $e(E', S')^{\bf
t}$), as a smooth map from $\widehat{\cal M}^{\bf
t}(S)\subset\widehat V$ (resp. $\widehat{\cal M}^{\bf
t}(S')\subset\widehat V'$) to $X$ (resp. $X'$) satisfy
\begin{equation}\label{e:2.137}
f\circ e(E, S)^{\bf t}=e(E', S')^{\bf t}\circ\hat f.
\end{equation}
\end{description}
\end{proposition}

The proof of Proposition~\ref{prop:2.59} can be completed by
slightly modifying the proof of Proposition~\ref{prop:4.15}. We omit
them.

\begin{remark}\label{rm:2.60}
{\rm We conjecture that there exists some kind $\Q$-value
topological degree ${\rm deg}(f)$ of $f$ so that (\ref{e:2.137})
gives
$$
f_\ast e(E, S)={\rm deg}(f) e(E', S')
$$
in the case $d=0$ and both $e(E, S)$ and $e(E', S')$ exist.}
\end{remark}

Finally, corresponding to Proposition~\ref{prop:1.15}, one may {\bf
expect} the following:

{\it  Let $(X, E, S)$ and $(X', E', S')$ be two Banach Fredholm
orbibundles with compact zero loci $Z(S)$ and $Z(S')$. $r$ and $r'$
are their indexes respectively.  Then the natural product $(X\times
X', E\times E', S\times S')$ is such a bundle of index $r''=r+ r'$
and with zero locus $Z(S\times S')=Z(S)\times Z(S')$, and there
exist virtual Euler chains $C$ of $(X, E, S)$, $C'$ of $(X', E',
S')$ such that the product $C''=C\times C'$ is a Euler chain of
$(X\times X', E\times E', S\times S')$. Furthermore, if
Assumption~\ref{ass:2.45} is also satisfied for $(X, E, S)$ and
$(X', E', S')$ (so that $e(E, S)$ and $e(E', S')$ exist),  then
 $e(E\times E', S\times S')$ also exists, and is equal to
$e(E, S)\times e(E', S')$.}

We cannot prove it here. The problem can be explained briefly as
follows: Let $(\widetilde W_i, \Gamma_i, \pi_i)$ be  the orbifold
charts on $X$ such that $Z(S)\subset\cup^n_{i=1}W_i$, and
$(\widetilde E_i,\Gamma_i,\Pi_i)$ be the corresponding orbifold
bundle charts on $E$ with $(\widetilde W_i, \Gamma_i, \pi_i)$,
$i=1,\cdots, n$; moreover,  let $(\widetilde W'_i, \Gamma'_i,
\pi'_i)$ be the orbifold charts on $X'$ such that
$Z(S')\subset\cup^n_{i=1}W'_i$, and $(\widetilde
E'_i,\Gamma'_i,\Pi'_i)$  be the corresponding orbifold bundle charts
on $E'$ with $(\widetilde W'_i, \Gamma'_i, \pi'_i)$, $i=1,\cdots,
n'$. For $(X\times X', E\times E', S\times S')$ we naturally choose
the orbifold charts $(\widetilde W_i\times\widetilde W'_k,
\Gamma_i\times\Gamma'_k, \pi_i\times\pi'_k)$ on $X\times X'$, and
the corresponding orbifold bundle charts $(\widetilde
E_i\times\widetilde E'_k, \Gamma_i\times\Gamma'_k,
\Pi_i\times\Pi'_k)$  on $E\times E'$ with $(\widetilde
W_i\times\widetilde W'_k, \Gamma_i\times\Gamma'_k,
\pi_i\times\pi'_k)$, $i=1,\cdots, n, k=1,\cdots, n'$. Clearly,
$Z(S\times S')\subset\cup^n_{i=1}\cup^{n'}_{k=1}W_i\times W'_k$.
When we start from $(\widetilde W_i\times\widetilde W'_k,
\Gamma_i\times\Gamma'_k, \pi_i\times\pi'_k)$ and $(\widetilde
E_i\times\widetilde E'_k, \Gamma_i\times\Gamma'_k,
\Pi_i\times\Pi'_k)$ and attempt to use the standard way to construct
a system of virtual Banach vector bundles, {\bf the first difficulty
we meet} is: not each index ${\bf I}\in{\cal N}\times{\cal N}'$ can
be expressed as form $I\times I'$, $I\in{\cal N}$ and $I'\in{\cal
N}'$. Here ${\cal N}$ and ${\cal N}'$ are respectively corresponding
sets of index with $\{\widetilde W_i\}^n_{i=1}$ and $\{\widetilde
W'_i\}^{n'}_{i=1}$.

 \pagebreak

\section{The Euler cycle of the partially smooth  Banach Fredholm
bundle}\label{sec:3}

\setcounter{equation}{0}

We begin with the introduction of the notion of partially smooth
Banach manifolds and Banach vector bundles and then extend
Theorem~\ref{th:1.5} to the partially smooth category. Here one
needs to overcome  another kind of difficulty. We shall abstract the
most essential properties which can guarantee the expected
conclusions. That is, Theorem~\ref{th:1.5} can be extended to a
class of  the partially smooth  Banach Fredholm bundles. From now on
we abbreviate ``{\bf partially smooth}" as ``{\bf PS}". Moreover, we
consider the boundary as a special case of the corner.

\begin{definition}\label{def:3.1}
{\rm A Hausdorff topological space $X$ is called a {\bf PS Banach
manifold (with corner)} if there exists a decomposition
 $X=\cup_{\alpha\in\Lambda} X_\alpha$ indexed by a set
$\Lambda$ with a partial order $\prec$ with $0$ as the largest
element such that the following two conditions are satisfied:
\begin{description}
\item[(i)] Any two $X_\alpha$ and $X_\beta$ with $\alpha\ne\beta$ are
disjoint, and each point of $X$ has a neighborhood that may only
intersect with finitely many $X_\alpha$.
\item[(ii)] Each $X_\alpha$, called a {\bf stratum}, is a smooth
Banach manifold with corner whose topology assures the inclusion
$X_\alpha\hookrightarrow X$ to be continuous.
\item[(iii)] The stratum $X_0$,  called the {\bf top stratum} of $X$,
 is a nonempty open subset in $X$. Other strata cannot contain any interior point of $X$,
  and are called {\bf lower strata}.
\item[(iv)]  For each $\alpha\in\Lambda$ the stratum
$X_\alpha$ is relatively open in $\cup_{\beta\prec\alpha}X_\beta$
with respect to the induced topology from $X$. For any
$\beta\in\Lambda$ with $\beta\prec\alpha$ and $\beta\ne\alpha$, the
stratum $X_\beta$ cannot contain any relative interior point of
$X_\alpha$ in $\cup_{\beta\prec\alpha}X_\beta$.
\end{description}
If only conditions (i) and (ii) are satisfied, $X$ is called a {\bf
weakly PS Banach manifold}. (In this case $\Lambda$ may not required
to have a partial order.) For $k\in\N\cup\{0\}$, a {\bf weakly
$k$-PS manifold of dimension $r$} is a weakly PS manifold whose
lower strata have dimension less than $r-k$.  }
\end{definition}

\begin{remark}\label{rem:3.2}
{\rm ($1^\circ$) Actually  one usually makes arguments near a
compact subset. So we may assume that our PS Banach manifold has
only finitely many strata in the following arguments.  The
condition (ii) implies that the Banach manifold topology of each
stratum is either the same as the induced topology from $X$
or stronger than the induced topology.\\
($2^\circ$) The condition (iii) implies that
$$
X_0^\infty:=\bigcap_{K\subset X_0\;{\rm compact}}Cl(X_0\setminus
K)
$$
is contained in the union of all lower strata of $X$, where
$Cl(X_0\setminus K)$ is the closure of $X_0\setminus K$ in $X$.
(The condition (iv) implies a similar remark for $X_\alpha$.) Let
us denote by
$$
\partial_t X\quad{\rm (resp.}\quad\partial_t X_0\,)
$$
 the {\bf topological boundary} of $X$ (resp. $X_0$), and by
 $$
\partial X\quad{\rm (resp.}\quad\partial X_0\,)
$$
 the {\bf boundary PS Banach manifold of $X$}
(resp. the boundary Banach manifold of $X_0$). That is, $\partial
X=\cup_\alpha\partial X_\alpha=\cup_\alpha\partial^1 X_\alpha$,
where $\partial^1 X_\alpha$ is the (codimension one) boundary Banach
submanifold of $X_\alpha$, see Definition~\ref{def:A.6} for the
precise sense of $\partial^1$.
 Then $\partial_t
X=\partial_t X_0$ consists of two parts,  the union of the lower
strata of $X$, and the union of all corner submanifolds of $X_0$.
Thus  $\partial_t X=\partial_t X_0=X_0^\infty$ if $X_0$ has no
corners (including boundary). {\it Warning:} a PS Banach manifold
without corner
might still have nonempty topological boundary.\\
($3^\circ$)  A Banach manifold with corner may be viewed as a PS
Banach manifold. The lower strata of it consist of the various
dimensional corners of $X$. However, by Definition~\ref{def:A.4} one
easily sees that a PS Banach manifold is not necessarily a Banach
manifold with corner. They are different two notions.}
\end{remark}

A topological subspace $Z$ of $X$ is called a {\bf PS Banach
submanifold} of $X$ if there exists a partial order subset
$\Lambda_0\subset\Lambda$ with maximal element $\alpha_0$ such that
(i) $Z=\cup_{\alpha\in\Lambda_0}Z_\alpha$, where  $Z_\alpha=Z\cap
X_\alpha$ is a Banach submanifold of $X_\alpha$ for each
$\alpha\in\Lambda_0$, (ii) $Z=\cup_{\alpha\in\Lambda_0}Z_\alpha$ is
also a PS Banach manifold with respect to the induced stratified
smooth structure. Clearly, if $X=\cup_{\alpha\in\Lambda}X_\alpha$ is
a PS Banach manifold, for each $\alpha\in\Lambda$ the subset
$\cup_{\beta\prec\alpha}X_\beta$ is a PS submanifold of $X$. A PS
Banach submanifold $Z=\cup_{\alpha\Lambda_0}Z_\alpha$ as above is
called {\bf completely neat} if for any $\alpha\in\Lambda_0$ and
$x\in Z_\alpha$ it holds that ${\rm ind}_{Z_\alpha}(x)={\rm
ind}_{X_\alpha}(x)$ and $B_kZ_\alpha=Z_\alpha\cap B_kX_\alpha$ for
all $k\in\N\cup\{0\}$ (cf. Appendix~\ref{app:A} for the definitions
of ${\rm ind}_{Z_\alpha}(x)$ and $B_kZ_\alpha$.)
 For $r\in\N\cup\{0\}$, a weakly PS
Banach manifold $X$ is said to be {\bf $r$-dimensional} if each
stratum of it has finite dimension and $r$ is the maximum of
dimensions of all strata of it. (In this case the stratum of
dimension $r$ is called its {\bf top stratum}.) In particular, we
say a PS manifold to be {\bf $r$-dimensional} if the top stratum
of it has dimension $r$. A (weakly) PS manifold of dimension $r$
is said to be {\bf oriented} if its top stratum is oriented. A
topological subspace $Z$ of a weakly PS manifold $X$ is called a
{\bf weakly PS submanifold} if it is a weakly PS manifold with
respect to the induced topology structure and stratified smooth
structure from $X$. For $k\in\N\cup\{0\}$, a {\bf weakly $k$-PS
submanifold of dimension $r$} is a weakly PS submanifold whose
lower strata have dimension less than $r-k$.  That is, this weakly
PS submanifold itself is a weakly $k$-PS manifold according to the
induced structure. It is obvious that
$$
 {\rm weakly}\;
k\hbox{-PS}\;{\rm submanifold}\Longrightarrow{\rm weakly}\;{\rm
PS}\;{\rm submanifold}.
$$
 A compact (weakly) PS manifold $X$  (resp. a compact weakly $k$-PS manifold
 $X$) of dimension $r$  is called
{\bf closed} if (i)  its top stratum is a manifold without corner
(including boundary), and (ii) all strata of dimension lower than
$r$ have at most dimensions $r-2$. Thus a compact weakly $1$-PS
manifold $X$ of dimension $r$ is closed if its top stratum has no
corner (including boundary). (When the top stratum of a compact
(weakly) $k$-PS submanifold has no corner (including boundary) we
also say it to be {\bf closed}.) Clearly, one has

\begin{claim}\label{cl:3.3}
{\rm Any closed weakly PS Banach manifold $X$ of dimension $r$
carries a fundamental class in $H_r(X,\Z_2)$, or $H_r(X,\Z)$ if $X$
is also oriented.}
\end{claim}

 In particular, a compact (weakly) PS manifold of
dimension zero is a finite set. A (weakly) PS Banach manifold $X$ is
called {\bf separable} if each stratum $X_\alpha$ is separable.

 For $x\in X$,
 if $X_\alpha$ is the stratum of $X$ containing it, we write $T_x X=T_x X_\alpha$.
 For $0<k\le\infty$, a {\bf PS $C^k$-map} $f$ from the PS Banach manifold
$X$ to $Y$ is a continuous map $f:X\to Y$ such that the
restriction of $f$ to any stratum $X_\alpha$,
$f_\alpha=f|_\alpha$, is a $C^k$-smooth map into a certain stratum
$Y_\beta$ of $Y$. If $f$ is a topological embedding and the
restriction of it to any stratum $X_\alpha$ is a smooth embedding
into some stratum $Y_\beta$ of $Y$, then $f$ is called a {\bf PS
embedding}. In particular, if the image of $f$ is also open  in
$Y$  we say $f$ to be a {\bf PS open embedding}. When $f$ is a
homeomorphism and both $f$ and $f^{-1}$ are PS embeddings we call
$f$ a {\bf PS-diffeomorphism} from $X$ to $Y$. In particular the
PS-diffeomorphism from $X$ to itself is called a {\bf PS
automorphism} of $X$. Denote by ${\rm Aut}(X)$ the group of
automorphisms of $X$.

\begin{definition}\label{def:3.4}
{\rm Let both $X$ and $E$ be two PS Banach manifolds. A surjective
PS map $p:E\to X$ is called a {\bf PS Banach vector bundle} if for
the stratum decomposition $X=\bigcup_{\alpha\in\Lambda}X_\alpha$ one
has:
\begin{description}

\item[(i)] The inverse image $p^{-1}(X_\alpha)$ of each stratum
$X_\alpha$ of $X$ is exactly a stratum of $E$.

\item[(ii)] The restriction $p|_{p^{-1}(X_\alpha)}:
p^{-1}(X_\alpha)\to X_\alpha$ is a smooth Banach vector bundle.

\item[(iii)] $E$ has the stratum decomposition
$E=\bigcup_{\alpha\in\Lambda}p^{-1}(X_\alpha)$.
\end{description}}
\end{definition}

Naturally, one may define a {\bf weakly PS Banach bundle} over a
weakly PS Banach manifold.

 Corresponding with smooth bundle maps we may define a notion of
a {\bf PS bundle map}. In particular, for a PS bundle map $(f, F):
(E, p, X)\to (E^\prime, p^\prime, X^\prime)$, if both $f:X\to
X^\prime$ and $F:E\to E^\prime$ are PS open embeddings and for each
$x\in X$ the restriction $F|_{E_x}: E_x\to E^\prime_{f(x)}$ is a
Banach space isomorphism, we call $(f, F)$ a {\bf PS bundle open
embedding} from $(E, p, X)$ to $(E^\prime, p^\prime, X^\prime)$.
Similarly, we have also the notions of  {\bf PS bundle
isomorphisms}, {\bf PS bundle automorphisms}, {\bf PS subbundles},
and so on.  In particular, for an open subset $U\subset X$, a PS
finite rank equi-rank subbundle of $E|_U$ is called a {\bf local PS
finite rank equi-rank subbundle of} $E|_U$. A {\bf PS  section} of
$E$ is a continuous section $S:X\to E$ whose restriction on each
stratum $X_\alpha$ of $X$ gives a smooth section
$S|_{X_\alpha}:{X_\alpha}\to E_\alpha=E|_{X_\alpha}$. Such a section
is called  {\bf transversal to the zero section} if the restriction
$S|_{X_\alpha}$ of it to each stratum $X_\alpha$ is transversal to
the zero section. (In this case, we can only assure that the zero
locus $S^{-1}(0)$ is a weakly PS Banach manifold, rather than a PS
Banach manifold. This is a reason that we introduce a weakly PS
Banach manifold!) Denote by $\Gamma(E)$ the space of all PS sections
of $E$. It might be empty when $X$ is not a Banach manifold.

\begin{definition}\label{def:3.5}
{\rm  Let $E\to X$ be a PS bundle.  A PS section $S:X\to E$ is
called {\bf weakly Fredholm} if the restriction of it to each
stratum of $X$ is a strongly Fredholm section in the sense above
Lemma~\ref{lem:1.1}. The {\bf index} of such a section is defined as
the maximum of Fredholm indexes of restrictions of it to all strata.
For an integer $k\ge 0$, a weakly Fredholm PS section $S:X\to E$ is
said to be  a {\bf weakly $k$-Fredholm PS section} if
$$
 {\rm Index}(S|_{X_\alpha})+ k<{\rm
Index}(S|_{X_0})
$$
for any  $\alpha\in \Lambda$  with $\alpha\ne 0$. Furthermore, this
weakly Fredholm PS section is called a {\bf $k$-Fredholm PS section}
if
$$
 {\rm Index}(S|_{X_\beta})+ k<{\rm
Index}(S|_{X_\alpha})
$$
for any two different $\alpha,\beta\in\Lambda$ with
$\beta\prec\alpha$. More generally, for a family of negative
integers ${\bf
k}:=\{k_{\beta\alpha}\,|\,\beta\prec\alpha,\,\beta,\alpha\in\Lambda\}$,
the above weakly Fredholm PS section $S$ is called a ${\bf k}$-{\bf
Fredholm PS section} if
$$
 {\rm Index}(S|_{X_\beta})+ k_{\beta\alpha}<{\rm
Index}(S|_{X_\alpha})
$$
for any two different $\alpha,\beta\in\Lambda$ with
$\beta\prec\alpha$.}
\end{definition}

By the definition, a $k$-Fredholm PS section is a weakly
$k$-Fredholm PS section, and the latter is a weakly Fredholm PS
section. Moreover, it is also clear that  the index of the
restriction of a (weakly) $k$-Fredholm PS section to the top stratum
is at least greater $k$ than the indexes of the restrictions of it
to the lower strata.
 It may be considered as a generalization of the notion of the
strongly Fredholm section in Section~\ref{sec:1}.  In the case of
Gromov-Witten invariants the corresponding section is $1$-Fredholm.

\begin{remark}\label{rem:3.6}
{\rm In Definition~\ref{def:3.1}, one may define a PS Banach
manifold by the conditions (i)-(iii) and does not require
$\Lambda$ to be a partial order set. In this case we can define
the notion of the $k$-Fredholm section as the above definition of
the weakly $k$-Fredholm section. According to such definitions all
arguments in this section (and next one after some corresponding
changes have been made) still hold. However, such a definition of
a $k$-Fredholm section can not give a relation among indexes of
the restrictions to different lower strata, which is needed for
some possible problems. Perhaps, in this case the following
definition of a $k$-Fredholm section is a good choice. Under the
above representatives, it holds that ${\rm Index}(S|_{X_\beta})+
k<{\rm Index}(S|_{X_\alpha}) $ if $\alpha,\beta\in\Lambda$ are
such that $X_\alpha$ is relatively open in $X_\alpha\cup X_\beta$
and $X_\beta$ does not contain any relative interior points of
$X_\alpha\cup X_\beta$. }
\end{remark}

\begin{definition}\label{def:3.7}
{\rm  A {\bf weakly PS Banach Fredholm bundle} (of index $r$) is a
triple $(X, E, S)$ consisting of a weakly PS Banach bundle $p:E\to
X$ and a weakly Fredholm PS section $S:X\to E$ (of index $r$).
Clearly, each stratum of it, $(X, E, S)_\alpha=(X_\alpha,
E|_{X_\alpha}, S|_{X_\alpha})$ is a Banach Fredholm bundle. If ${\rm
Index}(S_\alpha)=n$ we call $(X, E, S)_\alpha$ {\bf the stratum of
codimension $r-n$} of $(X, E, S)$.  A {\bf (weakly)\footnote{In this
paper we often write ``(A) B''. It means two cases ``A" and ``A
B''.} PS Banach $k$-Fredholm bundle} (of index $r$) is a triple $(X,
E, S)$ consisting of a PS Banach bundle $p:E\to X$ and a (weakly)
$k$-Fredholm PS section $S:X\to E$ (of index $r$). With the same way
we define their strata of codimension $r-n$. These three kinds of
bundles are said to be {\bf oriented} if their strata of codimension
zero, still denoted by $(X, E, S)_0$ without occurrence of
confusions,  are oriented  (as Banach Fredholm bundles). }
\end{definition}

Clearly, the stratum of codimension zero of a (weakly) PS Banach
$k$-Fredholm bundle $(X, E, S)$ is exactly $(X_0, E|_{X_0},
S|_{X_0})$ for any $k\ge 0$. Moreover, every (weakly) PS Banach
$k$-Fredholm bundle $(X, E, S)$ with $k\ge 1$ has no strata of
codimension one.

If the above word ``Banach'' are all replaced by ``Hilbert'' one
gets the {\bf partially smooth Hilbert manifolds and Hilbert
 bundles}.

As said at the beginning of this section our purpose is to
generalize Theorem~\ref{th:1.5} to the PS Banach Fredholm bundle.
However carefully checking the proof of Theorem~\ref{th:1.5} it is
easy to see that in the present case one cannot necessarily extend
an element $v\in E_x$ into a local PS section of the PS bundle $E\to
X$ near $x$. This motivates us to introduce:

\begin{definition}\label{def:3.8}
 {\rm A class of PS sections ${\cal A}$ of $E$ is said to be {\bf rich near $x\in X$}
if for any $v\in E_x$, $\epsilon>0$ and  given open neighborhood
of $x$, $U\subset X$, there exists a PS section $S_v:X\to E$ of
class ${\cal A}$ such that $\|S_{v}(x)-v\|<\epsilon$ and ${\rm
supp}(S_v)\subset U$.\footnote{The second condition shall be used
in the proofs of Theorem~\ref{th:3.14} and
Proposition~\ref{prop:3.16} and below (\ref{e:4.4}).} It is called
{\bf rich near} a subset $K\subset X$ if it is rich near each
point of $K$. }
\end{definition}

Clearly, if $v\ne 0$ we may choose $\epsilon<\frac{1}{2}\|v\|$ so
that $\|S_v(x)\|>\frac{1}{2}\|v\|$. It follows that $S_v(y)\ne 0$
for any $y$ sufficiently close to $x$. More generally, we have:

\begin{lemma}\label{lem:3.9}
 Let ${\cal A}$ be a class of PS sections of $E$, which is rich near $x$.
Then for linearly independent elements $v_1,\cdots, v_n$ in $E_x$
 and $0<\epsilon<<\frac{1}{2}\min_{1\le i\le n}|v_i|$
 there exist an open neighborhood $U$ of $x$ and PS sections
$S_1,\cdots, S_n$ of $E$ of class ${\cal A}$ such that
$\|S_i(x)-v_i\|<\epsilon$, $i=1,\cdots, n$, and that the elements
$S_1(y),\cdots, S_n(y)$ are linearly independent in $E_y$ for each
$y\in U$.
\end{lemma}

\noindent{\bf Proof.}\quad By the definition there exist  PS
sections $S_1,\cdots, S_n$ of $E$ of class ${\cal A}$ such that
$\|S_i(x)-v_i\|<\epsilon$ for $i=1,\cdots, n$. The remains is to
prove that $S_1(y),\cdots, S_n(y)$ are linearly independent for
each $y$ in a small neighborhood $U$ of $x$. Otherwise there exist
a sequence $\{y_m\}\subset X$ and ${\bf
t}^{(m)}=(t_1^{(m)},\cdots, t_n^{(m)})\in S^{n-1}$ such that
$\sum^n_{k=1}t_k^{(m)}S_k(y_m)=0$ and $y_m\to x$ as $m\to\infty$.
Since $S^{n-1}$ is compact we may assume that ${\bf
t}^{(m)}\to{\bf t}=(t_1,\cdots,t_n)\in S^{n-1}$ as $m\to\infty$.
So $\sum^n_{k=1}t_kS_k(x)=0$. This contradiction implies the
desired result. \hfill$\Box$\vspace{2mm}

It should be remarked that  we do not assume that our PS Banach
manifolds must admit PS cut-off functions. (In fact, in actual
examples one needs such a condition to check that a class of PS
sections is rich near a compact subset in the sense of
Definition~\ref{def:3.8}). The assumption of richness will take the
place of the assumptions about cutoff functions used in the previous
sections. It is easily checked that for a PS Banach bundle having
rich PS sections of class ${\cal A}$ near $x\in X$ and  a PS map
$\varphi:Y\to X$ from a PS Banach manifold $Y$ to $X$, if $Y$ admits
PS cut-off functions with arbitrarily small support near any point
then the pull-back PS Banach bundle $\varphi^\ast E\to Y$ has also
rich PS sections of class ${\cal A}$ near any point $y\in
f^{-1}(x)$.

Let $(X, E, S)$ be any one of the three kinds of bundles in
Definition~\ref{def:3.7}. It  is said to have {\bf rich
 PS sections of class ${\cal A}$} if the (weakly) PS Banach bundle $E\to
X$ has such PS sections near the zero locus $Z(S)$ of $X$. For such
a bundle $(X, E, S)$ and each $x\in Z(S)$ we have finitely many
elements $v_1,\cdots, v_n$ in $E_x$ such that ${\rm Im}(DS(x))+ {\rm
span}(\{v_1,\cdots, v_n\})=E_x$. It follows from
Lemmas~\ref{lem:1.2} and \ref{lem:3.9} that there exist PS sections
$\sigma_1,\cdots,\sigma_k$ of $E$ of class ${\cal A}$ such that the
linear operator
$$
T_xX\times\R^k\to E_x,\;(\xi,{\bf t})\mapsto
DS(x)+\sum^k_{i=1}t_i\sigma_i(x)
$$
is surjective. Furthermore, if $X_\alpha$ is the stratum
containing $x$ we can find a small neighborhood $U$ of $x$ in $X$
such that the linear operator
$$
 L(S,y):T_yX\times\R^k\to
E_y,\;(\xi,{\bf t})\mapsto DS(y)+\sum^k_{i=1}t_i\sigma_i(y)
$$
is surjective for any $y\in X_\alpha\cap U\cap Z(S)$. However, if
$Z(S)$ cannot be contained a single stratum then it is impossible
using the present assumptions  to derive that $L(S,y)$ is
surjective for $y\in U\cap Z(S)$ no matter how  $U$ is small. The
reasons are that different strata are only related continuously
and thus one cannot compare differentials of different Fredholm
sections $S|_{X_\alpha}$. Furthermore, even if $L(S,y)$ is
surjective for any $y\in U\cap Z(S)$ we can only deduce that the
restriction of the section
$$
 X\times\R^k\to\Pi_1^\ast E,\;(y,{\bf
t})\mapsto S(y)+\sum^k_{i=1}t_i\sigma_i(y)
 $$
 to $U\times\R^k$ is
transverse to the zero section along $(U\cap Z(S))\times\{0\}\subset
U\times\R^k$, and cannot prove that the restriction of the section
to a small neighborhood of $(x,0)$ in $X\times\R^k$ is transverse to
the zero section yet. Here $\Pi_1$ is the projection from
$X\times\R^k$ to $X$. This means that the local transversality
cannot be obtained from the present assumptions simply. One needs to
apply more assumptions on the PS Banach Fredholm bundle $(X, E, S)$.

\begin{definition}\label{def:3.10}
{\rm Let $(X, E, S)$ be a weakly PS Fredholm bundle of index $r$
(resp. a (weakly) PS $k$-Fredholm bundle of index $r$) as in
Definition~\ref{def:3.7}, and have rich PS sections of class
${\cal A}$. $S$ is called {\bf quasi-transversal  (resp. (weakly)
$k$-quasi-transversal) to the zero section relative to class
${\cal A}$ near} $x\in Z(S)$ if for any PS sections
$\sigma_1,\cdots, \sigma_m$ of $E$ of class ${\cal A}$ the
equality
$$
 DS(x)(T_xX)+ {\rm
span}(\{\sigma_1(x),\cdots, \sigma_m(x)\})=E_x
$$
(including the case $DS(x)(T_xX)=E_x$) implies that there exist a
small open neighborhood $U(x)$ of $x$ in $X$ and $\eta>0$ such that
the section
\begin{eqnarray}\label{e:3.1}
&&\Phi:U(x)\times B_\eta(\R^m)\to \Pi_1^\ast(E|_{U(x)}),\nonumber\\
&&\qquad (x; t_1,\cdots, t_m)\mapsto
S(x)+\sum^m_{i=1}t_i\sigma_i(x),
\end{eqnarray}
is  transversal to the zero section, and is also a weakly Fredholm
PS section of index $m+ r$ (resp. a (weakly) $k$-Fredholm PS
section of index $m+r$) in the sense of  Definition~\ref{def:3.5}.
Here $\Pi_1: U(x)\times B_\eta(\R^m)\to U(x)$ is the natural
projection. We say $S$ to be {\bf quasi-transversal (resp.
(weakly) $k$-quasi-transversal)
 to the zero section relative to class ${\cal A}$} if
it is  quasi-transversal (resp. (weakly) $k$-quasi-transversal) to
the zero section relative to class ${\cal A}$ near each $x\in Z(S)$.
Furthermore, $S$ is called {\bf smoothly quasi-transversal (resp.
smoothly (weakly) $k$-quasi-transversal) to the zero section
relative to class ${\cal A}$ near} $x\in Z(S)$ if the zero locus of
$\Phi$ in (\ref{e:3.1}), $\Phi^{-1}(0)$ is a true smooth manifold of
dimension $m+r$.
 Naturally,
$S$ is said to be {\bf smoothly quasi-transversal (resp. smoothly
(weakly) $k$-quasi-transversal) to the zero section relative to
class ${\cal A}$} if it is smoothly quasi-transversal (resp.
smoothly (weakly) $k$-quasi-transversal) to the zero section
relative to class ${\cal A}$. (In this case an orientation on the
top stratum of $\Phi^{-1}(0)$ may naturally induce ones on the lower
strata.)}
\end{definition}

\begin{remark}\label{rem:3.11}
{\rm (i) Clearly, when $S$ is quasi transversal to the zero section
  relative to class ${\cal A}$ near $x\in Z(S)$, if $DS(x)$ is also onto
  then the restriction
of $S$ to a small neighborhood $U(x)$ of $x$ is transversal to the
zero section.\\
(ii) Note that $U(x)\cap Z(S)\times\{0\}\subset Z(\Phi)$ and for any
$y\in U(x)\cap Z(S)$ and $(\xi,{\bf u})\in T_yX\times\R^m$,
$$
D\Phi(y,0)(\xi,{\bf u})=DS(y)(\xi)+ \sum^m_{i=1}u_i\sigma_i(y).
$$
That the section $\Phi$ in (\ref{e:3.1}) is transversal to the zero
section implies that
$$
DS(y)(T_yX)+ {\rm span}(\{\sigma_1(y),\cdots, \sigma_m(y)\})=E_y
\;\forall y\in U(x)\cap Z(S).
$$
(iii) When the section $\Phi$ in (\ref{e:3.1}) is  a weakly Fredholm
PS section of index $m+ r$ (resp. a (weakly) $k$-Fredholm PS section
of index $m+r$), for any ${\bf t}\in B_\eta(\R^m)$, it follows from
Lemma~\ref{lem:1.3}(ii) that the section $\Phi_{\bf t}:U(x)\to E,
y\mapsto \Phi(y, {\bf t})$ is also a weakly Fredholm PS section of
index $r$ (resp. a (weakly) $k$-Fredholm PS section of index $r$).
In fact, let $\Phi^\alpha$ be the restriction of $\Phi$ to the
stratum $U(x)_\alpha\times B_\eta(\R^m)$ of $U(x)\times
B_\eta(\R^m)$. It is given by
$$
U(x)_\alpha\times B_\eta(\R^m)\to \Pi_1^\ast(E|_{U(x)_\alpha}),\,
(y; t_1,\cdots, t_m)\mapsto
S^\alpha(y)+\sum^m_{i=1}t_i\sigma^\alpha_i(y),
$$
where $S^\alpha$ and $\sigma_i^\alpha$ are the restrictions of $S$
and $\sigma_i$ to $U(x)_\alpha$, $i=1,\cdots,m$. Assume that
$E|_{U(x)_\alpha}$ is trivializable over an open subset $V\subset
U(x)$ and that $S^{\alpha t},\sigma_i^{\alpha t}: V\mapsto E_z$ for
some $z\in V$ are trivialization representatives of $S^\alpha$ and
$\sigma^\alpha_i$ over $V$, $i=1,\cdots,m$. By the definition of the
strongly Fredholm section the map
$$
\Phi^{\alpha t}:V\times B_\eta(\R^m)\to E_z,\, (y; t_1,\cdots,
t_m)\mapsto S^{\alpha t}(y)+\sum^m_{i=1}t_i\sigma^{\alpha t}_i(y),
$$
is  Fredholm. So for any $(y,{\bf t})\in V\times B_\eta(\R^m)$ the
tangent map
\begin{eqnarray*}
&&d\Phi^{\alpha t}(y,{\bf t}):T_yV\times \R^m\to E_z, \\
&&(\xi;u_1,\cdots, u_m)\mapsto dS^{\alpha
t}(y)(\xi)+\sum^m_{i=1}t_id\sigma^{\alpha t}_i(y)(\xi)+
\sum^m_{i=1}u_i\sigma^{\alpha t}_i(y),
 \end{eqnarray*}
 is Fredholm. By Lemma~\ref{lem:1.3}(ii) the linear map
$$
d\Phi^{\alpha t}_{\bf t}(y):T_yV\to E_z,\; \xi\mapsto dS^{\alpha
t}(y)(\xi)+\sum^m_{i=1}t_id\sigma^{\alpha t}_i(y)(\xi),
$$
 is Fredholm, where $\Phi^{\alpha t}_{\bf t}:V\to E_z, y\mapsto\Phi^{\alpha t}(y,{\bf
 t})$. It follows that $\Phi^{\alpha t}_{\bf t}$ and thus
 $\Phi^{\alpha}_{\bf t}$, $\Phi_{\bf t}$ are all Fredholm. The conclusions on index are easily derived from
 Lemma~\ref{lem:1.3}(ii).\\
 (iv)  In applications to
 Floer homology one needs stronger conclusions because the strata of
 codimension one must exist. This is the reason that we introduce the smooth quasi-transversality.}
\end{remark}

Unlike in Section~\ref{sec:1} the corresponding section with one in
(\ref{e:3.1}) is Fredholm and thus is locally proper,  a weakly
Fredholm section might not be locally proper. We need to introduce
the following notion.

\begin{definition}\label{def:3.12}
{\rm Let $(X, E, S)$ be any one of the three kinds of bundles in
Definition~\ref{def:3.7}, and have rich PS sections of class ${\cal
A}$. $S$ is called {\bf locally proper relative to class ${\cal A}$
near $Z(S)$} if for any $x\in Z(S)$ and PS sections
$\sigma_1,\cdots,\sigma_m$ of class ${\cal A}$ there exists a
neighborhood $U(x)$ and $\eta>0$ such that  the closure of the set
$$
\bigl\{(x,{\bf t})\in U(x)\times B_\varepsilon(\R^m)\,|\,
(S+\sum^m_{i=1}t_i\sigma_i)(x)=0\bigr\}
$$
is compact in $X\times\R^m$, where ${\bf t}=(t_1,\cdots,
t_m)\in\R^m$.}
\end{definition}

\begin{definition}\label{def:3.13}
{\rm Let $(X, E, S)$ be a weakly PS Fredholm bundle of index $r$
(resp. a (weakly) PS $k$-Fredholm bundle of index $r$) as in
Definition~\ref{def:3.7}, and have rich PS sections of class ${\cal
A}$.  It  is called {\bf good} (resp. ({\bf weakly}) $k$-{\bf good})
{\bf relative to class} ${\cal A}$ if the following conditions are
satisfied:
\begin{description}
\item[(a)] $S$ is quasi (resp. (weakly) $k$-quasi) transversal to the zero section relative to class
${\cal A}$;
\item[(b)] $S$ is locally proper relative to class ${\cal A}$ near
      $Z(S)$.
\end{description}
If $S$ is   smoothly quasi-transversal (resp. smoothly (weakly)
$k$-quasi-transversal) to the zero section relative to class
${\cal A}$, the above bundle $(X, E, S)$ is called {\bf smoothly
good} (resp. {\bf smoothly} ({\bf weakly}) $k$-{\bf good}) {\bf
relative to class} ${\cal A}$. }
\end{definition}

 We now are in a position to give and prove
the analogue of Theorem~\ref{th:1.5}.

\begin{theorem}\label{th:3.14} Let $X$ be a separable PS Banach manifold $X$
with corner, and $(X, E, S)$ be a (weakly) PS Banach $k$-Fredholm
bundle of index $r$.  Then one has:
\begin{description}
\item[(i)] $Z(S)$ is a weakly $k$-PS submanifold of $X$ of dimension $r$ and with
corners provided that  $S$ is transversal to the zero section.
Moreover, an orientation of $(X, E, S)$ determines one of $Z(S)$.
(If $(X, E, S)$ is only a weakly PS Banach Fredholm bundle of index
$r$ then $Z(S)$ is a weakly PS submanifold of dimension $r$.)

 \item[(ii)] In general, if $(X, E,  S)$ has {\rm rich PS sections of
 class} ${\cal A}$ and  {\rm compact zero locus} $Z(S)$, and is also
 {\rm (weakly) $k$-good} then one  can find an open neighborhood ${\cal W}$ of $Z(S)$
in $X$ satisfying:
\begin{description}
 \item[(ii.1)] There exists a sufficiently small
 perturbation PS section $S^\prime$ of $E$, which is a  sum of finite PS
  sections of class ${\cal A}$,  such that the restriction
  $$(S+ S^\prime)|_{\cal W}: {\cal W}\to E|_{\cal W}$$
 is a (weakly) $k$-Fredholm PS section which is
 transversal to the zero section and has also a compact zero locus $Z((S+ S^\prime)|_{\cal W})$.
   In particular, if $S$ is transversal to the zero section one may
take $S^\prime=0$.

\item[(ii.2)]  $S^\prime$ can be chosen so small that
$Z((S+S^\prime)|_{\cal W})$ is contained in a given small
neighborhood of $Z(S)$. Specially, one can require all
$k$-boundaries $\partial^k(Z((S+S^\prime)|_{\cal W}))$ to be equal
to $Z((S+S^\prime)|_{\cal W})\cap\partial^k X$ for $k=1,2,\cdots$.

\item[(ii.3)] If $(X, E, S)$ is oriented  then
$Z((S+S^\prime)|_{\cal W})$ is an
 oriented weakly $k$-PS manifold of dimension $r$.

\item[(ii.4)] Any two such sufficiently small perturbations
$S^\prime$ and $S^{\prime\prime}$  yield
 cobordant $r$-dimensional weakly $k$-PS manifolds $Z((S+ S^\prime)|_{\cal W})$
 and $Z((S+S^{\prime\prime})|_{\cal W})$.
\end{description}
Moreover, if $(X, E, S)$ is  smoothly good (resp.  smoothly (
weakly) $k$-good) relative to class ${\cal A}$ then all words
``weakly $k$-PS manifold(s)'' in (ii.1)-(ii.4) above will be
replaced by ``manifold(s)''.
 \end{description}
 \end{theorem}

Similar to the convention below Theorem~\ref{th:1.5}, in this and
next sections we make the following {\bf convention}: when using
Sard theorem or Sard-Smale theorem to get an empty set or a PS
manifold of dimension lower than expected one {\it we shall not
mention that case} ! \vspace{2mm}

\noindent{\bf Proof of Theorem~\ref{th:3.14}.}\quad  Noting the
convention just, by the definition of the weakly PS Banach manifold
in Definition~\ref{def:3.1}, the definitions above
Claim~\ref{cl:3.3} and Definition~\ref{def:3.7} we may apply
Theorem~\ref{th:1.5} to the restriction of $S$ to each stratum of
$X$ to arrive at the desired (i).

 To prove (ii), note that $S$ is stratawise Fredholm.
 For each $x\in Z(S)$ there exist $v_1,\cdots, v_p\in E_x$ such that
 $DS(x)(T_xX)+ {\rm Span}(\{v_1,\cdots, v_p\})=E_x$. Since $(X, E, S)$ has
 rich PS sections of class ${\cal A}$ it follows from Lemma~\ref{lem:1.2} and
 Lemma~\ref{lem:3.9} that there exist  PS sections $\sigma_1,\cdots, \sigma_p$
of class ${\cal A}$ such that
$$
DS(x)(T_xX)+ {\rm Span}(\{\sigma_1(x),\cdots, \sigma_p(x)\})=E_x.
$$
Since $S$ is (weakly) $k$-quasi transversal to the zero section
relative to class ${\cal A}$ near any $x\in Z(S)$,
 by Remark~\ref{rem:3.11}(ii) and the compactness of $Z(S)$  one can find
finitely many PS sections of class ${\cal A}$, say $S_1,\cdots, S_m$
\footnote{These sections can be required to have small supports by
Definition~\ref{def:3.8}.}, such that
\begin{equation}\label{e:3.2}
DS(z)(T_zX)+ {\rm span}(\{S_1(z),\cdots, S_m(z)\})=E_z
\end{equation}
  for any $z\in Z(S)$. Using the facts that
  $S$ is (weakly) $k$-quasi transversal to the zero section relative to class
${\cal A}$ and also locally proper relative to class ${\cal A}$
again,  it follows from the compactness of $Z(S)$ that there exist
an open neighborhood ${\cal W}$ of $Z(S)$ and $\varepsilon>0$ such
that:
\begin{description}
\item[$(\ref{e:3.2})_1$] The PS section
$$
 \Phi:{\cal
W}\times B_\varepsilon(\R^m)\to\Pi_1^\ast(E|_{\cal W}),\;(y,{\bf
t})\mapsto S(y)+\sum^m_{i=1}t_iS_i(y)$$
  is a (weakly) $k$-Fredholm PS section of index $m+ r$ which is transversal to the
zero  section. Here $\Pi_1$ is
  the projection from ${\cal
W}\times B_\varepsilon(\R^m)$ to ${\cal W}$. So $\Phi^{-1}(0)$ is a
weakly $k$-PS submanifold of dimension $m+r$( if $Z(S)\cap
X_0\ne\emptyset$).

\item[$(\ref{e:3.2})_2$] The closure of $\Phi^{-1}(0)$ in $X\times\R^m$ is compact.
 It follows that for any small neighborhood
 ${\cal V}$ of $Z(S)$ in $X$ there exists an $\epsilon\in
 (0,\varepsilon)$ such that for each ${\bf t}\in
 B_\epsilon(\R^m)$,
 \begin{equation}\label{e:3.3}
 Cl({\cal W})\cap\bigl(S+\sum^m_{i=1 }t_iS_i\bigr)^{-1}(0)\subset {\cal V}.
\end{equation}
\end{description}
  Let $P_2$ be the restriction of the natural projection
  ${\cal W}\times B_\varepsilon(\R^m)\to B_\varepsilon(\R^m)$ to $\Phi^{-1}(0)$.
  For each stratum ${\cal W}_\alpha$ with
$\Phi^{-1}(0)\cap({\cal W}_\alpha\times
B_\varepsilon(\R^m))\ne\emptyset$, as in the proof of
Theorem~\ref{th:1.5} we deduce that the restriction of $P_2$ to
$\Phi^{-1}(0)\cap({\cal W}_\alpha\times B_\varepsilon(\R^m))$ is
Fredholm and has the same index as $S|_{{\cal W}_\alpha}$. Moreover,
${\bf t}\in B_\varepsilon(\R^m)$ is a regular value of this
restriction if and only if the section $\Phi_{\bf t}|_{{\cal
W}_\alpha}: {\cal W}_\alpha\to E|_{{\cal W}_\alpha}$ is transversal
to the zero section. (cf. Proposition~\ref{prop:A.14}). Here
\begin{equation}\label{e:3.4}
 \Phi_{\bf t}:{\cal W}\to E|_{{\cal
W}},\;y\mapsto\Phi(y, {\bf t}).
\end{equation}
 Applying the
Sard-Smale theorem to the restriction of $P_2$ to each stratum of
$\Phi^{-1}(0)$ we get a residual subset
$B_\varepsilon(\R^m)_{reg}\subset B_\varepsilon(\R^m)$ such that for
each ${\bf t}\in B_\varepsilon(\R^m)_{reg}$ the section  $\Phi_{\bf
t}$ in (\ref{e:3.4})  is transversal to the zero section, and is
also a (weakly) $k$-Fredholm PS section of index $r$ by
Remark~\ref{rem:3.11}(iii). Take an open neighborhood
 ${\cal V}$ of $Z(S)$ such that $Cl({\cal V})\subset{\cal W}$.
By $(\ref{e:3.2})_2$ we may find an $\epsilon\in
 (0,\varepsilon)$ such that (\ref{e:3.3}) holds for each ${\bf t}\in
 B_\epsilon(\R^m)$.
It follows that $\Phi_{\bf t}^{-1}(0)\subset{\cal V}$ for any ${\bf
t}\in B_\epsilon(\R^m)$. So each ${\bf t}\in
B_\epsilon(\R^m)_{reg}:=B_\epsilon(\R^m)\cap
B_\varepsilon(\R^m)_{reg}$ gives a desired perturbation section
$S'=\sum^m_{i=1}t_iS_i$  for which $Z((S+S')|_{\cal W})=\Phi_{\bf
t}^{-1}(0)$ is a compact weakly $k$-PS manifold of dimension $r$.
Then the conclusions in (ii.1) and (ii.2) may follow directly.

The proof of (ii.3) is almost repeat of that of
Theorem\ref{th:1.5}(B.2). Indeed, the map
$$
{\mathcal W}\times B_\varepsilon(\R^m)\times [0,1]\to {\cal W}\times
B_\varepsilon(\R^m),\;(x, {\bf t}, t)\mapsto (x, t{\bf t})
$$
is stratawise Fredholm and the restriction to each stratum has index
$1$. So the composition of it with $\Phi$ gives a (weakly)
$k$-Fredholm PS section
$$
 {\cal
W}\times B_\varepsilon(\R^m)\times [0,1]\to\Pi_1^\ast(E|_{\cal
W}),\;(y,{\bf t}, t)\mapsto S(y)+ t\sum^m_{i=1}t_iS_i(y),
$$
where $\Pi_1$ is also the projection to the first factor. As in
Remark~\ref{rem:3.11}(iii) we get a family of homotopic (weakly)
$k$-Fredholm PS sections
$$
 \Psi^t: {\cal
W}\times B_\varepsilon(\R^m)\to\Pi_1^\ast(E|_{\cal W}),\;(y,{\bf t},
t)\mapsto S(y)+ t\sum^m_{i=1}t_iS_i(y)
$$
for any $t\in [0,1]$. As in the proof of (ii.3) of
Theorem\ref{th:1.5}(B.2), we can derive
   from a given continuous nowhere zero  section of the restriction of
    ${\rm det}(DS)$ to  $Z(S)\cap X_0$
   that the weakly $k$-PS manifold
   of dimension $m+r$,
    $\Phi^{-1}(0)=\{(y,{\bf t})\in
    {\cal W}\times B_\varepsilon(\R^m)\,|\,\Phi(y,{\bf t})=0\}$ carries an
   orientation on its stratum of dimension $r+ m$. This orientation naturally
   induces one on the top stratum of each weakly $k$-PS submanifold $\Phi_{\bf
   t}^{-1}(0)$ of $X$.

Obviously, if $(X, E, S)$ is  smoothly good (resp.  smoothly (
weakly) $k$-good) relative to class ${\cal A}$, $\Phi_{\bf
   t}^{-1}(0)$ is a manifold of dimension $r$ because $\Phi_{\bf
   t}^{-1}(0)=P_2^{-1}({\bf t})$ and that ${\bf t}$ is a common regular value of the restriction of $P_2$ to
each stratum $\Phi^{-1}(0)\cap({\cal W}_\alpha\times
B_\varepsilon(\R^m))$ implies that ${\bf t}$ is a regular value of
$P_2$ as a smooth map from the smooth manifolds $\Phi^{-1}(0)$ to
$B_\varepsilon(\R^m)$. Note also that each  stratum
$\Phi^{-1}(0)\cap({\cal W}_\alpha\times B_\varepsilon(\R^m))$ is a
submanifold of $\Phi^{-1}(0)$.

 Finally we prove (ii.4) in two steps.
The first step is completely similar to \S1. The method in the
second step can also be applied to \S1. For completeness we will
give the detailed proof.
\\
{\bf Step 1}. For ${\bf t}=(t_1,\cdots, t_m)\in
B_\varepsilon(\R^m)_{reg}$ and ${\bf r}=(r_1,\cdots, r_m)\in
B_\varepsilon(\R^m)_{reg}$ let ${\cal P}^l({\bf t}, {\bf r})$
denote the space of all $C^l$-smooth paths $\gamma:[0,1]\to
B_\varepsilon(\R^m)$ with $\gamma(0)={\bf t}$ and $\gamma(1)={\bf
r}$. Consider the section
\begin{eqnarray*}
&&\Psi: {\cal W}\times [0, 1]\times {\cal P}^l({\bf t}, {\bf
r})\to
\Pi_1^\ast (E|_{\cal W}),\\
&&\qquad (x, \tau,\gamma)\mapsto \Phi(x,\gamma(\tau))=
S(x)+\sum^m_{i=1}\gamma_i(\tau)S_i(x),
\end{eqnarray*}
where $\Pi_1: {\cal W}\times [0, 1]\times {\cal P}^l({\bf t}, {\bf
r})\to X$ is the natural projection and $\gamma=(\gamma_1,\cdots,
\gamma_m)\in {\cal P}^l({\bf t}, {\bf r})$. (Note that
$\Psi^{-1}(0)\cap({\cal W}_0\times [0, 1]\times {\cal P}^l({\bf t},
{\bf r}))\ne\emptyset$ for ${\cal W}_0:={\cal W}\cap X_0$ if
$Z(S)\cap X_0\ne\emptyset$.) For any zero $(x,\tau,\gamma)$ of
$\Psi$ one has the (vertical) differential
\begin{eqnarray*}
&&D\Psi(x,\tau,\gamma): T_xX\times\R\times T_\gamma{\cal P}^l({\bf
t}, {\bf r})\to E_x,\\
&&\qquad (\xi, e, \alpha)\to D\Phi(x, \gamma(\tau))(\xi,
\alpha(\tau)+ e\dot\gamma(\tau)).
\end{eqnarray*}
For $\tau=0$ and $\tau=1$ it is surjective as a function of $\xi$
alone since the sections $\Psi(\cdot, 0,\gamma)=\Phi(\cdot, {\bf
t})$ and $\Psi(\cdot, \tau,\gamma)=\Phi(\cdot, {\bf r})$ are
transversal to the zero section. For $\tau\in [0,
1]\setminus\{0,1\}$ we can choose $\alpha(\tau)$ arbitrarily and
thus derive that this operator is also surjective as a function of
$\xi$ and $\alpha$. Hence $\Psi$ is transversal to the zero section.
Let $P$ be the projection from $\Psi^{-1}(0)$ to the third factor.
Applying the arguments in Step 3 of the proof of
Theorem~\ref{th:1.5} to each stratum $\Psi^{-1}(0)\cap{\cal
W}_\alpha$, we get that $\dim {\rm Ker}(DP(x,\tau,\gamma))={\rm
Ind}(S|_{{\cal W}_\alpha})+ 1$ at any $(x,\tau,\gamma)\in
\Psi^{-1}(0)\cap{\cal W}_\alpha$. Hence the restriction of $P$ to
$\Psi^{-1}(0)\cap{\cal W}_\alpha$ is a Fredholm operator of index
${\rm Ind}(S|_{{\cal W}_\alpha})+ 1$. Moreover, $\gamma\in {\cal
P}^l({\bf t}^{(1)}, {\bf t}^{(2)})$ is a regular value of the
restriction of $P$ to $\Psi^{-1}(0)\cap{\cal W}_\alpha$ if and only
if the restriction of  the section
$$
\Psi_\gamma: {\cal W}\times [0, 1]\to \Pi_1^\ast E, (x, \tau)\mapsto
S(x)+\sum^m_{i=1}\gamma_i(\tau)\sigma_i(x),
$$
to ${\cal W}_\alpha\times [0, 1]$ is  transversal to the zero
section, which implies that this restriction is a (weakly)
$k$-Fredholm section of index ${\rm Ind}(S|_{{\cal W}_\alpha})+1$.
Since ${\cal W}$ has only finitely many strata we may use the
Sard-Smale theorem to get a residual subset ${\cal P}^l_{reg}({\bf
t}, {\bf r})\subset {\cal P}^l({\bf t}, {\bf r})$
  such that each $\gamma\in {\cal P}^l_{reg}({\bf
t}, {\bf r})$ gives a compact weakly $k$-PS manifold
$\Psi^{-1}_\gamma(0)$ which forms a weakly $k$-PS manifold cobordsim
 between $\Phi_{\bf t}^{-1}(0)$ and $\Phi_{\bf
r}^{-1}(0)$.

\noindent{\bf Step 2}.\quad Assume that $S_1',\cdots, S_{m'}'$ are
another group of PS sections satisfying
\begin{equation}\label{e:3.5}
DS(z)(T_zX)+ {\rm span}(\{S'_1(z),\cdots, S'_{m'}(z)\})=E_z
\end{equation}
  for any $z\in Z(S)$. Then we have an open neighborhood ${\cal
W}'$ of $Z(S)$ and $\varepsilon'>0$ such that:
\begin{description}
\item[$(\ref{e:3.5})_1$] The  PS section
$$
 \Phi': {\cal W}'\times
B_{\varepsilon'}(\R^{m'})\to\Pi_1^\ast(E|_{{\cal W}'}),\;(y,{\bf
t})\mapsto S(y)+\sum^{m'}_{i=1}t'_iS'_i(y)$$
  is  a
(weakly) $k$-Fredholm PS section of index $m'+ r$ that is
transversal to the zero  section. Here $\Pi_1$ is
  the projection from ${\cal
W}'\times B_\varepsilon(\R^{m'})$ to ${\cal W}$. So
$(\Phi')^{-1}(0)$ is a weakly $k$-PS submanifold of dimension $m'+r$
(if $Z(S)\cap X_0\ne\emptyset$).

 \item[$(\ref{e:3.5})_2$]  The closure of the zero locus $(\Phi')^{-1}(0)$ in
 $X\times\R^{m'}$ is compact.
 \end{description}
Let  $B_{\varepsilon'}(\R^{m'})_{reg}\subset B_\varepsilon(\R^{m'})$
be the corresponding residual subset such that for each ${\bf
t}'=(t'_1,\cdots, t'_{m'})\in B_{\varepsilon'}(\R^{m'})_{reg}$ the
section
$$
\Phi'_{{\bf t}'}: {\cal W}'\to E|_{{\cal W}'},\;y\to
\Phi'(\cdot,{\bf t}')
$$
is a (weakly) $k$-Fredholm PS section of index $r$ that is
transversal to the zero section.

By (\ref{e:3.2}) and (\ref{e:3.5}), for any $z\in Z(S)$ one has
$$
 DS(z)(T_zX)+ {\rm
span}(\{S_1(z),\cdots, S_m(z), S'_1(z),\cdots, S'_{m'}(z)\})=E_z.
$$
 As before we can take an open neighborhood ${\cal
W}^\ast\subset{\cal W}\cap{\cal W}'$ of $Z(S)$ and
$0<\varepsilon^\ast<\min\{\varepsilon,\varepsilon'\}$ such that:
\begin{description}
\item[$(\ref{e:3.5})_3$] The  PS section
\begin{eqnarray*}
 {\cal G}: {\cal W}^\ast\times
B_{\varepsilon^\ast}(\R^{m})\times
B_{\varepsilon^\ast}(\R^{m'})\to\Pi_1^\ast (E|_{\cal W}),\\
(y,{\bf t}, {\bf t}')\mapsto S(y)+ \sum^{m}_{i=1}t_iS_i(y)+
\sum^{m'}_{i=1}t'_iS'_i(y) \end{eqnarray*}
  is a (weakly) $k$-Fredholm PS section of index $m+ m'+ r$
  that is  transversal to the zero  section. Here $\Pi_1$
is the projection from ${\cal W}^\ast\times
B_{\varepsilon^\ast}(\R^{m})\times B_\varepsilon(\R^{m'})$ to ${\cal
W}^\ast$. So ${\cal G}^{-1}(0)$ is a weakly $k$-PS submanifold of
dimension $m+ m'+r$( if $Z(S)\cap X_0\ne\emptyset$).

\item[$(\ref{e:3.5})_4$] The closure of the zero locus ${\cal
G}^{-1}(0)$ in $X\times R^m\times\R^{m'}$ is compact. (One has also
${\cal G}^{-1}(0)\cap(X_0\times R^m\times\R^{m'})\ne\emptyset$ if
$Z(S)\cap X_0\ne\emptyset$.)
 \end{description}
Let $(B_{\varepsilon^\ast}(\R^{m})\times
B_{\varepsilon^\ast}(\R^{m'}))_{reg}\subset
B_{\varepsilon^\ast}(\R^{m})\times B_{\varepsilon^\ast}(\R^{m'})$ be
the corresponding residual subset such that for each $({\bf t}, {\bf
t}')\in (B_{\varepsilon^\ast}(\R^{m})\times
B_{\varepsilon^\ast}(\R^{m'}))_{reg}$ the  section
$$
 {\cal G}_{({\bf t}, {\bf t}')}:
{\cal W}^\ast\to E|_{{\cal W}^\ast},\;y\to {\cal G}(\cdot, {\bf t},
{\bf t}')
$$
is a (weakly) $k$-Fredholm PS section of index $r$ that is
transversal to the zero section.

Take an open neighborhood ${\cal V}$ of $Z(S)$ such that $Cl({\cal
V})\subset{\cal W}^\ast$. By $(\ref{e:3.5})_4$ we have
$\epsilon\in (0,\varepsilon^\ast]$ such that  the set
\begin{equation}\label{e:3.6}
\cup_{|{\bf t}|\le\epsilon}
 \cup_{|{\bf t}'|\le\epsilon}Cl({\cal
 W}^\ast)\cap ({\cal G}_{({\bf t}, {\bf t}')})^{-1}(0)\subset{\cal V},
\end{equation}
where ${\bf t}=(t_1,\cdots,
 t_{m})\in\R^m$ and ${\bf t}'=(t'_1,\cdots,
 t'_{m'})\in\R^{m'}$.

 By Lemma~\ref{lem:1.7} we can take a small
 $$(\hat {\bf t},
\hat {\bf t}')\in (B_{\varepsilon^\ast}(\R^{m})\times
B_{\varepsilon^\ast}(\R^{m'}))_{reg}\cap
B_{\epsilon}(\R^{m})\times B_{\epsilon}(\R^{m'})
$$
such that $\hat {\bf t}\in B_{\varepsilon}(\R^m)_{reg}$ and $\hat
{\bf t}'\in B_{\varepsilon'}(\R^{m'})_{reg}$. Let ${\cal B}^l(\hat
{\bf t}, \hat {\bf t}')$ denote the space of all $C^l$-smooth
paths $\zeta:[0,1]\to B_{\epsilon}(\R^m)\times
B_{\epsilon}(\R^{m'})$ with $\zeta(0)=(\hat{\bf t}, 0)$ and
$\zeta(1)=(0, \hat{\bf t}')$. Consider the section
\begin{eqnarray*}
&&\Theta: {\cal W}^\ast\times [0, 1]\times {\cal B}^l(\hat {\bf
t}, \hat {\bf t}')\to \Pi_1^\ast(E|_{{\cal W}^\ast}),\\
&&\qquad (x, \tau,\zeta)\mapsto
S(x)+\sum^m_{i=1}\zeta_i(\tau)S_i(x)+\sum^{m'}_{j=1}\zeta'_j(\tau)S'_j(x),
\end{eqnarray*}
where $\zeta=(\zeta_1,\cdots, \zeta_m; \zeta'_1,\cdots,
\zeta'_{m'} )$ belongs to ${\cal B}^l(\hat {\bf t}, \hat {\bf
t}')$ and $\Pi_1: {\cal W}^\ast\times [0, 1]\times {\cal B}^l(\hat
{\bf t}, \hat {\bf t}')\to {\cal W}^\ast$ is the projection.
Clearly, $\Theta(\cdot, 0,\zeta)=\Phi(\cdot, \hat {\bf t})$ and
$\Theta(\cdot, 1,\zeta)=\Phi'(\cdot, \hat {\bf t}')$. Since both
sections $\Phi(\cdot, \hat {\bf t})$ and $\Phi'(\cdot, \hat {\bf
t}')$ are transversal to the zero section it follows from
$(\ref{e:3.5})_3$ that at any zero $(x,\tau,\zeta)$ of $\Theta$
the (vertical) differential
$$
D\Theta(x,\tau,\zeta): T_xX\times\R\times T_\zeta{\cal B}^l(\hat
{\bf t}, \hat {\bf t}')\to E_x
$$
is surjective. Hence $\Theta$ is transversal to the zero section.
(Clearly, $\Theta^{-1}(0)\cap({\cal W}_0^\ast\times [0, 1]\times
{\cal B}^l(\hat {\bf t}, \hat {\bf t}'))\ne\emptyset$ if  $Z(S)\cap
X_0\ne\emptyset$ and ${\cal W}^\ast_0:={\cal W}^\ast\cap X_0$). As
in Step 1 we have a residual subset ${\cal B}^l(\hat {\bf t}, \hat
{\bf t}')_{reg}\subset {\cal B}^l(\hat {\bf t}, \hat {\bf t}')$ such
that any path $\zeta\in {\cal B}^l(\hat {\bf t}, \hat {\bf
t}')_{res}$ yields a compact weakly $k$-PS manifold
$(\Theta_\zeta(\cdot,\cdot))^{-1}(0)$ which forms a weakly $k$-PS
cobordism between $(\Theta_\zeta(\cdot, 0))^{-1}(0)= (\Phi_{\hat
{\bf t}})^{-1}(0) $ and $(\Theta_\zeta(\cdot,
1))^{-1}(0)=(\Phi'_{\hat {\bf t}'})^{-1}(0)$ because (\ref{e:3.6})
implies that $\Theta_\zeta^{-1}(0)\subset{\cal V}\times [0,1]$. Here
$\Theta_\zeta:{\cal W}^\ast\times [0,1]\to \Pi_1^\ast(E|_{{\cal
W}^\ast})$ is given by $\Theta_\zeta(x, \tau)=\Theta(x,\tau,\zeta)$.

Combing Step 1 and Step 2 we complete the proof of (ii.4).

 Finally, if $(X, E, S)$ is  smoothly good (resp.  smoothly ( weakly)
$k$-good) relative to class ${\cal A}$ it is easily checked that the
argument techniques above Step 1 may also apply to Step 1 and Step
2. \hfill$\Box$\vspace{2mm}

In general $Z(S+S^\prime)$ is called a {\bf Euler chain} of
$(X,E,S)$. Under the assumptions of Theorem~\ref{th:3.14}, if  $S$
is a (weakly) $1$-Fredholm PS section and the top stratum $X_0$ of
$X$ has no boundary (including corner) as a Banach manifold then it
follows from Theorem~\ref{th:3.14} (ii) that $Z(S+S^\prime)$ is a
compact (weakly) $1$-PS manifold of dimension $r$ and without strata
of codimension one, and, in particular, is a closed weakly PS
manifold of dimension $r$. Thus it determines a well-defined
homology class $[Z(S+S^\prime)]$ in $H_d(X,\Z_2)$, or a class
$[Z(S+S^\prime)]$ in $H_d(X,\Z)$ if $(X, E, S)$ is oriented, which
is independent of choices of generic $S^\prime$. In this case this
class is called the {\bf Euler class} of the triple $(X, E,S)$, and
denoted by $e(E, S)$. Correspondingly, $Z(S+ S^\prime)$ is called a
{\bf Euler cycle} of $(X, E,S)$. Later we always consider the
oriented case without special statements.

Before discussing furthermore properties we give a relation between
finitely dimensional closed weakly PS manifolds and pseudo-cycles.
Recall on the page 90 of [McSa1] that
 a $m$-dimensional {\bf
pseudo-cycle} in a (smooth) manifold $V$ is a smooth map $g:M\to
V$ defined on an oriented  $m$-dimensional manifold $M$(possibly
noncompact) such that  the boundary
$$g(M^\infty)=\bigcap_{K\subset M\,compact}\overline{g(M\setminus K)}$$
of $g(M)$ has dimension at most $m-2$ in the sense that there
exists a smooth manifold $W$ of dimension at most $m-2$ and a
smooth map $h:W\to V$ such that $g(M^\infty)\subset h(W)$. When
$g(M)$ has compact closure in $V$, $g$ was called a {\bf strong
pseudo-cycle} in [Lu1]. Clearly,  these two notions are equivalent
in a compact manifold. Note that  the identity map $V\to V$ is not
a strong pseudo-cycle in the noncompact manifold $V$ according to
the definition. It was proved in \cite{McSa} and \cite{Sch}: (i)
every strong pseudo-cycle in $V$ represents an integral singular
homology class of $V$, (ii) each integral singular homology class
$\alpha\in H_\ast(V,\Z)$ may be represented by {\bf smooth cycles}
(a kind of special strong pseudo-cycles) and all smooth cycles
representing $\alpha$ are cobordant. Actually Thom showed in
\cite{Th} that a suitable integer multiple of any integral
singular homology class $\alpha\in H_\ast(V,\Z)$ may be
represented by a compact oriented submanifold of $X$ without
boundary. Two pseudo-cycles $e:P\to V$ and $g:Q\to V$ are called
{\bf transverse} if either $e(P)\cap g(Q)=\emptyset$ or
$e(P^\infty)\cap\overline{g(Q)}=\emptyset$, $\overline{e(P)}\cap
g(Q^\infty)=\emptyset$ and $T_xV={\rm Im}(de(p))+{\rm Im}(dg(q))$
whenever $e(p)=g(q)=x$. However, for two transverse pseudo-cycles
$e$ and $g$ as above, if $\Delta(e,g):=\{(p,q)\in P\times Q|
e(p)=g(q)\}\ne\emptyset$ then it is  a compact manifold of
dimension ${\rm dim}P+{\rm dim}Q-{\rm dim}V$ provided that only
one of $e$ and $g$ is strong.  Specially, $\Delta(e,g)$ is a
finite set if $P$ and $Q$ are of complementary dimension. By
Remark~\ref{rem:3.2}($2^\circ$) and the definition of the closed
weakly PS Banach manifold below Remark~\ref{rem:3.2} we
immediately get:

\begin{proposition}\label{prop:3.15}
Let $N$ be a closed oriented weakly PS manifold of dimension $n$.
Then for any smooth manifold $Y$, the restriction of any PS map
$f:N\to Y$ to the top stratum of $N$ is a strong pseudo-cycle of
dimension $n$ in $Y$, and therefore determines an integral singular
homology class in $H_n(Y,\Z)$, denoted by $f_\ast([N])$.
\end{proposition}

It follows that a closed oriented (weakly) PS manifold of finite
dimension behaves as a closed oriented manifold from the point of
view of topology.

As in Remark~\ref{rem:1.9} we can start from a compact subset
$\Lambda\subset Z(S)$  to construct a family of relatively compact
(in $X$) weakly $k$-PS manifolds of dimension $r$ in any given open
neighborhood $U(\Lambda)$,
\begin{equation}\label{e:3.7}
\{(\Phi^\Lambda_{\bf t})^{-1}(0)\,|\,{\bf t}\in
B_\epsilon(\R^l)_{res}\}
\end{equation}
with the following properties:
\begin{description}
\item[(i)] Each intersection $\Lambda\cap(\Phi^\Lambda_{\bf
t})^{-1}(0)$ is compact. (Such $(\Phi_{\bf t}^\Lambda)^{-1}(0)$ is
also called a {\bf local Euler chain} of $(X, E, S)$ near
$\Lambda$ later.)

\item[(ii)] Any two  ${\bf t},   {\bf t}'\in
B_\epsilon(\R^l)_{res}$ give a
  weakly $k$-PS manifold cobordism between $(\Phi^\Lambda_{{\bf t}})^{-1}(0)$ and $(\Phi^{\Lambda}_{{\bf
  t}'})^{-1}(0)$ in the sense that there exist generic paths $\gamma:[0,1]\to
B_\epsilon(\R^l)$ with $\gamma(0)={\bf t}$ and $\gamma(1)={\bf
t}'$ such that
$$(\Phi^\Lambda)^{-1}(\gamma):=\cup_{t\in
[0,1]}\{t\}\times(\Phi^\Lambda_{\gamma(t)})^{-1}(0)
$$
is a weakly $k$-PS manifold with boundary
$\{0\}\times(\Phi^{\Lambda}_{\bf
t})^{-1}(0)\cup(-\{1\}\times(\Phi^\Lambda_{{\bf t}'})^{-1}(0))$
which is relative compact in $[0,1]\times Cl(U(\Lambda))$.
\end{description}

 Now let $(X,E,S)$ be an oriented (weakly) PS Banach $k$-Fredholm
bundle of index $r$ and with $k\ge 1$.
 Furthermore, we also make the following assumptions: \\
 $\bullet$  The  zero locus $Z(S)$ is compact.\\
 $\bullet$ It has  rich PS sections of class ${\cal A}$
 and  is (weakly) $k$-good relative to the class ${\cal A}$;\\
 $\bullet$ The top stratum $X_0$
of $X$ has no boundary (including corners) as a Banach manifold
(near $Z(S)$).\\
 Then by Theorem~\ref{th:3.14}(ii), we get the Euler
class $e(E,S)=[\Phi_{\bf t}^{-1}(0)]\in H_r(X,\Z)$, given by  a
family of cobordant Euler cycles of $(X,E,S)$,
$$
\{\Phi_{\bf t}^{-1}(0)\,|\, {\bf t}\in
B_\varepsilon(\R^m)_{reg}\},
$$
 where each $\Phi_{\bf t}^{-1}(0)$ is
actually an oriented closed weakly $k$-PS manifold of dimension $r$.
If $f:X\to P$ is a PS map from $X$ to a manifold $P$ of finite
dimension,  by Proposition~\ref{prop:3.15} we get a class
$$
 f_\ast e(E,S)=f_\ast([\Phi_{\bf
t}^{-1}(0)])\in H_r(P,\Z).
$$
 Let $\alpha\in H^r(P, \R)$ and $\alpha^\ast$ be a
differential form representative of it. One easily checks that for
any ${\bf t}\in B_\varepsilon(\R^m)_{reg}$,
\begin{equation}\label{e:3.8}
\langle f^\ast\alpha, e(E,S)\rangle=\langle\alpha, f_\ast
e(E,S)\rangle =\int_{\Phi_{\bf t}^{-1}(0)}f^\ast\alpha^\ast,
\end{equation}
where the integration in the right side is understood as follows:
Multiply the form $f^\ast\alpha^\ast$ by  a PS function to $[0,1]$
which is $1$ outside a small neighborhood of $\Phi_{\bf
t}^{-1}(0)\setminus(\Phi_{\bf t}^{-1}(0))_{0}$ and $0$ inside a
smaller neighborhood and then pass to a limit as these
neighborhoods get smaller. Here $(\Phi_{\bf
t}^{-1}(0))_{0}=\Phi_{\bf t}^{-1}(0)\cap X_{0}$. Then almost
completely repeating the proof of Proposition~\ref{prop:1.11} we
get:

\begin{proposition}\label{prop:3.16}
({\bf First localization  formula}). If $\alpha$ has a
representative form $\alpha^\ast$ such that $f^\ast\alpha^\ast$
has support ${\rm supp}(f^\ast\alpha^\ast)$ contained in a compact
subset $\Lambda\subset Z(S)$ then for the family
$\{(\Phi^\Lambda_{\bf t})^{-1}(0)\,|\,{\bf t}\in
B_\epsilon(\R^l)_{res}\}$ in (\ref{e:3.7}) there exists a residual
subset $B_\epsilon(\R^l)^\star_{res}\subset
B_\epsilon(\R^l)_{res}$ such that
$$\langle e(E,S), f^\ast\alpha\rangle=
\int_{(\Phi^\Lambda_{\bf t})^{-1}(0)}f^\ast\alpha^\ast
$$
for ${\bf t}\in B_\epsilon(\R^l)^\star_{res}$.
\end{proposition}

The arguments above Proposition~\ref{prop:1.12} can also move to
here almost directly. Let $(X, E, S)$ be as in
Proposition~\ref{prop:3.16}.
 Recall that $X=\cup_{\alpha\in\Lambda}X_\alpha$. Denote
 by $r_\alpha$ the index of the restriction $S$ to $X_\alpha$.
We are interested in the following two classes of PS submanifolds of
$X$:\\
$\bullet$ $X^\alpha:=\cup_{\beta\prec\alpha}X_\beta$, where
$\alpha\in\Lambda$ and $\alpha\ne 0$.\\
$\bullet$ A PS submanifold
$Y=\cup_{\alpha\in\Lambda_0}Y_\alpha\subset X$ indexed by a partial
order subset $\Lambda_0\subset\Lambda$ with the maximal element
$\alpha_0$ is said to be of {\bf finite codimension $n$} if the top
stratum $Y_{\alpha_0}$ of $Y$ is a Banach submanifold of codimension
$n$ in $X_{\alpha_0}$ and other strata of $Y$ have at least finite
codimension $n$ in the corresponding strata of $X$.

  By definition~\ref{def:3.5} it is easily
proved that the (weakly) $k$-Fredholm PS section $S:X\to E$
restricts to a (weakly) $k$-Fredholm section
$S|_{X^\alpha}:X^\alpha\to E|_{X^\alpha}$. So $(X^\alpha,
E|_{X^\alpha}, S|_{X^\alpha})$ is also a (weakly) PS Banach
$k$-Fredholm bundle of index $r_\alpha$.

 By Proposition~\ref{prop:C.1}, the (weakly) $k$-Fredholm PS section
$S:X\to E$ restricts to a Fredholm section  $S|_{Y_{\alpha_0}}$ of
index $r_{\alpha_0}-n$.  Moreover, for other strata $Y_\alpha$,
since $\alpha\prec\alpha_0$ it follows from
Definition~\ref{def:3.5} that
$$
{\rm Index}(S|_{Y_{\alpha}})+ k\le r_\alpha-n+ k<r_{\alpha_0}-n
={\rm Index}(S|_{Y_{\alpha_0}}).
$$
These show that $(Y, E|_Y, S|_Y)$ is only a {\bf weakly} PS Banach
$k$-Fredholm bundle of index $r_{\alpha_0}-n$. Clearly, if $Y$ is
closed in $X$  then $Z(S|_Y)$ is also closed and thus compact in
$Y$ because of the compactness of $Z(S)$.

In the following we always assume that the PS submanifold  $Y$
belong to the above two classes without special statements, and also
assume:\\
$\bullet$ $Y$ is  closed in $X$. \\
Then $S|_Y$ is locally proper relative to class ${\cal A}|_Y$ near
$Z(S|_Y)$. Here ${\cal A}|_Y$ is the class of PS sections of
$E|_Y$ obtained by restricting the PS sections of $E$ of class $A$
to $Y$. One also easily prove that the bundle $(Y, E|_Y,
S|_Y)$ has also rich PS sections of class ${\cal A}|_Y$. We conclude:\\
$\bullet$ $S|_Y$ is  weakly $k$-quasi transversal to the zero
section relative to class ${\cal A}|_Y$. \\
In fact, for any $x\in
Z(S|_Y)$, if PS sections $\sigma_1|_Y,\cdots, \sigma_m|_Y$ of $E|_Y$
of class ${\cal A}|_Y$ satisfies
$$
D(S|_Y)(x)(T_xY)+ {\rm span}(\{\sigma_1(x),\cdots,
\sigma_m(x)\})=E_x
$$
then the relation  $T_xY\subset T_xX$ implies
$$
D(S)(x)(T_xX)+ {\rm span}(\{\sigma_1(x),\cdots,
\sigma_m(x)\})=E_x.
$$
Since $S$ is (weakly) $k$-quasi transversal to the zero section
relative to class ${\cal A}$ there exist a small open neighborhood
$U(x)$ of $x$ in $X$ and $\eta>0$ such that
$$
U(x)\times B_\eta(\R^m)\to \Pi_1^\ast(E|_{U(x)}),\;
 (x; t_1,\cdots, t_m)\mapsto S(x)+\sum^m_{i=1}t_i\sigma_i(x),
$$
is a (weakly) $k$-Fredholm PS section of index $r+m$ that is
transversal to the zero section, where $\Pi_1: U(x)\times
B_\eta(\R^m)\to U(x)$ is the natural projection. So the restriction
section
\begin{eqnarray*}
&&U(x)\cap Y\times B_\eta(\R^m)\to \Pi_1^\ast(E|_{U(x)\cap Y}),\\
&& (x; t_1,\cdots, t_m)\mapsto
S|_Y(x)+\sum^m_{i=1}t_i\sigma_i|_Y(x)
\end{eqnarray*}
is transversal to the zero section as well. Note that this section
is  (weakly) $k$-Fredholm (resp. weakly $k$-Fredholm ) if $Y$ is
of form $X^\alpha$ (resp. a PS submanifold of finite codimension
as the above second class). In summary, under the above
assumptions we have proved:
\begin{center}
{\bf $(Y, E|_Y, S|_Y)$ is a (weakly) PS Banach $k$-Fredholm bundle
which is (weakly) $k$-good if $Y=X^\alpha$, and a weakly PS Banach
$k$-Fredholm bundle which is weakly $k$-good if $Y$ is a PS
submanifold of finite codimension that is closed in $X$.}
\end{center}
The following is a generalization of Proposition~\ref{prop:1.12}.

\begin{proposition}\label{prop:3.17}
({\bf Second localization  formula}). Let $(X,E,S)$ be an oriented
(weakly) PS Banach $k$-Fredholm bundle of index $r$ and with
compact zero locus $Z(S)$ and $k\ge 1$, and the closed subset
$Y\subset X$ be one of the above two classes of PS submanifolds.
Assume also that $(X,E,S)$ has  rich PS sections of class ${\cal
A}$ and  is (weakly) $k$-good relative to the class ${\cal A}$.
Then there exist Euler chains $M$ of $(X, E, S)$ and $N$ of $(Y,
E|_Y, S|_Y)$ such that $M\cap Y=N$. Moreover, if  the top strata
of $X$ and $Y$ have no boundary (including corner), then
 the above  $M$ and $N$ can be chosen as  closed weakly $k$-PS
manifolds.  In particular, for a smooth map $f$ from $X$ to an
oriented smooth manifold $P$ of finite dimension, assume that it
is transversal to an oriented closed submanifold $Q\subset P$ and
$Y:=f^{-1}(Q)\subset X$ is a PS Banach submanifold of finite
codimension whose top stratum is contained in $X_0$. Then when
both $e(E,S)$ and $e(E|_Y, S|_Y)$ exist, it holds that
 $$
 I_{Q\ast}(a)\cdot_P f_\ast(e(E,S))
 = a\cdot_Q
f_\ast(e(E|_Y, S|_Y))
$$
 for any class $a\in H_\ast(Q, \R)$ of codimension $r$,
where $I_Q:Q\to P$  is the inclusion map, $\cdot_P$ (resp.
$\cdot_Q$) is the  intersection product in $P$ (resp. $Q$).
\end{proposition}

\noindent{\bf Proof.} Under the present assumptions it follows from
the proof of Theorem~\ref{th:3.14} that  there exist an open
neighborhood ${\cal W}$ of $Z(S)$ in $X$, $\varepsilon>0$ and PS
sections $S_i:X\to E$, $i=1,\cdots,m$, such that:
\begin{description}
\item[(i)]  $\bullet$ the PS section
$$
\Phi: {\cal W}\times B_\varepsilon(\R^m)\to \Pi_1^\ast(E|_{{\cal
W}}),\; y\mapsto S(y)+\sum^m_{i=1}t_i\sigma_i(y)
$$
is a (weakly) $k$-Fredholm PS section of index $r+m$ that is
transversal to the zero section;

 $\bullet$ the PS section
$$\Phi^Y: {\cal W}\cap Y\times B_\varepsilon(\R^m)\to \Pi_1^\ast(E|_{{\cal
W}\cap Y}),\; y\mapsto S|_Y(y)+\sum^m_{i=1}t_i\sigma_i|_Y(y)
$$
is transversal to the zero section and is also a weakly
$k$-Fredholm PS section (of index $r_\alpha+ m$ in the case
$Y=X^\alpha$ or $r_{\alpha_0}-n+ m$ in another case).

\item[(ii)] Both zero sets $\Phi^{-1}(0)\subset Cl({\cal W}\times
B_\varepsilon(\R^m))$ and $(\Phi^Y)^{-1}(0)\subset Cl({\cal W}\cap
Y\times B_\varepsilon(\R^m))$ are compact.
\end{description}
 Let $B_\varepsilon(\R^m)_{res}\subset B_\varepsilon(\R^m)$
be the set of the common regular values of both projections
\begin{eqnarray*}
&&\Phi^{-1}(0)\to B_\varepsilon(\R^m),\;(y,{\bf t})\mapsto{\bf
t}\quad{\rm and}\\
&& (\Phi^Y)^{-1}(0)\to B_\varepsilon(\R^m),\;(y,{\bf
t})\mapsto{\bf t}.
\end{eqnarray*}
 Then as in the arguments below (\ref{e:3.4}),
 for each ${\bf t}\in B_\varepsilon(\R^m)_{res}$
$$
\Phi_{\bf t}: {\cal W}\to E,\; y\mapsto
S(y)+\sum^m_{i=1}t_i\sigma_i(y)
$$
is a (weakly) $k$-Fredholm PS section that is transversal to the
zero section, and
$$\Phi^Y_{\bf t}: {\cal W}\cap Y\to E|_Y,\; y\mapsto
S|_Y(y)+\sum^m_{i=1}t_i\sigma_i|_Y(y)
$$
is a weakly $k$-Fredholm section that is transversal to the zero
section. Since ${\cal W}\cap Y$ is an open neighborhood of
$Z(S|_Y)=Z(S)\cap Y$ in $Y$, it is clear that
\begin{equation}\label{e:3.9}
 \Phi_{\bf t}^{-1}(0)\cap Y=(\Phi^Y_{\bf
t})^{-1}(0).
\end{equation}
This proves the first claim. The second claim is also obvious. Other
conclusions may be easily proved as in Proposition~\ref{prop:1.12}
as well. \hfill$\Box$\vspace{2mm}

Now let us state a generalization of Proposition~\ref{prop:1.13}.
Let $(X, E', S')$ (resp. $(X, E'', S'')$) be an oriented (weakly) PS
Banach $k$-Fredholm bundle with rich PS sections of class ${\cal
A}'$ (resp. ${\cal A}''$) and compact zero locus $Z(S')$ (resp.
$Z(S'')$). Consider the direct sum $E=E'\oplus E''$, $S=S'\oplus
S'':X\to E$ and
 ${\cal A}={\cal A}'\oplus{\cal A}''$. Clearly,
$Z(S)=Z(S')\cap Z(S'')$, $Z(S)=Z(S'|_{Z(S'')})$, and $(X, E, S)$
has rich PS sections of class ${\cal A}$ near $Z(S)$. Suppose that
both $(X, E', S')$ and $(X, E'', S'')$ are (weakly) $k$-good. We
conclude:\\
$\bullet$ $(X, E, S)$ is also (weakly) $k$-good.\\
  In fact, for $x\in Z(S)$ let PS sections
$\sigma_1,\cdots,\sigma_m$ of $E$ of class ${\cal A}$ be such that
$$
DS(x)(T_xX)+ {\rm span}(\{\sigma_1(x),\cdots,\sigma_m(x)\})=E_x.
$$
Since $E_x=E'_x\oplus E''_x$, $DS(x)(v)=DS'(x)(v)\oplus
DS''(x)(v)$ for $v\in T_xX$ and
$\sigma_i=\sigma'_i\oplus\sigma''_i$, $i=1,\cdots,m$, the above
equality implies
\begin{eqnarray*}
DS'(x)(T_xX)+ {\rm span}(\{\sigma'_1(x),\cdots,\sigma'_m(x)\})=E'_x,\\
DS''(x)(T_xX)+ {\rm
span}(\{\sigma''_1(x),\cdots,\sigma''_m(x)\})=E''_x.
\end{eqnarray*}
 Hence that $(X, E', S')$ and $(X, E'', S'')$ are (weakly) $k$-good implies that
 there exist an open neighborhood $U(x)$ and $\eta>0$ such that:
\begin{description}
\item[(i)]  $\bullet$ both sections
$$
\Phi': U(x)\times B_\varepsilon(\R^m)\to \Pi_1^\ast(E'|_{U(x)}),\;
y\mapsto S'(y)+\sum^m_{i=1}t_i\sigma'_i(y)
$$
and
$$
\Phi'': U(x)\times B_\varepsilon(\R^m)\to
\Pi_1^\ast(E''|_{U(x)}),\; y\mapsto
S''(y)+\sum^m_{i=1}t_i\sigma''_i(y)
$$
are (weakly) $k$-Fredholm PS sections of indexes ${\rm Index}(S')+
m$ and ${\rm Index}(S'')+ m$ respectively  that are transversal to
the zero section.

\item[(ii)] Both zero sets $(\Phi')^{-1}(0)\subset Cl(U(x)\times
B_\varepsilon(\R^m))$ and $(\Phi'')^{-1}(0)\subset Cl(U(x)\times
B_\varepsilon(\R^m))$ are compact.
\end{description}
Consider the PS section
\begin{eqnarray*}
&& \Phi: U(x)\times B_\varepsilon(\R^m)\to
\Pi_1^\ast(E|_{U(x)}),\\
&&\qquad y\mapsto S'(y)\oplus
S''(y)+\sum^m_{i=1}t_i\sigma'_i(y)\oplus\sigma''_i(y).
\end{eqnarray*}
 Clearly, it is  (weakly) $k$-Fredholm and transversal to the zero
section (because its zero locus $Z(\Phi)$ is equal to
$Z(\Phi')\cap Z(\Phi'')$). It follows that $S$ is (weakly)
$k$-quasi transversal to the zero section relative to class ${\cal
A}$. Moreover, by the compactness of $Z(\Phi)=Z(\Phi')\cap
Z(\Phi'')$ we get that  $S=S'\oplus S''$ is also proper relative
to class ${\cal A}$ near $Z(S)$. Hence $(X, E, S)$ is  (weakly)
$k$-good.

\begin{proposition}\label{prop:3.18}
Let $(X, E', S')$ and $(X, E'', S'')$ be as above, and have indexes
$r'$ and $r''$ respectively. Furthermore, assume that $DS''|_{Z(S)}$
is surjective and $Z(S)\cap X_0\ne\emptyset$. Then
\begin{description}
\item[(a)] For some small open neighborhood ${\cal U}$ of $Z(S)$
in $X$ the intersection  $Z(S'')^\star:=Z(S'')\cap{\cal U}$ is a
weakly $k$-PS manifold of dimension $r''$. (In particular, the top
stratum of $Z(S'')^\star$ is contained in $X_0$.)

\item[(b)] $D(S'|_{Z(S'')}):(TZ(S''))|_{Z(S)} \to E'|_{Z(S)}$ is a
PS Fredholm bundle map of index ${\rm Ind}(D(S'|_{Z(S'')})={\rm
Ind}(DS)$ (which precisely means
$${\rm
Ind}(D(S'|_{Z(S'')})(x))={\rm Ind}(DS(x))\;\forall x\in Z(S)\;).
$$
Moreover, if  $(Z(S'')^\star, E'|_{Z(S'')^\star},
S'|_{Z(S'')^\star})$ is (weakly) $k$-good relative to the class
${\cal A}'|_{Z(S'')^\star}$ then there exist Euler chains $N$ of
it and $M$ of $(X, E,S)$ such that $M=N$. Furthermore, if both
$e(E, S)$ and $e(E'|_{Z(S'')^\star}, S'|_{Z(S'')^\star})$ exist
then $M$ and $N$ can be chosen as closed weakly $k$-PS manifolds,
and therefore
$$
e(E,S) =(i_{Z(S'')^\star})_\ast e(E'|_{Z(S'')^\star},
S'|_{Z(S'')^\star}).
$$
Here $(i_{Z(S'')^\star})_\ast$ is the homomorphism between the
homology groups induced by the inclusion $i_{Z(S'')^\star}:
Z(S'')^\star\hookrightarrow X$.
\end{description}
\end{proposition}

\noindent{\bf Proof.}\quad
  (a) Since $(X, E'', S'')$ is (weakly) $k$-good, $S''$ is (weakly) $k$-quasi transversal to the
zero section relative to class ${\cal A}''$. Note that
$DS''|_{Z(S)}$ is surjective. By Remark~\ref{rem:3.11}(i), for each
$x\in Z(S)$ there exists a small open neighborhood $U(x)$ of $x$ in
$X$ such that $S'':U(x)\to E''|_{U(x)}$ is transversal the zero
section. It follows from the compactness of $Z(S)$ that there exists
a small open neighborhood ${\cal U}$ of $Z(S)$ in $X$ such that
$S'':{\cal U}\to E''|_{\cal U}$ is transversal the zero section.
Moreover, our assumption $Z(S)\cap X_0\ne\emptyset$ implies
$Z(S')\cap X_0\ne\emptyset$ and $Z(S'')\cap X_0\ne\emptyset$. So
$Z(S'')^\star:=Z(S'')\cap{\cal U}$ is a weakly $k$-PS manifold of
dimension $r''$ and its top stratum is also intersecting with
$Z(S'|_{Z(S'')})=Z(S)$.

 (b)  Applying the first claim in Proposition~\ref{prop:1.13}(b) to each stratum
of $Z(S'')^\star$ we immediately get the first conclusion.

To prove the second claim note that $Z(S'')^\star$ has finite
dimension and $DS''|_{Z(S)}$ is surjective. We can find finitely
many PS sections $S'_{j}:X\to E'$, $j=1,\cdots, m$ such that for
each $y\in Z(S)$,
\begin{equation}\label{e:3.10}
\left.\begin{array}{ll}
D(S'|_{Z(S'')})(y)(T_yZ(S'')^\star)+\\
 \qquad\qquad\qquad {\rm span}(\{S'_{1}(y),\cdots,
S'_{m}(y)\}=E'_y,\\
DS'(y)(T_yX)+ {\rm span}(\{S'_{1}(y),\cdots,
S'_{m}(y)\}=E'_y,\\
DS(y)(T_yX)+ {\rm span}(\{S'_{1}(y),\cdots, S'_{m}(y)\}=E_y.
\end{array}
\right\}
\end{equation}
Since $(Z(S'')^\star, E'|_{Z(S'')^\star}, S'|_{Z(S'')^\star})$ is
(weakly) $k$-good relative to the class ${\cal A}'|_{Z(S'')^\star}$,
$(X, E', S')$ is (weakly) $k$-good relative to the class ${\cal A}'$
and $(X, E, S)$ is (weakly) $k$-good relative to the class ${\cal
A}\supset{\cal A}'$ it follows from (\ref{e:3.10}) that
 there exist an open
neighborhood ${\cal W}$ of $Z(S)$ in ${\cal U}$ and
$\varepsilon>0$ such that:
\begin{description}
\item[($1^\circ$)]  The PS sections
\begin{eqnarray*}
&&\Phi: {\cal W}\times B_\varepsilon(\R^m)\to E,\; y\mapsto
S(y)+\sum^m_{i=1}t_iS'_i(y),\nonumber\\
&&\Phi': {\cal W}\times B_\varepsilon(\R^m)\to E',\; y\mapsto
S'(y)+\sum^m_{i=1}t_i S'_i(y),\nonumber\\
&&\Psi: {\cal W}\cap Z(S'')^\star\times B_\varepsilon(\R^m)\to
E'|_{{\cal W}\cap Z(S'')^\star},\nonumber\\
&&\qquad y\mapsto S'|_{Z(S'')}(y)+\sum^m_{i=1}t_i
S'_i|_{Z(S'')}(y),\nonumber
\end{eqnarray*}
are (weakly) $k$-Fredholm and transversal to the zero section.
\item[($2^\circ$)] The zero sets $\Phi^{-1}(0)\subset Cl({\cal
W}\times B_\varepsilon(\R^m))$, $(\Phi')^{-1}(0)\subset Cl({\cal
W}\times B_\varepsilon(\R^m))$ and $(\Psi)^{-1}(0)\subset Cl({\cal
W}\cap Z(S'')^\star\times B_\varepsilon(\R^m))$ are all compact
subsets.
\end{description}
Let $B_\varepsilon(\R^m)_{res}\subset B_\varepsilon(\R^m)$ be the
set of the common regular values of the following three
projections
\begin{eqnarray*}
&&\Phi^{-1}(0)\to B_\varepsilon(\R^m),\;(y,{\bf t})\mapsto{\bf t},\\
&&(\Phi')^{-1}(0)\to B_\varepsilon(\R^m),\;(y,{\bf t})\mapsto{\bf
t},\\
&&\Psi^{-1}(0)\to B_\varepsilon(\R^m),\;(y,{\bf t})\mapsto{\bf t}.
\end{eqnarray*}
 Then by Remark~\ref{rem:3.11}(iii), for each ${\bf t}\in B_\varepsilon(\R^m)_{res}$ the PS sections
\begin{eqnarray*}
&&\Phi_{\bf t}: {\cal W}\to E|_{{\cal W}}, \quad x\mapsto
S(x)+\sum^m_{i=1}t_iS'_i(x),\nonumber\\
&&\Phi'_{{\bf t}}: {\cal W}\to E'|_{{\cal W}},
 \quad x\mapsto S'(x)+\sum^m_{i=1}t_iS'_i(x),\nonumber\\
&&\Psi_{\bf t}: {\cal W}\cap Z(S'')^\star\to E'|_{{\cal W}\cap
Z(S'')^\star}, \quad x\mapsto S'(x)+\sum^m_{i=1}t_iS'_i(x),\nonumber
\end{eqnarray*}
are all (weakly) $k$-Fredholm and transversal to the zero section.
Let $R_{\bf t}=\Phi'_{\bf t}|_{\cal N}$.  Note that ${\rm
P}_{E'}\circ\Phi_{\bf t}=\Phi_{\bf t}'$, ${\rm
P}_{E''}\circ\Phi_{\bf t}=S''$ and ${\rm P}_{E'}\circ\Phi_{\bf
t}+{\rm P}_{E''}\circ\Phi_{\bf t}=\Phi_{\bf t}$. As in the proof of
Proposition~\ref{prop:1.13}  we can derive that $
 \Phi_{\bf t}^{-1}(0)=R_{\bf
t}^{-1}(0)$ for ${\bf t}\in B_\varepsilon(\R^m)_{res}$ small enough.
By the assumptions, $\Phi$, $\Phi'$ and $\Psi$ have all zero loci
intersecting with the top strata of their domains. For ${\bf t}\in
B_\varepsilon(\R^m)_{res}$ small enough, $\Phi_{\bf t}^{-1}(0)$,
$(\Phi'_{\bf t})^{-1}(0)$ and $\Psi_{\bf t}^{-1}(0)$ have nonempty
intersections with the top strata of ${\cal W}$ and ${\cal W}\cap
Z(S'')^\star$ as well. The desired result follows.
 \hfill$\Box$\vspace{2mm}

Corresponding to Proposition~\ref{prop:1.15} we easily get:

\begin{proposition}\label{prop:3.19}  Let $(X_i, E_i, S_i)$ be
 (weakly) PS Banach $k$-Fredholm bundles of index $r_i$ and
with compact zero locus $Z(S_i)$, $i=1,2$. Assume that $(X_1, E_1,
S_1)$ (resp. $(X_2, E_2, S_2)$) is (weakly) $k$-good relative to the
class ${\cal A}_1$ (resp. ${\cal A}_2$) and that $Z(S_1)$ (resp.
$Z(S_2)$) has a nonempty intersection with the top stratum  of $X_1$
(resp. $X_2$). Then the natural product $(X_1\times X_2, E_1\times
E_2, S_1\times S_2)$ is a (weakly) PS Banach $k$-Fredholm bundle of
index $r=r_1+ r_2$ that is (weakly) $k$-good relative to the class
${\cal A}_1\times{\cal A}_2$, and there exist Euler chains $M_i$ of
$(X_i, E_i,S_i)$, $i=1,2$, such that the product weakly $k$-PS
manifold $M=M_1\times M_2$ is a Euler chain of $(X_1\times X_2,
E_1\times E_2, S_1\times S_2)$. Moreover, if $e(X_1, E_1)$ and
$e(X_2, E_2)$ exist, then $e(E_1\times E_2, S_1\times S_2)$ also
exists and $M_i$ can be chosen as closed weakly $k$-PS manifolds,
$i=1,2$.
\end{proposition}

Similar to Proposition~\ref{prop:1.16} we  also have:

\begin{proposition}\label{prop:3.20}
({\bf Homotopy})  Let $(X, E^{(i)}, S^{(i)})$ be oriented (weakly)
PS Banach $k$-Fredholm bundles of index $r$ and with compact zero
loci, and be also (weakly) $k$-good relative to the classes ${\cal
A}_i$, $i=0,1$. They are said to be {\bf oriented homotopic} if
there exists an oriented  (weakly) PS Banach $k$-Fredholm bundle
$(X\times [0, 1], E, S)$ of index $r+1$ and with compact zero locus
that is (weakly) $k$-good relative to some class ${\cal A}$, such
that $(X\times \{i\}, E|_{X\times\{i\}}, S|_{X\times\{i\}})=(X,
E^{(i)}, S^{(i)})$ and ${\cal A}|_{X\times\{i\}}={\cal A}_i$, $i=0,
1$.
 In this case, if $r$-dimensional Euler classes
$e(E^{(0)}, S^{(0)})$ and $e(E^{(1)}, S^{(1)})$ exist then these two
Euler classes are same.
 \end{proposition}

We can also define the notion of a {\bf morphism} from oriented
(weakly) PS Banach $k$-Fredholm bundles $(X, E, S)$ to $(X', E',
S')$ and get a corresponding generalization of
Proposition~\ref{prop:1.17}. A {\bf PS Fredholm map of index $d$}
from PS Banach manifolds $Y$ to $X$ is a PS map whose restriction to
the top stratum (resp. each lower stratum) is a Fredholm map of
index $d$ (resp. index less than or equal to $d$).

The following proposition is a generalization of
Proposition~\ref{prop:1.18}. Its proof is still given for reader's
conveniences.

\begin{proposition}\label{prop:3.21} Let $(X, E, S)$ be a
(weakly) PS Banach $k$-Fredholm bundle of index $r$ and with compact
zero locus $Z(S)$.  Assume that $Z(S)\cap X_0\ne\emptyset$ and $(X,
E, S)$  is also (weakly) $k$-good relative some class ${\cal A}$.
Let $Y$ be another PS Banach manifold admitting PS cut-off functions
with arbitrarily small support near any point, and $f:Y\to X$ be a
proper PS Fredholm map of index $d$ from $Y$ to $X$. Assume that $f$
is a PS submersion and $f^{-1}(X_0\cap Z(S))\cap Y_0\ne\emptyset$.
Then the natural pullback $(Y, f^\ast E, f^\ast S)$ is a weakly  PS
Banach $k$-Fredholm bundle of index $r+d$ and with compact zero
locus $Z(f^\ast S)$   that is weakly $k$-good relative to the class
$f^\ast{\cal A}$, and there exists a Euler chain $M$ (resp. $N$) of
$(X, E, S)$ (resp. $(Y, f^\ast E, f^\ast S)$)  such that
$f^{-1}(M)=N$. Furthermore, if both Euler classes $e(E,S)$ and
$e(f^\ast E, f^\ast S)$ exist then $M$ and $N$ can be chosen as
closed weakly $k$-PS manifolds.
\end{proposition}

As noted below Proposition~\ref{prop:1.18}, if $d=0$, $X_0$ is
connected and the Euler classes $e(E,S)$ and $f_\ast(e(f^\ast E,
f^\ast S))$ exist, then the final claim in
Proposition~\ref{prop:3.21} implies
$$
f_\ast(e(f^\ast E, f^\ast S))={\rm deg}(f) e(E,S)
$$
for some kind topological degree ${\rm deg}(f)$ of
$f$.\vspace{2mm}

\noindent{\bf Proof.}\quad Since  $Z(S)$ is compact, so is $Z(f^\ast
S)=f^{-1}(Z(S))$ by the properness of $f$. Clearly, $f^{-1}(X_0\cap
Z(S))\cap Y_0\ne\emptyset$ implies that $Z(f^\ast S)\cap
Y_0\ne\emptyset$.
  From the proof of
Proposition~\ref{prop:1.18} it is not hard to see that pullback $(Y,
f^\ast E, f^\ast S)$ has index $r+d$.  One also easily checks that
$(Y, f^\ast E, f^\ast S)$ has rich local PS sections of class
$f^\ast{\cal A}$. Moreover, for $x\in Z(S)$ and PS sections of $E$,
$S_1,\cdots, S_m$, assume that there exist an open neighborhood
$U(x)$ of $x$ in $X$ and $\delta>0$ such that
$$
\bigcup_{|{\bf t}|\le\delta}Cl(U(x))\cap(S+
t_i\sum^m_{i=1}S_i)^{-1}(0)
$$
is compact. Then
\begin{eqnarray*}
&&\quad\bigcup_{|{\bf t}|\le\delta}f^{-1}(Cl(U(x)))\cap(f^\ast S+
t_i\sum^m_{i=1}f^\ast S_i)^{-1}(0)\\
&&=\bigcup_{|{\bf t}|\le\delta}f^{-1}\Bigl(Cl(U(x))\cap (S+
t_i\sum^m_{i=1} S_i)^{-1}(0)\Bigr)
\end{eqnarray*}
is also compact because $f$ is proper. So $f^\ast S$ is  locally
proper relative to the class $f^\ast{\cal A}$.

We conclude that {\bf $f^\ast S$ is weakly $k$-quasi transversal to
the zero section relative to the class $f^\ast{\cal A}$.} In fact,
for any $y\in Z(f^\ast S)$ assume that there exist PS sections
$f^\ast\sigma_1,\cdots, f^\ast\sigma_p$ of $f^\ast E$ of class
$f^\ast{\cal A}$ such that
$$
D(f^\ast S)(y)(T_yY)+ {\rm span}(\{f^\ast\sigma_1(y),\cdots,
f^\ast\sigma_p(y)\})=(f^\ast E)_y=E_{f(y)}.
$$
Then $x=f(y)\in Z(S)$ and
$$
DS(x)(T_xX)+ {\rm span}(\{f^\ast\sigma_1(y),\cdots,
f^\ast\sigma_p(y)\})=(f^\ast E)_y=E_{f(y)}.
$$
because $D(f^\ast S)(y)=DS(x)\circ df(y)$ and $DS(x)\circ
df(y)(T_yY)\subset DS(x)(T_xX)$. Note that $S$ is (weakly) $k$-quasi
transversal to the zero section. There exist an open neighborhood
$U(x)$ of $x$ and $\delta>0$ such that the PS section
$$
\Psi: U(x)\times B_\delta(\R^p)\to \Pi_1^\ast(E|_{U(x)}),\; (x, {\bf
t})\mapsto S(x)+\sum^p_{i=1}t_i\sigma_i(x)
$$
 is  (weakly) $k$-Fredholm and also transversal to the zero section.
Take a small open neighborhood $V(y)$ of $y$ so that $f(V(y))\subset
U(x)$. Consider the section
\begin{eqnarray*}
f^\ast\Psi: V(y)\times B_\delta(\R^p)\to
\Pi_1^\ast(f^\ast E|_{V(y)}),\\
 (y,{\bf t})\mapsto f^\ast
S(y)+\sum^p_{i=1}t_i (f^\ast \sigma_i)(y).
\end{eqnarray*}
 For any zero $(y, {\bf t})$ of $f^\ast\Psi$,
one easily checks that $(f(y),{\bf t})\in Z(\Psi)$ and the vertical
differentials of $f^\ast\Psi$ at $(y,{\bf t})$ is given by
$$D(f^\ast\Psi)(y,{\bf t})(\xi, {\bf
u})=D\Bigl(S+\sum^p_{i=1}t_i\sigma_i\Bigr)(f(y))\circ df(y)(\xi)+
\sum^p_{i=1}u_i\sigma_i(f(y)).
$$
Let us prove that $D(f^\ast\Psi)(y,{\bf t})$ is surjective.  For any
given $\zeta\in (f^\ast E)_y=E_{f(y)}$. Since $D\Psi(x,{\bf
t}):T_xX\times\R^p\to E_x$ is surjective we have $(\eta, {\bf v})\in
T_xX\times\R^p$ such that $D\Psi(x,{\bf t})(\eta,{\bf v})=\zeta$,
i.e.
$$
D\Bigl(S+\sum^p_{i=1}t_iS_i\Bigr)(x)(\eta)+
\sum^p_{i=1}v_i\sigma_i(x)=\zeta.
$$
Since $f$ is a PS submersion there exists $\xi\in T_yY$ such that
$df(y)(\xi)=\eta$. It follows  that $D(f^\ast\Psi)(y,{\bf t})(\xi,
{\bf v})=\zeta$. Consequently, {\bf  $f^\ast\Psi$ is transversal to
the zero section.}

Let $S^{\alpha t}, \sigma_i^{\alpha t}: W\subset U(x)\to E_z$ be the
local representatives of the restrictions of $S$ and $\sigma_i$ to
the stratum $U(x)_\alpha$, $i=1,\cdots,p$. Then the restriction of
$\Psi$ to $U(x)_\alpha\times B_\delta(\R^p)$ has a local
representative
$$
\Psi^{\alpha t}: W\times B_\delta(\R^p)\to E_z,\; (u, {\bf
t})\mapsto S^{\alpha t}(u)+\sum^p_{i=1}t_i\sigma^{\alpha t}_i(u).
$$
Correspondingly, the restriction of $f^\ast\Psi$ to
$V(y)_\alpha\times B_\delta(\R^p)$ has a local representative
$$
(f^\ast\Psi)^{\alpha t}: f^{-1}(W)\times B_\delta(\R^p)\to E_{z},
 (y,{\bf t})\mapsto
S^{\alpha t}(f(v))+\sum^p_{i=1}t_i \sigma^{\alpha t}_i(f(y)).
$$
Clearly, this representative may be viewed a composition of
$\Psi^{\alpha t}$ with the Fredholm map
$$
f^{-1}(W)\times B_\delta(\R^p)\to W\times B_\delta(\R^p),\,(u,{\bf
t})\mapsto (f(u),{\bf t})
$$
of index ${\rm Ind}(f)$. As in the arguments of
Remark~\ref{rem:3.11}(iii), one easily proves that
$(f^\ast\Psi)^{\alpha t}$ is a Fredholm map of index
$$
p+ {\rm
Ind}(S|_{U(x)_\alpha})+ d_\alpha\le p+ {\rm Ind}(S|_{U(x)_\alpha})+
d.
$$
It follows that $f^\ast\Psi$ is a weakly $k$-Fredholm PS section.
Summarizing up the above arguments, we get that {\bf $(Y, f^\ast E,
f^\ast S)$ is weakly $k$-good relative to the class $f^\ast{\cal
A}$.}

In order to prove other conclusions, assume that
$$
\Phi: {\cal W}\times B_\varepsilon(\R^m)\to \Pi_1^\ast(E|_{\cal
W}),\; (x, {\bf t})\mapsto S(x)+\sum^m_{i=1}t_iS_i(x)
$$
 is a (weakly) $k$-Fredholm PS section of index $r+m$ that is transversal to the zero section.
 Since $Z(S)\cap X_0\ne\emptyset$,  $\Phi^{-1}(0)$ is intersecting with the top strata of
 ${\cal W}\times B_\varepsilon(\R^m)$.
Let ${\cal W}^\ast=f^{-1}({\cal W})$. Consider the section
\begin{eqnarray*}
f^\ast\Phi: {\cal W}^\ast\times B_\varepsilon(\R^m)\to
\Pi_1^\ast(f^\ast E|_{{\cal W}^\ast}),\\
 (y,{\bf t})\mapsto f^\ast
S(y)+\sum^m_{i=1}t_i (f^\ast S_i)(y).
\end{eqnarray*}
 For any $(y, {\bf t})\in Z(f^\ast\Phi)$,
one easily checks that $(f(y),{\bf t})\in Z(\Phi)$ and the vertical
differentials of $f^\ast\Phi$ at $(y,{\bf t})$ is given by
$$D(f^\ast\Phi)(y,{\bf t})(\xi, {\bf
u})=D\Bigl(S+\sum^m_{i=1}t_iS_i\Bigr)(f(y))\circ df(y)(\xi)+
\sum^m_{i=1}u_iS_i(f(y)).
$$
As above we can prove that $D(f^\ast\Phi)(y,{\bf t})$ is surjective,
and thus that {\bf $f^\ast\Phi$ is transversal to the zero section.}

Let $B_\varepsilon(\R^m)_{res}$ be the set of the common regular
values of the projections
\begin{eqnarray*}
&&\Phi^{-1}(0)\to B_\varepsilon(\R^m),\;(x,{\bf t})\to {\bf
t}\quad{\rm and}\\
&& (f^\ast\Phi)^{-1}(0)\to B_\varepsilon(\R^m),\;(y,{\bf t})\to
{\bf t}.
\end{eqnarray*}
Then it is a residual subset in $B_\varepsilon(\R^m)$, and for each
${\bf t}\in B_\varepsilon(\R^m)_{res}$ both the (weakly) $k$-PS
sections
 $\Phi_{\bf t}:{\cal W}\to E|_{\cal W}$ and
 $(f^\ast\Phi)_{\bf t}=f^\ast(\Phi_{\bf t}):{\cal
W}^\ast\to (f^\ast E)|_{{\cal W}^\ast}$ are transversal to the zero
sections.  (Note that $Z(f^\ast S)\cap Y_0\ne\emptyset$ implies that
$(f^\ast\Phi)^{-1}(0)$ and the top strata of ${\cal W}^\ast\times
B_\varepsilon(\R^m)$ have a nonempty intersection.)  So the sets
\begin{eqnarray*}
&&\Phi_{\bf t}^{-1}(0):=\{x\in{\cal W}\,|\,\Phi(x, {\bf
t})=0\}\quad{\rm and}\\
&& (f^\ast\Phi)_{\bf t}^{-1}(0):=\{y\in{\cal
W}^\ast\,|\,f^\ast\Phi(x, {\bf t})=0\}
\end{eqnarray*}
 are respectively compact weakly $k$-PS
manifolds of dimensions $r$ and $r+d$ for small ${\bf t}\in
B_\varepsilon(\R^m)_{res}$. Now as in the proof of
Proposition~\ref{prop:1.16} we can prove that $f^{-1}(\Phi_{\bf
t}^{-1}(0))=(f^\ast\Phi)_{\bf t}^{-1}(0)$ for ${\bf t}\in
B_\varepsilon(\R^m)_{res}$ small enough. Finally, it is clear that
both $\Phi_{\bf t}^{-1}(0)$ and $(f^\ast\Phi)_{\bf t}^{-1}(0)$ are
closed weakly $k$-PS manifolds when the Euler classes $e(E, S)$ of
dimension $r$ and $e(f^\ast E, f^\ast S)$ of dimension $r+d$ exist.
 \hfill$\Box$\vspace{2mm}

\begin{remark}\label{rem:3.22}
{\rm The weakly $k$-Fredholm PS section may be defined for the
weakly PS Banach Fredholm bundle. One only need to require that the
restrictions of it to those strata which are not open in $X$ have
Fredholm indexes less at least $k$ than the Fredholm index of the
restriction of it to the (open) top stratum. In this case unless
Propositions~\ref{prop:3.17},~\ref{prop:3.18} need be suitably
modified other results in this section still hold.}
\end{remark}

In the constructions of Liu-Tian's virtual module cycles for Floer
homology and Gromov-Witten invariants the uniformizer charts of PS
Banach orbifolds  are PS Banach manifolds.  In
Appendix~\ref{app:D} we shall prove the existence of PS cut-off
functions on such PS Banach manifolds.

\pagebreak

\section{Virtual Euler chains and cycles of PS Banach Fredholm
orbibundles}\label{sec:4}

\setcounter{equation}{0}

As said in Introduction our motivations in this paper come from the
study of Floer homology and Gromov-Witten invariants theory. Our
abstract versions of this section exactly contain applications
needed in these two cases. We first give a general framework and
then discuss two kinds of more concrete ones, Framework I and
Framework II. They are respectively used to construct Gromov-Witten
invariants and  Floer homology.

\subsection{A general framework}\label{sec:4.1}

Let $X$ be a (weakly) PS Banach manifold and ${\rm Aut}(X)$ be the
group of PS automorphisms of $X$. From now on saying {\bf a finite
group $G$ acting on $X$} always means that $G$ acts on $X$ by PS
automorphisms of $X$, i.e. a homomorphism representation $G\to {\rm
Aut}(X)$.

In Definition~\ref{def:2.1}, if $\widetilde U$ is required to be a
connected open subset of a (weakly) PS Banach manifold (with corner)
one gets the notion of a {\bf (weakly) PS Banach orbifold chart} (or
simply {\bf (weakly) PS B-orbifold chart}) on $X$. Such a chart
$(\widetilde U, \Gamma_U,\pi_U)$ is called {\bf effective} if
$\Gamma_U$ acts on  $\widetilde U$ effectively. Replacing the words
``smooth'', ``diffeomorphism'', ``embedding'' and ``bundle open
embedding'') by ``PS'', ``PS diffeomorphism'', ``PS embedding'' and
``PS bundle open embedding'') in the other part of
Section~\ref{sec:2.1} one obtains all corresponding notions in the
PS category, e.x. a {\bf PS injection} from a (weakly) PS B-orbifold
chart $(\widetilde U, \Gamma_U,\pi_U)$ to another one $(\widetilde
V, \Gamma_V,\pi_V)$,  a {\bf (weakly) PS Banach orbifold}, and a
{\bf (weakly) PS Banach orbibundle} etc. In particular, a (weakly)
PS Banach orbifold $(X, [{\mathcal A}])$ is called an {\bf effective
(weakly) PS Banach orbifold} if $[{\mathcal A}]$ contains an atlas
${\mathcal A}$ consisting of effective (weakly) PS Banach orbifold
charts.

  Before
continuing our definitions we discuss stratifications of a PS Banach
orbifold. Let $(\widetilde U, \Gamma_U,\pi_U)$ and $(\widetilde V,
\Gamma_V,\pi_V)$ be two charts on $X$ with $U\subset V$, and
$\theta_{UV}=(\tilde\theta_{UV},\gamma_{UV})$ be an injection from
$(\widetilde U, \Gamma_U,\pi_U)$ to $(\widetilde V,
\Gamma_V,\pi_V)$. Assume that $\widetilde
U=\cup_{a\in\Lambda_u}\widetilde U_a$ and $\widetilde
V=\cup_{a\in\Lambda_v}\widetilde V_a$ are respectively
stratification decompositions of $\widetilde U$ and $\widetilde V$.
Correspondingly we have stratification decompositions
$U=\cup_{a\in\Lambda_u}U_a$ and $V=\cup_{a\in\Lambda_v}V_a$ of $U$
and $V$ respectively. Since $\widetilde\theta_{UV}$ is a PS open
embedding from $\widetilde U$ into $\widetilde V$, it must induce an
injection $\iota_{uv}:\Lambda_u\to\Lambda_v$ preserving order such
that for any $a\in\Lambda_u$ the restriction of
$\widetilde\theta_{UV}$ to the stratum $\widetilde U_a$ is an open
embedding into $\widetilde V_{\iota_{uv}(a)}$. Clearly,
$\cup_{a\in\Lambda_v\setminus\iota_{uv}(\Lambda_u)}V_a$ is disjoint
with $U$. Now we take a collection of PS Banach orbifold charts
$\Omega$ on $X$ such that  supports of all charts in $\Omega$ forms
an open cover of $X$. Note that each chart $(\widetilde U,
\Gamma_U,\pi_U)$ in $\Omega$ gives a partially order set $\Lambda_u$
via the stratification decomposition of $U$,
$U=\cup_{a\in\Lambda_u}U_a$. Consider the disjoint union
$$
\widehat\Lambda=\bigcup_{(\widetilde U,
\Gamma_U,\pi_U)\in\Omega}\Lambda_u,
$$
and introduce an equivalence in it as follows: For $(\widetilde U,
\Gamma_U,\pi_U)$ and $(\widetilde V, \Gamma_V,\pi_V)$ in $\Omega$ we
call $a\in\Lambda_u$ and $b\in\Lambda_v$ {\bf directly equivalence}
if $U_a\cap V_b\ne\emptyset$, and {\bf equivalence} if there exists
a chain of charts in $\Omega$, $(\widetilde U_i,
\Gamma_{U_i},\pi_{U_i})$, $i=0,\cdots,n$ such that $(\widetilde U,
\Gamma_U,\pi_U)=(\widetilde U_0, \Gamma_{U_0},\pi_{U_0})$,
$(\widetilde V, \Gamma_V,\pi_V)=(\widetilde U_n,
\Gamma_{U_n},\pi_{U_n})$ and there exist $a_i\in\Lambda_{u_i}$,
$i=0,\cdots,n$ such that $a_0=a$, $a_n=b$ and $a_i$ is directly
equivalent to $a_{i+1}$ for $i=0,\cdots, n-1$. It is easily checked
that this indeed defines an equivalence relation $\sim$ in
$\widehat\Lambda$ and that the equivalence relation is also
compatible with the partial order of each $\Lambda_u$. Thus the
partial orders in all these $\Lambda_u$ induce a partial order in
$\Lambda=\widehat\Lambda/\sim$ (with a maximal element $0$). It is
not hard to prove that this partial order set $\Lambda$ is unique up
to a bijection preserving order, i.e. independent choices of chart
covers of $X$. For $\alpha\in\Lambda$ let
$$
X_\alpha=\bigcup_{(\widetilde U,
\Gamma_U,\pi_U)\in\Omega}\,\bigcup_{a\in\Lambda_u, a\in\alpha}U_a
$$
Then $X=\cup_{\alpha\in\Lambda}X_\alpha$ is the stratification
decomposition of $X$ indexed by $\Lambda$, and for each chart
$(\widetilde U, \Gamma_U,\pi_U)$ with stratification decomposition
$\widetilde U=\cup_{a\in\Lambda_u}\widetilde U_a$ we have an
injection preserving order, $\iota_u:\Lambda_u\to\Lambda$ to satisfy
$U\cap X_{\iota_u(a)}=U_a$ for any $a\in\Lambda_u$. It is not hard
to check that these $X_\alpha$
satisfy:\\
$\bullet$  Any two $X_\alpha$ and $X_\beta$ with $\alpha\ne\beta$
are disjoint, and each point of $X$ has a neighborhood that may only
intersect with finitely many $X_\alpha$.\\
$\bullet$  Each $X_\alpha$, called a {\bf stratum}, is a smooth
Banach orbifold with corner whose topology assures the inclusion
$X_\alpha\hookrightarrow X$ to be continuous.\\
$\bullet$ The stratum $X_0$,  called the {\bf top stratum} of $X$,
 is a nonempty open subset in $X$. Other strata cannot contain any interior point of $X$,
  and are called {\bf lower strata}.\\
$\bullet$ For each $\alpha\in\Lambda$ the stratum $X_\alpha$ is
relatively open in $\cup_{\beta\prec\alpha}X_\beta$ with respect to
the induced topology from $X$. For any $\beta\in\Lambda$ with
$\beta\prec\alpha$ and $\beta\ne\alpha$, the stratum $X_\beta$
cannot contain any relative interior point of $X_\alpha$ in
$\cup_{\beta\prec\alpha}X_\beta$.

The above arguments show that the stratification decomposition of
$X$,  $X=\cup_{\alpha\in\Lambda} X_\alpha$, is intrinsic. Similarly,
we have a natural stratification decomposition of a weakly PS Banach
orbifold indexed by a finite set which is not necessarily a partial
order set.

Many corresponding terminology with ones in Section~\ref{sec:3} can
be naturally introduced. (They will be omitted and be directly used
when we need because one easily understands them.) In particular, we
point out that a {\bf closed (weakly) PS Banach orbifold dimension
$n$} is a compact (weakly) PS Banach orbifold whose lower strata
have at least codimension two and whose top stratum is a
$n$-dimensional Banach orbifold without corners (including
boundary).  Moreover, a topological subspace $Z$ of a (weakly) PS
Banach orbifold $X$ is called a {\bf (weakly) PS Banach suborbifold}
of $X$ if it is a (weakly) PS Banach orbifold with respect to the
induced PS Banach orbifold structure obtained as follows: For each
(weakly) PS Banach orbifold chart $(\widetilde U, \Gamma_U,\pi_U)$
of $X$ with $U\cap Z\ne\emptyset$ there exists a (weakly) PS Banach
submanifold $\widetilde Z_U\subset\widetilde U$ that is not only
stable under $\Gamma_U$ but also compatible with the inclusion map,
such that the restriction $(\widetilde Z_U, \Gamma_U|_{\widetilde
Z_U},\pi_U|_{\widetilde Z_U})$ is a (weakly) PS Banach orbifold
chart for $Z$.

Let $p:X\to E$ be a (weakly) PS Banach orbibundle whose definition
 is obtained by replacing the words ``Banach orbifold'', ``Banach
 orbifold chart'' and ``smooth'' in Definition~\ref{def:2.7} with
the words ``PS Banach orbifold'', ``PS Banach
 orbifold chart'' and ``PS'' respectively.   For a PS
B-orbifold chart $(\widetilde W,\Gamma_W,\pi_W)$ on $X$ let
$({\widetilde E}_W, \Gamma_W,\Pi_W)$ be the corresponding PS Banach
orbibundle chart on $E$ with it, and $\tilde p_W:{\widetilde E}_W\to
\widetilde W$ be the PS Banach bundle projection with
$\pi_W\circ{\widetilde p}_W=p\circ\Pi_W$ on $E_W$. If $\widetilde
W=\cup_{a\in\Lambda_w}\widetilde W_a$ is the stratification
decomposition of $\widetilde W$, then ${\widetilde
E}_W=\cup_{a\in\Lambda_w}(\widetilde p_W)^{-1}(\widetilde W_a)$ is
the stratification decomposition of $\widetilde E_W$. It follows
that $E=\cup_{\alpha\in\Lambda}p^{-1}(X_\alpha)$ is exactly the
stratification decomposition of $E$ which is obtained as the above
construction of the stratification decomposition of $X$. Denote by
$E_\alpha=p^{-1}(X_\alpha)$ and by $p_\alpha$ the restriction of $p$
to $E_\alpha$. Then each $p_\alpha:E_\alpha\to X_\alpha$ is a Banach
orbibundle.
  Having these we have essentially
same two methods to define a {\bf weakly Fredholm PS section} and a
{\bf weakly $k$-Fredholm PS section}. The first is to replace the
words ``PS bundle'' by ``PS orbibundle'' in
Definition~\ref{def:3.5}. To describe the second note that the
restriction of each continuous section $S:X\to E$ to
$W=\pi_W(\widetilde W)$ has a unique $\Gamma_W$-equivariant lift
${\widetilde S}_W:\widetilde W\to {\widetilde E}_W$, called a {\bf
local lift} of $S$. A continuous section $S:X\to E$ is called {\bf
weakly Fredholm PS section} (resp. {\bf weakly $k$-Fredholm PS
section}, resp. {\bf $k$-Fredholm PS section})  if each local lift
${\widetilde S}_W:\widetilde W\to {\widetilde E}_W$ as above is such
a section in the sense of Definition~\ref{def:3.5}. In particular,
it restricts to a Fredholm section $S_\alpha: X_\alpha\to
E|_{X_\alpha}$ on each stratum $X_\alpha$ of $X$. The {\bf index} of
$S$ is defined as that of the top stratum restriction
$S|_{X_0}:X_0\to E|_{X_0}$. Correspondingly, the notions of the {\bf
weakly PS Banach Fredholm orbibundle} of index $r$, {\bf weakly PS
Banach $k$-Fredholm orbibundle} and {\bf PS Banach $k$-Fredholm
orbibundle} may be defined naturally. For example, a triple $(X,
E,S)$ consisting of a PS Banach orbibundle $p:E\to X$ and its weakly
Fredholm section $S:X\to E$ (of index $r$) is called a {\bf weakly
PS Banach Fredholm orbibundle} (of index $r$). The Banach Fredholm
bundle $(\widetilde W, {\widetilde E}_W, {\widetilde S}_W)$ is
called a {\bf local lift} of $(X, E, S)$ over $W$.  A weakly PS
Banach Fredholm orbibundle  (resp. (weakly) PS Banach $k$-Fredholm
orbibundle) $(X, E, S)$ is called {\bf oriented} if its top stratum
$(X, E, S)_0=(X_0, E|_{X_0}, S|_{X_0})$, as a Banach Fredholm
orbibundle, is oriented.

Since there is no  corresponding slice theorem for action of Lie
groups on PS Banach manifolds, a {\bf $k$-regular PS Banach
orbifold} is defined as a PS Banach orbifold $X$ on which each chart
$(\widetilde W,  \Gamma_W,\pi_W)$ of it is {\bf $k$-regular} in the
sense that {\it the fixed point set ${\rm Fix}(g)$ of each
$g\in\Gamma_U\setminus{\rm Ker}(\Gamma_W,\widetilde W)$ and the top
stratum $\widetilde W_0$ of $\widetilde W$ intersects at the union
of finitely many submanifolds of $\widetilde W_0$ of codimension at
least $k$ locally, i.e. ${\rm Fix}(g)\cap \widetilde W_0$ is locally
 the union of finitely many submanifolds of $\widetilde W_0$ of codimension at
least $k$.} {\bf Remark} that in Definition~\ref{def:3.1} of the PS
Banach manifold we do not require that the top stratum is dense. So
$\widetilde W_0$ and thus the regular point set of $\widetilde W_0$,
$\widetilde W_0^\circ$ is not necessarily dense in $\widetilde W$.
(The dense condition is needed if we make our construction in the
category of effective PS Banach Fredholm orbibundles).

In the following, for simplicity, as in
Sections~\ref{sec:2.3}-\ref{sec:2.8} we {\it always assume}:\\
$\clubsuit 1$  {\bf all PS Banach orbifolds are effective and
$1$-regular,}\\
$\clubsuit 2$ {\bf all PS Banach orbifolds (or orbibundles) have
dense
top strata.}\\
 (Of course, the second condition will be
moved out if our arguments are not in the effective category, but in
arbitrarily one.)

 Then $\widetilde W_0$ and thus the regular point set of $\widetilde W_0$,
$\widetilde W_0^\circ$, is open and dense in $\widetilde W$. It is
sufficient for proving the action of $\Gamma_I$ on $\widetilde W_I$
to be effective as in Claim~\ref{cl:2.16}.

 For a compact subset $K\subset X$, one only need replace
the words ``Banach'', ``diffeomorphism'' and ``smooth'' by ``PS
Banach'', ``PS diffeomorphism'' and ``PS'' respectively in the
arguments in \S~\ref{sec:2.2} and \S~\ref{sec:2.3}, and then repeats
word by word the constructions therein to get a system of virtual PS
Banach bundles
$$
\bigl(\widetilde{\mathcal E}(K), \widetilde W(K)\bigr)=
\bigl\{(\widetilde E_I, \widetilde W_I), \pi_I, \Pi_I, \Gamma_I,
\pi^I_J, \Pi^I_J,\,\lambda^I_J\,\bigm|\, J\subset I\in{\mathcal
N}\bigr\}
$$
corresponding with (\ref{e:2.15}), and its resolution system of PS
Banach bundles
$$
 \bigl(\widehat{\mathcal E}(K),
\widehat W(K)\bigr)= \bigl\{\bigl(\widehat E_I, \widehat
W_I\bigr), \hat\pi_I, \hat\Pi_I, \Gamma_I,
\hat\pi^I_J,\hat\Pi^I_J, \lambda^I_J\,\bigm|\, J\subset
I\in{\mathcal N}\bigr\},
$$
Then take the pairs of open sets $W^j_i\subset\subset U^j_i$,
$j=1,\cdots, n-1$ such that $U^1_i\subset W^2_i\subset\subset
U^2_i\subset\subset\cdots\subset\subset W^{n-1}_i\subset\subset
U^{n-1}_i\subset\subset W^n_i=W_i$, $i=1,\cdots, n$, and that
$K\subset\cup^n_{i=1}W^1_i$. As in \S 2.3.2 we can use them to get a
renormalization of $(\widehat{\mathcal E}(K), \widehat W(K))$,
another system of PS Banach bundles
$$
 \bigl(\widehat{\mathcal F}(K),
\widehat V(K)\bigr)= \bigl\{\bigl(\widehat F_I, \widehat
V_I\bigr), \hat\pi_I, \hat\Pi_I, \Gamma_I,
\hat\pi^I_J,\hat\Pi^I_J, \lambda^I_J\,\bigm|\, J\subset
I\in{\mathcal N}\bigr\}.
$$
 From the local expressions in
(\ref{e:2.39}) and (\ref{e:2.52}) one easily sees that both maps
$\hat\pi^I_J$ and $\hat\Pi^I_J$ are {\bf partially smooth} in the
present case. As in Lemma~\ref{lem:2.30} we can prove that every PS
section $S$ of $E\to X$ may yield a global PS sections of the above
PS Banach bundle systems,  denoted by $\widehat S=\{\widehat
S_I\,|\, I\in{\mathcal N}\}$. The section $\widehat S$ is
equivariant in the sense that each $\widehat S_I$ is
$\Gamma_I$-equivariant. If $S$ is weakly Fredholm (resp. (weakly)
$k$-Fredholm), so is $\widehat S$, and ${\rm Ind}(S)={\rm
Index}(\widehat S)$.

Note that our PS Banach orbifolds might not have PS cut-off
functions even if each stratum (as a Banach orbifold) admits smooth
cut-off functions. From now on we make the following:
\begin{assumption}\label{ass:4.1}
{\rm Each $x\in X$ admits PS cut-off functions with arbitrarily
small supports.}\footnote{It can be satisfied in the cases of
Floer homology and GW-invariants, see Appendix~\ref{app:D}.}
\end{assumption}
So we easily obtain:

\begin{claim}\label{cl:4.2}
Under Assumption~\ref{ass:4.1}, for any PS B-orbifold chart
$(\widetilde W,\Gamma_W,\pi_W)$ on $X$, any point $\tilde
x\in\widetilde W$ and any neighborhood $\widetilde O$ of it there
exists a PS cut-off function $\delta$ such that ${\rm
supp}(\delta)\subset\widetilde O$ and $\delta\equiv 1$ near
$\tilde x$.
\end{claim}

A weakly PS Banach Fredholm orbibundle $(X, E, S)$ is called to have
{\bf rich local PS sections of class ${\cal A}$} if each local lift
$(\widetilde W, {\widetilde E}_W, {\widetilde S}_W)$ of it, as a
weakly PS Banach Fredholm bundle, has rich PS sections of class
${\cal A}(W)$.

Let $(X, E, S)$  have  rich local PS sections of class ${\cal A}$.
If $(\widetilde W, {\widetilde E}_W, {\widetilde S}_W)$ and
$(\widetilde V, {\widetilde E}_V, {\widetilde S}_V)$ are any two
local lifts of it with $W\cap V\ne\emptyset$, then for any $x\in
W\cap V$, $\tilde x_w\in (\pi_W)^{-1}(x)$ and $\tilde x_v\in
(\pi_V)^{-1}(x)$, as in Lemma~\ref{lem:2.19} it follows from
Lemma~\ref{lem:2.4} that there exist small open connected
neighborhoods $\widetilde O(\tilde x_w)$ of $\tilde x_w$ and
$\widetilde O(\tilde x_v)$ of $\tilde x_v$, a PS diffeomorphism
$\lambda_{WV}:\widetilde O(\tilde x_w)\to \widetilde O(\tilde x_v)$
that maps $\tilde x_w$ to $\tilde x_v$, and a group isomorphism
${\cal A}_{WV}:\Gamma_W(\tilde x_w)\to \Gamma_V(\tilde x_v)$ such
that
$$
 \lambda_{WV}\circ\phi={\cal
A}_{WV}(\phi)\circ\lambda_{WV}
$$
for any $\phi\in\Gamma_W(\tilde x_w)$. Here $\Gamma_W(\tilde x_w)$
is the isotropy group of $\Gamma_W$ at $\tilde x_w$. Note also that
the PS diffeomorphism $\lambda_{WV}$ is unique up to composition
with elements in $\Gamma_W(\tilde x_w)$ and $\Gamma_V(\tilde x_v)$.
Shrinking $\widetilde O(\tilde x_w)$ and $\widetilde O(\tilde x_v)$
if necessary it follows from the properties of PS orbifold bundles
that there exists a PS bundle isomorphism
$$
 \Lambda_{WV}:{\widetilde E}_W|_{{\widetilde
O}(\tilde x_w)}\to{\widetilde E}_V|_{{\widetilde O}(\tilde x_v)}
$$
 which is a lifting of $\lambda_{WV}$, to satisfy
$$
\Lambda_{WV}\circ\Phi=\boxed{{\mathcal
A}_{WV}(\phi)}\circ\Lambda_{WV}
$$
 for any $\phi\in\Gamma_W(\tilde x_w)$. Here
$$
 \Phi: {\widetilde
E}_W|_{{\widetilde O}(\tilde x_w)}\to{\widetilde E}_W|_{{\widetilde
O}(\tilde x_w)}\; \bigl({\rm resp.}\; \boxed{{\cal
 A}_{WV}(\phi)}: {\widetilde
E}_V|_{{\widetilde O}(\tilde x_v)}\to{\widetilde E}_V|_{{\widetilde
O}(\tilde x_v)}\bigr)
$$
is the PS bundle isomorphism lifting of $\phi$ (resp. ${\cal
A}_{WV}(\phi)$) produced in the definition of PS orbifold bundles.

Now let $\tilde\sigma_v:\widetilde V\to {\widetilde E}_V$ be a PS
section of class ${\cal A}(V)$. Then
$$
 \lambda_{WV}^\ast(\sigma):
{\widetilde O}(\tilde x_w)\to {\widetilde E}_W|_{{\widetilde
O}(\tilde x_w)},\; \tilde y\mapsto
\Lambda_{WV}^{-1}\bigl(\sigma_v(\lambda_{WV}(\tilde y))\bigr)
$$
 is also a PS section. By Claim~\ref{cl:4.2} we may take a PS cut-off function
$\delta: \widetilde W\to [0,1]$ such that ${\rm
supp}(\delta)\subset{\widetilde O}(\tilde x_w)$ and $\delta\equiv
1$ near $\tilde x_w$. Then
\begin{equation}\label{e:4.1}
\delta\cdot\lambda_{WV}^\ast(\sigma): {\widetilde W}\to
{\widetilde E}_W,
\end{equation}
is a PS section that equals $\lambda_{WV}^\ast(\sigma)$ near
$\tilde x_w$ and has support contained in ${\widetilde O}(\tilde
x_w)$. Adding all such PS sections of $\widetilde E_W$ to the
class ${\cal A}(W)$ we get a new class of PS sections of
$\widetilde E_W$, denoted by
\begin{equation}\label{e:4.2}
\widetilde{\cal A}(W).
\end{equation}
Corresponding with Definition~\ref{def:3.13} we introduce:

\begin{definition}\label{def:4.3}
{\rm  Let $(X,E,S)$ be a weakly PS Banach Fredholm orbibundle (resp.
a (weakly) PS Banach $k$-Fredholm orbibundle) with rich local PS
sections of class ${\cal A}$. It is called {\bf good} (resp. {\bf
(weakly) $k$-good relative to class ${\cal A}$})
  if the PS Banach orbifold $X$ satisfies
 Assumption~\ref{ass:4.1} and  each local lift $(\widetilde W, {\widetilde E}_W,
{\widetilde S}_W)$ of $(X, E, S)$ is  good (resp. (weakly) $k$-good
relative to class $\widetilde{\cal A}(W)$ in the sense of
Definition~\ref{def:3.13}, i.e. the following two conditions  are
also satisfied:
\begin{description}
\item[(i)] For each local lift $(\widetilde W, {\widetilde E}_W,
{\widetilde S}_W)$ of $(X, E, S)$ the PS section ${\widetilde S}_W$
is quasi (resp. (weakly) $k$-quasi) transversal to the zero section
relative to class $\widetilde{\cal A}(W)$.

\item[(ii)] For each local lift $(\widetilde W, {\widetilde E}_W,
{\widetilde S}_W)$ of $(X, E, S)$ the PS section ${\widetilde
S}_W$ is locally proper relative to class $\widetilde{\cal A}(W)$
near $Z({\widetilde S}_W)$.
\end{description}
Moreover, $(X,E,S)$  is called {\bf smoothly good} (resp. {\bf
smoothly (weakly) $k$-good ) relative to class ${\cal A}$} if  the
above condition
 (i) is replaced by\\
${\bf (i)'}$ For each local lift $(\widetilde W, {\widetilde E}_W,
{\widetilde S}_W)$ of $(X, E, S)$ the PS section ${\widetilde
S}_W$ is smoothly quasi-transversal (resp. smoothly (weakly)
$k$-quasi-transversal)  to the zero section relative to class
$\widetilde{\cal A}(W)$. }
\end{definition}

In the applications to Gromov-Witten invariants and Floer theory,
 the condition (i) or ${\bf (i)'}$ in Definition~\ref{def:4.3} shall be
checked by complicated gluing arguments.

{\bf From now on} we assume that {\bf $S$ is a (weakly) $k$-Fredholm
PS section of index $r$ and with compact zero locus $Z(S)$}. As
before we always write
$$
\bigl(\widehat{\mathcal E}, \widehat W\bigr)
=\bigl(\widehat{\mathcal E}(Z(S)), \widehat W(Z(S))\bigr),\quad
\bigl(\widehat{\mathcal F}, \widehat V\bigr)
=\bigl(\widehat{\mathcal F}(Z(S)), \widehat V(Z(S))\bigr).
$$
 As in Remark~\ref{rem:2.34}, for each
$i=1,\cdots,n$ let us take pairs of open subsets
$W^{+j}_i\subset\subset U^{+j}_i$, $j=1,\cdots, n-1$ such that
$$
W^j_i\subset\subset W^{+j}_i\subset\subset U^{+j}_i\subset\subset
U^j_i,\;j=1,\cdots, n-1.
$$
Then one gets another renormalization  of the PS bundle system
$(\widehat {\mathcal E}, \widehat W)$,
$$
  (\widehat {\mathcal F}^+,
\widehat V^+)= \bigl\{\bigl(\widehat F^+_I, \widehat V^+_I\bigr),
\hat\pi^+_I, \hat\Pi^+_I, \Gamma_I, \hat\pi^{+I}_J,\hat\Pi^{+I}_J,
\lambda^I_J\,\bigm|\, J\subset I\in{\mathcal N}\bigr\},
$$
such that $Cl(V_I)\subset V_I^+$ and thus $Cl(\widehat V_I)\subset
\widehat V_I^+$ for each $I\in{\cal N}$ and that $\hat\pi^+_I$,
$\hat\Pi^+_I$,  $\hat\pi^{+I}_J$ and $\hat\Pi^{+I}_J$ respectively
restrict to $\hat\pi_I$, $\hat\Pi_I$, $\hat\pi^{I}_J$ and
$\hat\Pi^{I}_J$.

Take the $\Gamma_i$-invariant PS cut-off functions\footnote{Under
Assumption~\ref{ass:4.1} one might suspect the existence of such PS
cut-off functions. If necessary, we can always shrink $W_i$ and make
suitable adjustments for completing arguments.} $\gamma_i:\widetilde
W_i\to [0, 1]$ with support in $\widetilde W_i^1$, $i=1,\cdots, n$
such that
\begin{equation}\label{e:4.3}
Z(S)\subset\cup^n_{i=1}U^0_i,
\end{equation}
where $U_i^0=\pi_i(\widetilde U_i^0)$ and  $\widetilde
U_i^0=\{\tilde x\in\widetilde W_i\,|\,\gamma_i(\tilde x)>0\}$. Let
$\widetilde S_i:\widetilde W_i\to\widetilde E_i$ be the unique
$\Gamma_i$-equivariant lifts of $S|_{W_i}$, $i=1,\cdots, n$. Since
$Cl(U^0_i)\subset W^1_i\subset\subset W_i$ and $Z(\widetilde
S_i)\cap Cl(\widetilde W^1_i)\subset\widetilde W_i$ are compact,
 and as in \S~\ref{sec:3} we can find finitely many
PS sections of the PS Banach bundles $\widetilde E_i\to\widetilde
W_i$ of class ${\cal A}(W_i)$, $\tilde s_{ij}$, $j=1,\cdots, m_i$
such that
$$
D\widetilde S_i(\tilde z)(T_{\tilde z}\widetilde W_i)+ {\rm
span}(\{\tilde s_{i1}(\tilde z),\cdots, \tilde s_{im_i}(\tilde
z)\})=(\widetilde E_i)_{\tilde z}
$$
  for any $\tilde z\in Z(\widetilde S_i)\cap Cl(\widetilde W^1_i)$,
  $i=1,\cdots,n$.
({\bf Note}: the second condition in Definition~\ref{def:3.8} is
used again!). Set  $\tilde\sigma_{ij}=\tilde\gamma_i\cdot\tilde
s_{ij}$, $j=1,\cdots, m_i$. Then
\begin{equation}\label{e:4.4}
D\widetilde S_i(\tilde y)+ {\rm span}(\{\tilde\sigma_{i1}(\tilde
y),\cdots,\tilde\sigma_{im_i}(\tilde y)\})=(\widetilde
E_i)_{\tilde y}
\end{equation}
for any $\tilde y\in Z(\widetilde S_i)\cap U^0_i$ and $i=1,\cdots,
n$. By Lemma~\ref{lem:2.33} each $\tilde\sigma_{ij}$ yields a global
PS section $\hat\sigma_{ij}=\{(\hat\sigma_{ij})_I\,|\, I\in{\cal
N}\}$ of the system of PS Banach bundles $(\widehat {\mathcal F},
\widehat V)$.

As in (\ref{e:2.82}) and (\ref{e:2.83}) let $m=m_1+\cdots+m_n$ and
consider the PS Banach bundle system
\begin{eqnarray}
&&\bigl({\bf P}_1^\ast\widehat {\mathcal F}, \widehat
V\times\R^m\bigr)\label{e:4.5}\\
&&=\bigl \{({\bf P}_1^\ast\widehat F_I, \widehat V_I\times\R^m),
\hat\pi_I,  \hat\pi^I_J, \hat\Pi_I, \hat\Pi^I_J, \hat
p_I,\Gamma_I\,\bigm|\, J\subset I\in{\mathcal N}\bigr\}, \nonumber
\end{eqnarray}
and its (weakly) $k$-Fredholm PS section $\Upsilon=\{\Upsilon_I\,|\,
I\in{\cal N}\}$,
\begin{eqnarray}
&&\Upsilon_I: \widehat V_I\times\R^m\to {\bf P}_1^\ast\widehat
F_I,\label{e:4.6}\\
&&\qquad \Bigl(\hat x_I, \{t_{ij}\}_{\substack{
 1\le j\le m_i\\
 1\le i\le n}}\Bigr)\mapsto \hat S_I(\hat x_I)+
 \sum^n_{i=1}\sum^{m_i}_{j=1}t_{ij}(\hat\sigma_{ij})_I(\hat
 x_I)\nonumber\\
&&\hspace{35mm} =\hat S_I(\hat x_I)+
 \sum_{i\in I}\sum^{m_i}_{j=1}t_{ij}(\hat\sigma_{ij})_I(\hat
 x_I).\nonumber
 \end{eqnarray}

 Then $Cl(\widehat V_I)\subset \widehat
V^+_I$ for each $I\in{\cal N}$, and $\Upsilon=\{\Upsilon_I\,|\,
I\in{\cal N}\}$ is naturally viewed as a PS section of $({\bf
P}_1^\ast\widehat {\mathcal F}^+, \widehat V^+\times\R^m)$.  Let
$\hat u_I\in Cl(\widehat V_I)\cap Z(\widehat S_I)$. Then
$\Upsilon_I(\hat u_I, 0)=0$ and it follows from (\ref{e:4.3}) that
$u_I=\hat\pi_I(\hat u_I)\in\cup^n_{i=1}U^0_i$. So $u_I$ sits in
$U^0_{i_q}\subset W_{i_q}$ for some $i_q\in I$ because of
(\ref{e:2.69}). As in (\ref{e:2.89}), in the natural pullback charts
of those in (\ref{e:2.87}) and (\ref{e:2.88}) we have the following
local expression of $\Upsilon_I$,
\begin{eqnarray}
&&\Upsilon_{Iq}:\widetilde O_{i_q}\times\R^m\to{\bf
P}_1^\ast(\widetilde E_{i_q}|_{\widetilde
O_{i_q}}),\label{e:4.7}\\
&&\Bigl(\tilde x, \{t_{ij}\}_{\substack{
 1\le j\le m_i\\
 1\le i\le n}}\Bigr) \mapsto\tilde
S_{i_q}(\tilde x) + \sum_{i\in
I}\sum^{m_i}_{j=1}t_{ij}(\Phi_i\circ\Lambda_{i_qi})^{-1}
\bigl(\tilde\sigma_{ij}(\phi_i\circ\lambda_{i_qi}(\tilde
x))\bigr)\nonumber\\
&&\hspace{28mm} =\tilde S_{i_q} + \sum^{m_{i_q}}_{j=1}t_{i_qj}
\tilde\sigma_{i_qj}+ \sum_{i\in
I\setminus\{i_q\}}\sum^{m_i}_{j=1}t_{ij}\tilde\tau^q_{ij}.\nonumber
\end{eqnarray}
Here $\tilde\tau^q_{ij}(\tilde x)
=(\Phi_i\circ\Lambda_{i_qi})^{-1}
\bigl(\tilde\sigma_{ij}(\phi_i\circ\lambda_{i_qi}(\tilde
x))\bigr)$. Note that multiplying $\tilde\tau^q_{ij}$ by a PS
cut-off function near $\tilde u_{i_q}\in\widetilde O_{i_q}$ gives
a PS section of the PS Banach bundle $\widetilde
E_{i_q}\to\widetilde W_{i_q}$ of class $\widetilde{\cal
A}(W_{i_q})$. So by Definition~\ref{def:4.3} we get a small open
neighborhood $O(\tilde u_{i_q})$ of $\tilde u_{i_q}$ in
$\widetilde O_{i_q}$ and $\varepsilon>0$ such that:
\begin{description}
\item[$(\ref{e:4.7})_1$]  The restriction of the PS section
$\Upsilon_{Iq}$ to $O(\tilde u_{i_q})\times B_\varepsilon(\R^m)$
  is a (weakly) $k$-Fredholm PS section that is transversal to the zero  section.

\item[$(\ref{e:4.7})_2$] The closure of the set $(O(\tilde
u_{i_q})\times B_\varepsilon(\R^m))\cap\Upsilon_{Iq}^{-1}(0)$  in
$Cl(O(\tilde u_{i_q}))\times\R^m\subset\widetilde
W_{i_q}\times\R^m$ is compact.
\end{description}

 It follows that there exists a neighborhood $U(\hat u_I)$ of $\hat u_I$ in
 $\widehat V_I^+$ such that:
\begin{description}
\item[$(\ref{e:4.7})_3$]  The restriction of the PS section
$\Upsilon_I$ to $U(\hat u_I)\times B_\varepsilon(\R^m)$
  is a (weakly) $k$-Fredholm PS section that is transversal to the zero  section.

   \item[$(\ref{e:4.7})_4$]  The closure of the
   set $(U(\hat u_I)\times B_\varepsilon(\R^m))\cap\Upsilon_I^{-1}(0)$
  in $Cl(U(\hat u_I))\times\R^m\subset\widehat V_I^+\times\R^m$ is compact.
\end{description}

As before using the compactness of $Cl(\widehat V_I)\cap
Z(\widehat S_I)\subset\widehat V_I^+$ we can find a small
neighborhood ${\cal U}_I$ of $Cl(\widehat V_I)\cap Z(\widehat
S_I)$ in $\widehat V_I^+$ such that:
\begin{description}
 \item[$(\ref{e:4.7})_5$]  The
restriction of the PS section $\Upsilon_I$ to ${\cal U}_I\times
B_\varepsilon(\R^m)$
  is a (weakly) $k$-Fredholm PS section that is transversal to the zero  section.

 \item[$(\ref{e:4.7})_6$]  The set $({\cal U}_I\times
B_\varepsilon(\R^m))\cap\Upsilon_I^{-1}(0)$ has a compact closure
in   $Cl({\cal U}_I)\times\R^m\subset\widehat V_I^+\times\R^m$.
\end{description}

Take an open neighborhood ${\cal W}^\ast$ of $Z(S)$ such that
$$
{\cal W}^\ast\subset (\cup_{I\in{\cal N}}{\cal
U}_I)\cap(\cup^n_{i=1}U_i^0).
$$
Then as before we have a system of PS Banach bundles
$$
 (\widehat {\mathcal
F}^\ast, \widehat V^\ast)=\bigl \{(\widehat F_I^\ast, \widehat
V_I^\ast), \hat\pi_I,
 \hat\pi^I_J, \hat\Pi_I, \hat\Pi^I_J, \hat p_I,\Gamma_I\,\bigm|\, J\subset I\in{\mathcal N}\bigr\},
$$
 called the {\bf restriction} of $(\widehat {\mathcal F},
\widehat V)$ to the open subset ${\cal W}^\ast$. Clearly,
$\widehat V^\ast_I\subset\widehat V_I$ and thus $Cl(\widehat
V^\ast_I)\subset\widehat V^+_I$ for any $I\in{\cal N}$. Note that
each global section $\hat \sigma =\{(\hat\sigma)_I\,|\,I\in{\cal
N}\}$ of $(\widehat{\cal F},\widehat V)$ restricts to a global
section of $(\widehat {\mathcal F}^\ast, \widehat V^\ast)$, still
denoted by $\hat\sigma$ without confusions.  Since ${\cal N}$ is
finite by shrinking $\varepsilon>0$ and ${\cal W}^\ast$ we get:

\begin{theorem}\label{th:4.4}
There exists a small open neighborhood ${\cal W}^\ast$ of $Z(S)$ in
$X$ and $\varepsilon>0$ such that for the restriction bundle system
$(\widehat{\cal F}^\ast,\widehat V^\ast)$ of $(\widehat{\cal
F},\widehat V)$ to ${\cal W}^\ast$ the restriction of the global
section $\Upsilon=\bigl\{\Upsilon_I\,|\, I\in{\cal N}\bigr\}$ to
$({\bf P}_1^\ast\widehat {\mathcal F}^\ast, \widehat V^\ast\times
B_\varepsilon(\R^m))$  is a (weakly) $k$-Fredholm PS section that is
transversal to the zero section. Consequently, for each $I\in{\cal
N}$ the set
$$
\widehat\Omega_I(S):=\{(\hat x_I, {\bf t})\in\widehat
V_I^\ast\times B_\varepsilon(\R^m)\,|\, \Upsilon_I(\hat x_I,{\bf
t})=0\}
$$
is a weakly $k$-PS submanifold of dimension ${\rm Ind}(S)+m$ and has
compact closure in $Cl(\widehat V^\ast_I)\times\R^m\subset\widehat
V_I^{+}\times\R^m$. They are compatible in the sense that for any
$J\subset I\in{\cal N}$,
$$
\hat\pi^I_J\bigl((\hat \pi^I_J)^{-1}(\widehat V_J^\ast\times
B_\varepsilon(\R^m))\cap\widehat\Omega_I(S)\bigr)={\rm
Im}(\hat\pi^I_J)\cap\widehat\Omega_J(S)\quad{\rm and}
$$
that for any $(\hat x_J, {\bf t})\in{\rm
Im}(\hat\pi^I_J)\cap\widehat\Omega_J(S)$ all
$$\frac{|\Gamma_I|-|\Gamma_I({\hat
x_I})|+1}{|\Gamma_J|-|\Gamma_J({\hat x_J})|+1}=\prod_{i\in
I\setminus J}(|\Gamma_i|-|\Gamma_J({\hat x_J})|^{1/|J|}+1)$$
elements of the inverse image $(\hat\pi^I_J)^{-1}((\hat x_J, {\bf
t}))$ sit in $(\hat \pi^I_J)^{-1}(\widehat V_J^\ast\times
B_\varepsilon(\R^m))\cap\widehat\Omega_I(S)$
 for any $(\hat x_I,{\bf t})\in(\hat\pi^I_J)^{-1}((\hat x_J, {\bf t}))$.
So for
\begin{eqnarray*}
&&\widehat W_J^\circ:=\{\hat x_I\in \widehat W_J(S)\,|\,
\Gamma_I(\hat x_I)=\{1\}\}\quad{\rm and}\\
&&\widehat\Omega_J(S)^\circ:=\widehat\Omega_J(S)\cap(\widehat
W_J^\circ\times B_\varepsilon(\R^m)),
\end{eqnarray*}
the restriction of the projection
$$
\hat\pi^I_J:(\hat \pi^I_J)^{-1}(\widehat V_J^\ast\times
B_\varepsilon(\R^m))\cap\widehat\Omega_I(S)\to {\rm
Im}(\hat\pi^I_J)\cap\widehat\Omega_J(S)
$$
to $(\hat
\pi^I_J)^{-1}(\widehat V_J^\ast\times
B_\varepsilon(\R^m))\cap\widehat\Omega_I(S)^\circ$
 is a $|\Gamma_I|/|\Gamma_J|$-fold (regular) PS covering to
 ${\rm Im}(\hat\pi^I_J)\cap\widehat\Omega_J(S)^\circ$. Moreover,
if $(X,E,S)$ is smoothly good (resp. smoothly (weakly) $k$-good)
then the above each $\widehat\Omega_I(S)$ is a smooth manifold of
dimension ${\rm Ind}(S)+m$ (though it is not necessarily equal to
$\widehat\Omega_I(S)^\circ$).
\end{theorem}

The remarks and conventions below Theorem~\ref{th:3.14} are also
effective in this section.

{\bf From now on} we also assume that $X$ is {\bf separable}, i.e.
each stratum of it is separable.  Since ${\cal N}$ is a finite set
 Theorem~\ref{th:4.4} and the Sard-Smale theorem immediately
give:

\begin{corollary}\label{cor:4.5}
There exists a residual subset $B_\varepsilon(\R^m)_{res}\subset
B_\varepsilon(\R^m)$ such that for each ${\bf t}\in
B_\varepsilon(\R^m)_{res}$ the global section
 $\Upsilon^{({\bf t})}=\{\Upsilon^{({\bf t})}_I\,|\, I\in{\cal
N}\}$ of the PS Banach bundle system $\bigl(\widehat {\mathcal
F}^\ast, \widehat V^\ast\bigr)$ is  a (weakly) $k$-Fredholm PS
section that is transversal to the zero section and  the set
$\widehat{\cal M}^{\bf t}_I(S):= (\Upsilon^{({\bf t})}_I)^{-1}(0)$
is a weakly $k$-PS submanifold of dimension ${\rm Ind}(S)$ and with
compact closure in $Cl(\widehat V^\ast_I)\subset\widehat
V_I^{+}\subset\widehat W_I^{+}$.
Here ${\bf t}=\{t_{ij}\}_{\substack{ 1\le j\le m_i\\
 1\le i\le n}}$ and
$$
 \Upsilon^{({\bf t})}_I: \widehat
V^\ast_I\to \widehat F^\ast_I,\; \hat x_I \mapsto
\Upsilon_{I}(\tilde x, {\bf t})=\hat S_I(\hat x_I)+
 \sum^n_{i=1}\sum^{m_i}_{j=1}t_{ij}(\hat\sigma_{ij})_I(\hat x_I).
 $$
Moreover, the family
\begin{equation}\label{e:4.8}
\widehat{\cal M}^{\bf t}(S)=\bigl\{\widehat{\cal M}^{\bf
t}_I(S)\,:\, I\in{\cal N}\bigr\}
\end{equation}
is compatible in the sense that for any $J\subset I\in{\cal N}$,
$$
 \hat\pi^I_J\bigl((\hat
\pi^I_J)^{-1}(\widehat V_J^\ast)\cap\widehat{\cal M}^{\bf
t}_I(S)\bigr)={\rm Im}(\hat\pi^I_J)\cap\widehat{\cal M}^{\bf
t}_J(S)\quad{\rm and}
$$
that for any $\hat x_J\in{\rm Im}(\hat\pi^I_J)\cap\widehat{\cal
M}^{\bf t}_J(S)$ the inverse image $(\hat\pi^I_J)^{-1}(\hat x_J)$
contains
$$\frac{|\Gamma_I|-|\Gamma_I({\hat
x_I})|+1}{|\Gamma_J|-|\Gamma_J({\hat x_J})|+1}=\prod_{i\in
I\setminus J}(|\Gamma_i|-|\Gamma_J({\hat x_J})|^{1/|J|}+1)$$
elements and all sit in $(\hat \pi^I_J)^{-1}(\widehat
V_J^\ast)\cap\widehat{\cal M}^{\bf t}_I(S)$
 for any $\hat x_I\in(\hat\pi^I_J)^{-1}(\hat x_J)$.
So for each stratum $\widehat V_{J\alpha}^{\ast}$ of $\widehat
V_{J}^{\ast}$ and
\begin{equation}\label{e:4.9}
\left.\begin{array}{ll}
 \widehat{\cal M}^{\bf t}_J(S)_\alpha:=\widehat{\cal M}^{\bf t}_J(S)\cap\widehat
V_{J\alpha}^{\ast},\\
 \widehat{\cal M}^{\bf t}_J(S)^\circ:=\widehat{\cal M}^{\bf
t}_J(S)\cap\widehat V_J^{\ast\circ}\quad{\rm and}\\
\widehat{\cal M}^{\bf t}_J(S)^{sing}:=\widehat{\cal M}^{\bf
t}_J(S)\cap\widehat V_J^{\ast sing}
\end{array}\right\}
\end{equation}
it also holds that
\begin{equation}\label{e:4.10}
\left.\begin{array}{ll} \hat\pi^I_J\bigl((\hat
\pi^I_J)^{-1}(\widehat V_J^\ast)\cap\partial^k\widehat{\cal
M}^{\bf t}_I(S)_\alpha\bigr)={\rm
Im}(\hat\pi^I_J)\cap\partial^k\widehat{\cal M}^{\bf
t}_J(S)_\alpha,\\
 \hat\pi^I_J\bigl((\hat
\pi^I_J)^{-1}(\widehat V_J^\ast)\cap\widehat{\cal M}^{\bf
t}_I(S)^\circ\bigr)={\rm
Im}(\hat\pi^I_J)\cap\widehat{\cal M}^{\bf t}_J(S)^\circ,\\
\hat\pi^I_J\bigl((\hat \pi^I_J)^{-1}(\widehat
V_J^\ast)\cap\widehat{\cal M}^{\bf t}_I(S)^{sing}\bigr)={\rm
Im}(\hat\pi^I_J)\cap\widehat{\cal M}^{\bf t}_J(S)^{sing}
\end{array}\right\}
\end{equation}
where $k\in\N$ and $\partial^k\widehat V_{J\alpha}^{\ast}$ is the
$k$-boundary of $\widehat V_{J\alpha}^{\ast}$ as defined in
Definition~\ref{def:A.6}. Moreover, the restriction of the
projection
\begin{equation}\label{e:4.11}
 \hat\pi^I_J: (\hat \pi^I_J)^{-1}(\widehat
V_J^\ast)\cap\widehat{\cal M}^{\bf t}_I(S)\to {\rm
Im}(\hat\pi^I_J)\cap\widehat{\cal M}^{\bf t}_J(S)
\end{equation}
to $(\hat \pi^I_J)^{-1}(\widehat V_J^\ast)\cap\widehat{\cal
M}^{\bf t}_I(S)^\circ$
 is a $|\Gamma_I|/|\Gamma_J|$-fold
(regular) PS covering to ${\rm Im}(\hat\pi^I_J)\cap\widehat{\cal
M}^{\bf t}_J(S)^\circ$. Note, as in \S2, that each $\hat\pi_I:
\widehat{\cal M}^{\bf t}_J(S) ({\rm resp}.\;\widehat{\cal M}^{\bf
t}_J(S)^\circ)\to X$ is also proper. Furthermore, if $(X,E,S)$ is
smoothly good (resp. smoothly (weakly) $k$-good) then the above
each $\widehat{\cal M}^{\bf t}_I(S)$ is a smooth manifold of
dimension ${\rm Ind}(S)$ (though it is not necessarily equal to
$\widehat{\cal M}^{\bf t}_I(S)^\circ$), and each stratum
$\widehat{\cal M}^{\bf t}_I(S)_\alpha$ is a submanifold of it.
\end{corollary}

For the final claim the reader may refer to the proof of
Theorem~\ref{th:3.14}. Corresponding to
Propositions~\ref{prop:2.40},~\ref{prop:2.44} we also have:

\begin{proposition}\label{prop:4.6}
\begin{description}
\item[(i)]  For any given small open neighborhood ${\cal U}$ of
$Z(S)$ in $\cup^n_{i=1}W_i$ there exists $\epsilon\in
(0,\varepsilon)$ such that $\hat\pi_I(\widehat{\cal M}^{\bf
t}_I(S))\subset {\cal U}$ for any ${\bf t}\in B_\epsilon(\R^m)$.

\item[(ii)] If $(X, E, S)$ is oriented then each weakly $k$-PS manifold
$\widehat{\cal M}^{\bf t}_I(S)$ in (\ref{e:4.8}) may be oriented
naturally and the map $\hat\pi^I_J$ in (\ref{e:4.11}) is
orientation preserving for any $J\subset I\in{\cal N}$.

\item[(iii)] For any two ${\bf t},\, {\bf t}'\in
B_\varepsilon(\R^m)_{res}$  the families of weakly $k$-PS manifolds
$\widehat{\mathcal M}^{{\bf t}}(S)=\{\widehat{\mathcal M}_I^{{\bf
t}}(S) | I\in{\cal N}\}$ and $\widehat{\mathcal M}^{{\bf
t}'}(S)=\{\widehat{\mathcal M}_I^{{\bf t}'}(S) | I\in{\cal N}\}$ are
cobordant in the sense that there exist generic paths
$\gamma:[0,1]\to B_\varepsilon(\R^m)$
from ${\bf t}$ to ${\bf t}'$ such that:\\
 (a) for
$I\in{\cal N}$ the set
$$\Upsilon_I^{-1}(\gamma):=\bigcup_{t\in
[0,1]}(\Upsilon_I^{\gamma(t)})^{-1}(0)\times\{t\}
$$
is a weakly $k$-PS manifold of dimension ${\rm Ind}(S)+1$ and with
boundary $\{0\}\times\widehat{\mathcal M}_I^{{\bf
t}}(S)\cup(-\{1\}\times\widehat{\mathcal M}_I^{{\bf t}'}(S))$
({\rm Here  $(X, E, S)$ has been assumed to be oriented. Otherwise
the negative sign before $\{1\}$ should be removed}.)\\
 (b) the
family $\Upsilon^{-1}(\gamma)=\{\Upsilon_I^{-1}(\gamma)\,|\,
I\in{\cal N}\}$ is compatible in the similar sense to
(\ref{e:4.11});\\
 (c) $\Upsilon_I^{-1}(\gamma)$ has a compact
closure in $Cl(\widehat V_I^\ast)\times[0,1]\subset\widehat
V_I^+\times [0,1]$.

\item[(iv)] Assume $\delta\in (0,\varepsilon]$ so that
 the closure of $\cup_{I\in{\mathcal N}}
\hat\pi_I(\widehat{\mathcal M}_I^{{\bf t}}(S))$ in
$\cup^n_{i=1}W_i$ may be contained in $\cup_{I\in{\cal
N}}V_I^\ast$ for any ${\bf t}\in B_\delta(\R^m)$. (This is always
possible by (i) and Theorem~\ref{th:4.4}.) Then for any ${\bf
t}\in B_\delta(\R^m)_{res}$ the set $\cup_{I\in{\mathcal N}}
\hat\pi_I(\widehat{\mathcal M}_I^{{\bf t}}(S))$ is compact.
 Consequently, the
family $\widehat{\cal M}^{\bf t}(S)=\bigl\{\widehat{\cal M}^{\bf
t}_I(S)\,|\, I\in{\cal N}\bigr\}$ is ``like'' an open cover of a
compact weakly $k$-PS manifold.  Specially, each $\widehat{\mathcal
M}_I^{{\bf t}}(S)$ is a finite set provided that ${\rm Ind}(S)=0$.

\end{description}
\end{proposition}

\noindent{\bf Proof.}\quad (i) The claim can be proved by completely
repeating the proof of Proposition~\ref{prop:2.40}(i).

(ii) Combining the proof of Proposition~\ref{prop:2.44} with that of
Theorem~\ref{th:3.14}(ii.3) one can complete the proof.

 (iii) Let ${\cal P}^l({\bf t},{\bf t}')$
denote the space of all $C^l$-smooth paths $\gamma:[0,1]\to
B_\varepsilon(\R^m)$ from ${\bf t}$ to ${\bf t}'$. It is a Banach
manifold. Consider the obvious pullback PS Banach bundle system
\begin{eqnarray*}
&&\bigl({\bf P}_1^\ast\widehat {\mathcal F}^\ast, \widehat
V^\ast\times{\cal P}^l({\bf t},{\bf t}')\times [0, 1]\bigr)\\
&&= \bigl\{\bigl({\bf P}_1^\ast\widehat F_I^\ast, \widehat
V_I^\ast\times{\cal P}^l({\bf t},{\bf t}')\times [0,1]\bigr),
\hat\pi_I, \hat\pi^I_J, \hat\Pi_I,\hat\Pi^I_J, \hat
p_I,\Gamma_I\,\bigm|\, J\subset I\in{\mathcal N}\bigr\}
\end{eqnarray*}
and its global section $\Psi=\{\Psi_I\,|\, I\in{\cal N}\}$ given
by
\begin{eqnarray*}
&&\Psi_I: \widehat V_I^\ast\times{\cal P}^l({\bf t},{\bf
t}')\times [0,1]\to {\bf P}_1^\ast\widehat F_I^\ast,\\
&& (\hat x_I,\gamma, s)\mapsto \hat S_I(\hat x_I)+
\sum^n_{i=1}\sum^{m_i}_{j=1}\gamma(s)_{ij}(\hat\sigma_{ij})_I(\hat
x_I).
\end{eqnarray*}
Here $\gamma(s)=\{\gamma(s)_{ij}\}_{\substack{
 1\le j\le m_i\\
 1\le i\le n}}$, ${\bf P}_1$ is the projection to the first factor.
As in Step 1 in the proof of Theorem~\ref{th:3.14}(ii.4)
 we can prove that $\Psi_I$ is  transversal to the zero section. Since
${\cal N}$ is finite,  the standard arguments yields a residual
subset ${\cal P}^l({\bf t},{\bf t}')_{reg}\subset{\cal P}^l({\bf
t},{\bf t}')$ such that for each $I\in{\cal N}$ and $\gamma\in{\cal
P}^l({\bf t},{\bf t}')_{reg}$ the PS section
\begin{eqnarray*}
&&\Psi_I^\gamma: \widehat V_I^\ast\times [0,1]\to {\bf
P}_1^\ast\widehat F_I,\\
&& (\hat x_I, s)\mapsto \hat S_I(\hat x_I)+
\sum^n_{i=1}\sum^{m_i}_{j=1}\gamma(s)_{ij}(\hat\sigma_{ij})_I(\hat
x_I)
 \end{eqnarray*}
is a (weakly) $k$-Fredholm PS section of index ${\rm Ind}(S)+ 1$
that is transversal to the zero section. That is,  for each
$\gamma\in{\cal P}^l({\bf t},{\bf t}')_{reg}$ the global section
$\Psi^\gamma=\{\Psi^\gamma_I\,|\, I\in{\cal N}\}$ of the pullback PS
Banach bundle system
\begin{eqnarray*}
&&\bigl({\bf P}_1^\ast\widehat {\mathcal F}^\ast, \widehat
V^\ast\times [0, 1]\bigr)\\
&&= \bigl\{({\bf P}_1^\ast\widehat F_I^\ast, \widehat
V_I^\ast\times [0,1]), \hat\pi_I, \hat\pi^I_J,
\hat\Pi_I,\hat\Pi^I_J, \hat p_I,\Gamma_I\,\bigm|\, J\subset
I\in{\mathcal N}\bigr\},
\end{eqnarray*}
 is a (weakly) $k$-Fredholm PS section of index ${\rm Ind}(S)+ 1$
that is transversal to the zero section.  It follows that for each
$\gamma\in{\cal P}^l({\bf t},{\bf t}')_{reg}$ and $I\in{\cal N}$ the
weakly $k$-PS manifold with boundary (and corner)
\begin{equation}\label{e:4.12}
(\Psi^\gamma_I)^{-1}(0)=\bigl\{(\hat x_I, s)\in \widehat
V_I^\ast\times [0,1]\,|\, \Psi_I(\hat x_I,\gamma,s)=0\bigr\}
\end{equation}
 forms a cobordism between the weakly $k$-PS manifolds $\widehat{\mathcal
M}_I^{{\bf t}}(S)$ and $\widehat{\mathcal M}_I^{{\bf t}'}(S)$. This
completes the proof of (a) in (iii).

The conclusion of (b) in (iii) is obtained naturally. In order to
prove (c) in (iii), note that
$\Upsilon_I^{\gamma(s)}(x)=\Psi_I^{\gamma}(x,s)$ for
$(x,s)\in\widehat V_I^\ast\times [0, 1]$. Let $(\hat x_{Ik},
s_k)\in(\Psi^\gamma_I)^{-1}(0)$ and $s_k\to s_0\in [0,1]$. Then
$$\hat S_I(\hat x_{Ik})+
\sum^n_{i=1}\sum^{m_i}_{j=1}\gamma(s_k)_{ij}(\hat\sigma_{ij})_I(\hat
x_{Ik})=0
$$
for $k=1,2,\cdots$, and $\gamma(s_k)_{ij}\to \gamma(s_0)_{ij}$ in
$\R^m$ as $k\to\infty$. By Theorem~\ref{th:4.4} the sequence
$\{(\hat x_{Ik}, \gamma(s_k))\}$ has a convergent subsequence in
 $Cl(\widehat V_I^\ast)\times\R^m\subset\widehat V_I^+\times\R^m$.
 It follows that $\{\hat x_{Ik}\}$ has a convergence subsequence in
 $Cl(\widehat V_I^\ast)\subset\widehat V_I^+$. Hence
$\Upsilon_I^{-1}(\gamma)=(\Psi^\gamma_I)^{-1}(0)$ has a compact
closure in $Cl(\widehat V_I^\ast)\times[0,1]\subset\widehat
V_I^+\times [0,1]$. (iii) is proved.

(iv) The proof can be obtained by  replaceing the word ``manifold''
therein by ``PS manifold'', and Corollary~\ref{cor:2.39} by
Corollary~\ref{cor:4.5} in the proof of
Proposition~\ref{prop:2.40}(iii). \hfill$\Box$\vspace{2mm}

As in Remark~\ref{rem:2.43}, if $Z(S)\subset{\rm Int}(X)$) then PS
manifold $(\Psi^\gamma_I)^{-1}(0)$ in (\ref{e:4.12}) has only
boundary but no corner for  $\varepsilon>0$  sufficiently small. In
this case the family $\widehat{\cal M}^{\bf
t}(S)=\bigl\{\widehat{\cal M}^{\bf t}_I(S)\,|\, I\in{\cal N}\bigr\}$
is ``like'' an open cover of a compact weakly $k$-PS manifold. If
$Z(S)\cap\partial X\ne\emptyset$  then the PS manifold
$(\Psi^\gamma_I)^{-1}(0)$ might have only boundary but corner.

Clearly, if $S$ is transversal to the zero section so is the section
$\widehat S=\{\widehat S_I:\,I\in{\cal N}\}$. In this case
$S^{-1}(0)$ is a compact PS orbifold of dimension ${\rm Ind}(S)$ and
${\bf t}={\bf 0}\in B_\varepsilon(\R^m)_{res}$. It follows from
Proposition~\ref{prop:4.6} that $\cup_{I\in{\mathcal N}}
\hat\pi_I(\widehat{\mathcal M}_I^{{\bf 0}}(S))=S^{-1}(0)$ and that
$\widehat{\mathcal M}^{{\bf 0}}(S)=\{\widehat{\mathcal M}_I^{{\bf
0}}(S):\,I\in{\cal N}\}$ is cobordant to $\widehat{\mathcal M}^{{\bf
t}}(S)=\{\widehat{\mathcal M}_I^{{\bf t}}(S):\,I\in{\cal N}\}$ for
any ${\bf t}\in B_\varepsilon(\R^m)_{res}$. Actually,
$\widehat{\mathcal M}^{{\bf t}}(S)$ is exactly a resolution of the
PS orbifold $S^{-1}(0)$.

\subsection{Framework I}\label{sec:4.2}

This framework is designed for the application to the (closed
string) Gromov-Witten invariants. We shall develop our theory in a
weaker framework though what is needed in GW-theory as in
\cite{LiT} is the case that $(X,E,S)$ is a (weakly) PS
$1$-Fredholm orbibundle which is smoothly (weakly) $1$-good and
satisfies Assumptions~\ref{ass:4.1},~\ref{ass:4.7}.

From now on, besides the assumptions in Section~\ref{sec:4.1} we
also assume that
$k\ge 1$ and \\
$\bullet$ $(X,E,S)$ is {\bf oriented},\\
$\bullet$ all PS Banach orbifolds satisfy
Assumption~\ref{ass:4.1}.\\
As before let $X_0$ be the top strata of $X$. Denote by
\begin{eqnarray*}
&&\widehat{\cal M}^{\bf t}_I(S)_0=\widehat{\cal M}^{\bf
t}_I(S)\cap(\hat\pi_I)^{-1}(X_0),\\
&&\partial^i\widehat{\cal M}^{\bf t}_I(S)_{0}=\widehat{\cal M}^{\bf
t}_I(S)\cap(\hat\pi_I)^{-1}(\partial^i X_{0})\quad\forall i\in\N.
\end{eqnarray*}
In this subsection we also write $\partial\widehat{\cal M}^{\bf
t}_I(S)_{0}=\partial^1\widehat{\cal M}^{\bf t}_I(S)_{0}$.
 By (\ref{e:4.10}) it holds that for any $J\subset I\in{\cal N}$ and
 $i\in\N$,
\begin{equation}\label{e:4.13}
\left.\begin{array}{ll}
 \hat\pi^I_J\bigl((\hat
\pi^I_J)^{-1}(\widehat V_J^\ast)\cap\widehat{\cal M}^{\bf
t}_I(S)_0\bigr)={\rm
Im}(\hat\pi^I_J)\cap\widehat{\cal M}^{\bf t}_J(S)_0,\\
\hat\pi^I_J\bigl((\hat \pi^I_J)^{-1}(\widehat
V_J^\ast)\cap\partial^i\widehat{\cal M}^{\bf t}_I(S)_0\bigr)={\rm
Im}(\hat\pi^I_J)\cap\partial^i\widehat{\cal M}^{\bf t}_J(S)_0.
\end{array}\right\}
\end{equation}

In order to use the family $\widehat{\mathcal M}^{{\bf t}}(S),{\bf
t}\in B_\delta(\R^m)_{res}$ constructing the virtual Euler class we
need to make a stronger assumption than $1$-regularity of $X$:

\begin{assumption}\label{ass:4.7}
{\rm  The top stratum $X_0$ of $X$  is contained in the regular part
$X^\circ$ of $X$.}
\end{assumption}

 If we only discuss
the results on the level of chains, this assumption is not actually
needed. Under Assumption~\ref{ass:4.7} the top stratum
$\widehat{\cal M}^{\bf t}_I(S)_0$ is contained in the regular part
$\widehat{\cal M}^{\bf t}_I(S)^\circ$ of $\widehat{\cal M}^{\bf
t}_I(S)$ and other strata of $\widehat{\cal M}^{\bf t}_I(S)$ have at
most dimension $r-2$. As in Section~\ref{sec:2.6} we can use the
family $\{\widehat{\cal M}^{\bf t}_I(S)_0\,|\, I\in{\cal N}\}$ to
construct a family of cobordant rational singular chains of
dimension ${\rm Ind}(S)$ in $X$, denoted by
\begin{equation}\label{e:4.14}
 e(E, S)^{\bf t}:=``\sum_{I\in{\cal
N}}"\frac{1}{|\Gamma_I|}\bigl\{\hat\pi_I:\widehat{\cal M}^{\bf
t}_I(S)_0\to X\bigr\}\;\forall{\bf t}\in B_\delta(\R^m)_{res}.
\end{equation}
 Then as before the boundary of $e(E, S)^{\bf t}$ is
given by
\begin{equation}\label{e:4.15}
\partial e(E, S)^{\bf t}=``\sum_{I\in{\cal
N}}"\frac{1}{|\Gamma_I|}\bigl\{\hat\pi_I: \partial\widehat{\cal
M}^{\bf t}_I(S)_0\to X\bigr\}.
\end{equation}
Corresponding with Claim~\ref{cl:2.46} we have:

\begin{claim}\label{cl:4.8}
$\partial e(E, S)^{\bf t}$ is a rational singular cycle in $X$ and
 the class of it, $[\partial e(E, S)^{\bf t}]$, is zero in $H_{r-1}(X,\Q)$.
\end{claim}

 As in Remark~\ref{rem:2.43}(iii) let
$\sim$ be the equivalence relation in the disjoint union
$\coprod_{I\in{\cal N}}\widehat{\mathcal M}_I^{{\bf t}}(S)_{0}$
generated by setting $\hat y_I\sim\hat z_J$ if $J\subset I$ and
$\hat\pi^I_J(\hat y_I)=\hat z_J$. Denote by
\begin{equation}\label{e:4.16}
\overline{\mathcal M}^{{\bf t}}(S)_{0}:= \coprod_{I\in{\cal
N}}\widehat{\mathcal M}_I^{{\bf t}}(S)_{0}/\sim
\end{equation}
and by $\hat q_I:\widehat{\mathcal M}_I^{{\bf t}}(S)_{0}\to
\overline{\mathcal M}^{{\bf t}}(S)_{0}$ the restriction of the
obvious quotient map to $\widehat{\mathcal M}_I^{{\bf t}}(S)_{0}$.
Since $\widehat{\mathcal M}_I^{{\bf
t}}(S)_{0}\subset\widehat{\mathcal M}_I^{{\bf t}}(S)^\circ$ we
have a well-defined label function
\begin{equation}\label{e:4.17}
\lambda: \overline{\mathcal M}^{{\bf t}}(S)_{0}\to\Q,\;\bar
x\mapsto\frac{|\hat q_I^{-1}(\bar x)|}{|\Gamma_I|}\quad{\rm if}\;
\bar x\in {\rm Im}(\hat q_I),
\end{equation}
where $|\hat q_I^{-1}(\bar x)|$ is the number of elements in the
set $\hat q_I^{-1}(\bar x)$ and $\hat x_I$ is any element in $\hat
q_I^{-1}(\bar x)$. With the same way we also define
\begin{equation}\label{e:4.18}
\partial\overline{\mathcal M}^{{\bf t}}(S)_{0}:=
\coprod_{I\in{\cal N}}\partial\widehat{\mathcal M}_I^{{\bf
t}}(S)_{0}/\sim
\end{equation}
Then the label function above naturally restricts to it, and
 all $\hat\pi_I$, $I\in{\cal N}$, can be glued into a
map $\hat\pi$ from $\overline{\mathcal M}^{{\bf t}}(S)_{0}$ to
$X$. It is not hard to check that as rational singular chains in
$X$,
\begin{equation}\label{e:4.19}
 e(E, S)^{\bf t}=[\hat\pi: \overline{\mathcal M}^{{\bf t}}(S)_{0}\to X].
\end{equation}

{\bf If  the top stratum $X_0$ is a Banach orbifold without
 boundary and corner}, then for ${\bf t}$ sufficiently small $e(E, S)^{\bf t}$ is a
rational singular cycle in $X$, called  a {\bf virtual Euler
cycle} of the triple $(X, E, S)$, and its homology class
\begin{equation}\label{e:4.20}
e(E, S)=[e(E, S)^{\bf t}]\in H_r(M,\Q)
 \end{equation}
 is said to be the {\bf virtual Euler class} of the triple.
 Combining the above arguments with ones
in Section~\ref{sec:2.6} we can also prove that the homology class
$e(E, S)$ is independent of choices of the PS section
$\tilde\sigma_{ij}$ of class ${\cal A}$ and open subsets $W_i$.

Now we begin to discuss the properties of the virtual Euler cycle
(class) of the (weakly) PS Banach $k$-Fredholm orbibundles. Their
proofs can be obtained by combining those of the corresponding
properties in \S~\ref{sec:2} and \S~\ref{sec:3}, and are left to the
interested readers. We only point out the places where some
differences shall occur.

Let $(X,E,S)$ be as above and $\Lambda\subset Z(S)$ be a nonempty
compact subset. As in Section~\ref{sec:2.7}, if one only requires: \\
 (i) $\{W_i\}^n_{i=1}$ is an open cover of $\Lambda$,\\
(ii) ${\cal W}^\ast$ is an open neighborhood of $\Lambda$, denoted by ${\cal W}^\ast(\Lambda)$,\\
then we get the corresponding systems of PS Banach bundles,
\begin{eqnarray}
&&\bigl(\widehat{\mathcal E}(\Lambda), \widehat W(\Lambda)\bigr):=
\bigl\{\bigl(\widehat E_I, \widehat W_I\bigr), \hat\pi_I,
\hat\Pi_I, \Gamma_I, \hat\pi^I_J,\hat\Pi^I_J,
\lambda^I_J\,\bigm|\, J\subset I\in{\mathcal N}\bigr\},\nonumber\\
&&\bigl(\widehat{\mathcal F}(\Lambda), \widehat V(\Lambda)\bigr):=
\bigl\{\bigl(\widehat F_I, \widehat V_I\bigr), \hat\pi_I,
\hat\Pi_I, \Gamma_I, \hat\pi^I_J,\hat\Pi^I_J,
\lambda^I_J\,\bigm|\, J\subset I\in{\mathcal
N}\bigr\},\nonumber\\
&&\bigl(\widehat{\mathcal F}^\ast(\Lambda), \widehat
V^\ast(\Lambda)\bigr):= \bigl\{\bigl(\widehat F_I^\ast, \widehat
V_I^\ast\bigr), \hat\pi_I, \hat\Pi_I, \Gamma_I,
\hat\pi^I_J,\hat\Pi^I_J, \lambda^I_J\,\bigm|\, J\subset
I\in{\mathcal N}\bigr\},\nonumber
\end{eqnarray}
where $V_I$, $\widehat V_I$, $\widehat F_I$ (resp. $V_I^\ast$,
$\widehat V_I^\ast$, $\widehat F_I^\ast$) are defined as in
(\ref{e:2.65}) and (\ref{e:2.66}) (resp.  (\ref{e:2.71})).
Corresponding with the sections $\Upsilon^{({\bf
t})}=\{\Upsilon^{({\bf t})}_I\,|\, I\in{\cal N}\}$ in
Corollary~\ref{cor:4.5}  we have a family of (weakly) $k$-Fredholm
PS sections of $\bigl(\widehat{\mathcal F}^\ast(\Lambda), \widehat
V^\ast(\Lambda)\bigr)$,
$$\Bigl\{\Upsilon^{({\bf
t})}(\Lambda)=\{\Upsilon^{({\bf t})}_I(\Lambda)\,|\, I\in{\cal
N}\}\Bigm|\;{\bf t}=\{t_{ij}\}_{\substack{ 1\le j\le m_i\\
 1\le i\le n}}\in B_\varepsilon(\R^m)\Bigr\},$$
$$
 \Upsilon^{({\bf t})}_I(\Lambda):
\widehat V^\ast_I\to \widehat F^\ast_I,\; \hat x_I \mapsto
\Upsilon_{I}(\hat x, {\bf t})=\hat S_I(\hat x_I)+
 \sum^n_{i=1}\sum^{m_i}_{j=1}t_{ij}(\hat\sigma_{ij})_I(\hat x_I),
 $$
 and a residual subset
$B_\varepsilon(\R^m)_{res}^\Lambda$ of $B_\varepsilon(\R^m)$ such that for each
$I\in{\cal N}$ the following are true.\\
(a) For each ${\bf t}\in B_\varepsilon(\R^m)_{res}^\Lambda$ the  PS
section  $\Upsilon_I^{({\bf t})}(\Lambda)$ is a (weakly)
$k$-Fredholm section that is transversal to the
zero section,\\
 (b) $(\Upsilon_I^{({\bf t})}(\Lambda))^{-1}(0)$ is a
PS manifold with compact closure in $Cl(\widehat
V_I^\ast)\subset\widehat V_I^+$,\\
 (c) For any ${\bf t}, {\bf
t}'\in B_\varepsilon(\R^m)_{res}^\Lambda$ the corresponding PS
manifolds $(\Upsilon_I^{({\bf t})}(\Lambda))^{-1}(0)$ and
$(\Upsilon_I^{({\bf
t}')}(\Lambda))^{-1}(0)$ are cobordant. \\
 So we get a family of cobordant rational chains near $\Lambda$ of dimension $r$,
  denoted by
$$
 e(E, S)^{\bf
t}_\Lambda:=``\sum_{I\in{\cal N}}"
\frac{1}{|\Gamma_I|}\{\hat\pi_I:\widehat{\cal M}^{\bf t}_I
(S,\Lambda)_{0}\to X\}\;\forall{\bf t}\in
B_\varepsilon(\R^m)_{res}^\Lambda,
$$
 and called {\bf virtual Euler chains relative to}
$\Lambda$ of the triple $(X, E, S)$. Here  $\widehat{\cal M}^{\bf
t}_I(S, \Lambda):= (\Upsilon^{({\bf t})}_I(\Lambda))^{-1}(0)$ is a
weakly $k$-PS manifold of dimension $r$ and with compact closure in
$Cl(\widehat V^\ast_I)\subset\widehat V_I^{+}$, and $\widehat{\cal
M}^{\bf t}_I(S, \Lambda)_{0}$ is the top stratum of $\widehat{\cal
M}^{\bf t}_I(S, \Lambda)$, i.e., $\widehat{\cal M}^{\bf t}_I(S,
\Lambda)_{0}=\widehat{\cal M}^{\bf t}_I(S,
\Lambda)\cap(\hat\pi_I)^{-1}(X_0)$.

 Let  $P$  be a smooth orbifold of finite
dimension and $f:X\to P$ be a PS map. If $\alpha\in H^r(P, \R)$ has
a differential form representative  $\alpha^\ast$ such that the
support ${\rm supp}(f^\ast\alpha^\ast)$ of $f^\ast\alpha^\ast$ is
contained in $\Lambda$ then we can define the integration
\begin{equation}\label{e:4.21}
\int_{e(E, S)^{\bf t}_\Lambda}f^\ast\alpha^\ast:=``\sum_{I\in{\cal
N}}"\frac{1}{|\Gamma_I|}\int_{\widehat{\cal M}^{\bf t}_I
(S,\Lambda)_0}(f\circ\hat\pi_I)^\ast\alpha^\ast.
\end{equation}
Here the integration in the right side is understand as follows:
Multiply the form $f^\ast\alpha^\ast$ by  a PS function to $[0,1]$
which is $1$ outside a small neighborhood of $\widehat{\cal
M}^{\bf t}_I (S,\Lambda)\setminus\widehat{\cal M}^{\bf t}_I
(S,\Lambda)_0$ and $0$ inside a smaller neighborhood and then pass
to a limit as these neighborhoods get smaller. As before,
$\sum_{I\in{\cal N}}$ is put in double quotation marks in
(\ref{e:4.21})  because  the repeating part is only counted once.
Corresponding with Proposition~\ref{prop:2.50} and
Proposition~\ref{prop:3.16} we have:

\begin{proposition}\label{prop:4.9}({\bf First localization formula}).
Let Assumption~\ref{ass:4.7} be satisfied so that the virtual
Euler class $e(E, S)$ is well-defined as in (\ref{e:4.14}). Then
 there exists a residual subset $B_\varepsilon(\R^m)^\star_{res}\subset
  B_\varepsilon(\R^m)_{res}^\Lambda$ such that
$$\langle e(E,S), f^\ast\alpha\rangle=\int_{e(E, S)^{\bf t}_\Lambda}f^\ast\alpha^\ast
$$
for any ${\bf t}\in B_\varepsilon(\R^m)^\star_{res}$.
\end{proposition}

As above Proposition~\ref{prop:2.51}, if $e(E,S)$ exists and
$\Lambda_i$, $i=1,\cdots, p$, are all connected components of
$Z(S)$, then one can start from these $\Lambda_i$ to construct
homology classes $e(E,S)_{\Lambda_i}\in H_r(X, \Q)$,
$i=1,\cdots,p$ such that
$$
e(E,S)=\sum^p_{i=1}e(E,S)_{\Lambda_i}.
$$

Let $(X,E,S)$  be as above and $X=\cup_{\alpha\in\Lambda}X_\alpha$
be the stratification decomposition of $X$. Denote
 by $r_\alpha$ the index of the restriction of $S$ to $X_\alpha$.
Then for each $\alpha\in\Lambda$,
$X^\alpha:=\cup_{\beta\prec\alpha}X_\beta$ is a PS Banach
suborbifold of $X$, and $(X^\alpha, E|_{X^\alpha}, S|_{X^\alpha})$
is also a (weakly) PS Banach $k$-Fredholm orbibundle of index
$r_\alpha$. Moreover, a PS suborbifold
$Y=\cup_{\alpha\in\Lambda_0}Y_\alpha\subset X$ indexed by a
partial order subset $\Lambda_0\subset\Lambda$ with the maximal
element $\alpha_0$ is said to have {\bf finite codimension $n$} if
the top stratum $Y_{\alpha_0}$ of $Y$ is a Banach suborbifold of
codimension $n$ in $X_{\alpha_0}$ and other strata of $Y$ have at
least finite codimension $n$ in the corresponding strata of $X$.
As arguments below Proposition~\ref{prop:3.16}, $(Y, E|_Y, S|_Y)$
is only a {\bf weakly} PS Banach $k$-Fredholm orbibundle of index
$r_{\alpha_0}-n$. When $Y$ is closed in $X$, $Z(S|_Y)$ is also
closed and thus  compact in $Y$ by the compactness of $Z(S)$ in
$X$. In the following we {\bf always consider these two kinds of
PS suborbifolds and assume that they are closed in $X$}. Then $(Y,
E|_Y, S|_Y)$ is also a PS Banach Fredholm orbibundle with rich PS
sections of ${\cal A}|_Y$, where ${\cal A}|_Y$ is the class of PS
sections of $E|_Y$ obtained by restricting the PS sections of $E$
of class ${\cal A}$ to $Y$.

Note that  each local lift $(\widetilde W, {\widetilde E}_W,
{\widetilde S}_W)$ of $(X, E, S)$ with $W\cap Y\ne\emptyset$
restricts to a natural lift of $(Y, E|_Y, S|_Y)$,
\begin{equation}\label{e:4.22}
(\widetilde W^y, {\widetilde E}^y_W, {\widetilde S}^y_W)=
\bigl(\pi_W^{-1}(W\cap Y), {\widetilde E}_W|_{\pi_W^{-1}(W\cap
Y)}, {\widetilde S}_W|_{{\widetilde S}_W}\bigr).
\end{equation}
Since $Y$ is closed in  $X$ one may easily prove that ${\widetilde
S}^y_W$ is locally proper relative to class $\widetilde{\cal
A}(W\cap Y)$ near $Z({\widetilde S}^y_W)$. Here $\widetilde{\cal
A}(W\cap Y)$ denotes the class of those sections consisting of the
restrictions of the PS sections in $\widetilde{\cal A}(W)$ to
$\widetilde W^y$.

Carefully checking arguments above Proposition~\ref{prop:3.17} and
ones between Claim~\ref{cl:4.2} and Definition~\ref{def:4.3} we can
prove that ${\widetilde S}^y_W$ is weakly $k$-quasi transversal to
the zero section relative to class $\widetilde{\cal A}(W\cap Y)$. So
we  get
\begin{center}
{\bf $(Y, E|_Y, S|_Y)$ is a (weakly) PS Banach $k$-Fredholm
orbibundle which is (weakly) $k$-good if $Y=X^\alpha$, and a
weakly PS Banach $k$-Fredholm orbibundle which is weakly $k$-good
if $Y$ is a PS suborbifold of finite codimension that is closed in
$X$.}
\end{center}
The following is a generalization of Proposition~\ref{prop:2.51} and
Proposition~\ref{prop:3.17}.

\begin{proposition}\label{prop:4.10}
 ({\bf The restriction principle}).  Let $(X, E, S)$ be an oriented
 (weakly) PS Banach $k$-Fredholm orbibundle of index $r$ and with compact zero locus,
  have rich PS sections  of class ${\cal A}$ and be (weakly) $k$-good relative to the class ${\cal A}$.
Suppose that  a closed subset $Y\subset X$ belongs to one of the
above two classes of  PS suborbifolds.
 \footnote{Actually we only need to assume that some closed
neighborhood of $Z(S|_Y)$ in $Y$ is flat in $X$.} Assume also that
$Y$ is {\bf flat} in $X$ in the sense that {\rm for each point $y\in
Y$ and each open neighborhood $U$ of it in $X$
  there exists another open neighborhood $V\subset U$ of $y$ in $X$ such that
  $Cl(V\cap Y)=Cl(V)\cap Y$}.
Then there exist a family of cobordant virtual Euler chains  of
the triple $(X, E, S)$,
$$
 e(E, S)^{{\bf
t}'}=``\sum_{I\in{\cal N}'}"
\frac{1}{|\Gamma_I|}\{\hat\pi_I:\widehat{\cal M}^{{\bf t}'}_I
(S)_0\to X\}\;\forall{\bf t}'\in B_\varepsilon(\R^{m'})_{res},
$$
and another family of cobordant virtual Euler chains of  $(Y,
E|_Y, S|_Y)$,
$$
 e(E|_Y, S|_Y)^{\bf
t}=``\sum_{I\in{\cal N}}"
\frac{1}{|\Gamma_I|}\{\hat\pi_I^y:\widehat{\cal M}^{\bf t}_I
(S|_Y)_0\to Y\}\;\forall{\bf t}\in B_\varepsilon(\R^m)_{res},
$$
such that for
$$
{\bf t}'=\{t_{ij}\}_{\substack{ 1\le j\le m_i\\
 1\le i\le n'}}\in B_\varepsilon(\R^{m'})_{res}\quad{\rm and}\quad
{\bf t}:=\{t_{ij}\}_{\substack{ 1\le j\le m_i\\
 1\le i\le n}}
$$
one has
$$
 \left.\begin{array}{ll}
\widehat{\cal M}^{{\bf t}'}_I(S)\cap
\hat\pi_I^{-1}(Y)=\emptyset\;\forall I\in{\cal N}'\setminus{\cal
N},\\
 \widehat{\cal M}^{\bf t}_I
(S|_Y)=\widehat{\cal M}^{{\bf t}'}_I(S)\cap
\hat\pi_I^{-1}(Y),\\
\widehat{\cal M}^{\bf t}_I (S|_Y)^\circ=\widehat{\cal M}^{{\bf
t}'}_I(S)^\circ\cap \hat\pi_I^{-1}(Y),\\
\widehat{\cal M}^{\bf t}_I (S|_Y)^{sing}=\widehat{\cal M}^{{\bf
t}'}_I(S)^{sing}\cap \hat\pi_I^{-1}(Y),\\
 \hat\pi^{yI}_J=\hat\pi^I_J|_{\widehat{\cal M}^{\bf t}_I (S|_Y)}
\end{array}
\right\}
$$
for any $I\in{\cal N}$ and $J\subset I$. It follows that
$$
 \overline{\cal M}^{\bf
t}(S|_Y)^\circ=\overline{\cal M}^{{\bf t}'}(S)^\circ\cap
\hat\pi^{-1}(Y)\quad{\rm and}\quad
\hat\pi^y=\hat\pi|_{\overline{\cal M}^{\bf t}(S|_Y)^\circ}.
$$
Roughly speaking,  $e(E|_Y, S|_Y)^{\bf t}$ is the intersection of
$e(E, S)^{{\bf t}'}$ with $Y$.
\end{proposition}

It should be noted that if $Y\cap X_0=\emptyset$ then
$\overline{\cal M}^{{\bf t}'}(S)^\circ\cap
\hat\pi^{-1}(Y)=\emptyset$. In this case even if both classes
$e(E|_Y, S|_Y)$ and $e(E, S)$ exist we cannot get any relation
between them. However, this restriction principle is convenient
for understanding relations among chains. \vspace{2mm}

\noindent{\bf Proof of Proposition~\ref{prop:4.10}.}
 This result can be obtained by almost
repeating the proof of Proposition~\ref{prop:2.51}. We only give
different arguments.  The first is to omit the second sentence ``
One can check that Assumption~\ref{ass:2.45} also holds for the
suborbifold $Y$.'' in the proof of Proposition~\ref{prop:2.51}.
Next, in the present case  the corresponding result with
(\ref{e:2.125}) cannot be obtained. The purpose that we assume $Y$
to be flat is for getting it.

Finally, let us replace the arguments below (\ref{e:2.127}) as
follows. Recall the arguments between (\ref{e:4.3}) and
(\ref{e:4.4}). We have the $\Gamma_i$-invariant PS cut-off functions
$\gamma_i:\widetilde W_i\to [0, 1]$ with support in $\widetilde
W_i^1$, $i=1,\cdots, n'$ such that
$Z(S)\subset\cup^{n'}_{i=1}U^0_i$, where $U_i^0=\pi_i(\widetilde
U_i^0)$ and  $\widetilde U_i^0=\{\tilde x\in\widetilde
W_i\,|\,\gamma_i(\tilde x)>0\}$. Let $\widetilde S_i:\widetilde
W_i\to\widetilde E_i$ be unique $\Gamma_i$-equivariant lifts of
$S|_{W_i}$, $i=1,\cdots, n'$. Since $Cl(U^0_i)\subset
W^1_i\subset\subset W_i$, using the compactness of $Z(\widetilde
S_i)\cap Cl(\widetilde W^1_i)\subset\widetilde W_i$ and the
assumption of richness we can find finitely many PS sections of the
PS Banach bundle $\widetilde E_i\to\widetilde W_i$ of class ${\cal
A}(W_i)$, $\tilde s_{ij}$, $j=1,\cdots, m_i$ such that
$$
D\widetilde S_i(\tilde z)(T_{\tilde z}\widetilde W_i)+ {\rm
span}(\{\tilde s_{i1}(\tilde z),\cdots, \tilde s_{im_i}(\tilde
z)\})=(\widetilde E_i)_{\tilde z}
$$
  for any $\tilde z\in Z(\widetilde S_i)\cap Cl(\widetilde W^1_i)$,
$i=1,\cdots,n'$.
 Set  $\tilde\sigma_{ij}=\tilde\gamma_i\cdot\tilde
s_{ij}$, $j=1,\cdots, m_i$. Then
$$
D\widetilde S_i(\tilde z)+ {\rm span}(\{\tilde\sigma_{i1}(\tilde
z),\cdots,\tilde\sigma_{im_i}(\tilde z)\})=(\widetilde E_i)_{\tilde
z}
$$
for any $\tilde z\in Z(\widetilde S_i)\cap Cl(U^0_i)$ and
$i=1,\cdots, n'$. By Lemma~\ref{lem:2.33} each $\tilde\sigma_{ij}$
yields a global PS section
$\hat\sigma_{ij}=\{(\hat\sigma_{ij})_I\,|\, I\in{\cal N}\}$ of the
system of PS Banach bundles $(\widehat {\mathcal F}, \widehat V)$.

 Note that for each $i=1,\cdots, n$, the restriction
$$\widetilde S_i^y=\widetilde S_i|_{\widetilde
W_i^y}:\widetilde W_i^y\to\widetilde E_i^y
$$
is exactly the lift of the restriction  $S|_{W_i^y}$, and the
restriction of $\gamma_i$ to $\widetilde W_i^{y1}$, denote by
$\gamma_i^y$, has support contained in $\widetilde W_i^{y1}$. Hence
\begin{eqnarray*}
&&\widetilde U^{y0}_i:=\{\tilde x\in\widetilde W_i^y\,|\,
\gamma^y_i(\tilde x)>0\}=\widetilde U^{0}_i\cap\widetilde
W_i^y\subset\subset\widetilde W_i^{y1} \quad{\rm and}\\
&&Z(S|_Y)\subset\bigcup^n_{i=1}U^{y0}_i=\bigcup^n_{i=1}\pi_i(\widetilde
U^{y0}_i).
\end{eqnarray*}
Let $\tilde s_{ij}^y$ denote the restriction of $\tilde s_{ij}$ to
$W_i^y$. Define sections
$$
\tilde\sigma_{ij}^y:=\gamma^y_i\cdot \tilde s_{ij}^y: \widetilde
W_i^y\to\widetilde E_i^y,\;j=1,\cdots, m_i.
$$
They are restrictions of $\tilde\sigma_{ij}$ to $W_i^y$. Since both
the sections $\widetilde S_i$ and $\widetilde S^y_i$ are weakly
$k$-Fredholm, by increasing $m_i$ we may assume
$$
D\widetilde S^y_i(\tilde z)+ {\rm span}(\{\tilde\sigma^y_{i1}(\tilde
z),\cdots,\tilde\sigma^y_{im_i}(\tilde z)\})=(\widetilde
E^y_i)_{\tilde z}
$$
for any $\tilde z\in Z(\widetilde S^y_i)\cap Cl(U^{y0}_i)$ and
$i=1,\cdots, n$.

By the assumptions, $(X, E, S)$ is (weakly) $k$-good relative to
class ${\cal A}$ and $(Y, E|_Y, S|_Y)$  is weakly $k$-good relative
to class ${\cal A}|_Y$ as showed above Proposition~\ref{prop:4.10}.
As the arguments below (\ref{e:4.5}) we use these to get a global PS
section  of $\bigl(\widehat{\mathcal F}^{\ast}, \widehat
V^{\ast}\bigr)$, $\Upsilon^{({\bf t}')}=\{\Upsilon^{({\bf
t}')}_I\,|\, I\in{\cal N}'\}$ given by
$$
 \Upsilon^{({\bf t}')}_I: \widehat
V^{\ast}_I\to \widehat F^{\ast}_I,\; \hat x_I \mapsto\hat S_I(\hat
x_I)+
 \sum^{n'}_{i=1}\sum^{m_i}_{j=1}t_{ij}(\hat\sigma_{ij})_I(\hat x_I),
$$
and  that of $\bigl(\widehat{\mathcal F}^{y\ast}, \widehat
V^{y\ast}\bigr)$, $\Upsilon^{y({\bf t})}=\{\Upsilon^{y({\bf
t})}_I\,|\, I\in{\cal N}\}$ given by
$$
\Upsilon^{y({\bf t})}_I: \widehat V^{y\ast}_I\to \widehat
F^{y\ast}_I,\; \hat x_I \mapsto\hat S_I^y(\hat x_I)+
 \sum^n_{i=1}\sum^{m_i}_{j=1}t_{ij}(\hat\sigma_{ij}^y)_I(\hat x_I).
 $$
Here ${\cal N}=\{I\in{\cal N}'\,|\, \max I\le n\}$. Furthermore we
can find a small $\varepsilon>0$ and a residual subset
$B_\varepsilon(\R^{m'})_{res}\subset B_\varepsilon(\R^{m'})$,
$m'=\sum^{n'}_{i=1}m_i$ such that for any
$$
{\bf t}'=\{t_{ij}\}_{\substack{ 1\le j\le m_i\\
 1\le i\le n'}}\in B_\varepsilon(\R^{m'})_{res}\quad{\rm and}\quad
{\bf t}:=\{t_{ij}\}_{\substack{ 1\le j\le m_i\\
 1\le i\le n}}
$$
 the above PS sections  $\Upsilon^{({\bf
t}')}_I$ and $\Upsilon^{y({\bf t})}_I$ are all weakly $k$-Fredholm
PS sections that are transversal to the zero section (by shrinking
${\cal W}^\ast$ and increasing $n$ and $n'$ if necessary). (Note
that when ${\bf t}'$ takes over $B_\varepsilon(\R^{m'})_{res}$ all
corresponding ${\bf t}:=\{t_{ij}\}_{\substack{ 1\le j\le m_i\\
 1\le i\le n}}$ form a residual subset
$B_\varepsilon(\R^{m})_{res}\subset B_\varepsilon(\R^{m})$.) Set
\begin{eqnarray*}
&& \widehat{\cal M}^{{\bf t}'}(S)=\bigl\{\widehat{\cal M}^{{\bf
t}'}_I(S)=(\Upsilon^{({{\bf t}'})}_I)^{-1}(0) \,|\,I\in{\cal N}'\bigr\}\quad{\rm and}\\
 &&\widehat{\cal M}^{\bf t}(S|_Y)=\bigl\{\widehat{\cal M}^{\bf
t}_{I}(S|_Y)=(\Upsilon^{({\bf t})}_{I})^{-1}(0) \,|\,I\in{\cal
N}\bigr\},
\end{eqnarray*}
and repeat the remnant arguments in the proof of
Proposition~\ref{prop:2.51} to arrive at the desired conclusions.
\hfill$\Box$\vspace{2mm}

Recall that Proposition~\ref{prop:2.51} is used in the proofs of
Propositions~\ref{prop:2.54}, \ref{prop:2.56} and \ref{prop:2.57}.
When proving the corresponding result with  \ref{prop:2.56}, i.e.,
the following Proposition~\ref{prop:4.14}, we need the flatness of
$Y=X\times\{i\}$ in $X\times [0, 1]$, $i=0, 1$. These hold
obviously. However, for the corresponding results with
Propositions~\ref{prop:2.54} and \ref{prop:2.57}, i.e., the
following Proposition~\ref{prop:4.12} and \ref{prop:4.15} we shall
add this condition.

\begin{corollary}\label{cor:4.11}
 ({\bf Second localization formula}).
 Let $f:X\to P$ be a PS map from $X$ to an oriented smooth orbifold of finite
 dimension, and $Q\subset P$ be an oriented
 closed suborbifold of dimension $m$ such
 that\\
  (i) $Y:=f^{-1}(Q)\subset X$ is a PS Banach
 suborbifold of codimension $\dim P-m$,\\
 (ii) $Y$ is flat in $X$ and satisfies Assumption~\ref{ass:4.7},\\
 (iii) the top stratum of $Y$, $Y_0$, is contained in the top stratum
 $X_0$ of $X$. \\
  Assume also that  $e(E,S)$ exists and that $(Y, E|_Y, S|_Y)$ is orientable and $e(E|_Y, S|_Y)$
 exists (after choosing an orientation of it).
 Then for any class
$\kappa\in H_\ast(Q,\R)$ of dimension $\dim P-r$ it holds that
$$\langle PD_Q(\kappa), f_\ast(e(E|_Y, S|_Y))\rangle= \langle PD_P(\kappa),
f_\ast(e(E,S))\rangle.
$$
Here $PD_Q(a)$ (resp. $PD_P(a)$) is the Poincar\'e duality of $a$
in $H^\ast(Q,\R)$ (resp. $H^\ast(P,\R)$).
\end{corollary}

Since $Y_0\subset X_0$, $(Y, E|_Y, S|_Y)$ has index $r-\dim P+ m$
and it follows from Proposition~\ref{prop:4.10} that
$\widehat{\cal M}^{\bf t}_I (S|_Y)_0=\widehat{\cal M}^{{\bf
t}'}_I(S)\cap \hat\pi_I^{-1}(Y_0)=\widehat{\cal M}^{{\bf
t}'}_I(S)_0\cap \hat\pi_I^{-1}(Y)$ for any $I\in{\cal N}$. The
desired equality can be obtained as in the proof of
Proposition~\ref{prop:1.12}. If the top stratum of $Y$ can be
contained a lower stratum of $X$, then $f_\ast(e(E|_Y, S|_Y))\in
H_\ast(P,\Q)$ has degree less than $r+m-\dim P$. In this case
$\langle PD_Q(\kappa), f_\ast(e(E|_Y, S|_Y))\rangle$ is
necessarily zero because $PD_Q(\kappa)\in H^{r+m-\dim P}(Q,\R)$.
But $\langle PD_P(\kappa), f_\ast(e(E,S))\rangle$ might be
nonzero. So the assumption $Y_0\subset X_0$ is reasonable and
necessary.

Now we shall give a common generalization of
Proposition~\ref{prop:2.54} and Proposition~\ref{prop:3.18}.

Let $X$ be a separable PS Banach orbifold satisfying
Assumption~\ref{ass:4.7}, and $(X, E^{(i)}, S^{(i)})$ be the
oriented (weakly) PS Banach $k$-Fredholm orbibundles of index
$r_i$ and with compact zero locus $Z(S^{(2)})$, $i=1,2$. Assume
also that $(X, E^{(i)}, S^{(i)})$ has rich local PS sections of
class ${\cal A}_i$ and is (weakly) $k$-good relative to the class
${\cal A}_i$, $i=1,2$.  Consider the direct sum $E:=E^{(1)}\oplus
E^{(2)}\to X$, $S:=S^{(1)}\oplus S^{(2)}:X\to E$ and ${\cal
A}:={\cal A}_1\oplus{\cal A}_2$.
 Clearly, $Z(S)=Z(S^{(1)})\cap Z(S^{(2)})$ and
$Z(S)=Z(S^{(1)}|_{Z(S^{(2)}})$. Then $(X, E, S)$ is an oriented
(weakly) PS Banach $k$-Fredholm orbibundle with rich local PS
sections of class ${\cal A}$. One also easily shows that $(X, E,
S)$ is also (weakly) $k$-good relative to the class ${\cal A}$
from the arguments above Proposition~\ref{prop:3.18} and ones
above Definition~\ref{def:4.3} (by showing that ${\cal A}(W)={\cal
A}_1(W)\oplus{\cal A}_2(W)$ and $\tilde{\cal A}(W)=\tilde{\cal
A}_1(W)\oplus\tilde{\cal A}_2(W)$, which are easily proved).

\begin{proposition}\label{prop:4.12}
({\bf Stability}). Let $(X, E^{(i)}, S^{(i)})$, $i=1,2$ be above.
Furthermore, assume: (i) $Z(S^{(1)})\cap X_0\ne\emptyset$ and
$Z(S^{(2)})\cap X_0\ne\emptyset$, (ii) $S^{(2)}$ is transversal to
the zero section at each point $x\in Z(S)$, (iii) $DS^{(2)}|_{Z(S)}$
is surjective. Then one has:
\begin{description}
\item[(a)] For some small open neighborhood ${\cal U}$ of $Z(S)$
in $X$ the intersection  $Z(S^{(2)})^\star:=Z(S^{(2)})\cap{\cal
U}$ is a weakly $k$-PS Banach orbifold of dimension $r_2$. (In
particular, the top stratum of $Z(S^{(2)})^\star$ is contained in
$X_0$.)

 \item[(b)]
$D(S^{(1)}|_{Z(S^{(2)})^\star}):(TZ(S^{(2)}))|_{Z(S)} \to
E^{(1)}|_{Z(S)}$ is a PS Fredholm bundle map of index ${\rm
Ind}(D(S^{(1)}|_{Z(S^{(2)})^\star})={\rm Ind}(DS)$ (which
precisely means
$$
{\rm Ind}\bigl(D(S^{(1)}|_{Z(S^{(2)})^\star})(x)\bigr)={\rm
Ind}(DS(x))\;\forall x\in Z(S)\;).
$$

\item[(c)] If $(Z(S^{(2)})^\star, E^{(1)}|_{Z(S^{(2)})^\star},
S^{(1)}|_{Z(S^{(2)})^\star})$ is (weakly) $k$-good relative to
${\cal A}_1|_{Z(S^{(2)})^\star}$ and some closed neighborhood of
$Z(S)$ in $Z(S^{(2)})^\star$ is flat in $X$, then
$(Z(S^{(2)})^\star, E^{(1)}|_{Z(S^{(2)})^\star},
S^{(1)}|_{Z(S^{(2)})^\star})$ and $(X, E,S)$ have respectively
virtual Euler chains $C_1$ and  $C$  such that $C=C_1$. In
particular, if both $e(E,S)$ and $ e(E^{(1)}|_{Z(S^{(2)})^\star},
S^{(1)}|_{Z(S^{(2)})^\star})$ exist,  then $C$ and $C_1$ can be
chosen as closed weakly $k$-PS manifolds (thus being cycles).
Consequently,
$$e(E,S) =(i_{Z(S^{(2)})^\star})_\ast e(E^{(1)}|_{Z(S^{(2)})^\star}, S^{(1)}|_{Z(S^{(2)})^\star}).$$
Here $(i_{Z(S^{(2)})^\star})_\ast$ is the homomorphism between
homology groups induced by the inclusion $i_{Z(S^{(2)})^\star}:
Z(S^{(2)})^\star\to X$.
\end{description}
\end{proposition}

\noindent{\bf Proof.}\quad The proof may be obtained by almost
completely repeating the proof of Proposition~\ref{prop:2.54}. The
main changes are:\\
$\bullet$ All words ``Banach orbifold'', ``Banach orbibundle'',
''bundle'',  ``suborbifold'' and ``section'' {\it are changed into
words} ``PS Banach orbifold'', ``PS Banach orbibundle'', ``PS
bundle'', ``PS suborbifold'' and ``PS
section''.\\
$\bullet$ {\it Delete}  ``As before we assume ...'' from the third
line above (\ref{e:2.129}) to the end of this paragraph, and {\bf
add}  ``Since the lift PS section $\widetilde S_i^{(1)}:\widetilde
W_i\to\widetilde E_i^{(1)}$ of $S^{(1)}|_{W_i}$  restricts to a
lift PS section $\widetilde S^{(1)z}_i=\widetilde
S^{(1)}_i|_{\widetilde W_i^z}:\widetilde W_i^z\to\widetilde
E_i^{(1)z}$ of $S^{(1)}|_{W_i^z}$, as in the proof of
Proposition~\ref{prop:1.13} we can show that on each stratum the
vertical differential $D\widetilde S^{(1)z}_i(x)$ at each zero
point $x$ of $\widetilde S^{(1)z}_i$ is Fredholm and has the same
index as $D\widetilde
S^{(1)}_i(x)$. This yields the first conclusion.''\\
$\bullet$ {\it Change} ``Note that $S^{(2)}$ is transversal to the
zero section..... are strongly transversal to the zero section.''
in the fourth pargraph of proof of Proposition~\ref{prop:2.54}
{\bf into}:
 ``Note that $S^{(2)}$ is transversal to the zero section
and that
 $(Z(S^{(2)})^\star, E^{(1)}|_{Z(S^{(2)})^\star},
S^{(1)}|_{Z(S^{(2)})^\star})$ is (weakly) $k$-good relative to
${\cal A}_1|_{Z(S^{(2)})^\star}$. As in the proof of
Proposition~\ref{prop:4.10}, for each $i=1,\cdots,n$ we can
construct the PS sections $\tilde\sigma_{ij}^{(1)}=\gamma_i\cdot
\tilde s_{ij}^{(1)}: \widetilde W_i\to\widetilde
E_i^{(1)}\subset\widetilde E_i$ and
$$
\tilde\sigma_{ij}^{(1)z}=\gamma_i\cdot \tilde s_{ij}^{(1)z}:
\widetilde W_i^z\to\widetilde E_i^{(1)z}
$$
such that for sufficiently small ${\cal W}^\ast$, $\varepsilon>0$
and each ${\bf t}=\{t_{ij}\}_{\substack{ 1\le j\le m_i\\
 1\le i\le n}}$ in a residual subset
$B_\varepsilon(\R^m)_{res}\subset B_\varepsilon(\R^m)$ the global PS
section  of the PS Banach bundle system $\bigl(\widehat {\mathcal
F}^\ast, \widehat V^\ast\bigr)$, $\Upsilon^{({\bf
t})}=\{\Upsilon^{({\bf t})}_I\,|\, I\in{\cal N}\}$,
$$
 \Upsilon^{({\bf t})}_I: \widehat
V^\ast_I\to \widehat F^\ast_I,\; \hat x_I \mapsto\hat S_I(\hat x_I)+
 \sum^n_{i=1}\sum^{m_i}_{j=1}t_{ij}(\hat\sigma^{(1)}_{ij})_I(\hat x_I),
 $$
and that of $\bigl(\widehat {\mathcal F}^{(1)z\ast}, \widehat
V^{z\ast}\bigr)$, $\Upsilon^{({\bf t})(1)z}=\{\Upsilon^{({\bf
t})(1)z}_I\,|\, I\in{\cal N}\}$,
$$
 \Upsilon^{({\bf t})(1)z}_I: \widehat
V^{z\ast}_I\to \widehat F^{(1)z\ast}_I,\; \hat x_I \mapsto\hat
S_I^{(1)z}(\hat x_I)+
 \sum^n_{i=1}\sum^{m_i}_{j=1}t_{ij}(\hat\sigma_{ij}^{(1)z})_I(\hat x_I),
$$
 are all weakly $k$-Fredholm PS sections that
 are transversal to the zero section.''
{\bf The remnant arguments} are the same as ones in the proof of
Proposition~\ref{prop:2.54}. (Note that in the present case our
proof seems to be simper  because many technically local arguments
are actually absorbed in the good assumption of PS orbibundles.)
\hfill$\Box$\vspace{2mm}

Let $(X, E^{(i)}, S^{(i)})$ be the oriented (weakly) PS Banach
$k$-Fredholm orbibundle of index $r$ and with compact zero locus
$Z(S^{(i)})$, $i=1,2$. Assume also that $(X, E^{(i)}, S^{(i)})$ has
rich local PS sections of class ${\cal A}_i$ and is (weakly)
$k$-good relative to the class ${\cal A}_i$, $i=1,2$.  As in
Section~\ref{sec:2.8} or Proposition~\ref{prop:3.20} we can
similarly define a homotopy between them and get  the following
result corresponding with Proposition~\ref{prop:2.56}.

\begin{proposition}\label{prop:4.13}
({\bf Homotopy}) Let $(X, E^{(i)}, S^{(i)})$, $i=0,1$, be as above,
and $(X\times [0, 1], E, S)$ be an oriented homotopy (weakly) PS
Banach $k$-Fredholm orbibundle between them which is (weakly)
$k$-good relative to class  ${\cal A}$. Assume that $\partial
X_0=\emptyset$. Then there exist virtual Euler chains of $(X\times
[0, 1], E, S)$, $(X, E^{(0)}, S^{(0)})$ and $(X, E^{(1)}, S^{(1)})$
respectively, $e(E, S)^{\bf t}$, $e(E^{(0)}, S^{(0)})^{{\bf t}_0}$
and $e(E^{(1)}, S^{(1)})^{{\bf t}_1}$ such that
$$
\partial e(E, S)^{\bf t}=e(E^{(0)}, S^{(0)})^{{\bf t}_0}\cup
(-e(E^{(1)}, S^{(1)})^{{\bf t}_1}).
$$
In particular $e(E^{(0)}, S^{(0)})=e(E^{(1)}, S^{(1)})$ because
they exist by Assumption~\ref{ass:4.7} and $\partial
X_0=\emptyset$.
 \end{proposition}

As pointed out below Proposition~\ref{prop:2.56},   for any
$m\in\N$ if the $m$-boundary $\partial^m X_0\ne\emptyset$ then we
have more general conclusion in the proposition above, i.e., after
neglecting the orientation the $m$-boundary of $e(E, S)^{\bf t}$,
$$
\partial^m e(E, S)^{\bf t}=\partial^{m-1}e(E^{(0)}, S^{(0)})^{{\bf t}_0}\cup
\partial^{m-1}e(E^{(1)}, S^{(1)})^{{\bf t}_1}\cup e(E, S)^{\bf t}_m,
$$
 where $e(E, S)^{\bf t}_m$ is defined by
$$
``\sum_{I\in{\cal
N}}"\frac{1}{|\Gamma_I|}\bigl\{\hat\pi_I:\widehat{\cal M}^{\bf
t}_I(S)\cap(\hat\pi_I)^{-1}(\partial^m X_0\times (0, 1))\to
X\times [0, 1]\bigr\}.
$$

Finally, we study the functoriality. Let $(X, E, S)$ and $(X', E',
S')$ be two oriented (weakly) PS Banach $k$-Fredholm orbibundles
with compact zero loci. Assume that $(X, E, S)$ (resp. $(X', E',
S')$) has rich PS sections of class ${\cal A}$ (resp. ${\cal A}'$)
and is also (weakly) $k$-good relative to class ${\cal A}$ (resp.
${\cal A}'$).
 A {\bf morphism} from $(X, E, S)$ to $(X', E', S')$ is a PS Banach
orbibundle embedding $(\psi, \Psi):E\to E'$ with following
properties:
\begin{description}
\item[(A)] $S\circ \psi=\Psi\circ S'$ and $Z(S')=\psi(Z(S))$;

\item[(B)] For any $x\in Z(S)$, the differential $d\psi(x):T_xX\to
T_{\psi(x)}X'$ and the above restriction $\Psi_x: E_x\to E'_{f(x)}$
induce isomorphisms
\begin{eqnarray*}
&&d\psi(x): {\rm Ker}(DS(x))\to {\rm Ker}(DS'(\psi(x)))\quad{\rm and}\\
&&[\Psi_x]: {\rm Coker}(DS(x))\to {\rm Coker}(DS'(\psi(x))),
\end{eqnarray*}
and the resulting isomorphism from ${\rm det}(DS)$ to ${\rm
det}(DS')$ is orientation preserving. (These are all understand on
levels of lifts.)
\end{description}

Clearly, $(X, E, S)$ and $(X', E', S')$ have the same index.
Moreover, (A) implies that $Z(S')\subset\psi(X)$. Later writing
``$(\psi,\Psi)^\ast{\cal A}'\supseteq{\cal A}$'' will be understood
as follows. (By the definition of the morphism we may assume: $X$ is
a PS Banach suborbifold of $X'$ , $E$ is a PS suborbibundle of
$E'|_X$ and $S=S'|_X$. $\psi$ and $\Psi$ become the natural
inclusions.) For a PS orbifold chart $(\widetilde W',\Gamma',\pi')$
on $X'$ and $W=W'\cap X$ let $\widetilde W\subset\widetilde W'$ be a
PS submanifold that is stable under $\Gamma$ and compatible with the
inclusion maps, such that the restriction $(\widetilde
W,\Gamma',\pi)=(\widetilde W, \Gamma'|_{\widetilde
W},\pi|_{\widetilde W})$ is a PS Banach orbifold chart on $X$ with
support $W$. Then ``$(\psi,\Psi)^\ast{\cal A}'\supseteq{\cal A}$''
means that each PS section of the PS bundle $\widetilde
E\to\widetilde W$ of class ${\cal A}(\widetilde W)$ can be extended
into a PS section of the PS bundle $\widetilde E'\to\widetilde W'$
of class ${\cal A}'(\widetilde W)$.

The following is corresponding result with
Proposition~\ref{prop:2.57}:

\begin{proposition}\label{prop:4.14}
Let $(X, E, S)$ and $(X', E', S')$ be as above, and $(\psi, \Psi)$
be a morphism from $(X, E, S)$ to $(X', E', S')$ satisfying
$(\psi,\Psi)^\ast{\cal A}'\supseteq{\cal A}$. Suppose that $X$ and
$X'$ satisfy Assumption~\ref{ass:4.7} and that some closed
neighborhood of $Z(S')=\psi(Z(S))$ in $\psi(X)$ is flat in $X'$.
 Then there exist a family of strongly cobordant virtual Euler chains of the
triple $(X', E', S')$,
$$
 e(E', S')^{{\bf
t}}=``\sum_{I\in{\cal N}}"
\frac{1}{|\Gamma'_I|}\{\hat\pi'_I:\widehat{\cal M}^{{\bf t}}_I
(S')_0\to X'\}\;\forall{\bf t}\in B_\varepsilon(\R^{m})_{res},
$$
and that of strongly cobordant virtual Euler chains of the triple
$(X, E, S)$,
$$
 e(E, S)^{{\bf
t}}=``\sum_{I\in{\cal N}}"
\frac{1}{|\Gamma'_I|}\{\hat\pi_I:\widehat{\cal M}^{{\bf t}}_I
(S)_0\to X\}\;\forall{\bf t}\in B_\varepsilon(\R^{m})_{res},
$$
and a family of PS embeddings $\{\widehat\psi_I\,|\, I\in{\cal
N}\}$, $\widehat\psi_I:\widehat W_I\to\widehat W'_I$, which are also
compatible with the projections $\hat\pi^I_J$ and $\hat\pi'^I_J$,
such that for any $I\in{\cal N}$ and $J\subset I$,
\begin{eqnarray*}
&& \widehat\psi_I(\widehat{\cal M}^{\bf
t}_I(S)^\circ)=\widehat{\cal M}^{\bf t}_I(S')^\circ\cap(\hat\pi'_I)^{-1}(\psi(X)),\\
&&\widehat\psi_I(\widehat{\cal M}^{\bf
t}_I(S)^{sing})=\widehat{\cal M}^{\bf t}_I(S')^{sing}\cap(\hat\pi'_I)^{-1}(\psi(X)),\\
&&\widehat\psi_J\circ\hat\pi^I_J
 =\hat\pi'^I_J|_{\widehat\psi_I(\widehat{\cal M}^{\bf t}_I (S))}\circ \widehat\psi_I,\\
&&\Bigl(\hat\pi'^I_J\bigl((\hat \pi'^I_J)^{-1}(\widehat
V'^\ast_J)\cap\widehat{\cal M}^{\bf t}_I(S')\bigr)\Bigr)\cap
(\hat\pi'_J)^{-1}(\psi(X))\\
&&=\widehat\psi_I\Bigl(\hat\pi^I_J\bigl((\hat
\pi^I_J)^{-1}(\widehat V_J^\ast)\cap\widehat{\cal M}^{\bf
t}_I(S)\bigr)\Bigr),\\
&&\Bigl({\rm Im}(\hat\pi'^I_J)\cap\widehat{\cal M}^{\bf
t}_J(S')\Bigr)\cap(\hat\pi'_J)^{-1}(\psi(X))\\
&&=\widehat\psi_I\bigl({\rm Im}(\hat\pi^I_J)\cap\widehat{\cal
M}^{\bf t}_J(S)\bigr).
\end{eqnarray*}
If $\overline W$ and $\overline W'$ are the spaces constructed
from $\{\widehat W_I\,|\, I\in{\cal N}\}$ and $\{\widehat
W'_I\,|\, I\in{\cal N}\}$ respectively then $\{\widehat\psi_I\,|\,
I\in{\cal N}\}$ induce a natural map $\bar\psi:\overline
W\to\overline W'$ such that
\begin{eqnarray*}
&&\bar\psi(\overline{\cal M}^{\bf t}(S)^\circ)=\overline{\cal
M}^{\bf t}(S')^\circ\cap
(\hat\pi')^{-1}(\psi(X)),\\
&&\bar\psi\circ\hat\pi=\hat\pi'|_{\bar\psi(\overline{\cal M}^{\bf
t}(S)^\circ)}\circ\bar\psi.
\end{eqnarray*}
 Consequently, if $Z(S)\cap X_0\ne\emptyset$, $Z(S')\cap X'_0\ne\emptyset$
  and $\psi(X_0)\subset X'_0$, $\partial X_0=\partial X'_0=\emptyset$ and hence
 $e(E,S)$ and $e(E', S')$ exist, then
$\langle\psi^\ast\alpha, e(E,S)\rangle= \langle\alpha, e(E',
S')\rangle$ for any $\alpha\in H^\ast(X', X'\setminus
\psi(X);\R)$. (In this case we say $\psi_\ast e(E, S)$ to be the
intersection  of $e(E', S')$ with $\psi(X)$.)
\end{proposition}

\noindent{\bf Proof.}\quad The proof may be obtained by almost
completely repeating the proof of Proposition~\ref{prop:2.57}. The
main changes are:\\
$\bullet$ The terms such as ``manifold'', ``orbifold'' and
``bundle'' are changed into corresponding terms in the PS
category.\\
$\bullet$ Below Claim~\ref{cl:2.58}, from ``For the chosen
trivializations....''  to the end of the next paragraph  ``.....are
strongly transversal to the zero section for all $I\in{\cal N}$.''
{\bf is changed into}: ``Note that $(\psi,\Psi)^\ast{\cal
A}'\supseteq {\cal A}$. By the explanation above
Proposition~\ref{prop:4.14}, for $\widetilde W_i'$ and $\widetilde
W_i$ as above each PS section of the PS bundle $\widetilde
E_i\to\widetilde W_i$ of class ${\cal A}(\widetilde W_i)$ can be
extended into a PS section of the PS bundle $\widetilde
E'_i\to\widetilde W'_i$ of class ${\cal A}'(\widetilde W'_i)$. Since
$(X,E,S)$ has rich PS sections of class ${\cal A}$, we may choose PS
sections $\tilde s_{ij}$, $j=1,\cdots, m_i$ of the PS bundle
$\widetilde E_i\to\widetilde W_i$ such that
$$
D\widetilde S_i(\tilde x)(T_{\tilde x}\widetilde W_i)+ {\rm
span}\{\tilde s_{i1}(\tilde x),\cdots, \tilde s_{im_i}(\tilde
x)\}=(\widetilde E_i)_{\tilde x}\quad\forall\tilde x\in Z(\widetilde
S_i).
$$
(This is always possible by shrinking $W_i$ if necessary.) Let
$\tilde s'_{ij}$, $j=1,\cdots, m_i$ be the PS sections of
$\widetilde E'_i\to\widetilde W'_i$ of class ${\cal A}'(\widetilde
W'_i)$ obtained by extending $\tilde s_{ij}$, $j=1,\cdots,m_i$.
 With the same proof as in Claim~\ref{cl:2.58} we get that
$$
D\widetilde S'_i(\tilde x)(T_{\tilde x}\widetilde W_i')+ {\rm
span}\{\tilde s'_{i1}(\tilde x),\cdots, \tilde s'_{im_i}(\tilde
x)\}=(\widetilde E'_i)_{\tilde x}\quad\forall\tilde x\in
Z(\widetilde S'_i)=Z(\widetilde S_i).
$$

Consider the sections $\sigma'_{ij}:=\gamma'_i\cdot\tilde s'_{ij}$
of $\widetilde E'_i$ with supports in $\widetilde W'^1_i$, and
$\sigma_{ij}:=\gamma'_i\cdot\tilde s_{ij}$ of $\widetilde E_i$ with
supports in $\widetilde W^1_i$, $j=1,\cdots, m_i$. Clearly,
$\sigma_{i1}=\sigma'_{i1}|_{\widetilde W_i}, \cdots,
\sigma_{im_i}=\sigma'_{im_i}|_{\widetilde W_i}$, $i=1,\cdots, n$.
Note that $Z(S')\subset X$. We can assume that the above small
neighborhood ${\cal W}'^\ast$ of $Z(S')$ in $X$ is contained in
$\cup^n_{i=1}U'^0_i$ (and thus ${\cal W}^{\ast}={\cal W}'^\ast\cap
X\subset\cup^n_{i=1}U^0_i$). Since $(X, E, S)$ (resp. $(X', E',
S')$) is  (weakly) $k$-good relative to class ${\cal A}$ (resp.
${\cal A}'$),  as before we can shrink ${\cal W}'^\ast$  and find a
small $\varepsilon>0$ such that for
 each ${\bf t}=\{t_{ij}\}_{\substack{ 1\le j\le m_i\\
 1\le i\le n}}$ in a residual subset
$B_\varepsilon(\R^m)_{res}\subset B_\varepsilon(\R^m)$ with
$m=\sum^n_{i=1}m_i$, the sections
\begin{eqnarray*}
 \Upsilon^{({\bf t})}_I: \widehat
V^\ast_I\to \widehat F^\ast_I,\; \hat x_I \mapsto\hat S_I(\hat x_I)+
 \sum^n_{i=1}\sum^{m_i}_{j=1}t_{ij}(\hat\sigma_{ij})_I(\hat x_I),\\
  \Upsilon'^{({\bf t})}_{I}: \widehat
V'^{\ast}_I\to \widehat F'^{\ast}_I,\; \hat x_I \mapsto\hat
S'_I(\hat x_I)+
 \sum^{n}_{i=1}\sum^{m_i}_{j=1}t_{ij}(\hat\sigma'_{ij})_I(\hat
 x_I),
\end{eqnarray*}
are all weakly $k$-Fredholm PS sections that are transversal to the
zero section for all $I\in{\cal N}$. The remnant arguments are the
same as  the proof of Proposition~\ref{prop:2.57}.
 \hfill$\Box$\vspace{2mm}

As above Proposition~\ref{prop:3.21}, a {\bf PS Fredholm map of
index $d$} from PS Banach orbifolds $Y$ to $X$ is a PS map whose
restriction to the top stratum (resp. each lower stratum) is a
Fredholm map of index $d$ (resp. index less than or equal to $d$).
Finally we give a generalization of Proposition~\ref{prop:2.59} and
prove it carefully.   Still using the notations therein  we have:

\begin{proposition}\label{prop:4.15}
Let $X$ and $X'$ be two separable PS Banach orbifolds satisfying
Assumptions~\ref{ass:4.1}, \ref{ass:4.7}, and $f:X\to X'$ be a
proper Fredholm map of index $d$.  Let $(X, E, S)$ (resp. $(X', E',
S')$) be the (weakly) PS Banach $k$-Fredholm orbibundle of index $r$
(resp. $r'$) with compact zero locus. Assume that $\bar f:(X, E,
S)\to (X', E', S')$ is a PS bundle map, i.e., a bundle map $\bar f:
E\to E'$ covering $f$ and satisfying $S'\circ f=\bar f\circ S$, and
that each local lifting of $\bar f$ is a Banach space isomorphism
when restricted to each fiber. (That is, $(X, E, S)$ is the
pull-back of $(X', E', S')$ via $f$). Furthermore assume that $f$ is
a PS submersion, $Z(S)\cap X_0\ne\emptyset$ and $f^{-1}(X_0\cap
Z(S))\cap Y_0\ne\emptyset$. (The latter implies $f(Y_0)\subset
X_0$). Then\\
{\bf (i)} $r=r'+d$ and $Z(S')=f^{-1}(Z(S))$ is compact.\\
Moreover, if $(X',E',S')$ has also rich PS sections of class ${\cal
A}'$ and
is (weakly) $k$-good relative to the class ${\cal A}'$, then\\
{\bf (ii)} $(X, E, S)$ is weakly $k$-good relative to the class
 $\bar f^\ast{\cal A}'$, and there exist a resolution $\widehat V$ of $X$
near $Z(S)$ and that $\widehat V'$ of $X'$ near $Z(S')$, and a
natural map $\hat f:\widehat V\to\widehat V'$ induced by $f$,  and
the virtual Euler chain $e(E, S)^{\bf t}$ (resp. $e(E', S')^{\bf
t}$), as a PS map from $\widehat{\cal M}^{\bf t}(S)\subset\widehat
V$ (resp. $\widehat{\cal M}^{\bf t}(S')\subset\widehat V'$) to $X$
(resp. $X'$) satisfy
\begin{equation}\label{e:4.23} f\circ e(E,
S)^{\bf t}=e(E', S')^{\bf t}\circ\hat f.
\end{equation}
In addition, for $i=0, 1, 2$, if $f$ sends the strata $(X, E, S)_i$
of $(X, E, S)$ of codimension $i$ to those $(X', E', S')_i$ of $(X',
E', S')$ of codimension $i$, $\hat f_I$ also sends the strata of
$\widehat{\cal M}^{\bf t}_I(S)$ of codimension $i$ to those of
$\widehat{\cal M}^{\bf t}_I(S')$ of codimension $i$.
\end{proposition}

As conjectured in Remark~\ref{rm:2.60} there should exists some kind
of $\Q$-valued topological degree ${\rm deg}(f)$ of $f$ so that
(\ref{e:4.23}) gives
$$
f_\ast e(E, S)={\rm deg}(f) e(E', S')
$$
in the case that $d=0$ and both $e(E, S)$ and $e(E', S')$
exist.\vspace{2mm}

\noindent{\bf Proof of Proposition~\ref{prop:4.15}.} Since each
local lifting of $f$ restricts to a Banach space isomorphism on each
fiber we have $Z(S')=f^{-1}(Z(S))$. This local lift is also proper
because of the properness of $f$.

{\bf Firstly, we prove (i)}. For $x\in Z(S)$ let  $(\widetilde U,
\Gamma_U,\pi_U)$ be the PS B-orbifold chart centered at $x$,  and
$(\widetilde U^\prime, \Gamma'_{U^\prime},\pi_{U^\prime})$ be one
centered at $f(x)$. As in Section~\ref{sec:2.1} assume that
$(\tilde{\bar f}_U, \tilde f_U,\varphi_U):(\widetilde E_U,
\widetilde U, \Gamma_U)\to (\widetilde E'_{U'},\widetilde U^\prime,
\Gamma'_{U^\prime})$ is a local representation of $\bar f$. Let
$\widetilde S_U:\widetilde U\to\widetilde E_U$ and $\widetilde
S'_{U'}:\widetilde U'\to\widetilde E'_{U'}$ be the local lifts of
$S$ and $S'$ in the above PS Banach charts respectively.

 By (\ref{e:2.4}), $\tilde{\bar f}_U\circ\widetilde
S_U= \widetilde S'_{U'}\circ\tilde f_U$.  It follows that
$$
 d\tilde{\bar f}_U(\widetilde S_U(\tilde x))\circ d\widetilde
S_U(\tilde x)=d\widetilde S'_{U'}(\tilde f_U(\tilde x))\circ
d\tilde f_U(\tilde x)
$$
for the lifting $\tilde x$ of $x$ in $\widetilde U$.  Note that
$\widetilde S_U(\tilde x)=(\tilde x, 0)$ and that $T_{(\tilde x,
0)}\widetilde E_U=T_{\tilde x}\widetilde U\oplus (\widetilde
E_U)_{\tilde x}$. One easily checks that
\begin{eqnarray*}
&&d\tilde{\bar f}_U(\tilde x, 0)(T_{\tilde x}\widetilde U)=d\tilde
f_U({\tilde x})(T_{\tilde x}\widetilde U)\quad{\rm and}\\
&&d\tilde{\bar f}_U(\tilde x,0)((\widetilde E_U)_{\tilde
x})=(\tilde{\bar f}_U)_{\tilde x}((\widetilde E_U)_{\tilde x}).
\end{eqnarray*}
So the vertical differentials $D\widetilde S_U(\tilde x)$ and
$D\widetilde S'_{U'}(\tilde f_U(\tilde x))$ satisfy:
\begin{equation}\label{e:4.24}
 (\tilde{\bar f}_U)_{\tilde x}\circ
D\widetilde S_U(\tilde x)=D\widetilde S'_{U'}(\tilde f_U(\tilde
x))\circ d\tilde f_U(\tilde x).
\end{equation}
 This leads to
$r=r'+ d$. (i) is proved.

{\bf Next we prove that $(X, E, S)$ is weakly $k$-good relative to
the class $\bar f^\ast{\cal A}'$.} For $\tilde y\in Z(\widetilde
S_U)$ denote by $\tilde y'=\tilde f_U(\tilde y)$. For PS sections of
$\widetilde E'_{U'}$, $\sigma'_1,\cdots,\sigma'_k$, we can define PS
sections $\sigma_i: \widetilde U\to \widetilde E_U$ by
\begin{equation}\label{e:4.25}
\sigma_i(\tilde y):=(\tilde{\bar f})^{-1}_{\tilde y}\bigl(
\sigma'_i(\tilde f_U(\tilde y))\bigr)
\end{equation}
because $(\tilde{\bar f}_U)_{\tilde z}$ is the Banach space
isomorphism for each $\tilde z$.
  Suppose that
\begin{equation}\label{e:4.26}
{\rm Im}(D\widetilde S_{U}(\tilde y))+ {\rm
span}(\{\sigma_1(\tilde y),\cdots, \sigma_k(\tilde
y)\})=(\widetilde E_U)_{\tilde y}.
\end{equation}
By (\ref{e:4.24}) we can easily get that
$$
{\rm Im}\bigl(D\widetilde S'_{U'}(\tilde y')\circ d\tilde f_U(\tilde
y')\bigr)+ {\rm span}(\{\sigma'_1(\tilde y'),\cdots,
\sigma'_k(\tilde y')\})=(\widetilde E'_U)_{\tilde y'}.
$$
Note that ${\rm Im}\bigl(D\widetilde S'_{U'}(\tilde y')\circ d\tilde
f_U(\tilde y)\bigr)\subset {\rm Im}\bigl(D\widetilde S'_{U'}(\tilde
y')\bigr)$. We get that
\begin{equation}\label{e:4.27}
{\rm Im}\bigl(D\widetilde S'_{U'}(\tilde y')\bigr)+ {\rm
span}(\{\sigma'_1(\tilde y'),\cdots, \sigma'_k(\tilde
y')\})=(\widetilde E'_U)_{\tilde y'}.
\end{equation}
Since $(X', E', S')$ is (weakly) $k$-good relative to the class
${\cal A}'$ there exist a small open neighborhood $O(\tilde y')$ of
$\tilde y'$ in $\widetilde U'$ and $\eta>0$ such that the section
\begin{eqnarray*}
&&O(\tilde y')\times B_\eta(\R^k)\to \Pi_1^\ast(\widetilde E'_{U'}|_{O(\tilde y')}),\\
&&\Phi': (\tilde z'; t_1,\cdots, t_k)\mapsto \widetilde
S'_{U'}(\tilde z')+\sum^k_{i=1}t_i\sigma'_i(\tilde z'),
\end{eqnarray*}
is a (weakly) $k$-Fredholm PS section that is transversal to the
zero section. Here $\Pi_1: O(\tilde y')\times B_\eta(\R^k)\to
O(\tilde y')$ is the natural projection. As in (\ref{e:4.24}),  we
can derive from (\ref{e:4.25}) that for any $\tilde z\in (\widetilde
S_{U}+\sum^k_{i=1}t_i\sigma_i\bigr)^{-1}(0)$,
$$
(\tilde{\bar f}_U)_{\tilde z}\circ D\bigl(\widetilde
S_{U}+\sum^k_{i=1}t_i\sigma_i\bigr)(\tilde z)= D\bigl(\widetilde
S'_{U'}+\sum^k_{i=1}t_i\sigma_i'\bigr)\circ d\tilde f_U(\tilde z).
$$
As in the proof of Proposition~\ref{prop:3.21}, all these imply that
the section
\begin{eqnarray*}
&&\tilde f^{-1}_U(O(\tilde y'))\times B_\eta(\R^k)\to \Pi_1^\ast(\widetilde E_{U}|_{\tilde f^{-1}_U(O(\tilde y'))}),\\
&&\Phi: (\tilde z; t_1,\cdots, t_k)\mapsto \widetilde S_{U}(\tilde
z)+\sum^k_{i=1}t_i\sigma_i(\tilde z),
\end{eqnarray*}
is  a weakly PS $k$-Fredholm section that is transversal to the zero
section. Hence this shows that $\widetilde S_U$ is weakly $k$-quasi
transversal to the zero section relative to class $\bar
f^\ast\widetilde{\cal A}(\widetilde U')$.

Note that for any ${\bf t}=(t_1,\cdots, t_k)\in\R^k$,
$$
\bigl(\widetilde S_{U}+\sum^k_{i=1}t_i\sigma_i\bigr)^{-1}(0)=
\tilde f^{-1}_U\Bigl[\bigl(\widetilde
S'_{U'}+\sum^k_{i=1}t_i\sigma_i'\bigr)^{-1}(0)\Bigr].
$$
One also easily proves that $\widetilde S_U$ is locally proper
relative to class $\bar f^\ast\widetilde{\cal A}(\widetilde U')$
near $Z(\widetilde S_U)$. So $(X, E, S)$ is weakly $k$-good relative
to the class $\bar f^\ast{\cal A}'$.

{\bf Finally we prove (ii)}.  Since $Z(S')$ is compact we can
choose finitely many such local representations of $\bar f$,
 \begin{equation}\label{e:4.28}
 (\tilde{\bar f}_{U_i}, \tilde f_{U_i},\varphi_{U_i}):(\widetilde E_{U_i},
 \widetilde U_i, \Gamma_{U_i})\to (\widetilde E'_{U'_i},\widetilde U'_i,
\Gamma'_{U'_i}),
\end{equation}
  the corresponding open subsets
$U'^0_i\subset W'^j_i\subset\subset U'^j_i\subset\subset
W'_i\subset\subset U'_i$, $i=1,\cdots, n, j=1,\cdots, n-1$, and
the PS sections
$$
\sigma'_{i1},\cdots,\sigma'_{im_i}:\widetilde U'_i\to \widetilde
E'_{U'_i}
$$
with support in $Cl(\widetilde U'^0_i)$,  such that
\begin{equation}\label{e:4.29}
\left.\begin{array}{ll} Z(S')\subset\cup^n_{i=1}U'^0_{i},\\
 {\rm
Im}(D\widetilde S'_{U'_i}(\tilde y'))+ {\rm
span}(\{\sigma'_{i1}(\tilde y'),\cdots, \sigma'_{im_i}(\tilde
y')\})=(\widetilde E'_{U'_i})_{\tilde y'}\\
\qquad\forall\tilde y'\in Cl(\widetilde U'_{i0})\cap Z(\widetilde
S'_{U'_i}),\;i=1,\cdots, n.
\end{array}
\right\}
\end{equation}
For $i=1,\cdots, n$, set $U^0_i=f^{-1}(U'^0_i)$,
$W^j_i=f^{-1}(W'^j_i)$, $U^j_i=f^{-1}(U'^j_i)$, $W_i=f^{-1}(W'_i)$
and $U_i=f^{-1}(U'_i)$, and
$$
\sigma_{il}(\tilde y):=(\tilde{\bar f})^{-1}_{\tilde y}\bigl(
\sigma'_{il}(\tilde f_U(\tilde y))\bigr),\;l=1,\cdots, m_i,
$$
then  $U^0_i\subset W^j_i\subset\subset U^j_i\subset\subset
W_i\subset\subset U_i$, $j=1,\cdots, n-1$, each $\sigma_{il}$ has
support in $Cl(\widetilde U^0_i)$, and it also holds that
\begin{equation}\label{e:4.30}
\left.\begin{array}{ll}
Z(S)\subset\cup^n_{i=1}U^0_{i},\\
{\rm Im}(D\widetilde S_{U_i}(\tilde y))+ {\rm
span}(\{\sigma_{i1}(\tilde y),\cdots, \sigma_{im_i}(\tilde
y)\})=(\widetilde E_{U_i})_{\tilde y}\\
\qquad\forall\tilde y\in Cl(\widetilde U_{i0})\cap Z(\widetilde
S_{U_i}), \;i=1,\cdots, n.
\end{array} \right\}
\end{equation}

 Let $m=m_1+\cdots+m_n$.  As in (\ref{e:4.5}) and (\ref{e:4.6})
 we have the PS Banach bundle systems
\begin{eqnarray*}
&&\bigl({\bf P}_1^\ast\widehat {\mathcal F}, \widehat
V\times\R^m\bigr)\\
&&=\bigl \{({\bf P}_1^\ast\widehat F_I, \widehat V_I\times\R^m),
\hat\pi_I,  \hat\pi^I_J, \hat\Pi_I, \hat\Pi^I_J, \hat
p_I,\Gamma_I\,\bigm|\, J\subset I\in{\mathcal N}\bigr\},
\nonumber\\
&&\bigl({\bf P}_1^\ast\widehat{\mathcal F}', \widehat
V'\times\R^m\bigr)\\
&&=\bigl \{({\bf P}_1^\ast\widehat F'_I, \widehat V'_I\times\R^m),
\hat\pi'_I,  \hat\pi'^I_J, \hat\Pi'_I, \hat\Pi'^I_J, \hat
p'_I,\Gamma'_I\,\bigm|\, J\subset I\in{\mathcal N}\bigr\} \nonumber
\end{eqnarray*}
and their global PS  sections $\Upsilon=\{\Upsilon_I\,|\, I\in{\cal
N}\}$ and $\Upsilon'=\{\Upsilon'_I\,|\, I\in{\cal N}\}$, where
\begin{eqnarray*}
&&\Upsilon_I: \widehat V_I\times\R^m\to {\bf P}_1^\ast\widehat
F_I,\\
&&\qquad \Bigl(\hat x_I, \{t_{ij}\}_{\substack{
 1\le j\le m_i\\
 1\le i\le n}}\Bigr)\mapsto \hat S_I(\hat x_I)+
 \sum^n_{i=1}\sum^{m_i}_{j=1}t_{ij}(\hat\sigma_{ij})_I(\hat
 x_I),\nonumber\\
 &&\Upsilon'_I: \widehat V'_I\times\R^m\to {\bf P}_1^\ast\widehat
F'_I,\\
&&\qquad \Bigl(\hat x_I, \{t_{ij}\}_{\substack{
 1\le j\le m_i\\
 1\le i\le n}}\Bigr)\mapsto \hat S'_I(\hat x_I)+
 \sum^n_{i=1}\sum^{m_i}_{j=1}t_{ij}(\hat\sigma'_{ij})_I(\hat
 x_I)\nonumber
 \end{eqnarray*}
 So each ${\bf t}=\{t_{ij}\}_{\substack{
 1\le j\le m_i\\
 1\le i\le n}}\in \R^m$ gives a global PS section
$\Upsilon^{({\bf t})}=\{\Upsilon_I^{({\bf t})}\,|\, I\in{\cal
N}\}$ (resp. ${\Upsilon'}^{({\bf t})}=\{{\Upsilon'}_I^{({\bf
t})}\,|\, I\in{\cal N}\}$) of the PS Banach bundle systems
$(\widehat V, \widehat F)$ (resp. $(\widehat V', \widehat F')$),
where
\begin{eqnarray*} &&\Upsilon_I^{({\bf t})}:\widehat
V_I\to\widehat F_I,\;\hat
x_I\mapsto\Upsilon_I(\hat x_I, {\bf t})\quad{\rm and}\\
&&{\Upsilon'}_I^{({\bf t})}:\widehat V'_I\to\widehat F'_I,\;\hat
x_I\mapsto\Upsilon'_I(\hat x'_I, {\bf t}).
\end{eqnarray*}

Carefully applying the resolution arguments in
Section~\ref{sec:2.3}, one can prove that for each $I\in{\cal N}$,
(\ref{e:4.28}) induces a natural group homomorphism $\varphi_I:
\Gamma_I\to\Gamma_I'$, a $\varphi_I$-equivariant PS map $\hat
f_I:\widehat V_I\to \widehat V_I'$ and a $\varphi_I$-equivariant PS
bundle map $\hat{\bar f}_I: \widehat F_I\to\widehat F_I'$ covering
$\hat f_I$ such that for any $J\subset I\in{\cal N}$ it holds that
\begin{equation}\label{e:4.31}
 \hat\pi'^I_J\circ\hat f_I=\hat
f_J'\circ\hat\pi^I_J\;{\rm (resp.}\;\hat\Pi'^I_J\circ\hat{\bar
f}_I=\hat{\bar f}_J'\circ\hat\Pi^I_J)
\end{equation}
on $(\hat\pi^I_J)^{-1}(\hat V_J)$ (resp. $(\hat\Pi^I_J)^{-1}(\hat
F_J)$). So they give a PS map
$$
(\hat f, \hat{\bar f})=\{(\hat f_I, \hat{\bar f}_I)\,|\, I\in{\cal
N}\}: (\widehat V, \widehat F)\to (\widehat V', \widehat F').
$$
 By the assumptions one can check that each $\hat f_I$ is a proper
 PS Fredholm map and sends the strata of $\widehat V_I$ of codimension
$i$ to those of $\widehat V_I'$ of codimension $i$, $i=0, 1, 2$.
Moreover, each $\hat{\bar f}$ restricts to a Banach space
isomorphism on each fiber. The sections $\Upsilon_I^{({\bf t})}$
and ${\Upsilon'}_I^{({\bf t})}$ above satisfy:
\begin{equation}\label{e:4.32}
\hat{\bar f}_I\circ \Upsilon_I^{({\bf t})}={\Upsilon'}_I^{({\bf
t})}\circ\hat f_I
\end{equation}
for any $I\in{\cal N}$.

Since $f$ is proper for each sufficiently small open neighborhood
${\cal W}'^\ast$ of $Z(S')$ in $\cup^n_{i=1}U'^0_i$ the inverse
image ${\cal W}^\ast:=f^{-1}({\cal W}'^\ast)$ is also very small
open neighborhood of $Z(S)$ in $\cup^n_{i=1}U_i^0$. Note that $(X,
E, S)$ and $(X', E', S')$ are $\bar f^\ast{\cal A}$-good and ${\cal
A}$-good respectively. As before, using (\ref{e:4.29}) and
(\ref{e:4.30}) we have a sufficiently small $\varepsilon>0$ and a
residual subset $B_\varepsilon(\R^m)_{res}$ in $B_\varepsilon(\R^m)$
such that for each ${\bf t}=\{t_{ij}\}_{\substack{
 1\le j\le m_i\\
 1\le i\le n}}\in B_\varepsilon(\R^m)_{res}$ both sections
$$
\Upsilon_I^{({\bf t})}:\widehat V_I^\ast\to\widehat F_I^\ast,\;\hat
x_I\mapsto\Upsilon_I(\hat x_I, {\bf t})
$$
and
$$
{\Upsilon'}_I^{({\bf t})}:\widehat V'^\ast_I\to\widehat
F'^\ast_I,\;\hat x_I\mapsto\Upsilon'_I(\hat x'_I, {\bf t})
$$
are (weakly) $k$-Fredholm PS section and also transversal to the
zero section.  Here
\begin{eqnarray*}
&&\widehat V'^\ast_I=(\hat\pi'_I)^{-1}(W'^\ast\cap
V_I'),\quad\widehat F'^\ast_I=\widehat F'_I|_{\widehat
V'^\ast_I},\\
 &&\widehat V^\ast_I=(\hat\pi_I)^{-1}(W^\ast\cap
V_I)\quad{\rm and}\quad \widehat F^\ast_I=\widehat F_I|_{\widehat
V^\ast_I}.
\end{eqnarray*}
As before we set
$$
\widehat{\cal M}^{\bf t}_I(S)=(\Upsilon_I^{({\bf
t})})^{-1}(0)\quad{\rm  and}\quad\widehat{\cal M}^{\bf
t}_I(S')=({\Upsilon'}_I^{({\bf t})})^{-1}(0).
$$
Then (\ref{e:4.32}) implies that for any $I\in{\cal N}$ and
$J\subset I\in{\cal N}$,
\begin{equation}\label{e:4.33}
\hat f_I(\widehat{\cal M}^{\bf t}_I(S))=\widehat{\cal M}^{\bf
t}_I(S')\quad{\rm and}\quad \hat\pi'^I_J\circ\hat f_I=\hat
f_J\circ\hat\pi^I_J.
\end{equation}
Since $f(Y_0)\subset X_0$, it easily follows that the following
commutative diagram holds:
\begin{center}\setlength{\unitlength}{1mm}
\begin{picture}(80,30)
\thinlines \put(10,25){$\overline{\cal M}^{\bf t}(S)_{0}$}
\put(35,25){\vector (1,0){37}} \put(45,27){$\hat f$}
\put(15,22){\vector(0,-1){16}} \put(10,13){$\hat\pi$}
\put(75,25){$\overline{\cal M}^{\bf t}(S')_{0}$}
\put(80,22){\vector(0,-1){16}} \put(71,13){$\hat\pi'$}
\put(13,0){$X$} \put(35,1){\vector(1,0){38}} \put(45,3){$f$}
\put(78,0){$X'$}
\end{picture}\end{center}
which is,  by (\ref{e:4.16}),  is equivalent to (\ref{e:4.23}).
(That is, $e(E, S)^{\bf t}$, as a PS map from $\widehat{\cal M}^{\bf
t}(S)_0\subset\widehat V$ to $X$, and $e(E', S')^{\bf t}$, as a PS
map from $\widehat{\cal M}^{\bf t}(S')_0\subset\widehat V'$ to $X'$,
satisfy (\ref{e:4.23}).) Here we also use $\hat f$ to denote the
induced natural map from $\overline{\cal M}^{\bf t}(S)_{0}$ to
$\overline{\cal M}^{\bf t}(S')_{0}$.

Under the final assumption, it is easily seen that $\hat f_I$ also
sends the strata of $\widehat{\cal M}^{\bf t}_I(S)$ of codimension
$i$ to those of $\widehat{\cal M}^{\bf t}_I(S')$ of codimension $i$,
$i=0, 1, 2$.  \hfill$\Box$\vspace{2mm}

\subsection{Framework II}\label{sec:4.3}

In the studies of the Floer homology and the (open string)
Gromov-Witten invariants one needs to consider  the oriented
(weakly) PS Banach $0$-Fredholm orbibundles with  the strata of
codimension $1$. (Actually, the bundle should also be smoothly
(weakly) $0$-good.) By Definition~\ref{def:3.7} there is a stratum
of $X$, denoted by $X_1$, such that $(X, E, S)_1=(X_1, E|_{X_1},
S|_{X_1})$ is a Banach Fredholm bundle of index ${\rm Ind}(S)-1$. As
before let $(X, E, S)_0=(X_0, E|_{X_0}, S|_{X_0})$ denote the top
stratum of $(X, E, S)$, i.e.  the stratum of codimension $0$. Denote
by $X_{01}=X_0\cup X_1$ and
\begin{eqnarray*}
&&\widehat{\cal M}^{\bf t}_I(S)_i=\widehat{\cal M}^{\bf
t}_I(S)\cap(\hat\pi_I)^{-1}(X_i),\;i=0, 1,\\
&&\widehat{\cal M}^{\bf t}_I(S)_{01}=\widehat{\cal M}^{\bf
t}_I(S)\cap(\hat\pi_I)^{-1}(X_{01}).
\end{eqnarray*}
Then by (\ref{e:4.10}) it holds that
\begin{equation}\label{e:4.34}
\left.\begin{array}{ll}
 \hat\pi^I_J\bigl((\hat
\pi^I_J)^{-1}(\widehat V_J^\ast)\cap\widehat{\cal M}^{\bf
t}_I(S)_i\bigr)={\rm
Im}(\hat\pi^I_J)\cap\widehat{\cal M}^{\bf t}_J(S)_i,\\
\hat\pi^I_J\bigl((\hat \pi^I_J)^{-1}(\widehat
V_J^\ast)\cap\partial^k\widehat{\cal M}^{\bf t}_I(S)_i\bigr)={\rm
Im}(\hat\pi^I_J)\cap\partial^k\widehat{\cal M}^{\bf t}_J(S)_i
\end{array}\right\}
\end{equation}
for any $i=0, 1$ and $k\in\N$.  In general, $\widehat{\cal M}^{\bf
t}_I(S)_{01}$ is a weakly $0$-PS manifold of dimension $r$. Since
the closure of $\widehat{\cal M}^{\bf t}_I(S)$  in $Cl(\widehat
V_I^\ast)$ is compact, the topological boundary of $\widehat{\cal
M}^{\bf t}_I(S)_{01}$ in $\widehat V_I^\ast$
\begin{equation}\label{e:4.35}
\partial_t\widehat{\cal M}^{\bf t}_I(S)_{01}\subseteq\partial\widehat{\cal M}^{\bf
t}_I(S)_{0}\cup\widehat{\cal M}^{\bf t}_I(S)_1\cup\widehat{\cal
M}^{\bf t}_I(S)_{\ge 2} ,
\end{equation}
where $\partial\widehat{\cal M}^{\bf t}_I(S)_{0}=\widehat{\cal
M}^{\bf t}_I(S)_{0}\cap(\hat\pi_I)^{-1}(\partial X_0)$ and $\partial
X_0$ is the boundary Banach orbifold of $X_0$, and $\widehat{\cal
M}^{\bf t}_I(S)_{\ge 2}$ is the union of all strata of
$\widehat{\cal M}^{\bf t}_I(S)$ of codimension more than $1$. We
first prove:

\begin{lemma}\label{lem:4.16}
If $r=1$ and ${\bf t}\in B_\delta(\R^m)_{res}$ then all sets
 $\widehat{\cal M}^{\bf t}_I(S)_1$ and $\partial\widehat{\cal M}^{\bf t}_I(S)_0$
 are finite.
\end{lemma}

\noindent{\bf Proof.}\quad Firstly,  $\cup_{I\in{\mathcal N}}
\hat\pi_I(\widehat{\mathcal M}_I^{{\bf t}}(S)_1)$ is a closed subset
in $\cup_{I\in{\mathcal N}} \hat\pi_I(\widehat{\mathcal M}_I^{{\bf
t}}(S))$.
 In fact, let $\hat\pi_I(\hat x_{Ik})\to x$
for some sequence $\{\hat x_{Ik}\}\subset\hat\pi_I(\widehat{\cal
M}^{\bf t}_I(S)_1)$. Since the set $\cup_{I\in{\mathcal N}}
\hat\pi_I(\widehat{\mathcal M}_I^{{\bf t}}(S))$ is compact by
Proposition~\ref{prop:4.6}(iv), so is $Cl(\widehat{\cal M}^{\bf
t}_I(S)_1)$. It follows that $\hat\pi_I^{-1}\bigl(Cl(\widehat{\cal
M}^{\bf t}_I(S)_1)\bigr)$ is compact in $Cl(\widehat V_I^\ast)$. By
Remark~\ref{rem:2.34}, after passing a subsequence (if necessary),
we may assume $\hat x_{Ik}\to\hat x_I\in\hat\pi_I^{-1}(x)\subset
Cl(\widehat V_I^\ast)\subset\widehat V_I^+$. This $\hat x_I$ must
belong to $\widehat V_I^+\cap\hat\pi_I^{-1}(X_1)$ since $r=1$ and
thus $X=X_0\cup X_1$. Moreover, the compactness of
$\cup_{I\in{\mathcal N}} \hat\pi_I(\widehat{\mathcal M}_I^{{\bf
t}}(S))$ implies $x\in\hat\pi_L(\widehat{\cal M}^{\bf t}_L(S))$ for
some $L\in {\cal N}$, where $L\subset I$ or $I\subset L$. Let
$x=\hat\pi_L(\hat x_L)$ for some $\hat x_L\in\widehat{\cal M}^{\bf
t}_L(S)$. As at the beginning of the proof of Theorem~\ref{th:2.38},
$\Upsilon=\{\Upsilon_I\,|\, I\in{\cal N}\}$ can naturally extend to
a (weakly) PS $0$-Fredholm section of $({\bf P}_1^\ast\widehat
{\mathcal F}^+, \widehat V^+\times\R^m)$, denoted by
$\Upsilon_+=\{\Upsilon_{+_I}\,|\, I\in{\cal N}\}$. Let
$\Upsilon_{+I}^{({\bf t})}$ be the corresponding (weakly) PS
$0$-Fredholm section with $\Upsilon_{I}^{({\bf t})}$ in
Corollary~\ref{cor:4.5}. Then each $\Upsilon_{+I}^{({\bf t})}$
restricts to $\Upsilon_{I}^{({\bf t})}$ on $\widehat V_I^\ast$. By
the compatibility of the family $\{ (\Upsilon_{+I}^{({\bf
t})})^{-1}(0)\,|\, I\in{\cal N}\}$, we have $\hat x_I\in
(\Upsilon_{+I}^{({\bf t})})^{-1}(0)_1$ and thus $\hat
x_L\in\widehat{\cal M}^{\bf t}_L(S)_1$ and $x\in\hat\pi_L(
\widehat{\cal M}^{\bf t}_L(S)_1)$. The claim is proved.

Next, the compactness of $\cup_{I\in{\mathcal N}}
\hat\pi_I(\widehat{\mathcal M}_I^{{\bf t}}(S))$ implies that
$\cup_{I\in{\mathcal N}} \hat\pi_I(\widehat{\mathcal M}_I^{{\bf
t}}(S)_1)$ is compact. The desired conclusion easily follows from
the fact that each $\widehat{\cal M}^{\bf t}_I(S)_1$ is a manifold
of dimension zero.

Similarly, we can prove that $\cup_{I\in{\mathcal N}}
\hat\pi_I(\partial\widehat{\mathcal M}_I^{{\bf t}}(S)_0)$ is a
closed subset in $\cup_{I\in{\mathcal N}}
\hat\pi_I(\widehat{\mathcal M}_I^{{\bf t}}(S))$. Thus it is compact
and hence a finite set. \hfill$\Box$\vspace{2mm}

In contrast with Assumption~\ref{ass:4.7} in Framework I, {\bf from
now on} we make:

\begin{assumption}\label{ass:4.17}
{\rm  The union $X_{01}$  of $X_0$ and $X_1$  is contained in the
regular part $X^\circ$ of $X$. }
\end{assumption}

In contrast to (\ref{e:4.14}), for any ${\bf t}\in
B_\varepsilon(\R^m)_{res}$
 we formally write
\begin{equation}\label{e:4.36}
 e(E, S)^{\bf t}:=``\sum_{I\in{\cal
N}}"\frac{1}{|\Gamma_I|}\bigl\{\hat\pi_I:\widehat{\cal M}^{\bf
t}_I(S)_{01}\to X\bigr\}
\end{equation}
where $``\sum_{I\in{\cal N}}"$ is understand as in (\ref{e:2.111}).
It is a rational singular chain in $X$ of dimension $r={\rm
Ind}(S)$, called a {\bf virtual Euler chain} of the triple $(X, E,
S)$.

 As in Remark~\ref{rem:2.43}(iii) let
$\sim$ be the equivalence relation in the disjoint union
$\coprod_{I\in{\cal N}}\widehat{\mathcal M}_I^{{\bf t}}(S)_{01}$
generated by setting $\hat y_I\sim\hat z_J$ if $J\subset I$ and
$\hat\pi^I_J(\tilde y_I)=\hat z_J$. Denote by
\begin{equation}\label{e:4.40}
\overline{\mathcal M}^{{\bf t}}(S)_{01}:= \coprod_{I\in{\cal
N}}\widehat{\mathcal M}_I^{{\bf t}}(S)_{01}/\sim
\end{equation}
and by $\hat q_I:\widehat{\mathcal M}_I^{{\bf t}}(S)_{01}\to
\overline{\mathcal M}^{{\bf t}}(S)_{01}$ the restriction of the
obvious quotient map to $\widehat{\mathcal M}_I^{{\bf t}}(S)_{01}$.
Since $\widehat{\mathcal M}_I^{{\bf
t}}(S)_{01}\subset\widehat{\mathcal M}_I^{{\bf t}}(S)^\circ$ we have
a well-defined label function
\begin{equation}\label{e:4.41}
\lambda: \overline{\mathcal M}^{{\bf t}}(S)_{01}\to\Q,\;\bar
x\mapsto\frac{|\hat q_I^{-1}(\bar x)|}{|\Gamma_I|}\quad{\rm if}\;
\bar x\in {\rm Im}(\hat q_I),
\end{equation}
where $|\hat q_I^{-1}(\bar x)|$ is the number of elements in the set
$\hat q_I^{-1}(\bar x)$.

We hope that $\widehat{\cal M}^{\bf t}_I(S)_{01}$ is a smooth
manifold with boundary $\widehat{\cal M}^{\bf t}_I(S)_1$. In the
case that the bundle $(X, E, S)$ is also smoothly (weakly) $0$-good,
the final claim in Corollary~\ref{cor:4.5} tells us that each
$\widehat{\cal M}^{\bf t}_I(S)$ is a smooth manifold of dimension
$r$ and with corner, and each stratum $\widehat{\cal M}^{\bf
t}_I(S)_i$ is a submanifold.

\begin{proposition}\label{prop:4.18}
Under Assumption~\ref{ass:4.17}, if the bundle $(X, E, S)$ is also
smoothly (weakly) $0$-good, then
 each $\widehat{\cal M}^{\bf
t}_I(S)_{01}$ is a $r$-dimensional oriented manifold with boundary
$\partial\widehat{\cal M}^{\bf t}_I(S)_{0}\cup\widehat{\cal M}^{\bf
t}_I(S)_1$, i.e.
\begin{equation}\label{e:4.39}
\partial\widehat{\cal M}^{\bf t}_I(S)_{01}=\partial\widehat{\cal M}^{\bf t}_I(S)_{0}\cup\widehat{\cal M}^{\bf
t}_I(S)_1.
\end{equation}
In particular, if $X_0$ has no boundary as a Banach orbifold, then
$\partial\widehat{\cal M}^{\bf t}_I(S)_{01}=\widehat{\cal M}^{\bf
t}_I(S)_1$ because $\partial\widehat{\cal M}^{\bf
t}_I(S)_{0}=\emptyset$.
\end{proposition}

\noindent{\bf Proof}.\quad By (\ref{e:4.35}),
$\partial_t\widehat{\cal M}^{\bf
t}_I(S)_{01}\subseteq\partial\widehat{\cal M}^{\bf
t}_I(S)_{0}\cup\widehat{\cal M}^{\bf t}_I(S)_1\cup\widehat{\cal
M}^{\bf t}_I(S)_{\ge 2}$. In the present case we conclude
 \begin{equation}\label{e:4.40}
\partial_t\widehat{\cal M}^{\bf t}_I(S)_{01}=
\partial\widehat{\cal M}^{\bf t}_I(S)_{0}\cup\widehat{\cal M}^{\bf t}_I(S)_1\cup\widehat{\cal
M}^{\bf t}_I(S)_{\ge 2}.
\end{equation}
In fact, for any $x\in\widehat{\cal M}^{\bf
t}_I(S)_1\cup\widehat{\cal M}^{\bf t}_I(S)_{\ge 2}$, if there exists
a small open neighborhood ${\cal U}(x)$ of $x$ in $\widehat
V_I^\ast$ such that ${\cal U}(x)\cap\widehat{\cal M}^{\bf
t}_I(S)_{0}=\emptyset$, then the intersection ${\cal
U}(x)\cap\widehat{\cal M}^{\bf t}_I(S)={\cal U}(x)\cap\widehat{\cal
M}^{\bf t}_I(S)_{\ge 1}$ is a weakly PS manifold of dimension less
than $r$ and with finitely many strata. This contradicts to the fact
that ${\cal U}(x)\cap\widehat{\cal M}^{\bf t}_I(S)$ is a nonempty
open subset of the $r$-dimensional manifold $\widehat{\cal M}^{\bf
t}_I(S)$. Hence $\widehat{\cal M}^{\bf t}_I(S)_1\cup\widehat{\cal
M}^{\bf t}_I(S)_{\ge 2}\subseteq\partial_t\widehat{\cal M}^{\bf
t}_I(S)_{01}$ and thus (\ref{e:4.40}) holds true.

 Recall that a boundary manifold (or manifold boundary) of
the manifold $\widehat{\cal M}^{\bf t}_I(S)_{01}$ is a submanifold
of codimension one on which each point is a topological boundary
point. Hence (\ref{e:4.39}) follows immediately.
\hfill$\Box$\vspace{2mm}

Under the assumption of Proposition~\ref{prop:4.18}, the boundary of
$e(E, S)^{\bf t}$ is defined as
\begin{equation}\label{e:4.41}
\partial e(E, S)^{\bf t}:=``\sum_{I\in{\cal
N}}"\frac{1}{|\Gamma_I|}\bigl\{\hat\pi_I: \partial\widehat{\cal
M}^{\bf t}_I(S)_{01}\to X\bigr\}.
\end{equation}

Corresponding with Claim~\ref{cl:4.8} we still have:

\begin{claim}\label{cl:4.19}
Under the assumptions of Proposition~\ref{prop:4.18},  then
$\partial e(E, S)^{\bf t}$ is a rational singular cycle in $X$ and
its homology class $[\partial e(E, S)^{\bf t}]\in H_{r-1}(X,\Q)$ is
zero.
\end{claim}

 If $r=1$, by Proposition~\ref{prop:4.18}
 each $\widehat{\cal M}^{\bf t}_I(S)=\widehat{\cal
M}^{\bf t}_I(S)_{01}$ is a $1$-dimensional oriented manifold with
boundary $\partial\widehat{\cal M}^{\bf t}_I(S)_{0}\cup\widehat{\cal
M}^{\bf t}_I(S)_1$. By the classification theorem of the
1-dimensional manifold each connected component of $\widehat{\cal
M}^{\bf t}_I(S)$ is diffeomorphic to either the circle $S^1$ or one
of three intervals $[0, 1]$, $[0, 1)$ and $(0, 1)$. For each $\hat
x_I\in\partial\widehat{\cal M}^{\bf t}_I(S)$, as usual the
orientation of $\widehat{\cal M}^{\bf t}_I(S)$ induces an one of
$\hat x_I$, i.e., a number $o(\hat x_I)\in\{1, -1\}$.
 More precisely, since there must exist
  an embedding $\varphi:(0, 1)\to {\rm Int}(\widehat{\cal M}^{\bf t}_I(S))$ such that
$\varphi(1/n)\to\hat x_I$ as $n\to\infty$, $o(\hat x_I)$ is defined
as $1$ (resp. $-1$) if $\varphi$ is orientation preserving (resp.
reversing). From this definition and (\ref{e:4.34}) it immediately
follows that
\begin{equation}\label{e:4.42}
o(\hat\pi^I_J(\hat x_I))=o(\hat x_I)
\end{equation}
 for any $J\subset I\in{\cal N}$ and $\hat x_I\in (\hat
\pi^I_J)^{-1}(\widehat V_J^\ast)\cap\partial\widehat{\cal M}^{\bf
t}_I(S)$. Moreover, under Assumption~\ref{ass:4.17}, if $r=1$ then
 for any $J\subset I\in{\cal N}$ the projection
$$
 \hat\pi^I_J: (\hat \pi^I_J)^{-1}(\widehat
V_J^\ast)\cap\widehat{\cal M}^{\bf t}_I(S)\to {\rm
Im}(\hat\pi^I_J)\cap\widehat{\cal M}^{\bf t}_J(S) $$
 is a $|\Gamma_I|/|\Gamma_J|$-fold
(regular) PS covering. So  using (\ref{e:4.34}) and (\ref{e:4.42})
 we can derive that for any $J\subset I\in{\cal N}$,
\begin{equation}\label{e:4.43}
 \sum_{\hat x\in (\hat \pi^I_J)^{-1}(\widehat
V_J^\ast)\cap\partial\widehat{\cal M}^{\bf
t}_I(S)}\frac{1}{|\Gamma_I|}o(\hat x)=\sum_{\hat y\in {\rm
Im}(\hat\pi^I_J)\cap\partial\widehat{\cal M}^{\bf
t}_J(S)}\frac{1}{|\Gamma_J|}o(\hat y).
\end{equation}
We call the rational number
\begin{equation}\label{e:4.44}
 \sharp \partial e(E, S)^{\bf t}:=``\sum_{I\in{\cal
N}}" \sum_{\hat x\in\partial\widehat{\cal M}^{\bf t}_I(S)}
\frac{1}{|\Gamma_I|}o(\hat x)
\end{equation}
the {\bf oriented number} of $\partial e(E, S)^{\bf t}$, where as
before $``\sum_{I\in{\cal N}}"$ means that for the parts
$$
\sum_{\hat x\in (\hat \pi^I_J)^{-1}(\widehat
V_J^\ast)\cap\partial\widehat{\cal M}^{\bf
t}_I(S)}\frac{1}{|\Gamma_I|}o(\hat x)\quad{\rm in}\quad
\sum_{\tilde x\in\partial\widehat{\cal M}^{\bf t}_I(S)}
\frac{1}{|\Gamma_I|}o(\hat x)
$$
and
$$\sum_{\hat y\in {\rm
Im}(\hat\pi^I_J)\cap\partial\widehat{\cal M}^{\bf
t}_J(S)}\frac{1}{|\Gamma_J|}o(\hat y)\quad{\rm in}\quad
\sum_{\tilde y\in\partial\widehat{\cal M}^{\bf t}_J(S)}
\frac{1}{|\Gamma_J|}o(\tilde y)
$$
we only count them once because they are same by (\ref{e:4.43}).

By (\ref{e:4.42}), for  $\bar x\in\partial\overline{\mathcal
M}^{{\bf t}}(S)=\partial\overline{\mathcal M}^{{\bf t}}(S)_{01}$ we
can define
\begin{equation}\label{e:4.45}
o(\bar x):=o(\hat x)
\end{equation}
with any representative of the equivalence class $\bar x$, and then
using the label function $\lambda$ in (\ref{e:4.41}) define the {\bf
oriented number} of $\partial\overline{\mathcal M}^{{\bf t}}(S)$ by
\begin{equation}\label{e:4.46}
\sharp\partial\overline{\mathcal M}^{{\bf t}}(S):=\sum_{\bar
x\in\partial\overline{\mathcal M}^{{\bf t}}(S)}\lambda(\bar
x)o(\bar x).
\end{equation}
One easily proves that
\begin{equation}\label{e:4.47}
\sharp \partial e(E, S)^{\bf t}=\sharp\partial\overline{\mathcal
M}^{{\bf t}}(S).
\end{equation}
It follows from Claim~\ref{cl:4.19} that

\begin{theorem}\label{th:4.20}
Under the assumptions of Proposition~\ref{prop:4.18}, if $r=1$ then
the oriented number $\sharp\partial e(E, S)^{\bf t}=0$
 for any  ${\bf t}\in B_\varepsilon(\R^m)_{res}$.
\end{theorem}

In the Floer homology, though the case $r=0$ is sufficient for
defining the boundary map,  the case that $r=1$ and $X_0$ has no
boundary as a Banach orbifold will be used when one proves the
boundary map to be indeed a boundary operator. The case that $r=1$
and $X_0$ has a boundary is also needed for one proving the
associated chain homomorphisms with different Hamiltonian homotopies
are chain homotopy equivalent.

Finally,  after replaceing Assumption~\ref{ass:4.7} by
Assumption~\ref{ass:4.17} (if necessary) in Prop.\ref{prop:4.10}
(the restriction principle), Prop.\ref{prop:4.12} (stability),
Prop.\ref{prop:4.13} (homotopy), Prop.\ref{prop:4.14} and
Prop.\ref{prop:4.15} (functoriality) all corresponding conclusions
about the virtual Euler chains still hold.

The deep studies of different product operators in Floer theory
require us to consider the case $r>1$. For a smoothly (weakly)
$0$-good PS Banach $0$-Fredholm orbibundle, our theory and
techniques can easily be further developed to satisfy different
needs. We here do not present them.

\appendix
\section{Appendix:\quad Banach manifolds with
corners}\label{app:A}

For convenience of the readers we review some basic materials on
the Banach manifolds, cf. \cite{MaOu}. Let $B_1,\cdots, B_n$ be
real Banach spaces. For $(x_1,\cdots,x_n)\in B_1\times\cdots\times
B_n$,
$$
\|(x_1,\cdots,
 x_n)\|=\max\{\|x_1\|,\cdots,\|x_n\|\}
 ,$$
 define one of three typical equivalent norms on the real linear space
$B_1\times\cdots\times B_n$.
 Their topology is the product topology of
the spaces $B_1,\cdots,B_n$. For a  real Banach space $B$ denote
by  ${\mathcal L}(B,\R)$  the space of all continuous real linear
functionals on $B$. Let $\Lambda=\{\lambda_1,\cdots,\lambda_n\}$
be a linearly independent system of elements of ${\mathcal
L}(B,\R)$. The set
$$B^+_\Lambda:=\{x\in B\,|\, \lambda_1(x)\ge
0,\cdots,\lambda_n(x)\ge 0\}
$$
is called the $\Lambda$-{\bf quadrant of order $n$ } of $B$. If
$\Lambda=\{\lambda\}$ we also denote $B^+_\Lambda$ by $B^+_\lambda$.
Then it is clear that $B^+_\Lambda=\cap^n_{i=1}E^+_{\lambda_i}$ and
$B^+_\Lambda=B$ if $\Lambda=\emptyset$. We also denote
$\cap^n_{i=1}{\rm Ker}(\lambda_i)$ by $B^0_\Lambda$ or $B^0_\lambda$
if $\Lambda=\{\lambda\}$. Then we have $B^0_\Lambda=\cap^n_{i=1}
B^0_{\lambda_i}$ and $B^0_\Lambda=B$ if $\Lambda=\emptyset$.

\begin{definition}\label{def:A.1}
{\rm Let $B$ and $F$ be two real Banach spaces, $B^+_\Lambda$ a
$\Lambda$-quadrant of $B$ and $U$ an open set of $B^+_\Lambda$.
For $p\in\N$ we say that a map $f:U\to F$ is of class $p$ on $U$
if $f$ is $p$ times differentiable on $U$  and the map $D^pf: U\to
{\mathcal L}^p_s(B,F)$ is continuous. If $f$ is of class $p$ on
$U$ for all $p\in\N$ we say that $f$ is of class $\infty$ on $U$.}
\end{definition}

\begin{proposition}\label{prop:A.2}
Let $f:U\to F$ be as in Definition~\ref{def:A.1}, and
$p\in\N\cup{\infty}$.
\begin{enumerate}
\item[(i)] If the real Banach space $B$ admits differentiable
partitions of unity of class $p$, then $f$ is of class $p$ on $U$
if and only if for any $x\in U$ there exist an open neighborhood
$V(x)$ of $x$ in $B$ and a map $f_x:V(x)\to F$ of class $p$ (in
the sense of the ordinary Differential Calculus) such that $f_x$
and $f$ coincide on $V(x)\cap U$.

\item[(ii)] If $p\in\N$, then $f$ is of class $p$ on $U$ if and
only if there exist an open set $G$ of $B$ and a map $\bar f:G\to
F$ of class $p$ (in the sense of the ordinary Differential
Calculus) such that $G\cap B^+_\Lambda=U$ and $f=\bar f$ on $U$.
\end{enumerate}
\end{proposition}

Let $X$ be a Hausdorff topological space. A  {\bf Banach chart} on
$X$ is a $4$-tuple $(U,\varphi, B,\Lambda)$ consisting of an open
subset $U$  of $X$,  a real Banach space $B$, a linearly independent
system $\Lambda$ of elements of ${\mathcal L}(B,\R)$ and a
homeomorphism $\varphi$ from $U$ to an open subset $\varphi(U)$ of
$B^+_\Lambda$. The space $B$ is called the {\bf model} of the chart,
and $U$ the {\bf domain} of the chart.  Two Banach charts
$(U,\varphi, B,\Lambda)$ and $(U^\prime,\varphi^\prime,
B^\prime,\Lambda^\prime)$ on $X$ are said to be {\bf compatible } of
class $p$ if $\varphi(U\cap U^\prime)$ and $\varphi^\prime(U\cap
U^\prime)$ are open subsets in $B^+_\Lambda$ and $B^{\prime
+}_{\Lambda^\prime}$ respectively, and
$$
\varphi^\prime\circ\varphi^{-1}:\varphi(U\cap U^\prime)\to
\varphi^\prime(U\cap U^\prime)\quad{\rm and}\quad
\varphi\circ\varphi^{\prime-1}:\varphi^\prime(U\cap U^\prime)\to
\varphi(U\cap U^\prime)
$$
are of class $p$ (and hence $C^p$-diffeomorphisms).

\begin{definition}\label{def:A.3}
{\rm A collection ${\mathcal A}$ of charts on a Hausdorff space
$X$ is called an {\bf atlas} of class $p$ on $X$ if the domains of
the charts of ${\mathcal A}$ cover $X$ and any two of them are
compatible of class $p$. Such two atlases ${\mathcal A}$ and
${\mathcal A}^\prime$ on $X$ is said to be {\bf equivalent} of
class $p$ if ${\mathcal A}\cup{\mathcal A}^\prime$ is also an
atlas of class $p$ on $X$. It may be proved that this is indeed a
equivalence relation over the atlases of class $p$ on $X$ for each
$p\in \N\cup\{\infty\}$. Denote by $[{\mathcal A}]$ the equivalent
class of the atlas of class $p$ on $X$. It is called {\bf
differentiable structure} of class $p$ on $X$ and the pair $(X,
[{\mathcal A}])$ is called $C^p$-Banach manifold, usually only
denoted by $X$.}
\end{definition}

It should be remarked that for two charts  $(U,\varphi,
B,\Lambda)$ and $(U^\prime,\varphi^\prime,
B^\prime,\Lambda^\prime)$  of a Banach manifold $X$ at $x$ the
Banach spaces $B$ and $B^\prime$ are not linearly homeomorphic.

\begin{definition}\label{def:A.4}
{\rm Let $B^+_\Lambda$ be a $\Lambda$-quadrant of order $n$ in the
real Banach space $B$ as above, and $U$ an open set of
$B^+_\Lambda$. For $x\in B^+_\Lambda$, ${\rm
ind}_\Lambda(x):=\sharp(\{i\,|\,\lambda_i(x)=0\}$ is called
$\Lambda$-{\bf index} of $x$. The sets $\partial_\Lambda U:=\{x\in
U\,|\, {\rm ind}_\Lambda(x)\ge 1\}$ and ${\rm int}_\Lambda
U:=\{x\in U\,|\, {\rm ind}_\Lambda(x)=0\}$ are called
$\Lambda$-{\bf boundary} of $U$ and  $\Lambda$-{\bf interior} of
$U$, respectively.}
\end{definition}

\begin{proposition}\label{prop:A.5}
 Let $B^+_\Lambda$ be a $\Lambda$-quadrant of order $n$ in the
real Banach space $B$ as above, and $U$ an open set of
$B^+_\Lambda$. Then we have:
\begin{description}
\item[(i)] $U-\partial_\Lambda U={\rm int}_\Lambda U=U\cap{\rm
int}_\Lambda B^+_\Lambda$;

\item[(ii)] $U=(\partial_\Lambda U)\cup{\rm int}_\Lambda U$ and
$(\partial_\Lambda U)\cap{\rm int}_\Lambda U\ne\emptyset$;

\item[(iii)] If $\Lambda=\emptyset$ then $\partial_\Lambda
U=\emptyset$ and ${\rm int}_\Lambda U=U$;

\item [(iv)] ${\rm int}_\Lambda U$ is the topological interior of
$U$ in $B$.

\item[(v)] If $\partial_\Lambda U\ne\emptyset$ then
$\Lambda\ne\emptyset$ and ${\rm int}_\Lambda U\ne\emptyset$.
\end{description}
Moreover, if $F^+_\gamma$ is a $\Gamma$-quadrant of order $m$ in
another Banach space $F$,  $V$ an open set of $F^+_\gamma$ and
$f:U\to V$ a diffeomorphism of class $p\ge 1$, then
\begin{description}
\item[(vi)] ${\rm ind}_\Lambda(x)={\rm ind}_\Gamma(f(x))$ for all
$x\in U$;

\item[(vii)] $\partial_\Lambda U\ne\emptyset$ if and only if
$\partial_\Gamma V\ne\emptyset$;

\item[(viii)] $f({\rm int}_\Lambda U)={\rm int}_\Gamma V$ and
$f(\partial_\Lambda U)=\partial_\gamma v$;

\item[(ix)] $f|_{{\rm int}_\Lambda U}:{\rm int}_\Lambda U\to {\rm
int}_\Gamma V$ is a diffeomorphism of class $p$.
\end{description}
\end{proposition}

It follows that for any Banach manifold $X$ of class $p\ge 1$, if
$x\in X$ is contained in the domains of two Banach charts
$(U,\varphi, B,\Lambda)$ and $(U^\prime,\varphi^\prime,
B^\prime,\Lambda^\prime)$ of $X$ one  has ${\rm
ind}_\Lambda(\varphi(x))={\rm
ind}_{\Lambda^\prime}(\varphi^\prime(x))$. Thus we may define
${\rm ind}_\Lambda(\varphi(x))$ to be {\bf index} of $x$, denoted
by ${\rm ind}(x)$, where $(U,\varphi, B,\Lambda)$ is any Banach
chart of $X$ whose domain contains $x$. One can prove that for any
$x\in X$ there exists a Banach chart $(U,\varphi, B,\Lambda)$ of
$X$ such that $x\in U$ and $\varphi(x)=0$. Such a chart is called
{\bf centered chart} at $x$.

\begin{definition}\label{def:A.6}
{\rm Let $X$ be a Banach manifold of class $p\ge 1$. For any
$k\in\N\cup\{0\}$ we define
$$\partial^k X:=\{x\in X\,|\,{\rm ind}(x)\ge k\}\quad{\rm
and}\quad B_k X=\{x\in X\,|\,{\rm ind}(x)=k\}.
$$
The former is called $k$-{\bf boundary} of $X$. In particular
$\partial^1 X$ is called {\bf boundary} of $X$ and denoted by
$\partial X$. Clearly, $\partial^0 X=X$. $B_0 X$ is called {\bf
interior} of $X$ and denoted by ${\rm Int}(X)$.}
\end{definition}

\begin{proposition}\label{prop:A.7}
Let $X$ be a Banach manifold of class $p\ge 1$. Then we have:
\begin{description}
\item[(a)] ${\rm Int}(X)$ is a dense open set of $X$.

\item[(b)] $\partial^k X$ is a closed set of $X$ for  any
$k\in\N$.

\item[(c)] For every $k\in\N\cup\{0\}$ there exists a unique
differentiable structure of class $p$ on $B_k X$ such that for all
$x\in B_k X$ and all charts $(U,\varphi, B,\Lambda)$ of $X$ centered
at $x$, the triplet $(U\cap B_k X,\varphi|_{U\cap B_k X},
B^0_\Lambda)$ is a Banach chart of that structure. Furthermore $B_k$
has no boundary, i.e. $\partial(B_k X)=\emptyset$, and the topology
of the manifold $B_k X$ is the topology induced by $X$.

\item[(d)] There is a unique differentiable structure of class $p$
on ${\rm Int}(X)$ such that for all $x\in{\rm Int}(X)$ and all
charts $(U,\varphi, B,\Lambda)$ of $X$ centered at $x$ (hence
$\Lambda=\emptyset$, $U\subset{\rm Int}(X)$ and $B^0_\Lambda=B$) the
triplet $(U,\varphi, B)$ is a chart of ${\rm Int}(X)$. Furthermore
${\rm Int}(X)$ has no boundary and its topology is the topology
induced by $X$.

\item[(e)] If $\partial^2 X=\emptyset$, there is a unique
differentiable structure of class $p$ on $\partial X$ such that
for all $x\in\partial X$ and all charts $(U,\varphi, B,\Lambda)$
of $X$ centered at $x$ (hence $\Lambda=\{\lambda\}$), it happens
that $(U\cap\partial X,\varphi|_{U\cap\partial X}, B^0_\lambda)$
is a chart of $\partial X$. Furthermore $\partial X$ has boundary
and its topology is the topology induced by $X$.
\end{description}
\end{proposition}

\begin{proposition}\label{prop:A.8}
Let $X$ and $X^\prime$ be two Banach manifolds of class $p\ge 1$,
and $f:X\to X^\prime$ a diffeomorphism of class $p$. Then we have:
\begin{description}
\item[(1)] ${\rm ind}(x)={\rm ind}(f(x))$ for all $x\in X$;

\item[(2)] $f(\partial^k X)=\partial^k X^\prime$ and $f(B_k X)=B_k
X^\prime$ for all $k\in\N\cap\{0\}$;

\item[(3)] For all $k\in\N\cap\{0\}$, $f|_{B_k X}: B_k X\to B_k
X^\prime$ is a diffeomorphism of class $p$, where $B_k X$ and $B_k
X^\prime$ are the manifolds described in
Proposition~\ref{prop:A.7}(c). In particular, if $\partial^2
X=\emptyset$, $f$ is a diffeomorphism of class $p$ of $\partial
X=B_1 X$ onto $\partial X^\prime=B_1 X^\prime$.
\end{description}
\end{proposition}

\begin{proposition}\label{prop:A.9}
Let $X$ and $Y$ be two Banach manifolds of class $p$. Then there
is a unique structure of differentiable manifold of class $p$,
$[{\mathcal A}]$, in $X\times Y$ such that for every chart
$(U,\varphi,B,\Lambda)$ of $X$ and every chart $(V,\psi,
F,\Gamma)$, $(U\times V,\varphi\times\psi, B\times F, \Lambda\circ
p_1\cup\Gamma\circ p_2))$ is a chart of $(X\times Y,[{\mathcal
A}])$. The pair $(X\times Y,[{\mathcal A}])$ is called the {\bf
product manifold} of $X$ and $Y$. Then we have:
\begin{description}
\item[(1)] The topology of the product manifold $X\times Y$ is the
product of the topologies of $X$ and $Y$;

\item[(2)] For every $(x,y)\in X\times Y$, ${\rm ind}(x,y)={\rm
ind}(x)+{\rm ind}(y)$;

\item[(3)] $\partial(X\times Y)=(\partial x\times
Y)\cup(X\times\partial Y)$, more generally, for all
$k\in\N\cap\{0\}$,
$$\partial^k(X\times Y)=\bigcup_{{i+j=k}\atop
{i\ge 0, j\ge 0}}(\partial^i x\times\partial ^j Y); $$

\item[(4)] ${\rm Int}(X\times Y)={\rm Int}(X)\times{\rm Int}(Y)$,
more generally, for all $k\in\N\cap\{0\}$,
$$B_k(X\times
Y)=\bigcup_{{i+j=k}\atop {i\ge 0, j\ge 0}}(B_i x\times B_j Y); $$
and these $B_i X\times B_j Y$ are pairwise disjoint open sets of
$B_k(X\times Y)$.
\end{description}
\end{proposition}

\begin{theorem}\label{th:A.10}
({\bf Local Inverse Mapping Theorem for Open Sets with Corners}).
Let $U$ be an open subset of a quadrant $B^+_\Lambda$,
$F^+_\Gamma$ a quadrant of a Banach space $F$, $f:U\to F^+_\Gamma$
a map of class $p$ ($p\in\N\cup\{\infty\}$) such that
$f(\partial_\Lambda(U))\subset\partial F^+_\Gamma$. For $x\in U$,
suppose that $DF(x): B\to F$ is a linear homeomorphism. Then there
exist an open neighbourhood $U_1(x)$ of $x$ in $U$ and an open
neighborhood $U^\prime$ of $f(x)$ in $F^+_\Gamma$ such that
$f|_{U_1(x)}$ is a diffeomorphism of class $p$ from $U_1(x)$ onto
$U^\prime$.
\end{theorem}

\begin{theorem}\label{th:A.11}
({\bf The Implicit Function Theorem for Manifolds  with Corners}).
Let $X, Y, Z$ be differentiable manifolds of class $p$
($p\in\N\cup\{\infty\}$), $f:X\times Y\to Z$ a map of class $p$
and $(a,b)\in X\times Y$. Suppose that $D_2f(a,b): T_bY\to
T_{f(a,b)}Z$ is a linear homeomorphism and suppose that there are
an open neighborhood $V(a)$ of $a$ and an open neighborhood $V(b)$
of $b$ such that $f(V(a)\times(V(b)\cap\partial y))\subset\partial
Z$.

Then there are an open neighborhood $W(a)$ of $a$, an open
neighborhood $W(b)$ of $b$ and a unique map $g:W(a)\to W(b)$ such
that:
\begin{description}
\item[(1)] $f(x, g(x))=f(a,b)$ for $x\in W(a)$;

\item[(2)] $g(a)=b$ and $g$ is of class $p$ on $W(a)$;

\item[(3)] For every $x\in W(a)$, $D_2f(x,g(x))$ is a linear
homeomorphism from $T_x W(a)=T_xX$ to $T_{g(x)}W(b)=T_{g(x)}Y$,
and $Dg(x)=-(D_2f(x,g(x)))^{-1}\circ D_1f(x,g(x))$.
\end{description}
\end{theorem}

\begin{proposition}\label{prop:A.12}({\bf Topological properties of the differential
manifolds}). Let $X$ be a Banach manifold. Then
\begin{description}
\item[(1)]$X$ is locally connected and locally path-connected;

\item[(2)] $X$ satisfies the first axiom of countability;

\item[(3)] $X$ satisfies the second axiom of countability if and
only if one of the following two conditions hold:
\begin{enumerate}
\item[$(3.1)$] $X$ satisfies the Lindel\"of property and has an
atlas whose charts are modelled over separable real Banach spaces;

\item[$(3.2)$] $X$ has a countable atlas whose charts are modelled
over separable real Banach spaces;
\end{enumerate}

\item[(4)] $X$ is locally compact if and only if it is locally of
finite dimension;

\item[(5)] $X$ is a Baire space;

\item[(6)] $X$ is a metrizable space if and only if $X$ is a
paracompact space.

\end{description}
\end{proposition}

Let $X$ and $Y$ be $C^p$ Banach manifolds with $p\ge 1$. A
$C^1$-map $f:X\to Y$ is called {\bf Fredholm} at $x\in X$ if
$df(x):T_xX\to T_{f(x)}Y$ is a Fredholm operator. In this case the
Fredholm index of $f$ at $x$, ${\rm ind}_x(f)$, is defined as the
Fredholm index of $df(x)$. $f$ is said to be a Fredholm map if it
is Fredholm at each point of $X$. A point $x\in X$ is called {\bf
singular} if $df(x)$ is not surjective. The set of singular points
of $f$ is denoted by $S(f)$. Each point in $X\setminus S(f)$ is
called a {\bf regular point} of $f$. A point in $Y\setminus
f(S(f))$ is called a {\bf regular value} of $f$.

\begin{theorem}\label{th:A.13} ({\bf Smale and Quinn Density
Theorem}). Let $X$ and $Y$  be $C^p$ Banach manifolds. For a $C^p$
Fredholm map  $f:X\to Y$ such that $p\ge 1+{\rm ind}_x(f)$ for
every $x\in X$ one has:\\
({\bf Smale}) If $X$ and $Y$ are $T_3$ topological spaces and $X$
is Lindel\"of then $Y\setminus f(S(f))$ is a residual set in $Y$
and therefore it is dense in $Y$.\\
({\bf Quinn}) If $X$ is a $T_3$ space and $f$ is also
$\sigma$-proper then $Y\setminus f(S(f))$ is a residual set in $Y$
and therefore it is dense in $Y$. Moreover, if $f$ is proper then
$Y\setminus f(S(f))$ is open and dense in $Y$.
\end{theorem}

\begin{proposition}\label{prop:A.14}
Let $X$ and $Y$ be separable Banach manifolds, and $Z$ be a Banach
space. Suppose that $\Phi:X\times Y\to Z$ is transversal to a point
$q\in Z$. Then ${\cal M}(\Phi, q):=\{(x,y)\in X\times Y\,|\,
\Phi(x,y)=q\}$ is a separable Banach manifold and $y\in Y$ is a
regular value of the restriction of the projection ${\bf P}_Y$ (to
the second factor $Y$) to ${\cal M}(\Phi,q)$  if and only if the map
$$\Phi_y: X\to Z,\; x\mapsto \Phi(x,y)$$
is transversal to $q$. In particular, if $\Phi$ is equal to $f\oplus
g:X\times Y\to Z,\;(x,y)\mapsto f(x)+ g(y)$, where $f:X\to Z$ is a
$C^l$ Fredholm map of and $g:Y\to Z$ a $C^l$-smooth map, then the
restriction of the projection ${\bf P}_Y$ to ${\cal M}(\Phi,q)$ is a
Fredholm map of index ${\rm Ind}(f)$.
\end{proposition}

\noindent{\bf Proof.}\quad The first conclusion is standard.  To
arrive at the second note that $y\notin {\bf P}_Y({\cal M}(\Phi,q))$
if and only if $\Phi(x,y)\ne q$ for any $x\in X$. When $y\in {\bf
P}_Y({\cal M}(\Phi,q))$ and $x\in(\Phi_y)^{-1}(q)$, i.e.
$\Phi(x,y)=q$, the differential
$$
d\Bigl({\bf P}_Y\Bigm|_{{\cal M}(\Phi,q)}\Bigr)(x,y): T_{(x,y)}{\cal
M}(\Phi,q)\to T_yY,\,(\xi,\eta)\mapsto\eta
$$
is surjective if and only if for each $\eta\in T_yY$ there exists
$\xi\in T_xX$ such that $(\xi,\eta)\in{\rm Ker}(\Phi(x,y))$, i.e.
$d_x\Phi(x,y)(\xi)+ d_y\Phi(x,y)(\eta)=d\Phi_y(x)(\xi)+
d_y\Phi(x,y)(\eta)=0$. If $\Phi_y$ is transversal to $q$ at this
$x$, then $d\Phi_y(x):T_xX\to $ is onto, and thus for any $\eta\in
T_yY$ one has $\xi\in T_xX$ satisfying
$d\Phi_y(x)(\xi)=-d_y\Phi(x,y)(\eta)$. Conversely, for any $z\in Z$,
since $d\Phi(x,y)$ is surjective we have $(\zeta_1,\eta)\in
T_xX\times T_yY$ such that $d_x\Phi(x,y)(\zeta_1)+
d_y\Phi(x,y)(\eta)=z$. Moreover, $d\Bigl({\bf P}_Y\Bigm|_{{\cal
M}(\Phi,q)}\Bigr)(x,y)$ is surjective one has also $\zeta_2\in T_xX$
such that $(\zeta_2, -\eta)\in T_{(x,y)}{\cal M}(\Phi,q)$, i.e
$d_x\Phi(x,y)(\zeta_2)+ d_y\Phi(x,y)(-\eta)=0$. It follows that
$d_x\Phi(x,y)(\zeta_1+\zeta_2)=z$. That is, $d_x\Phi_y(x)$ is
surjective.

Finally, if $\Phi=f\oplus g$, then at any $(x,y)\in {\cal
M}(\Phi,q)$ the tangent space $T_{(x,y)}{\cal M}(\Phi,q)={\rm
Ker}(\Phi(x,y))=\{ (\xi,\eta)\in T_xX\times T_yY\,|\,df(x)(\xi)+
dg(y)(\eta)=0\}$. Since the  differential
$$d\Phi(x,y): T_{(x,y)}(X\times Y)\to T_qZ=Z, (\xi,\eta)\mapsto
df(x)(\xi)+ dg(y)(\eta)
$$
 is onto, and $df(x)$ is Fredholm,
 by Lemma~\ref{lem:1.3}(i)  the projection $P_2:T_xX\times T_yY\to T_xY$ to
${\rm Ker}(d\Phi(x,y))$ (which is equal to the restriction of the
differential $d{\bf P}_Y(x,y)$ to ${\rm Ker}(d\Phi(x,y))$) is
Fredholm and has index ${\rm Index}(df(x))$.
\hfill$\Box$\vspace{2mm}

\section{Appendix:\quad Slice theorem for actions of Lie groups on
Banach manifolds}\label{app:B}

This appendix is written in terms of Section 2 in Appendix.B of
\cite{GuGiKa} by  Guillemin-Ginzburg-Karshon. Let $X$ be a
connected Banach manifold, and a compact Lie group $G$ smoothly
acts on $X$ by automorphisms of $X$. The {\bf isotropy group} (or
{\bf stabilizer}) of a point $x\in X$ is
$$
G_x=\{g\in G\,|\, g\cdot x=x\},
$$
and the {\bf orbit} of a point $x\in X$ is
$$
G\cdot x=\{g\cdot x\,|\, g\in G\}.
$$
The {\bf orbit type} of a point $x\in X$ is defined to be the
conjugacy class of its isotropy group $G_x$ in $G$. Since $G_{g\cdot
x}=gG_xg^{-1}$ the conjugacy class is constant along an orbit. For
any subgroup $H\subset G$ let $(H)$ denote the conjugacy class of
$H$ in $G$. We set
\begin{eqnarray*}
&&X_H=\{x\in X\,|\, G_x=H\},\\
&&X_{(H)}=\{x\in X\,|\, G_x\,{\rm is}\,{\rm conjugate}\,{\rm to}\,
H\}.
\end{eqnarray*}
The fixed point set of a non-trivial subgroup $H\subset G$  is
defined by
$$X^H:=\{x\in X\,|\, g\cdot x=x\,\forall g\in H\}.$$

\begin{theorem}\label{th:B.1}
 For any $x\in X^G$ there exists a $G$-equivariant diffeomorphism
from a neighborhood of the origin in $T_xX$ onto a neighborhood of
$x$ in $X$.
\end{theorem}

\noindent{\bf Proof.} For any small $G$-invariant neighborhood $U$
of $x$ we can take a smooth map $f:U\to T_xX$ such that its
differential at $x$ is the identity map on $T_xX$. With the  Haar
measure $dg$ of $G$ we define the average of $f$ over $G$ as
$F:U\to T_xX$,
$$
F(y)=\int_Gg_\ast f(g^{-1}y)dg.
$$
(If $G$ is a finite group $F$ is defined by
$F(y)=\frac{1}{|G|}\sum_{g\in G} g_\ast f(g^{-1}\cdot y)$). Then
$dF(x)=1_{T_xX}$ and $F(h\cdot y)=h_\ast(F(y))$ for any $y\in U$
and $h\in G$. By the implicit function theorem there exists a
small $G$-invariant neighborhood $V\subset U$ of $x$ such that
$F|_V:V\to F(V)$ is a diffeomorphism. Its inverse is a desired
diffeomorphism.\hfill$\Box$\vspace{2mm}

Theorem~\ref{th:B.1} is often called the local linearization
theorem.

The tangent space of an orbit $G\cdot x$ at $x$, $T_x(G\cdot x)$, is
a $G_x$-invariant subspace of $T_xX$. Since $T_x(G\cdot x)$ is
finite dimensional,  elementary functional analysis arguments imply
that there exists a direct sum decomposition $T_xX=T_x(G\cdot
x)\oplus Y$. However $Y$ is not necessarily $G_x$-invariant. Note
that the direct sum decomposition $T_xX=T_x(G\cdot x)\oplus Y$
corresponds with bounded projection $P:T_xX\to T_x(G\cdot x)$ in a
one-to-one way. We can average a projection by the $G_x$-action to
get an $G_x$-invariant complement subspace $Y$ of $T_x(G\cdot x)$,
$T_xX=T_x(G\cdot x)\oplus Y$. Let $W$ be an open $G_x$-invariant
neighborhood of the origin in $Y$, and $G\times_{G_x}W$ be the
associated bundle.

\begin{theorem}\label{th:B.2}{\rm (Slice theorem)}
 For small $W$ there exists a $G$-equivariant diffeomorphism
 $\Phi$  from the bundle $G\times_{G_x}W$ onto a neighborhood of the orbit
 $G\cdot x$ in $X$ such that the restriction of it to the zero
 section $G\times_{G_x}\{0\}=G/G_x$ is given by
 $$
 G/G_x\to X,\quad gG_x\mapsto g\cdot x.
 $$
\end{theorem}

\noindent{\bf Proof.} By Theorem~\ref{th:B.1} there exists a
$G$-equivariant diffeomorphism $\phi$ from a neighborhood ${\cal
U}$ of the origin in $T_xX$ onto a $G$-invariant neighborhood
${\cal V}$ of $x$ in $X$. Shrink $W$ so that $W\subset {\cal U}$,
then the map
$$
\Phi:G\times_{G_x}W\to X,\quad [g, \xi]\mapsto g\cdot\phi(\xi),
$$
is well defined and is a local diffeomorphism at the point $[e,
0]$. Here $e$ is the unit in $G$. The $G$-equivariance implies
that $\Phi$ is a local diffeomorphism at all points of the form
$[e,0]$. Thus one only need to prove that $\Phi$ is one-to-one for
sufficiently small $W$. The latter can be easily obtained by a
contradiction argument because $G$ is compact.
\hfill$\Box$\vspace{2mm}

Using the slice theorem we easily derive:

\begin{corollary}\label{cor:B.3}
 For any subgroup $H\subset G$ the set $X_{(H)}$ is a submanifold
of $X$. If $H$ is nontrivial then its fixed point set $X^H$ is a
disjoint union of closed submanifolds of $X$. Moreover, if the
action of $G$ is effective, $X^\circ:=\{x\in X\,|\, G_x=\{1\}\}$
is an open and dense subset of $X$, and $X^{sing}:=\{x\in X\,|\,
G_x\ne\{1\}\}$ is equal to the union $\cup_{H\ne\{1\}}X^H$.
\end{corollary}

\begin{corollary}\label{cor:B.4}
If the action of $G$ on $X$ is locally free, i.e., each point
$x\in X$ has finite isotropy group $G_x$, then the quotient $X/G$
is a Banach orbifold.
\end{corollary}

\section{Appendix:\quad Three results on linear functional analysis}\label{app:C}

All Banach spaces are real unless otherwise stated.

\begin{proposition}\label{prop:C.1}
Let $A:X\to Y$ be a bounded linear Fredholm operator of index $r$
from Banach spaces $X$ to $Y$. Assume that $H\subset X$ is a
subspace of finite codimension $r'\ge 1$. Then the restriction
$A|_H: H\to Y$ is also a Fredholm operator and ${\rm
Index}(A|_H)=r-r'$.
\end{proposition}

\noindent{\bf Proof.} {\it Case 1}. {\bf ${\rm Ker}(A)\subset H$.}
We have direct sum decompositions $H={\rm Ker}(A)\oplus H_1$ and
$X=H\oplus X_1$. Then $A(X)=A(H_1)\oplus A(X_1)$. Since $\dim
Y/A(X)$ is finite there exists a direct sum decomposition
$Y=Y_1\oplus A(X)=Y_1\oplus A(H_1)\oplus A(X_1)$. Note that
$A|_H(H)=A(H_1)$. We have $Y/{\rm Im}(A|_H)\cong Y_1\oplus A(X_1)$
and thus
\begin{eqnarray*}
{\rm Index}(A|_H)\!\!\!\!\!\!&&=\dim{\rm Ker}(A|_H)-\dim Y/{\rm Im}(A|_H)\\
&&= \dim{\rm Ker}(A)-\dim (Y_1\oplus A(X_1))\\
&&= \dim{\rm Ker}(A)-\dim Y_1-\dim A(X_1)\\
&&=r-r'.
\end{eqnarray*}

{\it Case 2.} {\bf ${\rm Ker}(A)\cap H=\{0\}$.} We have direct sum
decompositions $X=H\oplus{\rm Ker}(A)\oplus H_2$ and $Y=Y_1\oplus
A(X)=Y_1\oplus A(H_2)\oplus A(H)$. So $Y/A(H)\cong Y_1\oplus
A(H_2)$. It follows that
\begin{eqnarray*}
{\rm Index}(A|_H)\!\!\!\!\!\!&&=\dim{\rm Ker}(A|_H)-\dim Y/{\rm Im}(A|_H)\\
&&= 0-\dim (Y_1\oplus A(H_2))\\
&&= 0-\dim{\rm Coker}(A)-\dim H_2\\
&&=-\dim{\rm Coker}(A)-(r'-\dim{\rm Ker}(A))\\
&&=\dim{\rm Ker}(A)-\dim{\rm Coker}(A)-r'\\
&&= r-r'.
\end{eqnarray*}

{\it Case 3}. {\bf ${\rm Ker}(A)\cap H\ne\{0\}$ and ${\rm Ker}(A)$
can not be contained in $H$.} We have direct sum decompositions
${\rm Ker}(A)={\rm Ker}(A)\cap H\oplus Z$, $H={\rm Ker}(A)\cap
H\oplus H_3$ and $X=H\oplus Z\oplus X_2$ and $Y=Y_1\oplus
A(X)=Y_1\oplus A(X_2)\oplus A(H)$. So $Y/A(H)\cong Y_1\oplus
A(X_2)$ and
\begin{eqnarray*}
{\rm Index}(A|_H)\!\!\!\!\!\!&&=\dim{\rm Ker}(A|_H)-\dim Y/{\rm Im}(A|_H)\\
&&= \dim{\rm Ker}(A|_H)-\dim (Y_1\oplus A(X_2))\\
&&= \dim{\rm Ker}(A)\cap H-\dim{\rm Coker}(A)-\dim A(X_2)\\
&&= \dim{\rm Ker}(A)\cap H-\dim{\rm Coker}(A)-\dim X_2\\
&&= \dim{\rm Ker}(A)\cap H-\dim{\rm Coker}(A)-(r'-\dim Z)\\
&&= \dim{\rm Ker}(A)\cap H +\dim Z -\dim{\rm Coker}(A)-r'\\
&&=\dim{\rm Ker}(A)- \dim{\rm Coker}(A)-r'\\
&&=r-r'.
\end{eqnarray*}
\hfill$\Box$\vspace{2mm}

Let $X$ be a Banach space and $M\subset X$ be a closed subspace.
${\rm codim}M:=\dim(X/M)$ is called codimension of $M$ (in $X$).
We have the following elementary fact, which is still proved for
completeness.

\begin{proposition}\label{prop:C.2}
If the subspace $M\subset X$ has codimension $k\in [2, +\infty]$
then for any natural number $l\le k$ there exists a closed subspace
$X_0\subset X$ of codimension $l$ such that $M\subset X_0$.
\end{proposition}

\noindent{\bf Proof.} Since $M\subset X$ is a proper closed
subspace there exists a nonzero $x_1\in X$ such that $d_1={\rm
dist}(x_1, M)>0$. By the Hahn-Banach theorem we have $f_1\in
X^\ast$ such that $f_1(x_1)=1$ and $f_1(M)=0$.  Let $M_1$ be the
subspace spanned by $M$ and $x_1$. It is a closed subspace of
codimension $k-1$. If $l=1$ or $k$ we may take $X_0={\rm
Ker}(f_1)$ or $M$. If $1<l<k-1$ We can continue to choose a
nonzero $x_2\in X$ such that $d_2={\rm dist}(x_2, M_1)>0$. Using
the Hahn-Banach theorem again we get a $f_2\in X^\ast$ such that
$f_2(x_2)=1$ and $f_2(M_1)=0$. Note that $f_1$ and $f_2$ are
linearly independent. Thus ${\rm Ker}(f_1)\cap{\rm Ker}(f_2)$ is a
closed subspace of codimension $2$ and contains $M$. If $l=2$ we
may take $X_0={\rm Ker}(f_1)\cap{\rm Ker}(f_2)$.      If $l>2$ we
may continue this procedure. After $l$ steps we get $l$ linearly
independent vectors $x_1,\cdots, x_l$ and functionals $f_i\in
X^\ast$, $i=1,\cdots, l$ such that $f_i(x_i)=1$ and
$f_i(M_{i-1})=0$, where $M_i$ is the subspace spanned by $x_i$ and
$M_{i-1}$, $i=1,\cdots, l$ and $M_0=M$. Clearly, $f_1,\cdots, f_l$
are linearly independent. Then $X_0=\cap^l_{i=1}{\rm Ker}(f_i)$ is
a desired closed subspace because any finite elements in
$\{f_1,\cdots, f_l\}$ are also linearly independent.
\hfill$\Box$\vspace{2mm}

\begin{proposition}\label{prop:C.3}
Let $G$ and $H$ be two Banach spaces, and  $A:G\to H$ be a bounded
linear Fredholm operator. Suppose that $H=H_1\oplus H_2$ is a direct
sum decomposition of Banach spaces, and $P_i:H\to H_i$ are natural
projections and $A_i=P_i\circ A$ for $i=1,2$. If $A_2: G\to H_2$ is
onto, then
$$\tilde A:{\rm Ker}(A_2)=A^{-1}(H_1)\to H_1,\; x\mapsto A(x)$$
is a bounded linear Fredholm operator and ${\rm Ind}(\tilde A)={\rm
Ind}(A)$.
\end{proposition}

\noindent{\bf Proof.} Denote by $A_i=P_i\circ A$, $i=1,2$. We first
prove that $\tilde A$ has a closed range. Let $\{Ax_n\}\subset{\rm
Im}(\tilde A)$ converge to $y\in H$. Since ${\rm Im}(A)$ is closed,
there is a $x\in G$ such that $Ax=y$. Note that $H_1$ is closed and
$Ax_n\in H_1$ for $n=1,2,\cdots$. We get $Ax=y\in H_1$, and thus
$x\in A^{-1}(H_1)$. Hence $y\in {\rm Im}(\tilde A)$. Namely, ${\rm
Im}(\tilde A)$ is closed.

Next it is clear that ${\rm Ker}(\tilde A)={\rm Ker}(A)$ since ${\rm
Ker}(A)\subset A^{-1}(H_1)$.

Finally, we prove ${\rm Coker}(\tilde A)={\rm Coker}(A)$. Since
${\rm Im}(\tilde A)\subset{\rm Im}(A)$ we have a well-defined linear
map
$$
\varphi: {\rm Coker}(\tilde A)=H_1/{\rm Im}(\tilde A)\to {\rm
Coker}(A)=H/{\rm Im}(A),\; [y]_0\mapsto [y],
$$
where $[y]_0=y+ {\rm Im}(\tilde A)$ and $[y]=y+ {\rm Im}(A)$. So it
suffices to prove that $\varphi$ is a linear isomorphism. Assume
that $\varphi([y]_0)=0$ for some $[y]_0\in {\rm Coker}(\tilde A)$.
Then $[y]=0$ and thus $y\in{\rm Im}(A)$, i.e., $y=Ax$ for some $x\in
G$. Since $y\in H_1$, we get $A_2x=0$ or $x\in{\rm
Ker}(A_2)=A^{-1}(H_1)$. Hence $y=Ax\in{\rm Im}(\tilde A)$, i.e.,
$[y]_0=0$. This shows that $\varphi$ is injective. To see that
$\varphi$ is also surjective, we take any $[y]\in {\rm Coker}(A)$.
Decompose $y$ into $y_1 + y_2$, where $y_i\in H_i$, $i=1,2$. Since
$A_2$ is onto, we have $y_2=A_2x$ for some $x\in G$. Set
$z=y_1-A_1x\in H_1$. Then
$z-y=y_1-A_1x-y_1-y_2=-A_1x-A_2x=-Ax\in{\rm Im}(A)$. Hence
$\varphi([z]_0)=[z]=[y]$. It follows that $\tilde A$ is Fredholm and
has index ${\rm Ind}(\tilde A)={\rm Ind}(A)$ because $\varphi$ is a
linear isomorphism. \hfill$\Box$\vspace{2mm}

\section{Appendix:\quad Existence of partially smooth cut-off functions in GW-theory}\label{app:D}
\setcounter{equation}{0}

 In this appendix we shall show how to
construct desired PS cut-off functions. To limit the length of this
paper we directly adopt the constructions in \S 2.4 of \cite{Lu3}.
Let $[{\bf f}]\in\overline{\mathcal M}_{g,m}(M, J, A)$ and ${\bf
f}=(f;\Sigma,\bar{\bf z})$ be a representative of it. Let $\bar{\bf
y}=\{y_1,\cdots, y_l\}$ be the set of all points added to $\Sigma$
for stabilizing $(\Sigma,\bar{\bf z})$. Then
$\sigma:=(\Sigma,\bar{\bf z}\cup\bar{\bf y})$ is a genus $g$ stable
curve with $m+l$ marked points. By (2.5) in \cite{Lu3} we have a
local deformation of $\sigma$:
$$
 V_{deform}(\sigma)\times
V_{resolve}(\sigma)\to\overline{\mathcal M}_{g,m+l},\;
(u,v)\mapsto [\Sigma_{(u,v)},\bar{\bf z}\cup\bar{\bf y} ],
$$
 where $V_{deform}(\sigma)\times V_{resolve}(\sigma)$ is a
small open neighborhood of the origin in some vector spaces of
finite dimension. In \S2.4.2 of \cite{Lu3},  for each $(u,v)\in
V_{deform}(\sigma)\times V_{resolve}(\sigma)$ we constructed a map
$f_{(u,v)}:\Sigma_{(u,v)}\to M$ satisfying
 $$
 {\rm diam}_\mu(f_{(u,v)}(\Sigma_{(u,v)}))\le {\rm
  diam}_\mu (f(\Sigma))+ 4|v|^{1/4}
  $$
with respect to a fixed Riemannian distance $d_M$ on $M$. By
Remark 2.12 in \cite{Lu3} we can choose a $\delta_{\bf f}>0$ so
small that the open ball $V_{\delta_{\bf f}}$ with center $0$ and
radius $\delta_{\bf f}$ is contained in $V_{deform}(\sigma)\times
V_{resolve}$, and that  each ${\bf f}_{(u,v)}=(f_{(u,v)},
\Sigma_{(u,v)},\bar{\bf
 z})$, $(u,v)\in V_{\delta_{\bf f}}$, is a stable $L^{k,p}$-map and has the homology class $A$ and the energy
 $E({\bf f}_{(u,v)})\le\omega(A)+1$ since $E({\bf f})=\omega(A)$.
For $\epsilon>0$ let $\widetilde{\bf U}_\epsilon(\delta_{\bf f})$
by the set of all tuples $(g_{(u,v)},\Sigma_{(u,v)},\bar{\bf
z}\cup\bar{\bf y})$ satisfying the following conditions: (i)
  $(u,v)\in V_{\delta_{\bf f}}$, (ii) $g_{(u,v)}:\Sigma_{(u,v)}\to M$ is a $L^{k,p}$-map,
  (iii) $\|g_{(u,v)}-f_{(u,v)}\|_{k,p}<\epsilon$, where the norm
 $\|\cdot\|_{k,p}$ is measured with respect to the metrics
 $\tau_{(u,v)}$ and $g_J=\mu$. For each $0<\delta\le\delta_{\bf f}$ we denote by
 $$
 \widetilde{\bf U}_{\delta,\epsilon}({\bf f}):=\{(g_{(u,v)},\Sigma_{(u,v)},\bar{\bf
z}\cup\bar{\bf y})\in\widetilde{\bf U}_\epsilon(\delta_{\bf
f})\,|\, (u,v)\in V_{\delta_{\bf f}},\,|(u,v)|<\delta\}
$$
and by
$$
{\bf U}_{\delta,\epsilon}({\bf
f}):=\{[g_{(u,v)},\Sigma_{(u,v)},\bar{\bf
z}]\,|\,(g_{(u,v)},\Sigma_{(u,v)},\bar{\bf z}\cup\bar{\bf
y})\in\widetilde{\bf U}_{\delta,\epsilon}({\bf
f})\}\subset{\mathcal B}^M_{g,m,A}.
$$
The latter is the desired neighborhood of $[{\bf f}]\in {\mathcal
B}^M_{g,m,A}$, and the $\widetilde{\bf U}_{\delta,\epsilon}({\bf
f})$ is the uniformizer of ${\bf U}_{\delta,\epsilon}({\bf f})$.

For $0<\delta_0\ll\delta$ and $0<\epsilon_0\ll\epsilon$ let us
choose two smooth cut-off functions
$$
\eta_1: V_{\delta_{\bf f}}\to [0,1]\quad{\rm and}\quad \eta_2:[0,
+\infty)\to [0,1]
$$
such that
$$
\eta_1((u,v))=\left\{\begin{array}{ll}
  1&\;{\rm if}\; |(u,v)|\le \delta_0,\\
  0&\;{\rm if}\; |(u,v)|\ge 2\delta_0,
  \end{array}\right.
\eta_2(t)=\left\{\begin{array}{ll}
  1&\;{\rm if}\; 0\le t\le\epsilon^p_0,\\
  0&\;{\rm if}\; t\ge 2\epsilon^p_0.
  \end{array}\right.
$$
Let us define $\eta:\widetilde{\bf U}_{\delta,\epsilon}({\bf
f})\to \R$ by
$$
\eta(g_{(u,v)})=\eta_1((u,v))\eta_2(\|g_{(u,v)}-f_{(u,v)}\|^p).
$$
It is not hard to check that $\eta$ is a partially smooth cut-off
function on $\widetilde{\bf U}_{\delta,\epsilon}({\bf f})$. By
shrinking $\delta_0>0$ and $\epsilon_0>0$ we can also require that
its support is contained in a given small neighborhood of ${\bf
f}$. As in \S 2.4.4 of \cite{Lu3} let
$$
\widetilde{\bf U}_{\delta,\epsilon}({\bf f}, {\rm H}):=\{
(g_{(u,v)},\Sigma_{(u,v)},\bar{\bf z}\cup\bar{\bf y})\in
\widetilde{\bf U}_{\delta,\epsilon}({\bf f})\,|\,
g_{(u,v)}(y_j)\in H_j\;\forall y_j\in\bar{\bf y}\}
$$
for $0<\delta\le\delta_{\bf f}$ and small $\epsilon>0$. Clearly, the
above $\eta$  restricts to a natural PS cut-off function on
$\widetilde{\bf U}_{\delta,\epsilon}({\bf f}, {\rm H})$. By the
obvious averaging process we can get a $Aut({\bf f})$-invariant PS
cut-off function on $\widetilde{\bf U}_{\delta,\epsilon}({\bf f},
{\rm H})$ from it.

\end{document}